\documentclass[english]{amsart}
\usepackage[T1]{fontenc}
\usepackage[utf8]{inputenc}
\usepackage{geometry}
\usepackage{babel}
\usepackage{mathrsfs}
\usepackage{amsmath}
\usepackage{amsthm}
\usepackage{amssymb}
\usepackage{graphicx}
\usepackage{bm}
\usepackage[all]{xy}
\usepackage{color}
\usepackage[usenames,dvipsnames,svgnames,table]{xcolor}
\usepackage{hyperref}
\usepackage{graphicx}
\usepackage{csquotes}
\usepackage[inline,shortlabels]{enumitem}
\usepackage[capitalize]{cleveref}
\usepackage{tikz-cd}
\usepackage{stmaryrd}

\tikzset{isom/.style={every to/.append style={edge node={node [sloped, allow upside down, auto=false, above]{$\sim$}}}}}
\tikzset{isomb/.style={every to/.append style={edge node={node [sloped, allow upside down, auto=false, below]{$\sim$}}}}}

\makeatletter
\theoremstyle{plain}
\newtheorem{thm}{\protect\theoremname}[section]
\theoremstyle{plain}
\newtheorem{lem}[thm]{\protect\lemmaname}
\theoremstyle{plain}
\newtheorem{cor}[thm]{\protect\corollaryname}
\theoremstyle{plain}
\newtheorem{conjecture}[thm]{\protect\conjecturename}
\theoremstyle{plain}
\newtheorem{prop}[thm]{\protect\propositionname}
\theoremstyle{definition}
\newtheorem{defn}[thm]{\protect\definitionname}
\newtheorem{setting}[thm]{Setting}
\theoremstyle{definition}
\newtheorem{exmpl}[thm]{\protect\examplename}

\theoremstyle{definition}
\newtheorem{rmk}[thm]{\protect\remarkname}
\newtheorem{rmks}[thm]{\protect\remarksname}

\newlist{thmenum}{enumerate}{1}
\setlist[thmenum]{label=(\roman*), ref=\thethm.(\roman*)}
\crefalias{thmenumi}{thm}

\newlist{lemenum}{enumerate}{1}
\setlist[lemenum]{label=(\roman*), ref=\thelem.(\roman*)}
\crefalias{lemenumi}{lem}

\newlist{propenum}{enumerate}{1}
\setlist[propenum]{label=(\roman*), ref=\theprop.(\roman*)}
\crefalias{propenumi}{prop}

\newlist{corenum}{enumerate}{1}
\setlist[corenum]{label=(\roman*), ref=\thecor.(\roman*)}
\crefalias{corenumi}{cor}

\newlist{defenum}{enumerate}{1}
\setlist[defenum]{label=(\alph*), ref=\thedefn.(\alph*)}
\crefalias{defenumi}{defn}

\newlist{exmplsenum}{enumerate}{1}
\setlist[exmplsenum]{label=(\alph*), ref=\theexmpls.(\alph*)}
\crefalias{exmplsenumi}{exmpl}

\newlist{rmksenum}{enumerate}{1}
\setlist[rmksenum]{label=(\roman*), ref=\thermks.(\roman*)}
\crefalias{rmksenumi}{rmk}

\usepackage{xy}
\usepackage{mathrsfs}
\usepackage{color}

\usepackage{xcolor}
\hypersetup{colorlinks=true,linkcolor=purple,citecolor=magenta}

\numberwithin{thm}{subsection}

\DeclareMathOperator{\Ind}{Ind}
\DeclareMathOperator{\Pro}{Pro}
\DeclareMathOperator{\Hom}{Hom}
\DeclareMathOperator{\Ext}{Ext}
\DeclareMathOperator{\End}{End}
\DeclareMathOperator{\Map}{Map}
\DeclareMathOperator{\Rep}{Rep}
\DeclareMathOperator{\Isom}{Isom}

\DeclareMathOperator{\RHom}{RHom}
\DeclareMathOperator{\intHom}{\underline\Hom}
\DeclareMathOperator{\intMap}{\underline{\smash\Map}}
\DeclareMathOperator{\intRep}{\underline{\smash\Rep}}
\DeclareMathOperator{\intIsom}{\underline{\Isom}}
\DeclareMathOperator{\Fun}{Fun}
\DeclareMathOperator{\Ani}{Ani}
\DeclareMathOperator{\Adm}{Adm}
\DeclareMathOperator{\Eis}{Eis}
\DeclareMathOperator{\CT}{CT}
\DeclareMathOperator{\GL}{GL}
\DeclareMathOperator{\SL}{SL}
\newcommand{\et}{{\mathrm{et}}}
\newcommand{\nuc}{{\mathrm{nuc}}}
\DeclareMathOperator{\Div}{Div}
\DeclareMathOperator{\Sat}{Sat}
\newcommand{\spec}{{\mathrm{spec}}}
\newcommand{\coarse}{{\mathrm{crs}}}
\newcommand{\oc}{{\mathrm{oc}}}
\newcommand{\sm}{{\mathrm{sm}}}
\newcommand{\disc}{{\mathrm{disc}}}
\DeclareMathOperator{\cusp}{cusp}
\DeclareMathOperator{\cts}{cts}
\DeclareMathOperator{\Coh}{Coh}
\DeclareMathOperator{\Perf}{Perf}
\DeclareMathOperator{\QCoh}{QCoh}
\DeclareMathOperator{\IndCoh}{IndCoh}
\DeclareMathOperator{\PQCoh}{PQCoh}
\DeclareMathOperator{\Bun}{Bun}
\DeclareMathOperator{\Par}{Par}
\DeclareMathOperator{\Spec}{Spec}
\newcommand{\fin}{\mathrm{fin}}
\DeclareMathOperator{\adm}{adm}
\DeclareMathOperator{\ren}{ren}
\newcommand{\mot}{\mathrm{mot}}
\newcommand{\eff}{\mathrm{eff}}
\newcommand{\restr}{\mathrm{restr}}
\newcommand{\bdd}{{\mathrm{bdd}}}
\DeclareMathOperator{\Qellbar}{\overline{\mathbf{Q}}_{\ell}}
\DeclareMathOperator{\D}{\mathcal{D}}
\DeclareMathOperator{\Dadm}{\mathbf D_{adm}}
\DeclareMathOperator{\Dgs}{\mathbf D_{GS}}
\DeclareMathOperator{\Dtwgs}{\mathbf D_{tw.GS}}
\DeclareMathOperator{\Dtwadm}{\mathbf D_{tw.adm}}
\DeclareMathOperator{\Dbz}{\mathbf D_{BZ}}
\DeclareMathOperator{\Dverd}{\mathbf D_{Verd}}

\DeclareMathOperator{\FF}{FF}

\DeclareMathOperator{\LocSys}{LocSys}
\DeclareMathOperator{\Perfd}{Perfd}

\newcommand{\op}{{\mathrm{op}}}
\newcommand{\co}{{\mathrm{co}}}
\newcommand{\solid}{{\scalebox{0.5}{$\square$}}}

\newcommand{\noloc}{\nobreak\mskip6mu plus1mu\mathpunct{}\nonscript\mkern-\thinmuskip{:}\mskip2mu\relax} 

\DeclareMathOperator{\id}{id}
\DeclareMathOperator{\Fin}{Fin}
\DeclareMathOperator{\Ring}{Ring}

\DeclareMathOperator{\Op}{Op}
\DeclareMathOperator{\Alg}{Alg}
\DeclareMathOperator{\CAlg}{CAlg}
\newcommand{\rel}{\mathrm{rel}}

\DeclareMathOperator{\Cond}{Cond}
\DeclareMathOperator{\Mod}{Mod}

\DeclareMathOperator{\Shv}{Shv}

\DeclareMathOperator{\Stk}{Stk}
\DeclareMathOperator{\Sch}{Sch}
\DeclareMathOperator{\fib}{fib}
\DeclareMathOperator{\cofib}{cofib}

\DeclareMathOperator{\Spa}{Spa}
\DeclareMathOperator{\Spd}{Spd}
\DeclareMathOperator{\Corr}{Corr}
\DeclareMathOperator{\Ar}{Ar}
\DeclareMathOperator{\Conv}{Conv}
\DeclareMathOperator{\Cat}{Cat}

\DeclareMathOperator{\Der}{Der}
\DeclareMathOperator{\Top}{Top}
\DeclareMathOperator{\Ab}{Ab}
\DeclareMathOperator{\Mon}{Mon}
\DeclareMathOperator{\Grp}{Grp}
\DeclareMathOperator{\CMon}{CMon}
\DeclareMathOperator{\Enr}{Enr}
\DeclareMathOperator{\Ass}{Ass}
\DeclareMathOperator{\Sym}{Sym}
\DeclareMathOperator{\Lie}{Lie}

\DeclareMathOperator{\locHck}{\mathcal{H}ck}
\DeclareMathOperator{\glbHck}{Hck}
\DeclareMathOperator{\Fil}{Fil}
\DeclareMathOperator{\Gr}{Gr}
\DeclareMathOperator{\Pair}{Pair}
\DeclareMathOperator{\Ideals}{Ideals}
\DeclareMathOperator{\ev}{ev}
\DeclareMathOperator{\ins}{ins}
\DeclareMathOperator{\Nm}{Nm}
\DeclareMathOperator{\Tw}{Tw}
\DeclareMathOperator{\vStk}{vStk}
\DeclareMathOperator{\PrL}{Pr^L}
\DeclareMathOperator{\Sp}{Sp}
\DeclareMathOperator{\PrR}{Pr^R}
\DeclareMathOperator{\dimtrg}{dimtrg}
\DeclareMathOperator{\DAC}{DAC}
\newcommand{\tensor}{\otimes}

\newcommand{\aff}{\mathrm{aff}}

\newcommand{\dbl}{\mathrm{dbl}}
\newcommand{\lis}{\mathrm{lis}}
\newcommand{\rev}{\mathrm{rev}}
\newcommand{\lax}{\mathrm{lax}}

\DeclareMathOperator{\Gal}{Gal}
\DeclareMathOperator{\Ad}{Ad}

\DeclareMathOperator{\IPCoh}{IPCoh}
\DeclareMathOperator{\cind}{c-ind}

\DeclareMathOperator{\ICoh}{ICoh}
\DeclareMathOperator{\PCoh}{PCoh}

\DeclareMathOperator{\FilgeRing}{Fil^{\ge0}Ring}
\DeclareMathOperator{\RingMod}{RingMod}
\DeclareMathOperator{\Rees}{Rees}

\DeclareMathOperator{\pr}{pr}
\DeclareMathOperator{\orb}{orb}

\newcommand{\C}{\mathbf C}
\newcommand{\R}{\mathbf R}
\newcommand{\Q}{\mathbf Q}
\newcommand{\Z}{\mathbf Z}
\newcommand{\N}{\mathbf N}
\newcommand{\F}{\mathbf F}
\newcommand{\calO}{\mathcal O}
\newcommand{\pitch}{\mathbin{\rotatebox[origin=c]{180}{\bm{$\pitchfork$}}}}

\newcommand{\cat}[1]{\mathcal{#1}}

\newcommand{\isom}{\cong}
\newcommand{\from}{\mathrel{\leftarrow}}
\newcommand{\xto}[1]{\mathbin{\xrightarrow{#1}}} 
\newcommand{\xfrom}[1]{\mathbin{\xleftarrow{#1}}} 
\newcommand{\isoto}{\xto{\smash{\raisebox{-.4ex}{\ensuremath{\scriptstyle\sim}}}}}
\newcommand{\injto}{\mathrel{\hookrightarrow}}
\newcommand{\injfrom}{\mathrel{\hookleftarrow}}
\newcommand{\surjto}{\mathrel{\twoheadrightarrow}}

\newcommand{\pc}{\mathrm{pc}}
\newcommand{\qc}{\mathrm{qc}}
\newcommand{\qcqs}{\mathrm{qcqs}}
\newcommand{\abs}[1]{|#1|}
\newcommand{\cl}{\mathrm{cl}}
\newcommand{\red}{\mathrm{red}}
\newcommand{\afp}{\mathrm{afp}}
\newcommand{\lafp}{\mathrm{lafp}}
\newcommand{\comp}{\mathbin{\circ}}
\newcommand{\fpj}{{\mathrm{fpj}}}

\DeclareMathOperator{\SD}{SD}

\makeatother

\providecommand{\conjecturename}{Conjecture}
\providecommand{\corollaryname}{Corollary}
\providecommand{\lemmaname}{Lemma}
\providecommand{\propositionname}{Proposition}
\providecommand{\theoremname}{Theorem}
\providecommand{\definitionname}{Definition}
\providecommand{\examplename}{Example}
\providecommand{\examplesname}{Examples}
\providecommand{\remarkname}{Remark}
\providecommand{\remarksname}{Remarks}

\begin{document}
\title{The categorical local Langlands conjecture}
\author{David Hansen and Lucas Mann}
\begin{abstract}We formulate a program to prove the categorical local Langlands conjecture (CLLC) of Fargues-Scholze, for all quasisplit $p$-adic groups where the Fargues-Scholze $L$-parameters agree with the semisimplification of a known ``automorphic'' local Langlands parametrization. A key working hypothesis - which we expect to prove elsewhere jointly with Hamann - is the compatibility of the \emph{enhanced Whittaker coefficient functor} $c_\psi$ with Eisenstein series. For $\mathrm{GL}_n$, we show that this hypothesis alone implies the full CLLC. For more general groups $G$, we prove an induction principle which reduces CLLC for $G$ to CLLC for all proper Levi subgroups together with a very small amount of information about $G$. This principle applies unconditionally to many classical groups with current technology.

Along the way, we establish many foundational results. In particular:
\begin{itemize}
    \item We prove a very strong finiteness theorem for spectral constant term functors.
    \item We prove a spectral analogue of Bernstein's finite global dimension theorem for $p$-adic Hecke algebras.
    \item We introduce and develop the theory of \emph{admissible} ind-coherent sheaves and admissible duality on derived stacks.
    \item We prove a duality theorem for the spectral action.
\end{itemize}
Using all of these results, we unconditionally define a new and explicit functor $t_{\psi}$ from the spectral side to the automorphic side, which is defined on enough ind-coherent sheaves to control the entire conjecture.

\end{abstract}

\maketitle
\tableofcontents{}
\section{Introduction}

\subsection{Main results}

Fix a finite extension $F/\Q_p$ and a quasisplit connected reductive group $G$ over $F$. The \emph{local Langlands conjecture}, in its broadest form, seeks to parametrize the irreducible smooth representations of $G(F)$ in Galois-theoretic terms. A huge amount of effort has been expended on making this conjecture precise, and proving it in many cases. The work of Fargues--Scholze \cite{FS} puts the conjecture on new grounds by categorifying both sides of the sought-after correspondence.
\begin{itemize}
    \item On the side of smooth representations \cite{FS} introduce the category $\D(\Bun_G) = \D_\lis(\Bun_G, \Qellbar)$, which roughly parametrizes families of smooth representations for all inner forms of Levis of $G$.

    \item On the side of Galois representations one considers the category $\IndCoh(\Par_G)$ of ind-coherent sheaves on the algebraic stack $\Par_G$ of $L$-parameters. An $L$-parameter for $G$ is a continuous 1-coycle $W_F \to \hat G(\Qellbar)$, where $W_F$ denotes the Weil group of $F$ and $\hat G$ denotes the dual group of $G$ over $\Qellbar$ (equipped with its natural action of $W_F$ by pinned automorphisms).
\end{itemize}
In \cite[Conjecture~X.1.4]{FS} Fargues--Scholze conjecture the existence of a natural equivalence of categories
\begin{align*}
     \mathbf L_\psi\colon \D(\Bun_G) \isoto \IndCoh(\Par_G)
\end{align*}
that only depends on a choice of Whittaker datum $(B, \psi)$ for $G$ and enjoys various compatibilities with the basic structures on both sides. This hoped-for equivalence has come to be known as the \emph{categorical local Langlands conjecture}. Following \cite{Beijing}, we observe that the constructions in \cite{FS} allow one to pin down the equivalence uniquely, making the above conjecture precise (see \cref{prop:CLLCuniqueintro} and \cref{conj:CLLCintro} for a precise statement). The vague intuition for the equivalence is that an irreducible smooth $G(F)$-representation (viewed as an element in $\D(\Bun_G)$ by $!$-extension) corresponds to the skyscraper sheaf on $\Par_G$ concentrated at the corresponding $L$-parameter.

In this paper we study this conjecture from first principles, and develop a program to prove it for many groups. In particular, we prove the following theorem.

\begin{thm} \label{thm:intro-main-result}
Let $G = \GL_n$ with fixed Whittaker datum $(B, \psi)$. Assume that:
\begin{itemize}
    \item[(Eis)] For every standard parabolic subgroup $P$ of every standard Levi of $G$, the functor of enhanced Whittaker coefficient is compatible with Eisenstein functors along $P$ (see \cref{conj:Eiscpsiintro} for a precise statement).
\end{itemize}
Then the categorical local Langlands conjecture holds for $\GL_n$, i.e. there is a natural equivalence
\begin{align*}
    \mathbf L_\psi\colon \D(\Bun_{\GL_n}) \isoto \IndCoh(\Par_{\GL_n}).
\end{align*}
\end{thm}

Our proof makes crucial use of the known local Langlands correspondence for $\GL_n$ \cite{HarrisTaylor, Henniart}. The methods of proof also yield similar results for other classical groups $G$ for which the classical conjecture is known (see \cref{sec:intro.classical-groups} below). The Eisenstein compatibility (Eis) is joint work in progress with Linus Hamann and was announced recently for $\GL_2$ by Hamann-Imai.

While we draw some key inspirations from the recent proof of the geometric Langlands conjecture (in the setting of complex geometry), the categorical local Langlands conjecture has its own unique difficulties, and we are forced to introduce some essential new ideas. The arguments here also build on our earlier papers \cite{HanBasic,HHS,HKW,HJnote,heyer-mann-6ff,mann-p-adic-6-functors,mann-nuclear}.

In \cref{sec:summary} below we give a brief overview of the key ideas which go into the proof of \cref{thm:intro-main-result}, which should serve the reader as a rapid entry point into the main notations, constructions and results of this paper. We make no attempt to give complete proof sketches or thorough references in this overview. In the remainder of the introduction, we give a much more detailed discussion of all our results; the structure of the rest of the introduction is explained at the end of \Cref{sec:summary}. 


\subsection{Summary of the proof} \label{sec:summary}

In the following we sketch the ideas of the proof of Theorem \ref{thm:intro-main-result} and give a survey of the necessary tools that go into it. A much more elaborate discussion of all these ideas and many other results is given in the following subsections of the introduction.

The proof can be viewed as proceeding in six steps: The first five steps construct various related functors between different incarnations of both sides of the correspondence, and the last step concludes the desired equivalence. Throughout, we fix a quasisplit group $G$ and a Whittaker datum $(B, \psi)$.

\subsubsection{Step 1: The functor \texorpdfstring{$a_\psi$}{a_psi}}

One of the main results of Fargues--Scholze is the construction of a natural action
\begin{align*}
    \QCoh(\Par_G) \times \D(\Bun_G) \to \D(\Bun_G), \qquad (\mathcal F, A) \mapsto \mathcal F * A,
\end{align*}
called the \emph{spectral action}. Here $\QCoh(-)$ denotes the (derived) category of quasi-coherent sheaves on the algebraic stack $\Par_G$.

We can use the spectral action to define a functor $a_\psi$ as follows. Recall that there is an open immersion $i_!\colon */G(F) \injto \Bun_G$ and $!$-push along this maps embeds smooth $G(F)$-representations fully faithfully into $\D(\Bun_G)$. We furthermore denote by $W_\psi := \cind_{U(F)}^{G(F)} \psi$ the Whittaker representation associated to the datum $(B,\psi)$.

\begin{defn}
We define
\begin{align*}
    a_\psi\colon \QCoh(\Par_G) \to \D(\Bun_G), \qquad \mathcal F \mapsto a_\psi(\mathcal F) := \mathcal F * i_{1!} W_\psi.
\end{align*}
\end{defn}

By construction, $a_\psi$ is a colimit-preserving $\QCoh(\Par_G)$-linear functor. We call $G$ ($\psi$-)\emph{reasonable} if $a_\psi$ preserves compact objects, i.e.\ sends $\Perf^\qc(\Par_G)$ to $\D(\Bun_G)^\omega$, where $\Perf^\qc(\Par_G) \subseteq \QCoh(\Par_G)$ denotes the full subcategory spanned by the perfect complexes with quasicompact support. In \cref{prop:apsifinitebasic} we characterize reasonable groups and show in particular that this condition holds if the Fargues--Scholze $L$-parameter map has finite fibers; this is known for many groups, e.g. $\GL_n$.

\subsubsection{Step 2: The functor \texorpdfstring{$c_\psi$}{c_psi}}

We now pass to the right adjoint of $a_\psi$.

\begin{defn}
Let
\begin{align*}
    c_\psi\colon \D(\Bun_G) \to \QCoh(\Par_G)
\end{align*}
be the right adjoint of $a_\psi$.
\end{defn}

By abstract nonsense we see that if $G$ is reasonable then $c_\psi$ preserves colimits, and it is also not hard to see that it is then $\QCoh(\Par_G)$-linear.

It is hard to access the functor $c_\psi$ in general, and to make any serious analysis of it we need to impose some assumptions on $G$. We call $G$ \emph{well-understood} if there is a known local Langlands correspondence for $G$ and its standard Levis which is compatible with the Fargues--Scholze construction of $L$-parameters, see \cref{def:wellunderstoodG} for the precise condition. Note that this hypothesis implies that $G$ is reasonable (see \cref{rslt:well-understood-implies-reasonable}). In \cref{rslt:c-psi-sends-compact-to-coherent} we prove:

\begin{thm} \label{rslt:intro-c-psi-sends-compact-to-coherent}
Assume that $G$ is well-understood and satisfies the Eisenstein compatibility (Eis) from \cref{thm:intro-main-result}. Then $c_\psi$ sends compact sheaves to coherent sheaves with quasicompact support.
\end{thm}

We sketch the idea, so fix $A \in \D(\Bun_G)^\omega$. If $c_\psi(A)$ is supported on the supercuspidal locus of $\Par_G$ then the claim is easy (and unconditional on $G$). The hard part is to handle the case that $c_\psi(A)$ is supported away from the supercuspidal locus. In this case, we show that $A$ can be built via Eisenstein functors from compact sheaves $A_0 \in \D(\Bun_M)^\omega$ for various proper Levis $M$; by induction on $M$ and (Eis), we can conclude. The fact that $A$ can be built as claimed is highly non-trivial and uses the main results from \cite{HKW} (which requires that $G$ is well-understood) and \cite{HHS}.

\subsubsection{Step 3: The functor \texorpdfstring{$\mathbf L_\psi$}{L_psi}} With \cref{rslt:intro-c-psi-sends-compact-to-coherent} in hand we can \emph{renormalize} $c_\psi$ in the following manner.

\begin{defn} \label{defn:Langlandsfunctorintro}
Assume that $G$ is well-understood and satisfies (Eis). We define
\begin{align*}
    \mathbf L_\psi\colon \D(\Bun_G) \to \IndCoh(\Par_G)
\end{align*}
to be the ind-extension of the functor $c_\psi\colon \D(\Bun_G)^\omega \to \Coh^\qc(\Par_G)$. Here $\Coh^\qc(\Par_G) \subseteq \QCoh(\Par_G)$ denotes the full subcategory spanned by the coherent sheaves with quasicompact support.
\end{defn}

It follows immediately from the construction and from the properties of $c_\psi$ that $\mathbf L_\psi$ is $\QCoh(\Par_G)$-linear, preserves colimits and compact objects, and commutes with Eisenstein functors in the following sense: For every standard parabolic $P \subseteq G$ with opposite $\overline P$ and Levi $M$, there is a commuting square
\begin{equation*}\begin{tikzcd}
    \D(\Bun_M) \arrow[r,"\mathbf L_{\psi_M}"] \arrow[d,"\Eis_P"] & \IndCoh(\Par_M) \arrow[d,"\Eis_{\overline P,!}^\spec"] \\
    \D(\Bun_G) \arrow[r,"\mathbf L_\psi"] & \IndCoh(\Par_G)
\end{tikzcd}\end{equation*}
Here $\Eis_{\overline P,!}$ and $\Eis_P^\spec$ are defined using pull-push to $\Bun_{\overline P}$ and $\Par_P$, respectively. We refer the reader to \cref{rslt:intro-properties-of-L-psi} for details.

\subsubsection{Step 4: The functor \texorpdfstring{$\mathbf R_\psi$}{R_psi}}

Our goal is to prove that $\mathbf L_\psi$ is an equivalence, for which we look at its right adjoint:

\begin{defn}
Assume that $G$ is well-understood and satisfies (Eis). Define
\begin{align*}
    \mathbf R_\psi\colon \IndCoh(\Par_G) \to \D(\Bun_G)
\end{align*}
to be the right adjoint of $\mathbf L_\psi$.
\end{defn}

Since $\mathbf L_\psi$ preserves compact objects, it follows immediately that $\mathbf R_\psi$ preserves all small colimits. The next result, proved in \cref{thm:rightadjointcompactness}, gives us more control on $\mathbf R_\psi$:

\begin{thm} \label{rslt:intro-Rpsi-preserves-compact-objects}
Assume that $G$ is well-understood and satisfies (Eis). Then $\mathbf R_\psi$ preserves compact objects.
\end{thm}

This is one of the hardest results in this paper and requires us to develop a considerable amount of technical tools on the spectral side. Roughly, by analyzing the structure of compact objects in $\D(\Bun_G)$ and playing around with adjunctions and induction over Levi subgroups, one reduces the claim to certain key properties of coherent sheaves on $\Par_G$, which we prove unconditionally. In particular, we show the following difficult theorem.

\begin{thm}[\cref{rslt:miracle-bound}] \label{rslt:intro-miracle-bound}
For any $\mathcal F, \mathcal G \in \Coh^\qc(\Par_G)$, the complex $\RHom(\mathcal F, \mathcal G) \in \D(\Qellbar)$ is bounded and all of its cohomology groups are finitely generated $\mathfrak Z_G^\spec$-modules.
\end{thm}

Here $\mathfrak Z_G^\spec$ denotes the ring of global functions on the stack $\Par_G$.

Let us sketch the idea of the proof. The claim is easy if $\mathcal F$ is perfect. Now by a result of Takaya, any $\mathcal F \in \Coh^\qc(\Par_G)$ can be built from sheaves of the form $\Eis_P^\spec(\mathcal F_0)$ for parabolics $P = MU \subseteq G$ and sheaves $\mathcal F_0 \in \Perf^\qc(\Par_M)$. This reduces the claim to showing that $\RHom(\Eis^\spec_P \mathcal F_0, \mathcal G) = \RHom(\mathcal F_0, \CT_P^\spec (\mathcal G))$ has the claimed finiteness properties, where $\CT_P^\spec$ is the constant term functor (right adjoint to $\Eis_P^\spec$). This then furthermore reduces to showing that $\CT_P^\spec$ preserves coherent sheaves, at least on graded pieces for the central grading. This is a weak form of Second Adjointness on $\Par_G$, which we prove by establishing a subtle version of the geometric lemma for ind-coherent sheaves on $\Par_G$ (see \cref{rslt:spectral-geometric-lemma}).

\subsubsection{Step 5: The functor \texorpdfstring{$t_\psi$}{t_psi}}

In order to show that $\mathbf L_\psi$ and $\mathbf R_\psi$ are inverse to each other, we need even more control on $\mathbf R_\psi$. A key insight in this paper, and the starting point for the whole project, is an \emph{alternative} construction of $\mathbf R_\psi$. It relies on the expected compatibility of the Langlands correspondence with dualities. Let us briefly recall the most important finiteness conditions and dualities in our context:
\begin{enumerate}[(a)]
    \item The compact objects in $\D(\Bun_G)$ together with Bernstein-Zelevinsky duality $\Dbz$.
    \item The ULA sheaves in $\D(\Bun_G)$ together with Verdier duality $\Dverd = \intHom(-, \Qellbar)$.
    \item The coherent sheaves (with quasicompact support) in $\IndCoh(\Par_G)$ together with Grothendieck-Serre duality $\Dgs = \intHom(-, \omega_{\Par_G})$.
\end{enumerate}
By construction, if $\mathbf L_\psi$ is an equivalence then it matches the objects in (a) and (c). It is furthermore expected that it is also compatible with the respective dualities, up to a twist by the Chevalley involution. It is a natural question whether there is a finiteness condition together with duality on $\IndCoh(\Par_G)$ that matches with (b). To approach this question we recall from \cite[\S 4.4]{heyer-mann-6ff} that in every 6-functor formalism there are two natural finiteness conditions, both coming with duality: suaveness and primness. Now (a) and (c) correspond to prim objects and prim duality in their respective 6-functor formalism, while (b) corresponds to suave objects and suave duality. This motivates:

\begin{defn}
A sheaf $\mathcal F \in \IndCoh(\Par_G)$ is called \emph{admissible} if it is suave in the $\IndCoh$-formalism. We call the associated duality \emph{admissible duality} and denote it by $\Dadm$.
\end{defn}

Admissible sheaves are discussed in \cref{sec:adm-sheaves}, where in particular it is shown that $\mathcal F$ is admissible if and only if $\RHom(\mathcal G, \mathcal F) \in \D(\Qellbar)$ is perfect for all $\mathcal G \in \Coh^\qc(\Par_G)$. Admissibility is a somewhat strange condition on ind-coherent sheaves, in particular there are lots of examples of very \enquote{big} admissible sheaves. It turns out that, assuming the Langlands correspondence holds, then sheaves in (b) will match with admissible sheaves on $\Par_G$ and this matching is compatible with duality (see \cref{prop:Admdualityintro}).

In \cref{thm:spectraltemperizationintro} we make a surprising observation: If a sheaf $\mathcal F \in \IndCoh(\Par_G)$ is coherent and admissible, then $\Dadm(\mathcal F)$ lies in the full subcategory $\QCoh(\Par_G) \subseteq \IndCoh(\Par_G)$. This allows the following definition:

\begin{defn}\label{def:cohfinhardapproach}
Denote by $\Coh(\Par_G)_\fin \subseteq \IndCoh(\Par_G)$ the full subcategory spanned by the sheaves that are both coherent (with quasicompact support) and admissible. We define
\begin{align*}
    t_\psi\colon \Coh(\Par_G)_\fin \to \D(\Bun_G), \qquad \mathcal F \mapsto t_\psi(\mathcal F) := \Dverd a_{\psi^{-1}} c^* \Dadm(\mathcal F),
\end{align*}
where $c^*$ denotes the involution on $\Par_G$ induced by twisting an $L$-parameter by the Chevalley involution of $\phantom{}^L G$.
\end{defn}

By construction $a_\psi$ is left adjoint to $\mathbf L_\psi$ on $\QCoh(\Par_G)$. By some careful manipulations using \cref{thm:spectraltemperizationintro} together with duality properties of the spectral action (which we prove in this paper), one deduces that $t_\psi$ is \emph{right} adjoint to $\mathbf L_\psi$ on $\Coh(\Par_G)_\fin$ -- the key observation is that because the dualities are contravariant, they swap the directions of certain adjunctions. In particular, we get the following result.

\begin{prop} \label{rslt:intro-relation-of-t-psi-and-R-psi}
Suppose that $G$ is well-understood and satisfies (Eis), so that $\mathbf R_\psi$ is defined. Then
\begin{align*}
    \mathbf R_\psi|_{\Coh(\Par_G)_\fin} = t_\psi.
\end{align*}
\end{prop}
The great advantage of this description is that $t_\psi$ is built out of very explicit functors over which we have a lot of control.

\subsubsection{Step 6: Conclusion}

So far we have constructed the adjunction
\begin{align*}
    \mathbf L_\psi\colon \D(\Bun_G) \rightleftarrows \IndCoh(\Par_G) \noloc \mathbf R_\psi
\end{align*}
and we wish to prove that these functors are inverse to each other. By abstract nonsense, this is equivalent to showing that $\mathbf R_\psi$ is fully faithful and $\mathbf L_\psi$ is conservative. Let us discuss the full faithfulness of $\mathbf R_\psi$. With a bit of work, one can deduce from the definition of $t_\psi$ and from \cref{rslt:intro-relation-of-t-psi-and-R-psi} that
\begin{align*}
    \text{$a_\psi$ is fully faithful} \implies \text{$t_\psi$ is fully faithful} \implies \text{$\mathbf R_\psi$ is fully faithful}.
\end{align*}
These implications are the main reason we introduced $t_\psi$! Now, in order to show that $a_\psi$ is fully faithful, we use the $\QCoh(\Par_G)$-linearity of $a_\psi$ and $c_\psi$ to observe
\begin{align*}
    \text{$a_\psi$ is fully faithful} &\iff \text{$\calO_{\Par_G} \to c_\psi a_\psi \calO_{\Par_G}$ is an isomorphism}\\
    &\iff \text{$c_\psi i_{1!} W_\psi$ is a line bundle}.
\end{align*}
We have thus reduced the full faithfulness of $\mathbf R_\psi$ to a simple property of the functor $c_\psi i_{1!}$. Finally, using a subtle induction over Levis and compatibility with Eisenstein functors, we also manage to reduce the conservativity of $\mathbf L_\psi$ to a property of $c_\psi i_{1!}$, leading to the following key result.

\begin{thm}[\cref{thm:DreamInductionOnLevisTheorem}] \label{thm:intro-dream-induction-result}
Assume that $G$ is well-understood and satisfies (Eis). Assume furthermore:
\begin{enumerate}[(1)]
    \item For every proper standard Levi $M \subsetneq G$, the functor $\mathbf L_{\psi_M}$ is an equivalence of categories.

    \item The functor $\mathbf L_\psi i_{1!}\colon \D(G(F), \Qellbar) \to \IndCoh(\Par_G)$ sends $W_\psi$ to a line bundle, and is conservative on compact supercuspidal representations and sends them to \emph{cuspidal} coherent sheaves in the sense of Definition \ref{defn:basiccoherent}.(2).
\end{enumerate}
Then $\mathbf L_\psi$ is an equivalence.
\end{thm}

We have thus reduced the desired equivalence to an understanding of the functor $\mathbf{L}_\psi i_{1!}$. Now for $G = \GL_n$, Ben-Zvi--Chen--Helm--Nadler \cite{BCHN} already constructed a fully faithful functor
\begin{align*}
    \mathscr{S}_{\GL_n}\colon \D(\GL_n(F), \Qellbar) \injto \IndCoh(\Par_{\GL_n}),
\end{align*}
over which they have quite some control. Using Eisenstein compatibility and some basic algebra, we are able to prove that $\mathbf L_\psi i_{1!} = \mathscr{S}_{\GL_n}$, which in turn implies that this functor satisfies condition (2) in \cref{thm:intro-dream-induction-result}. This proves \cref{thm:intro-main-result}.

\subsubsection{Structure of the introduction}

In the remainder of this introduction, we give a much more leisurely and thorough discussion of all our results, structured roughly as follows.
\begin{itemize}
    \item In \cref{sec:history} we provide some historical background on the development of the local Langlands program and put our results in this broad context. 

    \item In \Cref{sec:intro-FS} we fix important notation and summarize the main constructions from the Fargues--Scholze program \cite{FS}. We also draw some first conclusions from this machinery and introduce the functors $a_\psi$, $c_\psi$ and $\mathbf L_\psi$. In particular we carefully formulate the precise version of the conjecture and explain why Definition \ref{defn:Langlandsfunctorintro} is the \emph{only} option for the correct functor.

    \item In \Cref{sec:firstproperties} we explore some first consequences of the conjecture, especially its compatibility with finiteness conditions and dualities. In particular, we state an important new duality theorem for the spectral action. We also explain how we arrived at the formula for $t_{\psi}$ given above. 

    \item In \Cref{sec:wellunderstood} we give a more detailed account of the well-understood condition and its consequences, roughly corresponding to (most of) Steps 4-6 above. 

    \item In \Cref{sec:GLnintro} we explain how the results in \cite{BCHN} allow us to verify the conditions of \cref{thm:intro-dream-induction-result} for $\GL_n$, completing the proof of \cref{thm:intro-main-result}.

    \item In \Cref{sec:intro.classical-groups} we discuss how to apply \cref{thm:intro-dream-induction-result} to various classical groups. This is much more subtle than the case of general linear groups.

    \item In \Cref{ss:spectralintro} we give a detailed discussion of \cref{rslt:intro-miracle-bound} and its application in Step 4 above.

    \item In \Cref{sec:introindcoherent} we survey the long \Cref{sec:alggeo}, which gives a detailed development, essentially from scratch, of quasicoherent and ind-coherent sheaves on derived algebraic stacks. 

    \item In \Cref{sec:applications} we mention (without proof) some first applications of the categorical local Langlands conjecture.

    \item In \Cref{sec:otherpeople} we carefully analyze the similarities and differences between our arguments and the recent proof of the geometric Langlands conjecture. We also discuss some other recent works on categorical local Langlands and their connection with this paper.

    \item Finally, in \Cref{sec:future} we highlight a few things which are \emph{not} done in this paper, and mention some possible next steps.
\end{itemize}

\subsection{Historical background} \label{sec:history}

As mentioned above, the local Langlands conjecture seeks to parametrize the irreducible smooth representations of $G(F)$ in Galois-theoretic terms. In its modern form, the conjecture proposes that irreducible smooth representations of $G(F)$ should be parametrized by (suitable) pairs $(\phi,\rho)$ where $\phi\colon W_{F} \times \SL_{2} \to \phantom{}^{L}G$ is an L-parameter and $\rho$ is an irreducible algebraic representation of a group closely related to the centralizer group $S_{\phi}=\mathrm{Cent}_{\hat{G}}(\phi)$. Moreover, this parametrization should depend only on a choice of Whittaker datum for $G$, and should satisfy many good properties. We refer to \cite{KalethaNQS} for a precise formulation of this conjecture in modern terms.

The local Langlands conjecture was originally motivated by the needs of the global Langlands program. For $\mathbf{G}$ a connected reductive group over a number field $K$, fundamental conjectures of Arthur \cite{ArthurUnipC, ArthurUnipG} predict a decomposition of the automorphic discrete spectrum $L^2_{\mathrm{disc}}(\mathbf{G}(K)\backslash \mathbf{G}(\mathbf{A}_K))$ indexed by suitable global analogues of $L$-parameters. In fact, given an abstract irreducible admissible $\mathbf{G}(\mathbf{A}_K)$-representation $ \pi = \otimes_{v}' \pi_v$, Arthur gave a precise conjectural formula for the multiplicity of $\pi$ in the discrete automorphic spectrum of $\mathbf{G}$. However, this formula is contingent on the local Langlands conjecture for the local groups $\mathbf{G}(K_v)$.

In the 1980s, the Langlands program split off in an unforseen \emph{geometric} direction. Motivated by Grothendieck's function-sheaf dictionary and Weil's formula
\[\Bun_G(X)(\mathbf{F}_q) = G(K) \backslash G(\mathbf{A}_K) / G(\mathbf{O}_K)\]
for $G/\mathbf{F}_q$ a connected reductive group and $X/\mathbf{F}_q$ a smooth projective geometrically connected curve with function field $K$, Drinfeld realized that the notion of automorphic function could be geometrized into a suitable notion of \emph{automorphic sheaf} \cite{DrinfeldAJM}. This idea, and related insights about moduli spaces of shtukas, turned out to be incredibly fruitful and spawned a huge amount of work which among other things led to a resolution of many of Langlands's original conjectures over function fields \cite{LLafforgue, VLafforgue}. Moreover, the notion of automorphic sheaf on $\Bun_G(X)$ can be interpreted in many different sheaf theories, and in particular makes perfectly good sense for $X/\mathbf{C}$ if we take the sheaf theory of \emph{D-modules}. In this setting, a precise categorical geometric Langlands conjecture was formulated by Arinkin-Gaitsgory \cite{AGsingular}, building on earlier work of Beilinson-Drinfeld: they predicted an equivalence of categories
\[\operatorname{D-mod}(\Bun_G(X)) \simeq \IndCoh_{\mathrm{Nilp}}(\mathrm{LocSys}_{\hat{G}}(X))\]
compatible with the spectral action. This conjecture is now a theorem, through the recent combined efforts of many people \cite{FR,GLC1,GLC2,GLC3,GLC4,GLC5}. For a wonderful survey of these developments, we refer the reader to Scholze's Bourbaki report \cite{ScholzeGLC}.

Until 2014, however, there remained a significant divide between the geometric Langlands program and the original arithmetic aspects of Langlands's vision. The situation then changed entirely in late 2014, with Scholze's discovery of diamonds and moduli spaces of local shtukas in mixed characteristic \cite{scholze-berkeley-lectures}, and Fargues's formulation of (an initial form of) his geometrization conjecture \cite{fargues-overview}. It quickly became clear with these ideas that a significant portion of the machinery from the geometric Langlands program should be adaptable to mixed characteristic, with the central object on the geometric side being the stack $\Bun_G$ of $G$-bundles on the Fargues-Fontaine curve. Turning this vision into precise mathematics took seven years of hard effort from Fargues and Scholze, and resulted in their magnificent paper \cite{FS}. At the end of their paper, they formulate a rigorous \emph{categorical} form of the local Langlands conjecture: as discussed above, they predict an equivalence
\[\D(\Bun_G) \simeq \IndCoh(\Par_G)\]
for any quasisplit $G/F$, compatible with the spectral action and normalized by a choice of Whittaker datum.\footnote{We implicitly take $\Qellbar$-coefficients here, as in the discussion below, where the condition of nilpotent singular support turns out to be automatic. Fargues-Scholze also formulate a version of this conjecture with integral coefficients, in which case one must write $\IndCoh_{\mathrm{Nilp}}$ on the spectral side.}  Note the clear and extremely strong analogy with the geometric Langlands conjecture formulated above. Already in \cite{FS}, they explained how this conjecture is perfectly compatible with the Kottwitz conjecture and the expected properties of supercuspidal $L$-packets. 

More generally, from this optic the mysterious pairs $(\phi,\rho)$ appearing in the classical local Langlands correspondence acquire an extremely clean geometric meaning: An $L$-parameter $\phi$ gives a point in the stack $\Par_G$, and the stabilizer group of this point is exactly $S_{\phi}$, so $\rho$ determines a vector bundle on the residual gerbe of $\phi$, which can then be extended to a quasicoherent sheaf on $\Par_G$. For generic $L$-parameters, the matching sheaf on $\D(\Bun_G)$ is expected to precisely encode the local Langlands correspondent of $(\phi,\rho)$ (see \cite[Section 3.1]{Beijing} for an exact statement).

\subsection{Overview of Fargues-Scholze and the main conjecture} \label{sec:intro-FS}

In this section we recall the main constructions in \cite{FS} and give a precise formulation of the categorical local Langlands conjecture. A more detailed account of many of the claims made here can be found in \cref{sec:automorphic-results}. We keep the setup from before, so we fix a finite extension $F/\mathbf{Q}_p$ and a quasisplit connected reductive group $G/F$. Fix also a prime $\ell\neq p$ and an isomorphism $\overline{\mathbf{Q}_\ell}\simeq \mathbf{C}$. 

On the geometric side, the main player is the stack $\Bun_G$ of $G$-bundles on the Fargues-Fontaine curve, and the associated category $\D(\Bun_G) = \D_{\mathrm{lis}}(\Bun_G,\Qellbar)$ of $\ell$-adic \'etale sheaves on $\Bun_G$. By a theorem of Fargues \cite{FarguesTorsors}, there is a canonical identification $|\Bun_G| = B(G)$, where $B(G)$ is the Kottwitz set of $F$-isocrystals with $G$-structure. By a theorem of Viehmann \cite{Viehmann}, this is even a homeomorphism of topological spaces if we equip $B(G)$ with the Newton point partial order topology. As such, the stack $\Bun_G$ is stratified into locally closed substacks $\Bun_{G}^{b}$ for $b \in B(G)$. On each stratum there is a canonical equivalence $\D(\Bun_{G}^{b}) \cong \D(G_b(F),\Qellbar)$ where the right side is the usual derived category of smooth representations. Here as usual $G_b$ denotes the stabilizer group of the isocrystal $b$, so $G_b$ is an inner form of a Levi subgroup in $G$. 

For any $b \in B(G)$, there are natural functors $i_{b}^{\ast}:\D(\Bun_G) \to \D(G_b(F),\Qellbar)$ and $i_{b!}: \D(G_b(F),\Qellbar) \to \D(\Bun_G)$ which admit right adjoints $i_{b \ast}$ and $i_{b}^{!}$ as usual. More surprisingly, the functor $i_{b}^{\ast}$ also admits a left adjoint $i_{b \sharp}$. In practice, it is better to work with the \emph{renormalized} versions of these functors \cite[Definition 1.1.1]{Beijing}, which are denoted $i_{b!}^{\ren}$, $i_{b}^{\ast\ren}$, etc. When $b$ is basic, the standard and renormalized functors agree, but for non-basic $b$ the renormalized functors are cleaner.

With this notation in hand, the following theorem summarizes the main features of the geometric side. Essentially all of this theorem is proved in \cite{FS}, except some portions of item (iii) which are proved in \cite{HHS}.

\begin{thm}\label{thm:BunGreminderintro} Keep the notation as above.
\begin{thmenum}
\item The category $\D(\Bun_G)$ is compactly generated. A sheaf $A\in \D(\Bun_G)$ is compact if and only if $i_{b}^{\ast}A$ is compact for all $b$ and zero for all but finitely many $b$. 

\item A sheaf $A \in \D(\Bun_G)$ is ULA (for the structure map) if and only if $i_{b}^{\ast}A$ is ULA for all $b$, if and only if $\RHom(B,A)$ lies in $\Perf(\Qellbar)$ for all compact $B$.

\item The functors $i_{b\sharp}$, $i_{b}^{\ast}$, $i_{b!}$, and $i_{b}^{!}$ preserve compact objects, and all five functors preserve ULA objects. If $A \in \D(G_b(F),\Qellbar)$ is compact, then $i_{b\ast}A$ has compact stalks.

\item The objects $i_{b\sharp}\cind_{K}^{G_b(F)} \Qellbar$ give a set of compact generators for $\D(\Bun_G)$. 

\item On compact sheaves, there is a canonical involutive anti-equivalence $\Dbz \circlearrowright \D(\Bun_G)^\omega$ called Bernstein-Zelevinsky duality, characterized by the identity
\[ \RHom(\Dbz A,B) = \Gamma_c(\Bun_G, A \otimes B)\]
for all $A \in \D(\Bun_G)^\omega$ and all $B \in \D(\Bun_G)$.
\end{thmenum}
\end{thm}

The main goal of the categorical local Langlands conjecture is a description of the category $\D(\Bun_G)$ in \emph{spectral} terms, i.e. in terms of coherent sheaves on the stack of $L$-parameters. For this, we need to recall the spectral side. Here the main geometric player is the stack $\Par_G$ of $\ell$-adically continuous $L$-parameters valued in $\phantom{}^L G(\Qellbar)$. It is a little bit subtle to give the precise moduli problem, but this can be done in several ways, and the resulting geometric object turns out to be extremely nice. The following result (in a more refined form) is due to Dat-Helm-Kurinczuk-Moss \cite{DHKM-ParG}, except for the result on the dualizing sheaf, which is sketched in \cite{Zhu} and proved in detail in \cref{rslt:ParG-is-lci-with-trivial-dualizing-sheaf}.

\begin{thm}\label{thm:ParGreminderintro}
\begin{thmenum}
\item The stack $\Par_G$ is an Artin stack locally of finite type over $\Spec \Qellbar$, which is reduced and equidimensional of dimension zero. Moreover, $\Par_G$ is a local complete intersection, and the dualizing sheaf of $\Par_G$ is isomorphic to the structure sheaf.

\item The stack $\Par_G$ has a countably infinite set of connected components, and every connected component is of the form $X_i/\hat{G}$ for some affine variety $X_i$. In particular, $\QCoh(\Par_G)$ is compactly generated, with compact objects the perfect complexes on $\Par_G$ with quasicompact support.

\item \label{thm:ParGcoarsereminderintro} The coarse moduli space $X_{G}^{\spec}$ of $\Par_G$ is a disjoint union of affine varieties, each of which is isomorphic to the (coarse) quotient of a torus by a finite group, and the natural map
\[q: \Par_{G} \to X_{G}^{\spec}\] is a good moduli space in the sense of Alper. In particular, the ring 
\[\mathfrak{Z}_{G}^{\spec} \overset{\mathrm{def}}{=} \mathcal{O}(X_{G}^{\spec}) \cong \mathcal{O}(\Par_G)\]
is a countable product of reduced finite type Cohen-Macaulay $\Qellbar$-algebras. Moreover, the closed points of $X_{G}^{\spec}$ are in canonical bijection with semisimple $L$-parameters.
\end{thmenum}
\end{thm}

The starting point of our work is the following very deep theorem of Fargues-Scholze \cite[Corollary X.1.3]{FS}.

\begin{thm} There is a canonical tensor action of $\QCoh(\Par_G)$ on $\D(\Bun_G)$, called the \emph{spectral action}.
\end{thm}
To fix notation, we write $\mathcal{F} \ast A \in \D(\Bun_G)$ for the spectral action of $\mathcal{F} \in \QCoh(\Par_G)$ on $A \in \D(\Bun_G)$. We briefly recall the essence of this theorem. Using a suitable version of the geometric Satake equivalence, Fargues-Scholze construct Hecke operators $T_V$ acting on $\D(\Bun_G)$ for all $V \in \Rep \hat{G}$ such that the association $V \rightsquigarrow T_V$ is monoidal and such that $T_V(A)$ carries a natural $W_F$-action. On the other hand, the parameter stack comes with a canonical map $\Par_G \to B \hat{G}$ induced by the framing of the universal parameter, so any $V \in \Rep \hat{G} \subset \Perf(B \hat{G})$ pulls back to a $W_F$-equivariant vector bundle $V$ on $\Par_G$. To first approximation, the spectral action is uniquely characterized by the existence of identifications $T_V(-)= V \ast (-)$ which are monoidal and compatible with the $W_F$-equivariant structure on both sides; we refer the reader to \cref{sec:spectral-action} for a detailed account of this construction.

Now, the most optimistic guess would be that there is some equivalence $\D(\Bun_G) \simeq \QCoh(\Par_G)$, compatible with the $\QCoh(\Par_G)$-action on both sides. While this is in fact true when $G$ is a torus \cite{ZouTori}, it fails for all other groups. As in classical geometric Langlands, the main issue is that one expects such an equivalence to be compatible with Eisenstein series functors on both sides. However, the Eisenstein functor on $\D(\Bun_G)$ preserves compact objects \cite[Theorem 1.2.1]{HHS}, but the Eisenstein functor on $\QCoh(\Par_G)$ does not preserve compacts.

To fix this, one instead considers the category $\IndCoh(\Par_G)$ of ind-coherent sheaves. We note that on general stacks $\IndCoh$ is defined by descent and is not very explicit, but on $\Par_G$ we have a canonical equivalence $\IndCoh(\Par_G) = \Ind(\Coh^\qc(\Par_G))$, where by definition $\Coh^\qc(\Par_G) \subseteq \QCoh(\Par_G)$ denotes bounded complexes with coherent cohomologies and quasicompact support. Ind-completing the natural tensor action of $\Perf^{\qc}$ on $\Coh^\qc$, there is a natural tensor action of $\QCoh(\Par_G)$ on $\IndCoh(\Par_G)$. We are then led to hope for some equivalence of categories $\D(\Bun_G) \simeq \IndCoh(\Par_G)$ which is $\QCoh(\Par_G)$-linear.

This is still too vague. However, note that if we do have such an equivalence, there will be some object in $\D(\Bun_G)$ which matches the structure sheaf on $\Par_G$. We rigidify the conjecture by specifying this object. For this we actually need to make a small choice. Precisely, we need to choose a \emph{Whittaker datum} for $G$, i.e. a pair $(B,\psi)$ where $B \subset G$ is a Borel subgroup and $\psi:U(F) \to \Qellbar^\times$ is a nondegenerate character on the unipotent radical of $B$. If $G=\GL_n$, or more generally if $G$ is split with connected center, there is a unique isomorphism class of Whittaker data, but in general there are finitely many isomorphism classes of Whittaker data which are permuted simply transitively by the finite abelian group $G_{\mathrm{ad}}(F)/G(F)$. For example, $\mathrm{SL}_{2}$ admits four isomorphism classes of Whittaker data if $p>2$. We note that all isomorphism classes of Whittaker data can be obtained by fixing $B$ and varying $\psi$; as such we will usually notate dependence on a Whittaker datum by the symbol $\psi$ alone. 

In any case, given a Whittaker datum $(B,\psi)$, we write $W_{\psi} = \cind_{U(F)}^{G(F)}\psi$ for the associated Whittaker representation. This is a smooth $G(F)$-representation with excellent properties: among other things, it is projective, and its projection to any Bernstein component is compact. Using the simple-minded functor $i_{1!}$, we can extend $W_{\psi}$ to a sheaf on $\Bun_G$. We now arrive at the following slightly refined conjecture.

\begin{conjecture}[Categorical local Langlands conjecture, first version] \label{conj:cllcvague} After choosing a Whittaker datum $(B,\psi)$, there is some equivalence of categories
\[\mathbf{L}_{\psi}:\D(\Bun_G) \simeq \IndCoh(\Par_G)\]
which is colimit-preserving and $\QCoh(\Par_G)$-linear, and such that $\mathbf{L}_{\psi}(i_{1!}W_{\psi}) \simeq \mathcal{O}_{\Par_G}$.
\end{conjecture}
This still looks perhaps somewhat vague. However, let us analyze some consequences of this form of the conjecture. For this, we recall that the natural inclusions $\Perf^{\qc}(\Par_G) \hookrightarrow \Coh^\qc(\Par_G) \hookrightarrow \QCoh(\Par_G)$ ind-extend to colimit-preserving functors 
\[\QCoh(\Par_G) \overset{\Xi}{\to} \IndCoh(\Par_G) \overset{\Psi}{\to} \QCoh(\Par_G)\]
which turn out to have excellent properties. In particular, $\Xi$ is fully faithful, and in fact $\Xi(\mathcal{F})=\mathcal{F} \otimes \mathcal{O}_{\Par_G}$ using that $\Par_G$ has trivial dualizing sheaf. Slightly more surprisingly, it is also true that $\Psi$ is the right adjoint to $\Xi$. Now, suppose given a $\QCoh(\Par_G)$-linear equivalence $\mathbf{L}_{\psi}$ as in Conjecture \ref{conj:cllcvague}. Then for any $\mathcal{F} \in \QCoh(\Par_G)$, we compute that
\begin{align*}
    \mathbf{L}_{\psi}(\mathcal{F} \ast i_{1!}W_{\psi}) & = \mathcal{F} \otimes \mathbf{L}_{\psi} (i_{1!}W_{\psi}) \\
    & = \mathcal{F} \otimes \mathcal{O}_{\Par_G} \\
    & = \Xi(\mathcal{F}).
\end{align*}
Here we used the linearity in the first line, and the stipulated value of $\mathbf{L}_{\psi}$ on the Whittaker sheaf in the second line. This suggests introducing the functor
\begin{align*}
    a_{\psi}: \QCoh(\Par_G) & \to \D(\Bun_G) \\
    \mathcal{F} & \mapsto \mathcal{F} \ast i_{1!} W_{\psi}.
\end{align*}
By the basic properties of the spectral action, this functor preserves colimits and is $\QCoh(\Par_G)$-linear. In this notation, the preceding discussion shows that if Conjecture \ref{conj:cllcvague} is true, the given equivalence automatically makes the diagram
\[
\xymatrix{\mathrm{QCoh}(\mathrm{Par}_{G})\ar[r]^{a_{\psi}}\ar[dr]^{\Xi} & \D(\mathrm{Bun}_{G})\ar[d]_{\wr}^{\mathbf{L}_{\psi}}\\
 & \mathrm{IndCoh}(\mathrm{Par}_{G})
}
\]
commute. 

Next, observe that since $a_{\psi}$ preserves colimits, it has an abstractly defined right adjoint \[c_{\psi}: \D(\Bun_G) \to \QCoh(\Par_G).\] 
In classical geometric Langlands, the functor $c_{\psi}$ is usually called the \textbf{enhanced Whittaker coefficient}, and this gadget will play a central role in our analysis. We caution the reader that a priori, $c_{\psi}$ is not known to preserve colimits or be $\QCoh(\Par_G)$-linear, although it is $\Perf(\Par_G)$-linear. From its definition as a right adjoint, it is an easy exercise to check that
\[\RHom(i_{1!}W_{\psi},T_{V}A) = \Gamma(\Par_G,V \otimes c_{\psi}(A))\]
for any $V \in \Rep \hat{G}$ and any automorphic sheaf $A$. Thinking of the left side as the $V$th Whittaker coefficient of $A$, we see that $c_{\psi}(A)$ captures these coefficients for all $V$ simultaneously.

In any case, passing to the right adjoint of the preceding diagram, we get a commutative diagram
\[
\xymatrix{\D(\mathrm{Bun}_{G})\ar[d]_{\wr}^{\mathbf{L}_{\psi}}\ar[r]^{c_{\psi}} & \mathrm{QCoh}(\mathrm{Par}_{G})\\
\mathrm{IndCoh}(\mathrm{Par}_{G})\ar[ur]^{\Psi}
}
\]
Since $\mathbf{L}_{\psi}$ is an equivalence of compactly generated categories, it will restrict to an equivalence from $\D(\Bun_G)^\omega$ towards $\Coh^\qc(\Par_G)$, so we can enlarge the preceding diagram to a diagram
\[
\xymatrix{\D(\Bun_G)^{\omega}\ar[d]_{\wr}^{\mathbf{L}_{\psi}}\ar[r] & \D(\mathrm{Bun}_{G})\ar[d]_{\wr}^{\mathbf{L}_{\psi}}\ar[r]^{c_{\psi}} & \mathrm{QCoh}(\mathrm{Par}_{G})\\
\Coh^\qc(\Par_G)\ar[r] & \mathrm{IndCoh}(\mathrm{Par}_{G})\ar[ur]^{\Psi}
}
\]
where the unlabelled horizontal arrows are the evident full inclusions. Now we make an absolutely crucial observation: tautologically from its definition, the functor $\Psi$ restricted to $\Coh^\qc(\Par_G) \subset \IndCoh(\Par_G)$ is just the usual inclusion $\Coh^\qc(\Par_G) \subset \QCoh(\Par_G)$. Equivalently, $\Psi$ restricts to the identity functor from $\Coh^\qc(\Par_G) \subset \IndCoh(\Par_G)$ towards $\Coh^\qc(\Par_G) \subset \QCoh(\Par_G)$. Under this identification, the preceding diagram shows that $\mathbf{L}_{\psi}(A) = c_{\psi}(A)$ for all $A \in \D(\Bun_G)^\omega$. In other words, after restricting to compact sheaves, $\mathbf{L}_{\psi}$ and $c_{\psi}$ must give \emph{the same functor} towards $\Coh^\qc(\Par_G)$. However, unlike the conjectural functor $\mathbf{L}_{\psi}$, $c_{\psi}$ is completely canonical and defined unconditionally for all groups! Moreover, in restricting to compact sheaves we haven't lost any information at all, since the full equivalence $\mathbf{L}_{\psi}$ can be recovered as the ind-completion of its restriction to compact objects. Summarizing this discussion, we have arrived at a proof of the following result.

\begin{prop}\label{prop:CLLCuniqueintro}There is \emph{at most one} equivalence $\mathbf{L}_\psi : \D(\Bun_G) \overset{\sim}{\to} \IndCoh(\Par_G)$ satisfying either of the following equivalent conditions:
\begin{propenum}
    
\item The diagram
\[
\xymatrix{\mathrm{QCoh}(\mathrm{Par}_{G})\ar[r]^{a_{\psi}}\ar[dr]^{\Xi} & \D(\mathrm{Bun}_{G})\ar[d]_{\wr}^{\mathbf{L}_{\psi}}\\
 & \mathrm{IndCoh}(\mathrm{Par}_{G})
}
\]
commutes.

\item The diagram
\[
\xymatrix{\D(\mathrm{Bun}_{G})\ar[d]_{\wr}^{\mathbf{L}_{\psi}}\ar[r]^{c_{\psi}} & \mathrm{QCoh}(\mathrm{Par}_{G})\\
\mathrm{IndCoh}(\mathrm{Par}_{G})\ar[ur]^{\Psi}
}
\]
commutes.
\end{propenum}
Moreover, such an equivalence exists if and only if $c_\psi$ restricts to an equivalence of categories
\[c_\psi : \D(\Bun_G)^{\omega} \overset{\sim}{\to} \Coh^\qc(\Par_G),\]
in which case $\mathbf{L}_{\psi}$ is obtained from this equivalence by ind-completion.
\end{prop}

In other words, the equivalence of categories postulated in Conjecture \ref{conj:cllcvague}, which initially looked perhaps somewhat underdetermined, is in fact extremely rigid: there is at most one such equivalence, and if it does exist it is completely determined by the unconditional functor $c_{\psi}$. Putting this functor at center stage, we can now state the official form of the categorical local Langlands conjecture.

\begin{conjecture}[Categorical local Langlands conjecture, official version] \label{conj:CLLCintro}For any quasisplit $G$ equipped with a Whittaker datum $(B,\psi)$, the functor $c_{\psi}$ restricts to an equivalence of categories $c_{\psi}:\D(\Bun_G)^{\omega} \overset{\sim}{\to} \Coh^\qc(\Par_G)$. Equivalently, the functor $c_\psi$ carries $\D(\Bun_G)^{\omega}$ into $\Coh^\qc(\Par_G)$, and the resulting ind-completed functor
\[\mathbf{L}_{\psi}: \D(\Bun_G) \to \IndCoh(\Par_G)\]
is an equivalence of categories.
\end{conjecture}

Our goal in this paper is to study Conjecture \ref{conj:CLLCintro} from first principles, and to prove it for many groups. Notably, we prove the conjecture for $\GL_n$ (assuming the compatibility of $c_{\psi}$ with Eisenstein series).

\subsection{First properties of the conjecture} \label{sec:firstproperties}

At present, it seems hopeless to attack Conjecture \ref{conj:CLLCintro} in all generality, because it simply contains far too much deep arithmetic information. As an instructive example, we mention the following result (a proof of this statement will appear elsewhere).
\begin{prop}Let $G$ be a quasisplit group with a Whittaker datum $(B,\psi)$, and suppose Conjecture \ref{conj:CLLCintro} is true. Then:
\begin{propenum}
    \item The Fargues-Scholze map $\pi \mapsto \varphi_{\pi}$ from smooth irreducible $G(F)$-representations towards semisimple $L$-parameters is surjective and has finite fibers.

    \item (Johansson) For every semisimple parameter $\varphi$, there is a unique irreducible $G(F)$-representation $\pi$ with $\varphi_{\pi} = \varphi$ such that $\Hom(W_{\psi},\pi) \neq 0$.
\end{propenum}
\end{prop}
For general groups, both the surjectivity of the Fargues-Scholze map and the finiteness of its fibers are completely open, while (ii) is an even stronger form of surjectivity as it exhibits a distinguished section of the map $\pi \mapsto \varphi_{\pi}$. We emphasize that this result was chosen somewhat at random, and Conjecture \ref{conj:CLLCintro} implies many other things which seem very difficult to access directly. Nevertheless, it is still very enlightening to contemplate various properties and implications of Conjecture \ref{conj:CLLCintro}, and the goal of this section is to record some results along these lines. Later, we will turn many of these results around and use them as guidance on our path towards the proof.

First, we observe that the categorical equivalence must match certain finiteness conditions on sheaves on both sides. To explain this, recall that the spectral action induces a natural map $\mathfrak{Z}_{G}^{\spec} \to \mathfrak{Z}(\D(\Bun_G))$, where $\mathfrak Z_G^\spec = \calO(X^\spec_G)$ is the spectral Bernstein center and $\mathfrak{Z}(\mathcal{C})$ denotes the center of the category $\mathcal{C}$. In particular, for any sheaf $A \in \D(\Bun_G)$ there is a natural ring map $\mathfrak{Z}_{G}^{\spec} \to \End(A)$. If $A$ is compact, the radical of the kernel of this map cuts out a quasicompact closed subscheme $\mathrm{spec.supp}(A) \subset X_{G}^{\spec}$ which we call the \emph{spectral support} of $A$ (see \cref{def:spectral-support} for the general definition). When $A=i_{1!}\pi$ for some irreducible representation $\pi$, the spectral support is just the closed point corresponding to the Fargues-Scholze parameter of $\pi$. If the equivalence $\mathbf{L}_{\psi}$ exists, it is completely formal to see that for any compact $A$, the support of the coherent sheaf $\mathbf{L}_{\psi}(A)$ is contained in $q^{-1}(\mathrm{spec.supp}(A))$, where $q:\Par_G \to X_{G}^{\spec}$ is the map to the coarse moduli as in \cref{thm:ParGcoarsereminderintro}.

\begin{defn}Let $G$ be any reductive group.
\begin{defenum}
\item A sheaf $A \in \D(\Bun_G)$ is \emph{finite} if it compact and $\mathrm{spec.supp}(A)$ has finite length, or equivalently the ring map $\mathfrak{Z}_{G}^{\spec} \to \End(A)$ factors over an Artinian quotient. We write $\D(\Bun_G)_{\fin}$ for the resulting full subcategory.
\item A sheaf $\mathcal{F} \in \Coh(\Par_G)$ is \emph{finite} if it is supported set-theoretically on finitely many closed fibers of the map $q\colon \Par_G \to X_G^\spec$. We write $\Coh(\Par_G)_{\fin}$ for the resulting full subcategory.\footnote{We will quickly resolve the apparent notational clash with \cref{def:cohfinhardapproach}, see $(\dagger)$ and \cref{thm:spectralboundedintro} below.}
\end{defenum}
\end{defn}
It is clear that in both cases, finite sheaves form a thick triangulated subcategory. Moreover, it is immediate from the discussion above that if the equivalence $\mathbf{L}_{\psi}$ exists, it automatically restricts to an equivalence $\D(\Bun_G)_{\fin}\overset{\sim}{\to} \Coh(\Par_G)_{\fin}$. This alone is not particularly interesting, but it becomes interesting once we realize that finite sheaves on $\Bun_G$ have several equivalent characterizations. More precisely, we prove in Proposition \ref{prop:finite4ways} that $A\in \D(\Bun_G)$ is finite in the above sense if and only if $A$ is compact and ULA, if and only if $A$ can be obtained by finitely many shifts and cones from sheaves of the form $i_{b!}\pi$ where $\pi$ is an irreducible $G_b(F)$-representation. This last condition intuitively says that $A$ has finite length in the strongest possible sense. On the other hand, the second condition suggests that we should try to figure out which subcategory of $\IndCoh(\Par_G)$ matches with $\D(\Bun_G)^{\mathrm{ULA}}$. Here Theorem \ref{thm:BunGreminderintro}.(ii) gives the key clue: since compact sheaves on $\Bun_G$ match with $\Coh^\qc(\Par_G)$, a direct translation across the equivalence $\mathbf{L}_{\psi}$ shows that $A \in \D(\Bun_G)$ is ULA if and only if $\RHom(\mathcal{G},\mathbf{L}_{\psi}(A))$ is a perfect complex of $\Qellbar$-vector spaces for all $\mathcal{G} \in \Coh^\qc(\Par_G)$. We now turn this into a definition.

\begin{defn}Let $X$ be a disjoint union of QCA stacks over a characteristic zero field $k$. An ind-coherent sheaf $\mathcal{F} \in \IndCoh(X)$ is \emph{admissible} if $\RHom(\mathcal{G},\mathcal{F})$ is a perfect complex of $k$-vector spaces for all $\mathcal{G} \in \Coh^\qc(X)$. We write $\Adm(X) \subset \IndCoh(X)$ for the resulting full subcategory.
\end{defn}

We emphasize that while admissibility is a kind of finiteness condition on ind-coherent sheaves (in fact, it coincides with the suaveness condition in the 6-functor formalism of ind-coherent sheaves, see \cref{sec:adm-sheaves}), it is typically quite far from the familiar finiteness condition of coherence. For example, if $X=BG$ is the classifying stack of a positive-dimensional reductive group over $k$ (so $\IndCoh(X)=\QCoh(X)$) and $\mathcal{O}(G)$ is the ind-coherent sheaf given by the regular representation of $G$, then $\mathcal{O}(G)$ is admissible but not coherent. On the other hand, if $X$ is a separated finite type $k$-scheme, it is an amusing exercise to check that there is a containment $\Coh(X) \subseteq \Adm(X)$ \emph{if and only if} $X$ is smooth and proper over $k$.

Anyway, in this notation, the above discussion shows that an equivalence $\mathbf{L}_{\psi}$ automatically restricts to an equivalence $\D(\Bun_G)^{\mathrm{ULA}}\overset{\sim}{\to} \Adm(\Par_G)$. Now, finite sheaves on $\Bun_G$ are exactly the compact sheaves which are also ULA. Translating this across the equivalence, we arrive at the following conjecture, which can be unconditionally formulated for any group. 
\begin{itemize}
    \item[($\dagger$)] Inside $\IndCoh(\Par_G)$, there should be an identification of subcategories
    \[\Coh(\Par_G)_{\fin} = \Coh^\qc(\Par_G) \cap \Adm(\Par_G).\]
\end{itemize}

It's not hard to prove that any coherent admissible sheaf is finite, so the essential content of this conjecture is that every finite coherent sheaf is admissible. Concretely, this boils down to proving that for any $\mathcal{G} \in \Coh^\qc(\Par_G)$ and any $\mathcal{F} \in \Coh(\Par_G)_{\fin}$, $\RHom(\mathcal{G},\mathcal{F})$ is a perfect complex of $\Qellbar$-vector spaces. It is not hard to prove that $\RHom(\mathcal{G},\mathcal{F})$ is left-bounded and that all cohomology groups are finite-dimensional vector spaces, but the right-boundedness is severely non-obvious due to the singularities of $\Par_G$. 

In fact, contemplating this difficulty further, one realizes that the finiteness condition on $\mathcal{F}$ is a red herring. More precisely, if $A,B$ are any compact sheaves on $\Bun_G$, it is not hard to prove that $\RHom(A,B)$ is a bounded complex of $\Qellbar$-vector spaces, using Theorem \ref{thm:BunGreminderintro}.(i) together with Bernstein's fundamental finite global dimension theorem for Hecke algebras of $p$-adic groups. Translating this back to the spectral side, one is finally led to guess the following result, which settles $(\dagger)$ by the preceding discussion.

\begin{thm}\label{thm:spectralboundedintro}For any reductive group $G$ and any $\mathcal{F},\mathcal{G} \in \Coh^\qc(\Par_G)$, $\RHom(\mathcal{G},\mathcal{F})$ is a bounded complex.
\end{thm}

Again, we emphasize that the essential difficulty here is that $\Par_G$ is singular. Note that for any singular \emph{variety}, no boundedness property like this can possibly hold for all coherent sheaves, so its truth for $\Par_G$ indicates that the stackiness of $\Par_G$ is somehow negating the singularities. This is a difficult theorem, and we will sketch the proof later in the introduction (the impatient reader can skip to Theorem \ref{thm:spectral-main}).

The next key observation is that if Conjecture \ref{conj:CLLCintro} is true, the equivalence $\mathbf{L}_{\psi}$ enjoys an automatic compatibility with duality. To explain this, recall that compact sheaves on $\Bun_G$ admit the involutive Bernstein-Zelevinsky self-duality $\Dbz$ recalled above. On the spectral side, $\Coh^\qc(\Par_G)$ carries the usual involutive Grothendieck-Serre duality $\Dgs$. However, it is clear that these can't match immediately under the conjectural equivalence, because any functor which matches $\Dbz$ must change the support of the underlying coherent sheaf. For instance, if $G$ is semisimple and $\pi$ is an irreducible supercuspidal with Fargues-Scholze parameter $\varphi$, $\mathbf{L}_{\psi}(i_{1!}\pi)$ will be supported on the fiber of $q$ over $\varphi$, while $\mathbf{L}_{\psi}(\Dbz i_{1!}\pi)= \mathbf{L}_{\psi}( i_{1!}\pi^\vee)$ will be supported on the fiber of $q$ over $c \circ \varphi$. Here $c$ is the Chevalley involution on $\phantom{}^L G$, and we also write $c:\Par_G \overset{\sim}{\to} \Par_G$ for the involutive automorphism induced by composing an $L$-parameter with the Chevalley involution. Noting that $c$ moves the fiber over $\varphi$ to the fiber over $c \circ \varphi$, we are led to contemplate the functor $\Dtwgs \overset{\mathrm{def}}{=} c^{\ast} \Dgs$ of \emph{twisted} Grothendieck-Serre duality. Since $\Par_G$ has trivial dualizing complex, $c^{\ast}$ commutes with $\Dgs$, so $\Dtwgs$ is still an involutive self-equivalence on $\Coh^\qc(\Par_G)$.

\begin{thm}[Theorem \ref{thm:automaticduality}.(i)]\label{thm:GSdualityintro}Let $G$ be a quasisplit group with a Whittaker datum $(B,\psi)$. Suppose Conjecture \ref{conj:CLLCintro} is true for $G$ equipped with both Whittaker data $(B,\psi^{\pm 1})$. Then there is a canonical equivalence of functors
\[\Dtwgs \mathbf{L}_{\psi} \cong \mathbf{L}_{\psi^{-1}} \Dbz : \D(\Bun_G)^{\omega} \to \Coh^\qc(\Par_G). \]
\end{thm}
The necessity to assume Conjecture \ref{conj:CLLCintro} for both the datum $(B,\psi)$ and the dual datum $(B,\psi^{-1})$ is a slight annoyance, and we will comment more generally on the dependence of this conjecture on the Whittaker datum in Section \ref{ss:changeofWhittaker}.

The proof of this theorem is somewhat tricky. The key ingredient is a duality theorem for $a_{\psi}$, which we can state much more generally. More precisely, if $G$ is a quasisplit group with a Whittaker datum $(B,\psi)$, we say $G$ is ($\psi$-)\emph{reasonable} if the functor $a_{\psi}$ preserves compact objects. This is equivalent to a rather concrete condition on $W_{\psi}$ (Proposition \ref{prop:apsifinitebasic}), and is also equivalent to requiring that $c_{\psi}$ be colimit-preserving (Proposition \ref{prop:decentequivalents}). Moreover, this is implied by (and much weaker than) finiteness of fibers for the Fargues-Scholze map $\pi \mapsto \varphi_{\pi}$, so many groups are reasonable with our current knowledge, including $\GL_n$ and all quasisplit classical groups of types A, B and D.

\begin{thm}[Theorem \ref{thm:apsiduality}] If $G$ is $\psi$-reasonable, it is also $\psi^{-1}$-reasonable, and there is a canonical equivalence of functors
\[\Dbz a_{\psi} \cong a_{\psi^{-1}} \Dtwgs :\Perf^{\qc}(\Par_G) \to \D(\Bun_G)^{\omega}.\]
\end{thm}

This theorem in turn depends on two key ingredients. One of them (Theorem \ref{thm:Whittakerbasics}.(iv)) is a duality theorem for $W_{\psi}$ projected to a Bernstein component, which was already stated in \cite{Beijing}. The other ingredient, which is new, is the following extremely useful duality theorem for the spectral action.
\begin{thm}[Theorem \ref{thm:spectraldualityuseful}]For all $A \in \D(\Bun_G)^\omega$ and $\mathcal{F} \in \Perf(\Par_{G})$, we have a functorial isomorphism
    \[ \Dbz(\mathcal{F} \ast A) = \Dtwgs\mathcal{F} \ast \Dbz A\]
    in $\D(\Bun_G)^{\omega}$. 
\end{thm}
In fact, with future applications in mind, we prove a more general result which works with integral coefficients as well. The proof of Theorem \ref{thm:spectraldualityuseful} ultimately boils down to a certain compatibility of the geometric Satake equivalence with the Chevalley involution, which was already proved in Fargues-Scholze \cite[Proposition VI.12.1]{FS}. However, in reducing the claimed duality theorem to this result we must (unsurprisingly) dig into the construction of the spectral action, and we found ourselves confused by certain aspects of the construction of the spectral action which are given with rather light detail in \cite{FS} and \cite{FS-motivic}. In the end, this resulted in the somewhat technical Appendices \ref{sec:extended-geom-setups}-\ref{sec:convolution-stacks} and \ref{sec:abstract-Hecke-action}, whose content can be summarized informally in the following theorem.
\begin{thm}Let $G$ be any reductive group, and let $\D:\vStk \to \Pr^L$ be a sheaf theory underlying a reasonable six-functor formalism on small v-stacks. For any finite set $I$, let  $\locHck_{G}^{I}$ be the local Hecke stack over $\Div^I$.
\begin{thmenum}
    \item There is a natural subcategory $\D^{\bdd} (\locHck_{G}^{I})\subset \D(\locHck_{G}^{I})$ of bounded sheaves equipped with a $\D(\Div^I)$-linear convolution monoidal structure which is functorial in $I$.
    \item There is a natural $\D(\Div^I)$-linear action of $\D^{\bdd} (\locHck_{G}^{I})$ on $\D(\Bun_G \times \Div^I)$, induced by a natural $\D(\Div^I)$-linear monoidal functor
    \[\D^{\bdd} (\locHck_{G}^{I}) \to \D(\Bun_G \times \Bun_G \times \Div^I)\]
    which is functorial in $I$.
\end{thmenum}
\end{thm}
We emphasize that this is an informal statement, and the point of the appendices is to state and prove a precise form of this theorem which explicitly accounts for all the implicit higher coherences. We also make the technical change of rebuilding the spectral action using the category of overconverget motivic sheaves defined by Scholze \cite{ScholzeBerkovich}, or more precisely using its base change $\D_{\rel}$ along the $\ell$-adic realization functor over a geometric base point as considered in \cite{HHS}: by forthcoming joint work of Gleason and the first author \cite{GH}, on $\Bun_G$ this category is canonically equivalent to $\D_{\mathrm{lis}}(\Bun_G)$ compatibly with Hecke operators. This eliminates the need to use the dual embedding on the Satake category when defining Hecke operators, and gives us access to a full six-functor formalism in the arguments, which is lacking in the setting of solid sheaves. We refer the reader to \cref{sec:def-of-D-rel} for an introduction to $\D_\rel$.

This is not the end of the duality story. Recall that aside from the Bernstein-Zelevinsky duality on $\D(\Bun_G)$, we also have the contravariant functor of Verdier duality $\Dverd$, which restricts to an involutive self-equivalence on ULA sheaves. This suggests the question of whether there is some contravariant duality functor on $\IndCoh(\Par_G)$ which matches Verdier duality under the categorical conjecture. Here the key clue comes from the tight relationship between Bernstein-Zelevinsky duality and Verdier duality: there is a ``duality exchange'' formula
\[\RHom(A, \Dverd B) = \RHom(\Dbz A, B)^\vee \]
valid for all $A \in \D(\Bun_G)^\omega$ and $B\in \D(\Bun_G)$. The matching notion on the spectral side turns out to be the following operation, which we call \emph{admissible duality}.

\begin{prop}Let $X$ be a disjoint union of QCA stacks over a characteristic zero field $k$. There is a natural contravariant functor $\Dadm: \IndCoh(X) \to \IndCoh(X)$ defined by the formula
\[\RHom(\mathcal{G},\Dadm \mathcal{F}) = \Gamma^{\IndCoh}(X,\mathcal{F} \otimes^{!} \mathcal{G})^\vee\]
for all $\mathcal{F},\mathcal{G} \in \IndCoh(X)$.
This functor is also characterized by the duality exchange formula
\[\RHom(\mathcal{G},\Dadm \mathcal{F}) = \RHom(\Dgs \mathcal{G},\mathcal{F})^\vee\]
valid for all $\mathcal{G} \in \Coh^\qc(X)$ and $\mathcal{F} \in \IndCoh(X)$.

An ind-coherent sheaf $\mathcal{F}$ is admissible if and only if $\Dadm \mathcal{F}$ is admissible, and $\Dadm$ restricts to an involutive self-equivalence on admissible ind-coherent sheaves.
\end{prop}

Here $\Gamma^{\IndCoh}(X,-)$ denotes the (colimit-preserving) functor $\Gamma_{!}(X,\Psi(-))$, where $\Gamma_!(X,-):\QCoh(X) \to \D(k)$ is the unique colimit-preserving functor agreeing with $\Gamma(X,-)$ on complexes with quasicompact support. We note that the defining formula of $\Dadm$ immediately implies its existence: for any fixed $\mathcal{F}$, the functor $\mathcal{G} \mapsto \Gamma^{\IndCoh}(X,\mathcal{F} \otimes^{!} \mathcal{G})^\vee$ converts colimits into limits, and therefore is representable by some unique object $\Dadm \mathcal{F}$. In fact, as we discuss in \cref{sec:adm-sheaves}, $\Dadm$ is simply the suave duality functor in the 6-functor formalism on $\IndCoh$. The remaining parts of this proposition are not particularly difficult. Moreover, using Theorem \ref{thm:GSdualityintro} along with the perfect parallel between the duality exchange formulas on both sides, it is straightforward to prove the following matching of dualities.

\begin{prop}[Theorem \ref{thm:automaticduality}.(ii)]\label{prop:Admdualityintro}Let $G$ be a quasisplit group with a Whittaker datum $(B,\psi)$. Suppose Conjecture \ref{conj:CLLCintro} is true for $G$ equipped with both Whittaker data $(B,\psi^{\pm 1})$. Then there is a canonical equivalence of functors
\[\Dtwadm \mathbf{L}_{\psi} \cong \mathbf{L}_{\psi^{-1}} \Dverd : \D(\Bun_G) \to \IndCoh(\Par_G) \]
where $\Dtwadm \overset{\mathrm{def}}{=} c^{\ast} \Dadm$ denotes the functor of \emph{twisted} admissible duality.
\end{prop}

Now we come to a great surprise. To set the stage for this, we point out that unlike the mild and well-understood Grothendieck-Serre duality, the functor of admissible dual typically destroys coherence. For a simple-minded example, if $X=\Spec R$ is a smooth affine $k$-variety with (invertible) dualizing complex $\omega_X \in \Perf(X)$, and $\mathcal{F}=\widetilde{M}$ is the coherent sheaf associated with some finitely generated $R$-module $M$, it is not difficult to check that $\Dadm \mathcal{F} = \omega_X \otimes \widetilde{N}$ where $N=\mathrm{Hom}_{k}(M,k)$ is a typically huge $R$-module. On algebraic stacks, admissible duals are typically still huge, and much harder to compute. However, admissible duality turns out to have the following astonishing property.

\begin{thm}[Theorem \ref{rslt:spectral-temperization}]\label{thm:spectraltemperizationintro}Let $X$ be a disjoint union of QCA stacks over a characteristic zero field $k$. Assume that $X$ is Gorenstein. Then for any $\mathcal{F} \in \Coh(X) \subset \IndCoh(X)$, the sheaf $\Dadm \mathcal{F}$ lies in the essential image of the fully faithful functor $\Xi: \QCoh(X) \to \IndCoh(X)$.
\end{thm}
When $X$ is smooth this theorem has no content, but for singular $X$ there is no philosophical reason for it to hold. We discovered this result experimentally, by computing admissible duals on the stack $(\Spec k[x,y]/(xy))/\mathbf{G}_m$ where $\mathbf{G}_m$ acts trivially on $x$ and acts with weight two on $y$. This stack is a local model for the singularities of $\Par_{\mathrm{PGL}_2}$, and in doing these computations we relied heavily on some calculations with coherent sheaves on this parameter stack due to Bertoloni Meli and Koshikawa. Later we realized that Theorem \ref{thm:spectraltemperizationintro} is roughly analogous with a result of Beraldo, who proved that Verdier duality on (classical) $\Bun_G$ sends compact $D$-modules to \emph{tempered} $D$-modules.

Let us contemplate this theorem in conjunction with our analysis of Conjecture \ref{conj:CLLCintro} so far. If $\mathcal{F} \in \IndCoh(\Par_G)$ is any sheaf, then Proposition \ref{prop:Admdualityintro} implies that the inverse functor $\mathbf{L}_{\psi}^{-1}$ satisfies $\Dverd \mathbf{L}_{\psi}^{-1}\mathcal{F} \cong \mathbf{L}_{\psi^{-1}}^{-1} \Dtwadm \mathcal{F}$. If $\mathcal{F}$ is coherent, Theorem \ref{thm:spectraltemperizationintro} shows that $\Dtwadm \mathcal{F}$ lies in the essential image of the fully faithful functor $\Xi$, or equivalently that the natural map $\Xi \Psi \Dtwadm \mathcal{F} \to \Dtwadm \mathcal{F}$ is an isomorphism. Combining these observations shows that for coherent $\mathcal{F}$, we have a canonical isomorphism
\[\Dverd \mathbf{L}_{\psi}^{-1}\mathcal{F} \cong \mathbf{L}_{\psi^{-1}}^{-1} \Xi \Psi \Dtwadm \mathcal{F}.\]
Moreover, by \ref{prop:CLLCuniqueintro}.(i), $\mathbf{L}_{\psi^{-1}}^{-1} \Xi$ agrees with the functor $a_{\psi^{-1}}$, so combining this with the previous formula we deduce that
\[\Dverd \mathbf{L}_{\psi}^{-1}\mathcal{F} \cong  a_{\psi^{-1}} \Psi \Dtwadm \mathcal{F}\]
for coherent $\mathcal{F}$. Now typically we cannot eliminate the Verdier dual on the left here, but if $\mathbf{L}_{\psi}^{-1}\mathcal{F}$ is ULA, or equivalently if $\mathcal{F}$ is admissible, we \emph{can} eliminate it by applying Verdier duality again on both sides and using Verdier biduality on ULA sheaves. In summary, this analysis shows that if Conjecture \ref{conj:CLLCintro} is true, then for all $\mathcal{F}$ which are coherent and admissible, or equivalently finite by our previous analysis, there is an equivalence
\[ \mathbf{L}_{\psi}^{-1}\mathcal{F} \cong \Dverd a_{\psi^{-1}} \Psi \Dtwadm \mathcal{F}.\]
While this formula should only be true for $\mathcal{F} \in \Coh(\Par_G)_{\fin}$, which is a rather thin subcategory of $\IndCoh(\Par_G)$, there is an enormous advantage to this formula: the right-hand side is defined unconditionally for any group.

\begin{defn}For $G$ quasisplit with a fixed Whittaker datum $(B,\psi)$, we define the functor
\begin{align*}
    t_{\psi}: \Coh(\Par_G)_{\fin} & \to \D(\Bun_G) \\
    \mathcal{F} & \mapsto \Dverd a_{\psi^{-1}} \Psi \Dtwadm \mathcal{F}.
\end{align*}
\end{defn}

Again, we emphasize that this functor is defined unconditionally in all generality. The following conjecture follows from Conjecture \ref{conj:CLLCintro} by the discussion above, but deserves to be stated separately.
\begin{conjecture}\label{conj:CLLCrestrictedintro}For $G$ quasisplit with a fixed Whittaker datum, the functor $t_{\psi}$ induces an equivalence of categories
\[t_{\psi}:\Coh(\Par_G)_{\fin} \overset{\sim}{\to} \D(\Bun_G)_{\fin}.\]
\end{conjecture}

This is a form of categorical local Langlands with ``restricted variation'' (see below for some more comments on this terminology). Of course, this conjecture is probably completely out of reach in general. However, it turns out we can prove a tight relationship between $c_{\psi}$ and $t_{\psi}$ under a very mild assumption on $G$.

\begin{thm}[Theorem \ref{thm:tpsiadjunction}] \label{thm:ctadjunctionintro}If $G$ is reasonable, there is a natural isomorphism
\[\RHom(c_{\psi}A,\Psi \mathcal{F}) \cong \RHom(A, t_{\psi}\mathcal{F})\]
for all $A\in \D(\Bun_G)$ and $\mathcal{F} \in \Coh(\Par_G)_{\fin}$.
\end{thm}

In other words, ignoring the functor $\Psi$ for the moment, $t_{\psi}$ gives a partially defined right adjoint to the functor $c_{\psi}$. This result is not a formality; after all, $c_{\psi}$ is itself defined as an abstract right adjoint, but here we are trying to compute maps out of it, rather than maps towards it. The proof is a somewhat elaborate calculation; the key ingredients are Theorem \ref{thm:apsiduality}, the duality exchange formulas on both sides, and Theorem \ref{thm:spectraltemperizationintro}.

This concludes our analysis of the formal properties and consequences of Conjecture \ref{conj:CLLCintro}. However, we still haven't touched the key question of how to attack this conjecture. In the next section we take up this task.

\subsection{Well-understood groups} \label{sec:wellunderstood}
As noted above, Conjecture \ref{conj:CLLCintro} seems completely out of reach in general. However, for some groups we already have a solid understanding of the classical local Langlands conjecture and its interaction with Fargues-Scholze $L$-parameters, and here the conjecture is more approachable. We give the following definition in a slightly imprecise form, referring to the main text for the precise definition (Definition \ref{def:wellunderstoodG}). 

\begin{defn}A quasisplit group $G$ with a Whittaker datum $(B,\psi)$ is \emph{well-understood} if there is a known Whittaker-normalized $B(G)_{\mathrm{basic}}$ local Langlands correspondence for $G$ and all of its standard Levi subgroups which satisfies the expected parametrizations of discrete series $L$-packets and supercuspidal $L$-packets, and the endoscopic character identities for them, and whose semisimplification agrees with the Fargues-Scholze construction of $L$-parameters.
\end{defn}

At present, many groups are well-understood, including $\mathrm{GL}_n$ and all classical groups of types A, B, and D, along with $\mathrm{GSp}_4$ and some groups related to these by products, central isogenies, etc. For a more thorough discussion and precise references, we refer to Section \ref{ss:wellunderstood}.

In this section we explain a strategy to prove Conjecture \ref{conj:CLLCintro} for well-understood groups. However, even giving ourselves this knowledge of the group still doesn't seem to be enough leverage, because the functor $c_{\psi}$ is so hard to access from first principles. Most of our remaining results will depend on the compatibility of this functor with parabolic induction. More precisely, we will need the following.

\begin{conjecture}[Conjecture \ref{conj:cpsiEiscompatible}] \label{conj:Eiscpsiintro}
Let $G$ be any quasisplit group with a Whittaker datum $(B,\psi)$ and maximal torus $T \subset B$, so every standard Levi subgroup $M \subset G$ inherits its own Whittaker datum $(B\cap M, \psi_M)$. Then for every standard parabolic $P\subset G$ with Levi $M$, there is a commutativity datum for the diagram
\[
\xymatrix{\D(\mathrm{Bun}_{M})^{\omega}\ar[r]^{c_{\psi_{M}}}\ar[d]^{\mathrm{Eis}_{P^{-}!}} & \mathrm{QCoh}(\mathrm{Par}_{M})\ar[d]^{\Eis_{P}^{\spec,\coarse}}\\
\D(\mathrm{Bun}_{G})^{\omega}\ar[r]^{c_{\psi}} & \mathrm{QCoh}(\mathrm{Par}_{G})
}
\]
i.e. there is an equivalence $c_{\psi} \circ \Eis_{P^{-}} \simeq \Eis_{P}^{\spec,\coarse} \circ c_{\psi_M}$ of functors $\D(\Bun_M)^{\omega} \to \QCoh(\Par_G)$.
\end{conjecture}

Here $\Eis_{P}^{\spec,\coarse}$ is the naive push-pull functor on quasicoherent sheaves induced by the diagram $\Par_{M} \leftarrow \Par_P \rightarrow \Par_G$, while $\Eis_{P^{-}!}$ is the geometric Eisenstein functor defined in \cite{HHS}. Informally, this conjecture asserts that $c_{\psi}$ is compatible with parabolic induction. \textbf{In the remainder of this section, we will assume the truth of Conjecture \ref{conj:Eiscpsiintro} unless stated otherwise.}

Our first main result is that for well-understood groups, this compatibility already allows us to construct the (only) correct candidate for the functor $\mathbf{L}_{\psi}$. To explain this, we observe the following general result, whose proof is a simple variant of the discussion preceding Conjecture \ref{conj:CLLCintro}.

\begin{prop}\label{prop:Lpsiexistsintro}For a quasisplit $G$ with a Whittaker datum $(B,\psi)$, the following are equivalent.
\begin{propenum}
\item There is a colimit-preserving functor $\mathbf{L}_{\psi}:\D(\Bun_G) \to \IndCoh(\Par_G)$ which preserves compact objects and such that the diagram \[
\xymatrix{\D(\mathrm{Bun}_{G})\ar[dr]^{c_{\psi}}\ar[r]^{\mathbf{L}_{\psi}} & \mathrm{IndCoh}(\mathrm{Par}_{G})\ar[d]^{\Psi}\\
& \mathrm{QCoh}(\mathrm{Par}_{G})
}
\]
commutes.
\item The functor $c_{\psi}$ carries $D(\Bun_G)^{\omega}$ into $\Coh^\qc(\Par_G)$.
\end{propenum}
If $\mathbf{L}_{\psi}$ satisfying (i) exists, it necessarily agrees with the ind-completion of $c_{\psi}:\D(\Bun_G)^{\omega} \to \Coh^\qc(\Par_G)$. In particular, such a functor $\mathbf{L}_{\psi}$ is automatically unique and $\QCoh(\Par_G)$-linear.
\end{prop}

If it exists, we call the (unique) functor $\mathbf{L}_{\psi}$ satisfying (i) the Langlands functor. By this proposition, proving the existence of the Langlands functor is equivalent to showing a strong finiteness property of $c_{\psi}$. Unfortunately, this finiteness seems completely out of reach in general.

\begin{thm} \label{rslt:intro-properties-of-L-psi}
For a well-understood group $G$ with a Whittaker datum $(B,\psi)$, the Langlands functor $\mathbf{L}_{\psi}$ exists. Moreover, it is compatible with parabolic induction in the sense that for any standard Levi $M \subset G$, there is a commutative diagram
\[
\xymatrix{\D(\mathrm{Bun}_{M})\ar[r]^{\mathbf{L}_{\psi_{M}}}\ar[d]^{\mathrm{Eis}_{P^{-}!}} & \mathrm{IndCoh}(\mathrm{Par}_{M})\ar[d]^{\Eis_{P}^{\spec}}\\
\D(\mathrm{Bun}_{G})\ar[r]^{\mathbf{L}_{\psi}} & \mathrm{IndCoh}(\mathrm{Par}_{G})
}
\]
\end{thm}
Here the second part is well-posed by noting that tautologically, all standard Levis in a well-understood group are well-understood. The functor $\Eis_{P}^{\spec}$ is the true spectral Eisenstein functor, given by a suitalble push-pull on ind-coherent sheaves.

To prove this, we verify the condition of Proposition \ref{prop:Lpsiexistsintro}.(ii) by induction on the semisimple rank. A key observation is that for any compact $A \in \D(\Bun_G)^\omega$ whose spectral support lies in the open-closed locus of irreducible (i.e., supercuspidal) $L$-parameters, we can prove by a direct computation that $c_{\psi}(A) \in \Coh^\qc(\Par_G)$ (Theorem \ref{thm:cpsisupercuspidalcoherent}). This result holds unconditionally for any group. For well-understood $G$, we are able to prove a complementary result: if $A$ is compact with spectral support in the complement of the irreducible locus, then $A$ can be obtained from sheaves of the form $\Eis_{P^{-}!}B$ for some proper parabolics $P=MU\subset G$ and compact sheaves $B \in \D(\Bun_M)^\omega$ by finitely many cones and retracts. This uses the full extent of the definition of well-understood groups, along with the main results of \cite{HKW} and the orthogonal decomposition of $\D(\Bun_G)$ into cuspidal and Eisenstein subcategories proved in \cite[Theorem 1.3.2.(iii)]{HHS}. In fact, what we really prove (Theorem \ref{thm:wellunderstoodSameOrthogonal}) is that the two natural orthogonal decompositions of $\D(\Bun_G)$ agree, as conjectured in \cite[\S 1.3]{HHS}. Since $\Eis_{P}^{\spec,\coarse}$ preserves $\Coh^\qc$ and $c_{\psi_M}$ carries compact sheaves towards $\Coh^\qc(\Par_M)$ by the induction hypothesis, we get the desired result. The compatibility of the Langlands functor with parabolic induction is then a formal consequence of the assumed compatibility for $c_{\psi}$ on compact objects.

To proceed further, it would be very useful to have an adjoint pair of functors. As a first observation, note that if the Langlands functor $\mathbf{L}_{\psi}$ exists, it automatically has a right adjoint $\mathbf{R}_{\psi}$ by the usual adjoint functor theorems. For easy formal reasons, $\mathbf{R}_{\psi}$ turns out to be colimit-preserving and $\QCoh(\Par_G)$-linear, carrying admissible sheaves towards ULA sheaves, and the compatibility of $\mathbf{L}_{\psi}$ with Eisenstein series induces a compatibility of $\mathbf{R}_{\psi}$ with constant term functors. 

Beyond these general properties, the right adjoint is not very explicit or easy to compute with. However, we have two more key pieces of leverage. Firstly, it is more or less a formal consequence of Theorem \ref{thm:ctadjunctionintro} that $\mathbf{R}_{\psi}|_{\Coh(\Par_G)_{\fin}} \cong t_{\psi} $. This is already very useful, as it identifies the inexplicit right adjoint with the very explicit functor $t_{\psi}$, at least after restricting to $\Coh(\Par_G)_{\fin}$. Even more importantly, we also have the following deep result, stated also in \cref{rslt:intro-Rpsi-preserves-compact-objects} above.

\begin{thm}\label{thm:rightadjcompactintro}For a well-understood group $G$ with a Whittaker datum $(B,\psi)$, the right adjoint functor $\mathbf{R}_{\psi}$ preserves compact objects.
\end{thm}
This is one of the hardest theorems in this paper, and we will discuss its proof in Section \ref{ss:spectralintro} below.

In any case, we now have a pair of adjoint functors 
\[\mathbf{L}_{\psi}:\D(\Bun_G) \rightleftarrows \IndCoh(\Par_G) : \mathbf{R}_{\psi}\]
between compactly generated categories which are both colimit-preserving and $\QCoh(\Par_G)$-linear, and both preserve compact objects. To conclude the full conjecture, it would be enough to see that $\mathbf{R}_{\psi}$ is fully faithful and that $\mathbf{L}_{\psi}$ is conservative on compact objects. Towards the first of these goals, we have the following crucial reduction.

\begin{thm}\label{thm:aRffintro}Let $G$ be a quasisplit group with a Whittaker datum $(B,\psi)$.
\begin{thmenum}
    \item If $G$ is reasonable, then $a_{\psi}$ is fully faithful if and only if $c_{\psi}(i_{1!}W_{\psi})$ is a line bundle on $\Par_G$.

    \item If $G$ is well-understood and the equivalent conditions of (i) hold, then $\mathbf{R}_{\psi}$ is fully faithful.
\end{thmenum}
\end{thm}
The first part here is totally formal (Proposition \ref{prop:fullyfaithfuleasy}), and doesn't require compatibility with parabolic induction. However, the second part is much more interesting, since after all we constructed $\mathbf{R}_{\psi}$ from $a_{\psi}$ in a very indirect manner (by passing to the right adjoint $c_{\psi}$, carefully lifting this along $\Psi$ to get $\mathbf{L}_{\psi}$, and then passing to another right adjoint). The argument here crucially relies on the functor $t_{\psi}$, and in fact, we separately prove the implications ``$a_{\psi}$ fully faithful $\Rightarrow$ $t_{\psi}$ fully faithful'' (Theorem \ref{thm:atff}) and ``$t_{\psi}$ fully faithful $\Rightarrow$ $\mathbf{R}_{\psi}$ fully faithful'' (Theorem \ref{thm:tRff}). The first of these follows quickly from the good properties of $a_{\psi}$ and $t_{\psi}$ established so far, but the second is more subtle, and the argument here borrows some ideas from ``geometric Langlands with restricted variation'' as in \cite{AGKRRV}. 

Let us highlight one pleasant aspect of this argument, which works unconditionally for all groups. Let $\mathcal{O}^{\restr} \in \QCoh(\Par_G)$ be the direct sum, over all closed points $\varphi \in X_{G}^{\spec}$, of the local cohomology sheaves of $\mathcal{O}_{\Par_G}$ along the fibers $q^{-1}(\varphi)$. It is easy to see that there is a natural map $\mathcal{O}^{\restr}\to \mathcal{O}_{\Par_G}$ inducing an idempotent coalgebra structure on $\mathcal{O}^{\restr}$. If $\mathcal{C}$ is any presentable stable $\infty$-category with a tensor action of $\QCoh(\Par_G)$, we write $\mathcal{C}^{\restr}$ for the full subcategory of ``restricted'' objects, which by definition is spanned by objects $A$ such that the natural map $\mathcal{O}^{\restr} \otimes A \to A$ is an isomorphism. Note that any restricted object $A$ has a canonical direct sum decomposition $A=\oplus_{\varphi} A_{\varphi}$ induced by the direct sum structure on $\mathcal{O}^{\restr}$. In general $\mathcal{C}^{\restr}$ is stable under colimits, and under mild conditions the inclusion $\mathcal{C}^{\restr} \subset \mathcal{C}$ has a colimit-preserving right adjoint given by $A\mapsto \mathcal{O}^{\restr} \otimes A$ (Proposition \ref{prop:restrnaiveadjoint}). The following result identifies the restricted objects in $\IndCoh(\Par_G)$ and $\D(\Bun_G)$ in very explicit terms.

\begin{thm}[Proposition \ref{prop:restrsheavesParG}, Theorem \ref{thm:restrsheavesBunG}, Theorem \ref{thm:ULAisrestricted}] Fix any reductive group $G$.
\begin{thmenum}
\item There is a natural identification
\[\IndCoh(\Par_G)^{\restr} = \Ind(\Coh(\Par_G)_{\fin})\]
as full subcategories of $\IndCoh(\Par_G)$.

\item There is a natural identification
\[\D(\Bun_G)^{\restr} = \Ind(\D(\Bun_G)_{\fin}) \]
as full subcategories of $\D(\Bun_G)$, and any ULA sheaf on $\Bun_G$ is restricted.
\end{thmenum}
In both cases, the right adjoint to the inclusion $\mathcal{C}^{\restr} \subset \mathcal{C}$ is given by the action of $\mathcal{O}^{\restr}$.
\end{thm}

The last statement here combined with the $\QCoh(\Par_G)$-linearity of the functors $\mathbf{L}_{\psi}$ and $\mathbf{R}_{\psi}$ implies that both functors are compatible with passage to and from the restricted subcategories. We also prove that in both cases, passing to restricted sheaves is not too lossy: if $A$ is a compact object, then $\mathcal{O}^{\restr} \otimes A$ is nonzero (Lemma \ref{lem:conservativetrick}). From here it is not very difficult to conclude the proof of Theorem \ref{thm:aRffintro}.(ii). Along the way, we also show that Conjecture \ref{conj:CLLCrestrictedintro} implies Conjecture \ref{conj:CLLCintro} for well-understood groups.

This analysis reduces the full conjecture to two properties of $c_{\psi}$: we need to show that $c_{\psi}(i_{1!}W_{\psi})$ is a line bundle, and that $c_{\psi}$ is conservative on compact sheaves. Unfortunately, both of these turn out to be extremely difficult questions. Let us discuss the conservativity first. Unwinding the definition of $c_{\psi}$, it is not hard to prove that conservativity on compact objects boils down to the following concrete problem.
\begin{conjecture}For varying irreducible $V_{\mu} \in \Rep \hat{G}$, the functors
\[c_{\psi}^{\mu}(-) \overset{\mathrm{def}}{=}\RHom(T_{V_{\mu}} i_{1!}W_{\psi},-) = \RHom(W_{\psi},i_{1}^{\ast}T_{V_{\mu^{-1}}}(-)):\D(\Bun_G)^\omega \to \D(\Qellbar)\]
give a conservative family.
\end{conjecture} 
In classical geometric Langlands, the exact analogue of this result is a recent breakthrough of Faergeman-Raskin \cite{FR}. While one can prove an analogue of the geometric Casselman-Shalika formula in our setting (which we pursue in our upcoming work with Hamann), this alone doesn't seem to help much, as the microlocal methods of \cite{FR} don't adapt at all.

Instead, we attack the conservativity of $c_{\psi}$ using the semiorthogonal decomposition of $\D(\Bun_G)^{\omega}$. Ideally, if we knew that $\Coh^\qc(\Par_G)$ carried a $B(G)$-indexed semiorthogonal decomposition compatible via $c_{\psi}$ with the canonical decomposition of $\D(\Bun_G)^\omega$, it would be enough to prove that the individual functors $c_{\psi}i_{b!}^{\ren}$ are conservative on compact representations of $G_b(F)$. We could then try to argue by induction on the semisimple rank of $G$, exploiting the fact that $i_{b!}^{\ren}$ is a fragment of an Eisenstein series for non-basic $b$, in conjunction with the assumed Eisenstein compatibility of $c_{\psi}$. Unfortunately, it seems extremely hard to construct the desired semiorthogonal decomposition on the spectral side directly. 

However, very fortunately for us, Bertoloni Meli-Koshikawa \cite{BMK} have proved some partial results towards such a decomposition. In particular, they define a canonical category of \emph{cuspidal} coherent sheaves on $\Par_G$, characterized by the vanishing of specific graded summands of spectral constant term functors, which (under the conjecture) match those sheaves on $\Bun_G$ which are $!$-extensions of compact supercuspidal representations from basic points $b \in B(G)_{\mathrm{bas}}$. As a consequence of their very detailed analysis of this category, they prove a deep vanishing theorem for $\RHom$'s between certain spectral Eisenstein series of cuspidal coherent sheaves (Theorem \ref{thm:BMKsemiorthogonalKey}), which is a spectral analogue of certain trivial vanishing theorems for groups of the form $\RHom(i_{b!}A,i_{b' !}A')$. Using their result together with some further tricks and ideas, we are (barely) able to push through an argument along the lines suggested in the previous paragraph. Ultimately this leads to the following theorem (stated as \cref{thm:intro-dream-induction-result} above), which is our main general result for well-understood groups.
\begin{thm}[Theorem \ref{thm:DreamInductionOnLevisTheorem}]Let $G$ be a well-understood group with a Whittaker datum $(B,\psi)$ and maximal torus $T \subset B$. Assume the following conditions.
\begin{thmenum}

    \item For each proper standard Levi $M \subsetneq G$, the functor $\mathbf{L}_{\psi_M}$ is an equivalence of categories. 

    \item The functor $c_{\psi} i_{1!}$ sends $W_{\psi}$ to a line bundle.

    \item The functor $c_{\psi} i_{1!}$ is conservative on compact supercuspidal representations and sends them towards cuspidal coherent sheaves.
\end{thmenum}
Then $\mathbf{L}_{\psi}$ is an equivalence of categories.
\end{thm}
We note that in (ii) and (iii), it is immaterial whether one argues with $c_{\psi}$ or $\mathbf{L}_{\psi}$, as they give the same functor on compact sheaves.

Let us comment on the conditions in this theorem. Since standard Levis in well-understood groups are well-understood tautologically, (i) is well-posed. Moreover, from the perspective of actually proving categorical Langlands, arguing by induction on Levis is an extremely natural technique. The condition (ii) is exactly what we need to conclude that $\mathbf{R}_{\psi}$ is fully faithful, as in Theorem \ref{thm:aRffintro} and the discussion thereafter. Finally, condition (iii) seems to be the bare minimum which (together with (i)) lets us prove the conservativity of $\mathbf{L}_{\psi}$ on compact sheaves by an argument based on the semiorthogonal decomposition as discussed above.  For more details, we refer the reader to the proof of Theorem \ref{thm:DreamInductionOnLevisTheorem}.

In any case, this theorem reduces the categorical local Langlands conjecture for well-understood groups to the same conjecture for proper Levis, together with two conjectural properties of the functor $c_{\psi}i_{1!}:\D(G(F),\Qellbar) \to \QCoh(\Par_G)$. Our next goal is to discuss what we can prove towards these properties.

\subsection{General linear groups}\label{sec:GLnintro}

We maintain the notations of the previous section, so $G$ is well-understood and equipped with a Whittaker datum $(B,\psi)$. We continue to assume Conjecture \ref{conj:Eiscpsiintro}. Our next task is to get a handle on the functor $c_{\psi}i_{1!}$. Note that compared to the full functor $c_{\psi}$, the domain of this functor is of much more classical nature: it is simply the derived category of smooth representations of $G(F)$. We remind the reader that this category admits a canonical decomposition into mutually orthogonal indecomposable subcategories, usually called Bernstein components or Bernstein blocks. These are indexed by the set $\mathfrak{B}(G)$ of primitive idemponents in the Bernstein center $\mathfrak{Z}_G$, and this decomposition is compatible with (and in fact refines) the decomposition induced by the set of primitive idempotents $\mathfrak{B}^{\spec}(G)$ in the ring $\mathfrak{Z}_{G}^{\spec}$ under the Fargues-Scholze map $\mathfrak{Z}_{G}^{\spec} \to \mathfrak{Z}_{G}$. Note that $\mathfrak{B}^{\spec}(G)$ is simply the set of connected components of $\Par_G$. Since $c_{\psi}i_{1!}$ is $\mathfrak{Z}^{\spec}_{G}$-linear via the Fargues-Scholze map, we are free to analyze this functor on individual spectral Bernstein components.

We first consider the case of $\mathrm{GL}_n$, or more generally \emph{groups of GL-type}, i.e. groups the form $G=\mathrm{GL}_{n_1}\times \cdots \times \mathrm{GL}_{n_i}$. This is the minimal class of groups containing $\GL_n$ and stable under passing to Levis. Note that the Whittaker datum for such a group is unique. We also have the important simplifying feature that for any group of GL-type the Fargues-Scholze map $\mathfrak{Z}_{G}^{\spec}\to \mathfrak{Z}_{G}$ is an isomorphism; this is a theorem of Helm. In what follows we will freely identify these two rings along the Fargues-Scholze map without further comment. 

Moreover, for groups of GL-type, there is already a known functor with the same domain and target as $c_{\psi}i_{1!}$, and with excellent properties. Precisely, we will need the following major result from \cite{BCHN}.

\begin{thm}[Ben--Zvi-Chen-Helm-Nadler]\label{thm:BCHNintro} For $G$ any group of GL-type, there is a natural fully faithful functor
\[\mathscr{S}_{G}: \D(G(F),\Qellbar)^{\omega} \to \Coh^\qc(\Par_G) \]
which is linear over the Bernstein center and compatible with parabolic induction. Moreover, its ind-completion sends $W_{\psi}$ to the structure sheaf $\mathcal{O}_{\Par_G}$.
\end{thm}

Unlike the horribly inexplicit functor we are reckoning with, the BCHN functor $\mathscr{S}_G$ is defined very explicitly. More precisely, for each Bernstein component $\mathfrak{s}$, Bernstein theory gives a canonical compact projective generator $\Pi_{\mathfrak{s}}$ for this component, whose endomorphism ring $\mathrm{End}(\Pi_{\mathfrak{s}})$ is a tensor product of type A Iwahori-Hecke algebras. On the other side, BCHN construct a coherent sheaf $\mathcal{S}_{\mathfrak{s}} \in \Coh^\qc(\Par_G)$ with two miraculous properties:

\begin{enumerate}
    \item There is an isomorphism $\mathrm{End}(\Pi_{\mathfrak{s}}) \simeq \mathrm{End}(\mathcal{S}_{\mathfrak{s}})$ of $\mathfrak{Z}_{G,\mathfrak{s}}$-algebras.

    \item For any $n\neq 0$, $\Hom(\mathcal{S}_{\mathfrak{s}},\mathcal{S}_{\mathfrak{s}}[n])=0$.
\end{enumerate}
Choosing isomorphisms as in (1), it is then formal nonsense that there is a unique $\mathfrak{Z}_{G}$-linear fully faithful functor $\mathscr{S}_G$ sending $\Pi_{\mathfrak{s}}$ to $\mathcal{S}_{\mathfrak{s}}$ and such that the induced map on endomorphism rings recovers the isomorphisms (1). The key connection with our strategy is the following result.

\begin{thm}[Theorem \ref{thm:GLnfullfaithndconsequences}] \label{thm:springervscpsiintro}For $G$ any group of GL-type, there is an equivalence of functors $c_{\psi}i_{1!} \simeq \mathscr{S}_G : \D(G(F),\Qellbar)^{\omega} \to \Coh^\qc(\Par_G)$.
\end{thm}

The essential gain here is that the functor $\mathscr{S}_G$ is extremely well-understood by Theorem \ref{thm:BCHNintro}, so this theorem immediately gives us the same knowledge for $c_{\psi}i_{1!}$. Combining this with a small further argument using Theorem \ref{thm:DreamInductionOnLevisTheorem}, we finally get the following theorem (stated as \cref{thm:intro-main-result} above).

\begin{thm}\label{thm:CLLCforGLnintro}For $G$ any group of GL-type, the functor $\mathbf{L}_{\psi}$ is an equivalence of categories.
\end{thm}

To prove Theorem \ref{thm:springervscpsiintro}, it's enough to check that $c_{\psi}i_{1!}\Pi_{\mathfrak{s}} \simeq \mathcal{S}_{\mathfrak{s}}$ and that the induced map on endomorphism rings is an isomorphism. For the first part, we exploit the fact that both $\Pi_{\mathfrak{s}}$ and $\mathcal{S}_{\mathfrak{s}}$ are parabolic inductions, together with the assumed compatibility of $c_{\psi}$ with parabolic induction. For the claim on endomorphism rings, we note that the source and target are abstractly isomorphic torsion-free $\mathfrak{Z}_{G,\mathfrak{s}}$-algebras. We use this along with a trick of commutative algebra to reduce to showing \emph{injectivity} of the given ring map. This is much easier, as it can be checked after specialization at any dense set of points in the Bernstein center. We now argue at the dense set of \emph{generous} points $\varphi \in \Spec \mathfrak{Z}_{G,\mathfrak{s}}$, where life becomes much simpler: if $\varphi$ is generous, there is a unique irreducible $G(F)$-representation $\pi_{\varphi}$ with semisimple parameter $\varphi$, and we have isomorphisms $\Pi_{\mathfrak{s}} \otimes \mathfrak{Z}_{G,\mathfrak{s}}/\mathfrak{m}_{\varphi} \simeq \pi_{\varphi}^{d}$ and $\mathcal{S}_{\mathfrak{s}} \otimes \mathfrak{Z}_{G,\mathfrak{s}}/\mathfrak{m}_{\varphi} \simeq k_{\varphi}^{d}$. Here $k_{\varphi}$ is the skyscraper at (the unique lift) $\varphi \in \Par_{G}$ and $d$ is some multiplicity which depends only on $\mathfrak{s}$. The essential point is then that $c_{\psi}(i_{1!}\pi_{\varphi}) \simeq k_{\varphi}$ by an explicit computation, again using the compatibility of $c_{\psi}$ with parabolic induction. For the details, we refer the reader to Sections \ref{sec:IwahoriGLn}-\ref{generalGLn}.

\subsection{Classical groups} \label{sec:intro.classical-groups}
In applying Theorem \ref{thm:DreamInductionOnLevisTheorem} to groups of GL-type, we took advantage of several special features of these groups: every Bernstein block is Whittaker generic, and any supercuspidal $G(F)$-representation automatically goes to a cuspidal coherent sheaf. However, to access condition (ii) in Theorem \ref{thm:DreamInductionOnLevisTheorem}, we don't need to understand $c_{\psi}i_{1!}$ on all Bernstein blocks: we just need to understand it on $\psi$-generic Bernstein blocks. Note that in general, full faithfulness in the categorical conjecture predicts that the composite ring map
\[\mathfrak{Z}_{G}^{\spec} \to \mathfrak{Z}_{G} \to \End(W_{\psi})\]
is an isomorphism. If this map were an isomorphism, in conjunction with Theorem \ref{thm:Whittakerbasics}.(iii) it would give a preferred section of the map $\mathfrak{B}(G) \to \mathfrak{B}^{\spec}(G)$ identifying the set of spectral Bernstein components with the set of $\psi$-generic Bernstein components, and inducing compatible identifications $\mathfrak{Z}_{G,\mathfrak{s}}^{\spec} = \mathfrak{Z}_{G,\mathfrak{s}}$ for all such $\mathfrak{s}$. 

Motivated by this observation, we say a well-understood group $G$ is \emph{very well-understood} if for $G$ and all its standard Levis, the ring map above is an isomorphism, and the component groups of centralizers of supercuspidal $L$-parameters are abelian. This condition includes most groups which are currently known to be well-understood \cite{DHKM2}, and in particular includes quasisplit classical groups of types A and B, as well as $\mathrm{GSp}_4$.\footnote{Quasisplit classical groups of type D should also be within reach, but strictly speaking they are not covered by \cite{DHKM2}.}

For very well-understood groups, we can prove several things. Firstly, we prove the full conjecture over the irreducible locus.
\begin{thm}[Theorem \ref{thm:CLLCoverirredparsVWE}] Let $G$ be a very well-understood group with a fixed Whittaker datum. Then $c_{\psi}$ induces an equivalence
\[\D(\Bun_G)^{\cusp} \overset{\sim}{\to} \QCoh(\Par_{G}^{\mathrm{irr}}).\]
\end{thm}
Unsurprisingly, this does not require the compatibility of $c_{\psi}$ with parabolic induction. Instead, we leverage the very well-understood condition via some tricks based on the main results of \cite{HJnote} and \cite{HKW}.
 
Using this result, we then formulate a precise conjectural analogue of Theorem \ref{thm:BCHNintro} for the $\psi$-generic Bernstein components of $G$, giving an explicit recipe for a functor
\[\mathscr{S}_{G,\psi}:\D(G(F),\Qellbar)^{\omega}_{\psi-\mathrm{gen}} \to \Coh^\qc(\Par_G)\]
defined similarly as in the case of general linear groups. We refer the reader to Conjecture \ref{conj:BCHNpsigenericdreamconjecture} for the exact statement. We then prove the following result, by adapting the arguments from the previous section.
\begin{thm}If $G$ is very well-understood with a fixed Whittaker datum and Conjecture \ref{conj:BCHNpsigenericdreamconjecture} is true, then $c_{\psi}i_{1!}|_{\D(G(F),\Qellbar)^{\omega}_{\psi-\mathrm{gen}}} \simeq \mathscr{S}_{G,\psi}$ and $c_{\psi}(i_{1!}W_{\psi})$ is a line bundle.
\end{thm}

For many groups which are known to be very well-understood, Conjecture \ref{conj:BCHNpsigenericdreamconjecture} is in fact a work in progress of Helm-Solleveld-Xu. As such, the above result gives a viable approach to computing $c_{\psi}(i_{1!}W_{\psi})$ for these groups and proving that it is a line bundle, which should be unconditional in the very near future.

At this point, we've verified the first two conditions of Theorem \ref{thm:DreamInductionOnLevisTheorem} for very well-understood groups. It remains to understand the final condition in that theorem. However, note that the above analogue of the BCHN functor is no longer enough, as there are plenty of supercuspidals which are not $\psi$-generic for any Whittaker datum. There is also the additional severe complication that supercuspidal representations no longer have supercuspidal $L$-parameters in general.

Instead, we argue with the right adjoint. For simplicity, we assume for this discussion that $G$ is semisimple. Then the supercuspidal part of $\D(G(F),\Qellbar)^{\omega}$ is just a direct sum of copies of $\Perf(\Qellbar)$, indexed by pairs $(\phi,\rho)$ where $\phi$ is a discrete $L$-parameter and $\rho$ is a cuspidal representation of the extended component group. In particular, there is a grading by discrete $L$-parameters, and the $\phi$-graded summand is a direct sum of a specific finite number of copies of $\Perf(\Qellbar)$, i.e. there is an equivalence $\D(G(F),\Qellbar)^{\omega}_{\phi}\simeq \Perf(\Qellbar)^{\oplus m_{\phi}^{\mathrm{aut}}}$ for some multiplicity $m_{\phi}^{\mathrm{aut}}\geq 0$. Note that $m_{\phi}^{\mathrm{aut}}$ is just the number of supercuspidal representations in the $L$-packet $\Pi_{\phi}(G)$.

Now we consider the other side. Here again we use some deep results of Bertoloni Meli-Koshikawa \cite{BMK}. As mentioned above, they have defined an intrinsic subcategory $\Coh(\Par_G)_{\cusp}$ of cuspidal coherent sheaves. We write $\Coh(\Par_G)^0$ for the direct summand spanned by sheaves on which $Z(\hat{G})^{W_F}$ acts trivially. If $G$ is semisimple, then Bertoloni Meli-Koshikawa prove unconditionally that $\Coh(\Par_G)_{\cusp}^{0}$ admits a grading by discrete $L$-parameters $\phi$, and for each graded summand there is an equivalence $\Coh(\Par_G)^{0}_{\cusp,\phi} \simeq \Perf(\Qellbar)^{\oplus m_{\phi}^{\spec}}$ for some multiplicity $m_{\phi}^{\spec}\geq 0$.  

Recall our goal is to see that $\mathbf{L}_{\psi} i_{1!}$ is conservative on compact supercuspidals and sends them towards cuspidal coherent sheaves. The key idea is that by its assumed constant term compatibility, the right adjoint $\mathbf{R}_{\psi}$ automatically carries $\Coh^\qc(\Par_G)_{\cusp}^{0}$ towards objects of the form $i_{1!}\pi$ where $\pi$ is compact and supercuspidal. Moreover, by our analysis above, we can assume that $\mathbf{R}_{\psi}$ is fully faithful. Combining these observations, we see that $\mathbf{R}_{\psi}$ induces a \emph{fully faithful} functor
\[i_{1}^{\ast} \mathbf{R}_{\psi}:\Coh^\qc(\Par_G)_{\cusp}^{0} \to \D(G(F),\Qellbar)^{\omega,\mathrm{sc}}\]
which is automatically compatible with the gradings by discrete parameters. Restricting to each $\phi$-graded summand, we get a fully faithful functor
\[i_{1}^{\ast} \mathbf{R}_{\psi}:\Perf(\Qellbar)^{\oplus m_{\phi}^{\spec}}\simeq \Coh(\Par_G)_{\cusp,\phi}^{0} \to \D(G(F),\Qellbar)_{\phi}^{\omega}\simeq \Perf(\Qellbar)^{\oplus m_{\phi}^{\mathrm{aut}}}.\]
The essential observation is that if we knew these functors were equivalences for all $\phi$, the full left adjoint $\mathbf{L}_{\psi}i_{1!}$ of $i_{1}^{\ast}\mathbf{R}_{\psi}$ would then restrict to an inverse equivalence from $\D(G(F),\Qellbar)^{\omega,\mathrm{sc}}$ towards $\Coh^\qc(\Par_{G})_{\cusp}^{0}$, which is more than enough to verify the final condition in Theorem \ref{thm:DreamInductionOnLevisTheorem}. However, by the subsequent (trivial) lemma, this boils down to a concrete counting problem.

\begin{lem}Let $F:\Perf(k)^{\oplus n} \to \Perf(k)^{\oplus m}$ be a fully faithful functor. Then $n\leq m$, and $F$ is an equivalence if and only if $n=m$.
\end{lem}

This reduces us to proving that $m_{\phi}^{\spec}\geq m_{\phi}^{\mathrm{aut}}$ for all discrete parameters. But this is an entirely manageable task! More precisely, for quasisplit classical groups Moeglin's work gives a completely elementary and concrete formula for $m_{\phi}^{\mathrm{aut}}$ in all cases, while Bertoloni Meli-Koshikawa are independently able to compute $m_{\phi}^{\spec}$ for classical groups, and they have assured us that the formulas match: $m_{\phi}^{\spec} = m_{\phi}^{\mathrm{aut}}$ for all discrete parameters of all quasisplit classical groups. Putting this whole discussion together, and giving some additional arguments to deal with the case of general reductive $G$, we finally prove the following theorem.

\begin{thm}[Theorem \ref{thm:ultimateVWO}]Let $G$ be a very well-understood group with a fixed Whittaker datum. Suppose that for $G$ and all its standard Levi subgroups, Conjecture \ref{conj:BCHNpsigenericdreamconjecture} is true for all $\psi$-generic Bernstein components and $m_{\phi}^{\mathrm{aut}} = m_{\phi}^{\spec}$ for all discrete $L$-parameters $\phi$. Then $\mathbf{L}_{\psi}$ is an equivalence of categories.
\end{thm}

In combination with the forthcoming work of Helm-Solleveld-Xu mentioned above, this Theorem implies the categorical local Langlands conjecture for all unramified odd special orthogonal and unitary groups.

\subsection{Spectral results}\label{ss:spectralintro}

In the discussion above, we punted on one crucial item, namely the proof of Theorem \ref{thm:rightadjcompactintro}. This requires several results on the spectral side. Our starting point is the following compactness criterion:

\begin{thm}[Theorem \ref{thm:compactcriterion}] Let $G$ be a well-understood group. Then a sheaf $A \in \D(\Bun_G)$ is compact if and only if it has quasicompact support, quasicompact spectral support, and for all compact sheaves $B \in \D(\Bun_G)^{\omega}$, the complex $\RHom(B,A)$ is bounded and its cohomology groups are finitely generated $\mathfrak{Z}_{G}^{\mathrm{spec}}$-modules.
\end{thm}

We want to apply this to $A=\mathbf{R}_{\psi}(\mathcal{F})$ for some coherent $\mathcal{F}$. Here the condition on the spectral support is trivial, but the first and third conditions are highly nontrivial. However, using the defining adjunction of $\mathbf{R}_{\psi}$, the third condition is quickly reduced to the first part of the following theorem.

\begin{thm}\label{thm:spectral-main}Let $G$ be any quasisplit group.

\begin{thmenum}
    \item For any $\mathcal{F},\mathcal{G} \in \Coh^\qc(\Par_G)$, $\RHom(\mathcal{F},\mathcal{G})$ is a bounded complex of $\Qellbar$-vector spaces whose cohomology groups are finitely generated $\mathfrak{Z}_{G}^{\mathrm{spec}}$-modules.

    \item We have a coincidence of full subcategories $\Adm(\Par_G) \cap \Coh^\qc(\Par_G) = \Coh(\Par_G)_{\fin}$ inside $\IndCoh(\Par_G)$.

\end{thmenum}
\end{thm}

Note that (i) here is a strengthening of Theorem \ref{thm:spectralboundedintro} above, and (ii) follows easily from (i) by the discussion preceding that result. Part (i) is also the exact analogue of Bernstein's classic theorem that for any compact $A,B \in \D(G(F),\Qellbar)^{\omega}$, $\RHom(A,B)$ is a bounded complex whose cohomology groups are finitely generated $\mathfrak{Z}_{G}$-modules.

Using that $q:\Par_{G} \to X_{G}^{\spec}$ is a good moduli space, it is not hard to prove that the complex in (i) is left-bounded with individual cohomologies finitely generated over $\mathfrak{Z}_{G}^{\spec}$. The essential difficulty is to prove it is right-bounded. For this we argue by induction on the semisimple rank of $G$. The first ingredient which makes the induction work is the following result of Takaya \cite{Takaya}.

\begin{thm}[Takaya] For any $G$, $\Coh^\qc(\Par_G)$ is generated under finite colimits and retracts by objects of $\Perf^{\qc}(\Par_G^{\mathrm{irr}})$ together with objects of the form $\Eis_{P}^{\mathrm{spec}}(\mathcal{G})$ for $\mathcal{G} \in \Coh^\qc(\Par_M)$ with $P=MU \subsetneq G$ varying over all proper standard parabolics. 
\end{thm}

 By an easy induction, in the second part of this theorem we could equally well assume that $\mathcal{G}$ lies in $\Perf^{\qc}(\Par_{M}^{\mathrm{irr}})$. 
 
 Now when $\mathcal{G}$ is perfect with quasicompact support and $\mathcal{F}$ is coherent it is easy to see that $\RHom(\mathcal{G},\mathcal{F})$ is bounded, so Takaya's theorem reduces us to understanding such RHoms when $\mathcal{G}$ is a spectral Eisenstein series. By adjunction, this can be rewritten as a similar RHom on $\Par_M$ where $\mathcal{F}$ is a spectral constant term, so we now conclude by the first part of the following theorem. This theorem makes use of the central grading on $\IndCoh(\Par_M)$, for which we refer the reader to \cref{def:central-grading-on-ParG}.

\begin{thm}\label{thm:CTmagic}For any $P=MU \subset G$ and any $\chi \in X^\ast(Z(\hat{M})^{W_F})$, the functor $\mathrm{CT}_{P}^{\mathrm{spec},\chi}$ defined as the composition
\[\IndCoh(\Par_G) \overset{\mathrm{CT}_{P}^{\mathrm{spec}}}{\longrightarrow} \IndCoh(\Par_M) \to \IndCoh(\Par_M)^\chi
\]preserves $\Coh^\qc$.

Moreover, if $\mathcal{F} \in \Coh^\qc(\Par_G)$ is fixed, then for any $P$, $\mathrm{CT}_{P}^{\mathrm{spec},\chi}(\mathcal{F})=0$ for all but finitely many $P$-dominant $\chi$.
\end{thm}

The first part of this theorem is the exact analogue of \cite[Theorem 1.2.1.(2)]{HHS}.\footnote{This is an interesting example of a pair of parallel theorems, matching each other under CLLC, which are \emph{both} difficult to prove.} We also note that the second part of this theorem quickly implies the first condition in the compactness criterion above, using that $\mathbf{R}_{\psi}$ is compatible with constant terms together with the fact that non-basic stalks on $\Bun_G$ are given by (stalks of) specific constant terms. As such, Theorem \ref{thm:CTmagic} finishes the proof of Theorem \ref{thm:rightadjcompactintro}.

The proof of Theorem \ref{thm:CTmagic} proceeds in several steps. First, note that constant terms will kill any sheaf supported on the irreducible locus, so by Takaya's theorem it suffices to prove similar coherence and vanishing statement for the composition $\mathrm{CT}_{P}^{\mathrm{spec},\chi} \Eis_{Q}^{\spec}(\mathcal{G})$ where $Q=LV$ is another proper standard parabolic. In the spirit of the famous geometric lemma of Bernstein-Zelevinsky \cite[Lemma 2.12]{BZinduced}, we then prove an analogous \emph{spectral geometric lemma}, Theorem \ref{rslt:spectral-geometric-lemma}, which equips the functor $\mathrm{CT}_{P}^{\mathrm{spec}} \Eis_{Q}^{\spec}$ with a canonical $W_{M} \backslash W / W_L$-indexed filtration whose graded pieces $F_{w}$ can themselves be explicated as certain pullback-pushforward functors. However, unlike the usual geometric lemma where the graded pieces are of very simple nature, the graded pieces here are still very complicated and each admit their own $\mathbf{N}$-indexed filtrations, and only the graded pieces $F_{w}^{n}$ of these additional filtrations can be written somewhat explicitly. Ultimately, we are able to prove just enough about the functors $F_{w}^{n}$ to complete the proof of Theorem \ref{thm:CTmagic}. We refer to the body of the paper for a precise discussion of the spectral geometric lemma and its application here.

\subsection{Ind-coherent sheaves} \label{sec:introindcoherent}

The attentive reader may have noticed that in our discussion so far, we have only cited one result from \cref{sec:alggeo} (namely, Theorem \ref{rslt:spectral-temperization}), and yet this section is over 80 pages long. In the following we survey the contents of this section and its relation to the rest of the paper.

The goal of \cref{sec:alggeo} is to introduce the tools that we need on the spectral side of the Langlands correspondence, i.e.\ algebraic stacks and ind-coherent sheaves. These tools are used crucially in \cref{sec:spectral}, especially in the proof of the spectral geometric lemma (\cref{rslt:spectral-geometric-lemma}), which lies at the heart of results like \cref{rslt:intro-miracle-bound}. The spectral geometric lemma necessitates the use of \emph{derived} algebraic geometry, because for parabolic subgroups $P \subseteq G$ the stack $\Par_P$ is derived in general. There are several references on derived algebraic geometry, for example \cite{DAG-Toen}, \cite{gaitsgory-rozenblyum-vol1}, \cite{gaitsgory-rozenblyum-vol2}, \cite{DAG-Lurie} and \cite{lurie-SAG}, but we were unsatisfied with the treatment of ind-coherent sheaves in the literature: It relies heavily on the geometry over a characteristic $0$ base, uses complicated $(\infty,2)$-categories without proper definition in its construction, and contains errors (see \cref{rmk:Halpern-Leistner-argument-error} and \cref{rmk:Gaitsgory-error}). Moreover, we did not find a clean formulation of the excision sequence for ind-coherent sheaves in the literature. For these reasons, we decided to build the theory of derived algebraic stacks and ind-coherent sheaves from the ground up, with a focus towards the applications in this paper. We use modern techniques in higher categories and 6-functor formalisms, and in particular provide a completely new construction of the ind-coherent 6-functor formalism that works over any regular noetherian base.

\subsubsection{Basic definitions}

In this paper we work with the functor-of-points perspective and implicitly view all objects in the derived world. In particular, by a \emph{ring} we mean an animated/simplicial ring and by a \emph{module} we mean an object in the derived category of modules. We recall the basic definitions of rings and modules in \cref{sec:alggeo.rings-and-modules}. With these definitions at hand, a \emph{stack} over a ring $\Lambda$ is a functor $X\colon \Alg_\Lambda \to \Ani$ that descends along fppf covers of $\Lambda$-algebras, see \cref{sec:alggeo.stacks}; here $\Ani$ is the category of anima, i.e.\ \enquote{$\infty$-groupoids}. We denote the category of stacks over $\Lambda$ by $\Stk_\Lambda$ and also introduce the full subcategory $\Stk_\Lambda^\lafp \subseteq \Stk_\Lambda$ of \emph{locally almost finitely presented} (\emph{lafp}) stacks over $\Lambda$. For every $\Lambda$-algebra $A$ we denote by $\Spec A \in \Stk_\Lambda$ the functor corepresented by $A$, and call stacks of this form the \emph{affine schemes}. We then denote by $\Sch \subseteq \Stk_\Lambda$ the full subcategory spanned by the \emph{schemes} over $\Lambda$, i.e.\ those stacks that admit an open cover by affine schemes. Every stack $X$ has an underlying topological space $\abs X$ which is obtained by gluing the usual Zariski spectrum on affine schemes, see \cref{def:topology-on-stacks}; this topological space depends only on the underlying classical stack (i.e.\ the restriction of the stack to classical rings). We define a \emph{closed immersion} of stacks to be a map $i\colon Z \to X$ which after pullback to every cover of $X$ by affine schemes is of the form $\Spec B \to \Spec A$ such that the induced map $\pi_0 A \to \pi_0 B$ is surjective. We warn the reader that in the derived world, closed immersions are not injective (see \cref{rmk:closed-immersions-not-injective}).

In \cref{def:geometric-maps-of-stacks} we introduce the notion of \emph{geometric} maps of stacks, which is a higher analog of maps that are representable in Artin stacks. There are robust notions of \emph{smooth}, \emph{flat} and \emph{étale} geometric maps. A stack $X$ over $\Lambda$ is \emph{geometric} if the map $X \to \Spec \Lambda$ is geometric, e.g.\ every scheme and every classical Artin stack is geometric. We call a geometric stack \emph{reduced} if it admits a smooth cover by a reduced (classical) scheme. For every geometric stack $X$ and every closed subset $Z \subseteq \abs X$ there is an induced reduced stack $Z^\red$ together with a closed immersion $Z^\red \to X$ with image $Z$ (see \cref{rslt:reducedness-for-geometric-stacks}).

To every stack $X$ we get an associated symmetric monoidal category $\QCoh(X)$ of \emph{quasicoherent sheaves} on $X$ and to every map $Y \to X$ of stacks there is a symmetric monoidal \emph{pullback functor} $f^*\colon \QCoh(X) \to \QCoh(Y)$ together with a right adjoint $f_*$ (see \cref{sec:alggeo.qcoh}). If $X$ is geometric then $\QCoh(X)$ inherits a t-structure from the usual t-structure on modules over rings (see \cref{rslt:t-structure-on-QCoh}). For the construction of ind-coherent sheaves we make use of the following enhancement, discussed in \cref{sec:alggeo.solid}. To every stack $X$ we associate the category $\QCoh_\solid(X)$ of \emph{solid quasicoherent sheaves} on $X$, which are roughly \enquote{quasicoherent sheaves equipped with a complete topology}. This definition goes back to Clausen--Scholze \cite{condensed-mathematics} and is discussed in several places of the literature. The main advantage of $\QCoh_\solid$ is that it admits a full 6-functor formalism, with $!$-functors for every geometric map of stacks. We refer the reader to \cref{rslt:6ff-on-QCoh-solid} for a summary of the basic properties of this formalism. For every stack $X$ there is a canonical embedding $\QCoh(X) \subseteq \QCoh_\solid(X)$ by equipping a module with the discrete topology. We use the 6-functor formalism on $\QCoh(X)$ to introduce more properties of maps of stacks: We call a geometric map $Y \to X$ \emph{cohomologically smooth} resp.\ \emph{suave} if it is $\QCoh_\solid$-smooth resp.\ $\QCoh_\solid$-suave; see \cref{rslt:characterization-of-suave-maps} for a characterization of suaveness and see \cref{rslt:QCoh-solid-on-proper-and-smooth-maps} for the proof that all local complete intersections (and in particular all smooth maps) are cohomologically smooth. In the literature, \enquote{suave} often appears as \enquote{eventually coconnective} and \enquote{cohomologically smooth} often appears as \enquote{Gorenstein}.

\subsubsection{Ind-coherent sheaves}

With the above preparations at hand, we now discuss the construction of ind-coherent sheaves. Compared to other sources in the literature, our construction has the following characteristics:
\begin{enumerate}[(a)]
    \item It works over any regular noetherian base ring.
    \item It can be performed on the category of affine schemes first and then extended to stacks by descent. In particular, we make no use of any form of Nagata compactification.
\end{enumerate}
Part (b) might be surprising at first, as on affine schemes we cannot rely on compactifications to construct $!$-functors. Instead, we rely on the existing 6-functor formalism on $\QCoh_\solid$, as follows. For a ring $A$, we single out the full subcategory $\QCoh_{\solid,\pc}(A) \subseteq \QCoh_\solid(A)$ spanned by the pseudo-compact objects (i.e. those that admit a resolution by finite direct sums of compact generators). This subcategory is stable under $f^*$, $f_!$ and $\tensor$ on affine schemes and hence inherits a 3-functor formalism. The same is then true for $\PQCoh_{\solid,\pc}(A) := \Pro(\QCoh_{\solid,\pc}(A))$. But note that $\Coh(A) \subseteq \QCoh_{\solid,\pc}(A)$ and hence we get an embedding $\Pro(\Coh(A)) \injto \PQCoh_{\solid,\pc}(A)$. By formal nonsense this embedding admits a left adjoint
\begin{align*}
    \eta_A\colon \PQCoh_{\solid,\pc}(A) \to \Pro(\Coh(A)).
\end{align*}
Since $\eta_A$ is a localization, it is now a \emph{condition} to check that the 3-functor formalism on $\PQCoh_{\solid,\pc}(A)$ transfers to a 3-functor formalism on $\Pro(\Coh(A))$, which we readily do in \cref{rslt:3ff-for-PCoh-on-rings}. We finally pass to opposite categories and use the identification
\begin{align*}
    \Pro(\Coh(A))^\op = \Ind(\Coh(A)^\op) = \Ind(\Coh(A)) =: \ICoh(A),
\end{align*}
where the second identity is given by Grothendieck-Serre duality. Note that $\ICoh(A)$ is presentable, hence the 3-functor formalism on it is automatically a presentable 6-functor formalism. We then extend it to stacks via some by-now-standard descent procedures \cite{heyer-mann-6ff} in order to arrive at \cref{rslt:6ff-for-ICoh}. We state a simplified and less detailed version of this result here:

\begin{thm} \label{rslt:intro-ICoh-6ff}
Fix a classical regular noetherian ring $\Lambda$. There is a class $E$ of maps in $\Stk_\Lambda^\lafp$ and a presentable 6-functor formalism
\begin{align*}
    \ICoh\colon \Corr(\Stk_\Lambda^\lafp, E) \to \Cat
\end{align*}
with the following properties:
\begin{thmenum}
    \item If $X = \Spec A$ is affine then $\ICoh(X) = \Ind(\Coh(A))$.

    \item $E$ is stable under composition and pullback and if $f\colon Y \to X$ and $g\colon Z \to Y$ are maps in $\Stk^\lafp_\Lambda$ such that if $f, fg \in E$, then $g \in E$. Moreover, $E$ is local on the target and $!$-local (in particular, local for finite open covers) on source and target.

    \item For every $X \in \Stk_\Lambda^\lafp$ there is a natural functor
    \begin{align*}
        \gamma\colon \QCoh(X) \to \ICoh(X),
    \end{align*}
    compatible with tensor products and pullbacks. Moreover, if $f\colon Y \to X$ is a qcqs suave map representable in lafp schemes over $\Lambda$ then the natural map $\gamma f_* \isoto f_* \gamma$ is an isomorphism.
    
    \item For every afp qcqs scheme $X$ over $\Lambda$ there is a functor
    \begin{align*}
        \Psi\colon \ICoh(X) \to \QCoh(X)
    \end{align*}
    compatible with lower-! on $\ICoh$ and pushforward on $\QCoh$. If $f\colon Y \to X$ is a suave map of qcqs afp schemes over $\Lambda$ then $f_!\colon \ICoh(Y) \to \ICoh(X)$ admits a left adjoint $f^\natural$ and the natural map $f^* \Psi \isoto \Psi f^\natural$ is an isomorphism. Moreover, if $X$ is suave over $\Lambda$ then $\Psi$ admits a left adjoint
    \begin{align*}
        \Xi\colon \QCoh(X) \to \ICoh(X)
    \end{align*}
    and the natural map $\Xi f_* \isoto f_! \Xi$ is an isomorphism.

    \item $E$ contains every qcqs map $f\colon Y \to X$ that is representable in lafp schemes over $\Lambda$. Moreover, we have the following explicit descriptions:
    \begin{enumerate}[(a)]
        \item If $f$ is proper then $f^! \isom f^*$ on $\ICoh$.
        
        \item If $f$ is cohomologically smooth then $f_* = f_!(\gamma(\omega_f^\vee) \tensor -)$ on $\ICoh$, where $\omega_f$ denotes the dualizing sheaf of $f$. If $f$ is étale then $f_* = f_!$ on $\ICoh$.
    \end{enumerate}

    \item Let $E_0 \subseteq E$ be the class of qcqs maps that are representable in lafp schemes over $\Lambda$. Then there are sheafy 3-functor formalisms on $(\Stk_\Lambda^\lafp, E_0)$ together with morphisms of 3-functor formalisms
    \begin{align*}
        \ICoh \xfrom{\ \eta\ } \PQCoh_{\solid,\pc}^\op \injfrom \QCoh_{\solid,\pc}^\op \injto \QCoh_\solid^\op,
    \end{align*}
    which on affine schemes are defined as above.
\end{thmenum}
\end{thm}

Let us indicate the proof of this result. From what we have discussed above, we already get the 6-functor formalism $\ICoh$ satisfying (i) and (ii). Part (iii) reduces to the case of affine schemes $\Spec A$ by descent, where the functor $\gamma$ is given by the composition
\begin{align*}
    &\QCoh(A) = \Ind(\Perf(A)) \isoto \Ind(\Perf(A)^\op) \injto \Ind(\QCoh_{\solid,\pc}(A)^\op) \\ &\qquad= \PQCoh_{\solid,\pc}(A)^\op \xto{\eta} \ICoh(A).
\end{align*}
Here the second step uses naive duality $\intHom(-, A)$ and clearly all steps are symmetric monoidal. We emphasize that the construction of $\gamma$ (with all its functoriality) takes a considerable amount of effort in the existing literature. For part (iv) we again argue by descent to reduce the construction of $\Xi$ to the case of affine schemes, where it is simply the $\Ind$-extension of the inclusion $\Perf(A) \injto \Coh(A)$ (of course one has to work a bit in order to get the claimed functoriality). Part (vi) again reduces to affine schemes, where it is true by construction. Finally (v) formally reduces to (vi). Namely, if $f$ is cohomologically smooth then it is $\QCoh_\solid$-suave by definition, hence $\QCoh_\solid^\op$-prim. This property transfers along maps of 3-functor formalisms, so we deduce from (vi) that $f$ is $\ICoh$-prim, saying that $f_* = f_!(\delta_f^{\ICoh} \tensor -)$ on $\ICoh$, for a certain codualizing sheaf $\delta_f^{\ICoh}$. By comparing this construction with the above definition of $\gamma$, one arrives at $\delta_f^{\ICoh} = \gamma(\omega_f^\vee)$.

\begin{rmks}
\begin{rmksenum}
    \item Unlike other approaches to ind-coherent sheaves, our approach is not very explicit about the resulting functors in the 6-functor formalism. However, one can of course provide an explicit description in most cases of interest. For example, if $f\colon \Spec B \to \Spec A$ is a map of affine schemes, then the functor $f_!\colon \ICoh(B) \to \ICoh(A)$ is the $\Ind$-extension of the functor $f_*\colon \Coh(B) \to \QCoh(A) \to \ICoh(A)$, and similarly $f^*$ is the $\Ind$-extension of $f^!$ on $\Coh$ (cf.\ \cref{rslt:explicit-description-of-ICoh-functors-on-affines}). In practice, at least as far as this paper is concerned, the explicit description of these functors is not relevant -- it is enough to work with the above stated functoriality of $\gamma$, $\Xi$ and $\Psi$.

    \item Our notation for ind-coherent sheaves differs from the literature: We denote $\IndCoh$ by $\ICoh$, $f^{!,\ICoh}$ by $f^*$, $f_{*,\ICoh}$ by $f_!$ and $f^{*,\ICoh}$ by $f^\natural$. These conventions are chosen to represent the meaning of the functors in terms of the 6-functor formalism, but obfuscate how these functors are actually computed. We found neither notation completely satisfying.

    \item The morphism $\eta\colon \PQCoh_{\solid,\pc}^\op \to \ICoh$ explains why $\ICoh$ behaves somewhat strangely: It says that properties like smoothness and properness (which are reflected by $\QCoh_\solid$) translate to opposite notions in $\ICoh$.
\end{rmksenum}
\end{rmks}

We refer the reader to \cref{sec:alggeo.ICoh} for a detailed construction of $\ICoh$ and to \cref{sec:alggeo.properties-of-ICoh} for a collection of many more properties of $\ICoh$: We provide a large class of stacky maps along which $\ICoh$ admits $!$-functors (see \cref{rslt:QCA-maps-are-ICoh-fine}), we construct a t-structure on $\ICoh(X)$ for geometric stacks $X$ (see \cref{rslt:t-structure-on-ICoh}) and we show that $\ICoh(X) = \Ind(\Coh(X))$ on qcqs schemes and QCA stacks (see \cref{rslt:ICoh-equals-IndCoh-on-nice-stacks}).

\subsubsection{Excision}

One of the main tools developed in \cref{sec:alggeo} is an excision sequence for ind-coherent sheaves, which is crucially used in the proof of the spectral geometric lemma. Similar results are claimed in \cite{Halpern-Leistner} and \cite{gaitsgory-rozenblyum-vol2}, but the first reference contains an error and the second one does not state everything we need (and is hard to follow). To state our version of the result, we need the notion of the completion $X_{\hat Z}$ of a stack along a closed subset $Z \subseteq \abs X$: It is defined to be the substack whose points on a ring $A$ consist of those maps $\Spec A \to X$ that factor over $Z$ on underlying topological spaces (see \cref{def:completion-of-stack-at-closed-subset}). Moreover, for every closed immersion $i\colon Z \to X$ of stacks and integer $n \ge 1$, we provide a notion of \emph{$n$-th infinitesimal neighbourhood} $i_n\colon Z_n \to X$ in \cref{def:infinitesimal-neighbourhoods}; if $X = \Spec A$ and $Z = \Spec A/I$ then $Z_n = \Spec A/I^n$, for an appropriate definition of ideals in the derived context. With these definitions at hand, we prove the following:

\begin{thm}[\cref{rslt:ICoh-excision-all}] \label{rslt:intro-ICoh-excision}
Fix a classical regular noetherian ring $\Lambda$. Let $i\colon Z \to X$ be a closed immersion of lafp stacks over $\Lambda$ with completion $\hat i\colon X_{\hat Z} \injto X$ and open complement $j\colon U \injto X$. Let $M \in \ICoh(X)$.
\begin{thmenum}
    \item The map $\hat i$ is $\ICoh$-étale and the map $j$ is $\ICoh$-proper. Moreover, we have natural fiber sequences
    \begin{align*}
        \hat i_! \hat i^* M \to M \to j_* j^* M, \qquad j_* j^! M \to M \to \hat i_* \hat i^* M.
    \end{align*}

    \item For $n \ge 0$ let $i_n\colon Z_n \to X$ denote the $n$-th infinitesimal neighborhood of $i$. Then we have natural isomorphisms
    \begin{align*}
        \hat i_! \hat i^* M = \varinjlim_n i_{n!} i_n^* M, \qquad \hat i_* \hat i^* M = \varprojlim_n i_{n*} i_n^* M.
    \end{align*}

    \item Suppose that $L_{Z/X} \in \QCoh(Z)$ is perfect and denote $\Nm_{Z/X} := (L_{Z/X}[-1])^\vee$ the normal bundle of $i$. Then for all $n \ge 0$ we have
    \begin{align*}
        \cofib(i_{n!} i_n^* M \to i_{n+1,!} i_{n+1}^* M) = M \tensor i_! \gamma(\Sym^n_Z(\Nm_{Z/X})).
    \end{align*}
\end{thmenum}
\end{thm}

Note that in the last statement here, we do not put any strong conditions on $i$. For instance, $i$ can be any closed immersion such that the cotangent complexes $L_{X/\Lambda}$ and $L_{Z/\Lambda}$ are both perfect, which is the key case in our applications.

The proof of \cref{rslt:intro-ICoh-excision} relies heavily on a good understanding of the relation between completions and infinitesimal neighbourhoods, which we study in \cref{rslt:infinitesimal-neighborhoods}. Concretely, for every closed immersion $i\colon Z \to X$ we get a diagram
\begin{align*}
    Z = Z_1 \to Z_2 \to Z_3 \to \dots \to Z_n \to \dots \to X_{\hat Z} \injto X
\end{align*}
and in particular there is an induced map $\varinjlim_n Z_n \to X_{\hat Z}$. This map is not an isomorphism in general (only under strong assumptions on $i$, e.g.\ local complete intersections), but we prove that it is an isomorphism when evaluated on \emph{truncated} rings. Since $\ICoh$ only depends on truncated rings (see \cref{rslt:ICoh-only-depends-on-truncated-rings}) this is enough to deduce (i) and (ii) of \cref{rslt:intro-ICoh-excision}. Moreover, part (iii) follows from basic properties of infinitesimal neighborhoods, summarized in \cref{rslt:infinitesimal-neighborhoods}. Here $L_{Z/X}$ denotes the relative cotangent complex, i.e.\ a derived version of the sheaf $\Omega^1_{Z/X}$ of Kähler differentials (see \cref{sec:alggeo.cotangent}).

\begin{rmk}
There is a similar (yet less clean!) version of the excision sequence for $\QCoh$, see \cref{rslt:excision-for-QCoh}.
\end{rmk}

\subsubsection{Miscellany}

In the following we highlight some further results in \cref{sec:alggeo} that are used in this paper and may be of independent interest:
\begin{itemize}
    \item In \cref{sec:perf-and-coh} we observe that coherent sheaves admit a nice interpretation in terms of the 6-functor formalism $\QCoh_\solid$: They are exactly the suave sheaves over the base (see also \cite{Tang-MasterThesis}). This automatically proves many of the basic properties of coherent sheaves, e.g.\ stability under pushforward along proper maps. Similar results hold in other geometric settings (e.g. rigid or complex-analytic geometry) and give a new interpretation for coherent sheaves.

    \item In \cref{sec:adm-sheaves} we introduce and study the suave objects in the $\ICoh$-formalism, which we call the \emph{admissible sheaves}. After providing basic properties and examples, we prove in \cref{sec:spectral-temperization} the surprising result that the admissible dual of a coherent sheaf on a cohomologically smooth QCA stack $X$ lies in $\QCoh(X) \subseteq \ICoh(X)$.

    \item In \cref{sec:contraction} we study stacks $X$ with an action of $\mathbb A^1$ and prove a version of the contraction principle for quasicoherent sheaves. More concretely, if $X^0$ denotes the retract of the idempotent action of $0 \in \mathbb A^1$ on $X$ then we have a natural projection $\pi\colon X \to X^0$ with section $i\colon X^0 \to X$. The contraction principle seeks an understanding of the functor $\pi_*$ in terms of $i^*$, at least on $\mathbb G_m$-equivariant sheaves. In \cref{rslt:A1-retract-pushforward-in-terms-of-pullback} we show that for a perfect $\mathbb G_m$-equivariant sheaf $\mathcal F$ on $X$, the object $\pi_* \mathcal F$ can be built from $i^* \mathcal F$ with error terms of the form $(i^* i_* \calO_{X^0})^{\tensor n}$. This result is crucial to obtain the desired weight bounds in the spectral geometric lemma, and it also implies that $\pi_*$ sends perfect sheaves to coherent sheaves, at least on $\mathbb G_m$-graded pieces (see \cref{rslt:A1-retract-pushforward-preserves-perfect}).
\end{itemize}

\subsection{Applications} \label{sec:applications}

In this section we briefly sketch some applications of the categorical local Langlands conjecture. The details of these results will be given elsewhere.

\subsubsection{Computations of cohomology of local Shimura varieties}

Let $(G,\mu,b)$ be a local shtuka triple, where we assume for simplicity that $b$ is basic. Attached to this datum, we have a functor $R\Gamma_c(G,\mu,b)[-]$ from smooth $G_b(F)$-representations to smooth $G(F)$-representations. Classically, this functor sends $\rho$ to the $\rho$-isotypic part of the tower of local shtuka spaces attached to this triple. In the language of $\Bun_G$, however, this functor is incredibly simple: it sends $\rho$ to $i_{1}^\ast T_{V_{\mu}^{\vee}} i_{b!}\rho$. This reformulation already allowed Fargues-Scholze to prove strong finiteness properties of $R\Gamma_c(G,\mu,b)[-]$, which seem completely inaccessible from a classical point of view even for minuscule $\mu$. 

However, these functors are still extremely difficult to compute. When $\rho$ is an irreducible representation with supercuspidal $L$-parameter, Kottwitz gave a completely precise conjectural description of $R\Gamma_c(G,\mu,b)[\rho]$ in terms of the local Langlands correspondence, which is now a theorem in many cases. Outside the regime of supercuspidal parameters, there is not even a precise conjecture.

Nevertheless, there are a few scattered computations in the literature. In particular, let us take $G=\mathrm{PGL}_2$. Then highest weights $\mu=\mu_n$ are parametrized by nonnegative integers in the usual way, and $G_b$ is either $G$ or $PD^\times$ according to whether $n$ is even or odd; we write $b=1$ and $b_{1/2}$ for the two associated basic elements. Most classically, if $n=1$ then $(G,\mu_{1},b_{1/2})$ gives rise to the Lubin-Tate tower, and the cohomology $R\Gamma_c(G,\mu_1,b_{1/2})[\rho]$ is completely understood for irreducible $\rho$:
\begin{itemize}

\item  When $\dim \rho>1$, $\rho$ is the Jacquet-Langlands transfer of a supercuspidal $G(F)$-representation $\pi$, and 
\[R\Gamma_c(G,\mu_1,b_{1/2})[\rho] = \pi[0] \boxtimes \varphi_{\pi}\]
where $\varphi_{\pi}$ is the $L$-parameter of $\pi$. This is a famous theorem of Harris-Taylor.

\item When $\dim \rho=1$, $\rho=\mathbf{1}_{PD^\times}$ is (a twist of) the trivial representation, and ignoring the Weil action we have
\[R\Gamma_c(G,\mu_1,b_{1/2})[\mathbf{1}_{PD^\times}]=\mathrm{St}_G[0]\oplus \mathbf{1}_G[-1]. \]
This is an old calculation of Drinfeld, which follows immediately once one observes that the left-hand side identifies with the compactly supported \'etale cohomology of the Drinfeld space $\Omega^1 \subset \mathbf{P}^{1}_{\mathbf{C}_p}$.
\end{itemize}

Now suppose we replace the triple $(G,\mu_{1},b_{1/2})$ with the triple $(G,\mu_{2n+1},b_{1/2})$. For positive $n$, the resulting shtuka spaces are not rigid spaces, but it still makes perfectly good sense to try to compute $R\Gamma_c(G,\mu_{2n+1},b_{1/2})[\rho]$. For irreducible $\rho$ of dimension $>1$, this can be computed using the results in \cite{HKW, HanBasic}: in the notation of the example above, one gets
\[R\Gamma_c(G,\mu_{2n+1},b_{1/2})[\rho] = \pi[0] \boxtimes \mathrm{sym}^{2n+1}\varphi_{\pi}.\]
However, when $\rho=\mathbf{1}_{PD^\times}$ this cohomology seems impossible to compute directly on the automorphic side, but using the categorical conjecture, one can perform the following calculation.

\begin{prop}Forgetting the Weil group action, we have
\[R\Gamma_c(G,\mu_{2n+1},b_{1/2})[\mathbf{1}_{PD^\times}]=\mathrm{St}_{G}^{n+1}[0] \oplus \mathbf{1}_G[-1]\oplus \mathbf{1}_G[-3] \oplus \cdots \oplus \mathbf{1}_G[-2n-1] \]
for any $n\geq 0$.
\end{prop}

These sorts of calculations can be vastly generalized. We mention one striking result in this direction. Suppose $G=\mathrm{GL}_n$ and $\mu$ is any cocharacter, and $b$ is the (unique) basic element in $B(G,\mu)$. Then $R\Gamma_c(\GL_n,\mu,b)[\mathbf{1}]$ is the normalized compactly supported intersection cohomology of the basic Newton stratum $\mathrm{Gr}_{\GL_n, \leq \mu}^{\mathrm{bas}}$ in the Schubert cell $\mathrm{Gr}_{G, \leq \mu}$. When $\mu$ is minuscule, the Schubert cell is just a classical Grassmannian, and the basic locus is a standard example of a $p$-adic period domain, and in the most classical example where $\mu=(1,0,\dots,0)$ this is just the Drinfeld space inside $\mathbf{P}^{n-1}$. In that case the full cohomology was calculated by Schneider-Stuhler, but essentially no other examples were known beyond the Drinfeld case. Using the categorical conjecture along with further ideas, one can prove the following result.

\begin{thm}[Dai-H.] In the notation above, $R\Gamma_c(G,\mu,b)[\mathbf{1}]$ is a direct sum of shifted irreducible $\GL_n(F)$-representations, of total length $\dim V_{\mu}$, and it can be computed explicitly by an effective polynomial-time algorithm.
\end{thm}

This is a consequence of a close analysis of the categorical equivalence around the unramified semisimple parameter $\varphi$ sending Frobenius to $\mathrm{diag}(q^{(1-n)/2},q^{(1-n)/2-1},\dots,q^{(n-1)/2})$. Note that the stack of parameters is highly singular around the fiber over $\varphi$! A detailed explanation of these results will appear in \cite{DHstudy}.

\subsubsection{Hecke eigensheaves}

In his original conception of a geometrized form of the local Langlands correspondence, Fargues predicted the existence of canonical eigensheaves on $\Bun_G$ attached to $L$-parameters and enjoying many remarkable properties.

Using the categorical conjecture, it is very easy to build eigensheaves, since they are abundant on the spectral side. In fact, if $\phi$ is any $L$-parameter, corresponding to a map $i_{\phi}:\Spec \Qellbar \to \Par_G$, the sheaf $i_{\phi \ast} \Qellbar$ is formally seen to be an eigensheaf by an easy projection formula argument. Transporting this across the equivalence, we get a sheaf $\mathscr{F}_{\phi} = \mathbf{L}_{\psi}^{-1}(i_{\phi \ast} \Qellbar)$ on $\Bun_G$ which is formally a Hecke eigensheaf with eigenvalue $\phi$. We note that if $\phi$ lies in the smooth locus of the parameter stack, one can show that $i_{\phi \ast} \Qellbar$ is quasicoherent, so we are in the regime where $\mathbf{L}_{\psi}^{-1}= a_{\psi}$. For general $\phi$ this sheaf can be pathological, but under a mild genericity condition on $\phi$ we expect $\mathscr{F}_{\phi}$ to be a very beautiful object. The following example is representative.

\begin{conjecture}Let $G$ be any quasisplit group, and let $\phi_{\mathrm{St}}$ be the Steinberg parameter for $G$. Set $\mathscr{F}_{\mathrm{St}} = a_{\psi}(i_{\phi_{\mathrm{St}} \ast} \Qellbar)$. Then the sheaf $\mathscr{F}_{\mathrm{St}}$ enjoys the following properties.
\begin{propenum}
\item $\mathscr{F}_{\mathrm{St}}$ is a Hecke eigensheaf with eigenvalue $\phi_{\mathrm{St}}$.
\item $\mathscr{F}_{\mathrm{St}}$ is perverse and ULA, and its restriction to any connected component of $\Bun_G$ is indecomposable.
\item $\mathscr{F}_{\mathrm{St}}$ is Verdier self-dual.
\item For all $b$, we have
\begin{align*}
    i_{b}^{\ast \ren} \mathscr{F}_{\mathrm{St}} = & \mathrm{St}_{G_b} \otimes \delta_{b}^{-1/2} \\
    i_{b}^{! \ren} \mathscr{F}_{\mathrm{St}} = & \mathrm{St}_{G_b} \otimes \delta_{b}^{1/2}
\end{align*}
where $\mathrm{St}_{G_b}$ denotes the Steinberg representation of $G_b(F)$.
\end{propenum}
\end{conjecture}
Note that (i) is formal, but nothing else here is obvious, and it's not even clear that $\mathscr{F}_{\mathrm{St}}$ is nonzero. Using the methods of this paper together with some ideas from the work of Bertoloni Meli-Koshikawa, it is not hard to prove the following theorem.
\begin{thm}The Steinberg eigensheaf has all the desired properties for $\GL_n$.
\end{thm}

For a more detailed discussion of the expected properties of eigensheaves associated with nice $L$-parameters, see \cite[Section 3.1]{Beijing}. We note that Bertoloni Meli-Koshikawa have formulated a vast generalization of the conjectures outlined in \cite[Section 3.1]{Beijing}, which applies to $A$-parameters as well as $L$-parameters.

\subsubsection{Perverse t-exactness at generic parameters}

Like any good equivalence of categories, the categorical conjecture allows us to access deep properties on one side by moving to the other side. Here we mention just one example of this. It is well-known (see e.g. \cite[Section 1.2]{Beijing}) that $\D(\Bun_G)$ carries a natural perverse t-structure, and that perverse truncations preserve ULA sheaves. It is not typically true that Hecke operators are perverse t-exact, or that perverse truncation preserves finite sheaves. However, under a mild condition on the $L$-parameters, the categorical conjecture predicts both of these properties.

In the special case of $\GL_n$, we can be very concrete. A semisimple $L$-parameter $\varphi:W_F \to \GL_n(\Qellbar)$ is generic if $\Hom(\varphi,\varphi(1))=0$. This induces a direct factor decomposition \[\D(\Bun_{\GL_n})^{\mathrm{ULA}} = \D(\Bun_{\GL_n})^{\mathrm{ULA},\mathrm{gen}} \oplus \D(\Bun_{\GL_n})^{\mathrm{ULA},\mathrm{ngen}}\]
where the first resp. second summand is spanned by ULA sheaves $A$ whose $\varphi$-restricted summands $A_{\varphi}$ vanish for all non-generic resp. generic $\varphi$. This decomposition is stable under Hecke operators and perverse truncations, and induces an analogous direct factor decomposition on finite sheaves.

The following result will be proved in the forthcoming UCSD Ph.D. thesis of Nikolas Castro.

\begin{thm}The Hecke action on $\D(\Bun_{\GL_n})^{\mathrm{ULA},\mathrm{gen}}$ is perverse t-exact, and $\D(\Bun_{\GL_n})^{\mathrm{gen}}_{\fin}$ is preserved by perverse truncations.
\end{thm}

We refer to \cite[Sections 2.2-2.4]{Beijing} for a detailed overview of how the categorical conjecture can be used to prove results of this form.

\subsection{The anxiety of influence} \label{sec:otherpeople}
In this section we discuss the relationship between this paper and some related recent works in the literature.

First, let us mention the recent work of Zou \cite{Zou}, building on important ideas of Nguyen \cite{Nguyen}, which proves CLLC for $\GL_n$ after localizing onto a large explicit open substack of $\Par_{\GL_n}$. More precisely, writing $\Par_{\GL_n}^{\mathrm{gen}}$ for the open substack of \emph{generous} parameters, Zou proved that $a_{\psi}$ and $c_{\psi}$ induce mutually inverse equivalences of categories
\[\QCoh(\Par_{\GL_n}^{\mathrm{gen}}) \simeq \D(\Bun_{\GL_n})\otimes_{\QCoh(\Par_{\GL_n})} \QCoh(\Par_{\GL_n}^{\mathrm{gen}}).\]
Note that $\Par_{\GL_n}^{\mathrm{gen}}$ is smooth, so there is no difference between $\QCoh$ and $\IndCoh$ for this stack. Zou relies on ingenious tricks to deduce this inductively from the case of supercuspidal parameters, by an exhaustive analysis of the combinatorics of certain modifications of vector bundles. Moreover, his method also works integrally, and does not depend on any unverified compatibilities. On the other hand, conjectures of the first author \cite[\S 2.1]{Beijing} imply that categorical local Langlands should be much simpler over the generous locus, so it is perhaps not surprising that a direct attack is possible after localizing onto these parameters. In contrast, the essence of our work is to treat the whole conjecture in a uniform manner which accounts for the subtle phenomena associated with the singularities of $\Par_G$. It seems very unlikely to us that the methods of \cite{Zou} could be extended beyond the generous locus, or to groups other than $\GL_n$.

We also want to address two obvious questions.
\begin{itemize}
    \item What is the relationship between this work and the recent proof of the global unramified geometric Langlands conjecture \cite{GLC1,GLC2,GLC3,GLC4,GLC5}?

\item How is this work related to the work of Zhu \cite{ZhuTame} on categorical Langlands for the stack of $G$-isocrystals? 
\end{itemize}

We begin with the first question. Here the reader is perhaps still wondering:

\emph{Did we just copy the arguments of \cite{GLC1,GLC2,GLC3,GLC4,GLC5} but cross out ``projective curve'' and write ``Fargues-Fontaine curve'' everywhere?}

The answer to this is \textbf{definitely not}. Of course, in its broadest strokes, some aspects of the picture in \cite{GLC1,GLC2,GLC3,GLC4,GLC5} are parallel to the picture in the Fargues-Scholze setting: after all, we \emph{are} trying to do geometric Langlands over the Fargues-Fontaine curve. In particular, the spectral action, the Whittaker sheaf, the functor of enhanced Whittaker coefficient, and the Eisenstein series and constant term functors play broadly similar roles for us as in the proof of the GLC. Moreover, we construct the Langlands functor by suitably lifting the functor of enhanced Whittaker coefficient along the surjection $\Psi:\IndCoh(\Par_G)\twoheadrightarrow \QCoh(\Par_G)$, just as in \cite{GLC1}. However, as soon as one considers the details, there is a significant divergence of arguments at nearly every step.\footnote{Some portions of the following overlap with our discussion so far, but the point here is to contrast some key steps in our arguments against comparable steps in the proof of GLC.}

\begin{itemize}

 \item Many of the crucial intermediate ingredients in the proof of GLC, especially techniques related to opers and ``localization at critical level'', seem to have no analogue whatsoever over the Fargues-Fontaine curve.\footnote{It is possible that this is just a failure of imagination on our part, or at least that such ingredients do adapt to the setting of ``Fargues-Scholze with $p$-adic coefficients''.}

    \item In \cite{GLC1,GLC2,GLC3,GLC4,GLC5}, the compatibility of Eisenstein series with enhanced Whittaker coefficient is just one of many, many intermediate results which assemble into the complicated architecture of the proof. By contrast, one of the critical realizations for us was that for groups with a sufficiently well-understood local Langlands correspondence, this compatibility \textbf{alone} is the essential key to categorical local Langlands.

    \item We construct the Langlands functor by proving that $c_{\psi}$ sends compact automorphic sheaves towards coherent sheaves and ind-completing. In the proof of GLC, the fact that enhanced Whittaker sends compact automorphic sheaves to coherent sheaves is only deduced at the very end, as a \emph{consequence} of GLC. To construct the Langlands functor, \cite{GLC1} instead proves a t-structure bound on $c_{\psi}$ restricted to compact automorphic sheaves, uniformly for all groups, by $D$-module methods building on \cite{FR} which are totally unavailable to us. However, this means that during the intermediate steps of the arguments in the proof of the geometric Langlands conjecture, the Langlands functor is colimit-preserving but \emph{not} known to preserve compact objects. In our setting, preservation of compacts is baked into the construction.

    \item We work with the Langlands functor and its \emph{right} adjoint, which is colimit-preserving thanks to the previous step, and which we are able to show also preserves compact objects. By contrast, in the proof of GLC the Langlands functor isn't known to preserve compact objects until the end of the proof, so its right adjoint is a priori \emph{not colimit-preserving} which renders it of little use. Instead, \cite{GLC3} goes to great effort to construct a \emph{left} adjoint to the Langlands functor. This construction relies on the compatibility of the Langlands functor with constant terms, which is one of the hardest steps in the whole argument and relies critically on Kac-Moody localization and other techniques developed in \cite{GLC2} which are highly specific to the setting of $D$-modules. We emphasize that we do \emph{not} expect constant term compatibility to be available in the Fargues-Scholze setting in the near future!

    \item In both the proof of GLC and in the present paper, a key step is to prove that the enhanced Whittaker coefficient functor evaluated on the Whittaker sheaf yields a line bundle. In \cite{GLC1,GLC2,GLC3,GLC4,GLC5} this is proved by a very indirect argument, which ultimately boils down to a trick \cite{GLC5} relying on the fact that $\mathrm{LS}_{\hat{G}}$ has no nonconstant global functions. In our setting, the Whittaker sheaf is concentrated at $b=1$, so we need to compute $c_{\psi}i_{1!}$ on the Whittaker representation. Our strategy (roughly) is to identify \emph{the whole functor} $c_{\psi}i_{1!}$ with an explicit functor $\mathscr{S}_{G}: \D(G(F),\Qellbar)\to \QCoh(\Par_G)$ which provably has the correct value on the Whittaker representation. For $\GL_n$, the desired explicit functor was already constructed in \cite{BCHN}, and a critical realization for us was that compatibility of Eisenstein series with enhanced Whittaker coefficient already gives enough leverage to identify these functors for $\GL_n$.

    \item We deduce full faithfulness of the right adjoint from the previous step by piece of homological subterfuge, which relies on a coherent analogue (Theorem \ref{rslt:spectral-temperization}) of Beraldo's theorem that Verdier duality on $\Bun_G$ carries compact $D$-modules into \emph{tempered} $D$-modules. As a consequence of this subterfuge, we get an \emph{explicit} formula for the right adjoint on many coherent sheaves in terms of the spectral action and some explicit dualities. These arguments seem to have no analogue in the proof of GLC.

     \item The conservativity of enhanced Whittaker coefficient on compact $D$-modules was already proved in \cite{FR}, so in some sense the bulk of the effort in the proof of GLC goes into constructing the left adjoint and proving its full faithfulness. In our setting, the arguments of \cite{FR} are totally unavailable, and conservativity of the enhanced Whittaker coefficient on compact automorphic sheaves is one of the hardest steps in this paper. In fact, we do not prove it directly or for general groups, but only for well-understood groups by a subtle induction which assumes the truth of CLLC for all proper Levis and relies on a difficult vanishing result proved by Bertoloni Meli-Koshikawa.
\end{itemize}

We also want to comment on the relationship with the work of Zhu \cite{ZhuTame}. The real answer here is that \textbf{there is no direct relationship}. Recall that Zhu constructs an equivalence of categories $\mathbf{L}_{\psi}^{\mathrm{Zhu}}$ between ind-coherent sheaves on $\Par_{G}^{\mathrm{tame}}$ and a suitable category of depth-zero sheaves on the stack $\mathrm{Isoc}_G$ of $G$-isocrystals, under mild conditions on $G$. In extreme caricature, his strategy is to produce this equivalence directly in one stroke, by taking the Frobenius-twisted categorical trace of a suitable form of Bezrukavnikov's equivalence. Architecturally, this is extremely \emph{dissimilar} from the structure of the present paper (and the structure of the proof of GLC): there are no Hecke operators, no spectral action, no Langlands functor, and no Eisenstein or constant term functors. On the other hand, Zhu's approach has the evident advantages that it yields an equivalence of categories uniformly in the group and without any prior knowledge of local Langlands, and that it also works with integral or mod-$\ell$ coefficients once the appropriate forms of Bezrukavnikov's equivalence are available.

Of course, the work \cite{ZhuTame} should be compatible with the results in the present paper. In particular, ongoing work of Gleason-Hamann-Ivanov-Lourenco-Zou \cite{GHILZ} is expected to yield a canonical equivalence of categories $\pitch : D(\mathrm{Isoc}_G)\overset{\sim}{\to}\D(\Bun_G)$, and at least for well-understood groups one can ask whether the resulting diagram \[
\xymatrix{D(\mathrm{Bun}_{G})\ar[r]^{\mathbf{L}_{\psi}} & \mathrm{IndCoh}(\mathrm{Par}_{G})\\
D(\mathrm{Isoc}_{G})^{\mathrm{depth.0}}\ar[u]_{\pitch}\ar[r]_{\sim}^{\mathbf{L}_{\psi}^{\mathrm{Zhu}}} & \mathrm{IndCoh}(\mathrm{Par}_{G}^{\mathrm{tame}})\ar[u]
}
\]
commutes. This is presumably a very difficult problem.\footnote{The discussion in \cite[Section 2.7.11]{GaitsgoryICM} is a significant oversimplification of the situation.}

\subsection{The future}\label{sec:future}

In this section, we discuss some natural problems which are not addressed in this paper, and make some comments on how things might develop from here.

As noted several times, we do not prove Conjecture \ref{conj:cpsiEiscompatible} in this paper. One of our short-term goals is to prove this conjecture uniformly in all groups, in future joint work with Linus Hamann. In classical geometric Langlands, the analogue of this conjecture is a well-known theorem, with proofs discussed to varying degrees of detail in \cite{GLC3}, \cite{GaitsgoryRaskintrick}, and \cite{FHparabolic}. We expect that to first approximation, the known proofs in that setting should adapt without significant difficulty to the Fargues-Fontaine curve. In fact, Hamann-Imai have already announced a proof of this conjecture in the special case where $G=\GL_2$.

Next, let us comment on our choices regarding the base field and the coefficients. As noted above, we work over a mixed characteristic local field $F$ throughout. Most of our methods should adapt to the equal characteristic setting, but the literature on classical local Langlands is extremely thin in equal characteristic. Moreover, even the most basic finiteness properties of the Whittaker representation seem to be open questions here. It is plausible one could adapt our arguments at least to the case of $\GL_n$ over an equal characteristic local field, but other classical groups seem far out of reach.

Additionally, we work throughout with rational coefficients. While it is conceivable that some steps in our program adapt to integral and modular coefficients, these adaptations seem to require new ideas. Moreover, we are not sanguine about the prospects of proving Conjecture \ref{conj:cpsiEiscompatible} with integral coefficients in the near future. As such, we did not try too hard to see how far one can get with our methods integrally, and it would be interesting to understand this.

Finally, we do not have any idea how to prove categorical local Langlands without significant input from classical local Langlands. We very much hope that in the future, there will be a second-generation proof strategy which avoids such inputs. It's very conceivable that a future melding of our approach with Zhu's methods via the equivalence established in \cite{GHILZ} will lead to an improved strategy along these lines.

\subsection{Notation and conventions}

This paper lies at the interplay of several modern fields of mathematics and we borrow notation from all of them, as well as introduce our own. Especially in the setting of ind-coherent sheaves our notation differs a lot from the existing literature. In the following we provide a guideline for all the conventions used in this paper.

\subsubsection{Geomeric local Langlands}

We mostly stick with the notation in \cite{FS}, but denote our base field by $F$ instead of $E$ and usually assume it to live in mixed characteristic, i.e. to be a finite extension of $\Q_p$. We denote by $W_F$ the Weil group of $F$, with inertia subgroup $I_F$. In most parts of this paper we fix a connected reductive group $G$ over $F$ and denote by $\Bun_G$ the v-stack of $G$-bundles on the Fargues--Fontaine curve. On the spectral side for any $\Z_\ell$-algebra $\Lambda$ we denote by $\Par_{G,\Lambda}$ the stack of $L$-parameters for $G$ relative to $\Lambda$, which is denoted $Z^1(W_F, \hat G)_\Lambda / \hat G$ in \cite{FS}. We denote by $\D(\Bun_G, \Lambda) := \D_\lis(\Bun_G, \Lambda)$ the category of lisse solid sheaves of $\Lambda$-modules on $\Bun_G$. In the case $\Lambda = \Qellbar$ we omit $\Lambda$ from the above notation and simply write $\Par_G$ and $\D(\Bun_G)$. For a locally profinite group $H$ we denote by $\D(H, \Lambda)$ the (derived) category of smooth $H$-representations on $\Lambda$-modules; in some (few) parts of the paper this is also denoted by $\Rep_\Lambda^\sm(H)$. We denote by $\Div^1$ the stack of Cartier divisors of degree $1$ on the Fargues--Fontaine curve and for every finite set $I$ we denote $\Div^I := (\Div^1)^I$.

\subsubsection{Higher categories}

This paper is deeply rooted in the language of higher categories and we do not explicitly mention this fact in our definitions and results. In particular, we refer to $(\infty,1)$-categories simply as categories. We denote by $\Ani$ the category of anima, which is sometimes also called the category of spaces or $\infty$-groupoids. We denote by $\Cat$ the ($\infty$-)category of ($\infty$-)categories and by $\CMon = \CAlg(\Cat)$ the category of symmetric monoidal categories. We refer the reader to \cite[\S A, \S B]{heyer-mann-6ff} for a quick overview of a large part of the definitions and notation used in higher categories and higher algebra.

\subsubsection{Algebraic geometry}

Our spectral results make use of algebraic geometry and all our constructions are implicitly derived, e.g.\ by a ring we mean an animated ring, by a module we mean a complex of modules and by $\QCoh(X)$ we denote the derived ($\infty$-)category of quasicoherent sheaves on a stack $X$. We call an object \emph{static} if we want to indicate that it is of non-derived nature (e.g.\ a static ring is a classical commutative ring). We refer the reader to \cref{sec:introindcoherent} for a quick overview of most of the relevant definitions and notation used in this paper. Many of these notions have appeared in the literature before and often with different names, because our notation is inspired by the theory of 6-functor formalisms. We compile a list of \emph{translations} below.
\begin{itemize}
    \item What we call \emph{local complete intersection} maps are often called \enquote{quasi-smooth} maps in the literature.
    
    \item What we call \emph{suave} maps are known as \enquote{eventually coconnective} or \enquote{finite Tor amplitude} in the literature (cf.\ \cref{rslt:characterization-of-suave-maps}).

    \item What we call \emph{cohomologically smooth} maps are known as \enquote{Gorenstein} maps in the literature.

    \item For a geometric stack $X$ we denote by $\Coh(X) \subseteq \QCoh(X)$ the category of coherent sheaves on $X$, i.e. those that smooth locally can be written as a bounded complex with finitely generated cohomologies; in particular we do not enforce any boundedness on the support. In the context of the $L$-parameter stack, we denote by $\Coh^\qc(\Par_G) \subseteq \Coh(\Par_G)$ the full subcategory of those sheaves that have quasicompact support. We use the similar notation $\Perf^\qc(\Par_G) \subseteq \Perf(\Par_G)$ for perfect sheaves.

    \item In the context of ind-coherent sheaves, what we denote as $\ICoh(X)$, $f^*$, $f_!$, $f^\natural$, $\tensor$ and $\gamma$ are usually denoted by $\IndCoh$, $f^{!,\IndCoh}$, $f_{*,\IndCoh}$, $f^{*,\IndCoh}$, $\tensor^!$ and $\Upsilon$, respectively. Our use of the symbols $\Xi$ and $\Psi$ agrees with their common usage in the previous literature.

    \item In the context of a symmetric monoidal category $\cat C$ we often denote the tensor unit by $1$. For example, $1 \in \ICoh(X)$ denotes the tensor unit of $\ICoh(X)$, which is often denoted by $\omega_X$ in the literature (because in many geometric situations it corresponds to the dualizing sheaf on $X$ via $\Xi$).
\end{itemize}

\subsection{Acknowledgments}

First and foremost, we extend our deep and sincere gratitude to Alexander Bertoloni Meli and Teruhisa Koshikawa. Our perspective on categorical local Langlands was shaped by innumerable conversations with them over the last three years. Their calculations with coherent sheaves on the parameter stack for $\mathrm{PGL}_2$ gave the perfect laboratory for us to experiment in, and these experiments led directly to some of the key ideas in this paper, including the crucial realizations that something like Theorem \ref{rslt:spectral-temperization} might be true, and (one day later) that the functor $t_{\psi}$ could be a reasonable gadget to consider. Moreover, at a late stage in this project, their ongoing foundational work on the theory of cuspidal coherent sheaves entered into our arguments in an indispensable way (see Theorem \ref{thm:BMKsemiorthogonalKey} and the proof of Theorem \ref{thm:DreamInductionOnLevisTheorem}).

Special thanks go to Eugen Hellmann, who asked one of us (LM) in early 2023 what Verdier duality should correspond to on the spectral side. The resulting theory of admissible duality for ind-coherent sheaves (conceived in July 2023 on a high-speed train in China) was a crucial preliminary step towards this project, and plays a key role in our arguments. Eugen's paper \cite{Hel} was also deeply influential on our thinking.

Next, we want to thank Linus Hamann for a great many conversations about categorical local Langlands and Eisenstein series in particular, and for embarking with us on the serious task of proving the key Conjecture \ref{conj:cpsiEiscompatible}.

We thank Dennis Gaitsgory and Sam Raskin for many patient and helpful conversations about the geometric Langlands program. One of us (DH) began thinking seriously about these topics after watching an inspiring Zoom lecture by Sam in Spring 2022 \cite{RaskinTalk}.

We thank Wee Teck Gan, Zachary Gardner, Toby Gee, Ian Gleason, Linus Hamann, Jeroen Hekking, Eugen Hellmann, David Helm, Claudius Heyer, Pol van Hoften, Christian Johansson, Johan de Jong, Maxime Ramzi, Peter Scholze, Maarten Solleveld, Yifei Zhao, and Xinwen Zhu for some very helpful conversations.

We presented evolving versions of this project to audiences in M\"unster, Bonn, Oberwolfach, Hanoi, Chicago, Singapore, Bonn (again), Sydney, and Marseille, and we want to thank our colleagues at these events for their valuable feedback and questions, which in some cases helped us sharpen our results significantly.

Finally, particular thanks go to Joseph Bernstein, Vladimir Drinfeld, Wee Teck Gan, Michael Harris, Tasho Kaletha, Wieslawa Niziol, and Peter Scholze for their interest in this project and for offering invaluable encouragement during the writing process.

\section{Spectral preparations} \label{sec:alggeo}

In this section we survey derived algebraic geometry and develop the necessary tools to prove our results about ind-coherent sheaves on the $L$-parameter stack in \cref{sec:spectral}. This section is not meant as a comprehensive introduction to derived algebraic geometry, even though it does build a substantial amount of the theory from the ground up. Everything we prove and define in this subsection is in the pursuit of the main tools that we need in the next section: A good understanding of ind-coherent sheaves including the six operations on them (see \cref{rslt:6ff-for-ICoh}), a large class of stacks that have nice geometric properties (see \cref{sec:alggeo.QCA-stacks}), and an excision result for ind-coherent sheaves on these stacks (see \cref{rslt:ICoh-excision-all}). As auxiliary tools we develop, among others: A good definition of coherent sheaves on stacks including the expected stability properties, a definition of ideals in derived algebraic geometry, a good notion of reduced stacks, and basic properties of the cotangent complex and lci maps. Moreover, we use a novel approach to ind-coherent sheaves via the existing formalism of solid quasi-coherent sheaves. This approach has the advantage of giving a very quick definition of the six functors on ind-coherent sheaves (including all higher coherences) and also provides the expected properties, like the computation of dualizing complexes.

\subsection{Rings and modules} \label{sec:alggeo.rings-and-modules}

In the following we recall some basic definitions regarding rings and modules, including perfect, pseudo-coherent and flat modules, as well as faithfully flat descent. We refer the reader to \cite[Part~I]{RC.DeRhamStacksLecture} for a more detailed introduction of the $\infty$-categorical structures introduced below.

\begin{defn}
\begin{defenum}
    \item By a \emph{ring} we mean an animated ring, i.e. roughly a \enquote{simplicial commutative ring up to homotopy} (see \cite[\S25.1]{lurie-SAG}). Classical rings will be called \emph{static} rings. We denote by $\Ring$ the category of rings and by $\Alg_A := \Ring_{A/}$ the category of $A$-algebras, for a fixed ring $A$. For every ring $A$ we get an underlying static ring $\pi_0 A$ together with a sequence of static $\pi_0 A$-modules $\pi_i A$ for $i \ge 0$. For an integer $n \ge 0$, we say that $A$ is \emph{$n$-truncated} if $\pi_i A = 0$ for $i > n$.

    \item Given a ring $A$, we denote the category of $A$-modules by $\D(A)$ (denoted $\Mod_A$ in \cite[\S25.2]{lurie-SAG}). In the case that $A$ is static, this is just the usual derived category of $A$-modules. In general, $\D(A)$ is a stable (i.e. \enquote{triangulated}) presentable category with $t$-structure. We denote the homology and cohomology functors for this $t$-structure by $H_n$ and $H^n$ respectively, and call an $A$-module \emph{static} if it is concentrated in degree $0$. We denote by $\D^{\le n}(A), \D^{\ge n}(A) \subseteq \D(A)$ the full subcategories of bounded $A$-modules, where the bound is imposed in cohomological notation. We similarly denote $\tau^{\le n}$ and $\tau^{\ge n}$ the truncations functors.
    
    The category $\D(A)$ comes equipped with a tensor product $\tensor_A$. Note that for every $A$-algebra $B$, we can regard $B$ as an object of $\D(A)$ and then $\pi_i(B) = H_i(B)$.
\end{defenum}
\end{defn}

\begin{rmk}
The category $\Ring$ is the category of ring objects in the \enquote{derived algebraic context} $\D(\Z)$. We discuss this more general construction in \cref{sec:filtered-rings}, as we need it in order to define filtered and graded rings below.
\end{rmk}

The reader who is unfamiliar with animated rings may replace them with ordinary rings for the most part of this paper, especially in the later sections.

For a ring $A$ the category $\Alg_A$ of $A$-algebras is presentable (see \cite[Remark~25.1.1.2]{lurie-SAG}) and in particular has all small colimits. Filtered colimits (and more generally sifted colimits, i.e. filtered colimits plus \enquote{quotients}) are computed on the underlying $A$-modules, while the coproduct of two $A$-algebras $B$ and $C$ is given by $B \tensor_A C$. The initial object of $\Alg_A$ is $A$ and the initial object of $\Ring$ is $\Z$. Moreover, for every ring $A$ there is a polynomial ring $A[x]$ with the property
\begin{align*}
    \Hom_{\Alg_A}(A[x], B) = \Hom_{\D(A)}(A, B)
\end{align*}
for any $A$-algebra $B$. Restricting to $\pi_0$, we see in particular that an $A$-algebra map $A[x] \to B$ is the same as an element in $\pi_0(B)$; this fact will be used repeatedly throughout this section. Note furthermore that $A[x] = \Z[x] \tensor_\Z A$, where $\Z[x]$ is the usual polynomial ring.

We now introduce some finiteness conditions on modules and on maps of rings, generalizing the notions of \enquote{finitely presented module} and \enquote{finitely presented map}:

\begin{defn}
Let $A$ be a ring and $M \in \D(A)$ an $A$-module.
\begin{defenum}
    \item We say that $M$ is \emph{perfect} if it can be obtained via finite limits, colimits, and retracts from copies of $A$.
    
    \item \label{def:pseudocoherent-modules} We say that $M$ is \emph{pseudo-coherent} if for every integer $n$, $\tau^{\ge n} M$ is a compact object in $\D^{\ge n}(A)$. Equivalently $\Hom(M, -)$ commutes with colimits of uniformly left-bounded filtered diagrams in $\D(A)$.
\end{defenum}
We denote by
\begin{align*}
    \Perf(A) \subseteq \D_\pc(A) \subseteq \D(A)
\end{align*}
the full subcategories of perfect and pseudo-coherent $A$-modules, respectively.
\end{defn}

In the case that $A$ is a static ring, the notion of perfect and pseudo-coherent modules is classical, see e.g. \cite[Sections~0656,~064N]{stacks-project}. Recall that an $A$-module is perfect iff it is compact iff it is dualizable (cf. \cite[Example~B.1.17]{heyer-mann-6ff}). For pseudo-coherent modules we have the following characterizations:

\begin{lem} \label{rslt:criteria-for-pseudocoherent-modules}
Let $A$ be a ring and $M \in \D(A)$ an $A$-module. Then the following are equivalent:
\begin{lemenum}
    \item $M$ is pseudo-coherent.
    
    \item \label{rslt:approximation-of-pseudocoh-by-perfect} For every integer $n$ there is a perfect $A$-module $P$ and a map $P \to M$ such that $\fib(P \to M)$ sits in cohomological degrees $\le n$.
    
    \item \label{rslt:pseudo-coherent-implies-resolution-by-finite-free} There are an integer $m \in \Z$ and a simplicial diagram $(P_n)_{n\in\Delta}$ in $\D(A)$ such that all $P_n$ are finite free, the diagram is uniformly right-bounded, and $M[m] = \varinjlim_n P_n$.
\end{lemenum}
\end{lem}
\begin{proof}
The implication (i) $\implies$ (iii) is shown in \cite[Lemma~C.6.6.3]{lurie-SAG} (alternatively, see \cite[Proposition~7.2.4.11]{lurie-higher-algebra}). The other implications are easy: We first show that (iii) implies (ii), so assume $M = \varinjlim_{n\in\Delta} P_n$ as in (iii). By assumption all $P_n$ are right-bounded by a uniform bound and we can assume that bound to be $0$, i.e. $P_n \in \D^{\le0}(A)$ for all $n$. Now fix an integer $m \ge 0$ and let $M_m := \varinjlim_{n \in \Delta_{\le m}} P_n$ be the $m$-truncated colimit of $P_\bullet$. Then $M_m$ is perfect (being a finite colimit of perfect modules) and there is a natural map $M_m \to M$, whose fiber sits in cohomological degrees $\le -m$ by \cite[Proposition~1.2.4.5(4)]{lurie-higher-algebra}. This implies (ii). Finally, we show that (ii) implies (i). Fix some integer $n$ and let $P \to M$ be as in (ii). Then $\tau^{\ge n} P$ is compact in $\D^{\ge n}(A)$ and $\tau^{\ge n} P = \tau^{\ge n} M$. This implies (i).
\end{proof}

One should see condition (iii) in \cref{rslt:criteria-for-pseudocoherent-modules} as a generalization of the statement that a pseudo-coherent module admits a resolution by finite free modules to the case that $A$ is an animated ring. We also record the following stability properties of perfect and pseudo-coherent modules. Over classical rings, similar results are e.g. found in \cite[Sections~064N,~0656]{stacks-project}.

\begin{lem} \label{rslt:stability-of-perfect-and-pseudo-coherent}
Let $A$ be a ring.
\begin{lemenum}
    \item \label{rslt:perfect-and-pseudo-coherent-stable-under-tensor-colim-base-change} The subcategories of perfect and of pseudo-coherent $A$-modules are stable under tensor product, finite limits and colimits, and retracts. For every ring map $A \to B$, the base-change $- \tensor_A B\colon \D(A) \to \D(B)$ preserves perfect and pseudo-coherent modules. Additionally, pseudo-coherent $A$-modules are stable under uniformly right-bounded geometric realizations (i.e. colimits over $\Delta$).

    \item \label{rslt:stability-of-perfect-and-pseudo-coherent-under-forget} Let $A \to B$ be a map of rings and assume that $B$ is pseudo-coherent (resp. perfect) as an $A$-module. Let $N \in \D(B)$ be a $B$-module. If $N$ is pseudo-coherent (resp. perfect) as a $B$-module, then it is so as an $A$-module. The converse is true for pseudo-coherence.

    \item \label{rslt:pseudo-coherent-plus-fin-tor-dim-implies-perfect} An $A$-module $M$ is perfect if and only if it is pseudo-coherent and has finite Tor dimension (i.e. $- \tensor_A M$ sends $\D^{\ge0}(A)$ to $\D^{\ge n}(A)$ for some integer $n$).
\end{lemenum}
\end{lem}
\begin{proof}
The first two claims of (i) are clear. For the last part of (i), if $P = \varinjlim_{n\in\Delta} P_n$ with all $P_n$ pseudocoherent, then for any left-bounded $A$-module $M$, $\Hom(P, M)$ depends only on $P^{(m)} := \varinjlim_{n\in\Delta_{\le m}} P_n$ for some $m \ge 0$ depending only on the boundedness of $M$ (by \cite[Proposition~1.2.4.5(4)]{lurie-higher-algebra}). But $P^{(m)}$ is a finite colimit and hence pseudo-coherent.

The first part of (ii) follows easily from (i) by writing $N$ as a geometric realization (resp. finite colimit and retract) of copies of $B$. To prove the second part of (ii), we use that $N = \varinjlim_{n\in\Delta} N \tensor_A B^{\tensor_A n+1}$ (the right-hand side computes the relative tensor product $N \tensor_B B$). Now if $N$ is pseudo-coherent as an $A$-module then $N \tensor_A B^{\tensor_A n+1}$ is pseudo-coherent as a $B$-module for all $n$ and hence $N$ is pseudo-coherent as a $B$-module by (i).

Part (iii) is \cite[Proposition~7.2.4.23(4)]{lurie-higher-algebra}.
\end{proof}

Having discussed the finiteness properties of modules, let us now introduce a related finiteness property for maps of rings:

\begin{defn}
Let $f\colon A \to B$ be a map of rings.
\begin{defenum}    
    \item $f$ is \emph{almost finitely presented} (\emph{afp}) if there is a factorization $A \to A[x_1, \dots, x_n] \to B$ such that $B$ is pseudo-coherent as an $A[x_1, \dots, x_n]$-module. We denote by
    \begin{align*}
        \Alg_A^\afp \subseteq \Alg_A
    \end{align*}
    the full subcategory of almost finitely presented $A$-algebras.

    \item $f$ is \emph{finite} if the induced map $\pi_0 A \to \pi_0 B$ is finite, i.e. $\pi_0 B$ is a finitely generated $\pi_0 A$-module.
\end{defenum}
\end{defn}

\begin{lem}
\begin{lemenum}
    \item Almost finitely presented maps are stable under base-change and composition.

    \item \label{rslt:afp-same-as-truncated-compact} A ring map $A \to B$ is almost finitely presented if and only if for all $n \ge 0$ and all filtered systems $(C_i)_i$ of $n$-truncated $A$-algebras we have
    \begin{align*}
        \Hom_{\Alg_A}(B, \varinjlim_i C_i) = \varinjlim_i \Hom_{\Alg_A}(B, C_i).
    \end{align*}
    
    \item \label{rslt:afp-maps-are-cancellative} Given ring maps $f\colon A \to B$ and $g\colon B \to C$ with $f$ and $fg$ almost finitely presented, then $g$ is almost finitely presented.
    
    \item \label{rslt:finite-afp-map-equiv-pseudo-coherent} Let $f\colon A \to B$ be a ring map such that $\pi_0(A) \surjto \pi_0(B)$ is surjective. Then $f$ is almost finitely presented if and only if $B$ is pseudo-coherent as an $A$-module.
\end{lemenum}
\end{lem}
\begin{proof}
Part (i) follows easily from \cref{rslt:stability-of-perfect-and-pseudo-coherent}. We now prove (iv), so let $f\colon A \to B$ be given as in the claim. The \enquote{if} direction of the claim is clear from the definition, so we only need to prove the \enquote{only if} direction. Thus, assume that $f$ is almost finitely presented, so there is a map $g\colon A[x_1, \dots, x_n] \to B$ such that $B$ is pseudo-coherent as an $A[x_1, \dots, x_n]$-module. The map $g$ is equivalently given by elements $b_1, \dots, b_n \in \pi_0(B)$, and by assumption on $f$ we can choose preimages $a_1, \dots, a_n \in \pi_0(A)$. Then $g$ factors as $A[x_1, \dots, x_n] \to A \to B$. By \cref{rslt:stability-of-perfect-and-pseudo-coherent-under-forget} we are thus reduced to showing that $A$ is pseudo-coherent as an $A[x_1, \dots, x_n]$-module. By induction on $n$ we can reduce to the case $n = 1$, i.e. we want to see that for any element $a \in \pi_0(A)$ with induced map $A[x] \to A$, $A$ is pseudo-coherent as an $A[x]$-module. But $A = \cofib(A[x] \xto{x-a} A[x])$, as can be shown by base-change from the universal case $\Z[y] = \cofib(\Z[y,x] \xto{x-y} \Z[y,x])$. Thus $A$ is even perfect as an $A[x]$-module . This finishes the proof of (iv).

We next prove (iii), which by (i) and the argument in \cite[Lemma~2.1.5]{heyer-mann-6ff} reduces to showing that if $f\colon A \to B$ is an almost finitely presented map, then the same is true for the \enquote{diagonal} $\Delta_f\colon B \tensor_A B \to B$ given by multiplication. Using (i) one checks easily that this claim is stable under compositions in $f$, so we are reduced to the cases that $B = A[x]$ and that $B$ is pseudo-coherent as an $A$-module. In the second case we use the factorization $B \to B \tensor_A B \to B$ and \cref{rslt:stability-of-perfect-and-pseudo-coherent-under-forget} to deduce that $B$ is pseudo-coherent over $B \tensor_A B$, as desired. It remains to handle the case $B = A[x]$, in which case the diagonal is the map $A[x_1, x_2] \to A[x]$ sending both $x_1$ and $x_2$ to $x$. But $A[x]$ is perfect (and in particular pseudo-coherent) over $A[x_1, x_2]$, as can e.g. be reduced to $A = \Z$ by base-change, where it is clear.

·It remains to prove (ii) so let the ring map $A \to B$ be given. Let us say that $B$ is \emph{almost compact} over $A$ if it satisfies the colimit condition in (ii), so the claim is that $B$ is almost compact over $A$ if and only if it is almost finitely presented.

First assume that $A \to B$ is almost compact. By considering the case of $0$-truncated rings $C_i$ and using \cite[Lemma~00QO]{stacks-project} we see that the map $\pi_0 A \to \pi_0 B$ is finitely presented. Fix a map $g\colon A[x_1, \dots, x_m] \to B$ which is surjective on $\pi_0$. Note that for every $B$-algebra $C$ we have
\begin{align*}
    \Hom_{\Alg_{A[x_1, \dots, x_m]}}(B, C) = \Hom_{\Alg_A}(B, C) \times_{\Hom_{\Alg_A}(A[x_1, \dots, x_m], C)} \{ g \}.
\end{align*}
Since filtered colimits of anima commute with fiber products and both $A \to B$ and $A \to A[x_1, \dots, x_m]$ are almost compact, it follows that $A[x_1, \dots, x_m] \to B$ is almost compact. Replacing $A$ by $A[x_1, \dots, x_m]$ we can assume that $\pi_0 A \to \pi_0 B$ is surjective. By (iv) it is now enough to show that $B$ is pseudocoherent as an $A$-module. For fixed $n \ge 0$, it follows easily from the definition of almost compactness that there is a compact $A$-algebra $B'$ such that $\tau_{\le n} B \isom \tau_{\le n} B'$. But note that $B'$ is almost finitely presented over $A$: As this property is stable under pushouts, coproducts and retracts by (i), we reduce to the case of the compact generator $A[x]$, which is clear. Then by (iv) we see that $B'$ is pseudocoherent and as this is true for all $n$ we deduce that $B$ is pseudocoherent.

Conversely, assume that $A \to B$ is almost finitely presented. To show that it is almost compact, it is enough to construct for every integer $n \ge 0$ a compact $A$-algebra $B_n$ together with a map $B_n \to B$ such that for $K_n := \cofib(B_n \to B)$ we have $\tau_{\le n} K_n = 0$. To construct $B_n$ we argue by induction on $n$. In the case $n = 0$ we pick a map $A[x_1, \dots, x_n] \to B$ which is surjective on $\pi_0$ and set $B_0 := A[x_1, \dots, x_n]$. For the induction step, fix $n \ge 0$. The map $B_n \to B$ is surjective on $\pi_0$ and by (iii) it is almost finitely presented, hence by (iv) we deduce that $B$ is pseudocoherent as a $B_n$-module; the same is then true for $K_n$. In particular $\pi_{n+1} K_n$ is finitely generated as a $\pi_0 B_n$-module. Each generator gives an element in $\pi_n B_n$ together with a nullhomotopy in $B$. As in the proof of \cite[Proposition~3.2.18]{DAG-Lurie} we can use this data to attach finitely many $n$-cells to $B_n$ (i.e. form the pushout $B_n \tensor_{\Sym_A(A[n])} A$, where $\Sym_A(A[n])$ is the free $A$-algebra on the module $A[n]$) in order to arrive at the desired ring $B_{n+1}$.
\end{proof}

\begin{rmk}
\Cref{rslt:afp-same-as-truncated-compact} was explained to us by Ishan Levy and can be rephrased as saying that a ring map $A \to B$ is almost finitely presented if and only if $\tau_{\le n} B$ is compact in the category of $n$-truncated $A$-algebras for all $n \ge 0$. The related notion of $B$ being compact in $\Alg_A$ (without truncations) is often called \emph{locally finitely presented}.
\end{rmk}

Having discussed the relevant finiteness properties of rings and modules, we now come to the definition of flat, smooth, and étale morphisms of rings. Let us start with the discussion of flatness.

\begin{defn}
Let $A$ be a ring.
\begin{defenum}
    \item An $A$-module $M \in \D(A)$ is called \emph{flat} if the functor $- \tensor_A M\colon \D(A) \to \D(A)$ is t-exact. $M$ is called \emph{faithfully flat} if the functor $- \tensor_A M$ is additionally conservative, i.e. for all $A$-modules $N$, $N \tensor_A M = 0$ implies $N = 0$.

    \item An $A$-algebra $B$ is called \emph{(faithfully) flat} if $B$ is so as an $A$-module.

    \item A ring map $A \to B$ is called \emph{fppf} if it is faithfully flat and almost finitely presented.
\end{defenum}
\end{defn}

\begin{lem}
Let $A$ be a ring.
\begin{lemenum}
    \item \label{rslt:flat-module-in-terms-of-pi-n} An $A$-module $M$ is (faithfully) flat if and only if $H^0(M)$ is (faithfully) flat as a $\pi_0 A$-module and for all $n \in \Z$ the natural map
    \begin{align*}
        H^n(A) \tensor_{\pi_0 A} H^0(M) \isoto H^n(M)
    \end{align*}
    is an isomorphism of $\pi_0 A$-modules.

    \item \label{rslt:flat-map-in-terms-of-pi-n} An $A$-algebra $B$ is (faithfully) flat if and only if the map $\pi_0 A \to \pi_0 B$ is a (faithfully) flat map of classical rings and for all $n \ge 0$ the natural map
    \begin{align*}
        \pi_n A \tensor_{\pi_0 A} \pi_0 B \isoto \pi_n B
    \end{align*}
    is an isomorphism of $\pi_0 A$-modules.
\end{lemenum}
\end{lem}
\begin{proof}
Part (ii) is an immediate consequence of (i). To prove (i), let us first assume that $M$ is flat. Then $- \tensor_A M$ preserves $\D^{\le0}(A)$ and plugging in $A$ shows that $M \in \D^{\le0}(A)$. Using t-exactness of $- \tensor_A M$ and \cite[Corollary~7.2.1.23(2)]{lurie-higher-algebra} we deduce for all $A$-module $N$ and all $n \in \Z$:
\begin{align*}
    H^n(N \tensor_A M) = H^n(N) \tensor_A M = H^0(H^n(N) \tensor_A M) = H^n(N) \tensor_{\pi_0 A} H^0(M).
\end{align*}
Plugging in $N = A$ yields $H^n(A) \tensor_{\pi_0 A} H^0(M) = H^n(M)$, as desired. Conversely, if $H^0(M)$ is flat over $\pi_0 A$ and the above identity holds for all $n \in \Z$, then the Tor spectral sequence shows that $M$ is flat over $A$ (see \cite[Proposition~7.2.2.13]{lurie-higher-algebra}). The faithfully flat version of the claim follows easily by using the above identity for $H^n(N \tensor_A M)$.
\end{proof}

As in classical algebraic geometry, we have faithfully flat descent of modules over a ring, in the following sense:

\begin{prop}
Let $f\colon A \to B$ be a faithfully flat map of rings and let $f^\bullet\colon A \to B^\bullet$ be the induced cosimplicial object, i.e. $B^n = B^{\tensor_A (n+1)}$ for all $n \ge 0$.
\begin{propenum}
    \item \label{rslt:fpqc-descent-for-QCoh} $\D(-)$ descends along $f$, i.e. the natural functor
    \begin{align*}
        - \tensor_A B^\bullet\colon \D(A) \isoto \varprojlim_{n\in\Delta} \D(B^n)
    \end{align*}
    is an equivalence of categories. The same is true for the full subcategories of flat and of pseudo-coherent modules, respectively.

    \item \label{rslt:fpqc-site-is-subcanonical} For every ring $C$, the functor $\Hom(C, -)$ descends along $f$, i.e. the natural map
    \begin{align*}
        \Hom(C, A) \isoto \varprojlim_{n\in\Delta} \Hom(C, B^n)
    \end{align*}
    is an isomorphism of anima. Here $\Hom$ denotes the homomorphism anima in rings.
\end{propenum}
\end{prop}
\begin{proof}
For the first part of (i) see e.g. \cite[Corollary~D.6.3.3]{lurie-higher-algebra} (it ultimately reduces the claim to the case of classical rings, where it is a classical statement). It remains to show descent for flat and for pseudo-coherent modules. In both cases, by descent for $\D(-)$ it is enough to show that a module $M \in \D(A)$ is flat resp. pseudo-coherent if and only if this is true for $M \tensor_A B \in \D(B)$. For flatness this is obvious, and the ``only if'' part of pseudo-coherence is also clear. Now assume that $M \tensor_A B$ is pseudo-coherent as a $B$-module; we need to show that then $M$ is pseudo-coherent as an $A$-module. Given any $A$-module $N$, let $N^\bullet$ denote the cosimplicial $A$-module with $N^n = N \tensor_A B^n$. Then by descent we have $N = \varprojlim_{n\in\Delta} N^n$. Now let $N = \varinjlim_i N_i$ be a filtered and uniformly left-bounded colimit of $A$-modules. We compute
\begin{align*}
    &\Hom_A(M, \varinjlim_i N_i) = \varprojlim_{n\in\Delta} \Hom_{B^n}(M^n, \varinjlim_i N_i^n) = \varprojlim_{n\in\Delta} \varinjlim_i \Hom_{B^n}(M^n, N_i^n) =\\&\qquad\qquad= \varinjlim_i \varprojlim_{n\in\Delta} \Hom(M^n, N_i^n) = \varinjlim_i \Hom(M, N_i).
\end{align*}
Here the first step follows from descent, the second step follows from pseudo-coherence of $M^n$ as a $B^n$-module (note that $(N^n_i)_i$ is still uniformly left-bounded), and the third step follows from the fact that $\varprojlim_{n\in\Delta}$ commutes with uniformly left-bounded filtered colimits of spectra (because it is computed by a spectral sequence).

For (ii) we need to show that $A \isoto \varprojlim_{n\in\Delta} B^n$ is an isomorphism in the category of rings. Since the forgetful functor $\Ring \to \D(\Z)$ is conservative and preserves all small limits, we can check the desired isomorphism in $\D(\Z)$. But then it follows from (i), as (i) implies that the natural map $M \isoto \varprojlim_{n\in\Delta} M \tensor_A B^n$ is an isomorphism in $\D(A)$ (and hence in $\D(\Z)$) for every $A$-module $M$ (see \cite[Lemma~D.4.7(ii)]{heyer-mann-6ff}) and in particular for $M = A$.
\end{proof}

We now come to the definition of smooth and étale morphisms. Since both of these properties imply flatness, by \cref{rslt:flat-map-in-terms-of-pi-n} it is not surprising that they can be checked on $\pi_0$:

\begin{defn}
A ring map $A \to B$ is called \emph{smooth} (resp. \emph{étale}) if it is flat and the induced map $\pi_0 A \to \pi_0 B$ is smooth (resp. étale).
\end{defn}

\begin{rmk}
One can show that every smooth ring map is locally finitely presented (see e.g. \cite[Corollary~3.58]{derived-f-zips}) and in particular almost finitely presented.
\end{rmk}

We have discussed flat modules and perfect modules. In the following we introduce a few more basic finiteness properties and state some of their basic properties, to be used later.

\begin{defn}
Let $A$ be a ring and $P \in \D(A)$.
\begin{defenum}
    \item $P$ is called \emph{projective} if it is connective and for all $M \in \D^{<0}(A)$ we have $\pi_0 \Hom(P, M) = 0$. $P$ is called \emph{finite projective} if it is projective and $H^0(P)$ is finitely generated over $\pi_0 A$.

    \item Given numbers $a \le b \in \Z \cup \{ \pm\infty \}$, we say that $P$ has \emph{Tor amplitude} in $[a, b]$ if for all $M \in \D^\heartsuit(A)$ we have $P \tensor_A M \in \D^{[a,b]}(A)$. Note that $P$ is flat if and only if it has Tor amplitude in $[0,0]$.
\end{defenum}
\end{defn}

\begin{lem} \label{rslt:basic-properties-and-descent-for-finite-projective-modules}
Let $A$ be a ring and $P \in \D(A)$. Then the following are equivalent:
\begin{lemenum}
    \item $P$ is finite projective.
    \item $P$ is a retract of $A^n$ for some $n \ge 0$.
    \item $P$ is pseudo-coherent and flat.
    \item $P$ is flat and $H^0(P)$ is finite projective as a $\pi_0 A$-module.
    \item There exist elements $f_1, \dots, f_n \in \pi_0 A$ generating the unit ideal, such that $P \tensor_A A[1/f_i]$ is finite free as an $A[1/f_i]$-module for all $i$ (see \cref{def:standard-open-immersion} for the definition of $A[1/f_i]$).
\end{lemenum}
In particular, the property of being finite projective descends along faithfully flat covers of rings.
\end{lem}
\begin{proof}
It is clear that (ii) implies (i). For the converse, assume that $P$ is finite projective and pick a map $A^n \to P$ which is surjective on $H^0$. Note that $K := \cofib(A^n \to P)$ lies in $\D^{<0}(A)$, hence by projectivity of $P$ the canonical map $P \to K$ is zero. This implies that the identity $P \to P$ factors over $A^n \to P$, as desired. For the equivalence of (i) and (iv) see \cite[Proposition~7.2.2.18]{lurie-higher-algebra}. Clearly (iii) implies (iv) and (i) implies (iii).

It remains to prove the equivalence of (v) with the other four conditions. First assume (v) and set $B = \prod_i A[1/f_i]$. Then the map $A \to B$ is a faithfully flat and $P \tensor_A B$ is flat and pseudo-coherent as a $B$-module, so by \cref{rslt:fpqc-descent-for-QCoh} the same is true for $P$ as an $A$-module, i.e. $P$ satisfies (iii). Conversely, assume that $P$ is finite projective. By the classical theory (see e.g. \cite[Lemma~00NX]{stacks-project}), we can find $f_i \in A$ generating the unit ideal such that for $P_i := P \tensor_A A[1/f_i]$ we know that $H^0(P)$ is finite free as a module over $\pi_0 A_i$, for $A_i := A[1/f_i]$. For each $i$ we pick a map $A_i^m \to P$ which is an isomorphism on $H^0$. By flatness of $P$ and \cref{rslt:flat-module-in-terms-of-pi-n} we conclude that this map is an isomorphism $A_i^m \isom P$.
\end{proof}

\begin{lem} \label{rslt:Tor-amplitude-of-perfect-dual}
Let $A$ be a ring and $P \in \D(A)$ a perfect $A$-module. Pick $a \le b$ in $\Z \cup \{ \pm \infty \}$. If $P$ has Tor amplitude in $[a, b]$, then $P^\vee$ has Tor amplitude in $[-b, -a]$.
\end{lem}
\begin{proof}
This follows from the identity $\Hom(N \tensor_A P, M) = \Hom(N, M \tensor_A P^\vee)$ for all $M, N \in \D(A)$. E.g. suppose that $P$ has Tor amplitude in $[-\infty, a]$. To show that $P^\vee$ has Tor amplitude in $[-a, \infty]$, pick any $M \in \D^\heartsuit(A)$. We need to show that $M \tensor_A P^\vee \in \D^{\ge -a}(A)$. But this is equivalent to showing that for all $N \in \D^{<-a}(A)$, $\Hom(N, M \tensor_A P^\vee) = 0$. Now use the $\Hom$ identity above.
\end{proof}

In this paper we will almost exclusively work with a base ring $\Lambda$ which is a classical regular noetherian scheme, e.g. $\Z[1/p]$, $\Q_\ell$, $\Z_\ell$ or $\F_\ell$ for primes $p \ne \ell$. In this setting all lafp schemes over $\Lambda$ are locally noetherian, which simplifies the above finiteness conditions. More precisely, we get the following definition and results.

\begin{defn}
A ring $A$ is called \emph{noetherian} if $\pi_0 A$ is noetherian and for every $n \ge 0$, $\pi_n A$ is finitely generated as a $\pi_0 A$-module.
\end{defn}

\begin{lem} \label{rslt:properties-of-noetherian-rings}
Let $A$ be a noetherian ring.
\begin{lemenum}
    \item An $A$-module $M \in \D(A)$ is pseudocoherent if and only if $M$ is right-bounded and $H^n(A)$ is finitely generated over $\pi_0 A$ for all $n \in \Z$. In particular the t-structure on $\D(A)$ restricts to a t-structure on $\D_\pc(A)$.

    \item A ring map $A \to B$ is almost finitely presented if and only if $\pi_0 B$ is finitely generated as a classical $\pi_0 A$-algebra and for all $n \ge 0$, $\pi_n B$ is finitely generated as a $\pi_0 B$-module.

    \item If $A \to B$ is an almost finitely presented ring map then $B$ is noetherian.
\end{lemenum}
\end{lem}
\begin{proof}
Part (i) is a standard argument by inductively resolving $M$ using finite free $A$-modules, see \cite[Proposition~7.2.4.17]{lurie-higher-algebra}. For (ii) we factor $A \to B$ into a composition $A \to A[x_1, \dots, x_n] \to B$ such that $B$ is pseudocoherent over $A[x_1, \dots, x_n]$. Clearly $A[x_1, \dots, x_n]$ is noetherian, so by replacing it with $A$ we can w.l.o.g.\ assume that $B$ is pseudocoherent as an $A$-module. But then the claim follows immediately from (i). Finally, part (iii) is a direct consequence of (ii).
\end{proof}


\subsection{Stacks} \label{sec:alggeo.stacks}

In this subsection we explain how to glue rings to schemes and stacks. We employ the \enquote{functor of points} approach, as it provides a very quick and clean definition of schemes. We refer the reader to \cite[\S A.4]{heyer-mann-6ff} for a quick introduction to sites and sheaves.

\begin{defn}
Fix a ring $\Lambda$.
\begin{defenum}
    \item We equip the category $\cat C := \Alg_\Lambda^\op$ of $\Lambda$-algebras with the fppf site, which is defined as follows: A sieve $\cat C^0 \subseteq \cat C_{/A}$ is a covering sieve iff there are elements $B_1, \dots, B_n \in \cat C^0$ such that the induced ring map $A \to B_1 \times \dots \times B_n$ is fppf. By restriction we define the fppf site on the full subcategory $\Alg^{\afp,\op}_\Lambda \subseteq \Alg^\op_\Lambda$ of almost finitely presented $\Lambda$-algebras.

    \item We denote by
    \begin{align*}
        \Stk_\Lambda := \Shv(\Alg_\Lambda^\op) \subseteq \Fun(\Alg_\Lambda, \Ani)
    \end{align*}
    the category of sheaves for the fppf site, where we implicitly replace $\Alg_\Lambda$ by the full subcategory $\Alg_\Lambda^\kappa$ of $\kappa$-compact objects for a large regular cardinal $\kappa$ (to avoid set-theoretic issues, as $\Alg_\Lambda$ is a big category). The objects of $\Stk_\Lambda$ are called \emph{stacks} over $\Lambda$. We similarly define
    \begin{align*}
        \Stk_\Lambda^\lafp := \Shv(\Alg_\Lambda^{\afp,\op}) \subseteq \Stk_\Lambda
    \end{align*}
    and call its objects the \emph{locally almost finitely presented} (short \emph{lafp}) stacks over $\Lambda$. It embeds fully faithfully into $\Stk_\Lambda$ and is stable under small colimits and finite limits (see \cref{rslt:embedding-of-sheaves-on-subsite}).

    \item By \cref{rslt:fpqc-site-is-subcanonical} the site $\cat C$ is subcanonical, i.e. the Yoneda functor induces an embedding
    \begin{align*}
        \Alg_\Lambda^\op \injto \Stk_\Lambda, \qquad A \mapsto \Spec A.
    \end{align*}
    An object in the image will be called an \emph{affine scheme}. We denote the full subcategory of affine schemes by $\Sch^\aff_\Lambda \subseteq \Stk_\Lambda$.
\end{defenum}
\end{defn}

Recall that there is a good notion of injective maps and surjective maps of stacks, given by monomorphisms and effective epimorphisms, respectively (see e.g. \cite[Lemma~A.3.9]{mann-p-adic-6-functors} for an explicit description of effective epimorphisms in line with classical definitions of surjective maps of sheaves).

We can now introduce (derived) schemes in a similar way as in \cite[\S2.9]{mann-p-adic-6-functors}: They are the stacks that can be obtained by gluing affine schemes along open immersions.

\begin{defn}
Fix a ring $\Lambda$.
\begin{defenum}
    \item \label{def:standard-open-immersion} A ring map $f\colon A \to B$ is called a \emph{standard open immersion} if $B = A[1/a]$ for some $a \in \pi_0 A$. Here we define $A[1/a] := A \tensor_{\Z[x]} \Z[x, x^{-1}]$ for the map $\Z[x] \to A$ induced by $a$.

    \item A map $U \to X$ in $\Stk_\Lambda$ is called an \emph{open immersion} if after every base-change $U' \to X'$ to some affine scheme $X' = \Spec A'$, $U'$ can be covered by affine substacks $V_i = \Spec B_i \subset U'$ such that the induced maps $A' \to B_i$ are standard open immersions.

    \item A \emph{scheme} over $\Lambda$ is a stack $X \in \Stk_\Lambda$ such that there is an open cover of $X$ by affine schemes. If these affine schemes can be chosen of the form $\Spec A$ for $\Lambda \to A$ almost finitely presented, then we say that $X$ is \emph{locally almost finitely presented} (\emph{lafp}) over $\Lambda$. We denote by
    \begin{align*}
        \Sch_\Lambda^\lafp \subseteq \Sch_\Lambda \subseteq \Stk_\Lambda
    \end{align*}
    the full subcategories of schemes and lafp schemes.

    \item A map $f\colon Y \to X$ in $\Stk_\Lambda$ is called \emph{schematic} if for every (affine) scheme $X'$ with a map $X' \to X$ the fiber product $X' \times_X Y$ is a scheme. In this case we say that $f$ is \emph{locally almost finitely presented (lafp)}, resp. \emph{flat}, resp. \emph{smooth}, resp. \emph{étale} if after every pullback of $f$ to an affine scheme, there is an affine open cover of $Y$ such that the induced maps on rings have the respective property.

    \item A scheme $X$ is called \emph{quasi-compact} if it can be covered by finitely many affine open subschemes. A schematic map $f\colon Y \to X$ is called \emph{quasi-compact} if after pullback to every affine scheme, $Y$ is a quasi-compact scheme. A scheme $X$ is called \emph{quasi-separated} if the diagonal of $X$ is quasi-compact. A schematic map $f\colon Y \to X$ is \emph{quasi-separated} if after pullback to every affine scheme, $Y$ is quasi-separated. We write \emph{qcqs} for quasi-compact and quasi-separated.
\end{defenum}
\end{defn}

Via the Yoneda embedding, any reasonable definition of (derived) schemes is compatible with the above definition and we will use this fact without further mention. The next result shows that there is no ambiguity in the meaning of \enquote{lafp scheme over $\Lambda$}.

\begin{lem} \label{rslt:no-ambiguity-in-lafp-schemes}
Fix a ring $\Lambda$. Then $\Sch_\Lambda^\lafp = \Sch_\Lambda \cap \Stk_\Lambda^\lafp$ as full subcategories of $\Stk_\Lambda$. In other words, a scheme over $\Lambda$ admits an open cover by affine lafp schemes over $\Lambda$ if and only if the scheme is lafp as a stack over $\Lambda$.
\end{lem}
\begin{proof}
To see $\subseteq$, let $X$ be a lafp scheme over $\Lambda$; we need to see that $X$ is a lafp stack over $\Lambda$. Pick a covering $Y \surjto X$ by a disjoint union of lafp affine schemes over $\Lambda$. If $Y_\bullet \to X$ denotes the associated Čech nerve then $X = \varinjlim_{n\in\Delta} Y_n$. But all $Y_n$ are disjoint unions of open subsets of affine lafp schemes, so by \cref{rslt:embedding-of-sheaves-on-subsite} is easy to see that all $Y_n$ are lafp stacks over $\Lambda$. Thus the same is true for $X$.

To see $\supseteq$, let $X$ be a scheme over $\Lambda$ which is a lafp stack. Pick a (not necessarily open) cover $Y \surjto X$ by a disjoint union of lafp affine schemes. Let furthermore $Z = \Spec A \subset X$ be an affine open subset. Then $Z$ induces an open subset of $Y$, which implies that $Z$ can be covered by an affine lafp scheme. Altogether we see that there is an fppf cover $A \to B$ of $\Lambda$-algebras such that $B$ is lafp over $\Lambda$. We have to show that this implies that $A$ is lafp over $\Lambda$. Classically this is \cite[Lemma~02KK]{stacks-project} and in the derived world it is stated in \cite[Proposition~5.3.5]{DAG-Lurie}; see \cite[Lemma~4.20]{RC.DeRhamStacksLecture} for a quick reduction of the derived statement to the classical statement using the characterization of afp maps in terms of the cotangent complex (see \cite[Proposition~3.2.18]{DAG-Lurie}).
\end{proof}

\begin{rmk} \label{rslt:schematic-map-between-lafp-stacks-is-lafp}
By \cref{rslt:no-ambiguity-in-lafp-schemes} if $Y \to X$ is a schematic map of lafp stacks over $\Lambda$ then after pullback to every lafp scheme over $X$, $Y$ becomes a lafp scheme over $\Lambda$. We do not know if the converse holds, i.e. whether one can check that $f$ is schematic by only testing with lafp schemes over $\Lambda$.
\end{rmk}

In this paper we need to work with certain \enquote{nice} stacks, mostly with algebraic stacks in the sense of \cite{stacks-project}. In derived algebraic geometry, we use the following version of algebraic stacks, cf. \cite[Definition~3.2]{DAG-Toen}.

\begin{defn} \label{def:geometric-maps-of-stacks}
Fix a ring $\Lambda$ and a map $f\colon Y \to X$ of stacks over $\Lambda$.
\begin{defenum}
    \item $f$ is \emph{$0$-geometric} if it is schematic.
    
    \item Fix $n > 0$. Then $f$ is \emph{$n$-geometric} if for every $\Lambda$-algebra $A$ with a map $\Spec A \to X$ and pullback $Y' := Y \times_X \Spec A$ there is a smooth $(n-1)$-geometric cover $U \surjto Y'$ such that $U$ is a scheme over $\Lambda$. We say that $f$ is \emph{smooth} if $U$ can be chosen smooth over $\Spec A$.

    \item $f$ is called \emph{geometric} if it is $n$-geometric for some $n \ge 0$. In this case, $f$ is called \emph{locally almost finitely presented (lafp)}, resp. \emph{flat}, resp. \emph{smooth} if after pullback to every affine scheme over $\Lambda$, $Y$ admits a smooth cover $U \surjto Y$ by a scheme $U$ such that $U \to X$ has the respective property.

    \item $X$ is a \emph{geometric stack} over $\Lambda$ if the map $X \to \Spec \Lambda$ is geometric.
\end{defenum}
\end{defn}

Most of the stacks we care about are classical, i.e. can be defined as stacks on classical rings. We now discuss how to view classical stacks in the framework that was set up above.

\begin{defn}
Fix a ring $\Lambda$.
\begin{defenum}
    \item Let $\Alg_\Lambda^\cl = \Alg_{\pi_0\Lambda}^\cl$ denote the category of classical (i.e. static) $\Lambda$-algebras. We equip this category with the fppf topology (and implicitly choose a large cutoff cardinal $\kappa$ as above). A \emph{classical stack over $\Lambda$} is a sheaf of anima on $\Alg_\Lambda^{\cl,\op}$. We denote by
    \begin{align*}
        \Stk_\Lambda^\cl \subseteq \Stk_\Lambda
    \end{align*}
    the full subcategory of classical stacks, where the embedding is induced by \cref{rslt:embedding-of-sheaves-on-subsite}. It is stable under all small colimits.

    \item By \cref{rslt:embedding-of-sheaves-on-subsite} the embedding $\Stk_\Lambda^\cl \injto \Stk_\Lambda$ admits a colimit-preserving right adjoint
    \begin{align*}
        (-)^\cl\colon \Stk_\Lambda \to \Stk_\Lambda^\cl \subseteq \Stk_\Lambda
    \end{align*}
    sending a stack $X$ over $\Lambda$ to its \emph{underlying classical stack} $X^\cl$. Note that for every $\Lambda$-algebra $A$ we have $(\Spec A)^\cl = \Spec \pi_0 A$.
\end{defenum}
\end{defn}

\begin{lem}
Fix a ring $\Lambda$.
\begin{lemenum}
    \item A scheme over $\Lambda$ is classical if and only it admits an open cover by affine schemes of the form $\Spec A$ for $A$ a classical $\Lambda$-algebra.

    \item \label{rslt:flat-geometric-maps-stable-under-cl} Let $f\colon Y \to X$ is a flat $n$-geometric map of stacks over $\Lambda$. Then $f^\cl\colon Y^\cl \to X^\cl$ is a flat $n$-geometric map. Moreover, for every map of stacks $Z \to X^\cl$ we have $Z \times_X Y = Z \times_{X^\cl} Y^\cl$ and if $Z$ is classical then so is $Z \times_X Y$.

    \item \label{rslt:X-cl-preserves-geometric} For every ($n$-)geometric stack $X$ over $\Lambda$, $X^\cl$ is still ($n$-)geometric.  
\end{lemenum}
\end{lem}
\begin{proof}
Part (i) reduces by the same proof as in \cref{rslt:no-ambiguity-in-lafp-schemes} to the following two claims: (a) if $U \subseteq \Spec A$ is an open subset of a classical affine scheme then $U$ is classical, (b) if $A \to B$ is an fppf cover of $\Lambda$-algebras and $B$ is classical, then $A$ is classical. Part (b) follows straight from the definition of fppf covers. For part (a), cover $U$ by standard Zariski open subsets of $\Spec A$ (which are obviously classical). Then $U$ is the colimit of the Čech nerve of this cover, so it is enough to show that all terms in that Čech nerve are classical. But these are computed as fiber products over $U$ and equivalently over $\Spec A$, so the claim is clear.

We now prove (ii), so pick a flat $n$-geometric map $Y \to X$ of stacks over $\Lambda$. We first show the second claim, which reduces to showing that $X^\cl \times_X Y = Y^\cl$. We use induction on the number $n$. By writing $X$ as a colimit of affine schemes and using that $(-)^\cl$ preserves all small colimits, we reduce to the case that $X = \Spec A$ is an affine scheme. In particular $X^\cl = \Spec \pi_0 A$. Pick a smooth $(n-1)$-geometric cover $U \surjto Y$ by a scheme $U$ and let $U_\bullet \to Y$ be the associated Čech nerve. Then all maps $U_n \to X$ are $(n-1)$-geometric and flat, hence by induction we have $X^\cl \times_X U_k = U_k^\cl$ for all $k$. On the other hand, $Y^\cl = \varinjlim_{k\in\Delta} U_k^\cl$, which implies the induction step. It thus remains to handle the case $n = 0$, i.e. $f$ is schematic and hence $Y$ is a scheme. By picking a cover of $Y$ given by a disjoint union of affine open subsets of $Y$ we easily reduce to the case that $Y$ is an open subspace in an affine scheme. Using a further covering of $Y$ we reduce to the case that $Y = \Spec B$ is an affine scheme. But then the statement is simply that $\pi_0 A \tensor_A B = \pi_0 B$, which follows from flatness of $A \to B$.

Let us now prove the third part of (ii), so assume that $Z$ is classical. Since classical stacks are stable under colimits, we can assume that $Z$ is a classical affine scheme and then perform a similar induction step as above to deduce that $Z \times_X Y$ is classical. The first part of (ii) is easy: Pick a map $Z = \Spec A \to X^\cl$ from some affine scheme $Z$. Suppose first that $n > 0$. Since $f$ is $n$-geometric, there is a smooth $(n-1)$-geometric cover of $Z \times_X Y = Z \times_{X^\cl} Y^\cl$ by a scheme over $Z$, as desired. If $n = 0$, i.e. $f$ is schematic, one easily reduces to the case that the pullback of $Y$ is affine, where the claim is clear.

We now prove (iii), so let $X$ be an $n$-geometric stack over $\Lambda$. Pick a smooth $(n-1)$-geometric cover $U \surjto X$ by a scheme over $\Lambda$. By (ii) the map $U^\cl \to X^\cl$ is a smooth $(n-1)$-geometric map, and by the compatibility with fiber products in (ii) it follows easily that this map is still a cover (as $(-)^\cl$ commutes with the formation of the Čech nerve).
\end{proof}

\begin{rmk}
Via the Yoneda embedding we can identify the category of classical schemes over $\Lambda$ (in the above sense) with classical definitions of schemes over the classical ring $\pi_0 \Lambda$.
\end{rmk}

A useful tool when studying schemes (and more generally stacks) is the underlying topological space, which we introduce next. Recall that a faithfully flat map $A \to B$ of classical rings induces an open map on the associated Zariski spectra (see \cite[Proposition~00I1]{stacks-project}). Since it is also surjective, it is a quotient map. This allow us to make the following definition.

\begin{defn} \label{def:topology-on-stacks}
Fix a ring $\Lambda$. We denote by
\begin{align*}
    \abs{-}\colon \Stk_\Lambda \to \Top, \qquad X \mapsto \abs{X}
\end{align*}
the unique colimit-preserving functor sending $\Spec A \mapsto \abs{\Spec \pi_0 A}$ for every $\Lambda$-algebra $A$, where the right-hand side denotes the classical Zariski spectrum. Clearly $\abs X = \abs{X^\cl}$ for all $X$.
\end{defn}

\begin{lem}
Fix a ring $\Lambda$.
\begin{lemenum}
    \item For all maps $Y \to X \from Z$ of stacks over $\Lambda$, the map $\abs{Y \times_X Z} \surjto \abs Y \times_{\abs X} \abs Z$ is surjective.

    \item If $f\colon Y \to X$ is a lafp flat geometric map of stacks then $\abs f\colon \abs Y \to \abs X$ is open.
    
    \item \label{rslt:open-substacks-same-as-open-subsets} Let $X$ be a stack over $\Lambda$. Then the assignment $U \mapsto \abs U$ defines an equivalence between the set of open substacks $U \subseteq X$ and the set of open subsets of $\abs X$.
\end{lemenum}
\end{lem}
\begin{proof}
We start with (i), for which we first prove the following claim: For every stack $X$ over $\Lambda$, the set $\abs X$ can be described as the set of maps $\Spec K \to X$ from (classical) fields $K$ over $\Lambda$ under the relation that two maps $\Spec K, \Spec K' \to X$ are equivalent if there is a field $\Spec L$ inclusions $K, K' \subseteq L$ such that the induced maps $\Spec L \to \Spec K \to X$ and $\Spec L \to \Spec K' \to X$ are isomorphic (cf. \cite[Section~0X4E]{stacks-project}). Namely, pick any cover $Y \surjto X$ such that $Y$ is a disjoint union of affine schemes. Then $\abs Y \to \abs X$ is surjective (even a quotient map), which immediately implies that every point in $\abs X$ lies in the image of $* = \abs{\Spec K} \to \abs X$ for some field $K$ with a map $\Spec K \to X$. Suppose now that we have two maps $\Spec K, \Spec K' \to X$ mapping to the same point of $\abs X$. Without loss of generality we can assume that these two maps lift to $Y$. But $\abs X$ is the quotient of $\abs Y$ by the equivalence relation coming from $\abs{Y \times_X Y}$, which implies that there is a point $x \in \abs{Y \times_X Y}$ which maps to $\abs{\Spec K}$ and $\abs{\Spec K'}$ under the two projections $Y \times_X Y \to Y$. Pick a map $\Spec L \to Y \times_X Y$ for some field $L$ that hits the point $x$. Then we get induced maps $\Spec L \to \Spec K$ and $\Spec L \to \Spec K'$ which agree after projection to $X$. This proves the claimed description of $\abs X$. With this description at hand, (i) is straightforward (see e.g. \cite[Lemma~04XH]{stacks-project}).

We now prove (ii), so let $f\colon Y \to X$ be a flat $n$-geometric map for some $n \ge 0$. We argue by induction on $n$, the case $n = 0$ being classical (cf. \cite[Lemma~01UA]{stacks-project}). Pick a cover $g\colon X' \surjto X$ such that $X'$ is a disjoint union of affine schemes and consider the base-change diagram
\begin{equation*}\begin{tikzcd}
    Y' \arrow[r,"g'"] \arrow[d,"f'"] & Y \arrow[d,"f"]\\
    X' \arrow[r,"g"] & X
\end{tikzcd}\end{equation*}
Given an open subset $U \subseteq \abs Y$, we have $g^{-1}(f(U)) = f'(g'^{-1}(U))$: By (i) we can replace $\abs{Y'}$ by $\abs Y \times_{\abs X} \abs{X'}$ for this question, then the claim becomes obvious. Now since $\abs g\colon \abs{X'} \surjto \abs X$ is a quotient map, the openness of $\abs f$ reduces to the openness of $\abs{f'}$. Altogether this implies that from now on we can assume that $X$ is an affine scheme. But then there is a smooth $(n-1)$-geometric cover $U \surjto Y$ such that $U$ is a scheme which is lafp and flat over $X$. By induction $\abs{U} \to \abs{Y}$ is open and by the scheme case also $\abs U \to \abs X$ is open, which together imply that $\abs f$ is open.

We now prove (iii), so let $X$ be a stack over $\Lambda$. By writing $X$ as a colimit of affine schemes and using the fact that $\abs{-}$ preserves colimits, the claim reduces to the following: (a) the claim is true in the case that $X$ is an affine scheme and (b) if $f\colon Y \to X$ is a map of affine schemes and $U \subseteq X$ an open subscheme with pullback $V \subseteq Y$, then $\abs V = \abs f^{-1}(\abs U)$. But these are standard facts in algebraic geometry (note that the claims can be checked after passing to classical schemes).
\end{proof}

By \cref{rslt:open-substacks-same-as-open-subsets} there is no ambiguity in the phrase \enquote{open substack of a stack $X$}. Next we discuss closed immersions of stacks and their relation to the underlying topological space.

\begin{defn}
Fix a ring $\Lambda$ and a map $f\colon Y \to X$ of stacks over $\Lambda$.
\begin{defenum}
    \item We say that $f$ is \emph{affine} if for every affine scheme $X' = \Spec A$ with a map $X' \to X$, the fiber product $Y' = Y \times_X X'$ is affine, i.e. $Y' = \Spec B$ for some $A$-algebra $B$.

    \item We say that $f$ is a \emph{closed immersion} if it is affine and for each $A$ and $B$ as in (i), the map $A \to B$ induces a surjective map $\pi_0 A \surjto \pi_0 B$. In this case we say that $f$ is \emph{finitary} if the kernel of the map $\pi_0 A \to \pi_0 B$ is finitely generated.

    \item If $f$ is schematic then we say that it is \emph{proper} if after pullback to every affine scheme, the map $Y^\cl \to X^\cl$ is a proper map of classical schemes.\footnote{There is a more conceptual and slightly more general definition of properness directly in the derived world, see e.g. \cite[Definition~2.9.28(c)]{mann-p-adic-6-functors}.}
\end{defenum}
\end{defn}

Clearly every closed immersion is both affine and proper. Also, affine maps are schematic and in particular geometric. Hence every closed immersion of a geometric stack over $\Lambda$ is again a geometric stack over $\Lambda$.

\begin{rmk} \label{rmk:closed-immersions-not-injective}
Unlike in classical algebraic geometry, closed immersions are usually \emph{not} monomorphisms, hence the name \enquote{immersion} is a bit misleading. For example, the inclusion $\{ 0 \} \injto \mathbb A^1$ is not a monomorphism: The fiber product of $\{ 0 \}$ with itself over $\mathbb A^1$ is a derived affine scheme whose ring of functions has non-trivial $\pi_1$ (cf. \cref{rmk:intIsom-often-derived} below).
\end{rmk}

\begin{lem}
Fix a ring $\Lambda$.
\begin{lemenum}
    \item \label{rslt:affineness-and-closed-immersion-are-local-on-target} A map $f\colon Y \to X$ of stacks over $\Lambda$ is affine (resp. a closed immersion, resp. a finitary closed immersion) if and only if it is so on a cover of $X$.

    \item If $Z \to X$ is a closed immersion of stacks over $\Lambda$ then the induced map $\abs Z \injto \abs X$ is injective and its image is closed.

    \item \label{rslt:X-cl-to-X-is-closed-immersion} If $X$ is a geometric stack over $\Lambda$ then the map $X^\cl \to X$ is a finitary closed immersion.
\end{lemenum}
\end{lem}
\begin{proof}
We start with (i) and first prove the claim about affineness. This readily reduces to the following claim: Let $A \to B$ be an fppf cover of $\Lambda$-algebras and let $Y$ be a stack with a map $Y \to \Spec A$ whose base-change to $\Spec B$ becomes an affine scheme; then $Y$ is an affine scheme. Let $A \to B^\bullet$ denote the co-Čech nerve of $A \to B$, so that $B^n = B^{\tensor_A (n+1)}$. For every ring $R$ recall that $\Sch_R^\aff \subseteq \Stk_R$ denotes the full subcategory spanned by the affine schemes over $R$, so that $\Sch_R^\aff = \Alg_R^\op$ (we ignore the set-theoretic bound $\kappa$ here). We claim that the assignment $R \mapsto \Sch_R^\aff$ satisfies fppf descent, i.e. the functor
\begin{align*}
    \Sch_A^\aff \isoto \varprojlim_{n\in\Delta} \Sch_{B^n}^\aff
\end{align*}
is an isomorphism. This follows from descent for $\D(-)$ (as proved in \cref{rslt:fpqc-descent-for-QCoh}) because the above functor admits a natural right adjoint computed via totalizations, so the claim reduces to checking that the unit and counit maps are isomorphisms; but this can be checked on underlying modules. It is an easy standard fact about sites that $\Stk_A \isoto \varprojlim_{n\in\Delta} \Stk_{B^n}$ is also an isomorphism (i.e. \enquote{sheaves on a slice satisfy descent}), and by construction this isomorphism is compatible with the above isomorphism via the embedding of affine stacks into all stacks. Now $Y$ induces an object in $\varprojlim_{n\in\Delta} \Sch_{B^n}^\aff$ and hence must lie in $\Sch_A^\aff$, as desired.

We now prove the claim about closed immersions in (i), which again reduces to the case of an fppf cover $A \to B$ of rings and a stack $Y$ over $\Spec A$ whose pullback to $\Spec B$ becomes a closed immersion; we have to show that then $Y \to \Spec A$ is a closed immersion. By the above we know that $Y = \Spec A'$ is affine, so we only need to show that the map $\pi_0 A \to \pi_0 A'$ is surjective. But this can be verified after tensoring with $B$, as $B$ is faithfully flat over $A$.

Finally, the claim about finitary closed immersions in (i) follows easily from the fact that the property of being finitely presented is fppf local for classical rings (see \cite[Lemma~00QQ]{stacks-project}).

We now prove (ii). As in the proof of \cref{rslt:open-substacks-same-as-open-subsets} the claim reduces to the case of affine schemes as long as one shows an additional pullback compatibility. Hence we can assume that all involved schemes are affine; then the statement only depends on the underlying classical rings, where it is easily verified (see e.g. \cite[Lemma~01HQ]{stacks-project}).

It remains to prove (iii), so let $X$ be a geometric stack over $\Lambda$. Pick a smooth geometric cover $U \surjto X$ by a scheme $U$. By \cref{rslt:flat-geometric-maps-stable-under-cl} the map $U^\cl \to U$ is the base-change of $X^\cl \to X$ along $U \surjto X$, so by (i) it is enough to check that $U^\cl \to U$ is a finitary closed immersion. But this is easy: We can e.g. pass to another cover to assume that $U$ is a disjoint union of affine schemes, then reduce to the case that $U = \Spec A$ is affine and then note that in this case $U^\cl = \Spec \pi_0 A$, so the claim is clear.
\end{proof}

Given a stack $X$ and a closed subset $Z \subseteq \abs X$ it is a natural question to ask whether $Z$ comes as the image of a closed embedding to $X$. At least if $X$ is geometric, the answer is yes and there is even a canonical choice for a closed embedding, given by a \emph{reduced} substack:

\begin{defn}
Fix a ring $\Lambda$.
\begin{defenum}
    \item A $\Lambda$-algebra $A$ is \emph{reduced} if $A$ is classical and reduced in the classical sense (i.e. has no non-zero nilpotent elements). We denote by
    \begin{align*}
        \Alg_\Lambda^\red \subseteq \Alg_\Lambda^\cl \subseteq \Alg_\Lambda
    \end{align*}
    the full subcategory spanned by the reduced $\Lambda$-algebras. Note that $\Alg_\Lambda^\red = \Alg_{\pi_0\Lambda}^\red$.

    \item A scheme over $\Lambda$ is called \emph{reduced} if it admits an open cover by reduced affine schemes. A geometric stack over $\Lambda$ is called \emph{reduced} if it admits a smooth geometric cover by a reduced scheme.
\end{defenum}
\end{defn}

\begin{lem} \label{rslt:reducedness-for-geometric-stacks}
Fix a ring $\Lambda$.
\begin{lemenum}
    \item \label{rslt:reducedness-is-smooth-local} Let $f\colon Y \to X$ be a smooth geometric map of geometric stacks over $\Lambda$. If $X$ is reduced, then so is $Y$. If $f$ is surjective and $Y$ is reduced, then $X$ is reduced.

    \item Let $X$ be a geometric stack over $\Lambda$ and $Z \subseteq \abs X$ a closed subset. Consider the full subcategory $\cat C_{X,Z} \subseteq \Stk_{/X}$ spanned by the closed immersions $Y \to X$ such that $Y$ is reduced and the induced map $\abs Y \injto \abs X$ factors over $Z$. Then $\cat C_{X,Z}$ has a final object $Z^\red$ which is uniquely determined by the requirement $\abs{Z^\red} = Z$.
\end{lemenum}
\end{lem}
\begin{proof}
We first prove (i), so let $f\colon Y \to X$ be a smooth $n$-geometric map of $n$-geometric stacks over $\Lambda$. We prove the claim by induction on $n$. For $n = 0$ the claim is that reducedness of schemes is a smooth local property and this claim reduces immediately to the case of classical schemes, where it is found in \cite[Lemma~034E]{stacks-project}. Now assume $n > 0$ and pick a smooth $(n-1)$-geometric cover $U \surjto X$ by a scheme $U$ with base-change $V \surjto Y$ along $f$. Further pick a smooth $(n-1)$-geometric cover $V' \surjto V$ by a scheme $V'$, so we have the following diagram:
\begin{equation*}\begin{tikzcd}
    V' \arrow[r,two heads] & V \arrow[r,two heads] \arrow[d] & Y \arrow[d,"f"]\\
    & U \arrow[r,two heads] & X
\end{tikzcd}\end{equation*}
Now assume that $X$ is reduced. By induction, $U$ is reduced. Since $V' \to U$ is smooth and $0$-geometric, we deduce that $V'$ is reduced. By induction again, we deduce that $V$ is reduced, hence by definition $Y$ is reduced. Now assume that $f$ is surjective and $Y$ is reduced. By induction $V'$ is reduced. Since $V' \surjto U$ is a smooth surjective map of schemes, $U$ is reduced. Thus $X$ is reduced by definition.

We now prove (ii) so let $X$ and $Z$ be given. We induct on the integer $n$ such that $X$ is $n$-geometric. In the case $n = 0$ the stack $X$ is a scheme and the statement only depends on the underlying classical scheme, where it is well-known. Now assume $n > 0$ and pick a smooth $(n-1)$-geometric cover $U \surjto X$ with Čech nerve $U_\bullet \to X$. For each $n$ denote by $Z_{U_n} \subseteq \abs{U_n}$ the preimage of $Z$. As in every site, the natural map
\begin{align*}
    (\Stk_\Lambda)_{/X} \isoto \varprojlim_{n\in\Delta} (\Stk_\Lambda)_{/U_n}
\end{align*}
is an isomorphism. We note that by (i) and \cref{rslt:affineness-and-closed-immersion-are-local-on-target}, containment in $\cat C_{X,Z}$ can be checked after pullback to $U$. Thus by passing to these full subcategories in the above limit, we deduce that $\cat C_{X,Z} \isoto \varprojlim_{n\in\Delta} \cat C_{U_n,Z_{U_n}}$ is an isomorphism. Each $U_n$ is $(n-1)$-geometric, so by induction we deduce that all $\cat C_{U_n,Z_{U_n}}$ have a final object and by the explicit description of this final object we see that the transition maps in the diagram $(\cat C_{U_n,Z_{U_n}})_n$ preserve it. But this implies that the limit $\cat C_{X,Z}$ also has a final object, namely exactly the object that pulls back to a final object in $\cat C_{U_n,Z_{U_n}}$ for all $n$. This proves the induction step. 
\end{proof}

\begin{defn}
Fix a ring $\Lambda$. For a geometric stack $X$ over $\Lambda$ we denote $X^\red := \abs X^\red$ the \emph{underlying reduced geometric stack} of $X$. We have natural closed immersions $X^\red \to X^\cl \to X$.
\end{defn}

Given two points in a stack one can define the stack of maps between them, cf. \cite[Definition~02ZB]{stacks-project}. In the following we introduce a derived version. We will focus on the case where both points are equal, as this is the only case we need.

\begin{defn} \label{def:intIsom-stack}
Fix a ring $\Lambda$ and a stack $X$ over $\Lambda$. Given a map $x\colon U \to X$ from some stack $U$ we denote
\begin{align*}
    \intIsom_X(x) := [V \mapsto \Hom_{X(V)}(x|_V, x|_V)] \in \Stk_{/U}
\end{align*}
and call it the \emph{stack of isomorphisms} of $X$ at $x$. The precise construction (in the $\infty$-world) will be part of the next result.
\end{defn}

\begin{lem} \label{rslt:basic-properties-of-intIsom-stack}
Fix a ring $\Lambda$, a stack $X$ over $\Lambda$ and a map $x\colon U \to X$ of stacks.
\begin{lemenum}
    \item There is a natural isomorphism $\intIsom_X(x) = U \times_{X \times U} U$, where the two maps $U \to X \times U$ are given by $(x, \id_U)$.

    \item \label{rslt:intIsom-is-group-and-classifying-stack-embeds} $\intIsom_X(x)$ is naturally a stack of groups and there is a natural monomorphism $*/\intIsom_X(x) \injto X \times U$ of stacks over $U$.
\end{lemenum}
\end{lem}
\begin{proof}
We replace $X$ by $X \times U$ in order to assume that $X$ lies in the topos $\cat X = \Stk_{/U}$. Note that $U = *$ is the final object in $\cat X$. Consider the \v{C}ech nerve $G_\bullet$ of the map $x\colon * \to X$. Then $G_\bullet$ defines a group in $\cat X$ with underlying object $G := G_1 = * \times_X *$ (see \cite[Lemma~B.3.2(i)]{heyer-mann-6ff}). Furthermore, the quotient $*/G$ is by definition given as the colimit over the \v{C}ech nerve. This provides a natural map $*/G \to X$, which by \cite[Proposition~6.2.3.4]{lurie-higher-topos-theory} is a monomorphism.

To finish the proof, it remains to identify $G$ with $\intIsom_X(x)$. But this is straightforward: Given an anima $S$ with an object $x \in S$ and associated map $* \to S$, we have $\Hom_S(x, x) = * \times_S *$. Therefore, for $V \in \Stk_{/U}$ we have $G(V) = * \times_{X(V)} * = \Hom_{X(V)}(x|_V, x|_V)$.
\end{proof}

\begin{rmk} \label{rmk:intIsom-often-derived}
We warn the reader that $\intIsom_X(x)$ will often be a genuinely derived stack, even if $X$ is not. For example, suppose that $X = \mathbb A^1 := \Spec \Z[x]$ and let $x\colon \Spec \Z \to X$ be the inclusion of $0$. We obtain
\begin{align*}
    \intIsom_{\mathbb A^1}(0) = \{ 0 \} \times_{\mathbb A^1} \{ 0 \} = \Spec (\Z \tensor_{\Z[x]} \Z) = \Spec \Z[\varepsilon_1],
\end{align*}
where $\varepsilon_1$ is a free generator in homological degree $1$. The underlying classical scheme of $\intIsom_{\mathbb A^1}(0)$ is $\Spec \Z$ and agrees with the definition in \cite[Definition~02ZB]{stacks-project}. The difference in the derived world is explained as follows: Given a ring $A$, the $A$-points of $\intIsom_{\mathbb A^1}(0)$ parametrize the isomorphisms of $0$ in $\mathbb A^1(A) = \pi_* A$. If $A$ is classical then there is exactly one such isomorphism, but in general these isomorphisms are given by $\pi_1 A$. A ring map $\Z[\varepsilon_1] \to A$ is the same as an element in $\pi_1 A$.
\end{rmk}

\subsection{Quasicoherent sheaves} \label{sec:alggeo.qcoh}

In this subsection we introduce the category of quasicoherent sheaves on a stack by gluing the category of modules over a ring. We refer the reader to \cite[\S A.4]{heyer-mann-6ff} for basic definitions and properties of descent.

\begin{defn}
Fix a ring $\Lambda$.
\begin{defenum}
    \item By \cref{rslt:fpqc-descent-for-QCoh} the assignment $A \mapsto \D(A)$ is an fppf sheaf on $\Alg_\Lambda^\op$ and hence extends uniquely to a sheaf of symmetric monoidal categories
    \begin{align*}
        \QCoh\colon (\Stk_\Lambda)^\op \to \CMon
    \end{align*}
    by descent. Given a stack $X \in \Stk_\Lambda$ we call the objects in $\QCoh(X)$ the \emph{quasicoherent sheaves} on $X$. We will often denote the tensor unit in $\QCoh(X)$ by $\calO_X$. Given a morphism $f\colon Y \to X$ of stacks over $\Lambda$ we denote
    \begin{align*}
        f^*\colon \QCoh(X) \rightleftarrows \QCoh(Y) \noloc f_*
    \end{align*}
    the \emph{pullback} $f^*$ (encoded by the sheaf $\QCoh$) and the \emph{pushforward} $f_*$ (the right adjoint of $f^*$).

    \item \label{def:pseudo-coherent-sheaves-on-stacks} Fix a stack $X$ over $\Lambda$ and a sheaf $M \in \QCoh(X)$. Then $M$ is called \emph{pseudo-coherent} if it is so after pullback to any affine scheme. We denote by
    \begin{align*}
        \QCoh_\pc(X) \subseteq \QCoh(X)
    \end{align*}
    the full subcategory of pseudo-coherent modules. By \cref{rslt:fpqc-descent-for-QCoh} this category satisfies descent.
\end{defenum}
\end{defn}

We first discuss the following base-change result for quasicoherent sheaves. This has a classical incarnation known as \emph{flat base-change}, see e.g. \cite[Lemma~02KH]{stacks-project}. One benefit of working in derived algebraic geometry is that the flatness hypothesis is not needed.

\begin{lem} \label{rslt:qcqs-base-change-for-QCoh}
Let $\Lambda$ be a ring and let
\begin{equation*}\begin{tikzcd}
    Y' \arrow[r,"g'"] \arrow[d,"f'"] & Y \arrow[d,"f"]\\
    X' \arrow[r,"g"] & X
\end{tikzcd}\end{equation*}
be a cartesian diagram in $\Stk_\Lambda$. Assume that $f$ is schematic and qcqs. Then the natural map
\begin{align*}
    g^* f_* \isoto f'_* g'^*
\end{align*}
is an isomorphism of functors $\QCoh(Y) \to \QCoh(X')$.
\end{lem}
\begin{proof}
Suppose first that $X$ and $X'$ are affine schemes. Then $Y$ and $Y'$ are qcqs schemes, so we can pick a finite cover $Y = \bigcup_{i\in I} U_i$ by affine open subschemes $U_i \subseteq Y$. For every subset $J \subseteq I$ we denote $U_J := \bigcap_{i\in J} U_i$ and $j_J\colon U_J \injto X$ and $f_J\colon U_J \to X$ the induced maps. Then $f_* = \varprojlim_{\emptyset \ne J \subseteq I} f_{J*} j_J^*$ and a similar description holds for $f'_*$ (after pulling the affine cover to $Y'$). Since this limit is finite, $g^*$ commutes with it, so we reduce the claim to the case of $f_J$ in place of $f$. Thus we can assume that $Y$ is a qcqs open subset in an affine scheme. Repeating the same procedure again, we arrive at the case that $Y$ is affine (because now all the intersections $U_J$ are affine). But then all the schemes in the cartesian square are affine and the claim is clear.

We now handle the general case. By covering $X'$ with affine schemes we can reduce to the case that $X'$ is an affine scheme. Now write $X = \varinjlim_{i\in I} X_i$ where $I$ is the collection of all ($\kappa$-compact) $\Lambda$-algebras $A$ together with a map $\Spec A \to X$. In particular the map $X' \to X$ is part of $I$. Note that $Y = \varinjlim_i (Y \times_X X_i)$. Thus $\QCoh(X) = \varprojlim_i \QCoh(X_i)$ and $\QCoh(Y) = \varprojlim_i \QCoh(Y \times_X X_i)$. Now the base-change in the case that $X$ and $X'$ are affine schemes implies that $f_*$ is computed componentwise in this limit for $\QCoh$, so in particular it satisfies base-change with respect to pullback along the maps $X_i \to X$. Since $X' \to X$ is one of these maps, we are done.
\end{proof}

\begin{cor} \label{rslt:qcqs-projection-formula--and-cocont-for-QCoh}
Let $\Lambda$ be a ring and $f\colon Y \to X$ a qcqs schematic map of stacks over $\Lambda$. Then $f_*\colon \QCoh(Y) \to \QCoh(X)$ preserves all small colimits and satisfies the projection formula, i.e. for all $M \in \QCoh(X)$ and $N \in \QCoh(Y)$ the natural map
\begin{align*}
    f_* N \tensor M \isoto f_*(N \tensor f^* M)
\end{align*}
is an isomorphism in $\QCoh(X)$.
\end{cor}
\begin{proof}
By \cref{rslt:qcqs-base-change-for-QCoh} both claims can be checked after pullback to a cover of $X$ by affine schemes, so we can assume that $X$ itself is an affine scheme. Thus $Y$ is a qcqs scheme. As in the proof of \cref{rslt:qcqs-base-change-for-QCoh}, covering $Y$ with affine open subschemes we first reduce to the case that $Y$ is a qcqs open subscheme of an affine scheme and then to the case that $Y$ is affine. But then $f_*$ is just a forgetful functor and the claim is clear.
\end{proof}

An important feature of $\QCoh(X)$ is that it admits a nice t-structure, at least if $X$ is a geometric stack:

\begin{lem} \label{rslt:t-structure-on-QCoh}
Fix a ring $\Lambda$. For every geometric stack $X$ there is a uniquely determined t-structure on $\QCoh(X)$ with the following properties:
\begin{lemenum}
    \item If $X = \Spec A$ is an affine scheme then the t-structure on $\QCoh(X) = \D(A)$ agrees with the canonical one.

    \item If $f\colon Y \to X$ is a geometric map then $f^*\colon \QCoh(X) \to \QCoh(Y)$ is right t-exact and $f_*$ is left t-exact. If $f$ is flat, then $f^*$ is t-exact. If $f$ is affine then $f_*$ is t-exact.

    \item The t-structure on $\QCoh(X)$ is left and right complete and filtered colimits are t-exact.

    \item The closed immersion $i\colon X^\cl \to X$ induces an equivalence $i_*\colon \QCoh(X^\cl)^\heartsuit \isoto \QCoh(X)^\heartsuit$.
\end{lemenum}
\end{lem}
\begin{proof}
We induct on the number $n \ge 0$ such that $X$ and $f$ are $n$-geometric. Let us first handle the induction step, so assume $n > 0$. Pick a smooth $(n-1)$-geometric cover $U \surjto X$ by a scheme and let $U_\bullet \to X$ be the associated Čech nerve. Note that all $U_n$ and all maps in $(U_n)_n$ are $(n-1)$-geometric, so we can apply the induction hypothesis to $(U_n)_n$. Now in order to satisfy (ii) we are forced to define the t-structure on $\QCoh(X)$ by letting $\QCoh^{\le0}(X)$ and $\QCoh^{\ge0}(X)$ be the full subcategories of $\QCoh(X)$ spanned by those objects that lie in the respective full subcategory of $\QCoh(U)$ when pulled back along $U \to X$. By induction all maps in the diagram $(U_n)_n$ are t-exact and by descent we have $\QCoh(X) = \varprojlim_{n\in\Delta} \QCoh(U_n)$. From these two observations it is easy to deduce that the above definition indeed defines a t-structure on $\QCoh(X)$, i.e. satisfies the axioms in \cite[Definition~1.2.1.1]{lurie-higher-algebra}. Moreover, it follows that $\QCoh^{\le0}(X) = \varprojlim_{n\in\Delta} \QCoh^{\le0}(U_n)$ (and the same for $\QCoh^{\ge0}$), which easily implies by induction that the t-structure on $\QCoh(X)$ is left and right complete. It is clearly compatible with filtered colimits (by induction), hence (iii) is satisfied.

We now prove (ii). We first observe that if we had chosen a different smooth $(n-1)$-geometric cover $U' \surjto X$ then it would define the same t-structure on $\QCoh(X)$: This follows from the fact that pullback along the projection $U \times_X U' \to U$ is t-exact by induction. Now let $f\colon Y \to X$ be an $n$-geometric map. Pick a diagram containing $U$, $V$ and $V'$ as in the proof of \cref{rslt:reducedness-is-smooth-local}. We first show that $f^*$ is right t-exact. Since $V' \surjto Y$ is a smooth $(n-1)$-geometric cover, this can be checked after pullback along this map, i.e. we need to show that the composed map $V' \surjto Y \to X$ is right t-exact. This map factors as $V' \to U \surjto X$ and the latter map is t-exact, so we reduce to showing that pullback along $V' \to U$ is right t-exact. But this is now a map of schemes and in particular $0$-geometric, so we are done by induction. In the same way we deduce that flatness of $f$ implies t-exactness of $f^*$. Moreover, the right t-exactness of $f^*$ formally implies the left t-exactness of $f_*$ using adjunctions. Finally, to show that $f_*$ is t-exact in the case of affine $f$, we can by \cref{rslt:qcqs-base-change-for-QCoh} reduce to the case that $X$ is an affine scheme; then $Y$ is an affine scheme and the claim is clear.

We next prove (iv). Pick $U \surjto X$ and $U_\bullet \to X$ as above. By \cref{rslt:flat-geometric-maps-stable-under-cl} we get an induced Čech nerve $U_\bullet^\cl \to X^\cl$ by base-changing along $X^\cl \to X$. Letting $i_n\colon U_n^\cl \to U_n$ denote the associated closed immersions then by \cref{rslt:qcqs-base-change-for-QCoh} and the t-exactness of $i_{n*}$ we obtain a map of diagrams $(i_{\bullet*})\colon \QCoh(U_\bullet^\cl)^\heartsuit \to \QCoh(U_\bullet)^\heartsuit$ which is an equivalence by induction. But note that $\QCoh(X)^\heartsuit = \varprojlim_{n\in\Delta} \QCoh(U_n)^\heartsuit$ (and similarly for $X^\cl$), hence the claim follows.

It remains to handle the induction start, so now let $n = 0$. Thus $X$ and $Y$ are schemes. Similar to above we can then further reduce to the case that $X$ and $Y$ are open subsets of an affine scheme and then to the case that they are affine schemes. But then all claims are clear.
\end{proof}

\begin{defn}
Fix a ring $\Lambda$, stack $X$ over $\Lambda$ and a quasicoherent sheaf $M \in \QCoh(X)$.
\begin{defenum}
    \item We say that $M$ is \emph{connective} (resp. \emph{almost connective}) if the pullback of $M$ to every affine scheme lies in $\D^{\le0}$ (resp. $D^-$). We denote by $\QCoh^{\le0}(X)$ (resp. $\QCoh^-(X)$) the full subcategory of connective (resp. almost connective) sheaves in $\QCoh(X)$.

    \item \label{def:flat-qcoh-sheaf-and-vector-bundle} We say that $M$ is \emph{flat} if it is so after pullback to some cover of $X$ by affine schemes. We similarly define the \emph{Tor amplitude} of $M$. We say that $M$ is a \emph{vector bundle} if it is finite projective on some cover of $X$ by affine schemes.
\end{defenum}
\end{defn}

\begin{rmk}
If $X$ is a geometric stack over $\Lambda$, then the definition of $\QCoh^{\le0}(X)$ and $\QCoh^-(X)$ coincides with the one coming from the t-structure in \cref{rslt:t-structure-on-QCoh}.
\end{rmk}

\begin{lem}
Fix a ring $\Lambda$, a stack $X$ over $\Lambda$ and a quasicoherent sheaf $M \in \QCoh(X)$.
\begin{lemenum}
    \item $M$ being flat, having Tor amplitude in $[a, b]$, or being a vector bundle is stable under pullbacks and can be checked on any cover of $X$.

    \item Suppose that $X$ is geometric. Then $M$ is flat if and only if the functor $M \tensor -\colon \QCoh(X) \to \QCoh(X)$ is t-exact for the t-structure from \cref{rslt:t-structure-on-QCoh}. Similarly, $M$ has Tor amplitude in $[a, b]$ if and only if the functor $M \tensor -$ sends $\QCoh^\heartsuit(X)$ to $\QCoh^{[a,b]}(X)$.
\end{lemenum}
\end{lem}
\begin{proof}
Part (i) follows from faithfully flat descent of flat and of finite projective modules, see \cref{rslt:fpqc-descent-for-QCoh} and \cref{rslt:basic-properties-and-descent-for-finite-projective-modules} (and note that descent of Tor amplitude is shown in the same way as for flatness). Part (ii) then follows easily from (i) by passing to any smooth cover of $X$ by affines, using that the pullback along this cover is conservative and t-exact.
\end{proof}

One useful feature of $\QCoh^-(X)$ is that it only depends on the value of $X$ on \emph{truncated} rings. More precisely, we have the following result:

\begin{lem} \label{rslt:QCoh-minus-only-depends-on-truncated-rings}
Fix a ring $\Lambda$ and let $f\colon Y \to X$ be a map of stacks over $\Lambda$. Assume that for every truncated $\Lambda$-algebra $A$ the induced map $Y(A) \isoto X(A)$ is an isomorphism. Then $f^*$ induces an equivalence
\begin{align*}
    f^*\colon \QCoh^-(Y) \isoto \QCoh^-(X).
\end{align*}
\end{lem}
\begin{proof}
Note that $\QCoh^-\colon \Stk_\Lambda^\op \to \Cat$ satisfies descent (because containment in it can be checked fppf locally) and is thus the right Kan extension of the functor $\D^-\colon \Alg_\Lambda \to \Cat$ along the inclusion of affine schemes. To prove the claim, it is thus enough to show that $\D^-$ is itself the right Kan extension from its restriction to truncated $\Lambda$-algebras. By the pointwise formula for Kan extensions this amounts to showing the following: Given any $\Lambda$-algebra $A$ we have $\D^-(A) = \varprojlim_{A \to B} \D^-(B)$, where the limit runs over all maps from $A$ to a truncated $\Lambda$-algebra $B$. If $B$ is $n$-truncated then $A \to B$ factors uniquely over $\tau_{\le n} A$, from which we deduce that the maps $A \to \tau_{\le n} A$ are final among all maps $A \to B$. We are thus left to show that the natural functor
\begin{align*}
    F\colon \D^-(A) \isoto \varprojlim_{n\ge0} \D^-(\tau_{\le n} A), \qquad M \mapsto (M \tensor_A \tau_{\le n} A)_n,
\end{align*}
is an isomorphism. This functor has a right adjoint $G$ given by $(M_n)_n \mapsto \varprojlim_n M_n$. We note that $GF(M) = \varprojlim_n M \tensor_A \tau_{\le n} A = M$ for all $M \in \D^-(A)$, as this can be checked on cohomology and by the right-boundedness of $M$ each cohomology group in the limit stabilizes for large enough $n$. One similarly shows $FG = \id$ using that an element $(M_n)_n$ in $\varprojlim_{n\ge0} \D^-(\tau_{\le n} A)$ each cohomology group stabilizes for large enough $n$.
\end{proof}

\subsection{The cotangent complex} \label{sec:alggeo.cotangent}

In the previous subsection we introduced a general theory of geometric stacks over a ring $\Lambda$ and discussed some basic properties. We now turn our focus towards the infinitesimal geometry of these stacks, which will provide us with powerful tools of decomposing stacks into simpler pieces. We start by introducing the cotangent complex following Lurie's work \cite[\S3.2]{DAG-Lurie}, \cite[\S25.3]{lurie-SAG}. Let us first consider the affine case. We denote by $\RingMod$ the category of pairs $(A, M)$ consisting of a ring $A$ and an $A$-module $M$ (this category is denoted $\operatorname{SCRMod}$ in \cite[Notation~25.2.1.1]{lurie-SAG}).

\begin{defn}
\begin{defenum}
    \item Given a ring $A$ and a connective $A$-module $M \in \D^{\le0}(A)$ we denote by
    \begin{align*}
        A \oplus M \in \Alg_A
    \end{align*}
    the $A$-algebra which is informally obtained by requiring the product of elements in $M$ to be zero. This $A$-algebra is called the \emph{trivial square zero extension} of $A$ by $M$. More formally, the functor $\RingMod \to \Ring$, $(A, M) \mapsto A \oplus M$ commutes with sifted colimits and can hence be constructed on the ordinary category of pairs $(A, M)$ where $A$ is a polynomial ring over $\Z$ and $M$ is a finite free $A$-module (see \cite[Construction~25.3.1.1]{lurie-SAG}).

    \item Given a ring map $A \to B$ we denote
    \begin{align*}
        \Der_A(B, -)\colon \D^{\le0}(B) \to \Ani, \qquad M \mapsto \Hom_{(\Alg_A)_{/B}}(B, B \oplus M)
    \end{align*}
    and call $\Der_A(B, M)$ the anima of \emph{$A$-derivations of $B$ in $M$}.

    \item Given a ring map $A \to B$ the functor $\Der_A(B, -)$ is accessible and preserves small limits, hence by the adjoint functor theorem it is corepresentable by a $B$-module $L_{B/A} \in \D^{\le0}(B)$, called the \emph{relative cotangent complex of $B$ over $A$}. Thus for every connective $B$-module $M$ we have
    \begin{align*}
        \Hom_B(L_{B/A}, M) = \Der_A(B, M).
    \end{align*}
\end{defenum}
\end{defn}

\begin{rmk} \label{rmk:cotangent-complex-derived-version-of-Kaehler-diff}
If $A$ is a classical ring and $M$ is a classical $A$-module then one checks that the above definition of $\Der_A(B, M)$ coincides with the classical definition of $A$-derivations $B \to M$. Thus the above definition of the cotangent complex $L_{B/A}$ is a derived generalization of the classical definition of Kähler differentials $\Omega_{B/A}$. 
\end{rmk}

\begin{lem}
Fix a ring $A$. The assignment $B \mapsto (B, L_{B/A})$ defines a functor $L_{-/A}\colon \Alg_A \to \RingMod$ which is uniquely determined by the following two properties:
\begin{lemenum}
    \item \label{rslt:cotangent-complex-for-rings-preserves-sifted-colimits} $L_{-/A}$ preserves all sifted colimits.
    \item \label{rslt:cotangent-complex-of-polynomial-algebra} For $n \ge 0$ we have $L_{A[x_1, \dots, x_n]/A} = A[x_1,\dots,x_n]^n$.
\end{lemenum}
\end{lem}
\begin{proof}
Let $\RingMod^{\le0,\vee}_A \to \Alg_A$ be the cocartesian unstraightening of the functor $B \mapsto \D^{\le0}(B)^\op$, so that an object in $\RingMod^{\le0,\vee}_A$ is a pair $(B, M)$ consisting of an $A$-algebra $B$ and a connective $B$-module $M$, and a morphism $(B, M) \to (C, N)$ consists of a ring map $B \to C$ and a $B$-module map $N \to M$. The anima $\Der_A(B, M)$ is functorial in $(B, M)$ and hence defines a functor $\Der_A\colon \RingMod^{\le0,\vee}_A \to \Ani$, equivalently a section $\Alg_A \to \Fun_A(\RingMod^{\le0,\vee}_A, \Ani)$, where $\Fun_A$ denotes the relative functor category over $\Alg_A$. Now $L_{B/A}$ pointwise represents this section and thus by the relative Yoneda embedding induces a functor $\Alg_A \to \RingMod^{\le0}_A$, where $\RingMod^{\le0}_A \to \Alg_A$ is the unstraightening of $B \mapsto \D^{\le0}(B)$. We thus obtain the desired functor $L_{-/A}$.

We now prove the remaining claims. We start with (i), which amounts to the following: Given a sifted colimit $B = \varinjlim_i B_i$ in $\Alg_A$ and a $B$-module $M$, the map $\Der_A(B, M) \isoto \varprojlim_i \Der_A(B_i, M)$ is an isomorphism. But this follows easily from the definition and the observation $B_i \oplus M = B_i \times_B (B \oplus M)$. Part (ii) follows immediately from the definitions and the universal property of the polynomial ring. Finally, we observe that $\Alg_A$ is compact projectively generated by the polynomial rings over $A$ (this reduces to the case $A = \Z$, where it is true by construction). Thus (i) and (ii) indeed determine $L_{-/A}$ uniquely.
\end{proof}

A very important property of the cotangent complex is that it often allows inductive argument of over the truncations of a ring, by the following result:

\begin{lem} \label{rslt:truncation-of-rings-via-square-zero-extension}
Let $A$ be a ring. Then for every $n \ge 0$ there is a pullback square of rings
\begin{equation*}\begin{tikzcd}
    \tau_{\le n+1} A \arrow[r] \arrow[d] & \tau_{\le n} A \arrow[d]\\
    \tau_{\le n} A \arrow[r] & \tau_{\le n} A \oplus \pi_{n+1} A[n+2]
\end{tikzcd}\end{equation*}
Here the bottom hoizontal map is the canonical zero section.
\end{lem}
\begin{proof}
By \cite[Definition~3.1.1]{DAG-Lurie} we need to show that $\tau_{\le n+1} A$ is a small extension of $\tau_{\le n} A$ by $\pi_{n+1} A[n+1]$. By applying \cite[Proposition~3.3.3]{DAG-Lurie} to $B = \tau_{\le n+1} A$ and $I = \pi_{n+1} A$, we see that the map $\tau_{\le n+1} A \to \tau_{\le n} A$ is a square-zero extension of $\tau_{\le n} A$ by $\pi_{n+1} A[n+1]$. By \cite[Proposition~3.3.5]{DAG-Lurie} this implies that it is a small extension, as desired.
\end{proof}

\begin{rmk}
There is also a version of \cref{rslt:truncation-of-rings-via-square-zero-extension} for $\mathbb E_\infty$-rings instead of animated rings, see \cite[Remark~7.4.1.29]{lurie-higher-algebra}.
\end{rmk}

We need the following corollary of \cref{rslt:truncation-of-rings-via-square-zero-extension} below. It is a ring version of the standard fact in derived categories that if one composes $n+1$ maps of complexes in degrees $[0,n]$ which vanish on cohomology then one obtains the zero map.

\begin{cor} \label{rslt:composition-of-phantom-ring-maps-vanishes}
Fix an integer $n \ge 0$ and let $A_0 \to \dots \to A_n$ be maps of $n$-truncated rings such that for all $k \ge 1$ the induced map $\pi_k A_{k-1} \to \pi_k A_k$ vanishes. Then there is a map $\pi_0 A_0 \to A_n$ such that the composition $A_0 \to A_n$ factors as $A_0 \to \pi_0 A_0 \to A_n$.
\end{cor}
\begin{proof}
By the proof of \cref{rslt:truncation-of-rings-via-square-zero-extension} (more precisely, by \cite[Proposition~3.3.3]{DAG-Lurie} applied to $A = A_n$, $B = A_{n-1}$ and $I = \pi_n A_{n-1}$) the fact that the map $\pi_n A_{n-1} \to \pi_n A_n$ is zero implies that the map $A_{n-1} \to A_n$ factors as $A_{n-1} \to \tau_{\le n-1} A_{n-1} \to A_n$. We get the following commuting diagram of rings:
\begin{equation*}\begin{tikzcd}
    A_{n-2} \arrow[r] \arrow[d] & A_{n-1} \arrow[r] \arrow[d] & A_n\\
    \tau_{\le n-1} A_{n-1} \arrow[r] & \tau_{\le n-1} A_{n-1} \arrow[ur]
\end{tikzcd}\end{equation*}
By the same argument as before, the fact that the map $\pi_{n-1} A_{n-2} \to \pi_{n-1} A_{n-1}$ is zero implies that the map $\tau_{\le n-1} A_{n-2} \to \tau_{\le n-1} A_{n-1}$ factors over $\tau_{\le n-2} A_{n-2}$. Thus $A_{n-2} \to A_n$ factors over $\tau_{\le n-2} A_{n-2}$. Inductively we deduce that the composition $A_k \to \dots \to A_n$ factors over $\tau_{\le k} A_k$ for $k = n, \dots, 0$, as desired.
\end{proof}

We wish to extend the relative cotangent complex to stacks. It might not always exist (although we will see below that it exists quite generally), so we introduce the following terminology.

\begin{defn}
Fix a ring $\Lambda$ and a map $f\colon Y \to X$ of stacks over $\Lambda$.
\begin{defenum}
    \item \label{def:derivation-of-stacks} Given an affine scheme $Z = \Spec C$ with a map $\eta\colon Z \to Y$ and a connective $C$-module $M$ we define $\Der_X(Y, M)_\eta \in \Ani$ as the pullback
    \begin{equation*}\begin{tikzcd}
        \Der_X(Y, M)_\eta \arrow[r] \arrow[d] & * \arrow[d]\\
        Y(C \oplus M) \arrow[r] & Y(C) \times_{X(C)} X(C \oplus M)
    \end{tikzcd}\end{equation*}
    Here the right-hand vertical map is given by $\eta \in Y(C)$ and the image of $\eta$ under $Y(C) \to X(C) \to X(C \oplus M)$, where the second map is the canonical zero section. This definition is functorial in $M$ and $\eta$ and thus defines a functor
    \begin{align*}
        \Der_X(Y, -)\colon \QCoh^{\le0}_{/Y} \to \Ani,
    \end{align*}
    where $\QCoh^{\le0}_{/Y} \to \Sch^{\aff,\op}_{/Y}$ denotes the unstraightening of the functor $\QCoh^{\le0}(-)$, i.e. it is the category of pairs $(\eta, M)$ as above. By standard constructions, $\Der_X(Y, -)$ can equivalently be viewed as a section of $\Fun_{/Y}(\QCoh_{/Y}, \Ani) \to \Sch^{\aff,\op}_{/Y}$, by which we denote the unstraightening of the functor $\Fun(\QCoh(-), \Ani)$.

    \item We say that $f$ \emph{admits a cotangent complex} if it satisfies the following conditions for every $\eta\colon \Spec C \to Y$ as in (a):
    \begin{enumerate}[(i)]
        \item The functor $\Der_X(Y, -)_\eta\colon \D^{\le0}(C) \to \Ani$ is corepresentable by an object $(L_{Y/X})_\eta \in \D^-(C)$.
        \item For every map $\alpha\colon \Spec C' \to \Spec C$ the induced map $(L_{Y/X})_{\eta} \tensor_C C' \isoto (L_{Y/X})_{\eta\alpha}$ is an isomorphism.
    \end{enumerate}
    In this case, by basic properties of the Yoneda embedding and the section in (a), the assignment $\eta\colon (L_{Y/X})_\eta$ defines a uniquely determined element
    \begin{align*}
        L_{Y/X} \in \QCoh^-(Y)
    \end{align*}
    that pulls back to $(L_{Y/X})_\eta$ for all $\eta$. We call $L_{Y/X}$ the \emph{relative cotangent complex} of $f$.
\end{defenum}
\end{defn}

One checks that in the case that $X = \Spec A$ and $Y = \Spec B$ are affine then $L_{Y/X}$ exists and coincides with $L_{B/A}$. In the following we provide some basic properties of the relative cotangent complex and in particular show that it exists in suitable generality.

\begin{prop}
Fix a ring $\Lambda$ and a map $f\colon Y \to X$ of stacks over $\Lambda$.
\begin{propenum}
    \item \label{rslt:cotangent-complex-pullback} Let $g\colon X' \to X$ be a map of stacks and let $Y' := Y \times_X X'$ with projection $g'\colon Y' \to Y$. If $L_{Y/X}$ exists then so does $L_{Y'/X'}$ and we have
    \begin{align*}
        L_{Y'/X'} = g'^* L_{Y/X}.
    \end{align*}
    Moreover, if $g$ is surjective and $L_{Y'/X'}$ exists, then $L_{Y/X}$ exists.

    \item \label{rslt:fiber-sequence-for-cotangent-complex} Let $Z \to Y \to X$ be maps of stacks over $\Lambda$ such that $L_{Y/X}$ exists. Then if either of $L_{Z/Y}$ and $L_{Z/X}$ exists, so does the other and we have a fiber sequence
    \begin{align*}
        L_{Y/X}|_Z \to L_{Z/X} \to L_{Z/Y}
    \end{align*}
    in $\QCoh^-(Z)$.

    \item \label{rslt:cotangent-complex-for-geometric-map-and-smoothness} If $f$ is $n$-geometric then $L_{Y/X}$ exists and $L_{Y/X}[n]$ is connective. If $f$ is additionally lafp, then the following is true:
    \begin{enumerate}[(a)]
        \item $L_{Y/X}$ is pseudo-coherent.
        \item $f$ is smooth if and only if $L_{Y/X}$ is perfect and its dual is connective.
    \end{enumerate}

    \item \label{rslt:cotangent-complex-of-limit} Suppose $f$ is the limit $[Y \to X] = \varprojlim_{i \in I} [Y_i \to X_i]$ of some diagram $(Y_i \to X_i)_i$ of maps in $\Stk_\Lambda$. Let $\pi_i\colon Y \to Y_i$ denote the projection and assume that all $L_{Y_i/X_i}$ exist. Then
    \begin{align*}
        L_{Y/X} = \varinjlim_i \pi_i^* L_{Y_i/X_i}.
    \end{align*}
    is a cotangent complex for $f$, provided that it lies in $\QCoh^-$ (e.g. if all $L_{Y_i/X_i}$ lie in $\QCoh^{\le0}$).
\end{propenum}
\end{prop}
\begin{proof}
We note that $\Der_X(Y, -)$ is functorial in $f$, in the following sense. Let $\cat C \to \Fun([1], \Stk_\Lambda)$ denote the cocartesian unstraightening of the functor $[f\colon Y \to X] \mapsto \Sch^{\aff,\op}_{/Y}$ and let $\QCoh^{\le0}_{/-} \to \cat C$ denote the cocartesian unstraightening of the functor $(f, Z) \mapsto \QCoh^{\le0}(Z)$. Explicitly, an object in $\QCoh_{/-}$ is a triple $(f, \eta, M)$ consisting of a map $Y \to X$ of stacks over $\Lambda$, a map $\eta\colon \Spec C \to Y$ for some ring $C$, and a module $M \in \D^{\le0}(C)$. An object in $\cat C$ is a pair $(f, \eta)$ as above and the map $\QCoh^{\le0}_{/-} \to \cat C$ is the projection. Then we get a functor
\begin{align*}
    \Der\colon \QCoh^{\le0}_{/-} \to \Ani, \qquad (f, \eta, M) \mapsto \Der_X(Y, M)_\eta.
\end{align*}
Equivalently, we can view $\Der$ as a section $\cat C \to \Fun_{/-}(\QCoh^{\le0}_{/-}, \Ani)$, which provides us with all the implicit functoriality used in the arguments below.

We first prove (i), so let $g$ and $g'$ be given. First assume that $L_{Y/X}$ exists. For any $\eta'\colon \Spec C \to Y'$ let $\eta := g \eta'\colon \Spec C \to Y$. Then we naturally have $\Der_X(Y,-)_\eta = \Der_{X'}(Y',-)_{\eta'}$. Since this is true for all $\eta'$, we deduce that $L_{Y'/X'}$ exists and it has the claimed description. To prove the second part of (i), we now assume that $g$ is surjective and $L_{Y'/X'}$ exists. The claim reduces to the case that $X = \Spec A$ and $X' = \Spec A'$ are affine and $A \to A'$ is an fppf cover. Let $g_\bullet\colon X'_\bullet \to X$ denote the Čech nerve of $g$ with pullback $g'_\bullet\colon Y'_\bullet \to Y$. Let furthermore $\eta\colon \Spec C \to Y$ be arbitrary and let $\eta_\bullet\colon \Spec C'^\bullet \to Y'_\bullet$ denote the pullbacks. Note that for every $M \in \D^{\le0}(C)$ and all $n \ge 0$ map $C \oplus M \to C'^n \oplus (M \tensor_C C'^n)$ is an fppf cover and it is obtained as a base-change of $C \to C'^n$. Since $X$ and $Y$ are fppf sheaves, we deduce that the natural map
\begin{align*}
    \Der_X(Y, M)_\eta \isoto \varprojlim_{n\in\Delta} \Der_{X_n}(Y_n, M \tensor_C C'^n)_{\eta'_n}
\end{align*}
is an isomorphism for all $M$. By assumption all $\Der_{X_n}(Y_n,-)_{\eta'_n}$ are corepresentable and the corepresenting objects are compatible with pullback. Thus the above limit implies that $\Der_X(Y,-)_\eta$ is corepresentable by the unique descent of the corepresenting objects for $n \ge 0$. Clearly this corepresenting object is almost coconnective, as this can be checked after pullback along an fppf cover. This finishes the proof of (i).

Part (ii) is straightforward, see \cite[Proposition~3.2.12]{DAG-Lurie}.

For the first part of (iii) we refer the reader to the proof of \cite[Proposition~5.1.5]{DAG-Lurie} which works verbatim in our setting. The basic idea is as follows. We use induction on $n$ (the base case being covered by the existence of the cotangent complex for affine maps as noted above). By (i) we can reduce to the case that $X$ is affine, then pick a smooth geometric cover $U \surjto X$ and a map $\eta\colon \Spec C \to Y$. We need to see that $\Der_X(Y,-)_\eta$ is corepresentable compatibly with pullbacks. By a similar argument as in the proof of (i) this can be checked locally on $\Spec C$, so we may assume that $\eta$ factors over $U$. Now one uses the existence of $L_{U/X}$ and $L_{U/Y}$ and the expected fiber sequence in (ii) pulled back to $\Spec C$ (which is reflected by a similar fiber sequence on $\Der$) to deduce the existence of $L_{Y/X}$.

We now prove the claims (a) and (b) in (iii), so assume that $f$ is lafp and fix $n \ge 0$ such that $f$ is $n$-geometric. In the case that $f$ is affine, (a) is proved in \cite[Proposition~3.2.18]{DAG-Lurie} while (b) is proved in \cite[Proposition~3.4.9]{DAG-Lurie}. From this one easily deduces the case that $f$ is schematic, i.e. $0$-geometric. We now argue by induction on $n$. After base-change we may assume that $X$ is an affine scheme. Pick a smooth $(n-1)$-geometric cover $g\colon U \surjto Y$, so that by (ii) we get the fiber sequence $g^* L_{Y/X} \to L_{U/X} \to L_{U/Y}$ in $\QCoh(U)$. By induction we know that $L_{U/Y}$ is perfect and has connective dual. By the affine case we know that $L_{U/X}$ is pseudo-coherent, so altogether we see that $g^* L_{Y/X}$ is pseudo-coherent. By descent of pseudo-coherence (see \cref{def:pseudo-coherent-sheaves-on-stacks}) we deduce that $L_{Y/X}$ is pseudo-coherent, proving (a). To prove (b), we observe that $f$ is smooth iff $fg$ is smooth iff $L_{U/X}$ is perfect with connective dual (by the case $n = 0$) iff $g^*L_{Y/X}$ is perfect with connective dual (by the fiber sequence) iff $L_{Y/X}$ is perfect with connective dual (by flatness of $f$). This proof (b) and hence finishes the proof of (iii).

To prove (iv) we can first replace $Y_i \to X_i$ by its base-change along $X \to X_i$ and use (i) in order to assume that $X = X_i$ for all $i$. Now the claim follows easily by the universal property of $L_{Y/X}$ (cf. \cite[Proposition~3.2.9]{DAG-Lurie}).
\end{proof}

We record the following important special case of \cref{rslt:cotangent-complex-for-geometric-map-and-smoothness}, characterizing smooth maps of schemes in terms of the cotangent complex.

\begin{cor} \label{rslt:smooth-iff-cotangent-complex-is-vector-bundle}
Fix a ring $\Lambda$ and a schematic map $f\colon Y \to X$ of stacks over $\Lambda$. Then $f$ is smooth if and only if it is lafp and $L_{Y/X} \in \QCoh(Y)$ is a vector bundle.
\end{cor}
\begin{proof}
By \cref{rslt:cotangent-complex-pullback} the cotangent complex exists and is local on $X$, hence the claim is local on $X$. We can therefore assume that $X = \Spec A$ is an affine scheme. Then $Y$ is a scheme and we can further localize the claim on $Y$ in order to assume that $Y = \Spec B$ is also affine. Now $L_{Y/X} \in \D(B)$ is connective and by \cref{rslt:cotangent-complex-for-geometric-map-and-smoothness} we know that $f$ is smooth if and only if $L_{Y/X}$ is perfect and its dual is connective. We need to see that this is the case if and only if $L_{Y/X}$ is finite projective. The ``if'' part is easy. For the ``only if'' part, assume now that $L_{Y/X}$ is perfect and has connective dual. Since $L_{Y/X}$ is connective, it has Tor amplitude in $[-\infty,0]$, and the same is true for its dual. Then by \cref{rslt:Tor-amplitude-of-perfect-dual} $L_{Y/X}$ has Tor amplitude on $[0, \infty]$. Together we deduce that $L_{Y/X}$ has Tor amplitude in $[0,0]$, i.e. it is flat. Then by \cref{rslt:basic-properties-and-descent-for-finite-projective-modules} it follows that $L_{Y/X}$ is finite projective, as desired.
\end{proof}

With the general formalism of the cotangent complex at hand, it is time for examples. Let us first make the following comparison to classical algebraic geometry:

\begin{exmpl} \label{rslt:compare-cotangent-complex-to-Kaehler-diff-on-smooth-schemes}
If $f\colon Y \to X$ is a smooth map of classical schemes, then $L_{Y/X} = \Omega^1_{Y/X}$ is the usual sheaf of Kähler differentials. Indeed, this follows from \cref{rslt:smooth-iff-cotangent-complex-is-vector-bundle} and \cref{rmk:cotangent-complex-derived-version-of-Kaehler-diff}.
\end{exmpl}

Another example of interest to us is the cotangent complex of the classifying stack of an algebraic group. In the following we will show that it is given by the adjoint representation, so let us first define this representation.

\begin{defn} \label{def:Lie-alg-and-adjoint-rep}
Fix a ring $\Lambda$ and a smooth group scheme $G$ over $\Lambda$. We denote $* = \Spec \Lambda$.
\begin{defenum}
    \item Let $e\colon * \to G$ be the unit section. We define
    \begin{align*}
        \Lie(G) := e^* L_{G/X}^\vee \in \D(\Lambda)
    \end{align*}
    and call it the \emph{Lie algebra} of $G$. By \cref{rslt:smooth-iff-cotangent-complex-is-vector-bundle} it is a finite projective $\Lambda$-module.

    \item \label{def:adjoint-representation} Consider the stack $G/G$, where $G$ acts on itself via conjugation. The section $e$ from (a) upgrades to a map $e\colon */G \to G/G$. We denote
    \begin{align*}
        \Ad_G := e^* L_{(G/G) \, / \, (*/G)}^\vee \in \QCoh(*/G)
    \end{align*}
    and call it the \emph{adjoint representation} of $G$. By \cref{rslt:cotangent-complex-pullback} the underlying $\Lambda$-module of $\Ad_G$ is $\Lie(G)$.
\end{defenum}
\end{defn}

\begin{rmk}
In the setup of \cref{def:Lie-alg-and-adjoint-rep}, we have the following more classical description of $\Lie(G)$ and $\Ad_G$. Pick any map $f\colon Y \to X$ of smooth schemes over $\Lambda$. By \cref{rslt:fiber-sequence-for-cotangent-complex} there is a canonical map $f^* L_X \to L_Y$, which dualizes to a map $L_Y^\vee \to f^* L_X^\vee$. Let us call $T_X := L_X^\vee$ and $T_Y := L_Y^\vee$ the tangent bundles of $X$ and $Y$ (these are vector bundles by \cref{rslt:smooth-iff-cotangent-complex-is-vector-bundle}). We thus obtain a map
\begin{align*}
    df\colon T_Y \to f^* T_X,
\end{align*}
which we call the \emph{differential} of $f$. By \cref{rslt:compare-cotangent-complex-to-Kaehler-diff-on-smooth-schemes} it coincides with the classical definition of differentials in case that $\Lambda$ is a static ring.

With the above notation, we see that $\Lie(G) = T_e G := e^* T_G$ is the tangent space of $G$ at the unit. Moreover, given any $g \in G(\Lambda)$, the induced action $\Ad_G(g)\colon \Lie(G) \to \Lie(G)$ is obtained as follows: Via the conjugation action of $G$ on itself, $g$ induces an automorphism $\Ad_g\colon G \isoto G$ and hence an induced isomorphism $\Ad_g^* L_G \isoto L_G$; by dualizing and restricting along $e$ we obtain $\Ad_G(g)$. Therefore $\Ad_G(g) = (d\Ad_g)_e$, so the adjoint representation has the expected form.
\end{rmk}

We now compute the cotangent complex of the classifying stack of some smooth group scheme $G$. It is known to experts that this cotangent complex can be described using the adjoint representation, but we could not find a reference for this in the literature. The closest we found is \cite[Example~8.6.5]{Khan-DAG}, where the underlying $\Lambda$-module of the cotangent complex is computed to be $\Lie(G)$, but the $G$-action is missing.

\begin{prop} \label{rslt:cotangent-of-BG-is-adjoint-representation}
Fix a ring $\Lambda$ and a smooth group scheme $G$ over $\Lambda$. Denote $* = \Spec \Lambda$ and $BG = */G$. Then
\begin{align*}
    L_{BG} = \Ad_G^\vee[-1] \in \QCoh(*/G).
\end{align*}
\end{prop}
\begin{proof}
We start with a general observation. Suppose $f\colon Y \to X$ is a map of stacks over $\Lambda$ such that $L_{Y/X}$ exists and let $\Delta_f\colon Y \to Y \times_X Y$ be the diagonal. Then $L_{Y/Y \times_X Y}$ exists and there is a canonical isomorphism
\begin{align}
    L_{Y/X} = L_{Y/Y \times_X Y}[-1] \label{eq:cotangent-complex-via-diagonal}
\end{align}
in $\QCoh(Y)$. To prove this, note that by \cref{rslt:cotangent-complex-of-limit} we have $L_{Y \times_X Y / X} = \pi_1^* L_{Y/X} \oplus \pi_2^* L_{Y/X}$, where $\pi_i\colon Y \times_X Y \to Y$ are the projections; in particular the left-hand side exists. We apply \cref{rslt:fiber-sequence-for-cotangent-complex} to the chain of maps $Y \to Y \times_X Y \to X$ in order to arrive at the fiber sequence
\begin{align*}
    L_{Y/X} \oplus L_{Y/X} = \Delta_f^* L_{Y \times_X Y / X} \to L_{Y/X} \to L_{Y / Y \times_X Y}
\end{align*}
in $\QCoh(Y)$; in particular the last term exists. The map $L_{Y/X} \oplus L_{Y/X} \to L_{Y/X}$ is the sum, whose fiber identifies with $L_{Y/X}$. This proves \cref{eq:cotangent-complex-via-diagonal}.

We apply \cref{eq:cotangent-complex-via-diagonal} to the map $BG \to *$ and obtain $L_{BG} = L_{BG/BG^2}[-1]$. We then apply \cref{eq:cotangent-complex-via-diagonal} to $BG \to BG^2$, so we need to understand the diagonal of this map, i.e. we need to compute $BG \times_{BG^2} BG$. By standard computations (cf. \cref{exmpl:convolution-of-classfying-stacks}) this stack can be identified with $G \backslash G^2 / G$, where $G$ acts on $G^2$ diagonally from both sides. By identifying $G \backslash G^2 \isom G$ as schemes, we arrive at $BG \times_{BG^2} BG = G /G $, where $G$ acts via conjugation. Altogether we deduce that
\begin{align*}
    L_{BG} = L_{BG/BG^2}[-1] = L_{BG \, / \, (G/G)}[-2].
\end{align*}
It remains to compute the cotangent complex on the right. But by applying \cref{rslt:fiber-sequence-for-cotangent-complex} to the map $BG \to G/G \to BG$ we see that $L_{BG \, / \, (G/G)} = e^* L_{(G/G) \, / \, BG}[1]$, where $e\colon BG \to G/G$ is the unit section. By definition of $\Ad_G$, we conclude $L_{BG} = \Ad_G^\vee[-1]$.
\end{proof}

We can use the cotangent complex to define a slightly weaker notion than smoothness, which still enjoys many of its nice properties. The following definition is often also called \emph{quasi-smooth}, see e.g. \cite[Definition~3.4.15]{DAG-Lurie}.

\begin{defn}
Fix a ring $\Lambda$ and a map $f\colon Y \to X$ of stacks over $\Lambda$.
\begin{defenum}
    \item Suppose that $f$ is schematic. Then $f$ is a \emph{local complete intersection} if it is lafp and $L_{Y/X}$ has Tor amplitude in $[-1,0]$ (in particular it is perfect by \cite[Proposition~3.2.14]{DAG-Lurie} and \cref{rslt:pseudo-coherent-plus-fin-tor-dim-implies-perfect}).

    \item Suppose $f$ is geometric. Then $f$ is a \emph{local complete intersection} if after pullback to every affine scheme over $X$ there is some smooth cover $U \surjto Y$ by a scheme $U$ such that the map $U \to X$ is a local complete intersection.
\end{defenum}
\end{defn}

\begin{lem}
Fix a ring $\Lambda$.
\begin{lemenum}
    \item Local complete intersection morphisms are stable under composition and base-change. Moreover, among geometric morphisms, being a local complete intersection can be checked locally on the target and smooth locally on the source.

    \item \label{rslt:map-between-smooth-stacks-is-lci} Consider a diagram of geometric morphisms of stacks over $\Lambda$ of the form
    \begin{equation*}\begin{tikzcd}
        Y \arrow[rr,"g"] \arrow[dr,"h",swap] && X \arrow[dl,"f"]\\
        & S
    \end{tikzcd}\end{equation*}
    If $f$ is smooth and $h$ is a local complete intersection, then $g$ is a local complete intersection.
\end{lemenum}
\end{lem}
\begin{proof}
Part (i) follows easily from the definitions and the basic properties of the cotangent complex. For (ii), we can first use (i) to reduce to the case that $X$, $Y$ and $S$ are schemes. Then the claim follows easily from \cref{rslt:fiber-sequence-for-cotangent-complex}.
\end{proof}

The name ``local complete intersection'' is justified by the following alternative characterization that is often taken as the definition:

\begin{prop} \label{rslt:characterization-of-lci-maps-via-regular-closed-immersions}
Fix a ring $\Lambda$ and let $f\colon Y \to X$ be a map of schemes over $\Lambda$. Then $f$ is a local complete intersection if and only if after passing to some open cover of $Y$, there is a factorization of $f$ of the form
\begin{align*}
    Y \xto{i} \tilde X \xto{\tilde f} X
\end{align*}
where $\tilde f$ is smooth and $i$ is a closed immersion with the following property: There exists an integer $n \ge 0$, a map $\tilde X \to \mathbb A^n$ and a cartesian square
\begin{equation*}\begin{tikzcd}
    Y \arrow[r,"i"] \arrow[d] & \tilde X \arrow[d]\\
    \{ 0 \} \arrow[r] & \mathbb A^n
\end{tikzcd}\end{equation*}
\end{prop}
\begin{proof}
This is for example shown in \cite[Proposition~2.3.14]{KhanRydh-VirtualCartier} and essentially also in the paragraphs after \cite[Example~3.4.16]{DAG-Lurie}. We sketch the argument. The ``if'' part readily reduces to showing that the closed immersion $\{ 0 \} \to \mathbb A^n$ is a local complete intersection, which follows from \cref{rslt:map-between-smooth-stacks-is-lci}.

We now prove the ''only if'' part, so assume that $f$ is a local complete intersection. After passing to open covers we can assume that $X = \Spec A$ and $Y = \Spec B$ are affine. Pick elements $b_1, \dots, b_n \in \pi_0 B$ such that the induced map $A[x_1, \dots, x_n] \to B$ is surjective on $\pi_0$. By \cref{rslt:map-between-smooth-stacks-is-lci} we may replace $A$ by $A[x_1, \dots, x_n]$ to assume that $\pi_0 A \to \pi_0 B$ is surjective, i.e. $f$ is a closed immersion. Then $H^0(L_{B/A}) = 0$ (e.g. by \cite[Proposition~3.2.16]{DAG-Lurie}) and hence $L_{B/A}[-1]$ has Tor amplitude in $[0,0]$ and is therefore finite projective by \cref{rslt:basic-properties-and-descent-for-finite-projective-modules}. By localizing further we can assume that $L_{B/A}[-1] \isom B^n$ for some $n \ge 0$. Then a basis for $L_{B/A}[-1]$ induces a representation $A = B/(b_1, \dots, b_n)$, see the proof of \cite[Proposition~2.3.8]{KhanRydh-VirtualCartier} for details.
\end{proof}

A useful property of local complete intersections is that they provide a criterion in terms of dimensions for when a scheme is classical. This will be discussed next.

\begin{defn}
Fix a ring $\Lambda$ and a map $f\colon Y \to X$ of schemes over $\Lambda$ and a point $y \in Y$ with residue field $k$.
\begin{defenum}
    \item The \emph{relative dimension of $f$ at $y$} is the minimum dimension of $U$ for $U$ an open neighbourhood of $y$ in the fiber $Y_{f(y)}$. Note that the relative dimension only depends on the underlying topological spaces (and in particular only on the underlying reduced classical schemes).

    \item Suppose that $f$ is a local complete intersection. Let $L_y := L_{Y/X}|_y \in \D(k)$. This $k$-module has Tor amplitude in $[-1,0]$ and is hence concentrated in these two degrees. The \emph{relative virtual dimension} of $f$ at $y$ is defined as $\dim_k H^0(L_y) - \dim_k H^{-1}(L_y)$.
\end{defenum}
\end{defn}

In terms of \cref{rslt:characterization-of-lci-maps-via-regular-closed-immersions} one can also characterize the relative virtual dimension of $f$ at $y$ as follows. Locally around $y$, factor $f$ as $Y \to \tilde X \to X$ such that $\tilde X \to X$ is smooth of pure dimension $d$ and $Y \to \tilde X$ is a pullback of the map $\{ 0 \} \to \mathbb A^n$ for some $n \ge 0$. Then the virtual dimension of $f$ at $y$ is $d - n$.

\begin{prop} \label{rslt:classicality-result-for-lci-maps}
Fix a Cohen--Macauley classical noetherian ring $A$ and let $f\colon Y \to \Spec A$ be a map of schemes. Suppose that $f$ is a local complete intersection and its relative dimension is at most its relative virtual dimension at each $y \in Y$. Then $Y$ is classical.
\end{prop}
\begin{proof}
This is proved in \cite[Proposition~B.0.1]{Beijing}. We quickly provide an argument. The claim can be checked on an open cover of $Y$, so it is enough to show that it holds in an open neighborhood of every closed point $y$ of $Y$. In particular may assume that $Y = \Spec B$ is affine and by the proof of \cref{rslt:characterization-of-lci-maps-via-regular-closed-immersions} we can write $B = A[x_1, \dots, x_d]/(a_1, \dots, a_n)$ such that the relative virtual dimension of $f$ at $y$ is $d - n$. To prove that $B$ is classical in a neighborhood of $y$ it is enough to show that $(a_1, \dots, a_n)$ is a regular sequence in a neighborhood of $y$. Thus it is enough to show that it is a regular sequence in the local ring of $A[x_1, \dots, x_d]$ at $y$. This local ring is Cohen-Macaulay by \cite[Lemma~00ND]{stacks-project} and we may also replace $A$ by its local ring at the image of $y$ under $f$. By \cite[Lemmas~09CC,~0668]{stacks-project} the claim can be checked after modding out by any regular sequence in $A$. Since $A$ is Cohen-Macaulay, we can find a regular sequence in $A$ whose quotient has dimension $0$, so from now on we may assume that $A$ has dimension $0$. Then the relative dimension of $f$ at $y$ is the dimension of the local ring of $B$ at $y$. This dimension is at least $d - n$ and by our assumptions it is at most $d - n$, hence it is equal to $d - n$. We conclude by \cite[Proposition~00N6]{stacks-project}.
\end{proof}

\subsection{Completion and excision} \label{sec:alggeo.completion-and-excision}

We now come to the completion of a stack at a closed subset and discuss an excision sequence for quasicoherent sheaves. Let us start with a very general and abstract definition and its interpretation in terms of quasicoherent sheaves.

\begin{defn}
Fix a ring $\Lambda$, a stack $X$ over $\Lambda$ and a closed subset $Z \subseteq \abs X$.
\begin{defenum}
    \item \label{def:completion-of-stack-at-closed-subset} We denote by
    \begin{align*}
        X_{\hat Z} \subseteq X
    \end{align*}
    the full substack whose value at a $\Lambda$-algebra $A$ is the full subanima of $X(A)$ spanned by the maps $\Spec A \to X$ such that the induced map $\abs{\Spec A} \to \abs X$ factors over $Z$. We call $X_{\hat Z}$ the \emph{completion of $X$ at $Z$}. If $Z \to X$ is a closed immersion then we denote $X_{\hat Z} := X_{\widehat{\abs Z}}$.

    \item Let $j\colon U \injto X$ be the open complement of $Z$ in $X$. We denote by
    \begin{align*}
        \QCoh(X)_Z \subseteq \QCoh(X)
    \end{align*}
    the full subcategory of those $M \in \QCoh(X)$ such that $j^* M = 0$.
\end{defenum}
\end{defn}

In general, completions at closed subsets can be badly behaved, in the same way that completions of modules at an infinitely generated ideal do not work well. In the finitely generated case, one gets a good general theory, as we will explain in the following. We start with a basic observation that is also found in \cite[Lemma~5.1.5]{Lurie-DAG-XII} (for $\mathbb E_\infty$-rings) and \cite[Proposition~6.7.4]{Gaitsgory-Rozenblyum-indschemes}:

\begin{lem} \label{rslt:basic-properties-of-completion}
Let $\Lambda$ be a ring, $X$ a stack over $\Lambda$ and $Z \subseteq \abs X$ a closed subset.
\begin{lemenum}
    \item \label{rslt:completion-commutes-with-base-change} Let $f\colon X' \to X$ be any map of stacks over $\Lambda$ and let $Z' := f^{-1}(Z)$. Then $X'_{\hat Z'} = X_{\hat Z} \times_X X'$.

    \item \label{rslt:completion-of-affine-is-colimit} Suppose that $X = \Spec A$ is affine and that the open complement $U \subset X$ of $Z$ is quasicompact. Then there is a sequence of closed immersions
    \begin{align*}
        Z_1 \to Z_2 \to \dots \to Z_n \to \dots \to X_{\hat Z} \subseteq X
    \end{align*}
    with image $Z$ such that $X_{\hat Z} = \varinjlim_n Z_n$. Moreover, writing $Z_n = \Spec B_n$ then $B_n$ is perfect as an $A$-module.
\end{lemenum}
\end{lem}
\begin{proof}
We can assume $\Lambda = \Z$. Part (i) is obvious from the definition. We now prove (ii), so let $X = \Spec A$ be affine and $U \subset X$ quasicompact. Then $U$ is covered by finitely many standard open subsets of the form $\Spec A[\frac1{f_i}]$ for $f_1, \dots, f_n \in \pi_0 A$. Let $f\colon X \to \mathbb A^n := \Spec \Z[x_1, \dots, x_n]$ be the map induced by that $f_i$'s. Then $Z = f^{-1}(0) \subset \abs X$ (as one checks easily by looking at the open complements). Using (i) we see that (ii) is stable under base-change, so we can reduce to the case that $X = \mathbb A^n$ and $Z = \{ 0 \}$. One checks easily that the claim in (ii) is stable under finite products in $X$ and $Z$, so we can further reduce to the case $X = \mathbb A^1 = \Spec \Z[x]$. Then define $Z_n = \Spec \Z[x]/x^n$. It remains to verify $X_{\hat Z} = \varinjlim_n Z_n$, the other claims in (ii) are clear. Let us denote $X' := \varinjlim_n Z_n$. There clearly is a map $X' \to X$ and we need to check that for every ring $B$ the map
\begin{align*}
    \varinjlim_n \Hom(\Z[x]/x^n, B) \isoto X_{\hat Z}(B) \subseteq \Hom(\Z[x], B) = \pi_* B
\end{align*}
is an isomorphism, where we denote by $\pi_* B$ the underlying anima of $B$. We follow the argument in \cite[Lemma~5.1.5]{Lurie-DAG-XII} (stated for $\mathbb E_\infty$-rings). Fix an element $b \in \pi_0 B$ and let the fiber of $\Hom(\Z[x]/x^n, B)$ over $b$ be denoted by $P_n$, so that we have to show that $\varinjlim_n P_n = *$ if $b$ is nilpotent and $= \emptyset$ if $b$ is not. Note that $\Z[x]/x^n = \Z[x] \tensor_{x^n,\Z[x]} \Z$, where the map $\Z[x] \to \Z[x]$ is given by $x^n$. Thus $\Hom(\Z[x]/x^n, B) = \fib(\pi_* B \to \pi_* B) := \pi_* B \times_{\pi_* B} *$, where the map $\pi_* B \to \pi_* B$ is given by $b \mapsto b^n$ and the map $* \to \pi_* B$ is given by $0 \in B$. This implies that there is a cartesian square
\begin{equation*}\begin{tikzcd}
    P_n \arrow[r] \arrow[d] & \Hom(\Z[x]/x^n, B) \arrow[r] \arrow[d] & * \arrow[d,"0"]\\
    * \arrow[r,"b"] & \pi_* B \arrow[r,"(-)^n"] & \pi_* B
\end{tikzcd}\end{equation*}
Thus $P_n = * \times_{b^n, B, 0} *$. In particular $P_n = \emptyset$ if $b^n \ne 0$ in $\pi_0 B$. Now assume that there is some $p_n \in P_n$. By looking at the induced long exact sequence on homotopy groups we deduce that there is a canonical isomorphism $\pi_i(P_n, p_n) = \pi_{i+1} B$ for all $i \ge 0$. For $m \ge n$ let $p_m \in \pi_0 P_m$ be the image of $p_n$ under the map $P_n \to P_m$. Then the map $\pi_i(P_n, p_n) \to \pi_i(P_m, p_m)$ corresponds to the map $\pi_{i+1} B \to \pi_{i+1} B$ given by multiplication with $b^{m-n}$. Since $b$ is nilpotent (by the existence of $p_n$) we deduce that for $m$ large enough, the map $\pi_i(P_n, p_n) \to \pi_i(P_m, p_m)$ is zero. Thus if $p_\infty \in \varinjlim_n P_n$ denotes the image of $p_n$ we deduce
\begin{align*}
    \pi_i(\varinjlim_n P_n, p_\infty) = \varinjlim_{m\ge n} \pi_i(P_m, p_m) = 0
\end{align*}
for all $i \ge 0$. Hence $\varinjlim_n P_n = *$, as desired.
\end{proof}

\begin{rmk}
The claim $\varinjlim_n \Spec(\Z[x]/x^n) = \mathbb A^1_{\hat 0}$ can also be seen using condensed mathematics: Given a ring $B$, the value of the left-hand side on $B$ coincides with $\Hom_{\Z[x]}(\Z[[x]], B)$ and this easily shows that the map to $\mathbb A^1$ is injective: One just observes that $\Z[[x]]$ is idempotent in $\D_\solid(\Z[x])$.
\end{rmk}

We now come to the first main result about completions, which tells us how quasicoherent sheaves behave on the completion and in particular provides a useful excision sequence. Similar results have been obtained in \cite[Theorem~5.1.9]{Lurie-DAG-XII} (for $\mathbb E_\infty$-rings) and \cite[Proposition~7.1.3]{Gaitsgory-Rozenblyum-indschemes}.

\begin{prop} \label{rslt:hat-i-pullback-induces-equiv}
Let $\Lambda$ be a ring, $X$ a stack over $\Lambda$, $Z \subseteq \abs X$ a closed subset and $\hat i\colon X_{\hat Z} \injto X$ the completion. Assume that the open complement $j\colon U \injto X$ is a quasicompact map. Then $\hat i^*$ induces an equivalence
\begin{align*}
    \hat i^*\colon \QCoh(X)_Z \isoto \QCoh(X_{\hat Z}).
\end{align*}
\end{prop}
\begin{proof}
By \cref{rslt:completion-commutes-with-base-change} the right-hand side glues along covers of $X$ and one sees easily that the same is true for the left-hand side (i.e. that containment in $\QCoh(X)_Z$ can be checked on a cover of $X$). By writing $X$ as a colimit of affine schemes we can therefore reduce to the case that $X = \Spec A$ is affine. Let $i_n\colon Z_n = \Spec B_n \to X$ be the closed immersions from \cref{rslt:completion-of-affine-is-colimit}, so that we have $X_{\hat Z} = \varinjlim_n Z_n$. In particular $\QCoh(X_{\hat Z}) = \varprojlim_n \QCoh(Z_n)$. Since $B_n$ is perfect as an $A$-module, the functor $i_n^*$ preserve limits and hence admits a left adjoint $i_{n\natural}\colon \QCoh(Z_n) \to \QCoh(X)$. By \cite[Lemma~D.4.7(i)]{heyer-mann-6ff} it follows that $\hat i^*$ admits a left adjoint $\hat i_\natural$ and in particular $\hat i^*\colon \QCoh(X) \to \QCoh(X_{\hat Z})$ preserves all small limits.

Note that $\hat i_\natural$ satisfies base-change with respect to pullback along any qcqs schematic map: By passing to adjoints this reduces to \cref{rslt:qcqs-base-change-for-QCoh}. In particular we deduce that $j^* \hat i_\natural = 0$, i.e. $\hat i_\natural$ factors over $\QCoh(X)_Z$. We also deduce $i_n^* \, \hat i_\natural = i'^*_n$, where we denote by $i'_n\colon Z_n \to X_{\hat Z}$ the natural map. This implies that the unit $\id \isoto \hat i^* \, \hat i_\natural$ is an isomorphism (as this can be checked after pullback along all $i'_n$), so that $\hat i_\natural$ is fully faithful. We have thus constructed the adjunction
\begin{align*}
    \hat i_\natural\colon \QCoh(X_{\hat Z}) \rightleftarrows \QCoh(X)_Z \noloc \hat i^*
\end{align*}
and we already know that $\hat i_\natural$ is fully faithful. To finish the argument it remains to see that $\hat i^*$ is conservative on $\QCoh(X)_Z$, because then the above adjunction is an equivalence. It is enough to show that $i_1^*$ is conservative on $\QCoh(X)_Z$, i.e. we need to show that for $M \in \QCoh(X)$, if $j^* M = 0$ and $i_1^* M = 0$ then $M = 0$. By picking generators of the ideal cutting out $Z$ and looking at the induced map $X \to \mathbb A^n$ (as in the proof of \cref{rslt:completion-of-affine-is-colimit}) we can reduce to the case $X = \mathbb A^n$ and $Z = \{ 0 \}$. By considering the sequence $Z = \mathbb A^0 \subseteq \mathbb A^1 \subseteq \dots \subseteq \mathbb A^n$, where the inclusions are given by forcing coordinates to be $0$, we can reduce the claim to the case $X = \mathbb A^1$ (each of the above inclusions is a base-change of this case). Thus we are left with showing the following: Given $M \in \D(\Z[x])$ such that $M \tensor_{\Z[x]} \Z[x,x^{-1}] = 0$ and $M \tensor_{\Z[x]} \Z = 0$ then $M = 0$. But $M \tensor_{\Z[x]} \Z = \cofib(M \xto{x} M)$, hence multiplication by $x$ is an isomorphism on $M$. On the other hand $M \tensor_{\Z[x]} \Z[x,x^{-1}]$ is the colimit of multiplication by $x$ on $M$, so $M = \varinjlim_x M = 0$ and we conclude.
\end{proof}

In \cref{rslt:completion-of-affine-is-colimit} we saw that the completion of an affine scheme at a finitely generated ideal can be written as a colimit of closed immersions. This description involved the choice of generators and hence does not globalize easily. In the following we will provide a different description of the completion which does globalize and admits an explicit description in terms of the cotangent complex. Similar results are claimed in \cite[\S A]{Halpern-Leistner} using model theory and in \cite[\S9]{gaitsgory-rozenblyum-vol2} using an advanced 2-categorical machine (where they aim to prove a more general statement). In the following we provide a quick direct argument using filtered rings, which roughly implements the argument in \cite{Halpern-Leistner} (more precisely, a fix of that argument, see \cref{rmk:Halpern-Leistner-argument-error}) without referring to model categories. We learnt many of the following definitions from Gardner--Hekking and they also appear in \cite[\S3]{GardnerHekking-Ideals}.

\begin{defn}
\begin{defenum}
    \item We denote
    \begin{align*}
        \Pair := \Pair(\D(\Z)), \qquad \FilgeRing := \FilgeRing(\D(\Z)),
    \end{align*}
    where $\D(\Z)$ is the standard derived algebraic context and $\Ring(-)$ and $\Fil^{\ge0}(-)$ are defined in \cref{sec:filtered-rings}. Objects in $\Pair$ are called \emph{pairs} and objects in $\FilgeRing$ are called \emph{(non-negatively) filtered rings}.

    \item Given a ring $A$, an \emph{ideal} of $A$ is a pair of the form $(A, I)$. We denote by
    \begin{align*}
        \Ideals_A := \Pair \times_{\Ring} \{ A \}
    \end{align*}
    the category of ideals in $A$.
\end{defenum}
\end{defn}

Recall from \cref{sec:filtered-rings} that an object in $\Pair$ is a pair $(A, I)$ consisting of a ring $A$ and an ideal $I$ in $A$ (we take this as the \emph{definition} of the notion of an ideal in a ring). More explicitly, $\Pair$ is the category of (animated) rings in the symmetric monoidal category $\Fun([1]^\op, \D(\Z))$ of maps $I \to A$, where the tensor product is given by
\begin{align*}
    (I \to A) \tensor (I' \to A') = (I \tensor A' \sqcup_{I \tensor I'} A \tensor I' \to A \tensor A').
\end{align*}
As shown in \cref{rslt:equivalence-of-pairs-and-ring-maps}, to every pair $(A, I)$ there is an associated ring map $A \to A/I$ and this construction provides an equivalence between the category of pairs and the category of closed immersions of rings.

The objects of $\FilgeRing$ are sequences $\dots I^2 \to I^1 \to I^0 = A$ together with a coherent multiplicative structure that induces compatible maps $I^n \tensor I^m \to I^{n+m}$ for all $n, m \ge 0$. We often denote the objects by $(A, I^\bullet)$. As shown in \cref{rslt:quotients-from-filtered-ring}, to every filtered ring $(A, I^\bullet)$ we get an associated sequence $A \to \dots \to A/I^n \to \dots \to A/I^2 \to A/I$ of rings.

With the above definitions at hand we can now come to the definition of the Rees algebra of a pair $(A, I)$, which generalizes the classical Rees construction.

\begin{defn}
Restriction along $[1]^\op \injto \Z_{\ge0}^\op$ provides a lax morphism of derived algebraic contexts $\Fil^{\ge0}(\D(\Z)) \to \Fil^{[0,1]}(\D(\Z))$, which induces a functor $\FilgeRing \to \Pair$. This functor preserves limits as that can be checked on underlying objects (see \cref{def:Sym-and-forgetful-functor-on-DAG}). We denote by
\begin{align*}
    \Rees\colon \Pair \to \FilgeRing, \qquad (A, I) \mapsto \Rees(A, I)
\end{align*}
the left adjoint of this functor and call the filtered ring $\Rees(A, I)$ the \emph{Rees algebra} associated to the pair $(A, I)$.
\end{defn}

In order to study the Rees algebra we need to introduce one more piece of notation, namely the symmetric powers associated to a connective module.

\begin{defn}
\begin{defenum}
    \item Let $A$ be a ring and $M \in \D^{\le0}(A)$ a connective $A$-module. We denote by $\Sym_B^n(M)$ the $n$-th (derived) symmetric power of a module over a ring, as defined in \cite[Construction~25.2.2.1]{lurie-SAG} by deriving the usual symmetric power (cf. \cref{def:symmetric-power}).

    \item The assignment $M \mapsto \Sym_A^n(M)$ commutes with base-change in $A$ (this easily reduces to the case of polynomial rings over $\Z$, see \cite[Proposition~25.2.3.1]{lurie-SAG}), hence for every stack $X$ over some ring $\Lambda$ and every $M \in \QCoh^{\le0}(X)$ and $n \ge 0$ we get an associated quasicoherent sheaf
    \begin{align*}
        \Sym_X^n(M) \in \QCoh^{\le0}(X),
    \end{align*}
    the \emph{$n$-th symmetric power of $M$}.
\end{defenum}
\end{defn}

\begin{rmk} \label{rlst:Sym-preserves-pseudo-coherence}
The functor $\Sym^n_X$ preserves many properties of quasicoherent sheaves, see \cite[\S 25.2]{lurie-SAG}. For example, it preserves pseudo-coherence by \cite[Corollary~25.2.5.2]{lurie-SAG} (this reference proves it for rings, but it immediately follows for stacks by descent).
\end{rmk}

The next result is taken from \cite[Proposition~13.3]{scholze-complex-geometry} (which in turn is based on \cite[\S3.4]{Mao-derived-crystalline-cohom}) and provides a tight relation between the Rees algebra and the cotangent complex.

\begin{lem} \label{rslt:properties-of-Rees-algebra}
Let $(A, I)$ be a pair.
\begin{lemenum}
    \item There are natural isomorphisms $\Rees(A, I)^0 = A$ and $\Rees(A, I)^1 = I$.
    \item \label{rslt:Rees-algebra-computed-using-cotangent-complex} Denote $\Rees(A, I) = (A, I^\bullet)$ and $B = A/I$. Then there are natural isomorphisms
    \begin{align*}
        I^n/I^{n+1} = \Sym_B^n(L_{B/A}[-1]).
    \end{align*}
    for all $n \ge 0$. Here $I^n/I^{n+1} := \cofib(I^{n+1} \to I^n)$.
\end{lemenum}
\end{lem}
\begin{proof}
Note that all constructions in question commute with sifted colimits in the pair $(A, I)$ (for the cotangent complex use that $L_{B/A} = \cofib(B \tensor L_{A/\Z} \to L_{B/\Z})$ and \cref{rslt:cotangent-complex-for-rings-preserves-sifted-colimits}). Therefore all the claimed identities can be proved by providing functorial equivalences in the case that $(A, I)$ is a compact projective generator of $\Pair$, i.e. of the form $(A, I) = (\Z[x_1, \dots, x_n,y_1, \dots, y_m], (y_1, \dots, y_m))$. This is the symmetric algebra of the inclusion map $\Z^m \to \Z^{n+m}$ in $\Fil^{[0,1]}(\D(\Z))$. Since taking symmetric algebras is left adjoint to the forgetful functor, we see that $\Rees \comp \Sym$ is left adjoint to the functor $\FilgeRing \to \Fil^{[0,1]}(\D(\Z))$ sending $(A, I^\bullet) \mapsto (I \to A)$. But the left adjoint of the latter functor can also be computed as $(I \to A) \mapsto \Sym(\dots \to 0 \to 0 \to I \to A)$. Therefore
\begin{align*}
    \Rees(A, I) &= \Rees(\Z[x_\bullet, y_\bullet], (y_\bullet))\\
    &= \left[\dots \to (y_\bullet)^3 \to (y_\bullet)^2 \to (y_\bullet) \to \Z[x_\bullet,y_\bullet]\right].
\end{align*}
In other words, $\Rees(A, I)$ is the ring $A$ equipped with the $I$-adic filtration. We immediately deduce (i). Also $B = \Z[x_\bullet]$ and we see that $I^n/I^{n+1} = \Sym_B^n(I/I^2)$. To prove (ii) it remains to show that $I/I^2 = L_{B/A}[-1]$ in a functorial way. But this follows immediately from the short exact sequence $0 \to L_{B/A}[-1] \to B \tensor_A L_{A/\Z} \to L_{B/\Z} \to 0$ induced by the fiber sequence in \cref{rslt:fiber-sequence-for-cotangent-complex} together with \cref{rslt:cotangent-complex-of-polynomial-algebra}.
\end{proof}

\begin{defn}
Given a pair $(A, I)$ we denote $I^n := \Rees(A, I)^n$ and call it the \emph{$n$-th power of $I$}. By \cref{rslt:quotients-from-filtered-ring} $I^n$ is again an ideal in $A$ and hence there is an induced ring $A/I^n$.
\end{defn}

The next result records some useful properties of ideals and their powers. These results are not needed in the sequel, but we believe the other statements to be interesting on its own.

\begin{lem}
\begin{lemenum}
    \item \label{rslt:powers-of-ideal-become-connective} Let $I$ be an ideal of a ring $A$ such that $H^0 I = 0$. Then $H^{-k}(I^n) = 0$ for all $n \ge 0$ and $k < n$.
    
    \item Let $B$ be a ring, $M \in \D^{\le0}(B)$ a connective $B$-module and $A = \Sym_B(M)$ with the canonical map $A \to B$ induced by $M \to 0$. If $I$ is the ideal of this map then there are natural isomorphisms $I^n = \bigoplus_{k \ge n} \Sym_B^k(M)$ for all $n \ge 0$.
\end{lemenum}
\end{lem}
\begin{proof}
In (ii) both sides of the claimed isomorphism commute with sifted colimits of the pair $(B, M)$ and can thus be checked in the case $B = \Z[x_1, \dots, x_k]$ and $M = B^m$ for some $k, m \ge 0$. In this case the claim follows easily from the explicit description of $I^n$ as in the proof of \cref{rslt:properties-of-Rees-algebra}.

To prove (i) let $\Pair' \subseteq \Pair$ denote the full subcategory spanned by those pairs $(A, I)$ such that $H^0 I = 0$. Then $\Pair'$ is stable under sifted colimits in $\Pair$ (as they can be computed on underlying modules) and it is compact projectively generated by the elements $s_{X,Y} = \Sym(X[1] \to X[1] \oplus Y)$ for finite free $\Z$-modules $X$ and $Y$. In particular every element of $\Pair'$ is a sifted colimit of $s_{X,Y}$'s, so we can in (i) assume that $(A, I) = s_{X,Y}$ for some $X$ and $Y$. Since $s_{X,Y} = s_{X,0} \tensor s_{0,Y}$ we can further reduce to the cases that either $X = 0$ or $Y = 0$. In the first case we have $I = 0$ and thus $I^n = 0$ for all $n$. In the second case we have $A = \Sym_{\Z}(X[1])$ and $B = \Z$, hence by (ii) we deduce $I^n = \bigoplus_{k\ge n} \Sym_{\Z}^k(X[1])$. We conclude by the connectivity of symmetric powers (see \cite[Proposition~25.2.4.1]{lurie-SAG}).
\end{proof}

We now extend \cref{rslt:properties-of-Rees-algebra} to stacks in order to obtain the desired construction of infinitesimal neighborhoods and with it a description of completions. A similar result is claimed in \cite[\S A]{Halpern-Leistner}, but the proof is flawed (see \cref{rmk:Halpern-Leistner-argument-error} below).

\begin{prop} \label{rslt:infinitesimal-neighborhoods}
Fix a ring $\Lambda$ and let $i\colon Z \to X$ be a closed immersion of stacks over $\Lambda$. There is a natural sequence
\begin{align*}
    \emptyset = Z_0 \to Z = Z_1 \to Z_2 \to \dots \to Z_n \to \dots \to X_{\hat Z} \subseteq X
\end{align*}
of closed immersions satisfying the following:
\begin{propenum}
    \item The formation of $Z_n$ is compatible with base-change in $X$. If $X = \Spec A$ is affine and $Z = \Spec A/I$, then $Z_n = \Spec A/I^n$.
    
    \item \label{rslt:fiber-of-infinitesimal-nbhds} Denoting $i_n\colon Z_n \to X$ the induced map then for all $n \ge 0$ there is a natural isomorphism
    \begin{align*}
        \fib(i_{n+1,*} \calO_{Z_{n+1}} \to i_{n,*} \calO_{Z_n}) = i_* \Sym_Z^n(L_{Z/X}[-1]).
    \end{align*}
    in $\QCoh(X)$.

    \item \label{rslt:completion-is-colimit-of-inf-nbhds-on-truncated} Suppose that $i$ is finitary. Then for all truncated rings $C \in \Alg_\Lambda$ the map
    \begin{align*}
        \varinjlim_n Z_n(C) \isoto X_{\hat Z}(C)
    \end{align*}
    is an isomorphism. Moreover, the functor
    \begin{align*}
        i_\bullet^*\colon \QCoh^-(X_{\hat Z}) \isoto \varprojlim_n \QCoh^-(Z_n)
    \end{align*}
    is an equivalence.
\end{propenum}
\end{prop}
\begin{proof}
Without loss of generality we can assume $\Lambda = \Z$. For every stack $X$ we denote by $\Stk_X^{\rm ci} \subseteq \Stk_{/X}$ the full subcategory spanned by the closed immersions. By \cref{rslt:affineness-and-closed-immersion-are-local-on-target} the assignment $X \mapsto \Stk_X^{\rm ci}$ satisfies descent, i.e. defines a sheaf of categories on $\Stk$. The same is then true for the assignment $X \mapsto \Fun(\Z_{\ge0}, \Stk_X^{\rm ci})$ because $\Fun(\Z_{\ge0}, -)$ preserves limits. We wish to construct a natural functor $\Stk_X^{\rm ci} \to \Fun(\Z_{\ge0}, \Stk_X^{\rm ci})$ mapping a closed immersion $Z \to X$ to the sequence $Z_0 \to Z_1 \to \dots \to X$. By sheafiness this functor can be constructed on the defining site (see \cite[Lemma~A.4.18]{heyer-mann-6ff}), i.e. on $\Ring$. Observe that the cocartesian unstraightening of the functor $\Ring \to \Cat$, $A \mapsto \Stk_A^{\rm ci}$ is the full subcategory of $\Fun([1], \Ring)$ spanned by the closed immersions, which by \cref{rslt:equivalence-of-pairs-and-ring-maps} is equivalent to $\Pair$. Moreover, the unstraightening of the functor $\Ring \to \Cat$, $A \mapsto \Fun(\Z_{\ge0}, \Stk_A^{\rm ci})$ is given by $\Fun(\Z_{\ge0}^{\triangleright,\op}, \Ring)$ (indeed, there is a natural map from the unstraightening to this category induced by evaluation and this map is clearly fiberwise an equivalence). Altogether the construction of the desired functor is now reduced to constructing the map
\begin{align*}
    \Pair \to \Fun(\Z_{\ge0}^{\triangleright,\op}, \Ring), \qquad (A, I) \mapsto [A \to \dots \to A/I^n \to \dots \to A/I^2 \to A/I \to 0]
\end{align*}
of cocartesian fibrations over $\Ring$. The above functor is given by the composition of $\Rees$ and \cref{rslt:quotients-from-filtered-ring}. It remains to show that this functor preserves cocartesian edges, which is the following statement: Given a ring map $A \to B$ and an ideal $I \in A$ with induced ideal $IB := I \tensor_A B$ then for all $n \ge 0$ we have $A/I^n \tensor_A B = B/(IB)^n$. For $n = 1$ this is true by definition of $IB$. For $n \ge 1$ we argue by induction using \cref{rslt:Rees-algebra-computed-using-cotangent-complex} and the fact that both $L_{(A/I)/A}$ and $\Sym_A^n$ commute with base-change along $A \to B$. This finishes the construction of $(Z_n)_n$.

Part (i) holds by construction. To prove (ii) we note that both sides of the claimed identity are functors $\Stk_X^{\rm ci} \to \QCoh(X)$ which assemble to morphisms of sheaves on $\Stk$. Hence the identities can be constructed on the site $\Ring$, where they were established in \cref{rslt:Rees-algebra-computed-using-cotangent-complex}.

We now prove (iii), so suppose that $i$ is finitary and let $C$ be a truncated ring. By writing $X$ as a colimit of affine schemes we can reduce to the case that $X = \Spec A$ is affine. Let $(A, I)$ be the pair associated to $Z$, i.e. $Z = \Spec A/I$. We need to show that the map
\begin{align}
    \varinjlim_n \Hom(A/I^n, C) \isoto X_{\hat Z}(C) \label{eq:infinitesimal-nbhd-colim}
\end{align}
is an isomorphism. Note that by \cref{rslt:powers-of-ideal-become-connective} this is clear if $H^0 I = 0$: In this case, for $n$ large enough we have $\Hom(A/I^n, C) = \Hom(A, C)$ because then $\tau^{<-n}(A/I^n) = \tau^{<-n} A$. This observation motivates the following claim about an arbitrary pair $(A, I)$:
\begin{itemize}
    \item[($*$)] Let $I \to J$ be a map of ideals of $A$ such that $H^0 I \isoto H^0 J$ is an isomorphism. Then the induced map $\varinjlim_n \Hom(A/J^n, C) \isoto \varinjlim_n \Hom(A/I^n, C)$ is an isomorphism.
\end{itemize}
We first observe that claim ($*$) is compatible with sifted colimits as follows. Fix an integer $k \ge 0$ such that $C$ is $k$-truncated, i.e. $\pi_i C = 0$ for all $i > k$. Then $\Hom(B, C)$ is $k$-truncated for every ring $B$, hence both sides of the claimed isomorphism in ($*$) are $k$-truncated. Now observe that filtered colimits of $k$-truncated anima commute with totalizations, because the category of $k$-truncated anima is a $(k+1)$-category and a functor between $(k+1)$-categories preserves totalizations as soon as it preserves finite limits (since totalizations in a $(k+1)$-category are the same as $(k+1)$-truncated totalizations, see the opposite of \cite[Lemma~1.3.3.10(2)]{lurie-higher-algebra}); but filtered colimits certainly preserve finite limits. Note also that the assignment $(A, I) \mapsto A/I^n$ preserves sifted colimits (this can be checked on underlying modules where it follows easily from the definitions). Altogether we see that if the diagram $A \to A/I \to A/J$ is obtained as a geometric realization $\varinjlim_{m\in\Delta} [A_m \to A_m/I_m \to A_m/J_m]$ then the claim ($*$) for $A \to A/I \to A/J$ reduces to the claim ($*$) for each $A_m \to A_m/I_m \to A_m/J_m$, because then
\begin{align*}
    &\varinjlim_n \Hom(A/J^n, C) = \varprojlim_{m\in\Delta} \varinjlim_n \Hom(A_m/J_m^n, C) = \varprojlim_{m\in\Delta} \varinjlim_n \Hom(A_m/I_m^n, C) \\&\qquad= \varinjlim_n \Hom(A/I^n, C).
\end{align*}
We will use this observation below.

Let $\cat C$ be the functor category of diagrams $A \to A' \to A''$ of rings and let $\cat C' \subseteq \cat C$ be the full subcategory spanned by those diagrams such that $A \to A'$ is $0$-connected and $A' \to A''$ is $1$-connected, i.e. $\pi_0 A \surjto \pi_0 A'$ and $\pi_1 A' \to \pi_1 A''$ are surjective and $\pi_0 A' \isoto \pi_0 A''$ is an isomorphism. By looking at long exact sequences we see that $\cat C'$ is exactly the category of maps $A \to A/I \to A/J$ in ($*$). Similarly to the proof of \cref{rslt:equivalence-of-pairs-and-ring-maps} we see that $\cat C'$ is compact projectively generated by the elements
\begin{align*}
    \Sym_{\Z}(X) \to \Sym_{\Z}(X) \to \Sym_{\Z}(X), \qquad \Sym_{\Z}(X) \to \Z \to \Z, \qquad \Z \to \Sym_{\Z}(X[1]) \to \Z
\end{align*}
for finite free abelian groups $X$. By \cite[Lemma~5.5.8.13]{lurie-higher-topos-theory} we can thus write $[A \to A/I \to A/J] = \varinjlim_{m\in\Delta} [A_m \to A_m/I_m \to A_m/J_m]$ such that each $A_m \to A_m/I_m \to A_m/J_m$ is a (possibly infinite) coproduct of the three generators above. Thus ($*$) reduces to the case that $A \to A/I \to A/J$ is such a coproduct. We split this coproduct into a coproduct of three coproducts, each consisting only of a single type of generator. We can then pull this finite coproduct out of both sides of the claimed isomorphism in ($*$) to reduce to the case that $A \to A/I \to A/J$ is a coproduct of a single type of generator above. Note that for the first two types of generators the map $A/I \to A/J$ is an isomorphism, hence ($*$) is clear. We are thus reduced to the third type of generator, i.e. we have reduced ($*$) to the following claim: Let $S$ be any set and let $[A \to A/I \to A/J] = [\Z \to \Sym_{\Z}(\Z[S][1]) \to \Z]$; then the natural map
\begin{align*}
    \varinjlim_n \Hom(A/I^n, C) \isoto \varinjlim_n \Hom(A/J^n, C) = \Hom(\Z, C) = *
\end{align*}
is an isomorphism. Picking $k \ge 0$ large enough so that $C$ is $k$-truncated, the above identity follows immediately if we can show that for each $n \ge 0$ and $m \gg n$ the map $\tau_{\le k}(A/I^m) \to \tau_{\le k}(A/I^n)$ factors over $\Z$. By \cref{rslt:composition-of-phantom-ring-maps-vanishes} this is implied if we can show that for all $i > 0$ and $n \ge 0$ the map $\pi_i(A/I^{n+1}) \to \pi_i(A/I^n)$ is zero. By writing $S$ as a filtered colimit of finite sets we can reduce to the case that $S$ is finite. Note that $A/I = \Z \tensor_{\Z[x_\bullet]} \Z$, where $\Z[x_\bullet]$ is the polynomial ring with variables parametrized by $S$ (because $\Sym$ preserves colimits). Thus $A \to A/I$ is the base-change of the map $\Z[x_\bullet] \to \Z$ along itself. By the compatibility of the formation of $A/I^n$ with base-change discussed in (i), we deduce
\begin{align*}
    A/I^n = \Z[x_\bullet]/(x_\bullet)^n \tensor_{\Z[x_\bullet]} \Z.
\end{align*}
Using the Koszul resolution of $\Z$ as a $\Z[x_\bullet]$-module we see that $A/I^n$ admits an explicit description as a complex whose terms are direct sums of copies of $\Z[x_\bullet]/(x_\bullet)^n$ and whose transition maps are matrices with entries being monomials of degree $1$. Moreover, the same complex with $\Z[x_\bullet]$ in place of $\Z[x_\bullet]/(x_\bullet)^n$ is exact everywhere except at $0$. This implies that every element in $\pi_i(A/I^n) = H^{-i}(A/I^n)$ comes from an element in a direct sum of $(x_\bullet)^{n-1}/(x_\bullet)^n \subseteq \Z[x_\bullet]/(x_\bullet)^n$. In particular this maps gets sent to zero under the map $A/I^n \to A/I^{n-1}$, as desired. This finishes the proof of ($*$).

With claim ($*$) at hand, we can now continue with the proof of (iii), so we come back to the pair $(A, I)$ associated to the finitary closed immersion $Z \injto X$ and we wish to prove the isomorphism in \cref{eq:infinitesimal-nbhd-colim}. By definition of finitary closed immersions, the image of $H^0 I$ in $\pi_0 A$ is finitely generated. Pick generators $f_1, \dots, f_k$ of this image and denote $A' := A \tensor_{\Z[x_1, \dots, x_k]} \Z$, where $\Z[x_1, \dots, x_k] \to A$ is the induced map. Then the induced map $A' \to A/I$ is an isomorphism on $\pi_0$. Let $A'' = A' \tensor_{\Z} \Sym_{\Z}(\pi_1(A/I)[1])$. Then the natural map $A'' \to A/I$ is still an isomorphism on $\pi_0$ and additionally surjective on $\pi_1$, i.e. $1$-connected, so by ($*$) we can replace the closed immersion $A \to A/I$ by $A \to A''$. Note that the projection $A'' \to A'$ (induced by the zero map $\pi_1(A/I)[1] \to 0$) is also $1$-connected, so by an additional application of ($*$) we can reduce to the map $A \to A'$. This map is a base-change of the map $\Z[x_\bullet] \to \Z$, so by (i) and \cref{rslt:completion-commutes-with-base-change} we can finally reduce to the case $(A, I) = (\Z[x_\bullet], (x_\bullet))$ (with finitely many polynomial variables). We now observe that the natural map $(\Z[x_\bullet]/(x_\bullet^n))_n \to (\Z[x_\bullet]/(x_\bullet)^n)_n$ is an isomorphism of sequences of rings because it splits: By the pigeon hole principle, for every $n \ge 0$ we have $(x_\bullet)^{kn} \subseteq (x_\bullet^n)$ (where $k$ is the number of polynomial variables) and hence there is a natural backward map $\Z[x_\bullet]/(x_\bullet)^{kn} \to \Z[x_\bullet]/(x_\bullet^n)$. This implies
\begin{align*}
    \varinjlim_n \Hom(A/I^n, C) = \varinjlim_n \Hom(\Z[x_\bullet]/(x_\bullet)^n, C) = \varinjlim_n \Hom(\Z[X_\bullet]/(x_\bullet^n), C) = X_{\hat Z}(C),
\end{align*}
where the last identity follows from (the proof of) \cref{rslt:completion-of-affine-is-colimit}. This finishes the proof of \cref{eq:infinitesimal-nbhd-colim}.

To finish the proof of (iii), it remains to prove the second part of the claim, i.e. the identity on $\QCoh^-$. But this follows immediately from the first part of (iii) and \cref{rslt:QCoh-minus-only-depends-on-truncated-rings}.
\end{proof}

\begin{rmks}
\begin{rmksenum}
    \item \label{rmk:Halpern-Leistner-argument-error} The results in \cref{rslt:infinitesimal-neighborhoods} are also claimed in \cite[\S A]{Halpern-Leistner} (under stronger hypotheses). However, the argument in the reference appears to be wrong: It relies crucially on the identity in line 4 of page 87, which fails in many cases, e.g. $A = \Z[x]$ and $I = (x)$. Our proof provides a different argument for the claims and also generalizes them slightly.

    \item A similar result to \cref{rslt:infinitesimal-neighborhoods} is proved in \cite[Theorem~5.4.7]{GardnerHekking-Ideals}, but under much stronger hypotheses: In the reference, the authors assume that the closed immersion $i\colon Z \to X$ is a complete intersection. In this case \cref{rslt:completion-is-colimit-of-inf-nbhds-on-truncated} holds without truncatedness assumption on $C$ by a similar argument as in \cref{rslt:completion-of-affine-is-colimit}.
\end{rmksenum}
\end{rmks}

\begin{defn} \label{def:infinitesimal-neighbourhoods}
Let $Z \to X$ be a closed immersion of stacks over some ring $\Lambda$. We call the stacks $Z_n$ constructed in \cref{rslt:infinitesimal-neighborhoods} the \emph{infinitesimal neighbourhoods of $Z$ in $X$}.
\end{defn}

Using \cref{rslt:infinitesimal-neighborhoods} we can now state the following version of the excision sequence for quasicoherent sheaves. We refer the reader to \cref{rslt:ICoh-excision-all} for a similar result with ind-coherent sheaves.

\begin{cor} \label{rslt:excision-for-QCoh}
Fix a ring $\Lambda$ and let $i\colon Z \to X$ be a closed immersion of stacks over $\Lambda$ such that the open complement $j\colon U \injto X$ is a quasicompact map. Let $\hat i\colon X_{\hat Z} \injto X$ denote the completion of $X$ at $Z$.
\begin{corenum}
    \item The functor $\hat i^*$ admits a fully faithful left adjoint $\hat i_\natural\colon \QCoh(X_{\hat Z}) \injto \QCoh(X)$ and the functor $j_*$ admits a right adjoint $j^?\colon \QCoh(X) \to \QCoh(U)$. Moreover, for every $M \in \QCoh(X)$ we have the following fiber sequences in $\QCoh(X)$:
    \begin{align*}
        \hat i_\natural \hat i^* M \to M \to j_* j^* M, \qquad j_* j^? M \to M \to \hat i_* \hat i^* M.
    \end{align*}

    \item Suppose that $X$ is geometric, $i$ is finitary and $M \in \QCoh^-(X)$. Let $i_n\colon Z_n \to X$ be the $n$-th infinitesimal neighbourhood of $i$ and $M_n := i_{n*} i_n^* M$. Then
    \begin{align*}
        \hat i_* \hat i^* M = \varprojlim_n M_n, \qquad \fib(M_{n+1} \to M_n) = M \tensor i_* (\Sym_Z^n(L_{Z/X}[-1]))
    \end{align*}
    for all $n$.
\end{corenum}
\end{cor}
\begin{proof}
Part (i) can be deduced from \cref{rslt:hat-i-pullback-induces-equiv} as follows. Namely, by that reference we know that $\hat i^*\colon \QCoh(X) \to \QCoh(X_{\hat Z})$ restricts an equivalence $\QCoh(X)_Z \isoto \QCoh(X_{\hat Z})$, so we can define $\hat i_\natural$ as the composition
\begin{align*}
    \hat i_\natural\colon \QCoh(X_{\hat Z}) \xto{(\hat i^*)^{-1}} \QCoh(X)_Z \injto \QCoh(X).
\end{align*}
Clearly $\hat i^* \hat i_\natural = \id$, which provides us with the unit of an adjunction between $\hat i^*$ and $\hat i_\natural$. To show that this is indeed an adjunction, we provide a different formula for $\hat i^*$. Consider the functor $F\colon \QCoh(X) \to \QCoh(X)_Z$ given by $F(M) = \fib(M \to j_* j^* M)$. Since $\hat i_\natural \hat i^*$ is clearly the identity on $\QCoh(X)_Z$, we have $\hat i_\natural \hat i^* F = F$. On the other hand we compute
\begin{align*}
    \hat i_\natural \hat i^* F = \fib(\hat i_\natural \hat i^* \to \hat i_\natural \hat i^* j_* j^*) = \hat i_\natural \hat i^*,
\end{align*}
where in the last step we used that $\hat i^* j_* = 0$ by \cref{rslt:qcqs-base-change-for-QCoh}. Altogether we see that $\hat i_\natural \hat i^* M = \fib(M \to j_* j^* M)$. Thus, for all $M \in \QCoh(X_{\hat Z})$ and $N \in \QCoh(X)$ we have
\begin{align*}
    &\Hom(M, \hat i^* N) = \Hom(\hat i_\natural M, \hat i_\natural \hat i^* N) = \Hom(\hat i_\natural M, \fib(N \to j_* j^* N)) =\\&\qquad= \fib(\Hom(\hat i_\natural M, N) \to \Hom(\hat i_\natural M, j_* j^* N)).
\end{align*}
On the other hand, $\Hom(\hat i_\natural M, j_* j^* N) = \Hom(j^* \hat i_\natural M, j^* N) = 0$, because $j^* \hat i_\natural = 0$ by definition of $\QCoh(X)_Z$. We deduce that $\hat i_\natural$ is left adjoint to $\hat i^*$ and we have also already established the first fiber sequence in (i). The second fiber sequence can be obtained in a similar way, or alternatively be deduced from the first one: For every $M \in \QCoh(X)$ we have the fiber sequence
\begin{align*}
    \intHom(j_* j^* 1, M) \to \intHom(1, M) \to \intHom(\hat i_\natural \hat i^* 1, M),
\end{align*}
which one easily evaluates to the desired fiber sequence.

We now prove (ii), so assume that $X$ is geometric and $i$ is finitary. Given $M \in \QCoh^-(X)$, denote $M_n$ and $i_n\colon Z_n \to X$ as in the claim. By \cref{rslt:infinitesimal-neighborhoods} and the projection formula for $i_n$ (see \cref{rslt:qcqs-projection-formula--and-cocont-for-QCoh}) we have
\begin{align*}
    &\fib(M_{n+1} \to M_n) = \fib(i_{n+1,*} i_{n+1}^* M \to i_{n*} i_n^* M) = M \tensor \fib(i_{n+1,*} 1 \to i_{n*} 1) =\\&\qquad= M \tensor i_* \Sym_Z^n(L_{Z/X}[-1]),
\end{align*}
proving the second part of (ii). To prove the first part, we have to show that for all $N \in \QCoh(X)$ the map $\Hom(N, M) \to \varprojlim_n \Hom(N, M_n)$ is an isomorphism. Since $X$ is geometric and hence $\QCoh(X)$ has a (left and) right complete t-structure (see \cref{rslt:t-structure-on-QCoh}), we can write $N = \varinjlim_m \tau^{\le m} N$. By pulling out this colimit we reduce to the case that $N \in \QCoh^-(X)$.
\end{proof}

\begin{defn} \label{def:local-cohomology-sheaf}
In the setting of \cref{rslt:excision-for-QCoh}, given any $M \in \QCoh(X)$ we call
\begin{align*}
    M_Z := \hat i_\natural \hat i^* M \in \QCoh(X)
\end{align*}
the \emph{local cohomology sheaf} of $M$ at $Z$.
\end{defn}

\begin{rmk} \label{rslt:computation-of-local-cohomology}
In the setting of \cref{rslt:excision-for-QCoh}, suppose that the closed immersion $Z \to X$ comes via base-change from the closed immersion $\Spec A/I \to \Spec A$ for some $\Lambda$-algebra $A$ and some finitely generated ideal $I = (f_1, \dots, f_n)$ of $A$. Then for every $M \in \QCoh(X)$ we have
\begin{align*}
    M_Z = M \tensor_A \left[A \to \bigoplus_i A[\frac1{f_i}] \to \bigoplus_{i<j} A[\frac1{f_if_j}] \to \dots \right],
\end{align*}
where by $M \tensor_A N$ we mean $M \tensor g^* N$ for $N \in \D(A)$ and $g\colon X \to \Spec A$ the implicit projection. Indeed, the above formula follows immediately from the excision sequence $M_Z \to M \to j_* j^* M$ and the fact that $j_* j^* M$ has the above explicit description as a complex (this can be checked via base-change from the universal case $A = \Z[x_1, \dots, x_n]$, $I = (x_1, \dots, x_n)$).
\end{rmk}

\subsection{QCA stacks and descendability} \label{sec:alggeo.QCA-stacks}

In this subsection we single out a particularly nice class of stacks in characteristic $0$, namely the QCA stacks. They were introduced in \cite[Definition~1.1.8]{DrinfeldGaitsgory-finiteness} where they are shown to have excellent finiteness properties with regard to quasi-coherent sheaves. We extend these finiteness results by proving a certain descendability property, which will be crucial for extending the six functors on ind-coherent sheaves to QCA stacks. Similar results have been obtained in \cite[Theorem~7.1]{Hall-descendability} on non-derived algebraic stacks.

\begin{defn}
Fix a $\Q$-algebra $\Lambda$. A map $f\colon Y \to X$ of stacks over $\Lambda$ is called \emph{QCA} if it satisfies the following conditions:
\begin{enumerate}[(i)]
    \item $f$ is geometric and quasicompact.
    \item The diagonal of $f$ is affine.
\end{enumerate}
We say that a stack $X$ over $\Lambda$ is \emph{QCA} if the map $X \to \Spec \Lambda$ is QCA.
\end{defn}

\begin{rmk}
In the definition of QCA stacks we require $\Lambda$ to be a $\Q$-algebra, because only then we can prove the desired finiteness and descendability results. There should be a version of QCA stacks that works over arbitrary rings where one adds an additional hypothesis on the stabilizer groups of $X$: They need to have finite cohomological dimension over $\Lambda$. In characteristic $p$, this means that they are linearly reductive, cf. \cite[Theorem~1.2]{HallRydh-compact-generation-of-alg-reps}. We do not pursue this approach here, as it only adds unnecessary complications without much benefit to us.
\end{rmk}

\begin{rmk} \label{rslt:stability-of-QCA-maps}
One checks easily that QCA maps are stable under composition and base-change, and that a map between QCA stacks is automatically QCA (apply \cite[Lemma~2.1.5]{heyer-mann-6ff} to the class of QCA maps).
\end{rmk}

Let us now come to the promised finiteness results on QCA stacks. They rely on the following structure result, which is based on the non-derived version in \cite{HallRydh-compact-generation-of-alg-reps}.

\begin{prop} \label{rslt:strata-decomposition-of-QCA-stack}
Let $\Lambda$ be a $\Q$-algebra and $X$ a QCA stack over $\Lambda$. Then there is a sequence of finitary closed immersions
\begin{align*}
    \emptyset = Z_0 \to Z_1 \to \dots \to Z_{r-1} \to Z_r = X
\end{align*}
such that for $i = 1, \dots, r$ the open substack $Z_i \setminus Z_{i-1} \subset Z_i$ is isomorphic to a stack of the form $Y_i/\GL_{n_i,\Lambda}$ for some $n_i \ge 0$, where $Y_i$ is a quasicompact open subset of a classical affine scheme over $\Lambda$.
\end{prop}
\begin{proof}
By \cref{rslt:X-cl-to-X-is-closed-immersion} the map $X^\cl \to X$ is a finitary closed immersion and in particular QCA, hence $X^\cl$ is still a QCA stack. We observe that the classical diagonal $X^\cl \to (X^\cl \times X^\cl)^\cl$ is affine. Namely, more generally if $Y' \to Y$ is an affine map of geometric stacks then $Y'^\cl \to Y^\cl$ is still affine: Pick a smooth cover $U \surjto Y$ by a disjoint union of affine schemes $U$, so that $U' := U \times_Y Y' \to U$ is is affine; then by \cref{rslt:flat-geometric-maps-stable-under-cl} we have $U'^\cl = U^\cl \times_{Y^\cl} Y'^\cl$ (as the right-hand side is classical) and clearly $U'^\cl \to U^\cl$ is still affine, so we conclude because affineness of $Y'^\cl \to Y^\cl$ can be checked after passing to the cover $U^\cl \to Y^\cl$ by \cref{rslt:affineness-and-closed-immersion-are-local-on-target}. Altogether we see that $X^\cl$ is a classical quasicompact algebraic stack over $\pi_0 \Lambda$ with affine diagonal (in the classical sense). In particular $X^\cl$ is quasiseparated.

Since $X^\cl$ has affine diagonal, it in particular has affine stabilizers. We now apply \cite[Proposition~2.6]{HallRydh-compact-generation-of-alg-reps} to $X^\cl$ in order to obtain a sequence $\emptyset = Z_0 \to Z_1 \to \dots \to Z_{r-1} \to X^\cl$ of classical finitely presented closed immersions of classical algebraic stacks. Note that these closed immersions are finitary closed immersions in the derived sense of this paper: By \cref{rslt:affineness-and-closed-immersion-are-local-on-target} this can be checked after passing to a smooth cover of $X^\cl$ by a classical affine scheme, and by \cref{rslt:flat-geometric-maps-stable-under-cl} passage to that smooth cover is compatible with the classical setting (i.e. the derived fiber product is classical). By construction, for all $i = 1, \dots, r-1$ the classical stack $Z_i \setminus Z_{i-1}$ admits a quasi-affine map to $*/\GL_{n_i,\Z}$ for some $n_i \ge 0$. Since this stack maps to $\Spec \Lambda$, we also get a quasi-affine map to $*/\GL_{n_i,\Lambda}$, which is the same as saying that $Z_i \setminus Z_{i-1}$ is of the form $Y_i/\GL_{n_i,\Lambda}$ as in the claim. Finally, we set $Z_r = X$ and note that the map $Z_{r-1} \to Z_r$ is still a finitary closed immersion (by \cref{rslt:X-cl-to-X-is-closed-immersion}) and the open complement is empty.
\end{proof}

With \cref{rslt:strata-decomposition-of-QCA-stack} at hand, we can now prove many results about QCA stacks by reducing them to the case $Y/\GL_n$ for a scheme $Y$. An important tool is the notion of descendable algebras from \cite[Definition~3.18]{akhil-galois-group-of-stable-homotopy}, so let us briefly recall it here:

\begin{defn} \label{def:descendable-algebra}
Let $\cat C$ be a stable symmetric monoidal category and $A \in \CAlg(\cat C)$ a commutative algebra in $\cat C$.
\begin{defenum}
    \item We denote by $\langle A \rangle \subseteq \cat C$ the full subcategory generated under tensor products, finite (co)limits and retracts by $A$.  We say that $A$ is \emph{descendable} if the tensor unit $1 \in \cat C$ lies in $\langle A \rangle$.

    \item We denote $K_A := \fib(1 \to A)$. We say that $A$ is \emph{descendable of index $\le d$} if the natural map $K_A^{\tensor d} \to 1$ is zero.
\end{defenum}
\end{defn}

\begin{defn} \label{def:descendable-map}
Let $f\colon Y \to X$ be a map of stacks over some ring $\Lambda$. We say that $f$ is \emph{descendable} (\emph{of index $d$}) if the algebra $f_* 1 \in \QCoh(X)$ is descendable (of index $d$). We also denote $K_f := \fib(1 \to f_* 1)$.
\end{defn}

By \cite[Lemma~11.20]{bhatt-scholze-witt} an algebra is descendable if and only if it is descendable of some index, so there is no ambiguity in the above definitions. More intuitively, $A \in \cat C$ is descendable if one can obtain $1$ from the tensor powers $A^{\tensor n}$ using only finitely many fiber sequences and retracts. This is a strong finiteness property on $A$ and we will see it in action below (see \cref{rslt:QCA-pushforward-properties}).

\begin{exmpl} \label{ex:fppf-cover-of-rings-is-descendable}
We follow \cite[Example~6.6]{Hall-descendability} in order to provide the following fundamental example of descendability: If $f\colon \Spec B \to \Spec A$ is an fppf cover of affine schemes then $f$ is descendable of index $\le 2$. In fact, even more is true: For every connective $A$-module $M \in \D^{\le0}(A)$ and every $i \ge 0$ we have
\begin{align*}
    H^2 \Hom(K_f^{\tensor i}[i], M) = 0.
\end{align*}
This implies the claimed descendability because it shows that there is no non-zero map $K_f^{\tensor 2} \to A$. In order to prove the $\Hom$-vanishing, note first that we can without loss of generality assume that $M$ is concentrated in degree $0$ (then for general $M$ use the Postnikov limit $M = \varprojlim_n \tau^{\le -n} M$). Then $M$ lies in the image of the forgetful functor $\D(\pi_0 A) \to \D(A)$, hence by adjunction we can replace $f$ by its pullback along $A \to \pi_0 A$, i.e. we can assume that $A$ is a classical ring; the same is then true for $B$. We observe that the classical $A$-module $K_f[1] = \cofib(A \to B) = B/A$ is flat: This can be checked after $- \tensor_A B$ (cf. \cite[Lemma~00HJ]{stacks-project}) where it becomes a direct summand of $B \tensor_A B$. Thus also $K_f^{\tensor i}[i]$ is a flat classical $A$-module for all $i \ge 0$. By \cite[Theorem~3.2]{Lazard-flatness} we see that this module has projective dimension $\le 1$, as desired.
\end{exmpl}

\begin{lem}
Fix a ring $\Lambda$ and let $f\colon Y \to X$ and $g\colon Z \to Y$ be qcqs schematic maps of stacks over $\Lambda$.
\begin{lemenum}
    \item If $f$ and $g$ are descendable, then so is $fg$.
    \item If $f$ is descendable, then so is any base-change of $f$.
    \item \label{rslt:descendable-implies-M-fin-from-Cech-nerve} Suppose that $f$ is descendable. Then every $M \in \QCoh(X)$ can be obtained using finitely many retracts and fiber sequences from the objects $f_{n*} f_n^* M$, where $f_\bullet\colon Y_\bullet \to X$ is the \v{C}ech nerve of $f$.
\end{lemenum}
\end{lem}
\begin{proof}
For (i), suppose that $f$ and $g$ are descendable. Then $g_* 1$ is a descendable algebra in $\QCoh(Y)$ and since $f_*$ is lax symmetric monoidal and preserves small colimits (see \cref{rslt:qcqs-projection-formula--and-cocont-for-QCoh}), it follows from \cite[Lemma~11.20]{bhatt-scholze-witt} that the map $f_* 1 \to g_*f_* 1$ is descendable in $\QCoh(X)$ (in the sense of \cite[Definition~11.18]{bhatt-scholze-witt}). By descendability of $f$ we know that $1 \to f_* 1$ is also descendable, hence by \cite[Lemma~11.17]{bhatt-scholze-witt} also $1 \to g_*f_* 1$ is descendable, as desired. Part (ii) follows easily from the fact that the formation of $f_* 1$ is compatible with any base-change by \cref{rslt:qcqs-base-change-for-QCoh}. Part (iii) is clear for $M = 1$ and then follows for arbitrary $M$ by noting that $f_{n*} f_n^* M = (f_* 1)^{\tensor n-1} \tensor M$ by the projection formula and base change (see \cref{rslt:qcqs-base-change-for-QCoh} and \cref{rslt:qcqs-projection-formula--and-cocont-for-QCoh}).
\end{proof}

Let us come back to QCA stacks. Our main result is that they always admit a descendable cover by an affine scheme, which we prove now. A similar result has been obtained in \cite[Theorem~7.1]{Hall-descendability}, but we provide an independent proof, relying only on \cref{rslt:strata-decomposition-of-QCA-stack}.

\begin{prop} \label{rslt:descendable-cover-of-QCA-stack}
Let $\Lambda$ be a $\Q$-algebra and $X$ a QCA stack over $\Lambda$. Then there is an integer $d \ge 0$ such that every flat geometric cover $f\colon U \to X$ by an affine scheme $U$ is descendable of index $\le d$.
\end{prop}
\begin{proof}
We roughly follow the strategy in \cite[Proposition~8.3.2]{Gaitsgory-1-affineness} but combine it with ideas from \cite[Theorem~7.1]{Hall-descendability}. Let us call a sheaf $M \in \QCoh(X)$ \emph{special} if for every map $f\colon \Spec A \to X$ and every $N \in \D^{\le0}(A)$ we have $H^2 \Hom(f^* M, N) = 0$. We then denote by $d(X) \in \Z \cup \{\infty\}$ the smallest number $d$ with the following property: For every special $M \in \QCoh(X)$ and every $N \in \QCoh^{\le-d}(X)$ we have $H^0 \Hom(M, N) = 0$. We first make the following observation:
\begin{itemize}
    \item[($*$)] If $d(X) < \infty$ then every flat geometric cover $f\colon U \surjto X$ by an affine scheme $U$ is descendable of index $\le d(X)$.
\end{itemize}
To prove ($*$), we note that $f$ is automatically affine because $X$ has affine diagonal. It follows easily that $K_f^{\tensor i}[i]$ is special for all $i \ge 0$: By \cref{rslt:qcqs-base-change-for-QCoh} this reduces to the case that $X$ is an affine scheme, where it is shown in \cref{ex:fppf-cover-of-rings-is-descendable}. Thus for $d = d(X)$, using the definition of $d(X)$ and the fact that $K_f^{\tensor d}[d]$ is special, we deduce that $H^0 \Hom(K_f^{\tensor d}, 1) = H^0 \Hom(K_f^{\tensor d}[d], 1[d]) = 0$. This immediately shows that $f$ is descendable of index $\le d$.

It is now enough to show that $d(X) < \infty$. In order to see this, we first prove the following intermediate claim:
\begin{itemize}
    \item[($**$)] Let $i\colon Z \to X$ be a finitary closed immersion with open complement $j\colon V \injto X$. If $d(Z) < \infty$ and $d(V) < \infty$ then $d(X) < \infty$.
\end{itemize}
To prove this claim, assume that both $d(Z)$ and $d(V)$ are finite. Let $\hat i\colon X_{\hat Z} \injto X$ denote the completion of $i$ and for $n \ge 0$ let $i_n\colon Z_n \to X$ denote the $n$-th infinitesimal neighborhood of $i$. We fix a special sheaf $M \in \QCoh(X)$. Then for any sheaf $N \in \QCoh(X)$, the excision sequence from \cref{rslt:excision-for-QCoh} together with adjunctions produce the following fiber sequence in $\D(\Lambda)$:
\begin{align*}
    \Hom(j^* M, j^? N) \to \Hom(M, N) \to \Hom(\hat i^* M, \hat i^* N)
\end{align*}
This reduces the proof of ($**$) to the following two subclaims:
\begin{enumerate}[(a)]
    \item If $N \in \QCoh^{\le d}(X)$ for $d \ll 0$ then $H^0 \Hom(\hat i^* M, \hat i^* N) = 0$.
    \item The functor $j^?\colon \QCoh(X) \to \QCoh(V)$ has finite cohomological dimension.
\end{enumerate}
Let us first prove (a). By \cref{rslt:completion-is-colimit-of-inf-nbhds-on-truncated} we have
\begin{align*}
    \Hom(\hat i^* M, \hat i^* N) = \varprojlim_n \Hom(i_n^* M, i_n^* N).
\end{align*}
Since $\varprojlim_n$ has cohomological dimension $\le 1$ in $\D(\Lambda)$ it suffices to prove that $H^0 \Hom(i_n^* M, i_n^* N) = 0$ if $N$ is connective enough (independent of $n$). More precisely, we claim that this is true if $N \in \QCoh^{\le-d(Z)}(X)$. The case $n = 1$ is true by definition of $d(Z)$. For $n \ge 2$ we rewrite the $\Hom$ as $H^0 \Hom(M, N_n)$ for $N_n := i_{n*} i_n^* N$ and by \cref{rslt:excision-for-QCoh} we have $\fib(N_{n+1} \to N_n) = i_* N'_n$ for $N'_n = i^* N \tensor \Sym^n_Z(L_{Z/X}[1])$. By induction we can assume that $H^0 \Hom(M, N_n) = 0$ and by definition of $d(Z)$ we know that $H^0 \Hom(M, i_* N'_n) = 0$. Hence also $H^0 \Hom(M, N_{n+1}) = 0$, as desired. This finishes the proof of subclaim (a).

It remains to prove subclaim (b). Since $j^? = j^* j_* j^?$ and $j^*$ is t-exact, it is enough to show that $j_* j^?$ is cohomologically bounded. By the excision sequence this reduces us to showing that $\hat i_* \hat i^*$ is cohomologically bounded. If $X$ is an affine scheme then this follows immediately from the formula $\hat i_* \hat i^* = \varprojlim_n i_{n*} i_n^*$ (from \cref{rslt:excision-for-QCoh}) because then $\varprojlim_n$ has cohomological dimension $\le 1$. In general, choose a flat geometric cover $f\colon X' \surjto X$ by an affine scheme $X'$ and denote $g\colon Y \to X_{\hat Z}$ the base-change of $f$ along $\hat i$ and $g_\bullet\colon Y_\bullet \to X_{\hat Z}$ its \v{C}ech nerve. By a similar argument as in the proof of ($*$) we deduce from (a) that $g$ is descendable. In particular, \cref{rslt:descendable-implies-M-fin-from-Cech-nerve} shows that for every $N \in \QCoh^{\le0}(X)$, $\hat i^* N$ can be constructed using retracts and finitely many fiber sequences from the objects $g_{n*} g_n^* \hat i^* N$. Thus it is enough to show that the functors $\hat i_* g_{n*} g_n^* \hat i^*$ are cohomologically bounded, reducing us to the case that $X$ is affine, handled above. This finishes the proof of ($**$).

With claim ($**$) at hand, we now use \cref{rslt:strata-decomposition-of-QCA-stack} to reduce the claim $d(X) < \infty$ to the case that $X = Y/\GL_{n,\Lambda}$ for a quasi-affine scheme $Y$ with an action by $\GL_{n,\Lambda}$. We can assume $\Lambda = \Q$. Let us denote $f\colon Y \surjto X$ the projection and also consider the projection $f_0\colon * \to */\GL_{n,\Q}$ of stacks over $\Q$, where $* = \Spec \Q$. Note that $\QCoh(*/\GL_{n,\Q})^\heartsuit$ is the category of algebraic representations of $\GL_{n,\Q}$ in the classical sense. In particular, the trivial representation is a direct summand of the regular representation, i.e. $1$ is a retract of $f_{0*} 1$ in $\QCoh(*/\GL_{n,\Q})$ (this is where we use the assumption that $\Lambda$ is a $\Q$-algebra!). Since $f$ is the base-change of $f_0$ along the natural map $Y/\GL_{n,\Q} \to */\GL_{n,\Q}$, we deduce from \cref{rslt:qcqs-base-change-for-QCoh} that $1$ is a retract of $f_* 1$ in $\QCoh(X)$ as well. By the projection formula (see \cref{rslt:qcqs-projection-formula--and-cocont-for-QCoh}) we deduce that for every $N \in \QCoh(X)$, $N$ is a retract of $f_* f^* N$. Hence for every $M \in \QCoh(X)$, $\Hom(M, N)$ is a retract of $\Hom(M, f_* f^* N) = \Hom(f^* M, f^* N)$. Thus if $d(Y) < \infty$ then $d(X) < \infty$.

We are now reduced to the case that $X$ is a quasicompact open subset of an affine scheme over $\Lambda$. Pick a finite open covering $X = \bigcup_{i\in I} U_i$ with affine open subsets $U_i \subset X$. Then the $\Hom$ in $\QCoh(X)$ is computed as the limit of the $\Hom$ in $\QCoh(U_{i_1} \cap \dots \cap U_{i_k})$ for all non-empty subsets $\{ i_1, \dots, i_k \} \subseteq I$. As this is a finite limit and hence has finite cohomological dimension, the claim $d(X) < \infty$ reduces to $d(U_{i_1} \cap \dots \cap U_{i_k}) < \infty$ for all non-empty subsets $\{ i_1, \dots, i_k \} \subseteq I$. But these are affine schemes, so we can assume that $X$ is affine. Then $d(X) \le 2$ by definition of special sheaves.
\end{proof}

\begin{rmk} \label{rmk:Gaitsgory-error}
In the proof of claim ($**$) in \cref{rslt:descendable-cover-of-QCA-stack} we followed the strategy of the proof of \cite[Proposition~8.3.2]{Gaitsgory-1-affineness}, but there appears to be a mistake in the proof of the reference: It is claimed in \cite[Lemma~8.3.5(a)]{Gaitsgory-1-affineness} that the functor $\hat i_* \hat i^*$ is right t-exact, but this is not clear, as it requires a bound on the cohomological dimension of $\varprojlim_n$ in $\QCoh(X)$. It is not clear that such a bound exists a priori (unleass $X$ is an affine scheme), but using descendability techniques we were able to salvage this claim (see the proof of subclaim (b) in \cref{rslt:descendable-cover-of-QCA-stack}).
\end{rmk}

With \cref{rslt:descendable-cover-of-QCA-stack} at hand, we can easily extend the base-change and projection formula from \cref{rslt:qcqs-base-change-for-QCoh} and \cref{rslt:qcqs-projection-formula--and-cocont-for-QCoh} to QCA maps. A similar discussion appears in \cite[Corollary~1.4.5]{DrinfeldGaitsgory-finiteness}.

\begin{prop} \label{rslt:QCA-pushforward-properties}
Fix a $\Q$-algebra $\Lambda$ and let $f\colon Y \to X$ be a QCA map of stacks over $\Lambda$.
\begin{propenum}
    \item \label{rslt:QCA-maps-have-base-change-and-projection-formula} The functor $f_*\colon \QCoh(Y) \to \QCoh(X)$ preserves all small colimits and satisfies base-change and the projection formula (in the sense of \cref{rslt:qcqs-base-change-for-QCoh} and \cref{rslt:qcqs-projection-formula--and-cocont-for-QCoh}).

    \item \label{rslt:QCA-pushforward-has-finite-cohomological-dimension} Suppose that $X$ is geometric (hence also $Y$), so that both $\QCoh(X)$ and $\QCoh(Y)$ have a t-structure. Then $f_*$ has finite cohomological dimension, i.e. sends $\QCoh^{\le0}(Y)$ to $\QCoh^{\le d}(X)$ for some $d \ge 0$.
\end{propenum}
\end{prop}
\begin{proof}
Let us first prove the base-change in (i). By the same argument as for \cref{rslt:qcqs-base-change-for-QCoh}, we are reduced to showing the base-change in the case that $X$ is affine and that we take the base-change along a map $g\colon X' \to X$ of affine schemes. Let us denote by $f'\colon Y' \to X'$ and $g'\colon Y' \to Y$ the base-changes of $f$ and $g$ and fix some $M \in \QCoh(Y)$. We have to show that the natural map $g^* f_* M \isoto f'_* g'^* M$ is an isomorphism. Pick a smooth cover $h\colon U \surjto Y$ by an affine scheme $U$ and let $h_\bullet\colon U_\bullet \to Y$ denote its \v{C}ech nerve. Note that $h$ is an affine map because the diagonal of $Y$ is affine. By \cref{rslt:descendable-cover-of-QCA-stack} $h$ is descendable and hence by \cref{rslt:descendable-implies-M-fin-from-Cech-nerve} we can reduce to the case that $M = h_{n*} M'$ for some $n \ge 0$ and some $M' \in \QCoh(U_n)$. But by \cref{rslt:qcqs-base-change-for-QCoh} we have base-change for $h_{n*}$ and we obviously have base-change for $(fh_n)_*$, so we easily arrive at the claim.

The remaining claims in (i) can be proved in a very similar manner, as follows. First note that by base-change for $f_*$ we can reduce to the case that $X$ is affine. We consider $h_\bullet\colon U_\bullet \to X$ as in the previous paragraph. Let us first check the projection formula, i.e. that for $N \in \QCoh(Y)$ and $M \in \QCoh(X)$ the natural map $f_* N \tensor M \isoto f_*(N \tensor f^* M)$ is an isomorphism. By \cref{rslt:descendable-implies-M-fin-from-Cech-nerve} we reduce to the case $N = h_{n*} N'$ for some $n \ge 0$ and $N' \in \QCoh(U_n)$. But then the claim follows easily from projection formula for $h_n$ and $fh_n$ (shown in \cref{rslt:qcqs-projection-formula--and-cocont-for-QCoh}). To finish the proof of (i) it remains to show that $f_*$ preserves colimits. Since it clearly preserves finite colimits (being exact), it is enough to show that $f_*$ preserves infinite direct sums. This follows again by the same strategy by using that $h_{n*} h_n^*$ and $(fh_n)_*$ both preserve direct sums (see \cref{rslt:qcqs-projection-formula--and-cocont-for-QCoh}). Part (ii) is proved analogously.
\end{proof}

\begin{rmk} \label{rmk:dualizability-of-QCoh-for-QCA-stacks}
By formal 6-functor arguments one can deduce from \cref{rslt:descendable-cover-of-QCA-stack} that for every QCA stack $X$ over a $\Q$-algebra $\Lambda$, the symmetric monoidal category $\QCoh(X)$ is rigid dualizable. We sketch the argument here but do not provide details, as we do not need the result. Similar results are for example obtained in \cite{Stefanich-Dualizability}.

First note that $\QCoh$ upgrades to a 6-functor formalism on stacks over $\Lambda$: On affine schemes, we declare all maps to be proper, i.e. $f_! = f_*$, which we can do by \cref{rslt:qcqs-base-change-for-QCoh} and \cref{rslt:qcqs-projection-formula--and-cocont-for-QCoh}. We then use \cite[Theorem~3.4.11]{heyer-mann-6ff} to extend this 6-functor formalism to stacks over $\Lambda$, for a certain class $E$ of $!$-able maps including all affine maps. Since affine maps are prim, it follows from \cite[Lemma~4.7.4]{heyer-mann-6ff} that descendable affine maps are universal $!$-covers. Thus \cref{rslt:descendable-cover-of-QCA-stack} implies that all QCA maps are $!$-able.

On the other hand, the 6-functor formalism on $\QCoh$ satisfies categorical Künneth, i.e. for two $!$-able stacks $X$ and $Y$ we have $\QCoh(X) \tensor_{\D(\Lambda)} \QCoh(Y) = \QCoh(X \times Y)$: This is easy to see for affine schemes and then follows by abstract nonsense for all $!$-able stacks (using $!$-colimit descent along a cover by affines). The categorical Künneth formula implies that the functor $\cat K_{\QCoh} \to \Mod_{\D(\Lambda)}(\PrL)$ from \cite[Remark~4.1.7]{heyer-mann-6ff} is fully faithful and symmetric monoidal. Since in $\cat K_{\QCoh}$ every object is canonically self-dual, the same must be true in its image, proving that $\QCoh(X)$ is dualizable. Moreover, if $X$ is QCA then rigidity of $\QCoh(X)$ follows by observing that both $X \to *$ and $\Delta_X\colon X \to X \times X$ are prim in the 6-functor formalism, so that $X$ is ``internally rigid'' in $\cat K_{\QCoh}$.

An even cleaner argument (still based on \cref{rslt:descendable-cover-of-QCA-stack}) can be obtained using Gestalten \cite{Scholze-Gestalten}.
\end{rmk}

\subsection{Solid quasicoherent sheaves} \label{sec:alggeo.solid}

In the previous subsections we introduced the category of (derived) schemes $X$ together with an associated category $\QCoh(X)$ of quasicoherent sheaves. However, one downside of $\QCoh$ is that it does not admit a 6-functor formalism that respects the geometric properties of maps of schemes. For example, the lower-$!$ functor along an open immersion of schemes should be the left adjoint of the pullback, but the pullback usually does not preserve limits and hence does not admit a left adjoint.

We will eventually fix the issue by passing to ind-coherent sheaves, as they are the ones that appear in categorical local Langlands. However, there is a more canonical coherent 6-functor formalism on schemes, which will play a crucial role in setting up ind-coherent sheaves, namely the 6-functor formalism for solid quasicoherent sheaves due to Clausen--Scholze. In the following we briefly introduce this 6-functor formalism following \cite{condensed-mathematics,scholze-analytic-spaces,scholze-complex-geometry,RC.SolidGeometry} and \cite[\S2]{mann-p-adic-6-functors}.

The main idea is to allow a certain kind of \enquote{topological} quasicoherent sheaves, i.e. sheaves that on an affine scheme $\Spec A$ look like \enquote{complete topological} $A$-modules. The precise terminology is \emph{solid} $A$-modules, which we introduce now.

\begin{defn}
\begin{defenum}
    \item A \emph{(light) profinite set} is a topological space that can be written as a (countable) inverse limit of finite discrete sets.
    
    \item Let $\Cond(\Z)$ denote the category of (light) \emph{condensed $\Z$-modules}, i.e. the derived category of the category of sheaves of abelian groups on the site of light profinite sets. For a ring $A$ we denote by $\Cond(A) := \Mod_A(\Cond(\Z))$ the category of \emph{condensed $A$-modules}. If $M$ is a condensed $A$-module then for every light profinite set $S$ we get an associated $A$-module $M(S) \in \D(A)$. We call $M(*)$ the \emph{underlying $A$-module} of $M$.

    The category $\Cond(A)$ is stable (i.e. \enquote{triangulated}) and presentable, and comes equipped with a $t$-structure. Moreover, it has a natural tensor product $\tensor_A$ induced by the usual tensor product of sheaves on sites.
\end{defenum}
\end{defn}

Intuitively a condensed $A$-module $M$ is a topological $A$-module, where $M(S)$ captures the $A$-module of continuous maps $S \to M$. Indeed, for $A = \Z$ this recipe defines a functor from topological abelian groups to (static) condensed abelian groups, which is fully faithful on metrizable compactly generated abelian groups; we will implicitly use this functor to identify topological abelian groups and rings with their condensed version. As an important special case we get:

\begin{defn}
Equipping an $A$-module with the discrete topology induces an embedding
\begin{align*}
    \D(A) \subset \Cond(A),
\end{align*}
for every ring $A$. For $M \in \D(A)$ the associated condensed $A$-module is given by sending $S \mapsto M(S) := \varinjlim_n \bigoplus_{S_n} M$ for a light profinite set $S = \varprojlim_n S_n$. We call these condensed $A$-modules \emph{discrete}.
\end{defn}

If $A$ is an adic ring (i.e. there is a finitely generated ideal $I$ of $A$ such that $A$ is $I$-adically complete) then there is an evident notion for an $A$-module $M$ to be complete: We can ask for it to be adically complete, for any ideal of definition $I$ of $A$. For a \emph{condensed} $A$-module $M$, Clausen--Scholze construct a related notion of \enquote{completeness} that takes into account the inherent topology of the module:

\begin{defn}
Let $A$ be a ring. A condensed $A$-module $M \in \Cond(A)$ is called \emph{solid} if for every finite-type static $\Z$-algebra $A_0$ with a map $A_0 \to A$, $M$ viewed as a condensed $A_0$-module lies in the full subcategory of $\Cond(A_0)$ generated under limits and colimits by $A_0$. We denote by
\begin{align*}
    \D_\solid(A) \subseteq \Cond(A)
\end{align*}
the full subcategory spanned by the solid $A$-modules. The inclusion $\D_\solid(A) \injto \Cond(A)$ admits a left adjoint denoted
\begin{align*}
    (-)^\solid\colon \Cond(A) \to \D_\solid(A), \qquad M \mapsto M^\solid
\end{align*}
The category $\D_\solid(A)$ is a stable category and the t-structure on $\Cond(A)$ restricts to a t-structure on $\D_\solid(A)$. Clearly every discrete $A$-module is solid, so we obtain fully faithful embeddings
\begin{align*}
    \D(A) \subseteq \D_\solid(A) \subseteq \Cond(A).
\end{align*}
For every light profinite set $S$, we denote by $A_\solid[S] := A[S]^\solid$ the associated free solid $A$-module on $S$. These are compact projective in $\D_\solid(A)$, in fact the single object
\begin{align*}
    A_\solid[\N_\infty] = \varinjlim_{A_0 \to A} \prod_\N A_0,
\end{align*}
generates $\D_\solid(A)$; here the colimit is taken over all finite-type $\Z$-algebras $A_0$ with a map to $A$. In particular, a condensed $A$-module is solid if and only if it can be written as a colimit of copies of $A_\solid[\N_\infty]$.
\end{defn}

\begin{defn}
Let $A$ be a ring.
\begin{defenum}
    \item The tensor product $\tensor_A$ on $\Cond(A)$ induces a tensor product on $D_\solid(A)$ via
    \begin{align*}
        M \tensor_A^\solid N := (M \tensor_A N)^\solid
    \end{align*}
    The embedding $\D(A) \injto \D_\solid(A)$ is symmetric monoidal, i.e. commutes with the tensor products on $\D(A)$ and $\D_\solid(A)$.

    \item Let $A \to B$ be a map of rings. Then we get an induced base-change functor
    \begin{align*}
        - \tensor_{A_\solid} B_\solid\colon D_\solid(A) \to \D_\solid(B), \qquad M \mapsto (M \tensor_A B)^\solid,
    \end{align*}
    where the solidification on the right is with respect to $B$. Since discrete modules are solid, the restriction of $- \tensor_{A_\solid} B_\solid$ to $\D(A)$ is the usual $- \tensor_A B$.
\end{defenum}
\end{defn}

We want to extend the assignment $A \mapsto \D_\solid(A)$ to stacks, for which we need to verify fppf descent. For static rings, this descent appeared first independently in \cite[Lemma~2.10.6]{mann-p-adic-6-functors} and \cite{mikami-fppf-descent1}, and was then extended to all (animated) rings in \cite{mikami-fppf-descent2}:

\begin{prop} \label{rslt:fppf-descent-for-QCoh-solid}
Let $f\colon A \to B$ be an fppf map of rings and let $f^\bullet\colon A \to B^\bullet$ be the induced cosimplicial object. Then $\D_\solid(-)$ descends along $f$, i.e. the natural functor
\begin{align*}
    - \tensor_{A_\solid} B_\solid^\bullet\colon \D_\solid(A) \isoto \varprojlim_{n\in\Delta} \D_\solid(B^n)
\end{align*}
is an equivalence of categories.
\end{prop}
\begin{proof}
See \cite[Theorem~0.5]{mikami-fppf-descent2}.
\end{proof}

\begin{rmk}
By \cite[Theorem~2.7.8]{mann-p-adic-6-functors} one can deduce that $\D_\solid(-)$ even satisfies ind-fppf descent. However, we warn the reader that faithfully flat descent is false in general.
\end{rmk}

\begin{defn} \label{def:QCoh-solid}
Fix a ring $\Lambda$. We denote by
\begin{align*}
    \QCoh_\solid\colon \Stk_\Lambda^\op \to \CMon, \qquad X \mapsto \QCoh_\solid(X)
\end{align*}
the sheaf of symmetric monoidal categories extended from $\QCoh_\solid(\Spec A) = \D_\solid(A)$ via descent (using \cref{rslt:fppf-descent-for-QCoh-solid}). For $X \in \Stk_\Lambda$, the objects in $\QCoh_\solid(X)$ are called the \emph{solid quasicoherent sheaves} on $X$. There is a natural fully faithful transformation
\begin{align*}
    \QCoh \injto \QCoh_\solid
\end{align*}
of sheaves of symmetric monoidal categories.
\end{defn}

We warn the reader that while there is a natural t-structure on $\QCoh(X)$ for every geometric stack $X$ (see \cref{rslt:t-structure-on-QCoh}), there is in general no good t-structure on $\QCoh_\solid(X)$, even though there is a t-structure on $\D_\solid(A)$ for every ring $A$. The problem is that for a flat map $A \to B$ the base-change $- \tensor_{A_\solid} B_\solid$ is rarely t-exact, even for open immersions. On the plus side, by \cite[Lemma~4.10]{mikami-fppf-descent2} the base-change $- \tensor_{A_\solid} B_\solid$ has Tor dimension $\le 1$ and hence preserves boundedness.

Let us now come to the promised 6-functor formalism. It first appeared in \cite{condensed-mathematics}, was fleshed out in \cite[\S2.9]{mann-p-adic-6-functors} and subsequent works by Clausen--Scholze and others, including \cite{RC.SolidGeometry}. In the following we use the notation of 6-functor formalisms from \cite[\S3, \S4]{heyer-mann-6ff}.

\begin{thm} \label{rslt:6ff-on-QCoh-solid}
Fix a ring $\Lambda$. There is a class $E$ of maps in $\Stk_\Lambda$ and a 6-functor formalism
\begin{align*}
    \QCoh_\solid\colon \Corr(\Stk_\Lambda, E) \to \Cat
\end{align*}
with the following properties:
\begin{thmenum}
    \item The underlying functor $\QCoh_\solid^*\colon \Stk_\Lambda^\op \to \CMon$ (where functoriality is given by pullback) coincides with the sheaf from \cref{def:QCoh-solid}.

    \item \label{rslt:QCoh-solid-on-proper-and-smooth-maps} $E$ contains every lafp geometric map $f\colon Y \to X$ of stacks over $\Lambda$. Moreover, we have the following description:
    \begin{enumerate}[(a)]
        \item If $f$ is proper and schematic then it is $\QCoh_\solid$-prim with codualizing sheaf $\delta_f \isom 1$, so we have an isomorphism $f_! \isom f_*$.

        \item If $f$ is flat then it is $\QCoh_\solid$-suave, so $f^! = f^* \tensor \omega_f$ for the dualizing complex $\omega_f := f^! 1 \in \QCoh_\solid(Y)$.
        
        \item If $f$ is a local complete intersection then it is $\QCoh_\solid$-smooth and $\omega_f = \det(L_f)$. If $f$ is étale then $\omega_f = 1$ and thus $f^! = f^*$.
    \end{enumerate}

    \item \label{rslt:stability-of-QCoh-solid-fine-maps} The class $E$ enjoys the following stability properties:
    \begin{enumerate}[(a)]
        \item $E$ is stable under composition and pullback. Moreover, if $f\colon Y \to X$ and $g\colon Z \to Y$ are maps in $\Stk_\Lambda$ such that $f, fg \in E$, then $g \in E$.
        
        \item The class $E$ is local on the target: If $f\colon Y \to X$ is a map in $\Stk_\Lambda$ and there is a cover of $X$ such that the pullback of $f$ to each element in that cover lies in $E$, then $f$ is in $E$.

        \item The class $E$ is $!$-local on source and target, i.e. containment of $f\colon Y \to X$ in $E$ can be checked after passing to a $!$-cover of $Y$ and $X$. An important class of $!$-covers are covers by $\QCoh_\solid$-suave maps, see \cite[Lemma~4.7.1]{heyer-mann-6ff}.
    \end{enumerate}
\end{thmenum}
\end{thm}
\begin{proof}
Most of this is done in \cite[Proposition~2.9.31]{mann-p-adic-6-functors}. In the following we sketch the construction and fill in the missing parts. As $f^*$ and $\tensor$ are already encoded in \cref{def:QCoh-solid}, the main task is to construct the functor $f_!$.

We first construct the 6-functor formalism on affine schemes. In this case the construction is purely algebraic (i.e. no Nagata compactifications are needed) if one is willing to replace the category of rings by the category of \emph{discrete Huber pairs}, i.e. pairs $(A, A^+)$ consisting of a $\Lambda$-algebra $A$ together with a map $A^+ \to A$ from some ring $A^+$ (more precisely, one considers such pairs up to a certain equivalence relation: A class of pairs only depends on the integral closure of the image of $\pi_0 A^+$ in $\pi_0 A$). One then defines $\D_\solid(A, A^+) := \Mod_A(\D_\solid(A^+)) \subseteq \Cond(A)$, so in particular we recover $\D_\solid(A) = \D_\solid(A, A)$. The advantage of working with discrete Huber pairs is that every map $f\colon (A, A^+) \to (B, B^+)$ of discrete Huber pairs factors into the maps $g\colon (A, A^+) \to (B, A^+)$ and $j\colon (B, A^+) \to (B, B^+)$. Now we deem $g$ a proper map and set $g_! := g_*$, while if $A \to B$ is afp then $j$ is an open immersion and thus we let $j_!$ be a left adjoint of $j^*$. We then define $f_! := g_! j_!$. This construction glues to a 6-functor formalism by \cite[Proposition~3.3.3]{heyer-mann-6ff}, the relevant conditions (i.e. projection formula and base-change for $g_*$ and $j_!$) are verified in the proof of \cite[Proposition~2.9.31]{mann-p-adic-6-functors}. By restricting to Huber pairs with $A^+ = A$ we obtain the 6-functor formalism
\begin{align*}
    \D_\solid\colon \Corr(\Alg_\Lambda^\op, E_0) \to \Cat,
\end{align*}
where $E_0$ is the class of almost finitely presented maps. In the next step we use \cite[Theorem~3.4.11]{heyer-mann-6ff} to extend $\D_\solid$ to the 6-functor formalism $\QCoh_\solid$ on $\Stk_\Lambda$. By construction it clearly satisfies (i).

We next verify (iii). Part (a) simply says that $E$ defines a geometric setup (see \cite[\S2.1]{heyer-mann-6ff} and particularly \cite[Lemma~2.1.5]{heyer-mann-6ff}), which is true by construction. Part (c) is a reformulation of \cite[Theorem~3.4.11(iii)]{heyer-mann-6ff}. It remains to prove (b), so let $f\colon Y \to X$ be a map in $\Stk_\Lambda$ and assume that there is a cover $(X_i \to X)_i$ of stacks such that each pullback $f_i\colon Y_i \to X_i$ lies in $E$. To show that $f \in E$, it is enough to do so after pullback to every affine scheme $X'$ with a map $X' \to X$ (see \cite[Theorem~3.4.11(ii)]{heyer-mann-6ff}), so we may w.l.o.g. assume that $X = \Spec A$ is affine. By \cite[Lemma~A.3.9]{mann-p-adic-6-functors} there is a flat cover $(Z_j \to X)_j$ in the site $\Alg_\Lambda^\op$ such that each $Z_j$ factors over some $X_i$. Then the pullback of $f$ to each $Z_j$ lies in $E$, so by (c) it is enough to verify that $(Z_j \to X)_j$ is a universal $!$-cover, which in turn is implied if each map $Z_j \to X$ is $\QCoh_\solid$-suave. This is true by (ii).(b), which will be proved below.

It remains to prove (ii). We will defer the proof of the $\QCoh_\solid$-suaveness of flat maps of rings to a more general result on suaveness below (see \cref{rslt:characterization-of-suave-maps}). Assuming this for now, then in particular every standard open immersion of rings is $\QCoh_\solid$-suave. Since suaveness is local on the target and suave-local on the source (see \cite[Lemmas~4.5.7,~4.5.8]{heyer-mann-6ff}) one easily deduces that every open immersion of stacks lies in $E$ and is $\QCoh_\solid$-suave. In particular, containment in $E$ can be checked on an open cover by (iii).(c), which immediately shows that lafp maps of schemes are in $E$ (as afp maps of affine schemes lie in $E$ by construction). Moreover, flat lafp maps of schemes are $\QCoh_\solid$-suave, as this can be checked suave locally on the source (see \cite[Lemma~4.5.8(i)]{heyer-mann-6ff}). One now inductively deduces that lafp $n$-geometric maps $f\colon Y \to X$ lie in $E$ and if $f$ is flat then it is $\QCoh_\solid$-suave: Both claims are local on $X$, so we can assume that $X$ is an affine scheme; then we can find a $\QCoh_\solid$-suave cover of $Y$ by a lafp scheme over $X$ and hence deduce that $f$ lies in $E$ by (iii).(c) and \cite[Lemma~4.7.1]{heyer-mann-6ff}. If $f$ is flat then suaveness is deduced from \cite[Lemma~4.5.8(i)]{heyer-mann-6ff}.

We now prove parts (a), (b) and (c) of (ii). Part (a) is shown in \cite[Proposition~2.9.31(ii)]{mann-p-adic-6-functors} using the more general version of $\QCoh_\solid$ on discrete adic spaces (the main advantage being that then properness can be checked locally). Part (b) was already shown above.

We now prove (iii).(c), so let $f\colon Y \to X$ be a geometric local intersection morphism. We first check that $f$ is suave. By \cref{rslt:characterization-of-lci-maps-via-regular-closed-immersions} $f$ is smooth locally on $Y$ given by a composition of a smooth map and a pullback of the map $i\colon \{ 0 \} \to \mathbb A^n$. Since smooth maps are flat and hence $\QCoh_\solid$-suave by (b) and $\QCoh_\solid$-suaveness is $\QCoh_\solid$-suave local on the source and stable under condition and base-change by \cite[Lemmas 4.5.7, 4.5.8, 4.5.9]{heyer-mann-6ff}, we reduce to showing that $i$ is suave. But this follows immediately from \cref{rslt:characterization-of-suave-maps} below, as $\Lambda$ is perfect as a $\Lambda[x_1, \dots, x_n]$-module (via the Koszul complex).

It remains to identify $\omega_f := f^! 1$ with $\det(L_f)$. This is a well-known folklore statement, but currently lacks a proof in the literature. In the case that $f$ is schematic and smooth, an argument is sketched in \cite[Theorem~13.6]{scholze-complex-geometry} using a deformation to the normal cone. There have since been generalizations of this construction to more general maps $f$, see e.g. \cite{HekkingKhanRydh-NormalCone}, which should lead to a general proof of the claimed identification of $\omega_f$, cf. \cite[\S0.4.3]{HekkingKhanRydh-NormalCone}. As we only need this statement in one particular example (the dualizing complex of $\Par_G$), we content ourselves with these references to future work.
\end{proof}

\begin{exmpl} \label{ex:solid-lower-shriek-for-A1}
Consider the map $f\colon \Spec \Z[x] \to \Spec \Z$. This map is smooth and its dualizing complex is given by $\Z[x][1] \in \D(\Z[x]) \subset \D_\solid(\Z[x])$. In particular $f^! = f^*[1]$. The functor $f_!$ is given by
\begin{align*}
    f_!\colon \D_\solid(\Z[x]) \to \D_\solid(\Z), \qquad M \mapsto M \tensor_{(\Z[x], \Z)}^\solid (\Z[x]_\infty/\Z[x])[-1].
\end{align*}
Here $\Z[x]_\infty := \Z((x^{-1}))$ is the Laurent series ring in $x^{-1}$ (viewed as a condensed $\Z[x]$-module via the natural topology) and $\tensor^\solid_{(\Z[x], \Z)}$ denotes the tensor product in $\D_\solid(\Z[x], \Z) = \Mod_{\Z[x]}(\D_\solid(\Z))$. Intuitively one should see $\Z[x]_\infty$ as the ring of functions that are defined in an infinitesimal neighborhood at infinity, so that
\begin{align*}
    (\Z[x]_\infty/\Z[x])[-1] = \fib(\Z[x] \to \Z[x]_\infty)
\end{align*}
is the space of those functions on $\Spec \Z[x]$ that \enquote{vanish in a neighborhood of $\infty$}. Observe that we can only make sense of this vanishing in the derived world, as the map $\Z[x] \to \Z[x]_\infty$ has trivial kernel.
\end{exmpl}

\begin{defn}
Fix a ring $\Lambda$. A map $f\colon Y \to X$ in $\Stk_\Lambda$ is \emph{$\QCoh_\solid$-fine} if it lies in the class $E$ from \cref{rslt:6ff-on-QCoh-solid}.
\end{defn}

Using the 6-functor formalism on $\QCoh_\solid$ we can now extend several of the notions from \cref{sec:alggeo.stacks} to more general maps of stacks, as well as introduce certain new properties.

\begin{defn}
Fix a ring $\Lambda$ and a $\QCoh_\solid$-fine map $f\colon Y \to X$ in $\Stk_\Lambda$.
\begin{defenum}
    \item \label{def:suave-and-coh-smooth-maps} $f$ is called \emph{suave}, resp. \emph{cohomologically smooth} if is it $\QCoh_\solid$-suave, resp. $\QCoh_\solid$-smooth.
    \item $f$ is called \emph{proper} if one of the following is true: Either on a cover of $X$ the map $f$ becomes a proper map of schemes, or the diagonal of $f$ is proper and $f$ is $\QCoh_\solid$-prim.
\end{defenum}
\end{defn}

\begin{lem}
Let $\Lambda$ be a ring.
\begin{lemenum}
    \item \label{rslt:stability-of-smooth-suave-and-proper} The properties of being suave, cohomologically smooth, and proper are stable under composition and base-change and can be checked on any cover of the target. The former two properties are additionally suave local, resp. cohomologically smooth local, on the source.

    \item If $f\colon Y \to X$ and $g\colon Z \to Y$ are maps of stacks over $\Lambda$ such that $f$ and $fg$ are proper, then $g$ is proper.
    
    \item If $f\colon Y \to X$ is a proper map in $\Stk_\Lambda$ then there is an isomorphism $f_! \isom f_*$ of functors $\QCoh_\solid(Y) \to \QCoh_\solid(X)$.
    
    \item \label{rslt:characterization-of-suave-maps} Let $f\colon A \to B$ be an almost finitely presented map of $\Lambda$-algebras. Then the following are equivalent:
    \begin{enumerate}[(a)]
        \item The map $\Spec B \to \Spec A$ is suave.
        \item $f$ has finite Tor dimension, i.e. the functor $- \tensor_A B\colon \D(A) \to \D(B)$ sends $\D^{\ge0}(A)$ to $\D^{\ge d}(B)$ for some integer $d$.
        \item There is a factorization $A \to A[x_1, \dots, x_n] \to B$ such that $B$ is perfect as a module over $A[x_1, \dots, x_n]$.
    \end{enumerate}
    In this case the dualizing complex of $f$ is given by the dual of $B$ in $\D(A[x_1, \dots, x_n])$. In particular this dualizing complex is discrete.
\end{lemenum}
\end{lem}
\begin{proof}
Parts (i) and (ii) follow easily from the respective properties of suave and prim maps, see \cite[Lemmas 4.5.7, 4.5.8, 4.5.9]{heyer-mann-6ff}, together with the stability properties of $\QCoh_\solid$-fine maps by \cref{rslt:stability-of-QCoh-solid-fine-maps}. The version for properness requires a little care and can be shown by inducting on the number of times one has to pass to the diagonal before the map becomes representable in proper maps of schemes. Part (iii) follows in the same way as in \cite[Lemma~4.6.4(ii)]{heyer-mann-6ff}.

We now prove (iv), so let $f\colon A \to B$ be given. Let us first check that (c) implies (a). For this we first show that the map $\Spec(A[x_1, \dots, x_n]) \to \Spec A$ is suave. In fact it is even cohomologically smooth: By (i) this reduces immediately to the case $\Spec \Z[x] \to \Spec \Z$, where it is an explicit computation (see e.g. \cite[Proposition~2.9.31(iv)]{mann-p-adic-6-functors}). Next we observe that if $B$ is perfect as an $A$ module then $\Spec B \to \Spec A$ is $\QCoh_\solid$-suave. Indeed, in this case the map is finite and hence $\QCoh_\solid$-prim by \cref{rslt:QCoh-solid-on-proper-and-smooth-maps}, so by \cite[Lemma~4.5.10]{heyer-mann-6ff} the map $\Spec B \to \Spec A$ is suave if and only if $B$ is perfect as an $A$-module. Altogether we see that (c) implies (a). We now show that (a) implies (c), so assume that $f$ is suave and pick a factorization $A \to A[x_1, \dots, x_n] \to B$ such that $B$ is pseudocoherent as a module over $A[x_1, \dots, x_n]$. By what we have seen above, the diagonal of $\Spec A[x_1, \dots, x_n] \to \Spec A$ is suave (because it is finite and perfect as a module), so from \cite[Lemma~4.5.9(ii)]{heyer-mann-6ff} we deduce that the map $\Spec B \to \Spec A[x_1, \dots, x_n]$ is suave. By the above argument this implies that $B$ is perfect as a module over $A[x_1, \dots, x_n]$, as desired.

It remains to prove that (b) is equivalent to (a) and (c). Clearly (c) implies (b). To prove the converse, assume now that $f$ has finite Tor dimension and pick a factorization $A \to A[x_1, \dots, x_n] \to B$ such that $B$ is pseudocoherent as a module over $A[x_1, \dots, x_n]$. Note that having finite Tor dimension is stable under composition and base-change and also the diagonal of $\Spec A[x_1, \dots, x_n] \to \Spec A$ has finite Tor dimension (one can see this e.g. in the universal case $A = \Z$ and further reduce to the case $n = 1$, where it is easy). Hence by the argument in \cite[Lemma~2.1.5]{heyer-mann-6ff} the finite Tor dimension of $A \to B$ implies the finite Tor dimension of $A[x_1, \dots, x_n] \to B$. But then \cref{rslt:pseudo-coherent-plus-fin-tor-dim-implies-perfect} shows that $B$ is perfect over $A[x_1, \dots, x_n]$, as desired.
\end{proof}

\begin{rmk} \label{rmk:suaveness-in-the-literature}
Using \cref{rslt:characterization-of-suave-maps} we see that suave maps appear as eventually coconnective maps in \cite[Definition~3.1.2]{gaitsgory-rozenblyum-vol1} (see \cite[Lemma~3.1.3]{gaitsgory-rozenblyum-vol1}) and as perfect morphisms in \cite[Section~0685]{stacks-project}. Our terminology is derived from the perspective of the quasi-coherent six-functor formalism, which also provides a robust notion of suaveness for stacks. Suave maps are abundant in algebraic geometry, e.g. any lci or smooth map is suave (even cohomologically smooth, see \cref{rslt:QCoh-solid-on-proper-and-smooth-maps}), and they enjoy particularly nice properties with respect to quasi-coherent and ind-coherent sheaves. However, suaveness is not stable under diagonals, and in particular \cref{rslt:afp-maps-are-cancellative} fails for suave maps.
\end{rmk}

\subsection{Perfect and relatively coherent sheaves} \label{sec:perf-and-coh}

One of the main goals of \cref{sec:alggeo} is the construction of the category of ind-coherent sheaves on a lafp stack, so in particular we need a good notion of coherent sheaves and the related notion of perfect sheaves. While both can easily defined by hand in the relevant settings, the present subsection provides a different perspective on the definition of coherent sheaves which we find rather enlightening and which explains (and proves!) all of the expected properties in a conceptual way. The same approach works almost verbatim in different geometric settings, e.g. formal schemes, rigid geometry and complex analytic geometry, so in particular provides an easy proof of the finiteness of coherent cohomology in each of these settings (see \cref{rmk:Coh-in-other-geometries} below). Without further ado, let us jump right into the main definitions:

\begin{defn}
Let $\Lambda$ be a ring, $X \in \Stk_\Lambda$ and $P \in \QCoh_\solid(X)$.
\begin{defenum}
    \item $P$ is \emph{perfect} if it is a perfect (discrete) module after pullback to every affine scheme mapping to $X$. We denote by
    \begin{align*}
        \Perf(X) \subseteq \QCoh(X) \subseteq \QCoh_\solid(X)
    \end{align*}
    the full subcategory of perfect sheaves.

    \item \label{def:relatively-coherent-sheaves} Let $f\colon X \to S$ a $\QCoh_\solid$-fine map of stacks over $\Lambda$. $P$ is called \emph{$f$-coherent} if it is $f$-suave (see \cite[Definition~4.4.1]{heyer-mann-6ff}) in the $\QCoh_\solid$ 6-functor formalism. We denote by
    \begin{align*}
        \Coh_S(X) := \Coh_f(X) \subseteq \QCoh_\solid(X)
    \end{align*}
    the full subcategory spanned by the $f$-coherent sheaves. In the case that $X = \Spec A$ and $S = \Spec R$ are affine, we also use the notation $\Coh_R(A) \subseteq \D_\solid(A)$.
\end{defenum}
\end{defn}

Perfect sheaves admit a more conceptual definition in terms of dualizable objects, which can be defined in any symmetric monoidal category (see e.g. \cite[Definition~B.1.15]{heyer-mann-6ff}):

\begin{lem} \label{rslt:Perf-same-as-dualizable}
Let $\Lambda$ be a ring, $X \in \Stk_\Lambda$ and $P \in \QCoh_\solid(X)$. Then $P$ is perfect if and only if it is dualizable.
\end{lem}
\begin{proof}
If we a priori know that $P$ is discrete, i.e. lies in $\QCoh(X) \subseteq \QCoh_\solid(X)$, then this is a classical result, see e.g. \cite[Section~0FPP]{stacks-project}. In the following we provide a quick general argument without this a priori assumption.

Note first that dualizable objects in a symmetric monoidal category satisfy the same descent as the ambient symmetric monoidal category (for example, this is a special case of \cite[Proposition~D.2.8]{heyer-mann-6ff} using \cite[Example~D.2.3]{heyer-mann-6ff}). This reduces the claim to the case that $X = \Spec A$ is affine, hence $P \in \D_\solid(A)$ is now a solid $A$-module. Since $A$ is dualizable and dualizability in $\QCoh_\solid(A)$ is clearly stable under finite limits/colimits and retracts, every perfect $A$-module is dualizable. Conversely, assume now that $P$ is dualizable. Note that then $\Hom(P, -) = \Hom(A, P^\vee \tensor -)$, showing that $P$ is compact. We now show that $P$ is discrete. Given a light profinite set $S = \varprojlim_n S_n$ we have
\begin{align*}
    & P(S) = \Hom(A_\solid[S], P) = \Hom(P^\vee, A_\solid[S]^\vee) = \Hom(P^\vee, \varinjlim_n A[S_n]) =\\&\qquad=\varinjlim_n \Hom(P^\vee, A[S_n]) = \varinjlim_n P(S_n),
\end{align*}
proving that $P$ is indeed discrete. Altogether we see that $P$ is a compact object in $\D(A)$. But since $\D(A)$ is compactly generated by $A$ this implies that $P$ is a retract of a finite colimit of copies of $A$, proving that $P$ is perfect.
\end{proof}

From the general properties of suave objects in a 6-functor formalism discussed in \cite[\S4.4]{heyer-mann-6ff} we can immediately deduce the following list of properties of perfect and relatively coherent sheaves.

\begin{prop} \label{rslt:properties-of-Coh-and-Perf}
Let $\Lambda$ be a ring and $f\colon X \to S$ a $\QCoh_\solid$-fine map in $\Stk_\Lambda$.
\begin{propenum}
    \item \label{rslt:colim-stability-and-self-duality-for-Perf-and-Coh} $\Perf(X)$ and $\Coh_S(X)$ are stable under finite limits and colimits (i.e. under \enquote{exact triangles}) and all retracts in $\QCoh_\solid$. Moreover, $\Perf(X)$ is self-dual under $\intHom(-, \calO_X)$ and $\Coh_S(X)$ is self-dual under $\intHom(-, f^! \calO_S)$.

    \item \label{rslt:descent-for-Coh-and-Perf} $\Perf(-)$ is stable under arbitrary pullback and descends along all covers. $\Coh_S(-)$ is stable under suave pullback and descends along suave covers.

    \item \label{rslt:Coh-stable-under-proper-pushforward} If $g\colon Y \to X$ is proper then $g_*$ preserves $S$-coherent objects and hence restricts to a functor $g_*\colon \Coh_S(Y) \to \Coh_S(X)$.

    \item \label{rslt:Coh-S-and-smoothness-diagonal} If $f$ is suave then $\Perf(X) \subseteq \Coh_S(X)$. If the diagonal of $f$ is suave then $\Coh_S(X) \subseteq \Perf(X)$. In particular $\Coh_X(X) = \Perf(X)$.

    \item There is a natural action of $\Perf(X)$ on $\Coh_S(X)$ induced by the tensor product in $\QCoh_\solid(X)$.
\end{propenum}
\end{prop}
\begin{proof}
In part (i) the claim for $\Perf(X)$ is easy and for $\Coh_S(X)$ is found in \cite[Lemmas~4.4.4, 4.4.5]{heyer-mann-6ff} and \cite[Corollary~4.4.13]{heyer-mann-6ff}. In part (ii) the claim for $\Perf(-)$ follows from \cref{rslt:Perf-same-as-dualizable} and the claim for $\Coh_S(-)$ follows from \cite[Lemma~4.5.16(i), Lemma~4.4.8(i)]{heyer-mann-6ff}. Part (iii) is \cite[Lemma~4.5.16(ii)]{heyer-mann-6ff}. Part (iv) is \cite[Lemma~4.5.17]{heyer-mann-6ff}. For (v) we write $\Perf(X) = \Coh_X(X)$ and observe that this is the category of left adjoint morphisms $X \to X$ over $X$ in the language of the category of kernels from \cite[\S4]{heyer-mann-6ff}. Via the functoriality of the category of kernels over different bases (see \cite[Theorem~4.2.4]{heyer-mann-6ff}) we can view every $\mathcal M \in \Coh_X(X)$ also as a left adjoint morphism $X \to X$ over $S$. Then the claimed action comes from composing this left adjoint morphism with the left adjoint morphism $X \to S$ represented by a sheaf in $\Coh_S(X)$.
\end{proof}

Using the above general properties we can now provide a more explicit description of relatively coherent sheaves and then also relate them to the notion of coherent sheaves used in \cite{gaitsgory-rozenblyum-vol1}.

\begin{cor}
Let $\Lambda$ be a ring and let $f\colon X \to S$ be a lafp geometric map of stacks over $\Lambda$.
\begin{corenum}
    \item Every $f$-coherent sheaf is discrete, i.e. $\Coh_S(X) \subset \QCoh(X)$.
    
    \item \label{rslt:Coh-on-smooth-stack-equals-Perf} If $f$ is smooth then $\Coh_S(X) = \Perf(X)$.

    \item \label{rslt:characterization-of-suave-sheaves-in-affine-case} If $X = \Spec A$ and $S = \Spec R$ are affine schemes, then $\Coh_R(A) \subseteq \D_\solid(A)$ consists precisely of those $P \in \D_\solid(A)$ which satisfy the following condition:
    \begin{itemize}
        \item[($*$)] For some (equivalently every) factorization $R \to R[x_1, \dots, x_n] \to A$ such that $A$ is finite over $R[x_1, \dots, x_n]$, $P$ is perfect as an $R[x_1, \dots, x_n]$-module.
    \end{itemize}

    \item \label{rslt:coh-over-regular-noetherian-base} Suppose $S$ is a classical regular noetherian scheme. Then $\Coh_S(X) = \Coh(X)$ is the category of those quasicoherent sheaves on $X$ which after pullback along every smooth geometric map $\Spec A \to X$ correspond to a bounded $A$-module whose cohomologies are finitely generated $\pi_0 A$-modules.
\end{corenum}
\end{cor}
\begin{proof}
Part (i) follows immediately from (iii) and part (iv) is an easy consequence of (iii) together with \cref{rslt:descent-for-Coh-and-Perf} and \cite[Lemma~066Z]{stacks-project}. It thus remains to prove (ii) and (iii).

We first prove (ii), so assume that $f$ is smooth. Then $f$ is suave and by \cref{rslt:map-between-smooth-stacks-is-lci} the diagonal of $f$ is lci and hence also suave (by \cref{rslt:QCoh-solid-on-proper-and-smooth-maps}). Thus the claim follows from \cref{rslt:Coh-S-and-smoothness-diagonal}.

We now prove (iii), so assume that everything is as in the claim. We can factor $f$ as $R \to R[x_1, \dots, x_n] \to A$ for some $n \ge 0$ such that the map $R[x_1, \dots, x_n] \to A$ is finite and hence proper. Then by \cite[Lemma~4.5.17]{heyer-mann-6ff} we have $\Coh_R(A) = \Coh_{R[x_1,\dots,x_n]}(A)$. We can therefore replace $R$ by $R[x_1,\dots,x_n]$ to assume that $f$ is finite. We then apply \cite[Lemma~4.4.15]{heyer-mann-6ff} with $Q_i = \calO_X$ to deduce that a solid $A$-module $M \in \D_\solid(A)$ lies in $\Coh_R(A)$ if and only if $M$ is perfect as an $R$-module, as desired.
\end{proof}

\begin{rmks}
\begin{rmksenum}
    \item \label{rmk:Coh-in-other-geometries} Similar definitions and results work in other geometric settings, like formal, rigid and complex-analytic geometry, once a good theory of \enquote{solid} quasicoherent sheaves together with a 6-functor formalism is established. This is one of the beautiful observations made by Clausen--Scholze and collaborators and it is one of the first applications of condensed mathematics. Our formulation in terms of suave sheaves is conceptualizing this idea a bit further.

    \item Given a reductive group $G$ over $\Q_p$ and a nice $\Z[1/p]$-algebra $\Lambda$, the categorical local Langlands conjecture predicts an equivalence of $\D(\Bun_G, \Lambda)$ and $\Ind(\Coh^\qc(\Par_G))$ (up to a singular support condition), where $\Par_G \in \Stk_\Lambda^\lafp$ denotes the $L$-parameter stack and $\Coh^\qc(-)$ denotes the category of coherent sheaves with quasi-compact support. If $G$ is a torus then $\Par_G$ is smooth and hence $\Coh^\qc(\Par_G) = \Perf^\qc(\Par_G)$, which is expected by the Langlands correspondence in this case. However, for general rings $\Lambda$ (e.g. $\Lambda = \Z/\ell^2$) this identity only holds if one replaces $\Coh$ by $\Coh_\Lambda$.
\end{rmksenum}
\end{rmks}

As shown in \cref{rslt:coh-over-regular-noetherian-base}, in the case of a classical regular noetherian base $S$, the category $\Coh_S(X) = \Coh(X)$ admits a direct description independent of $S$. The next result shows that $\Coh(X)$ is well-behaved with respect to t-structures:

\begin{lem} \label{rslt:Coh-duality-is-t-bounded}
Fix a classical regular noetherian ring $\Lambda$ and a lafp scheme $X$ over $\Lambda$. Then the t-structure on $\QCoh(X)$ restricts to a t-structure on $\Coh(X) = \Coh_\Lambda(X)$ and the suave duality functor
\begin{align*}
    \SD\colon \Coh(X)^\op \isoto \Coh(X), \qquad \mathcal M \mapsto \intHom(\mathcal M, \omega_X)
\end{align*}
is uniformly t-bounded by a constant depending only on $X$.
\end{lem}
\begin{proof}
By \cref{rslt:coh-over-regular-noetherian-base} it is clear that the t-structure on $\QCoh(X)$ restricts to $\Coh(X)$. To prove the claim on $\SD$, we can reduce to the case that $X = \Spec A$ is affine. Pick a map $\Lambda' := \Lambda[x_1, \dots, x_n] \to A$ which is surjective on $\pi_0$. Then the suave duality on $A$ over $\Lambda'$ differs from the suave duality over $\Lambda$ by a shift by $n$ (since the dualizing complex of $\Lambda'$ over $\Lambda$ is $\Lambda'[n]$), so we can replace $\Lambda$ by $\Lambda'$ to assume that $A$ is finite over $\Lambda$. But in this case $\omega_A = \intHom_\Lambda(A, \Lambda)$ and hence $\SD = \intHom_\Lambda(-, \Lambda)$. We conclude by observing that $\Lambda$ has finite global dimension (see \cite[Lemma~00OE]{stacks-project}).
\end{proof}

\subsection{The 6-functor formalism for ind-coherent sheaves} \label{sec:alggeo.ICoh}

We finally come to the definition of ind-coherent sheaves and the associated 6-functor formalism. In \cref{sec:perf-and-coh} we have defined relatively coherent sheaves over a base. In the following we will mostly choose a classical regular noetherian base ring $\Lambda$, which simplifies some things (see \cref{rslt:coh-over-regular-noetherian-base,rslt:properties-of-noetherian-rings}): A $\Lambda$-algebra $A$ is almost finitely presented if and only if $\pi_0 A$ is finitely generated as a classical $\Lambda$-algebra and $\pi_n A$ is a finitely generated $\pi_0 A$-module for all $n \ge 0$. Moreover, $\Coh_\Lambda(A) = \Coh(A) \subseteq \D(A)$ is the category of bounded complexes of $A$-modules all of whose cohomologies are finitely generated over $\pi_0 A$.

\begin{defn}
Let $\Lambda$ be a classical regular noetherian ring. For an almost finitely presented $\Lambda$-algebra $A$ we set
\begin{align*}
    \ICoh(A) := \Ind(\Coh(A))
\end{align*}
and call it the category of \emph{ind-coherent} $A$-modules. Here $\Coh(A) := \Coh_\Lambda(A)$ is as discussed above.
\end{defn}

By definition of $\Ind$-categories, an object in $\ICoh(A)$ is a filtered diagram $(P_i)_i$ of coherent $A$-modules. By \cref{rslt:properties-of-Coh-and-Perf}, $\Coh(-)$ is stable under pullback along suave maps (i.e. maps of finite Tor dimension, see \cref{rslt:characterization-of-suave-maps}) and under pushforward along finite maps. In particular, these two operations immediately extend to functors on $\ICoh(-)$. However, if we want to extend $\ICoh(-)$ to stacks over $\Lambda$ we need to construct a pullback functor along \emph{all} maps $A \to B$ of almost finitely presented $\Lambda$-algebras, and in order to construct a 6-functor formalism for $\ICoh$ we also need to construct lower-$!$ functors for all maps $A \to B$. These constructions are subtle and occupy a large part of \cite{gaitsgory-rozenblyum-vol1}, where the additional assumption that $\Lambda$ is a field of characteristic $0$ is required. In the following we show how $\QCoh_\solid$ can be employed to provide an alternative and quick construction that works for general (regular noetherian) $\Lambda$.

The basic idea is to extract the six functors for $\ICoh$ out of the six functors on $\QCoh_\solid$. To achieve this, we first need to pass to a suitable subcategory of $\QCoh_\solid$, which generalizes the definition of pseudo-coherent modules (see \cref{def:pseudocoherent-modules}):

\begin{defn} \label{def:pseudo-compact-modules}
Let $A$ be a ring. A solid $A$-module $M \in \D_\solid(A)$ is called \emph{pseudo-compact} if $\Hom(A, -)$ preserves filtered colimits of uniformly left-bounded objects in $\D_\solid(A)$. We denote by
\begin{align*}
    \D_{\solid,\pc}(A) \subseteq \D_\solid(A)
\end{align*}
the full subcategory of pseudo-compact objects. We will often also denote $\QCoh_\solid(A) = \D_\solid(A)$ and then $\QCoh_{\solid,\pc}(A) := \D_{\solid,\pc}(A)$.
\end{defn}

The results in \cref{rslt:criteria-for-pseudocoherent-modules} and \cref{rslt:stability-of-perfect-and-pseudo-coherent} hold verbatim for $\D_\solid(A)$ in place of $\D(A)$, \enquote{pseudo-compact} in place of \enquote{pseudo-coherent}, and \enquote{compact} in place of \enquote{perfect}, with the same proofs. Additionally, in \cref{rslt:pseudo-coherent-implies-resolution-by-finite-free} \enquote{finite free} needs to be replaced by \enquote{finite direct sum of copies of $A_\solid[\N_\infty]$}. We will use this observation freely below. The category of pseudo-compact solid modules enjoys the following additional properties:

\begin{lem} \label{rslt:basic-properties-of-pseudo-compact-solid-modules}
Fix a ring $\Lambda$.
\begin{lemenum}
    \item \label{rslt:3ff-on-pseudo-compact-modules} The 6-functor formalism from \cref{rslt:6ff-on-QCoh-solid} on $\QCoh_\solid$ restricts to a 3-functor formalism for pseudo-compact modules on $\Alg_\Lambda^\afp$, i.e. to a functor of categories
    \begin{align*}
        \D_{\solid,\pc}\colon \Corr(\Alg_\Lambda^{\afp,\op})^\tensor \to \Cat,
    \end{align*}
    where $\Corr(-)^\tensor$ denotes the operad from \cite[Definition~2.3.2]{heyer-mann-6ff} and where we use \cite[Example~B.1.9]{heyer-mann-6ff} to relate the above functor to a map of operads.
    
    \item \label{rslt:coherent-implies-pseudo-compact} For every almost finitely presented $\Lambda$-algebra $A$ we have $\Coh_\Lambda(A) \subseteq \D_{\solid,\pc}(A)$.
\end{lemenum}
\end{lem}
\begin{proof}
We first prove (i), for which we need to show that $\D_{\solid,\pc}$ is stable under $\tensor$, $f^*$ and $f_!$. All of these functors preserve colimits in $\D_\solid$, so by the analogs of \cref{rslt:perfect-and-pseudo-coherent-stable-under-tensor-colim-base-change} and \cref{rslt:pseudo-coherent-implies-resolution-by-finite-free} we are reduced to showing that each of these functors applied to a compact generator is pseudo-compact. Fix a map $f\colon A \to B$ and a light profinite set $S$. Note that $A_\solid[S] \tensor_A^\solid A_\solid[S] = A_\solid[S \times S]$ is compact (and in particular pseudo-compact), handling the case of tensor product. Similarly $f^* A_\solid[S] = B_\solid[S]$ is compact. It remains to show that $f_! B_\solid[\N_\infty] \in \D_\solid(A)$ is pseudo-compact. We factor $f$ as $A \to A[x_1, \dots, x_n] \to B$ in order to reduce the claim to the following two cases: Either $f$ is smooth or $B$ is pseudo-coherent as an $A$-module. In the first case $f_!$ preserves compact objects because its right adjoint $f^! = f^* \tensor \omega_f$ preserves colimits. In the second case $f_! = f_*$ is just the forgetful functor and hence the claim follows from the analog of \cref{rslt:stability-of-perfect-and-pseudo-coherent-under-forget}.

We now prove (ii), so let $A$ and $P \in \Coh_\Lambda(A)$ be given. Pick a factorization $\Lambda \to \Lambda[x_1, \dots, x_n] \to A$ such that $A$ is pseudo-coherent as a $\Lambda[x_1, \dots, x_n]$-module. By \cref{rslt:characterization-of-suave-sheaves-in-affine-case} $P$ is pseudo-coherent (even perfect) as a $\Lambda[x_1, \dots, x_n]$-module, so by the analog of \cref{rslt:stability-of-perfect-and-pseudo-coherent-under-forget} we deduce that $P$ is pseudo-coherent (and in particular pseudo-compact) as a $B$-module.
\end{proof}

Using \cref{rslt:basic-properties-of-pseudo-compact-solid-modules} we can now construct the desired functoriality for $\ICoh$. To minimize confusion we will first formulate everything in terms of \emph{pro}-coherent sheaves and later pass to $\ICoh$ via taking opposite categories (see \cref{def:6ff-for-ICoh-on-rings} below).

\begin{defn}
Fix a classical regular noetherian ring $\Lambda$ and an almost finitely presented $\Lambda$-algebra $A$.
\begin{defenum}
    \item \label{def:PCoh-and-PQCoh-solid-pc-on-rings} We denote
    \begin{align*}
        \PCoh(A) := \Pro(\Coh(A)), \qquad \PQCoh_{\solid,\pc}(A) := \Pro(\D_{\solid,\pc}(A))
    \end{align*}
    the pro-categories of $\Coh(A)$ and $\D_{\solid,\pc}(A)$. Their objects are cofiltered diagrams of coherent and pseudo-compact solid $A$-modules, respectively.

    \item The embedding $\Coh(A) \injto \D_{\solid,\pc}(A)$ from \cref{rslt:coherent-implies-pseudo-compact} induces an embedding $\iota_A$ on pro-categories via right Kan extension together with a left adjoint $\eta_A$, i.e. we have the following diagram of adjoint functors
    \begin{align*}
        \eta_A\colon \PQCoh_{\solid,\pc}(A) \rightleftarrows \PCoh(A) \noloc \iota_A.
    \end{align*}
\end{defenum}
\end{defn}

Our goal is to employ $\eta_A$ and the 3-functor formalism on $\D_{\solid,\pc}(-)$ to construct a 3-functor formalism on $\PCoh(-)$. Let us first sketch the idea before coming to the formal result. Given a map $f\colon A \to B$ of afp $\Lambda$-algebras, we obtain a functor $f^*\colon \D_{\solid,\pc}(A) \to \D_{\solid,\pc}(B)$ and hence also an induced functor on the associated pro-categories. We can thus define
\begin{align*}
    f^* := \eta_B f^* \iota_A\colon \PCoh(A) \to \PCoh(B).
\end{align*}
We can similarly define $f_!$ and $\tensor$ on $\PCoh(-)$. Now the hard part is to show that these definitions are compatible with compositions. This requires a good understanding of the functor $\eta_B$ and this is also the main reason why we currently have to restrict to the case that $\Lambda$ is a classical regular noetherian ring. Let us start with some preparations.

\begin{lem}
Fix a classical regular noetherian ring $\Lambda$ and an afp $\Lambda$-algebra $A$.
\begin{lemenum}
    \item \label{rslt:t-structure-on-PCoh} The t-structure on $\Coh(A)$ induces a left-complete t-structure on $\PCoh(A)$ such that $M \in \PCoh(A)$ lies in $\PCoh^{\le0}(A)$ (resp. $\PCoh^{\ge0}(A))$) if and only if $M$ is isomorphic to a cofiltered diagram of objects in $\Coh^{\le0}(A)$ (resp. in $\Coh^{\ge0}(A)$).
    
    \item \label{rslt:computing-eta-on-generator} For every light profinite set $S = \varprojlim_i S_i$ we have
    \begin{align*}
        \eta_A(A_\solid[S]) = \varprojlim_{n\ge0} \varprojlim_i \tau^{\ge -n} A[S_i]
    \end{align*}
    in $\PCoh(A)$.

    \item \label{rslt:inherent-functors-on-PCoh-commute-with-geom-realization} Let $f\colon A \to B$ be a suave (resp. finite) map of afp $\Lambda$-algebras. Then $f^*\colon \Coh(A) \to \Coh(B)$ (resp. $f_!\colon \Coh(B) \to \Coh(A)$) induces a functor
    \begin{align*}
        f^*\colon \PCoh(A) \to \PCoh(B) \qquad \text{(resp. $f_!\colon \PCoh(B) \to \PCoh(A)$)}
    \end{align*}
    that preserves all small limits and all uniformly right-bounded geometric realizations.
\end{lemenum}
\end{lem}
\begin{proof}
Part (i) is straightforward, see e.g. \cite[Lemma~C.2.4.3]{lurie-SAG}. We now prove (ii). Note that there is a natural map from left to right, so the equivalence can be checked after applying $\Hom(-, M)$ for all $M \in \PCoh(A)$. By writing $M = (P_j)_j$ as a cofiltered diagram of coherent $A$-modules and using $\Hom(-, M) = \varprojlim_j \Hom(-, P_j)$, we immediately reduce to the case that $M = P$ is coherent. By the definition of homomorphisms in a pro-category we have
\begin{align*}
    \Hom(\varprojlim_{n\ge0} \varprojlim_i \tau^{\ge-n} A[S_i], P) = \varinjlim_i \varinjlim_{n\ge0} \Hom(\tau^{\ge-n} A[S_i], P) = \varinjlim_i \Hom(A[S_i], P) = P(S).
\end{align*}
Here the $\Hom$ in the middle two terms is the $\Hom$ in $\D(A)$, the second step follows easily from the fact that $P$ is left-bounded and the third step follows from $P$ being discrete. On the other hand we have
\begin{align*}
    \Hom(\eta_A(A_\solid[S]), P) = \Hom(A_\solid[S], P) = P(S).
\end{align*}
Here the $\Hom$ in the second term is the $\Hom$ in $\D_\solid(A)$ and the second step follows by definition of the free solid generator $A_\solid[S]$. This finishes the proof of (ii).

We now prove (iii). Note first that $f^*$ and $f_!$ are defined on $\Coh(-)$ by \cref{rslt:descent-for-Coh-and-Perf} and \cref{rslt:Coh-stable-under-proper-pushforward}. By abstract properties of pro-categories these functor extend to $\PCoh(-)$, where they preserves all small limits. It remains to show that they preserve uniformly right-bounded geometric realizations. By left-completeness of $\PCoh(-)$ this can be checked on cohomologies, where it follows easily from the fact that each cohomology group of a uniformly right-bounded geometric realization depends only on a truncated geometric realization (i.e. colimit over $\Delta_{\le n}$ for some $n \ge 0$) and that all exact functors preserve truncated geometric realizations because they preserve finite colimits. The only non-formal input for this argument is that $f^*$ and $f_!$ are right-bounded, which we check now. For $f_!$ this is clear, because for finite $f$ we have $f_! = f_*$ and this is just the forgetful functor and in particular t-exact. On the other hand, the functor $f^*$ is always right t-exact (without any condition on $f$) because it is left adjoint to the t-exact forgetful functor.
\end{proof}

\begin{cor} \label{rslt:eta-compatible-with-certain-pull-and-push}
Let $\Lambda$ be a classical regular noetherian ring and $f\colon A \to B$ a map of afp $\Lambda$-algebras.
\begin{corenum}
    \item \label{rslt:suave-pullback-compatibility-of-PCoh-and-Pro-pc} If $f$ is suave then the natural map $\eta_B f^* \isoto f^* \eta_A$ is an isomorphism.

    \item \label{rslt:finite-pushforward-compatibility-of-PCoh-and-Pro-pc} If $f$ is finite then the natural map $\eta_B f_! \isoto f_! \eta_B$ is an isomorphism.
\end{corenum}
\end{cor}
\begin{proof}
Note first that $\eta_A$ and $\eta_B$ preserve all small limits and colimits. The latter is clear because they are left adjoint functors. The former follows from the fact that their right adjoints $\iota_A$ and $\iota_B$ preserve cocompact objects by construction (this is the $\op$-version of the statement that a functor between compactly generated stable categories preserve all small colimits if its left adjoint preserves compact objects). In particular both sides of the claimed isomorphisms in (i) and (ii) preserve all small limits and hence they can be checked on objects $P \in \D_{\solid,\pc}(A)$ resp.\ $Q \in \D_{\solid,\pc}(B)$. By the analog of \cref{rslt:pseudo-coherent-implies-resolution-by-finite-free} we can write a shift of $P$ resp. $Q$ as a geometric realization of copies of $A_\solid[S]$ resp. $B_\solid[S]$ for light profinite sets $S = \varprojlim_i S_i$. This geometric realization is in particular a geometric realization in $\PQCoh_{\solid,\pc}(-)$ by general properties of pro-categories. Since $\eta_A$ and $\eta_B$ preserve geometric realizations and $f^*$ resp. $f_!$ preserve uniformly right-bounded geometric realizations by \cref{rslt:inherent-functors-on-PCoh-commute-with-geom-realization}, we reduce to the case $P = A_\solid[S]$ resp. $Q = B_\solid[S]$. Altogether we have reduced the claims to showing that the natural maps
\begin{align*}
    \eta_B f^* A_\solid[S] \isoto f^* \eta_A A_\solid[S], \qquad \eta_B f_! B_\solid[S] \isoto f_! \eta_B B_\solid[S]
\end{align*}
are isomorphisms in $\PCoh(B)$ resp.\ $\PCoh(A)$. In the suave case we compute using \cref{rslt:computing-eta-on-generator}
\begin{align*}
    f^* \eta_A A_\solid[S] = f^*(\varprojlim_n \varprojlim_i \tau^{\ge-n} A[S_i]) = \varprojlim_i (\varprojlim_n f^*\tau^{\ge-n} A)[S_i] = \varprojlim_i \varprojlim_n \tau^{\ge-n} B[S_i] = \eta_B f^* A_\solid[S],
\end{align*}
where in the third step we used the identity $\varprojlim_n f^* \tau^{\ge-n} A = \varprojlim_n \tau^{\ge-n} B$ that can easily be checked on cohomologies (using that $f$ is right t-exact as shown in the proof of \cref{rslt:inherent-functors-on-PCoh-commute-with-geom-realization}). In the finite case, the isomorphism for $f_!$ follows very similarly (note that $f_! = f_*$ is just a forgetful functor).
\end{proof}

With the above preparations at hand, we can finally implement the above strategy and construct a 3-functor formalism for $\PCoh$:

\begin{prop} \label{rslt:3ff-for-PCoh-on-rings}
Fix a classical regular noetherian ring $\Lambda$. The 3-functor formalism from \cref{rslt:3ff-on-pseudo-compact-modules} induces 3-functor formalisms for $\PQCoh_{\solid,\pc}$ and $\PCoh$ on $\Alg_\Lambda^{\afp,\op}$ together with a natural transformation of 3-functor formalisms
\begin{align*}
    \eta\colon \PQCoh_{\solid,\pc} \to \PCoh.
\end{align*}
\end{prop}
\begin{proof}
We first note that the 3-functor formalism on $\D_{\solid,\pc}$ induces a 3-functor formalism on $\PQCoh_{\solid,\pc}$, i.e. that $\Pro(-)$ preserves the structure of the 3-functor formalism. By passing to opposite categories we reduce to the following observation: If $\D$ is a stable 3-functor formalism then there is a stable 3-functor formalism $\Ind \comp \D$. The only critical point is the tensor product, for which we first employ \cite[Remark~4.8.1.9]{lurie-higher-algebra} (with $\cat K$ the collection of finite simplicial sets and $\operatorname{Idem}$) to lift the 3-functor formalism from the cartesian product to a certain tensor product on stable categories, and then use \cite[Proposition~5.3.2.11(3)]{lurie-higher-algebra} to obtain the desired (presentable) 3-functor formalism. This establishes the existence of the 3-functor formalism on $\PQCoh_{\solid,\pc}$.

We now show that we also obtain a 3-functor formalism on $\PCoh(-)$. As explained above, for a map $f\colon A \to B$ of afp $\Lambda$-algebras we define $f^* := \eta_B f^* \iota_A$ and $f_! := \eta_A f_! \iota_B$, where the left-hand side denotes the desired functor on $\PCoh$; in fact we are forced to define the functors this way if we want $\eta$ to be a morphism of 3-functor formalisms. In order to show that this definition works, we apply \cite[Lemma~A.2.24]{heyer-mann-6ff} to the functor $\PQCoh_{\solid,\pc}\colon \Corr(\Alg_\Lambda^{\afp,\op})^\tensor \to \Cat$ and the full subcategories $\PCoh(A) \subseteq \PQCoh_{\solid,\pc}(A)$. More precisely, an object in $\Corr(\Alg_\Lambda^{\afp,\op})^\tensor$ is a tuple $(A_1, \dots, A_n)$ which gets sent to $\prod_{i=1}^n \PQCoh_{\solid,\pc}(A_i)$ and we pick the full subcategory $\prod_{i=1}^n \PCoh(A_i)$. Now \cite[Lemma~A.2.24]{heyer-mann-6ff} reduces the claim to checking the following condition: For every map $\alpha\colon A_\bullet \to B_\bullet$ in $\Corr(\Alg_\Lambda^{\afp,\op})^\tensor$ with induced functor $\alpha_\sharp$ on $\PQCoh_{\solid,\pc}(-)$ the natural map
\begin{align*}
    \eta_{B_\bullet} \alpha_\sharp \isoto (\eta_{A_\bullet} \alpha_\sharp \iota_{B_\bullet}) \eta_{B_\bullet}
\end{align*}
is an isomorphism. This property is stable under compositions in $\alpha$ and can be checked componentwise in a product of categories, which reduces the claim to showing that for every map $f\colon A \to B$ of afp $\Lambda$-algebras the following natural maps are isomorphisms:
\begin{enumerate}[(a)]
    \item $\eta_B f^* \isoto (\eta_B f^* \iota_A) \eta_A$,
    \item $\eta_A f_! \isoto (\eta_A f_! \iota_B) \eta_B$,
    \item $\eta_{A \tensor_\Lambda A} (- \boxtimes -) \isoto \eta_{A \tensor_\Lambda A}(\iota_A\eta_A \boxtimes \iota_A\eta_A)$.
\end{enumerate}
Let us start with (a). By using the usual factorization $A \to A[x_1, \dots, x_n] \to B$ we can reduce to the cases that $f$ is suave or finite. In the suave case we see that $f^* \iota_A = \iota_A f^*$, where the $f^*$ on the right-hand side is the one from \cref{rslt:inherent-functors-on-PCoh-commute-with-geom-realization}. Thus in this case the claim in (a) follows from \cref{rslt:suave-pullback-compatibility-of-PCoh-and-Pro-pc}. Now assume that $f$ is finite. Then $f^*$ is left adjoint to $f_!$ on $\PQCoh_{\solid,\pc}(-)$ and hence $\eta_B f^*$ is left adjoint to $f_! \iota_B = \iota_A f_!$, where on the right-hand side we use the functor from \cref{rslt:inherent-functors-on-PCoh-commute-with-geom-realization}. By restriction we see that $\eta_B f^* \iota_A$ is left adjoint to $f_!$. Thus by passing to right adjoints in (a) the claim reduces to the obvious identity $f_! \iota_B = \iota_A f_!$.

We now prove (b). The proof is analogous to (a), so we only sketch it. We first reduce to the case that either $B = A[x]$ or $A \to B$ is finite. In the latter case we finish by \cref{rslt:finite-pushforward-compatibility-of-PCoh-and-Pro-pc} and in the former case we use that $f_!$ is left adjoint to $f^*[1]$ in order to finish by passing to right adjoints.

It remains to prove (c). Note first that $\Coh(-)$ is stable under $\boxtimes$: One can either see this by explicit computation or using 2-categorical magic. The latter approach shows that suave objects over some fixed object are stable under $\boxtimes$ in any 3-functor formalism and the reason is that $\boxtimes$ is obtained from the symmetric monoidal structure $\cat K_{\D} \times \cat K_{\D} \to \cat K_{\D}$ encoded in \cite[Theorem~4.2.4]{heyer-mann-6ff}; now we simply observe that 2-functors preserve adjunctions. In any case, $\boxtimes$ passes to pro-categories to induce a functor $\boxtimes\colon \PCoh(A)^2 = \Pro(\Coh(A)^2) \to \Pro(A \tensor_\Lambda A)$ which is compatible via $\iota_A$ with $\boxtimes$ on $\PQCoh_{\solid,\pc}(A)$ by construction. The rest of the argument is as in the suave case of (a) using that \cref{rslt:inherent-functors-on-PCoh-commute-with-geom-realization} and \cref{rslt:suave-pullback-compatibility-of-PCoh-and-Pro-pc} hold also for $\boxtimes$ in place of $f^*$, with the same proofs.
\end{proof}

Let us now switch from $\PCoh$ to $\ICoh$ by passing to opposite categories. Since $(-)^\op$ is an endofunctor of $\Cat$ we retain the 3-functor formalism from \cref{rslt:3ff-for-PCoh-on-rings}. Moreover, as all transition functors preserve colimits and $\ICoh(-)$ is presentable, the 3-functor formalism even becomes a 6-functor formalism (cf. \cite[Lemma~3.2.5]{heyer-mann-6ff}).

\begin{defn} \label{def:6ff-for-ICoh-on-rings}
Fix a classical regular noetherian ring $\Lambda$. For every afp $\Lambda$-algebra $A$ we consider the identification
\begin{align*}
    \ICoh(A) = \Ind(\Coh(A)) = \Ind(\Coh(A)^\op) = \Pro(\Coh(A))^\op.
\end{align*}
In the second identity we used the suave duality functor
\begin{align*}
    \Coh(A)^\op \isoto \Coh(A), \qquad P \mapsto \intHom(P, \omega_A)
\end{align*}
in $\QCoh_\solid$, where $\omega_A = f^! \Lambda$ for $f\colon \Lambda \to A$ the structure map (cf. \cref{rslt:colim-stability-and-self-duality-for-Perf-and-Coh}). As explained above, the 3-functor formalism from \cref{rslt:3ff-for-PCoh-on-rings} induces a presentable stable 6-functor formalism on $\ICoh(-)$ together with a natural transformation of 6-functor formalisms $\PQCoh^\op_{\solid,\pc} \to \ICoh$. We will denote the six functors on $\ICoh$ in the usual 6-functor notation (this deviates from the notation in \cite{gaitsgory-rozenblyum-vol1}!).
\end{defn}

The transition from $\PCoh$ to $\ICoh$ creates some strange behavior that one needs to be careful with. Firstly, the passage to opposite categories flips the usual notions in the 6-functor formalism: Smooth maps behave like proper maps in $\ICoh$ while proper maps behave like smooth ones. Secondly, the suave duality functor used in the transition from $\PCoh$ to $\ICoh$ changes the tensor product. For example, if $A$ is truncated (i.e. suave over $\Lambda$) then the tensor unit in $\PCoh(A)$ is just $A$, while in $\ICoh(A)$ it is $\omega_A$.

We will postpone a more systematic study of the intricacies of $\ICoh$ to the next subsection and continue now with the construction of all the relevant structures. A crucial observation in \cite{gaitsgory-rozenblyum-vol1} is that $\ICoh$ has a tight relationship with $\QCoh$. Our construction provides a quick proof of that fact that avoids the construction of complicated 2-categories:

\begin{lem} \label{rslt:gamma-for-QCoh-and-ICoh-on-rings}
Fix a classical regular noetherian ring $\Lambda$.
\begin{lemenum}
    \item \label{rslt:construct-gamma-from-QCoh-to-ICoh-on-rings} There is a natural transformation $\gamma\colon \QCoh \to \ICoh$ of functors $\Alg_\Lambda^\afp \to \CMon$. In other words, for every afp $\Lambda$-algebra $A$ there is a symmetric monoidal functor
    \begin{align*}
        \gamma_A\colon \QCoh(A) \to \ICoh(A)
    \end{align*}
    which is compatible with pullbacks on both sides.

    \item \label{rslt:gamma-compatible-with-suave-pushforward-on-rings} If $f\colon A \to B$ is a suave map of afp $\Lambda$-algebras then the natural map
    \begin{align*}
        \gamma_A f_* \isoto f_* \gamma_B
    \end{align*}
    is an isomorphism of functors $\QCoh(B) \to \ICoh(A)$.
\end{lemenum}
\end{lem}
\begin{proof}
We first prove (i). Note that for every afp $\Lambda$-algebra $A$, every perfect $A$-module is in particular a pseudo-compact solid $A$-module, i.e. we have an embedding $\Perf(A) \subseteq \D_{\solid,\pc}(A)$ of symmetric monoidal categories. Passing to pro-categories produces a natural transformation of functors $\Pro(\Perf(-)) \to \PQCoh_{\solid,\pc}$ and composing this with $\eta$ yields a natural transformation $\Pro(\Perf(-)) \to \PCoh$. We now pass to opposite categories in order to obtain the natural transformation $\Ind(\Perf(-)^\op) \to \ICoh$ of functors $\Alg_\Lambda^{\afp,\op} \to \CMon$. Thus it only remains to observe that the duality functor $(-)^\vee\colon \Perf(A)^\op \isoto \Perf(A)$ is symmetric monoidal and induces an equivalence of functors $\Ind(\Perf(-)^\op) = \QCoh$. This is easy to check on objects, as it just amounts to the identities $(P \tensor Q)^\vee = P^\vee \tensor Q^\vee$ and $f(P)^\vee = f(P^\vee)$ for perfect $P$ and $Q$; we provide a formal construction of the necessary higher homotopies in the appendix, see \cref{rslt:functoriality-of-dualizable-dual}.

We now prove (ii), so let the suave map $f\colon A \to B$ be given. Since all functors in the claimed isomorphism commute with colimits, it is enough to check the isomorphism on $B$, i.e. we need to see that the natural map $\gamma_A f_* B \isoto f_* \gamma_B B$ is an isomorphism. To avoid suave duality trickery we work with the pro-perspective, i.e. pass to opposite categories and view $\gamma$ as a functor $\Pro(\Perf(-)) \to \PCoh(-)$ (note that $f_*$ now becomes a \emph{left} adjoint). By \cref{rslt:characterization-of-suave-maps} we can factor $f$ as $A \to A[x_1, \dots, x_n] \to B$ such that $B$ is perfect as an $A[x_1, \dots, x_n]$-module, which reduces the claim to the cases that either $B = A[x]$ or $B$ is perfect as an $A$-module.

Suppose first that $B = A[x]$. Then $f_* B = \prod_\N A$ and hence $\gamma_A f_* B = \prod_\N \gamma_A A$. We now compute, for every coherent $A$-module $P$,
\begin{align*}
    \Hom(f_* \gamma_B B, P) = \Hom(f_* \eta_B B, P) = \Hom(\eta_B B, f^* P) = \Hom(B, f^* P),
\end{align*}
where the last $\Hom$ is taken in $\PQCoh_{\solid,\pc}(B)$ or equivalently in $\QCoh(B)$. But then it computes as
\begin{align*}
    \Hom_B(B, f^* P) = \Hom_A(A, P \tensor_A A[x]) = \bigoplus_\N \Hom_A(A, P) = \Hom(\prod_\N \gamma_A A, P),
\end{align*}
where the last $\Hom$ is taken in $\PCoh(A)$. This shows that $f_* \gamma_B B = \prod_\N \gamma_A A$, as desired.

Now suppose that $B$ is perfect as an $A$-module. Then $\gamma_A f_* B = \eta_A B$ (here we use that $B$ is perfect as an $A$-module) and $f_* \gamma_B B = f_* \eta_B B = \eta_A f_* B = \eta_A B$, as desired. In the second computation we used the fact that $f_* = f_!$ on $\PQCoh_{\solid,\pc}$ and $\PCoh$ and that $\eta$ commutes with $f_!$ because it is a natural transformation of 3-functor formalisms.
\end{proof}

The functor $\gamma_A$ provides a compatibility of $\QCoh$ and $\ICoh$ with respect to the pullback and pushforward functors. We now construct a functor $\Psi_A$ that is compatible with lower-$!$ functors. In earlier constructions of ind-coherent sheaves, $\Psi_A$ is usually the first functor that is constructed (see e.g. \cite[Proposition~4.2.2.3]{gaitsgory-rozenblyum-vol1}). In our approach $\Psi_A$ does not play an important role and in fact is not really needed in the paper; we mainly construct it for the sake of completeness.

\begin{lem} \label{rslt:psi-for-ICoh-and-QCoh-on-rings}
Fix a classical regular noetherian ring $\Lambda$.
\begin{lemenum}
    \item \label{rslt:construct-psi-for-ICoh-and-QCoh-on-rings} There is a natural transformation $\Psi\colon \ICoh \to \QCoh$ of functors $\Alg_\Lambda^{\afp,\op} \to \Cat$. In other words, for every afp $\Lambda$-algebra $A$ there is a functor
    \begin{align*}
        \Psi_A\colon \ICoh(A) \to \QCoh(A)
    \end{align*}
    which is compatible with lower-! functors on the left and pushforward functors on the right. Explicitly, $\Psi_A$ is obtained as the ind-extension of the natural inclusion $\Coh(A) \subseteq \QCoh(A)$.

    \item \label{rslt:pullback-compatibility-of-psi-on-rings} If $f\colon A \to B$ is a suave map of afp $\Lambda$-algebras then $f_!\colon \ICoh(B) \to \ICoh(A)$ admits a left adjoint $f^\natural$ and the natural map
    \begin{align*}
        f^* \Psi_B \isoto \Psi_A f^\natural
    \end{align*}
    is an isomorphism of functors $\ICoh(A) \to \QCoh(B)$.

    \item \label{rslt:xi-for-QCoh-and-ICoh-on-rings} If $A$ is a suave $\Lambda$-algebra then $\Psi_A$ admits a fully faithful left adjoint
    \begin{align*}
        \Xi_A\colon \QCoh(A) \injto \ICoh(A)
    \end{align*}
    which is induced from the inclusion $\Perf(A) \subseteq \Coh(A)$ via passing to $\Ind$-categories. If $f\colon A \to B$ is a suave map of $\Lambda$-algebras then the natural map
    \begin{align*}
        \Xi_A f_* \isoto f_! \Xi_B
    \end{align*}
    is an isomorphism of functors $\QCoh(B) \to \ICoh(A)$.
\end{lemenum}
\end{lem}
\begin{proof}
To prove (i), we consider the natural transformation $\QCoh_{\solid,!}(A)^\op \to \QCoh_{\solid,*}(A)$, $M \mapsto \intHom(M, \omega_A)$ of functors $\Alg_\Lambda^{\afp,\op} \to \Cat$ from \cref{rslt:functoriality-of-suave-dual}. We compose this functor with the forgetful functor $\QCoh_\solid \to \QCoh$ that sends a solid $A$-module $M$ to its underlying discrete $A$-module $M(*)$. By restricting to pseudo-compact objects in the source and then passing to $\Pro$-categories, we obtain a natural transformation
\begin{align*}
    \Psi'\colon \PQCoh_{\solid,\pc}^\op \to \QCoh
\end{align*}
of functors $\Alg_\Lambda^{\afp,\op} \to \Cat$, where the transition maps in the source are given by lower-! functors and the transition maps in the target are pushforwards (i.e. forgetful functors). We claim that $\Psi'$ factors as $\PQCoh_{\solid,\pc}^\op \to \PCoh^\op \to \QCoh$, where the first map is $\eta$ and the second is the desired natural transformation $\Psi$. Since $\eta$ is adjoint to the fully faithful embedding $\PCoh^\op \injto \PQCoh_{\solid,\pc}^\op$, one checks that the desired factorization is unique if it exists, and it exists as soon as for any $A \in \Alg_\Lambda^\afp$ and any map $f\colon M \to N$ in $\PQCoh_{\solid,\pc}^\op$ such that $\eta(f)$ is an isomorphism, also $\Psi'(f)$ is an isomorphism (this can be proved similarly to \cite[Lemma~A.2.24]{heyer-mann-6ff} by passing to the cocartesian fibration of the functor $[1] \times \Alg_\Lambda^{\afp,\op} \to \Cat$ given by $\Psi'$ and then looking at the appropriate full subcategory of that fibration). To prove the required compatibility, suppose we are given $f\colon M \to N$ with $\eta(f)$ an isomorphism. To prove that $\Psi'(f)$ is an isomorphism, it is enough to do so after applying the conservative forgetful functor $\QCoh(A) \to \QCoh(\Lambda)$. Since $\eta$ and $\Psi'$ are natural transformations, this easily reduces the claim to the case $A = \Lambda$. Now $\Psi'_\Lambda = \Hom(-, \Lambda)$ and since $\eta(f)$ is an isomorphism and $\Hom(M, \Lambda) = \Hom(\eta(M), \Lambda)$, we deduce that $\Psi'(f)$ is indeed an isomorphism.

We have constructed a natural transformation $\Psi\colon \PCoh^\op \to \QCoh$. It remains to verify the explicit description at the end of (i) after using the identification $\PCoh(A)^\op = \ICoh(A)$. First note that $\Psi$ preserves all small colimits by construction, so we only need to check that its restriction to $\Coh(A)$ is the canonical embedding $\Coh(A) \injto \QCoh(A)$. By tracing through the construction and the identification of $\PCoh^\op$ with $\ICoh$, we see that $\Psi$ is given on $\Coh(A)$ by applying $\intHom(-, \omega_A)$; but this is the suave duality functor, hence its square is indeed the identity.

We now prove (ii), so let $f\colon A \to B$ be suave. Then $f$ is $\QCoh_\solid$-suave, hence $\PQCoh_{\solid,\pc}$-suave and since $\eta\colon \PQCoh_{\solid,\pc} \to \PCoh$ is a morphism of 3-functor formalisms it follows that $f$ is $\PCoh$-suave and hence $\ICoh$-prim (cf. the proof of descent for $\ICoh$ in \cref{rslt:fppf-descent-for-all-sheaf-theories}). In particular $f_!$ admits a left adjoint $f^\natural$ on $\ICoh$. It is now enough to prove desired identity on $\Coh(A) \subseteq \ICoh(A)$, because all functors in the claim commute with colimits. By passing to left adjoints in the identity $\eta_A f_! = f_! \eta_B$ we see that there is a natural isomorphism $f^\natural \iota_A = \iota_B f^\natural$ of functors $\ICoh(A) \to \PQCoh_{\solid,\pc}^\op(B)$. Now the functor $f^\natural$ on $\PQCoh_{\solid,\pc}^\op(B)$ is just the functor $f^!$ on $\PQCoh_{\solid,\pc}^\op(B)$. Since in the identification $\ICoh(A) = \PCoh(A)^\op$ we implicitly used suave duality, which swaps $f^!$ and $f^*$, we finally deduce that $f^\natural\colon \ICoh(A) \to \ICoh(B)$ restricts to the usual pullback $\Coh(A) \to \Coh(B)$, as desired.

We now prove (iii), so let $A$ be a suave $\Lambda$-algebra. Then $\Perf(A) \subseteq \Coh(A)$ (e.g. by \cref{rslt:Coh-S-and-smoothness-diagonal}), so we get a fully faithful functor $\Xi_A\colon \QCoh(A) \injto \ICoh(A)$ by applying $\Ind$. There is a natural isomorphism $\id \isoto \Psi_A \Xi_A$, as both sides commute with colimits and clearly agree on perfect $A$-modules. This induces a map $R \to \Psi_A$ from the right adjoint $R$ of $\Xi_A$ to $\Psi_A$. To show that this map is an isomorphism, first note that $R$ commutes with colimits because $\Xi_A$ preserves compact objects; hence it is enough to show that the map $R(M) \to \Psi_A(M)$ is an isomorphism for $M \in \Coh(A)$. This can be checked after applying $\Hom(P, -)$ for every perfect $A$-module $P$, but then both sides evaluate to $\Hom(P, M)$ in $\QCoh(A)$. This proves the existence and description of the left adjoint $\Xi_A$ of $\Psi_A$.

Now let $f\colon A \to B$ be suave. To show that the map $\Xi_A f_* \to f_! \Xi_B$ is an isomorphism, it is enough to show this when evaluated at $B$, because $B$ generate $\QCoh(B)$ under colimits. As in the proof of \cref{rslt:gamma-compatible-with-suave-pushforward-on-rings} we can reduce to the case that either $B$ is perfect as an $A$-module or $B = A[x]$. In the first case we see that $f_!\colon \ICoh(B) \to \ICoh(A)$ restricts to the forgetful functor $\Coh(B) \to \Coh(A)$, hence the claim follows immediately from the description of $\Xi_A$ and $\Xi_B$ (note that the perfectness of $B$ is crucial in order to guarantee that $f_*\colon \Perf(B) \to \QCoh(A)$ lands in $\Perf(A)$). It remains to deal with the case $B = A[x]$. Then $\Xi_A f_* B = \bigoplus_\N \Xi_A A$. One computes $f_! \Xi_A B = \bigoplus_\N \Xi_A A$ by going through the definitions: In terms of $\PCoh$ the statement reads $f_! \omega_B = \prod_\N \omega_A$. To verify this, it is enough to show that $\Hom(f_! \omega_B, P) = \bigoplus_\N \Hom(\omega_A, P)$ in $\PCoh(A)$, where $P \in \Coh(A)$ is arbitrary. Now the left-hand side can be computed in $\PQCoh_{\solid,\pc}$ and hence in $\QCoh_\solid$, and we obtain (the following notation is all in terms of $\QCoh_\solid$)
\begin{align*}
    &\Hom(f_! \omega_B, P) = \Hom(\omega_B, f^! P) = \Hom(f^* \omega_A \tensor_B \omega_f, f^* P \tensor_B \omega_f) = \Hom(f^* \omega_A, f^* P)=\\&\qquad= \Hom(\omega_A, f_* f^* P) = \Hom(\omega_A, \bigoplus_\N P) = \bigoplus_\N \Hom(\omega_A, P),
\end{align*}
as desired. Here the last step uses that $\omega_A$ is (pseudo-)coherent and $\bigoplus_\N P$ is a uniformly left-bounded direct sum.
\end{proof}

Our next goal is to extend the above constructions to stacks. Most of this is formal by \cite[Theorem~3.4.11]{heyer-mann-6ff}, but we need to verify fppf descent of the involved sheaf theories.


\begin{prop} \label{rslt:fppf-descent-for-all-sheaf-theories}
Fix a ring $\Lambda$. All of the following categories satisfy fppf descent on $\Alg_\Lambda$, i.e. they satisfy the limit property in \cref{rslt:fppf-descent-for-QCoh-solid}:
\begin{align*}
    \Perf, \qquad \QCoh, \qquad \QCoh_\solid, \qquad \QCoh_{\solid,\pc}, \qquad \PQCoh_{\solid,\pc}.
\end{align*}
Moreover, if $\Lambda$ is a classical regular noetherian ring, then additionally the following categories satisfy fppf descent on $\Alg_\Lambda^\afp$:
\begin{align*}
    \Coh, \qquad \PCoh, \qquad \ICoh.
\end{align*}
\end{prop}
\begin{proof}
For $\QCoh$ this was shown in \cref{rslt:fpqc-descent-for-QCoh}, for $\QCoh_\solid$ it was shown in \cref{rslt:fppf-descent-for-QCoh-solid} and for $\Perf$ and $\Coh$ it was shown in \cref{rslt:descent-for-Coh-and-Perf}. For $\PCoh$ it follows from the $\ICoh$-version because $(-)^\op$ preserves limits. We are thus left with showing fppf descent for $\QCoh_{\solid,\pc}$, $\PQCoh_{\solid,\pc}$ and $\ICoh$.

We start with the first two categories and fix an fppf map $f\colon A \to B$ of $\Lambda$-algebras. Denote by $f_\bullet\colon \Spec B^\bullet \to \Spec A$ the Čech nerve of $f$. We now make the following key observation, which is inspired by \cite[Proposition~4.6]{mikami-fppf-descent2}:
\begin{itemize}
    \item[($*$)] The identity functor on $\D_\solid(A)$ can be obtained using finite limits, colimits and retracts from the functors $f_{n!} f_n^!$ in $\Fun(\D_\solid(A), \D_\solid(A))$.
\end{itemize}
To prove ($*$), we first observe that $f_{n!} f_n^! = f_{n!} (f_n^* \tensor \omega_{f_n}) = - \tensor f_{n!} \omega_{f_n}$, so the claim reduces to showing that $A$ can be obtained using finite limits, colimits and retracts from the objects $f_{n!} \omega_{f_n}$. Note that $f$ is a suave cover, so by \cite[Lemma~4.7.1]{heyer-mann-6ff} it is a $!$-cover, which implies that $M = \varinjlim_{n\in\Delta} f_{n!} f_n^! M$ for every $M \in \D_\solid(A)$. Plugging in $M = A$ shows that the natural map
\begin{align*}
    \varinjlim_{n\in\Delta} f_{n!} \omega_{f_n} \isoto A
\end{align*}
is an isomorphism in $\D_\solid(A)$. We now claim that $f_{n!} \omega_{f_n} \in \D_\solid^{\le 1}(A)$ for all $n$. For this we have to show that for all $N \in \D^{>1}_\solid(A)$ we have $\Hom(f_{n!} \omega_{f_n}, N) = 0$. We compute
\begin{align*}
    \Hom(f_{n!} \omega_{f_n}, N) = \Hom(\omega_{f_n}, f_n^! N) = \Hom(B^n, f_n^* N),
\end{align*}
where in the last step we used the identity $f_n^* = \intHom(\omega_{f_n}, f_n^!)$ from \cite[Corollary~4.5.11(i)]{heyer-mann-6ff}. Now the claimed bound follows from the fact that $f_n$ has Tor dimension $\le 1$ by \cite[Lemma~4.10]{mikami-fppf-descent2}. For every integer $m \ge 0$ let us denote by $F_m$ the truncated geometric realization of $(f_{\bullet!} \omega_{f_\bullet})$, i.e. the colimit over $\Delta_{\le m}$. Then by the above isomorphism of the full geometric realization with $A$ we obtain a map $F_2 \to A$ which by \cite[Proposition~1.2.4.5]{lurie-higher-algebra} induces an isomorphism $\tau^{\ge0} F_2 \isoto \tau^{\ge0} A$. Let $R := \cofib(F_2 \to A)$, which thus lies in $\D_\solid^{<0}(A)$. But then the map $A \to R$ must be zero (as it corresponds to an element in $\pi_0 R = 0$), thus $A$ is a retract of $F_2$.


We now show descent for $\QCoh_{\solid,\pc}$. Since $\D_\solid(-)$ descends along $f$, we only need to show the following: Given $M \in \D_\solid(A)$ then $M$ is pseudo-compact if and only if $f_n^* M$ is pseudo-compact for all $n \ge 0$. The \enquote{only if} direction is clear, so it remains to prove the \enquote{if} direction. Thus assume that $f_n^* M$ is pseudo-compact for all $n$; we need to show that $M$ is pseudo-compact. Let us first observe that $f_{n!} f_n^! M$ is pseudo-compact for all $n$. Indeed, since $f_n$ is suave (see \cref{rslt:characterization-of-suave-maps}) we have $f_n^! = f_n^* \tensor \omega_{f_n}$, so by \cref{rslt:3ff-on-pseudo-compact-modules} we are reduced to showing that $\omega_{f_n} := f_n^! A \in \D_\solid(B^n)$ is pseudo-compact. This can be seen by factoring $f_n$ into $g_n\colon A \to A[x_1, \dots, x_m]$ and $h_n\colon A[x_1, \dots, x_m] \to B^n$ such that $B^n$ is perfect over $A[x_1, \dots, x_m]$ (see \cref{rslt:characterization-of-suave-maps}). Then
\begin{align*}
    \omega_{f_n} = h_n^! g_n^! A = h_n^! A[x_1, \dots, x_m][m] = \intHom_{A[x_1, \dots, x_m]}(B^n, A[x_1, \dots, x_m][m]).
\end{align*}
From this description it follows easily that $\omega_{f_n}$ is pseudo-coherent (even perfect) as an $A[x_1, \dots, x_m]$-module, hence by the analog of \cref{rslt:stability-of-perfect-and-pseudo-coherent-under-forget} it is also pseudo-coherent (in particular pseudo-compact) as a $B^n$-module. This finishes the proof that all $f_{n!} f_n^! M$ are pseudo-compact. But then ($*$) immediately implies that $M$ is pseudo-compact, as desired.

We next show descent for $\PQCoh_{\solid,\pc}$. Let us denote $\PQCoh_\solid(-) := \Pro(\QCoh_\solid(-))$ (we ignore the set-theoretic issues involved in this definition; they can be fixed by replacing this category with the filtered colimit of $\Pro(\QCoh_\solid(A)^\kappa)$ for increasing regular cardinals $\kappa$). We now show that $\PQCoh_\solid$ descends along $f$. For all $n \ge 0$ we have $\intHom(f_{n!} \omega_{f_n}, -) = f_{n*} \intHom(\omega_{f_n}, f_n^!) = f_{n*} f_n^*$, hence ($*$) implies that the identity endofunctor of $\D_\solid(A)$ is obtained using finite limits, colimits and retracts from the functors $f_{n*} f_n^*$. In other words, $f$ induces a \emph{descendable} map of analytic rings (this is also deduced in \cite[Theorem~4.15]{mikami-fppf-descent2}). By functoriality of $\Pro(-)$ the same is then true for $\PQCoh_\solid$. But this implies the desired descent (as in \cite[Proposition~2.6.3]{mann-p-adic-6-functors}). With the descent of $\PQCoh_\solid$ at hand, we can show the descent of $\PQCoh_{\solid,\pc}$ in the same way as the descent for $\QCoh_{\solid,\pc}$ (observing that ($*$) also holds in $\PQCoh_\solid$ by functoriality of $\Pro(-)$).

It remains to prove descent for $\ICoh$, so assume now that $\Lambda$ is classical regular noetherian and that $A$ and $B$ are afp over $\Lambda$. Recall that suaveness is transferred along any map of 3-functor formalisms by the functoriality of the category of kernels (see \cite[Proposition~4.2.1(i)]{heyer-mann-6ff}) and the fact that adjunctions are transferred along any 2-functor. Thus, since $f$ is suave (by \cref{rslt:characterization-of-suave-maps}), i.e. $\QCoh_\solid$-suave, it is in particular $\PQCoh_\solid$-suave. As the suave dualizing complex $\omega_f$ is pseudo-compact by the above argument, we deduce that $f$ is $\QCoh_{\solid,\pc}$- and in particular $\PQCoh_{\solid,\pc}$-suave. Thus $f$ is $\PQCoh^\op_{\solid,\pc}$-prim and via $\eta$ we see that $f$ is $\ICoh$-prim. By \cite[Lemma~4.7.4(i)]{heyer-mann-6ff} the descent for $\ICoh$ follows if we can show that $f_* 1$ is descendable in $\ICoh(A)$. But by \cref{rslt:gamma-for-QCoh-and-ICoh-on-rings} we have $f_* 1 = \gamma_A f_* 1$. Since $\gamma_A$ is symmetric monoidal, it preserves descendable algebras, so it is enough to show that $f_* 1 = B$ is descendable in $\QCoh(A)$. But this is a classical result by Mathew (see \cite[Proposition~3.31]{akhil-galois-group-of-stable-homotopy}) or alternatively follows from the above discussed stronger version for solid analytic rings.
\end{proof}

We can finally construct the 6-functor formalism for $\ICoh$ on stacks. For studying $\ICoh$ it is very useful to have the other 6-functor formalisms used above, so the following result includes them as well.

\begin{thm} \label{rslt:6ff-for-ICoh}
Fix a classical regular noetherian ring $\Lambda$. There is a class $E$ of maps in $\Stk_\Lambda^\lafp$ and a presentable 6-functor formalism
\begin{align*}
    \ICoh\colon \Corr(\Stk_\Lambda^\lafp, E) \to \Cat
\end{align*}
with the following properties:
\begin{thmenum}
    \item The underlying functor $\ICoh^*\colon \Stk_\Lambda^{\lafp,\op} \to \Cat$ is a sheaf and coincides with the functor from \cref{def:6ff-for-ICoh-on-rings} on affine schemes.

    \item \label{rslt:ICoh-on-lafp-schemes} $E$ contains every qcqs map $f\colon Y \to X$ that is representable in lafp schemes over $\Lambda$. Moreover, we have the following explicit description:
    \begin{enumerate}[(a)]
        \item If $f$ is proper then it is $\ICoh$-smooth with $\ICoh$-dualizing complex isomorphic to $1$, so we have an isomorphism $f^! \isom f^*$ on $\ICoh$.
        
        \item If $f$ is suave then it is $\ICoh$-prim, so $f_* = f_!(\delta_f^{\ICoh} \tensor -)$ for the $\ICoh$-codualizing complex $\delta_f^{\ICoh} \in \ICoh(Y)$.
        
        \item If $f$ is cohomologically smooth then $\delta_f^{\ICoh} = \gamma(\omega_f^\vee)$, where $\omega_f \in \QCoh(Y)$ is the $\QCoh_\solid$-dualizing complex. In particular we have $f_* = f_!(\gamma(\omega_f^\vee) \tensor -)$ on $\ICoh$, and if $f$ is étale then $f_* = f_!$ on $\ICoh$.
    \end{enumerate}

    \item \label{rslt:stability-of-ICoh-fine-maps} $E$ has the following stability properties:
    \begin{enumerate}[(a)]
        \item $E$ is stable under composition and pullback. Moreover, if $f\colon Y \to X$ and $g\colon Z \to Y$ are maps in $\Stk_\Lambda$ such that $f, fg \in E$, then $g \in E$.
        
        \item $E$ is local on the target, i.e. containment in $E$ can be checked after pullback to any cover of the target.

        \item $E$ is $!$-local on source and target, i.e. containment in $E$ can be checked after passing to any universal $\ICoh$-$!$-cover of the source and target.
    \end{enumerate}

    \item \label{rslt:gamma-from-QCoh-to-ICoh-on-stacks} There is a natural transformation $\gamma\colon \QCoh \to \ICoh$ of functors $\Stk_\Lambda^{\lafp,\op} \to \CMon$, i.e. for every lafp stack $X$ over $\Lambda$ there is a symmetric monoidal functor
    \begin{align*}
        \gamma\colon \QCoh(X) \to \ICoh(X)
    \end{align*}
    compatible with pullback. Moreover, if $f\colon Y \to X$ is a qcqs suave map representable in lafp schemes over $\Lambda$ then the natural map $\gamma f_* \isoto f_* \gamma$ is an isomorphism.

    \item \label{rslt:psi-from-ICoh-to-QCoh-on-stacks} There is a natural transformation $\Psi\colon \ICoh \to \QCoh$ of functors $\Sch_\Lambda^{\afp,\qcqs} \to \Cat$, i.e. for every afp qcqs scheme $X$ over $\Lambda$ there is a functor
    \begin{align*}
        \Psi\colon \ICoh(X) \to \QCoh(X)
    \end{align*}
    compatible with lower-! on $\ICoh$ and pushforward on $\QCoh$. If $f\colon Y \to X$ is a suave map of qcqs afp schemes over $\Lambda$ then $f_!\colon \ICoh(Y) \to \ICoh(X)$ admits a left adjoint $f^\natural$ and the natural map $f^* \Psi \isoto \Psi f^\natural$ is an isomorphism. Moreover, if $X$ is suave over $\Lambda$ then $\Psi$ admits a left adjoint
    \begin{align*}
        \Xi\colon \QCoh(X) \to \ICoh(X)
    \end{align*}
    and the natural map $\Xi f_* \isoto f_! \Xi$ is an isomorphism.

    \item \label{rslt:morphisms-of-3ffs-to-ICoh} Let $E_0 \subseteq E$ be the class of qcqs maps that are representable in lafp schemes over $\Lambda$. Then there are sheafy 3-functor formalisms on $(\Stk_\Lambda^\lafp, E_0)$ together with morphisms of 3-functor formalisms
    \begin{align*}
        \ICoh \xfrom{\ \eta\ } \PQCoh_{\solid,\pc}^\op \injfrom \QCoh_{\solid,\pc}^\op \injto \QCoh_\solid^\op,
    \end{align*}
    which on affine schemes are defined as in \cref{def:pseudo-compact-modules} and \cref{def:PCoh-and-PQCoh-solid-pc-on-rings}.
\end{thmenum}
\end{thm}
\begin{proof}
The construction of the 6-functor formalism $\ICoh$ on lafp stacks follows immediately by applying \cite[Theorem~3.4.11]{heyer-mann-6ff} to the 6-functor formalism from \cref{def:6ff-for-ICoh-on-rings} using the sheafiness of $\ICoh$ by \cref{rslt:fppf-descent-for-all-sheaf-theories}. This immediately proves (i). Part (iii) follows in the same way as in the proof of \cref{rslt:stability-of-QCoh-solid-fine-maps}. Namely, the only non-formal input is that fppf covers in $\Alg_\Lambda^\lafp$ are universal $!$-covers. But this holds for $\ICoh$ by the proof of \cref{rslt:fppf-descent-for-all-sheaf-theories}, where it was shown that an fppf-cover $f\colon A \to B$ is $\ICoh$-prim and $f_* 1$ is descendable, so we can apply \cite[Lemma~4.7.4(i)]{heyer-mann-6ff}.

We next show that every qcqs open immersion of lafp stacks over $\Lambda$ is $\ICoh$-prim (i.e. automatically $\ICoh$-proper). By \cite[Lemma~4.5.7]{heyer-mann-6ff} this can be checked locally on the target, so we reduce to the case of a qcqs open immersion $f\colon U \injto \Spec A$. We cover $U$ by finitely many open affines $\Spec B_i \subseteq U$. Then each $f_i\colon \Spec B_i \to \Spec A$ is $\ICoh$-prim because it was shown in the proof of \cref{rslt:fppf-descent-for-all-sheaf-theories} that suave maps are $\ICoh$-prim. But $g_i\colon \Spec B_i \to U$ is a base-change of $f_i$ and thus also $\ICoh$-prim. Thus the $g_i$ form a finite cover by $\ICoh$-prim subspaces and hence a universal $!$-cover for $\ICoh$ (by the same argument as in \cite[Lemma~4.8.3]{heyer-mann-6ff}). Hence $f$ lies in $E$ and by \cite[Lemma~4.5.8(ii)]{heyer-mann-6ff} (cf. \cite[Remark~4.4.11(iii)]{heyer-mann-6ff}) we deduce that it is $\ICoh$-prim.

We now observe that every qcqs map $f\colon Y \to X$ representable in lafp schemes over $\Lambda$ is in $E$. Indeed, since $E$ is local on the target we can assume that $X = \Spec A$ is an affine scheme, thus $Y$ is a qcqs scheme. Then as explained in the previous paragraph, a finite open cover of $Y$ by affine subschemes is a finite $\ICoh$-prim cover by monomorphisms and in particular a universal $!$-cover for $\ICoh$. By $!$-locality of $E$ we reduce to the case that $Y$ is affine, where the claim is clear.

We next prove (vi). All the sheaf categories in (vi) satisfy descent by \cref{rslt:fppf-descent-for-all-sheaf-theories}. The 6-functor formalism $\QCoh_\solid$ was constructed in \cref{rslt:6ff-on-QCoh-solid}, where it was shown that all maps in $E_0$ are $\QCoh_\solid$-fine. To get the 3-functor formalism on $\QCoh_{\solid,\pc}^\op \subset \QCoh_\solid$, we just need to show that this subcategory is stable under $\tensor$, $f^*$ and $f_!$. The first two are clear by \cref{rslt:3ff-on-pseudo-compact-modules} and descent. To show stability under $f_!$, fix some map $f\colon Y \to X$ in $E_0$; we need to show that $f_!\colon \QCoh_\solid(Y) \to \QCoh_\solid(X)$ preserves $\QCoh_{\solid,\pc}$. By descent for $\QCoh_{\solid,\pc}$ and base-change for $f_!$ we can assume that $X$ is an affine scheme. Then $Y$ is a qcqs scheme, so we can pick a finite cover $Y = \bigcup_{i\in I} U_i$ by affine open subschemes $U_i \subseteq Y$. For every subset $J \subseteq I$ denote $U_J = \bigcap_{i\in J} U_j \subseteq Y$ and let $g_J\colon U_J \to X$ be the induced map. Then $f_! = \varinjlim_{\emptyset \ne J \subseteq I} g_{J!} g_J^*$, so it is enough to show that each $g_{J!}$ preserves $\QCoh_{\solid,\pc}$. This reduces the claim to the case that $Y$ is a quasicompact open subset of an affine scheme. Repeating the above argument again reduces to the case that $Y$ is affine (because now all the intersections are affine), where it follows from \ref{rslt:3ff-on-pseudo-compact-modules}. The construction of the 3-functor formalism $\PQCoh_{\solid,\pc}^\op$ works in the same way as for $\ICoh$ (note that $\PQCoh_{\solid,\pc}^\op$ is presentable) and by the uniqueness of the extension we obtain the morphism $\eta\colon \PQCoh_{\solid,\pc}^\op \to \ICoh$ from its restriction to affine schemes, where it was constructed in \cref{rslt:3ff-for-PCoh-on-rings}. All maps in $E_0$ admit $!$-functors for $\PQCoh_{\solid,\pc}$ by the same argument as for $\ICoh$.

For the proof of (vi) it remains to construct the embedding $\QCoh_{\solid,\pc} \injto \PQCoh_{\solid,\pc}$ of 3-functor formalisms. Let $E'_0 \subseteq E_0$ be the subclass of maps that are representable in \emph{affine} schemes. Then clearly there is an embedding of 3-functor formalisms $\QCoh_{\solid,\pc} \injto \PQCoh_{\solid,\pc}$ on $(\Stk_\Lambda^\lafp, E'_0)$, because by the uniqueness of extension in \cite[Proposition~3.4.2]{heyer-mann-6ff} this can be checked on the site $\Alg_\Lambda^{\afp,\op}$, where it is clear by construction. We thus obtain a natural transformation of presentable 6-functor formalisms $\Ind \comp \QCoh_{\solid,\pc}^\op \to \PQCoh_{\solid,\pc}^\op$ on $(\Stk_\Lambda^\lafp, E'_0)$, which extends uniquely to $E_0$ by \cite[Proposition~3.4.8(ii)]{heyer-mann-6ff} and the above observation that finite open covers are universal $!$-covers (which holds in by the same argument as for $\ICoh$). It only remains to verify that the 6-functor formalism on $\Ind \comp \QCoh_{\solid,\pc}^\op$ restricts to a 6-functor formalism on $\QCoh_{\solid,\pc}^\op$, because by uniqueness of the extension it must be the same as the one coming from $\QCoh_\solid^\op$. But this follows in the same way as above, finishing the proof of (vi).

We now prove (ii), so fix a qcqs map $f\colon Y \to X$ representable in lafp schemes over $\Lambda$. We have already seen above that $f \in E$, so it remains to verify (a), (b) and (c). Parts (a) and (b) are true for $\QCoh_\solid^\op$ by \cref{rslt:QCoh-solid-on-proper-and-smooth-maps}. Moreover, the suave respectively prim dualizing complex in $\QCoh_\solid^\op(Y)$ is pseudo-compact: In case (a) this is clear (because it is $1$) and in case (b) we need to see that $\omega_f$ is pseudo-compact, which reduces immediately to the case that $X$ and $Y$ are affine schemes, where it is an easy computation (see the proof of \cref{rslt:fppf-descent-for-all-sheaf-theories}). By (vi) and the induced functoriality of the categories of kernels (see \cite[Proposition~4.2.1(i)]{heyer-mann-6ff}) we deduce that (a) and (b) hold for $\QCoh_{\solid,\pc}^\op$ and hence also for $\ICoh$. Moreover, in (b) we see that $\delta_f^{\ICoh} = \eta(\omega_f)$. Since $\gamma$ (as constructed in \cref{rslt:construct-gamma-from-QCoh-to-ICoh-on-rings} and then extended to sheaves) factors over $\eta$, we arrive at (c).

We now prove (iv). The existence of $\gamma$ follows immediately from \cref{rslt:construct-gamma-from-QCoh-to-ICoh-on-rings} by passing to sheaves, so we only need to prove the second part of the claim. Let $f\colon Y \to X$ be a qcqs suave map that is representable in lafp schemes over $\Lambda$. By (iii).(b) and \cite[Lemma~4.5.7]{heyer-mann-6ff} $f$ is $\ICoh$-prim and hence by \cite[Lemma~4.5.13]{heyer-mann-6ff} the functor $f_*$ on $\ICoh$ satisfies arbitrary base-change. By \cref{rslt:qcqs-base-change-for-QCoh} the functor $f_*$ on $\QCoh$ also satisfies arbitrary base-change. We can thus pass to a cover of $X$ by affine schemes in order to reduce to the case that $X$ itself is an affine scheme. Then $Y$ is a qcqs scheme and as in the proof of \cref{rslt:qcqs-base-change-for-QCoh} we can pass to a finite affine cover of $Y$ to first reduce to the case that $Y$ is an open subscheme of an affine scheme and then to the case that $Y$ itself is affine. But then the claim was shown in \cref{rslt:gamma-compatible-with-suave-pushforward-on-rings}.

It remains to prove (v). We already constructed the natural transformations $\Psi$ and $\Xi$ on affine schemes in \cref{rslt:psi-for-ICoh-and-QCoh-on-rings}. To construct $\Psi$ on qcqs afp schemes, we view both $\QCoh$ and $\ICoh$ as functors $\Sch_\Lambda^{\afp,\qcqs} \to \PrL$, with transition maps given by pushforward and lower-!, respectively. It is enough to show that both functors are left Kan extensions of their restriction to affine schemes. Since colimits in $\PrL$ are computed as limits in $\PrR$, we can equivalently consider $\QCoh$ and $\ICoh$ as functors $\Sch_\Lambda^{\afp,\qcqs,\op} \to \Cat$, where the transition maps are given by the right adjoints of pushforward (to be denoted upper-?) resp. upper-!, and ask whether these functors are the right Kan extension of their restriction to affine schemes. Note first that both of these functors satisfy fppf descent on affine schemes: In the case of $\ICoh$ this follows from \cite[Lemma~4.7.4(i)]{heyer-mann-6ff}, because as noted in the proof of \cref{rslt:fppf-descent-for-all-sheaf-theories}, a map $f\colon Y \to X$ of affine schemes is $\ICoh$-prim and $f_* 1$ is descendable. The same argument works for $\QCoh$, by noting that $\QCoh$ forms a 6-functor formalism on affine schemes with $f_! = f_*$ for all maps $f$. It remains to show that $\QCoh$ and $\ICoh$ (with upper-! pullbacks) descend along Zariski open covers of qcqs schemes. But these form finite $\ICoh$- and $\QCoh$-proper covers and thus universal !-covers (as discussed in the second paragraph of this proof). This finishes the construction of $\Psi$ on qcqs afp schemes.

We now prove the remaining claims in (v). Fix a suave map $f\colon Y \to X$ of qcqs afp schemes. By (ii).(b) $f$ is $\ICoh$-prim and in particular $f_!$ admits a left adjoint $f^\natural$. To prove the claimed isomorphism $f^* \Psi_X \isoto \Psi_Y f^\natural$, we pass to right adjoints: Denoting $R_X$ the right adjoint of $\Psi_X$, we need to see that the natural map $R_X f_* \to f_! R_Y$ is an isomorphism of functors $\QCoh(Y) \to \ICoh(X)$. Write $X = \varinjlim_i X_i$ where all $X_i$ are affine schemes (and the colimit is indexed over an arbitrary small category $\cat I$), and denote $Y_i := Y \times_X X_i$, so that $Y = \varinjlim_i Y_i$. Then $\ICoh(X) = \varprojlim_i \ICoh^!(X_i)$ and $\QCoh(X) = \varprojlim_i \QCoh^?(X_i)$ (where the transition maps are given by upper-! resp. upper-? functors) and $R_X$ is given componentwise by $(R_{X_i})_i$. Moreover, since $f$ is $\ICoh$-prim, $f_!$ satisfies base-change with respect to the upper-! functors in the above limit (see \cite[Lemma~4.5.13(ii)]{heyer-mann-6ff}) and is thus given componentwise by $(f_{i!})_i$ for $f_i\colon Y_i \to X_i$ the natural maps. A similar statement holds for $f_*$ on $\QCoh$ (using again that $\QCoh$ is itself a 6-functor formalism with $f_! = f_*$ on affines). Altogether we see that the desired isomorphism $R_X f_* \to f_! R_Y$ reduces to showing that $R_{X_i} f_{i*} \to f_{i!} R_{Y_i}$ is an isomorphism for all $i$, so we can reduce to the case that $X$ is affine. In an iterative procedure we can now deduce the claim first for affine open immersions (by reducing to \cref{rslt:pullback-compatibility-of-psi-on-rings}), then for general open immersions and then for $Y$ by covering it by passing to an open affine cover.

Now assume that $X$ is suave. We want to show that $\Psi_X$ admits a left adjoint $\Xi_X$. First assume that $X$ is an open subscheme of an affine scheme. Cover $X$ by affine open subsets $U_i \subseteq X$. Then all the intersections of the $U_i$'s are affine and we know that $\Xi$ exists on all these and commutes with pullback on $\QCoh$ and upper-$\natural$ on $\ICoh$ (because $\Psi$ is compatible with $f_*$ and $f_!$). By writing $\ICoh(X)$ and $\QCoh(X)$ as the limit of the respective categories over the cover (using upper-$\natural$ and upper-* transition maps; the limit holds by $\ICoh$-primness of the cover) we see that $\Xi$ assembles into a functor $\Xi_X\colon \QCoh(X) \to \ICoh(X)$ that is computed componentwise in that limit. By the previous paragraph we also know that $\Psi_X$ is computed componentwise in this limit. Since on each term of the limit, $\Psi$ and $\Xi$ are adjoint, we deduce that they are adjoint on $X$. This proves the existence of $\Xi_X$ in the case that $X$ is open inside an affine scheme. In general, again choose a cover of $X$ by affine opens; then all the intersections of these opens are open inside some affine scheme and we can argue as before, using that in this case we constructed $\Xi$.

Finally, we prove that the natural map $\Xi_X f_* \to f_! \Xi_Y$ is an isomorphism, where $f\colon Y \to X$ is suave. This can be checked after $\natural$-pulling back to an affine cover of $X$, and since $\Xi$ is compatible with $\natural$-pullback, we immediately reduce to the case that $X$ is affine. Now pass to a finite affine open cover of $Y$ (first in the case that $Y$ is itself an open inside an affine scheme, then for general $Y$), write $f_*$ as a finite limit over that cover and use that $\Xi_X$ preserves finite limits. This reduces the claim to the case that $Y$ is affine and hence to \cref{rslt:xi-for-QCoh-and-ICoh-on-rings}.
\end{proof}

\begin{defn}
Fix a classical regular noetherian ring $\Lambda$. A map $f\colon Y \to X$ of lafp stacks over $\Lambda$ is called \emph{$\ICoh$-fine} if it lies in the class $E$ of \cref{rslt:6ff-for-ICoh}, so that $f_!$ and $f^!$ exist on $\ICoh$.
\end{defn}

We have finally constructed the 6-functor formalism for ind-coherent sheaves. We postpone a more detailed study of this 6-functor formalism to the next subsection and finish the current subsection with several remarks.

\begin{rmk}
The only place where we used that the base ring $\Lambda$ is a classical regular noetherian ring is in the proof of \cref{rslt:eta-compatible-with-certain-pull-and-push}, because this proof relied on the t-structure on $\Coh$. It is an interesting question whether \cref{rslt:eta-compatible-with-certain-pull-and-push} is true without any assumption on $\Lambda$ (replacing $\Coh$ by $\Coh_\Lambda$), which would provide a theory of ind-coherent sheaves over any base ring. This theory would have some strange behavior though: Let $k$ be a field and consider it as an algebra over $\Lambda := k[\varepsilon]/\varepsilon^2$; then one checks that $\Coh_\Lambda(k) = 0$.
\end{rmk}

\begin{rmk} \label{rslt:explicit-description-of-ICoh-functors-on-affines}
One needs to be careful with the transition functors on $\ICoh$. Suppose $f\colon A \to B$ is a map of afp algebras over some fixed base $\Lambda$ (which is assumed to be classical regular noetherian, as always).
\begin{rmksenum}
    \item \label{rslt:ICoh-pullback-on-affines-is-upper-shriek} Assume that $f$ is suave. Then the pullback $f^*\colon \D_\solid(A) \to \D_\solid(B)$ restricts to a functor $f^*\colon \Coh(A) \to \Coh(B)$ and hence induces functors $f^*\colon \PCoh(A) \to \PCoh(B)$ and $\ICoh(A) \to \ICoh(B)$. However, while the pullback on $\PCoh$ is indeed the correct one, the naive pullback on $\ICoh$ is not the one we call $f^*$. Namely, the pullback on $\ICoh$ is computed, at least on coherent sheaves, as $f^*(P) = \SD(f^*\SD(P)) = f^! P$, where $\SD = \intHom(-,\omega)$ denotes suave duality and the $f^*$ appearing in the middle term is the canonical one on $\Coh$ from above. In other words, for suave $f$ the pullback on $\ICoh$ is given by the ind-extension of the functor $f^!\colon \Coh(A) \to \Coh(B)$. This is the reason why this pullback is often denoted $f^{!,\ICoh}$ in the literature, but we prefer to use the notation $f^*$ to emphasize its role in terms of the 6-functor formalism.

    \item \label{rslt:ICoh-lower-shriek-on-affines-is-pushforward} Assume that $f$ is finite. Then $f_*\colon \D_\solid(B) \to \D_\solid(A)$ restricts to a functor $f_*\colon \Coh(B) \to \Coh(A)$ and hence induces functors $\PCoh(B) \to \PCoh(A)$ and $\ICoh(B) \to \ICoh(A)$. These functors are the ones denotes $f_!$ above. In other words, $f_!\colon \ICoh(B) \to \ICoh(A)$ is the $\Ind$-extension of the forgetful functor and therefore often denoted $f_*^{\ICoh}$ in the literature.
\end{rmksenum}
\end{rmk}

\begin{rmk}
The 3-functor formalism on $\PCoh$ can be constructed in the same way as above on all noetherian schemes (i.e. without a specific base ring) if one defines $\Coh \subseteq \QCoh$ as the full subcategory of bounded complexes with finitely generated cohomologies. If one wants to pass to $\ICoh$ then a regular noetherian base is very helpful, though.
\end{rmk}

\subsection{Properties of ind-coherent sheaves} \label{sec:alggeo.properties-of-ICoh}

With the 6-functor formalism for $\ICoh$ constructed in the previous subsection (see \cref{rslt:6ff-for-ICoh}) we now come to the study of some of its basic properties. Our first goal is to provide a large class of $\ICoh$-fine maps. For this, we need one more computation (recall the definition of $\QCoh_{\solid,\pc}$ from \cref{def:pseudo-compact-modules}):

\begin{lem} \label{rslt:QCA-lower-shriek-preserves-pseudocompact}
Let $f\colon Y \to X$ be a lafp QCA map of stacks over some $\Q$-algebra $\Lambda$. Then $f_!\colon \QCoh_\solid(Y) \to \QCoh_\solid(X)$ preserves $\QCoh_{\solid,\pc}$.
\end{lem}
\begin{proof}
By descent (see \cref{rslt:fppf-descent-for-all-sheaf-theories}) this can be checked after base-change to some affine cover of $X$, so we can assume that $X = \Spec A$ is an affine scheme. By definition of QCA stacks, there is a smooth affine cover $g\colon U \surjto Y$. We denote its \v{C}ech nerve by $g_\bullet\colon U_\bullet \to Y$. Then for $M \in \QCoh_\solid(Y)$ we have $f_! M = \varinjlim_{n\in\Delta} (fg_n)_! g_n^! M$. Since $\QCoh_{\solid,\pc}(A)$ is stable under uniformly right-bounded geometric realizations, it is enough to show that the family $((fg_n)_! g_n^! M)$ is uniformly right-bounded.

For $n \ge 0$ we denote by $h_n\colon U_n = U^{\times_Y n+1} \to U$ the projection to the first factor in the fibrer product, set $M_0 := g^! M$ and denote $U_n = \Spec B^n$ and $B := B^0$. Then $(fg_n)_! g_n^! M = (fg)_! h_{n!} h_n^! M_0$. Now $(fg)_!\colon \QCoh_\solid(B) \to \QCoh_\solid(A)$ is right-bounded: By factoring this map into a polynomial ring and a closed immersion, we can reduce to these two cases; in the second case it is a forgetful functor and in the first case it reduces to $\mathbb A^1$, where it is a standard computation (see \cref{ex:solid-lower-shriek-for-A1}). Furthermore, $M_0 = g^* M \tensor \omega_g$ is pseudocompact and hence right-bounded. We are thus left to show that the functors $h_{n!} h_n^!$ are right-bounded, independent of $n$. Since $h_n$ is suave we have $h_{n!} h_n^! = h_{n\natural} h_n^*$, where $h_{n\natural}$ is the left adjoint of $h_n^*$ (this can be seen by direct computation or using the self-duality of the category of kernels). Now $h_n$ is right t-exact and since $h_n$ is flat, it has Tor dimension $\le 1$ by \cite[Lemma~4.10]{mikami-fppf-descent2}. This implies that $h_{n\natural}$ is right-bounded by $1$, so we conclude the same for $h_{n\natural} h_n^*$, as desired.
\end{proof}

With \cref{rslt:QCA-lower-shriek-preserves-pseudocompact} and our preparations in \cref{sec:alggeo.QCA-stacks} at hand, we easily obtain the following result:

\begin{thm}
Fix a classical regular noetherian $\Q$-algebra $\Lambda$ and let $f\colon Y \to X$ be a QCA map of lafp stacks over $\Lambda$. Then:
\begin{thmenum}
    \item \label{rslt:QCA-maps-are-ICoh-fine} $f$ is $\ICoh$-fine and hence admits $!$-functors in the $\ICoh$-formalism.
    
    \item \label{rslt:ICoh-codualizing-sheaf-for-QCA-maps} If $f$ is suave then it is $\ICoh$-prim, hence $f_* = f_!(\delta_f^{\ICoh} \tensor -)$ on $\ICoh$. If $\omega_f$ is invertible then $\delta_f^{\ICoh} = \gamma(\omega_f^\vee)$, where $\omega_f \in \QCoh(Y)$ is the $\QCoh_\solid$-dualizing complex.

    \item \label{rslt:gamma-commutes-with-pushforward-along-suave-QCA-maps} If $f$ is suave then the natural map $\gamma f_* \isoto f_* \gamma$ is an isomorphism of functors $\QCoh(Y) \to \ICoh(X)$.
\end{thmenum}
\end{thm}
\begin{proof}
We first prove (i). By \cref{rslt:stability-of-ICoh-fine-maps} we can check that $f$ is $\ICoh$-fine after base-change to any cover of $X$, so we may assume that $X$ is an affine scheme. In particular $Y$ is a QCA stack over $\Lambda$. Pick any smooth geometric cover $f\colon U \surjto Y$ by a scheme $U$. Note that $f$ is affine (because $Y$ has affine diagonal) and in particular $\ICoh$-fine. Clearly the map $U \to X$ is $\ICoh$-fine, so by \cref{rslt:stability-of-ICoh-fine-maps} it is enough to show that $f$ is a universal $\ICoh$-$!$-cover. By \cref{rslt:ICoh-on-lafp-schemes} $f$ is $\ICoh$-prim, so by \cite[Lemma~4.7.4]{heyer-mann-6ff} it is enough to show that $f_* 1 \in \Coh(Y)$ is descendable. On the other hand, $f_* 1 = \gamma(f_* 1)$ by \cref{rslt:gamma-from-QCoh-to-ICoh-on-stacks}, hence it is enough to show that $f_* 1 \in \QCoh(Y)$ is descendable. But this was shown in \cref{rslt:descendable-cover-of-QCA-stack}.

To prove (ii), we observe that as in the proof of \cref{rslt:morphisms-of-3ffs-to-ICoh}, the 6-functor formalism $\PQCoh_{\solid,\pc}^\op$ also extends to QCA maps (by the same argument as for $\ICoh$) and we get maps of 3-functor formalisms
\begin{align*}
    \ICoh \xfrom{\ \eta\ } \PQCoh_{\solid,\pc}^\op \injfrom \QCoh_{\solid,\pc}^\op \injto \QCoh_\solid^\op,
\end{align*}
which are now defined on all QCA maps; the only non-formal input is \cref{rslt:QCA-lower-shriek-preserves-pseudocompact}. Then the rest of (ii) follows as in the proof of \cref{rslt:ICoh-on-lafp-schemes}.

We now prove (iii). As in the proof of \cref{rslt:gamma-from-QCoh-to-ICoh-on-stacks} we can reduce to the case that $X = \Spec A$ is an affine scheme. Then $Y$ is a QCA stack and we can pick a smooth affine cover $g\colon U \surjto Y$ with \v{C}ech nerve $g_\bullet\colon U_\bullet \to Y$. All $U_n$ are affine schemes, hence we know that $\gamma$ commutes with the pushfoward along $f g_n$ for all $n$. Since $g$ is descendable by \cref{rslt:descendable-cover-of-QCA-stack}, the functor $f_*$ is a retract of a \emph{finite} limit of the functors $(f g_n)_* g_n^*$. As $\gamma$ commutes with finite limits, we conclude
\end{proof}

The next result is a version of \cref{rslt:QCoh-minus-only-depends-on-truncated-rings} for $\ICoh$, showing that $\ICoh(X)$ only depends on the value of $X$ on \emph{truncated} rings. This technical result will be used at different places below.

\begin{lem} \label{rslt:ICoh-only-depends-on-truncated-rings}
Fix a classical regular noetherian ring $\Lambda$ and let $f\colon Y \to X$ be a map of lafp stacks over $\Lambda$. Assume that for every truncated afp $\Lambda$-algebra $A$, the induced map $Y(A) \isoto X(A)$ is an isomorphism. Then $f^*$ induces an isomorphism
\begin{align*}
    f^*\colon \ICoh(Y) \isoto \ICoh(X).
\end{align*}
\end{lem}
\begin{proof}
As in the proof of \cref{rslt:QCoh-minus-only-depends-on-truncated-rings} we immediately reduce to showing that for every afp $\Lambda$-algebra $A$ with truncation maps $g_n\colon A \to \tau_{\le n} A$ the natural functor
\begin{align*}
    F\colon \ICoh(A) \isoto \varprojlim_{n\ge 0} \ICoh(\tau_{\le n} A), \qquad M \mapsto (g_n^* M)_n
\end{align*}
is an isomorphism. To see this, we use the argument in \cite[p.~72]{scholze-6ff}. By \cref{rslt:ICoh-on-lafp-schemes} each $g_n$ is $\ICoh$-smooth with trivial dualizing complex, hence the left adjoint of $g_n^*$ is given by $g_{n!}$. This implies that the above limit can also be computed as the colimit $\varinjlim_n \ICoh(\tau_{\le n} A)$ in $\Pr^L$, where the transition maps are now given by the $\ICoh$ lower-$!$ maps. By construction (see \cref{rslt:ICoh-lower-shriek-on-affines-is-pushforward}) these transition functors are the ind-extensions of the forgetful functors $\Coh(\tau_{\le n} A) \to \Coh(\tau_{\le m} A)$ for $n \ge m$. By \cite[Corollary~A.5.9]{heyer-mann-6ff} the colimit is given as $\Ind(\varinjlim_n \Coh(\tau_{\le n} A))$, so we are reduced to showing that the natural map
\begin{align*}
    \varinjlim_n \Coh(\tau_{\le n} A) \isoto \Coh(A)
\end{align*}
is an equivalence, where the filtered colimit is taken in $\Cat$. The functor is clearly essentially surjective, as every $M \in \Coh(A)$ is bounded and hence comes from some $\Coh(\tau_{\le n} A)$. To prove full faithfulness, fix $M, N \in \Coh(\tau_{\le n} A)$ for some $n \ge 0$. It is enough to show that $\Hom_{\tau_{\le m} A}(M, N) = \Hom_A(M, N)$ for $m \gg n$. But the right-hand side can be rewritten as $\Hom_{\tau_{\le m} A}(M \tensor_A \tau_{\le m} A, N)$. Since $N$ is bounded, it is enough to show that the map $\tau_{\le m} A \tensor_A M \to M$ is increasingly connective as $m$ increases. But this follows from $M$ being right-bounded, see the proof of \cref{rslt:QCoh-minus-only-depends-on-truncated-rings}.
\end{proof}

The comparison functor $\gamma\colon \QCoh \to \ICoh$ is a handy tool for studying $\ICoh$ and appears naturally when one tries to explicitly construct objects in $\ICoh$. The next result provides some additional properties of $\gamma$. 

\begin{prop}
Fix a classical regular noetherian ring $\Lambda$ and a lafp stack $X$ over $\Lambda$.
\begin{propenum}
    \item \label{rslt:gamma-fully-faithful-on-QCoh-minus} $\gamma$ induces a fully faithful symmetric monoidal embedding $\QCoh^-(X) \injto \ICoh(X)$. 
    
    \item \label{rslt:Perf-equals-dualizable-objects-in-ICoh} The essential image of $\Perf(X) \injto \ICoh(X)$ consists precisely of the dualizable objects in $\ICoh(X)$.

    \item \label{rslt:gamma-fully-faithful-for-suave} Suppose that $X$ is geometric and suave. Then $\gamma\colon \QCoh(X) \injto \ICoh(X)$ is fully faithful. If $X$ is smooth then $\gamma$ is an equivalence.
\end{propenum}
\end{prop}
\begin{proof}
We first prove (iii), so assume that $X$ is geometric and suave. Pick a smooth cover $U \to X$ by some scheme $U$. Then both $\QCoh(X)$ and $\ICoh(X)$ descend along the associated \v{C}ech nerve $U_\bullet \to X$, i.e. are computed as the limit of their values on all $U_m$. Both the full faithfulness and the isomorphism claim of $\gamma$ can therefore be checked on all $U_m$. By induction on the number $n$ such that $X$ is $n$-geometric, we can reduce to the case that $X$ is a scheme. By passing to an affine cover we can further reduce to the case that $X$ is quasi-affine and then that $X = \Spec A$ is affine. Since $A$ is suave, the tensor unit in $\ICoh(A)$ is given by the coherent $A$-module $\omega_A$ and the functor $\gamma$ is computed as $\gamma(P) = P \tensor \omega_A \in \Coh(A)$ on perfect $A$-modules $P$ (and then ind-extended to $\D(A)$). In particular $\gamma$ sends compact objects to compact objects and is fully faithful on them (because $\Hom_A(\omega_A, \omega_A) = \Hom_A(A, A)$ by suave duality), hence $\gamma$ is fully faithful. If $A$ is smooth then $\Coh(A) = \Perf(A)$ and $\omega_A$ is an invertible $A$-module -- this easily implies that $\gamma$ is an equivalence.

We next prove (i), so let $X$ be general again. Note that full faithfulness passes through limits, so by writing $X$ as a colimit of affine schemes we can reduce to the case that $X = \Spec A$ itself is an affine scheme. By the proofs of \cref{rslt:QCoh-minus-only-depends-on-truncated-rings} and \cref{rslt:ICoh-only-depends-on-truncated-rings} both $\ICoh(A)$ and $\QCoh^-(A)$ can be written as the limit of their values on the truncations of $A$, which further reduces us to the case that $A$ is truncated, i.e. suave. But then the claim follows from (iii).

It remains to prove (ii). Since dualizability descends along limits (e.g. by \cite[Proposition~D.2.16(ii)]{heyer-mann-6ff}), the same reduction steps as in the proof of (i) apply to reduce to the case that $X = \Spec A$ is for a truncated ring $A$. In this case the tensor unit in $\ICoh(A) = \Ind(\Coh(A))$ is coherent and hence compact, which implies that every dualizable object in $\ICoh(A)$ is compact and hence coherent. Recall that the tensor product on $\Coh(A)$ is defined via suave duality from the canonical tensor product, so by going through the definition of $\gamma$ the claim reduces to the following statement: The essential image of the canonical embedding $\Perf(A) \injto \Coh(A)$ consists precisely of the dualizable objects under the canonical tensor product on $\QCoh(A)$. But this claim is proved in \cref{rslt:Perf-same-as-dualizable}.
\end{proof}

Next we discuss the excision sequence for $\ICoh$, adapting the argument for $\QCoh$ in \cref{rslt:excision-for-QCoh}. As we will see, the statement for $\ICoh$ is much cleaner.

\begin{defn}
Fix a classical regular noetherian ring $\Lambda$, a lafp stack $X$ over $\Lambda$ and a closed subset $Z \subseteq \abs X$. We denote by
\begin{align*}
    \ICoh(X)_Z \subseteq \ICoh(X)
\end{align*}
the full subcategory spanned by those sheaves $M \in \ICoh(X)$ such that $j^* M = 0$, where $j\colon U \injto X$ is the open complement of $Z$.
\end{defn}

\begin{lem}
Fix a classical regular noetherian ring $\Lambda$ and let $i\colon Z \to X$ be a closed immersion of lafp stacks over $\Lambda$. Let $\hat i\colon X_{\hat Z} \injto X$ denote the completion of $Z$ in $X$ and let $i_n\colon Z_n \to X$ denote the $n$-th infinitesimal neighborhood.
\begin{lemenum}
    \item \label{rslt:hat-i-pullback-induces-equivavence-on-ICoh} The pullback functor $\hat i^*$ induces an equivalence
    \begin{align*}
        \hat i^*\colon \ICoh(X)_Z \isoto \ICoh(X_{\hat Z}).
    \end{align*}

    \item \label{rslt:fibers-of-inf-nbhds-on-ICoh} For all $n \ge 0$ there is some $L_n \in \ICoh(Z)$ such that for all $M \in \ICoh(X)$ there is a natural isomorphism
    \begin{align*}
        \cofib(i_{n!} i_n^* M \to i_{n+1,!} i_{n+1}^* M) = i_! (i^* M \tensor L_n).
    \end{align*}

    \item \label{rslt:pullback-from-completion-to-Z-is-conservative-on-ICoh} The functor $i^*\colon \ICoh(X)_Z \to \ICoh(Z)$ is conservative.
\end{lemenum}
\end{lem}
\begin{proof}
Let $j\colon U \to X$ denote the open complement of $i$ and for $n \ge 0$ let $i_n\colon Z_n \to X$ denote the $n$-th infinitesimal neighborhood of $Z$ in $X$ (as defined in \cref{def:infinitesimal-neighbourhoods}). Note that $i$ is automatically lafp by \cref{rslt:schematic-map-between-lafp-stacks-is-lafp} and in particular it is finitary. Thus by \cref{rslt:infinitesimal-neighborhoods} and \cref{rslt:ICoh-only-depends-on-truncated-rings} we deduce that the pullbacks induce an isomorphism
\begin{align*}
    \ICoh(X_{\hat Z}) \isoto \varprojlim_n \ICoh(Z_n).
\end{align*}
By \cref{rslt:ICoh-on-lafp-schemes}, each of the pullbacks $i_n^*$ admit the left adjoint $i_{n!}$ and the right adjoint $i_{n*}$. Hence by \cite[Lemma~D.4.7]{heyer-mann-6ff} the functor $\hat i^*$ admits both a left adjoint $\hat i_\natural$ and a right adjoint $\hat i_*$.

We now claim that $\hat i_\natural$ is fully faithful, i.e. that the unit $\id \isoto \hat i^* \hat i_\natural$ is an isomorphism. By passing to right adjoints, this is equivalent to showing that the counit $\hat i^* \hat i_* \isoto \id$ is an isomorphism. By the above limit formula for $\ICoh(X_{\hat Z})$ this can be checked after pullback along $i'_n\colon Z_n \to X_{\hat Z}$ for all $n$, which reduces the problem to showing that the natural maps $i_n^* \hat i_* \isoto i_n'^*$ are isomorphisms. But this follows easily from the fact that $\hat i_*$ satisfies base-change along $i_n^*$ by \cite[Lemma~4.5.13(i)]{heyer-mann-6ff}. One similarly sees that $j^* \hat i_\natural = 0$: After passing to right adjoints this becomes $\hat i^* j_* = 0$, which follows from base-change (see \cite[Lemma~4.5.13]{heyer-mann-6ff} and note that $j$ is $\ICoh$-prim by \cref{rslt:ICoh-on-lafp-schemes}). Thus $\hat i_\natural$ factors over $\ICoh(X)_Z$, so we obtain the adjunction
\begin{align*}
    \hat i_\natural\colon \ICoh(X_{\hat Z}) \leftrightarrow \ICoh(X)_Z \noloc \hat i^*.
\end{align*}
We have already seen above that $\hat i_\natural$ is fully faithful, so in order to show that the above adjunction consists of equivalences, it is enough to see that $\hat i^*$ is conservative. This can be checked after passing to a cover, so we can assume that $X = \Spec A$ is an affine scheme and in particular $Z$ induces a finitely generated ideal $(f_1, \dots, f_d) \subset \pi_0 A$. The $f_i$'s induce a map $X \to \mathbb A^d_\Lambda$ and we let $Z' \to X$ denote the pullback of the closed immersion $\{ 0 \} \to \mathbb A^d_\Lambda$. Then $Z$ and $Z'$ have the same image in $\abs X$ and thus $X_{\hat Z} = X_{\hat Z'}$. Therefore, for proving that $\hat i$ is conservative we can replace $Z$ by $Z'$. It is now enough to show that $i^*\colon \ICoh(X)_Z \to \ICoh(Z)$ is conservative, i.e. that the pair $(i^*, j^*)$ is conservative on $\ICoh(X)$. As in the last part of the proof of \cref{rslt:hat-i-pullback-induces-equiv} we can inductively reduce to the case $d = 1$, so that $Z = \Spec A/f$ and $U = \Spec A[1/f]$ for some element $f \in \pi_0 A$. Now suppose we have some $M \in \ICoh(X)$ such that $i^* M = 0$ and $j^* M = 0$. Note that both $i$ and $j$ are cohomologically smooth, hence $\ICoh$-prim by \cref{rslt:ICoh-on-lafp-schemes}. Therefore $0 = i_* i^* M = i_*1 \tensor M$ and $0 = j_* j^* M = j_* 1 \tensor M$. By \cref{rslt:gamma-from-QCoh-to-ICoh-on-stacks} we compute
\begin{align*}
    i_* 1 = \gamma(i_* 1) = \cofib(1 \xto{f} 1), \qquad j_* 1 = \gamma(j_* 1) = \varinjlim [1 \xto{f} 1 \xto{f} 1 \xto{f} \dots].
\end{align*}
Altogether we deduce that $\cofib(M \xto{f} M) = 0$ and $\varinjlim_f M = 0$, which implies $f = 0$. This finishes the proof of (i).

We now prove (ii) so let $Z$ and $X$ be general again. By the projection formula the claim immediately reduces to the case $M = 1$, i.e. we have to show that $\cofib(i_{n!} 1 \to i_{n+1,!} 1)$ lies in the image of $i_!$. Consider the morphism of 3-functor formalisms $\eta\colon \QCoh^\op_{\solid,\pc} \to \ICoh$ from \cref{rslt:morphisms-of-3ffs-to-ICoh}. It is thus enough to show the same claim in $\QCoh^\op_{\solid,\pc}$, i.e. we need to show that $\cofib(i_{n!} 1 \to i_{n+1,!} 1)$ lies in the image of $i_!$ in this 3-functor formalism. By removing the $(-)^\op$ and observing that $i_{n!} = i_{n*}$ in $\QCoh_{\solid,\pc}$ we are reduced to showing that $\fib(i_{n+1,*} 1 \to i_{n*} 1)$ lies in the image of $i_*$ in $\QCoh_{\solid,\pc}$. But this is shown in \cref{rslt:fiber-of-infinitesimal-nbhds}, where this fiber is explicitly computed as $i_* \Sym_Z^n(L_{Z/X}[-1])$ (note that $\Sym_Z^n(L_{Z/X}[-1])$ is pseudocoherent because $i$ is lafp, see \cref{rslt:cotangent-complex-for-geometric-map-and-smoothness} and \cref{rlst:Sym-preserves-pseudo-coherence}).

It remains to prove (iii). As the claim can be checked after passing to any cover, we can assume that $X = \Spec A$ is an affine scheme. Fix $M \in \ICoh(X)_Z$ such that $i^* M = 0$; we need to show that then $M = 0$. By the proof of (i) it is enough to show that $\hat i_\natural \hat i^* M = 0$. By the above limit $\ICoh(X_{\hat Z}) = \varprojlim_n \ICoh(Z_n)$ and \cite[Lemma~D.4.7(i)]{heyer-mann-6ff} we have $\hat i_\natural \hat i^* M = \varinjlim_n i_{n\natural} i_n^* M$, so it is enough to show that $i_{n\natural} i_n^* M = 0$ for all $n \ge 0$ -- here we note that $i_{n\natural} = i_{n!}$ by \cref{rslt:ICoh-on-lafp-schemes}. But this follows inductively from (ii).
\end{proof}

\begin{thm} \label{rslt:ICoh-excision-all}
Fix a classical regular noetherian ring $\Lambda$. Let $i\colon Z \to X$ be a closed immersion of lafp stacks over $\Lambda$ with completion $\hat i\colon X_{\hat Z} \injto X$ and open complement $j\colon U \injto X$. Let $M \in \ICoh(X)$.
\begin{thmenum}
    \item \label{rslt:ICoh-excision-triangles} The map $\hat i$ is $\ICoh$-étale (and in particular $\ICoh$-fine) and the map $j$ is $\ICoh$-proper. Moreover, we have natural fiber sequences
    \begin{align*}
        \hat i_! \hat i^* M \to M \to j_* j^* M, \qquad j_* j^! M \to M \to \hat i_* \hat i^* M.
    \end{align*}

    \item \label{rslt:ICoh-excision-filtration} For $n \ge 0$ let $i_n\colon Z_n \to X$ denote the $n$-th infinitesimal neighborhood of $i$. Then we have natural isomorphisms
    \begin{align*}
        \hat i_! \hat i^* M = \varinjlim_n i_{n!} i_n^* M, \qquad \hat i_* \hat i^* M = \varprojlim_n i_{n*} i_n^* M.
    \end{align*}

    \item \label{rslt:ICoh-excision-normal-bundle} Suppose that $L_{Z/X} \in \QCoh(Z)$ is perfect and denote $\Nm_{Z/X} := (L_{Z/X}[-1])^\vee$ the normal bundle of $i$. Then for all $n \ge 0$ we have
    \begin{align*}
        \cofib(i_{n!} i_n^* M \to i_{n+1,!} i_{n+1}^* M) = M \tensor i_! \gamma(\Sym^n_Z(\Nm_{Z/X})).
    \end{align*}
\end{thmenum}
\end{thm}
\begin{proof}
We first prove (i). By \cref{rslt:ICoh-on-lafp-schemes} $j$ is $\ICoh$-prim and hence $\ICoh$-proper (because it is a monomorphism). To show that $\hat i$ is $\ICoh$-étale, it is enough to show that it is $\ICoh$-suave (because it is a monomorphism). By \cref{rslt:ICoh-on-lafp-schemes} we know that $i$ is $\ICoh$-suave. Hence the same is true for the base-change of $i$ along $\hat i$, which is the map $i'\colon Z \to X_{\hat Z}$. By \cref{rslt:pullback-from-completion-to-Z-is-conservative-on-ICoh} $i'^*$ is conservative and hence by \cite[Lemma~4.7.1]{heyer-mann-6ff} $i'$ is a universal $\ICoh$-$!$-cover. By \cref{rslt:stability-of-ICoh-fine-maps} the factorization $Z \to X_{\hat Z} \injto X$ then implies that $\hat i$ is $\ICoh$-fine, and by \cite[Lemma~5.4.8(i)]{heyer-mann-6ff} the same factorization implies that it is $\ICoh$-suave.

To finish the proof of (i), it remains to show the claimed fiber sequences. Let $F := \fib(M \to j_* j^* M)$. Then clearly $F \in \ICoh(X)_Z$ (note that $j_* = j_!$ satisfies base-change). By \cref{rslt:hat-i-pullback-induces-equivavence-on-ICoh} we can identify $\ICoh(X)_Z$ with the image of the fully faithful functor $\hat i_!$. Hence we have $F = \hat i_! \hat i^* F$ and altogether
\begin{align*}
    F = \hat i_! \hat i^* F = \fib(\hat i_! \hat i^* M \to \hat i_! \hat i^* j_* j^* M) = \hat i_! \hat i^* M,
\end{align*}
where in the third step we used $\hat i^* j_* = 0$ by base-change. This proves the first fiber sequence. The second fiber sequence can be proved similarly or deduced from the first using $\intHom$, see the proof of \cref{rslt:excision-for-QCoh}.

To prove (ii) we observe that $\ICoh(X_{\hat Z}) = \varprojlim_n \ICoh(Z_n)$ by \cref{rslt:infinitesimal-neighborhoods} and \cref{rslt:ICoh-only-depends-on-truncated-rings} and so the claim follows easily from \cite[Lemma~D.4.7]{heyer-mann-6ff} (using that $i_{n!}$ is left adjoint to $i_n^*$ by \cref{rslt:ICoh-on-lafp-schemes}).

We now prove (iii). By \cref{rslt:fibers-of-inf-nbhds-on-ICoh} we already know that the claimed identity holds with $\gamma(\Sym^n(\Nm_{Z/X}))$ replaced by a certain element $L_n \in \ICoh$. By the proof of \cref{rslt:fibers-of-inf-nbhds-on-ICoh} we see that $L_n = \eta(\Sym^n(L_{Z/X}[-1]))$ for the morphism $\eta\colon \QCoh^\op_{\solid,\pc} \to \ICoh$ from \cref{rslt:morphisms-of-3ffs-to-ICoh}. We are thus reduced to the following observation: If $P \in \QCoh(Z)$ is perfect then $\gamma(P) = \eta(P^\vee)$ in $\ICoh(Z)$. To show this identity, note that it is an isomorphism between morphisms $\Perf^\op \to \ICoh$ of sheaves of categories, hence it can be checked on the site of affine schemes. But then it is true by construction of $\gamma$ in \cref{rslt:construct-gamma-from-QCoh-to-ICoh-on-rings}.
\end{proof}

Next we introduce a t-structure on $\ICoh(X)$ for geometric stacks $X$, similar to the one on $\QCoh(X)$. Note that if $X = \Spec A$ is an affine scheme then $\Coh(X) = \Coh(A)$ has a natural t-structure by restricting the one from $\D(A)$. By \cite[Lemma~C.2.4.3]{lurie-SAG} this induces a right-complete t-structure on $\ICoh(A) = \Ind(\Coh(A))$ such that $\ICoh^{\le0}(A)$ consists of those objects that can be written as a cofiltered system $(M_i)_i$ with $M_i \in \Coh^{\le0}(A)$ for all $i$ (and similarly for $\ICoh^{\ge0}(A)$). Unfortunately this t-structure is not compatible with $\ICoh$-pullback along flat maps, as these are computed via $!$-pullback on $\Coh$ (see \cref{rslt:ICoh-pullback-on-affines-is-upper-shriek}). We therefore use the following replacement:

\begin{lem} \label{rslt:ICoh-upper-natural-functors}
Fix a classical regular noetherian ring $\Lambda$. Let $\Stk_\Lambda^{\lafp,\mathrm{gsv}} \subset \Stk_\Lambda$ be the wide subcategory of lafp stacks with suave geometric transition maps. For every map $f\colon Y \to X$ in this category there is a functor $f^\natural\colon \ICoh(X) \to \ICoh(Y)$ with the following properties:
\begin{lemenum}
    \item \label{rslt:ICoh-upper-natural-descent} The assignment $f \mapsto f^\natural$ is natural, i.e. upgrades to a functor
    \begin{align*}
        \ICoh^\natural\colon (\Stk_\Lambda^{\lafp,\mathrm{gsv}})^\op \to \Cat, \qquad X \mapsto \ICoh(X).
    \end{align*}
    Moreover, $\ICoh^\natural$ descends along any smooth geometric cover.

    \item \label{rslt:ICoh-upper-natural-twist-of-upper-star} If $f\colon Y \to X$ is smooth then $f^\natural = f^* \tensor \gamma(\omega_f^\vee)$, where $\omega_f \in \QCoh(Y)$ is the dualizing complex for $f$ in the $\QCoh_\solid$-formalism.

    \item \label{rslt:f-upper-natural-left-adjoint-to-f-lower-shriek} If $f\colon Y \to X$ is qcqs schematic or QCA (so that $f_!$ is defined on $\ICoh$) and suave then $f^\natural$ is left adjoint to $f_!$.
\end{lemenum}
\end{lem}
\begin{proof}
If we restrict to the category of maps $f$ that are qcqs schematic or QCA, then the claim is easy: In this case we have a functor $f_!$ on $\ICoh$ and we simply define $f^\natural$ as the left adjoint, which exists because $f$ is $\ICoh$-prim by \cref{rslt:ICoh-on-lafp-schemes} and \cref{rslt:ICoh-codualizing-sheaf-for-QCA-maps}. This clearly satisfies (iii). The descent claim in (i) follows by applying \cite[Lemma~4.7.1]{heyer-mann-6ff} to $\ICoh^\op$ (where $f^\natural$ plays the role of $f^!$). To check (ii) we observe that on $\ICoh$ we have $f_! = f_* \intHom(\delta_f^{\ICoh}, -)$ by \cite[Corollary~4.5.11(ii)]{heyer-mann-6ff}, hence $f^\natural(-) = f^*(-) \tensor \delta_f^{\ICoh}$, so we conclude by the explicit computation of $\delta_f^{\ICoh}$ in the above references.

In the following we extend the definition of $f^\natural$ beyond the case where $f_!$ is defined on $\ICoh$ in order to avoid unnecessary technical restrictions in the results below. The reader is invited to skip this proof.

We consider the sheafy presentable 6-functor formalism $\IPCoh$ on $\Stk_\Lambda^\lafp$ which on affine schemes $X = \Spec A$ is given by $\IPCoh(X) = \Ind(\ICoh(A)^\op)$. To avoid set-theoretic issues one can instead take the colimit of $\Ind(\ICoh(A)^{\kappa,\op})$ for increasing regular cardinals $\kappa$, but we gloss over these details in the following. To define $\IPCoh$ on lafp stacks we implicitly use that it has fppf descent on lafp affine schemes, which follows in a similar way as for $\PQCoh_{\solid,\pc}$ using descendability (cf.\ \cref{rslt:fppf-descent-for-all-sheaf-theories}). By construction there is a natural fully faithful embedding $\ICoh^\op \injto \IPCoh$ compatible with pullbacks and tensor products. It upgrades to an embedding of 3-functor formalisms which admit $!$-functors for lafp affine maps.

We observe that suave affine maps are $\IPCoh$-suave: This can be checked on affine schemes, where it boils down to the observation that suave maps are $\ICoh$-prim (by \cref{rslt:ICoh-on-lafp-schemes}) and hence $\ICoh^\op$-suave. Hence, since suave covers are universal $!$-covers, we can inductively extend the $!$-able maps for $\IPCoh$ to the class of all geometric maps of lafp stacks and then every suave geometric map is $\IPCoh$-suave (cf.\ the proof of \cref{rslt:QCoh-solid-on-proper-and-smooth-maps}). In particular there is a functor
\begin{align*}
    \IPCoh^!\colon (\Stk_\Lambda^{\lafp,\mathrm{gsv}})^\op \to \Cat
\end{align*}
with transition maps given by upper-$!$ functors. We now claim that $\ICoh^\op \subseteq \IPCoh$ is stable under the upper-$!$ functors. Namely, for every suave geometric map $f\colon Y \to X$, by $\IPCoh$-suaveness of $f$ we have $f^! = f^* \tensor \omega_f$ as functors $\IPCoh(X) \to \IPCoh(Y)$ and by construction $\ICoh^\op$ is stable under $f^*$ and $\tensor$, so it is enough to verify that $\omega_f \in \ICoh^\op(Y)$. By descent we can reduce to the case that $X$ is an affine scheme. We now argue by induction on $n$ such that $f$ is $n$-geometric. Pick a smooth $(n-1)$-geometric cover $g\colon U \to Y$ by a scheme $U$. It is enough to show that $g^* \omega_f \in \ICoh^\op(U)$. But $g^* \omega_f = \omega_{fg} \tensor \omega_g^{-1}$ and by induction the right-hand side lies indeed in $\ICoh^\op(U)$. Altogether we naturally obtain functors $f^!$ on $\ICoh^\op$ for all maps in $\Stk_\Lambda^{\lafp,\mathrm{gsv}}$. We pass to opposite categories and rename $f^!$ to $f^\natural$ in order to arrive at the desired functor in (i).

It remains to check the claims (i), (ii) and (iii). In the same way as in the first paragraph of this proof, all of these statements hold true for $\IPCoh$. We deduce the statements for $\ICoh$, as follows. In (i) only the smooth descent remains, which now reduces to showing that containment in $\ICoh^\op \subseteq \IPCoh$ can be checked after $!$-pull along a smooth cover; but since $!$-pull is a twist of $*$-pull, this reduces to the case of $*$-pull, where it is true by construction. Part (ii) is clear (use a similar $\Ind$-extension of $\PQCoh_{\solid,\pc}$). For part (iii) we need to see that the embedding $\ICoh^\op \injto \IPCoh$ extends to an embedding of 3-functor formalisms on a class of $!$-able maps that contains the qcqs schematic or the QCA maps. But this follows in the same way as for the embedding of $\QCoh_{\solid,\pc} \injto \PQCoh_{\solid,\pc}$ in the proof of \cref{rslt:morphisms-of-3ffs-to-ICoh}.
\end{proof}

The functors $f^\natural$ often appear as $f^{*,\ICoh}$ in the literature. We can now use them to define the promised t-structure on $\ICoh$:

\begin{prop} \label{rslt:t-structure-on-ICoh}
Fix a classical regular noetherian ring $\Lambda$ and let $X$ be a lafp geometric stack over $\Lambda$. Then there exists a uniquely determined t-structure on $\ICoh(X)$ with the following properties:
\begin{propenum}
    \item If $X = \Spec A$ is an affine scheme then the t-structure on $\ICoh(X)$ is the one induced from the canonical t-structure on $\ICoh(A) = \Ind(\Coh(A))$ discussed above.

    \item \label{rslt:pullbacks-in-ICoh-t-structure} If $f\colon Y \to X$ is a suave geometric map then $f^\natural\colon \ICoh(X) \to \ICoh(Y)$ is right t-exact. If $f$ is flat then $f^\natural$ is t-exact.

    \item The t-structure on $\ICoh(X)$ is right complete and filtered colimits are t-exact.
\end{propenum}
\end{prop}
\begin{proof}
In the case that $X = \Spec A$ is affine then we equip $\ICoh(X) = \Ind(\Coh(A))$ with the $\Ind$-t-structure induced from the canonical t-structure on $\Coh(A)$. This t-structure is right complete and filtered colimits are exact (see \cite[Lemma~C.2.4.3]{lurie-SAG}). Suppose that $f\colon A \to B$ is a suave map of lafp $\Lambda$-algebras. We claim that in this case the functor $f^\natural\colon \ICoh(A) \to \ICoh(B)$ is induced by the canonical pullback functor on coherent modules via $\Ind$-extension. To see this, we pass to opposite categories and go through suave duality in order to reduce the claim to the following statement: Let $f^!\colon \PCoh(A) \to \PCoh(B)$ be the right adjoint of $f_!$; then $f^!$ is defined via $\Pro$-extension by the functor $f^!\colon \Coh(A) \to \Coh(B)$ coming from $\QCoh_\solid$. But this is straightforward to see: By construction we have $f_! \eta_B = \eta_A f_!$ as functors $\PQCoh_{\solid,\pc}(B) \to \PCoh(A)$, hence by passing to right adjoints we have $\iota_B f^! = f^! \iota_A$, which easily implies the claim. Altogether we deduce that $f^\natural\colon \ICoh(A) \to \ICoh(B)$ is right t-exact and if $f$ is flat then it is t-exact.

We have proved all claims in the case that $X$ is an affine scheme. One can now perform similar inductive steps as in the proof of \cref{rslt:t-structure-on-QCoh} in order to define the t-structure on $\ICoh$ in the claimed generality.
\end{proof}

On affine schemes $X = \Spec A$ we have $\ICoh(X) = \Ind(\Coh(X))$ by definition, so it is natural to ask if this statement is true for more general stacks $X$. Recall the morphism $\eta\colon \QCoh^\op_{\solid,\pc} \to \ICoh$ of 3-functor formalisms from \cref{rslt:morphisms-of-3ffs-to-ICoh}. Recall the definition of $\Coh(X) = \Coh_\Lambda(X) \subseteq \QCoh_{\solid,\pc}(X)$ from \cref{def:relatively-coherent-sheaves} and its explicit description in \cref{rslt:coh-over-regular-noetherian-base}. This leads to the following definition:

\begin{defn}
Let $\Lambda$ be a classical regular noetherian ring and $X$ a lafp geometric stack over $\Lambda$. We denote by
\begin{align*}
    \nu_X\colon \Ind(\Coh(X)) \to \ICoh(X)
\end{align*}
the ind-extension of the composition $\Coh(X) \isoto \Coh(X)^\op \injto \QCoh^\op_{\solid,\pc}(X) \to \ICoh(X)$, where the first map is suave duality and the last map is $\eta$. 
\end{defn}

\begin{rmk} \label{rmk:comparison-of-Coh-and-Perf-in-ICoh}
We warn the reader that the relation between $\nu\colon \Coh(X) \to \ICoh(X)$ and $\gamma\colon \Perf(X) \to \ICoh(X)$ is a bit subtle. While both functors are defined via $\eta$, they use the precompositions $\Coh(X) \isoto \Coh(X)^\op$ and $\Perf(X) \isoto \Perf(X)^\op$ respectively, and these two are not compatible even in the case that $\Perf(X)$ naturally embeds into $\Coh(X)$. The first one is suave duality, i.e. internal hom to the dualizing complex, while the second one is natural duality, i.e. internal hom to $\calO_X$.    
\end{rmk}

\begin{lem} \label{rslt:basic-properties-of-nu-on-ICoh}
Let $\Lambda$ be a classical regular noetherian ring and $X$ a lafp geometric stack over $\Lambda$. The restriction of $\nu_X$ to $\Coh(X)$ is fully faithful and hence induces an embedding
\begin{align*}
    \Coh(X) \subseteq \ICoh(X)
\end{align*}
with the following properties:
\begin{lemenum}
    \item An object $M \in \ICoh(X)$ lies in $\Coh(X)$ if and only if its pullback along every smooth geometric map $\Spec A \to X$ lies in $\Coh(A) \subseteq \ICoh(A)$.

    \item \label{rslt:Coh-in-ICoh-descends-along-suave-covers} $\Coh \subseteq \ICoh$ is stable under pullback along any suave map of lafp geometric stacks. Moreover, containment in $\Coh$ can be checked after pullback along any suave cover.

    \item \label{rslt:t-structure-on-ICoh-restricts-to-Coh} The t-structure on $\ICoh(X)$ restricts to a t-structure on $\Coh(X)$. Moreover, an object in $\ICoh^\heartsuit(X)$ lies in $\Coh^\heartsuit(X)$ as soon as it admits an injection into or a surjection from an object in $\Coh^\heartsuit(X)$.
\end{lemenum}
\end{lem}
\begin{proof}
Fix a suave map $f\colon Y \to X$ of lafp geometric stacks over $\Lambda$. Given $M \in \Coh(X) \subseteq \QCoh_{\solid,\pc}(X)$ we see by naturality of $\eta$ and the fact that $f^*$ on $\QCoh_{\solid,\pc}$ preserves $\Coh$ that
\begin{align*}
    f^* (\nu_X M) = \nu_Y(\SD_Y(f^* \SD_X(M))) = \nu_Y(f^! M),
\end{align*}
where $\SD_X$ denotes the suave duality functor $\SD = \intHom(-, \omega_X)$. This upgrades to a natural transformation $\nu\colon \Coh^! \to \ICoh^*$, where in the source we use the $\QCoh_\solid$ upper-$!$ transition maps and in the target we use $\ICoh$ upper-$*$ transition maps. Indeed, all the functors used in the definition of $\nu$ are natural in $X$ (for the suave duality functor use \cref{rslt:functoriality-of-suave-dual}).
Both $\Coh^!$ and $\ICoh^*$ descend along suave covers and hence the full faithfulness of $\Coh(X) \to \ICoh(X)$ can be reduced to the case that $X = \Spec A$ is an affine scheme, where it is clear. We have also already shown the first part of (ii), and the second part follows easily by descent of $\Coh^!$. Moreover, part (i) follows easily by (ii) and descent.

It remains to prove (iii). Since $\ICoh$-upper-$\natural$ along smooth maps is t-exact by \cref{rslt:pullbacks-in-ICoh-t-structure}, we can use (i) to reduce to the case that $X = \Spec A$ is an affine scheme (note that by \cref{rslt:ICoh-upper-natural-twist-of-upper-star}, (i) easily implies the same claim for $\natural$-pullbacks instead of pullbacks). Now fix some $M \in \ICoh^\heartsuit(A)$, i.e. $M = (M_i)_i$ is a filtered system of classical finitely generated $\pi_0 A$-modules. Suppose further that $M$ admits an injective map $\alpha\colon M \injto N$ in $\ICoh^\heartsuit(A)$ for some classical finitely generated $\pi_0 A$-module $N \in \Coh^\heartsuit(A)$. The map $\alpha$ is induced by a compatible family of maps $\alpha_i\colon M_i \to N$ and for each $i$ we let $M'_i$ be the image of $M_i$ in $N$. The $M'_i$ are finitely generated because $\pi_0 A$ is noetherian and hence they assemble into an object $M' \in \ICoh^\heartsuit(A)$ together with maps $M \to M' \injto N$. Since the composition of these maps is $\alpha$ and hence injective, it follows that $M \to M'$ is injective. It is also clearly surjective (as this can be checked termwise), so $M = M'$. We can therefore replace $M$ by $M'$ and from now on assume that all the maps $\alpha_i\colon M_i \injto N$ are injective. In other words, the $M_i$ form a filtered system of submodules of $N$. Since $\pi_0 A$ is noetherian, this filtered system is eventually constant and hence $M = M_i$ for large enough $i$. This shows that $M \in \Coh^\heartsuit(A)$, as desired. The claim about surjections follows formally from this by passing to the kernel.
\end{proof}

With the functor $\nu_X$ at hand we can now study when it is an equivalence. We start with the following descent property of this question.

\begin{lem} \label{rslt:ICoh-equal-IndCoh-descends}
Fix a classical regular noetherian ring $\Lambda$ and let $f\colon Y \surjto X$ be a qcqs suave schematic cover of lafp geometric stacks over $\Lambda$. Assume that
\begin{enumerate}[(a)]
    \item $f$ a universal $\ICoh$-$!$-cover.
    \item $\nu_Y$ is essentially surjective.
    \item $\nu_X$ is fully faithful.
\end{enumerate}
Then $\nu_X$ is an equivalence.
\end{lem}
\begin{proof}
Let $f_\bullet\colon Y_\bullet \to X$ be the \v{C}ech nerve of $f$. By (a) we have $\ICoh^!(X) = \varprojlim_{n\in\Delta} \ICoh^!(Y_n)$, where $\ICoh^!$ denotes the functor $\ICoh$ with upper-$!$ transition maps. By \cite[Lemma~D.4.7]{heyer-mann-6ff}, for every $M \in \ICoh(X)$ we have $M = \varinjlim_{n\in\Delta} f_{n!} f_n^! M$. We need to show that $M$ is in the image of $\nu_X$. Since $\nu_X$ is fully faithful and commutes with all colimits, this reduces to showing that $f_{n!} f_n^! M$ is in the image of $\nu_X$. By factoring $f_n$ over $f$ we further reduce to showing that for every $N \in \ICoh(Y)$, $f_! N$ is in the image of $\nu_X$. By (b) we can write $N$ as a colimit of objects in $\Coh(Y)$, which allows us to reduce to the case that $N \in \Coh(Y)$. In particular $N$ is bounded for the t-structure on $\ICoh$. By \cref{rslt:pullbacks-in-ICoh-t-structure} $f^\natural$ is t-exact and thus its right adjoint $f_!$ is left t-exact, so in particular $f_! N$ is left bounded. Since the t-structure on $\ICoh(X)$ is right complete we can write $f_! N$ as a colimit of its right truncations. These are bounded and can thus be decomposed using fiber sequences into static objects. Altogether we have reduced the claim to:
\begin{itemize}
    \item[($*$)] Any $M \in \ICoh^\heartsuit(X)$ lies in the image of $\nu_X$.
\end{itemize}
To prove ($*$), let $M$ be given. We follow the argument in \cite[Proposition~9.5.2.3]{lurie-SAG}. Let $N := f^\natural M$. By \cref{rslt:pullbacks-in-ICoh-t-structure} $N$ lies in $\ICoh^\heartsuit(Y)$. By (b) we can write $N = \varinjlim_i N_i$ as a filtered colimit of coherent sheaves $N_i \in \Coh(Y)$. Since the t-structure on $\ICoh(Y)$ is compatible with filtered colimits and restricts to a t-structure on $\Coh(Y)$ (see \cref{rslt:t-structure-on-ICoh-restricts-to-Coh}) we can assume that all $N_i$ lie in $\Coh^\heartsuit(Y)$. By \cref{rslt:t-structure-on-ICoh-restricts-to-Coh} we can replace each $N_i$ by its image in $N$ (inside the abelian category $\ICoh^\heartsuit(Y)$) in order to assume that all $N_i \injto N$ are injective. Altogether we have written
\begin{align*}
    N = \bigcup_i N_i \qquad \text{for $N_i \subseteq N$ with $N_i \in \Coh^\heartsuit(Y)$}.
\end{align*}
By \cref{rslt:pullbacks-in-ICoh-t-structure} and \cref{rslt:f-upper-natural-left-adjoint-to-f-lower-shriek} the functor $f^\natural$ is t-exact and left adjoint to $f_!$, so in particular $f_!$ is left t-exact. The adjunction unit provides a natural map $\alpha\colon M \to H^0 f_! N$ and the left t-exactness of $f_!$ shows that the map $H^0 f_! N_i \injto H^0 f_! N$ is injective for all $i$. We denote by $M_i \subset M$ the preimage of $H^0 f_! N_i$ under $\alpha$. Since $H^0 f_!$ preserves filtered colimits we have $M = \varinjlim_i M_i$, so in order to conclude it only remains to show that $M_i \in \Coh^\heartsuit(X)$ for all $i$. By \cref{rslt:Coh-in-ICoh-descends-along-suave-covers} it suffices to show that $f^\natural M_i$ lies in $\Coh^\heartsuit(Y)$. We have the following commuting diagram in the abelian category $\Coh^\heartsuit(Y)$:
\begin{equation*}\begin{tikzcd}
    f^\natural M_i \arrow[r] \arrow[d] & f^\natural H^0 f_! N_i \arrow[r] \arrow[d] & N_i \arrow[d] \\
    f^\natural M \arrow[r] & f^\natural H^0 f_! N \arrow[r] & N
\end{tikzcd}\end{equation*}
The composition of the bottom horizontal maps is an isomorphism and the left vertical map is injective. It follows that the composition of the top horizontal maps is injective, i.e. we have an injection $f^\natural M_i \injto N_i$. But $N_i$ is coherent, so by \cref{rslt:t-structure-on-ICoh-restricts-to-Coh} also $f^\natural M_i$ is coherent, as desired.
\end{proof}

We can now use \cref{rslt:ICoh-equal-IndCoh-descends} in order to show that $\nu_X$ is an isomorphism for a large class of stacks. Note that this behavior is different from $\QCoh$, for which it is harder to determine whether $\QCoh = \Ind(\Perf)$ (see \cite[\S3]{BenZviFrancisNadler-IntegralTransforms} for a basic treatment of when this is true).

\begin{prop} \label{rslt:ICoh-equals-IndCoh-on-nice-stacks}
Fix a classical regular noetherian ring $\Lambda$ and let $X$ be a lafp stack over $\Lambda$. Assume that one of the following is true:
\begin{enumerate}[(a)]
    \item $X$ is a qcqs scheme.
    \item $\Lambda$ is a $\Q$-algebra and $X$ is a QCA stack over $\Lambda$.
\end{enumerate}
Then $\nu_X$ induces an equivalence
\begin{align*}
    \ICoh(X) = \Ind(\Coh(X)).
\end{align*}
In particular $\Coh(X) = \ICoh(X)^\omega$ is the subcategory of compact objects. It coincides with the subcategory of prim objects over $\Lambda$.
\end{prop}
\begin{proof}
First assume we are in setting (b). Let $f\colon U \to X$ be a smooth cover by an affine scheme $U$. We want to apply \cref{rslt:ICoh-equal-IndCoh-descends} to $f$. Note that $f$ is suave and by \cref{rslt:descendable-cover-of-QCA-stack} it is descendable, hence $f$ is $\ICoh$-prim and a universal $\ICoh$-$!$-cover (as in the proof of \cref{rslt:QCA-maps-are-ICoh-fine}). This shows condition (a) of \cref{rslt:ICoh-equal-IndCoh-descends}. Moreover, $\nu_U$ is clearly an equivalence, because $U$ is an affine scheme, showing condition (b). It remains to check condition (c), i.e. that $\nu_X$ is fully faithful. By \cref{rslt:basic-properties-of-nu-on-ICoh} this is true on $\Coh$, so we are reduced to showing that all objects in $\Coh(X)$ are compact in $\ICoh(X)$. This can be directly deduced from descendability, but we choose a different argument because it also provides the last statement about prim objects. First note that by \cite[Lemma~4.4.18(ii)]{heyer-mann-6ff}, every prim object in $\ICoh(X)$ is compact, so it is enough to show that all objects in $\Coh(X)$ are prim. By \cite[Corollary~4.7.5(i)]{heyer-mann-6ff} primness can be checked after pullback along $f$, so we are reduced to showing that in the case that $X = \Spec A$ is an affine scheme, all objects in $\Coh(A)$ are prim. Now recall the morphism of 3-functor formalisms $\eta\colon \QCoh_{\solid}^\op \to \ICoh$ from \cref{rslt:morphisms-of-3ffs-to-ICoh}. Since every $P \in \Coh(A)$ lies in the image of $\eta_A$, it is enough to show that $\Coh(A)^\op \subseteq \QCoh_\solid(A)^\op$ consists of prim objects. But this follows immediately from the fact that $\Coh(A)$ are exactly the suave objects in $\QCoh_\solid(A)$ (by definition). This finishes the proof that $\ICoh(X) = \Ind(\Coh(X))$ in setting (ii).

In setting (a) we can argue in the same way: As established in the proof of \cref{rslt:6ff-for-ICoh}, finite open covers are universal $\ICoh$-$!$-covers (and even descendable, by the same argument as in the proof of \cite[Lemma~4.8.3]{heyer-mann-6ff}), so we can cover $X$ by a disjoint union of affine open subsets to reduce to the affine case.

It remains to prove the characterization of prim objects. But this is now obvious: It was already established that all prim objects are compact and all objects in $\Coh(X)$ are prim. By the equivalence $\ICoh(X) = \Ind(\Coh(X))$ we know that $\Coh(X)$ consists precisely of the compact objects in $\ICoh(X)$ (because $\Coh(X)$ is idempotent complete, i.e. closed under retracts). Thus $\Coh(X)$ consists of precisely the prim objects in $\ICoh(X)$.
\end{proof}

By \cref{rslt:ICoh-equals-IndCoh-on-nice-stacks}, for every lafp QCA stack or qcqs scheme $X$ we get a functor $\Psi_X\colon \ICoh(X) = \Ind(\Coh(X)) \to \QCoh(X)$ by ind-extending the natural inclusion $\Coh(X) \subseteq \QCoh(X)$. In \cref{rslt:psi-for-ICoh-and-QCoh-on-rings} we have already seen this functor on affine schemes, and by \cref{rslt:psi-from-ICoh-to-QCoh-on-stacks} it glues to lafp qcqs schemes. We now combine all these results to get the following description:

\begin{prop} \label{rslt:psi-and-xi-on-QCA-stacks}
Fix a classical regular noetherian ring $\Lambda$ and let $\cat C \subseteq \Stk_\Lambda^\lafp$ be the full subcategory spanned by
\begin{enumerate}[(a)]
    \item the qcqs schemes over $\Lambda$, or
    \item the QCA stacks over $\Lambda$ (in case $\Lambda$ is a $\Q$-algebra).
\end{enumerate}
Then there is a natural transformation $\Psi\colon \ICoh \to \QCoh$ of functors $\cat C \to \Cat$, i.e. for every $X \in \cat C$ there is a functor
\begin{align*}
    \Psi_X\colon \ICoh(X) \to \QCoh(X)
\end{align*}
that is compatible with lower-! functors on $\ICoh$ and pushforward functors on $\QCoh$. It satisfies the following properties:
\begin{propenum}
    \item Fix $X \in \cat C$. Using the equivalence $\ICoh(X) = \Ind(\Coh(X))$ from \cref{rslt:ICoh-equals-IndCoh-on-nice-stacks}, $\Psi_X$ is the ind-extension of the natural inclusion $\Coh(X) \injto \QCoh(X)$.
    
    \item In case (a), $\Psi$ agrees with the functor from \cref{rslt:psi-from-ICoh-to-QCoh-on-stacks}.

    \item If $f\colon Y \to X$ is a suave map in $\cat C$ then $f_!\colon \ICoh(Y) \to \ICoh(X)$ admits a left adjoint $f^\natural$ and the natural map $f^* \Psi_X \isoto \Psi_Y f^\natural$ is an isomorphism.

    \item If $X \in \cat C$ is suave then $\Psi_X$ admits a fully faithful left adjoint $\Xi_X\colon \QCoh(X) \injto \ICoh(X)$. If $f\colon Y \to X$ is a suave map then the natural map $\Xi_X f_* \isoto f_! \Xi_Y$ is an isomorphism.

    \item \label{rslt:coh-smooth-image-of-gamma-and-Xi-coincide} If $X \in \cat C$ is cohomologically smooth then the two functors $\gamma_X, \Xi_X\colon \QCoh(X) \injto \ICoh(X)$ have the same image.
\end{propenum}
\end{prop}
\begin{proof}
In case (a) $\Psi$ was already constructed in \cref{rslt:psi-from-ICoh-to-QCoh-on-stacks}. In case (b) we can construct it in a very similar way: The crucial observation is that both $\ICoh$ and $\QCoh$ descend along smooth covers of a QCA stack by an affine scheme, where the transition functors are upper-! maps (in the case of $\QCoh$ these are the right adjoints of the pushforward functors). But this follows from the descendability of smooth covers of QCA stacks (shown in \cref{rslt:descendable-cover-of-QCA-stack}) together with \cite[Lemma~4.7.4(i)]{heyer-mann-6ff}. Since $\Psi$ is constructed via a Kan extension, we see that in both (a) and (b) it is uniquely determined by its restriction to affine schemes.

We now prove the additional claims about $\Psi$. Part (ii) is clear by construction. Parts (iii) and (iv) are proved in \cref{rslt:psi-from-ICoh-to-QCoh-on-stacks} in case (a) and can be proved in a very similar way in case (b): Just replace ``finite open cover'' by ``affine smooth cover'' and observe that all the arguments still work because affine smooth covers of QCA stacks are descendable by \cref{rslt:descendable-cover-of-QCA-stack}. For part (v) we already know that $\Xi_X$ is fully faithful by (iv) and $\gamma_X$ is fully faithful by \cref{rslt:gamma-fully-faithful-for-suave}. To check that their images agree, pick a smooth cover $g\colon U \surjto X$ by an affine scheme $U$. Since $\gamma$ is compatible with pullbacks and $\Xi$ is compatible with pullback on $\QCoh$ and upper-$\natural$ on $\ICoh$ (by passing to left adjoints in the compatibility for $\Psi$) and since $g^\natural$ and $g^*$ differ by a twist and $\ICoh$ descends along both, the claim reduces to $U$ instead of $X$, i.e. we can assume that $X = \Spec A$ is an affine scheme. But then the claim reduces to the easy observation that $\gamma$ sends $\Perf(A)$ isomorphically to $\Perf(A) \subseteq \Coh(A) \subseteq \ICoh(A)$ (it is given by $- \tensor \omega_A$ and $\omega_A$ is invertible).

It remains to prove (i), so let $X \in \cat C$ be given. By construction $\Psi_X$ commutes with colimits, so it is enough to show that its restriction to $\Coh(X)$ coincides with the canonical embedding into $\QCoh(X)$. We prove this in setting (b), with (a) being very similar. Pick a smooth cover $U \surjto X$ by some affine scheme $U$ and let $g_\bullet\colon U_\bullet \to X$ denote the associated \v{C}ech nerve. Note that $\ICoh(X) = \varprojlim_{n\in\Delta} \ICoh^\natural(U_n)$ and $\QCoh(X) = \varprojlim_{n\in\Delta} \QCoh^*(U_n)$, where the transition maps are given by upper-$\natural$ and pullback functors, respectively (for $\ICoh$, the descent was shown in \cref{rslt:ICoh-upper-natural-descent}). By (iii) the functor $\Psi_X$ is given as
\begin{align*}
    \Psi_X = \Psi_{U_\bullet}\colon \varprojlim_{n\in\Delta} \ICoh^\natural(U_n) \to \varprojlim_{n\in\Delta} \QCoh^*(U_n).
\end{align*}
It is therefore enough to show that the restriction of the functor $\Psi_{U_n}$ to $\Coh(U_n)$ coincides with the natural embedding into $\QCoh(U_n)$, naturally in $n$. Since all $U_n$ are affine, the first part of this claim was shown in \cref{rslt:construct-psi-for-ICoh-and-QCoh-on-rings}, so we only need to prove functoriality in $n$. The embedding $\Coh \subset \ICoh$ is given by $\nu$ and it was shown in the proof of \cref{rslt:Coh-in-ICoh-descends-along-suave-covers} that $\nu$ is functorial with respect to $f^!$ on $\Coh$ and $f^*$ on $\ICoh$. To get a similar functoriality for $f^\natural$, note that since $f^\natural$ is left adjoint to $f_!$ on $\ICoh$, it corresponds to the usual $f^!$ on $\PQCoh_{\solid,\pc}$ via $\eta\colon \PQCoh_{\solid,\pc}^\op \to \ICoh$ (because the $(-)^\op$ swaps the direction of adjunction). Using the functorial equivalence $\Coh^\op = \Coh$ from suave duality (see the proof of \cref{rslt:Coh-in-ICoh-descends-along-suave-covers}) we arrive at the natural transformation $\nu\colon \Coh^* \injto \ICoh^\natural$ on $\cat C$, where we restrict to suave transition maps. A quick inspection of the construction in \cref{rslt:psi-for-ICoh-and-QCoh-on-rings} reveals that on affine schemes and suave maps, the composition $\Psi \comp \nu\colon \Coh^* \to \QCoh^*$ is indeed given by the natural inclusion, as desired.
\end{proof}

\begin{rmk}
We expect $\Psi$ and $\Xi$ to be $\QCoh$-linear, where the $\QCoh$-linear structure on $\ICoh$ is induced from $\gamma$. Moreover, in the setting of \cref{rslt:coh-smooth-image-of-gamma-and-Xi-coincide} we expect there to be a canonical equivalence $\gamma_X(-) = \Xi_X(-) \tensor \gamma(\omega_X)$. We do not need these statements in this paper and hence do not pursue them.
\end{rmk}

\subsection{Admissible ind-coherent sheaves} \label{sec:adm-sheaves}

In \cref{rslt:ICoh-equals-IndCoh-on-nice-stacks} we have characterized the $\ICoh$-prim sheaves on nice stacks: They are exactly the coherent ones. We now introduce the $\ICoh$-suave sheaves, which are much more exotic and -- to our knowledge -- have not been studied before.

\begin{defn}\label{def:admissibleIcoh}
Fix a classical regular noetherian ring $\Lambda$ and an $\ICoh$-fine lafp stack $X$ over $\Lambda$. We say that a sheaf $M \in \ICoh(X)$ is \emph{admissible} if it is suave over $\Spec \Lambda$ (in the sense of \cite[Definition~4.4.1]{heyer-mann-6ff}). We denote by
\begin{align*}
    \Adm(X) \subseteq \ICoh(X)
\end{align*}
the full subcategory spanned by the admissible sheaves. Denoting by $f\colon X \to \Spec \Lambda$ the structure map, we define the contravariant endofunctor
\begin{align*}
    \Dadm\colon \ICoh(X)^\op \to \ICoh(X), \qquad M \mapsto \intHom(M, f^! 1).
\end{align*}
and call it the \emph{admissible duality functor}. By \cite[Lemmas~4.4.4,~4.4.5]{heyer-mann-6ff} it restricts an equivalence $\Adm(X)^\op \isoto \Adm(X)$.
\end{defn}

\begin{defn}
Fix a classical regular noetherian ring $\Lambda$ and a lafp stack $X$ over $\Lambda$. We denote by
\begin{align*}
    \Dgs\colon \Coh(X)^\op \isoto \Coh(X)
\end{align*}
the suave duality functor in $\QCoh_\solid$, given as $\intHom(-, f^! 1)$. Note that if $X$ is as in \cref{rslt:ICoh-equals-IndCoh-on-nice-stacks} then by the proof of that result, $\Coh(X)$ identifies with the prim objects in $\ICoh(X)$ and $\Dgs$ is prim duality in $\ICoh(X)$. Recall from \cref{rslt:Coh-duality-is-t-bounded} that $\Dgs$ is uniformly t-bounded.
\end{defn}

We caution the reader that admissibility and coherence are often unrelated. We also emphasize that the formula $\Dadm = \intHom(-, f^! 1)$ uses the six functors of $\ICoh$, which are usually not explicit and are rarely considered in the literature. For example, recall that the functor $f^{!,\ICoh}$ from the literature is denoted $f^*$ in this paper. Let us start with a few examples:

\begin{prop}
Fix a field $k$ and an $\ICoh$-fine lafp stack $X$ over $k$.
\begin{propenum}
    \item Suppose $X = \Spec A$ for a classical finitely generated $k$-algebra $A$. Then $\Adm(X)$ consists precisely of the perfect complexes supported on a finite-length closed subscheme of $X$. Here we implicitly use the embedding $\gamma\colon \Perf(X) \injto \ICoh(X)$ from \cref{rslt:gamma-fully-faithful-on-QCoh-minus}.

    \item Suppose $X$ is a finite-type separated scheme over $k$. Then $X$ is smooth and proper if and only if $\Coh(X) \subseteq \Adm(X)$ inside $\ICoh(X)$.

    \item Suppose $k$ has characteristic $0$ and $X = */G$ for a smooth linear algebraic group $G$ over $k$. Then the regular representation $k[G]$ regarded as an object of $\QCoh(*/G) = \ICoh(*/G)$ is admissible. If $\dim G > 0$, it is not coherent.
\end{propenum}
\end{prop}
\begin{proof}
We first prove (i), so let $X = \Spec A$ be as in the claim. Let $M \in \Adm(X)$ be arbitrary. First note that the diagonal of $X$ is proper and hence $\ICoh$-suave by \cref{rslt:ICoh-on-lafp-schemes}. Thus by \cite[Corollary~4.5.18(ii)]{heyer-mann-6ff} $M$ is dualizable in $\ICoh(X)$ and hence lies in the image of $\Perf(X) \injto \ICoh(X)$ by \cref{rslt:Perf-equals-dualizable-objects-in-ICoh}; let us thus view $M \in \Perf(X)$ from now on. If $f\colon X \to *$ denotes the structure map of $X$ then $f$ is suave and hence by \cite[Lemma~4.5.16(ii)]{heyer-mann-6ff} $f_*\colon \ICoh(X) \to \ICoh(*)$ preserves admissible sheaves. Since $\gamma$ is compatible with pushforward (see \cref{rslt:gamma-from-QCoh-to-ICoh-on-stacks}) we altogether deduce that $M$ is a perfect $A$-module which is also perfect as a $k$-module. Equivalently, $M$ is a perfect complex on $X$ which is supported on a finite-length closed subscheme (see \cite{Tang-MasterThesis}).

The converse direction can be deduced from \cref{rslt:admissible-sheaf-criterion-via-coherent} below. 

For (ii), if $X$ is smooth and proper, the claimed containment follows easily from \cref{rslt:admissible-sheaf-criterion-via-coherent} below. For the converse direction, one uses some results from \cite{GLproperthing}. We leave the details to the reader.

To prove (iii), note that the regular representation is obtained as $k[G] = g_* 1$ for $g\colon * \to */G$ the projection (here we use implicitly that $\gamma$ is compatible with $g_*$ by \cref{rslt:gamma-from-QCoh-to-ICoh-on-stacks}). Since $G$ is smooth, $g$ is $\ICoh$-prim, hence $g_*$ preserves admissible objects by \cite[Lemma~4.5.16(ii)]{heyer-mann-6ff}. This shows that $g_* 1$ is indeed admissible, and it is also clear that it is not coherent if $\dim G > 0$.
\end{proof}





In this paper we are only interested in admissible sheaves on (disjoint unions of) QCA stacks in characteristic $0$, so let us focus on this case now. Recall that QCA stacks are $\ICoh$-fine by \cref{rslt:QCA-maps-are-ICoh-fine}. We get the following alternative characterization of admissibility, as well as several basic properties.

\begin{prop}\label{prop:admissiblefacts}Fix a lafp QCA stack $X$ over a classical regular noetherian $\Q$-algebra $\Lambda$.
\begin{propenum}
    \item \label{rslt:admissible-sheaf-criterion-via-coherent} A sheaf $\mathcal F \in \ICoh(X)$ is admissible if and only if for every coherent sheaf $\mathcal Q \in \Coh(X)$ the complex $\Hom(\mathcal Q, \mathcal F) \in \D(\Lambda)$ is perfect.

    \item \label{rslt:duality-exchange-formula} For all $\mathcal F \in \ICoh(X)$ and $\mathcal G \in \Coh(X)$ there is a natural isomorphism
    \begin{align*}
        \Hom(\mathcal G, \Dadm \mathcal F) = \Hom(\Dgs \mathcal G, \mathcal F)^\vee
    \end{align*}
    in $\D(\Lambda)$.
    
    \item \label{rslt:admissible-duality-reflects-suaveness} Suppose $\Lambda = k$ is a field. Then a sheaf $\mathcal F \in \ICoh(X)$ is admissible if and only if $\Dadm \mathcal F$ is admissible.

    \item \label{rslt:Dadm-commutes-with-lower-shriek-up-to-twist} If $f\colon Y \to X$ is a cohomologically smooth morphism of lafp QCA stacks with relative dualizing complex $\omega_f \in \Perf(X)$, there is a canonical isomorphism
    \begin{align*}
        \Dadm f_!(-) = f_!(\gamma(\omega_f)^\vee \tensor \Dadm(-))
    \end{align*}
    of functors $\ICoh(Y) \to \ICoh(X)$.

    \item For $f\colon Y \to X$ a proper schematic morphism of QCA stacks, there is a natural isomorphism
    \begin{align*}
        f^* \Dadm(-) = \Dadm(f^* -).
    \end{align*}
    of functors $\ICoh(X) \to \ICoh(Y)$.
\end{propenum}
\end{prop}
\begin{proof}
By \cref{rslt:ICoh-equals-IndCoh-on-nice-stacks} we know that $\Coh(X)$ are exactly the prim objects in $\ICoh(X)$ and that they form compact generators; the same is true on $X \times X$. Thus (i) follows from \cite[Corollary 4.4.15]{heyer-mann-6ff}.

For (ii), let $f_X\colon X \to *$ denote the projection. We apply \cite[Lemma~4.4.17(ii)]{heyer-mann-6ff} to $\Dgs \mathcal G$ and obtain
\begin{align*}
    \Hom(\Dgs \mathcal{G},\mathcal{F}) = f_!(\mathcal{G} \tensor \mathcal{F})
\end{align*}
in $\D(\Lambda)$. By passing to duals, this immediately implies the identity in (ii).

For (iii), we use the criterion in (i) to detect admissibility. Then the claim follows from (ii) together with the easy facts that $\Dgs$ is an equivalence on coherent sheaves and the naive $k$-linear dual $(-)^\vee$ preserves and detects perfect complexes of $k$-vector spaces.

For (iv), first note that $f$ is automatically a QCA map by \cref{rslt:stability-of-QCA-maps}. By easy manipulations we get $\Dadm f_!(-) = f_* \intHom(-, f_Y^! 1)$, where $f_Y\colon Y \to *$ is the projection. Now apply $f_* = f_!(\gamma(\omega_f^\vee) \tensor -)$ from \cref{rslt:ICoh-codualizing-sheaf-for-QCA-maps} to conclude.

Part (v) follows easily from the fact that $f$ is $\ICoh$-smooth with trivial dualizing complex (by \cref{rslt:ICoh-on-lafp-schemes}), i.e. $f^* = f^!$.
\end{proof}

\begin{rmk}
The formula in \cref{rslt:duality-exchange-formula} is extremely useful. We will often refer to it as the \emph{duality exchange formula}.
\end{rmk}

\subsection{Spectral temperization} \label{sec:spectral-temperization}


In this section we prove the following surprising theorem.\footnote{In the early stages of this project, we proved a weaker result which required the sheaf $\mathcal{F}$ to be coherent and admissible. We heartily thank Sam Raskin for suggesting that admissibility should not be needed.}

\begin{thm}\label{rslt:spectral-temperization}
Let $X$ be a lafp QCA stack over a classical regular noetherian $\Q$-algebra $\Lambda$. Assume that $X$ is cohomologically smooth. Then for every $\mathcal{F} \in  \ICoh^+(X)$, $\Dadm\mathcal{F}$ lies in the essential image of $\gamma\colon \QCoh(X) \injto \ICoh(X)$.
\end{thm}

\begin{rmks}
\begin{rmksenum}
    \item Recall from \cref{def:suave-and-coh-smooth-maps} and \cref{rmk:suaveness-in-the-literature}  that $X$ is cohomologically smooth if and only if it is truncated (i.e. admits a smooth cover by an affine scheme attached to a truncated ring) and has invertible dualizing sheaf. This condition also appears as \emph{Gorenstein} in the literature.
    
    \item It follows from \cref{rslt:basic-properties-of-nu-on-ICoh} that all coherent sheaves in $\ICoh(X)$ are bounded and in particular lie in $\ICoh^+(X)$, so \cref{rslt:spectral-temperization} applies to them.

    \item Recall the adjoint functors
    \begin{align*}
        \Xi_X\colon \QCoh(X) \leftrightarrows \ICoh(X) \noloc \Psi_X
    \end{align*}
    from \cref{rslt:psi-and-xi-on-QCA-stacks}. Under the cohomological smoothness assumption the essential images of $\gamma_X$ and $\Xi_X$ coincide (see \cref{rslt:coh-smooth-image-of-gamma-and-Xi-coincide}), so the conclusion is trivially equivalent to the claim that $\Dadm \mathcal{F}$ lies in the essential image of $\Xi$. This is the form of the result we will use later in the paper.
\end{rmksenum}
\end{rmks}

\begin{rmk}
We discovered \cref{rslt:spectral-temperization} experimentally, after proving it by a brute force calculation for the $L$-parameter stack for $\mathrm{PGL}_2$. There is a rough analogy with a theorem of Beraldo asserting that on $\Bun_G$ of a smooth projective curve, the functor $\mathrm{ps-id}_{\ast}$ sends compact $D$-modules into tempered $D$-modules. However, we only noticed this analogy after the fact, and we still don't see a ``philosophical'' reason why the general Theorem \ref{rslt:spectral-temperization} should be true.    
\end{rmk}


To prove \cref{rslt:spectral-temperization}, we will need the following general property of functors on $\Ind$-categories. For simplicity we restrict the result to stable categories.

\begin{lem} \label{rslt:adjunction-on-Ind}
Let $F\colon \cat C \injto \cat D$ be an exact fully faithful functor of stable categories and let $F'\colon \Ind(\cat C) \to \Ind(\cat D)$ be its $\Ind$-extension. Then $F'$ is fully faithful and an object $N \in \Ind(\cat D)$ lies in the essential image of $F'$ if and only if for all $Q \in \cat D$ the natural map 
\begin{align*}
    \varinjlim_{P \in (\cat C_{Q/})^{\op}} H^0\Hom(P,N) \isoto H^0\Hom(Q, N).
\end{align*}
is an isomorphism. Here all $\Hom$'s are considered in $\Ind(\cat D)$.
\end{lem}
\begin{proof}
We consider the natural embeddings $\Ind(\cat C) \subseteq \cat P(\cat C)$ and $\Ind(\cat D) \subseteq \cat P(\cat D)$, where $\cat P(-)$ denotes the Yoneda category. Then left Kan extension along $F\colon \cat C^\op \to \cat D^\op$ produces a functor $\cat D^\op \to \Ani$ out of a functor $\cat C^\op \to \Ani$, i.e. we get a map $\cat P(\cat C) \to \cat P(\cat D)$. The functor $F'$ is the restriction of this map to $\Ind(\cat C)$ (which automatically factors over $\Ind(\cat D)$ because it preserves all small colimits). By the pointwise formula for left Kan extensions, we deduce that for every $M \in \Ind(\cat C)$ and every $Q \in \cat D$
we have 
\begin{align*}
    \Hom(Q, F'(M)) = \varinjlim_{P \in (\cat C_{Q/})^\op} \Hom(P, M),
\end{align*}
where $\Hom$ on the left is in $\Ind(\cat D)$ and the $\Hom$ on the right is in $\Ind(\cat C)$; both $\Hom$'s are understood as anima. Here we implicitly identify $M$ with a functor $\cat C^\op \to \Ani$ sending $P$ to $\Hom(P, M)$.

Let us denote by $R\colon \cat P(\cat D) \to \cat P(\cat C)$ the restriction along $F$, which restricts to a functor $\Ind(\cat D) \to \Ind(\cat C)$ that is right adjoint to $F'$. By the above formula for $F'(M)$ we see immediately that $RF'(M) = M$ for all $M \in \Ind(\cat C)$. Moreover, an object $N \in \Ind(\cat D)$ lies in the essential image of $F'$ if and only if the natural map $F'R(N) \to N$ is an isomorphism. This can be checked after applying $\Hom(Q, -)$ for all $Q \in \cat D$, so we end up with the condition that 
\begin{align*}
    \varinjlim_{P \in (\cat C_{Q/})^{\op}} \Hom(P,N) \isoto \Hom(Q, N).
\end{align*}
is an isomorphism for all $Q \in \cat D$. Here all $\Hom$'s are understood as anima.

It remains to replace $\Hom$ with $H^0\Hom$ in the above condition. For this it is enough to show that the category $(\cat C_{Q/})^\op$ is filtered, because then taking $H^0$ commutes with this colimit (and by applying shifts one easily deduces the same isomorphism for $H^n\Hom$ for all $n$). To show that $(\cat C_{Q/})^\op$ is filtered, we need to verify that every finite diagram $(P_{i})_{i\in I}$ in $(\cat C_{Q/})^\op$ can be extended to a diagram over $I^{\triangleright}$ (see \cite[Definition~5.3.1.7]{lurie-higher-topos-theory}). For this it is enough to know that the diagram has a colimit, so we only need to show that $(\cat C_{Q/})^\op$ has finite colimits. Passing to opposite categories this reduces to showing that $\cat C_{Q/}$ has finite limits. This follows immediately from the fact that $\cat C$ has finite limits (see \cite[Proposition~1.2.13.8]{lurie-higher-topos-theory}).
\end{proof}

With \cref{rslt:adjunction-on-Ind} at hand we can now prove a version of \cref{rslt:spectral-temperization} on affine schemes in any characteristic.

\begin{prop} \label{rslt:spectral-temperization-on-affine-schemes}
Let $X = \Spec A$ be an afp affine scheme over a classical regular noetherian ring $\Lambda$. Assume that $X$ is cohomologically smooth over $\Lambda$. Then for every $\mathcal{F} \in  \ICoh^+(X)$, $\Dadm\mathcal{F}$ lies in the essential image of $\gamma\colon \QCoh(X) \injto \ICoh(X)$.
\end{prop}
\begin{proof}
Fix $\mathcal F \in \ICoh^+(X)$. The functor $\gamma\colon \QCoh(A) \injto \ICoh(A)$ is given by $P \mapsto P \tensor \omega_A$ for perfect $P$, where $\tensor$ denotes the usual tensor product in $\QCoh(A)$. We apply \cref{rslt:adjunction-on-Ind} to $\gamma$, which shows that $\Dadm \mathcal F$ lies in the essential image of $\gamma$ if and only if for any fixed $Q \in \Coh(A)$ the natural map
\begin{align*}
    \varinjlim_{Q \to P} H^0\Hom(P, \Dadm(\mathcal{F})) \overset{\sim}{\to} H^0\Hom(Q,\Dadm(\mathcal{F})). \tag{$*$}
\end{align*}
is an isomorphism of classical $\Lambda$-modules. Here the colimit on the left is taken over $(\Perf(A)_{Q/})^\op$ and all $\Hom$'s are computed in $\ICoh(A)$. Note that a direct application of \cref{rslt:adjunction-on-Ind} would involve the term $P \tensor \omega_A$ instead of $P$ everywhere, but since $\omega_A$ is invertible we can drop it from the claim. We fix $Q$ from now on. In order to prove ($*$), we make the following auxiliary observation:
\begin{itemize}
    \item[($**$)] For every $n \ge 0$ there is a perfect $A$-module $P_n \in \Perf(A)$ with a map $\alpha_n\colon Q \to P_n$ such that the cofiber $K_n := \cofib(Q \to P_n)$ lives in cohomological degrees $\ge n$.
\end{itemize}
To prove ($**$), first recall that $\Dgs\colon \Coh(A)^\op \isoto \Coh(A)$ is t-bounded (see \cref{rslt:Coh-duality-is-t-bounded}) so there is some constant $C \ge 0$ such that $\Dgs$ sends the heart to $\D^{[-C,C]}(A)$. We now apply \cref{rslt:approximation-of-pseudocoh-by-perfect} to the pseudocoherent $A$-module $Q$ in order to find a perfect $A$-module $P'_n$ together with a map $P'_n \to \Dgs Q$ whose fiber is concentrated in degrees $\le -n - C$. After applying $\Dgs$ we obtain the desired map $Q \to P_n$ for $P_n := \Dgs P_n'$. Note that $P_n$ is still perfect because $X$ is cohomologically smooth, so that $\omega_A$ is invertible.

With ($**$) at hand, we can now prove the isomorphism in ($*$). Let us first show that the map is surjective. For $n \ge 0$ choose $\alpha\colon Q \to P_n$ and $K_n = \cofib(Q \to P_n)$ as in ($**$). Then by duality exchange (see \cref{rslt:duality-exchange-formula}) we have
\begin{align*}
    \Hom(K_n,\Dadm(\mathcal{F})) = \Hom(\Dgs K_n,\mathcal{F})^\vee.
\end{align*}
in $\D(\Lambda)$. Since $\mathcal F$ is left-bounded and $\Dgs$ is t-bounded, we can pick $n \gg 0$ large enough so that $\Hom(\Dgs K_n, \mathcal F)$ lives in degrees $> 0$. Then $\Hom(K_n, \Dadm(\mathcal{F}))$ lives in degrees $< 0$, which implies that the natural map
\begin{align*}
    H^0\Hom(P_n, \Dadm(\mathcal{F})) \isoto H^0\Hom(Q,\Dadm(\mathcal{F}))
\end{align*}
is an isomorphism. This shows that the map in ($*$) is indeed surjective.

We now check that the map in ($*$) is injective. By the proof of \cref{rslt:adjunction-on-Ind} the colimit on the left-hand side of ($*$) is filtered, so the desired injectivity is equivalent to the claim that if we are given maps $Q \to P \to \Dadm(\mathcal{F})$ with $P$ perfect such that the composite map $Q \to \Dadm(\mathcal{F})$ is zero, then we can further refine the given data to a sequence of maps
\begin{align*}
    Q \to P' \to P \to \Dadm(\mathcal{F})
\end{align*}
where $P'$ is perfect and the composite map $P' \to \Dadm(\mathcal{F})$ is already zero. For some $n \ge 0$ we choose $\alpha_n\colon Q \to P_n$ as in ($**$). Then by the proof of the surjectivity of ($*$) we can choose $n$ large enough so that the (zero!) map $Q \to \Dadm(\mathcal{F})$ factors uniquely over the zero map $P_n \to \Dadm(\mathcal{F})$. It remains to show that the given map $g\colon Q \to P$ factors uniquely as $Q \to P_n \to P$ for large enough $n$. In fact, we claim that for $n \gg 0$ the map
\begin{align*}
    H^0\Hom(P_n, P) \isoto H^0\Hom(Q, P)
\end{align*}
is an isomorphism, equivalently that $\Hom(K_n, P)$ is concentrated in degrees $<0$. Since $\gamma$ is fully faithful, this $\Hom$ can be computed in $\QCoh$. Since $P$ is perfect, we see $\Hom(K_n, P) = K_n^\vee \tensor P$. But since $X$ is cohomologically smooth, $K_n^\vee$ differs from $\Dgs K_n$ only by a twist by the invertible $A$-module $\omega_A$ and thus the t-boundedness of $\Dgs$ implies the t-boundedness of $(-)^\vee$. Thus for $n \gg 0$ we see that $K_n^\vee$ lives in degrees $\ll 0$. This implies that indeed $K_n^\vee \tensor P$ lives in degrees $< 0$ for $n \gg 0$, as desired.
\end{proof}

With the affine case handled in \cref{rslt:spectral-temperization-on-affine-schemes} we can now use the strong descent statements for QCA stacks in order to deduce \cref{rslt:spectral-temperization}:

\begin{proof}[Proof of \cref{rslt:spectral-temperization}]
Let $U \to X$ be a surjective smooth map from an affine scheme, and let $U_n$ be the $n$-fold \v{C}ech nerve with $f_n\colon U_n \to X$ the evident smooth map. By \cref{rslt:descendable-cover-of-QCA-stack}, any $\mathcal{F} \in \ICoh^+(X)$ can be obtained by finite colimits and retracts from objects of the form $f_{n*} f_n^* \mathcal{F}$ for $n \ge 0$. Since $\gamma_X$ is fully faithful (see \cref{rslt:gamma-fully-faithful-for-suave}), it is enough to show that each $\Dadm(f_{n*} f_n^* \mathcal F)$ lies in the essential image of $\gamma$.

Since $f_n$ is smooth, we deduce from \cref{rslt:ICoh-on-lafp-schemes} that $f_{n*}$ and $f_{n!}$ differ only by a twist, and thus \cref{rslt:Dadm-commutes-with-lower-shriek-up-to-twist} reduces the claim to showing that $f_{n*} \Dadm(f_n^* \mathcal F)$ lies in the essential image of $\gamma$. By \cref{rslt:gamma-from-QCoh-to-ICoh-on-stacks} this further reduces to showing that $\Dadm(f_n^* \mathcal F)$ lies in the essential image of $\gamma$. By \cref{rslt:pullbacks-in-ICoh-t-structure} we see that $f^*$ is left t-bounded, hence $f_n^* \mathcal F \in \ICoh^+(U_n)$. On the other hand, since $X$ is QCA and cohomologically smooth, all $U_n$ are affine and cohomologically smooth, hence \cref{rslt:spectral-temperization-on-affine-schemes} applies to $U_n$. This shows that $\Dadm(f_n^* \mathcal F)$ does indeed lie in the essential image of $\gamma$, as desired.
\end{proof}

\begin{rmk}
Aside from its applications in this paper, \cref{rslt:spectral-temperization} can notably be used to compute admissible duals. For instance, suppose $X = U/G$ for a Gorenstein quasiprojective variety $U$ and a reductive group $G$, over a characteristic zero base field. Let $\kappa_X \in \QCoh(X)$ be the object representing the functor $\mathcal{F} \mapsto \Gamma(X,\mathcal{F})^{\vee}$ on $\QCoh(X)$. Then for any $\mathcal{F} \in \Coh(X)$,
\[\Dadm \mathcal{F} = \Xi(\kappa_X \otimes \Psi \Dgs \mathcal{F})\]
in $\ICoh(X)$. The details of this result and some related applications can be found in \cite{Hserrething}.
\end{rmk}

\subsection{A contraction principle} \label{sec:contraction}

The objective of the present subsection is to study quasi-coherent sheaves on stacks that come equipped with an action by the multiplicative monoid $\mathbb A^1$. The observations made here play an important role in the study of spectral Eisenstein and constant term functors in \cref{sec:spectral} and are independent from most of the subsections above. Let us start with the general setup:

\begin{defn}
Fix a ring $\Lambda$. An $\mathbb A^1$-retract over $\Lambda$ is a pair $(X, X^0)$ consisting of a stack $X \in \Stk_\Lambda$ equipped with an action by the multiplicative monoid $\mathbb A^1$, and the retract $X^0$ of the idempotent map $X \to X$ given by multiplication by $0 \in \mathbb A^1$. We denote by
\begin{align*}
    i\colon X^0 \rightleftarrows X \noloc \pi
\end{align*}
the canonical maps. We say that $(X, X^0)$ is \emph{QCA} if $\pi$ is QCA (and in particular $\Lambda$ is a $\Q$-algebra). Note that then $i$ is also QCA (cf.\ \cref{rslt:stability-of-QCA-maps}).
\end{defn}

\begin{rmk}
The QCA assumption is only needed to guarantee base-change and projection formula for pushforward of quasicoherent sheaves (see \cref{rslt:QCA-maps-have-base-change-and-projection-formula}). The following arguments thus apply to any situation where base-change and projection formula along $\pi$ and $i$ (and potentially their base-changes) are satisfied, e.g. if $\pi$ and $i$ are qcqs schematic.
\end{rmk}

Given an $\mathbb A^1$-retract $(X, X^0)$ there is a natural morphism $\pi_* = i^* \pi^* \pi_* \to i^*$ of functors $\QCoh(X) \to \QCoh(X^0)$ and the persistent theme of this subsection is the question how close this map is to being an isomorphism. Note that both $\pi$ and $i$ are $\mathbb A^1$-equivariant and in particular $\mathbb G_m$-equivariant maps, where $\mathbb A^1$ acts trivially on $X^0$; we will denote the induced maps
\begin{align*}
    i\colon X^0/\mathbb G_m \rightleftarrows X/\mathbb G_m \noloc \pi
\end{align*}
by the same letters. We observe that $X^0/\mathbb G_m = X^0 \times */\mathbb G_m$. The following result is a well-known special case of Cartier duality (see e.g. \cite[Proposition~3.1.3]{RC-CartierDuality} and \cite[Example~10.2]{Scholze-Gestalten}) and computes the category of quasi-coherent sheaves on $X^0 \times */\mathbb G_m$. In that result, by a split diagonalizable group $T$ over a ring $\Lambda$ we mean the base-change of a split diagonalizable group from $\Z$ to $\Lambda$. We denote by $X^*(T) = \Hom(T, \mathbb G_m)$ its character group.

\begin{lem} \label{rslt:Cartier-duality-for-torus}
Let $X$ be a stack over some ring $\Lambda$ and let $T$ be a split diagonalizable group over $\Lambda$. Then there is a canonical symmetric monoidal isomorphism
\begin{align*}
    \QCoh(X \times */T) = \prod_{\chi \in X^*(T)} \QCoh(X), \qquad M \mapsto (M_\chi)_\chi
\end{align*}
Here the symmetric monoidal structure on the left-hand side is induced by convolution with respect to the multiplication of $*/T$ and the symmetric monoidal structure on the right-hand side is given by componentwise tensor product.
\end{lem}
\begin{proof}
Denote $f\colon X \to X \times */T$ the natural projection. Then $f_* 1 = \bigoplus_\chi \calO(\chi)$ for line bundles $\calO(\chi) \in X \times */T$ whose underlying sheaf on $X$ is trivial and where $T$ acts via $\chi$: By base-change, this reduces to the case $\Lambda = \Z$ and $X = *$, in which case it is a well-known computation. Each $\calO(\chi)$ is an idempotent algebra for the convolution monoidal structure (which can again be checked via base-change). Thus acting on $\calO(\chi)$ defines a symmetric monoidal functor $\QCoh(X) \to \QCoh(X \times */T)$, whose right adjoint is automatically lax symmetric monoidal and one checks that it is even strictly symmetric monoidal. Combining these functors for all $\chi$ yields a symmetric monoidal functor
\begin{align*}
    \QCoh(X \times */T) \to \prod_{\chi \in X^*(T)} \QCoh(X).
\end{align*}
To check that this functor is an isomorphism, we consider its left adjoint $(M_\chi)_\chi \mapsto \bigoplus_\chi (M_\chi \tensor \mathcal O(\chi))$ and show that both unit and counit of the adjunction are an isomorphism. This can be checked on a cover of $X$, so we can assume that $X = \Spec A$ is an affine scheme. Via the forgetful functor we can further reduce to the case $X = \Spec \Z$, where it is an easy explicit computation.
\end{proof}

\begin{rmk} \label{rmk:sign-convention-for-torus-reps}
We warn the reader that there is a subtle sign change implicit in the proof of \cref{rslt:Cartier-duality-for-torus}. To explain this, for simplicity take $T = \mathbb G_m$, $X = *$ and $\Lambda = \Z$ and consider the algebraic $\mathbb G_m$-representation $f_* \Z = \bigoplus_{n\in\Z} \calO(n)$. We can explicitly describe the underlying abelian group as $\Z[t, t^{-1}] = \bigoplus_{n\in \Z} \Z t^n$. With this notation we have $\calO(n) = \Z t^{-n}$ for all $n \in \Z$, i.e. there is an inversion of sign. Indeed, by \cite[\S2.3.5]{Springer-LinAlgGrp} the $\mathbb G_m$-action on $f_! \Z = \calO(\mathbb G_m) = \Hom(\mathbb G_m, \mathbb A^1)$ is given by the formula $(gf)(x) = f(g^{-1} x)$ for $f \in \Hom(\mathbb G_m, \mathbb A^1)$ and $g, x \in \mathbb G_m$. This is consistent with $\calO(1)$ being the canonical line bundle on $*/\mathbb G_m$.
\end{rmk}

\begin{rmk} \label{rslt:compute-G-m-graded-pieces-using-retracts}
In the setting of \cref{rslt:Cartier-duality-for-torus} for $T = \mathbb G_m$, with $f\colon X \to X \times */\mathbb G_m$ the projection, one sees easily that the adjunction between $f^*$ and $f_*$ on $\QCoh$ is comonadic and hence identifies $\QCoh(X \times */\mathbb G_m)$ with the category of comodules under the coalgebra $\calO[t, t^{-1}]$ in $\QCoh(X)$. Thus, given some $\underline M \in \QCoh(X \times */\mathbb G_m)$ with underlying sheaf $M = \gamma^* \underline M \in \QCoh(X)$ there is a natural map
\begin{align*}
    M \to M \tensor \calO[t, t^{-1}] = \bigoplus_{n\in\Z} M t^n,
\end{align*}
which is part of a $\calO[t, t^{-1}]$-comodule structure on $M$. For each $n \in \Z$ we get an induced map $\alpha_n\colon M \to M t^n = M$ by composing the above map with the projection to the $n$-th factor. We claim that $\alpha_n$ is idempotent and the associated retract is the graded piece $\underline M_{-n}$ (note the sign!).

Indeed, the above map comes from the adjunction map $\underline M \to f_* f^* \underline M = \bigoplus_n M(n)$ by applying $f^*$; here we implicitly equip $M$ with the trivial $\mathbb G_m$-action and denote $M(n) = M \tensor \calO(n)$. Via the equivalence in \cref{rslt:Cartier-duality-for-torus} this adjunction map corresponds to the adjunction map $\underline M \to g^* g_\natural \underline M$, where $g\colon X \times \Z \to X$ is the projection and $g_\natural$ denotes the left adjoint of $g^*$. From this description it is clear that the induced maps $\underline M \to M(n)$ are given by the composition $\underline M \to \underline M_n \to M(n)$. After applying $f^*$ and observing that $M(n) = M t^{-n}$ (see \cref{rmk:sign-convention-for-torus-reps}) we arrive at the claimed description of $\alpha_n$.
\end{rmk}

By \cref{rslt:Cartier-duality-for-torus} we have a natural $\Z$-grading on $\QCoh(X^0/\mathbb G_m)$ for every $\mathbb A^1$-retract $(X, X^0)$. The next result studies the graded pieces of the sheaf $\pi_* 1$:

\begin{prop} \label{rslt:grading-on-cohomology-of-A1-retract}
Let $(X, X^0)$ be a QCA $\mathbb A^1$-retract over $\Q$. Then for the grading on the sheaf $\pi_* 1 \in \QCoh(X^0/\mathbb G_m)$ we have:
\begin{propenum}
    \item $(\pi_* 1)_n = 0$ for $n > 0$.
    \item $(\pi_* 1)_0 = 1$.
\end{propenum}
\end{prop}
\begin{proof}
Note that the QCA assumption guarantees that $\pi_*$ commutes with the pullback along $X^0 \to X^0/\mathbb G_m$ and $X \to X/\mathbb G_m$ (see \cref{rslt:QCA-maps-have-base-change-and-projection-formula}), which will be used implicitly in the proof (and is also implicitly used in the claim). The intuition behind the following proof is that $1 \in \QCoh(X)$ is not only $\mathbb G_m$-equivariant, but even $\mathbb A^1$-equivariant in an appropriate sense. The same is then true for $\pi_* 1$, resulting in the claimed description of its graded pieces. In fact, the proof works more generally for $\mathbb A^1$-equivariant sheaves on $X$.

More precisely, consider the following commuting diagram:
\begin{equation*}\begin{tikzcd}
    \mathbb G_m \times X \arrow[r,hook] \arrow[d,"\id \times \pi"] \arrow[rr,bend left,"g"] & \mathbb A^1 \times X \arrow[r,"a"] \arrow[d,"\id\times\pi"] & X \arrow[d,"\pi"] \\
    \mathbb G_m \times X^0 \arrow[r,hook] \arrow[rr,bend right,"g_0"] & \mathbb A^1 \times X^0 \arrow[r,"a_0"] & X^0
\end{tikzcd}\end{equation*}
Here the maps $g_0$ and $a_0$ are the projection to the second factor, while $g$ and $a$ are the action maps. Recall by \cref{rslt:compute-G-m-graded-pieces-using-retracts} that the graded pieces of $\pi_* 1$ are computed via the retracts of the idempotent maps coming from the natural map $\pi_* 1 \to \pi_* 1 \tensor g_{0*} 1 = \bigoplus_n (\pi_* 1) t^n$ by projecting to each factor. By the projection formula and base-change along the outer rectangle above (which is cartesian) we see that the map in question is the natural map $\pi_* 1 \to (\pi g)_* 1$. But this map factors over
\begin{align*}
    (\pi_ a)_* 1 = a_{0*} (\id \times \pi)_* 1 = a_{0*} 1 \tensor \pi_* 1 = \bigoplus_{n \ge 0} (\pi_* 1) t^n,
\end{align*}
which immediately shows $(\pi_* 1)_n = 0$ for $n > 0$ (mind the sign change in \cref{rslt:compute-G-m-graded-pieces-using-retracts}). It remains to compute $(\pi_* 1)_0$, which is the retract of the idempotent map
\begin{align*}
    \pi_* 1 \to (\pi a)_* 1 = \bigoplus_{n \ge 0} (\pi_* 1) t^n \xto{\pr_0} \pi_* 1.
\end{align*}
This map is induced by the projection $a_{0*} 1 \to 1$, which in turn comes from the closed immersion $X^0 = 0 \times X^0 \to \mathbb A^1 \times X^0$. Thus the above map $\pi_* 1 \to \pi_* 1$ is induced by the map $X = 0 \times X \to \mathbb A^1 \times X \xto{a} X$ by applying $\pi_*$. But by definition of $X^0$, this map is the same as $i \pi\colon X \to X^0 \to X$. Since $\pi i = \id_{X^0}$ we conclude that the map $\pi_* 1 \to \pi_* 1$ has retract $1$, as desired.
\end{proof}

An $\mathbb A^1$-retract $(X, X^0)$ lives in the category of $\mathbb A^1$-modules in stacks. By applying the tensor product $- \tensor_{\mathbb G_m} *$ we obtain the stacks $X/\mathbb G_m$ and $X^0/\mathbb G_m$, which are therefore both modules under the monoid $\mathbb A^1/\mathbb G_m$, with the module structure on $X^0$ being trivial. The following result provides a convenient way to reduce questions about the $\mathbb A^1/\mathbb G_m$-module $X/\mathbb G_m$ to the $\mathbb A^1/\mathbb G_m$-module $\mathbb A^1/\mathbb G_m \times X/\mathbb G_m$ on which $ \mathbb A^1/\mathbb G_m$ acts by multiplication on the first factor.

\begin{lem} \label{rslt:pushforward-along-action-of-A1-mod-Gm}
Let $(X, X^0)$ be a QCA $\mathbb A^1$-retract and denote $a\colon \mathbb A^1/\mathbb G_m \times X/\mathbb G_m \to X/\mathbb G_m$ the action map. Then $a$ is QCA and $a_* 1 = 1$.
\end{lem}
\begin{proof}
To show that $a$ is QCA it is enough to observe that its composition with $\pi\colon X/\mathbb G_m \to X^0/\mathbb G_m$ is QCA (see \cref{rslt:stability-of-QCA-maps}). But this composition is the same as the projection to $X/\mathbb G_m$ followed by $\pi$, hence is indeed QCA.

We now compute $a_* 1$. We rephrase the problem slightly and observe that the action of $\mathbb A^1/\mathbb G_m$ on $X/\mathbb G_m$ induces an action of $\QCoh(\mathbb A^1/\mathbb G_m)$ on $\QCoh(X/\mathbb G_m)$ via $F * M := a_*(F \boxtimes M)$ (in fact, for this proof we do not need the full action, but just the functor $- * -$). We then have $a_* 1 = 1 * 1$. Consider the open immersion $j\colon * = \mathbb G_m/\mathbb G_m \injto \mathbb A^1/\mathbb G_m$. Then $j_* 1$ is the identity for the convolution product on $\QCoh(\mathbb A^1/\mathbb G_m)$ and hence $(j_* 1) * M = M$ for all $M \in \QCoh(X/\mathbb G_m)$ (this identity can also easily be checked by hand). Denoting $K := \fib(1 \to j_* 1)$, our problem is now reduced to showing $K * 1 = 0$.

Denote $i'\colon 0/\mathbb G_m \to \mathbb A^1/\mathbb G_m$ the closed immersion. We observe that $K[1] = \Z[t, t^{-1}]/\Z[t] = \varinjlim_{n\ge 1} t^{-n} \Z[t]/\Z[t]$. The fiber of two consecutive terms in this colimit is isomorphic to $i'_* \calO(n)$ for $n \ge 1$ (cf. \cref{rmk:sign-convention-for-torus-reps}). Since $- * 1$ commutes with colimits (as $a$ is QCA), the claim $K * 1 = 0$ reduces to $(i'_* \calO(n)) * 1 = 0$ for all $n \ge 1$. Let us denote $a'\colon 0/\mathbb G_m \times X/\mathbb G_m \to X/\mathbb G_m$ the restricted action map. Then by the projection formula, for all $M \in \QCoh(X/\mathbb G_m)$ we have
\begin{align*}
    (i'_* \calO(n)) * M = a'_*(\calO(n) * M).
\end{align*}
Moreover, $a'$ factors as
\begin{align*}
    0/\mathbb G_m \times X/\mathbb G_m \xto{\id \times \pi} 0/\mathbb G_m \times X^0/\mathbb G_m \xto{a''} X^0/\mathbb G_m \xto{i} X/\mathbb G_m,
\end{align*}
where $a''$ is induced by the multiplication on $*/\mathbb G_m$, i.e. computes the convolution product on $X^0/\mathbb G_m = X^0 \times */\mathbb G_m$. From the symmetric monoidal equivalence in \cref{rslt:Cartier-duality-for-torus} we deduce that convolution with $\calO(n)$ spits out the $n$-graded piece. Altogether we conclude that
\begin{align*}
    (i'_* \calO(n)) * M = i_* (\pi_* M)_n.
\end{align*}
Using \cref{rslt:grading-on-cohomology-of-A1-retract} we conclude $(i'_* \calO(n)) * 1 = 0$, as desired.
\end{proof}

The next result is the key computation that goes into our main result on $\mathbb A^1$-retracts below. In the following we denote by $\hat{\mathbb A}^1 \subseteq \mathbb A^1$ the completion of $\mathbb A^1$ at $0$ (as defined in \cref{def:completion-of-stack-at-closed-subset}). Note that the $\mathbb G_m$-action on $\mathbb A^1$ restricts to a $\mathbb G_m$-action on the substack $\hat{\mathbb A}^1$.

\begin{lem} \label{rslt:hat-i-push-is-trivial-modulo-G-m}
Let $Y$ be a stack and denote by $\hat i\colon \hat{\mathbb A}^1/\mathbb G_m \times Y \injto \mathbb A^1/\mathbb G_m \times Y$ the canonical inclusion. Then $\hat i_* 1 = 1$.
\end{lem}
\begin{proof}
There is a natural map $1 \to \hat i_* 1$. Since the pushforward along $\pi\colon \mathbb A^1/\mathbb G_m \times Y \to */\mathbb G_m \times Y$ is conservative (as this map is affine), it is enough to show that the natural map $\pi_* 1 \to \pi_* \hat i_* 1$ is an isomorphism. The idea is now that $\pi_* \hat i_* 1$ wants to be the power series ring $1[[t]] = \varprojlim_n 1[t]/t^n$, while $\pi_* 1$ is the polynomial ring $1[t]$ -- but both of these rings agree on graded pieces, hence must agree in $\QCoh(*/\mathbb G_m \times Y)$.

To make this argument more precise, let $Z_n = \Spec \Z[t]/t^n$ and let $i_n\colon Z_n/\mathbb G_m \times Y \to \mathbb A^1/\mathbb G_m \times Y$ be the induced closed immersion. By \cref{rslt:completion-of-affine-is-colimit} we have $\hat{\mathbb A}^1 = \varinjlim_n Z_n$, hence the same is true after modding out by $\mathbb G_m$ and taking the product with $Y$, so we deduce $\hat i_* 1 = \varprojlim_n i_{n*} 1$. Note that $\pi_* i_{n*} 1 = 1[t]/t^n$ and the claim reduces to showing that the natural map
\begin{align*}
    1[t] \to \varprojlim_n 1[t]/t^n
\end{align*}
is an isomorphism in $*/\mathbb G_m \times Y$. This can be checked on graded pieces and taking graded pieces commutes with limits. Now for each fixed integer $m < 0$, the $m$-th graded piece of the right-hand side is constant $1$ for $n \ge -m$, and for $m > 0$ it is constant $0$ for all $n$. This proves the claim.
\end{proof}

We combine the above results in order to arrive at the following way to compute $\pi_* M$ in terms of $i^* M$ and some easy-to-control error terms:

\begin{prop} \label{rslt:perfect-module-is-retract-of-totalization-along-A1-retract}
Let $(X, X^0)$ be a QCA $\mathbb A^1$-retract and $M \in \Perf(X/\mathbb G_m)$. Denote $i_\bullet\colon X^\bullet \to X$ the \v{C}ech nerve of $i_0 = i\colon X^0 \to X$ and denote the maps modulo $\mathbb G_m$ by the same letters. Then the natural map
\begin{align*}
    M \to \varprojlim_{n\in\Delta} i_{n*} i_n^* M
\end{align*}
admits a splitting in $\QCoh(X/\mathbb G_m)$, i.e. $M$ is a retract of the righ-hand side.
\end{prop}
\begin{proof}
It is enough to prove the claim for $M = 1$, the general case then follows by tensoring with $M$ (observe that this preserves limits because $M$ is perfect). We consider the following commuting diagram:
\begin{equation*}\begin{tikzcd}
    */\mathbb G_m \times X/\mathbb G_m \arrow[rr,bend left,"i'"] \arrow[r,dashed] & Y^0 \arrow[r,"i''"] \arrow[d,"a'"] & \mathbb A^1/\mathbb G_m \times X/\mathbb G_m \arrow[d,"a"]\\
    & X^0/\mathbb G_m \arrow[r,"i"] & X/\mathbb G_m
\end{tikzcd}\end{equation*}
Here $Y^0$ is defined such that the square is cartesian, $a$ is the action map and $i'$ is induced by the closed immersion $\{ 0 \} \to \mathbb A^1$. Note that $ai'$ factors over $X^0/\mathbb G_m$, which provides the dashed arrow in the diagram. By \cref{rslt:pushforward-along-action-of-A1-mod-Gm} $a$ is QCA and $a_* 1 = 1$, hence the same is true for $a'$ by base-change and similarly for all transition maps $Y^n \to X^n$, where $i''_\bullet\colon Y^\bullet \to \mathbb A^1/\mathbb G_m \times X/\mathbb G_m$ denotes the \v{C}ech nerve of $i''$. Thus the map $1 \to \varprojlim_n i_{n*} 1$ in $\QCoh(X/\mathbb G_m)$ is obtained from the map $1 \to \varprojlim_n i''_{n*} 1$ in $\QCoh(\mathbb A^1/\mathbb G_m \times X/\mathbb G_m)$ by applying $a_*$, so it is enough to prove that the latter map admits a splitting. If $i'_\bullet$ denotes the \v{C}ech nerve of $i'$ then there is a natural map $i'_\bullet \to i''_\bullet$ of simplicial stacks, which induces maps
\begin{align*}
    1 \to \varprojlim_{n\in\Delta} i''_* 1 \to \varprojlim_{n\in\Delta} i'_* 1
\end{align*}
in $\QCoh(\mathbb A^1/\mathbb G_m \times X/\mathbb G_m)$. It is therefore enough to show that the composition of these two maps is an isomorphism. By \cref{rslt:hat-i-push-is-trivial-modulo-G-m} it is enough to show that
\begin{align*}
    \varprojlim_{n\in\Delta} i'_* 1 = \hat i_* 1,
\end{align*}
where $\hat i$ is defined as in \cref{rslt:hat-i-push-is-trivial-modulo-G-m}. To prove this identity, consider the map $f\colon */\mathbb G_m \times X/\mathbb G_m \to \hat{\mathbb A}^1/\mathbb G_m \times X/\mathbb G_m$. The \v{C}ech nerve of $f$ is the same as the \v{C}ech nerve of $i'$ (because $\hat{\mathbb A}^1 \injto \mathbb A^1$ is a monomorphism), so it is enough to show that $\QCoh$ descends along $f$. This is standard: By the dual of \cite[Corollary~4.7.5.3]{lurie-higher-algebra} it is enough to show that $f^*$ is conservative and preserves $f^*$-split totalizations. Under the equivalence in \cref{rslt:hat-i-pullback-induces-equiv} and after base-change to an affine cover, the conservativity of $f^*$ follows from the observation that for a ring $A$, an $A[t]$-module $M$ is zero as soon as $M/t$ and $M \tensor_{A[t]} A[t, t^{-1}]$ are zero. Similarly, the fact that $f^*$ preserves $f^*$-split totalizations reduces to the observation that, after passing to a cover by $\Spec A$, the functor $f_* f^*$ is given by $M \mapsto \fib(M \xto{t} M)$, which commutes with limlts.
\end{proof}

\begin{cor} \label{rslt:A1-retract-pushforward-in-terms-of-pullback}
Let $(X, X^0)$ be a QCA $\mathbb A^1$-retract and $M \in \Perf(X/\mathbb G_m)$. Then $\pi_* M \in \QCoh(X^0/\mathbb G_m)$ is obtained using countable limits from the objects
\begin{align*}
    i^* M \tensor (i^* i_* 1)^{\tensor n}, \qquad n \ge 0.
\end{align*}
Moreover, $i^* i_* 1$ is non-positively graded with $(i^* i_* 1)_0 = 1$.
\end{cor}
\begin{proof}
Since retracts can be computed via countable limits, by \cref{rslt:perfect-module-is-retract-of-totalization-along-A1-retract} we immediately reduce the claim to showing that for all $n \in \Delta$ the object $\pi_* i_{n*} i_n^* M$ lies in the claimed subcategory of $\QCoh(X^0/\mathbb G_m)$. Note that $i_n$ factors as
\begin{align*}
    i_n\colon X^n/\mathbb G_m \xto{f_n} X^0/\mathbb G_m \xto{i} X/\mathbb G_m,
\end{align*}
where $f_n$ is the projection onto the first factor. Using the identity $\pi i = \id$ and the projection formula (by \cref{rslt:QCA-maps-have-base-change-and-projection-formula}) we deduce
\begin{align*}
    \pi_* i_{n*} i_n^* M = \pi_* i_* f_{n*} f_n^* i^* M = f_{n*} f_n^* i^* M = f_{n*} 1 \tensor i^* M.
\end{align*}
Now observe that by base-change and projection formula we have $f_{n*} 1 \isom (f_{1*} 1)^{\tensor n} = (i^* i_* 1)^{\tensor n}$, where we implicitly use that $f_{n*} 1$ agrees with the similarly defined pushforwards along the projection $X^n/\mathbb G_m \to X^0/\mathbb G_m$ to any other factor (because they agree after further pushing along $i$ and hence after pushing along $\pi i = \id$).

It remains to prove the claim about the grading on $i^* i_* 1 = f_{1*} 1$. But note that $(X^0 \times_X X^0, X^0)$ is a QCA $\mathbb A^1$-retract and $f_1$ plays the role of $\pi$ for this retract, so we conclude by \cref{rslt:grading-on-cohomology-of-A1-retract}.
\end{proof}

\begin{rmk}
We view \cref{rslt:A1-retract-pushforward-in-terms-of-pullback} as an answer to the initial question of how much $\pi_*$ and $i^*$ differ: At least on perfect objects, one can quantify the error in terms of the sheaves $(i^* i_* 1)^{\tensor n}$.
\end{rmk}

From \cref{rslt:perfect-module-is-retract-of-totalization-along-A1-retract} we also obtain the following result about preservation of perfect sheaves along pushforward. The result may be far from optimal, but is strong enough for our applications.

\begin{cor} \label{rslt:A1-retract-pushforward-preserves-perfect}
Fix a classical regular noetherian $\Q$-algebra $\Lambda$ and let $(X, X^0)$ be a QCA $\mathbb A^1$-retract $\Lambda$. Assume that
\begin{enumerate}[(a)]
    \item $X^0$ is QCA and suave over $\Lambda$.
    \item For all $n \le 0$, the sheaf $(i^* i_* 1)_n \in \QCoh(X^0/\mathbb G_m)$ is perfect.
\end{enumerate}
Then for every perfect sheaf $M \in \Perf(X/\mathbb G_m)$ and all $n \in \Z$, the sheaf $(\pi_* M)_n$ is coherent on $X^0$.
\end{cor}
\begin{proof}
Fix $M \in \Perf(X/\mathbb G_m)$ and $n \in \Z$. By \cref{rslt:perfect-module-is-retract-of-totalization-along-A1-retract} the sheaf $\pi_* M$ is a retract of $\varprojlim_{m\in\Delta} M^m$ for $M^m := \pi_* i_{m*} i_m^* M$. Since passage to graded pieces commutes with limits, we see that $(\pi_* M)_n$ is a retract of $\varprojlim_{m\in\Delta} M^m_n$. We now claim that $M^m_n$ is perfect and uniformly bounded (for fixed $n$ and varying $m$) for the natural t-structure on $\QCoh(X^0)$ from \cref{rslt:t-structure-on-QCoh}. To see this, as in the proof of \cref{rslt:A1-retract-pushforward-in-terms-of-pullback} we compute
\begin{align*}
    &M^m_n = (\pi_* i_{m*} i_m^* M)_n = (i^* M \tensor (i^* i_* 1)^{\tensor m})_n = \bigoplus_{k\in\Z} (i^* M)_{n-k} \tensor ((i^* i_* 1)^{\tensor m})_k =\\&\qquad = \bigoplus_{k_1, \dots, k_m \in \Z} (i^* M)_{n - (k_1 + \dots + k_m)} \tensor (i^* i_* 1)_{k_1} \tensor \dots \tensor (i^* i_* 1)_{k_m}.
\end{align*}
Since $M$ is perfect, so is $i^* M$, and hence all $(i^* M)_{n-k}$ are perfect and only finitely many are non-zero. Moreover, since $i^* i_* 1$ is non-positively graded (see \cref{rslt:A1-retract-pushforward-in-terms-of-pullback}), the above terms can only be non-zero if $k_i \le 0$ for all $i$. Furthermore, $(i^* i_* 1)_0 = 1$, so if $k_i = 0$ then it can be ignored in the above tensor product. Altogether we see that there are only finitely many different terms that can appear in the above direct sum. Each term is perfect and hence bounded, because $X^0$ is suave and hence admits a smooth cover by $\Spec A$ for $A$ a suave $\Q$-algebra, which is truncated by \cref{rslt:characterization-of-suave-maps}. This finishes the proof that $M^m_n$ is perfect and uniformly bounded in $m$.

Since totalizations of uniformly left-bounded sheaves are computed via a spectral sequence, we observe that for every $k \in \Z$, the complex $\tau^{\le k} \varprojlim_{m\in\Delta} M^m$ only depends on a \emph{finite} limit over the $M^m$'s. This finite limit is perfect and hence coherent. Since coherent sheaves are stable under truncation (cf. \cref{rslt:coh-over-regular-noetherian-base}) we deduce that $\tau^{\le k} \varprojlim_{m\in\Delta} M^m$ is coherent for every $k \in \Z$. Now $(\pi_* M)_n$ is a retract of $\varprojlim_{m\in\Delta} M^m$, but since $\pi$ is QCA, $(\pi_* M)_n$ is also right-bounded (see \cref{rslt:QCA-pushforward-has-finite-cohomological-dimension}). We deduce that $(\pi_* M)_n$ is a retract of $\tau^{\le k} \varprojlim_{m\in\Delta} M^m$ for $k \gg 0$ and hence coherent.
\end{proof}

\section{Spectral results} \label{sec:spectral}

Throughout this section we fix a prime $p$ and a finite extension $F$ of $\Q_p$. We let $q$ be the size of the residue field $\calO_F/\mathfrak m_F$. We denote by $\Breve{F}$ the completed maximal unramified extension of $F$, obtained by adjoining all prime-to-$p$ roots of unity. We denote by $C \supset \Breve F$ a completed algebraic closure. We denote by $W_F$ the Weil group of $F$, i.e. $W_F = I_F \rtimes \Z$ where $I_F := \Gal(C/\Breve E)$ is the inertia subgroup and $\Z$ acts via the Frobenius element in $\Gal(\Breve E/E) \isom \hat \Z$. We also fix a prime $\ell \ne p$.

\subsection{Stacks of representations}

In the following we introduce and study representation stacks, which the L-parameter stack is a special case of. Our exposition is a variant of \cite[\S2]{Zhu}, but we generalize some aspects and provide a new approach to handling topological groups based on condensed sets.

Representation stacks are special types of mapping stacks, so let us first introduce the following notation.

\begin{defn}
Fix a ring $\Lambda$. Then $\Ani$ acts on $\Stk_\Lambda$ via $S \times X := \varinjlim_S X$ for all $S \in \Ani$ and $X \in \Stk_\Lambda$. We will often identify an anima $S$ with the stack $S \times *$. For fixed $S$, we denote by
\begin{align*}
    \intMap(S, -)\colon \Stk_\Lambda \to \Stk_\Lambda
\end{align*}
the right adjoint of $S \times -$.
\end{defn}

Note that $\intMap(S, X) = \intHom(S \times *, X) = \varprojlim_S X$, where $\intHom$ denotes the internal hom in $\Stk_\Lambda$. Explicitly, for any stack $X \in \Stk_\Lambda$ and any $\Lambda$-algebra $A$ we have
\begin{align*}
    \intMap(S, X)(A) = \Hom(S, X(A)),
\end{align*}
where on the right-hand side we have the usual homomorphism anima in $\Ani$. Indeed, note that on both sides one can pull out colimits in $S$ as limits, so the claim reduces to $S = *$, which is clear.

Given a stack $X$ over some ring $\Lambda$ and an anima $S$, we have $\QCoh(S \times X) = \Fun(S, \QCoh(X))$, as one checks by pulling out colimits in $S$ on both sides. In particular, if $S \to T$ is a map of anima and $f\colon S \times X \to T \times X$ is the induced map of stacks then $f^*$ admits a left adjoint
\begin{align*}
    f_\natural\colon \QCoh(S \times X) \to \QCoh(T \times X)
\end{align*}
given by left Kan extension. In the case $T = *$, $f_\natural$ is simply given as $\varinjlim_S$. Using these observations, we can compute the cotangent complex of $\intMap(S, X)$:

\begin{prop} \label{rslt:cotangent-complex-of-mapping-stack}
Fix a ring $\Lambda$, a stack $X \in \Stk_\Lambda$ and an anima $S$. Consider the correspondence
\begin{align*}
    X \xfrom{\ \ev \ } \intMap(S, X) \times S \xto{\ \pr \ } \intMap(S,X).
\end{align*}
Assume that $X$ admits a cotangent complex $L_{X/\Lambda} \in \QCoh(X)$. Then
\begin{align*}
    L_{\intMap(S,X)/\Lambda} := \pr_\natural \ev^* L_{X/\Lambda}
\end{align*}
is a cotangent complex for $\intMap(S, X)$, provided that it lies in $\QCoh^-$.
\end{prop}
\begin{proof}
We have $\intMap(S, X) = \varprojlim_S X$, so by \cref{rslt:cotangent-complex-of-limit} we have a natural isomorphism $L_{\intMap(S,X)/\Lambda} = \varinjlim_S L_{X/\Lambda}|_{\intMap(S,X)}$. On the other hand, $\pr_\natural$ is computed as $\varinjlim_S$, so we only need to verify that
\begin{align*}
    \ev^* L_{X/\Lambda} = (L_{X/\Lambda}|_{\intMap(S,X)})_S
\end{align*}
in $\Fun(S, \QCoh(\intMap(S, X)))$. This follows by going through the definitions: Note that $\ev$ corresponds to the canonical map $S \to \Hom(\intMap(S, X), X)$ and composing this map with $\QCoh$ provides a functor $\alpha\colon S \to \Fun(\QCoh(X), \QCoh(\intMap(S, X)))$. Now $L_{X/\Lambda}$ induces a (constant) functor $S \to \QCoh(X)$ and composing this with $\alpha$ yields $\ev^* L_{X/\Lambda}$. On the other hand, this also computes $(L_{X/\Lambda}|_{\intMap(S,X)})_S$ by definition.
\end{proof}

The next result shows that under mild assumptions, the mapping stack $\intMap(S, -)$ preserves closed and open immersions.

\begin{lem} \label{rslt:mapping-stack-preserves-locally-closed-immersions}
Fix an anima $S$ a ring $\Lambda$.
\begin{lemenum}
    \item The functor $\intMap(S, -)$ preserves closed immersions of stacks over $\Lambda$.
    \item Suppose that $\pi_0 S$ is finite. Then the functor $\intMap(S, -)$ preserves open immersions of stacks over $\Lambda$.
\end{lemenum}
\end{lem}
\begin{proof}
We first prove (i), so let $i\colon Z \to X$ be a closed immersion in $\Stk_\Lambda$. Since $\intMap(S, -) = \varprojlim_S (-)$, we need to show that the map $\varprojlim_S Z \to \varprojlim_S X$ is a closed immersion. Fix any ring $A$ together with a map $\Spec A \to \varprojlim_S X$. Since $i$ is a closed immersion, we know that $Z \times_X \Spec A$ is affine, say equal to $\Spec B$, and $\pi_0 A \to \pi_0 B$ is surjective. We now obtain the diagram
\begin{equation*}\begin{tikzcd}
    Y \arrow[r] \arrow[d] & \varprojlim_S \Spec B \arrow[r] \arrow[d] & \varprojlim_S Z \arrow[d]\\
    \Spec A \arrow[r] & \varprojlim_S \Spec A \arrow[r] & \varprojlim_S X
\end{tikzcd}\end{equation*}
The bottom right map is the one induced from $\Spec A \to \varprojlim_S X$ and clearly the right-hand square is cartesian. The bottom left horizontal map is the diagonal map and $Y$ is defined so that the left-hand square is also cartesian. Note that the composition of the bottom horizontal maps is the map $\Spec A \to \varprojlim_S X$ that we started with, and since this map was arbitrary, it is enough to show that the map $Y \to \Spec A$ is a closed immersion. This in turn reduces to showing that the middle vertical map is an isomorphism. Since the embedding of affine schemes into stacks preserves limits, we are reduced to showing that the map of rings $\varinjlim_S A \to \varinjlim_S B$ is surjective on $\pi_0$. Since $\pi_0$ is left adjoint and hence commuted with colimits, this question immediately reduces to the claim that classical colimits of classical rings preserve surjective maps. This is clear for filtered colimits and the initial map $\Z \to \Z$, hence reduces to pushouts. But pushouts of rings are computed via relative tensor products, which themselves are quotients of tensor products over $\Z$. We are thus left to show that for surjective maps $A_i \to B_i$, also $A_1 \tensor A_2 \to B_1 \tensor B_2$ is surjective; but this follows from right-exactness of tensor products, finishing the proof of (i).

We now prove (ii), so assume that $\pi_0 S$ is finite and let $j\colon U \injto X$ be an open immersion. Write $S = \bigsqcup_{i=1}^n S_i$ for connected anima $S_i$. Then $\intMap(S, -) = \prod_{i=1}^n \intMap(S_i, -)$ and since finite products of open immersions are again open immersions, we reduce to the case that $S$ is connected. Fix a surjective map $\rho\colon * \to S$. Using the identification $\intMap(*, -) = \id$, $\rho$ induces the vertical maps in the following commuting square:
\begin{equation*}\begin{tikzcd}
    \intMap(S, U) \arrow[r,hook,"j'"] \arrow[d] & \intMap(S, X) \arrow[d]\\
    U \arrow[r,hook,"j"] & X
\end{tikzcd}\end{equation*}
We claim that this square is cartesian, which immediately implies (ii). First observe that $j'$ is a monomorphism (i.e. its diagonal is an isomorphism), because $\intMap(S, -)$ preserves limits and in particular diagonals. The cartesianness can be checked on $A$-valued points for all rings $A$, where it becomes the claim that the full subanima $\Hom(S, U(A)) \subseteq \Hom(S, X(A))$ consists precisely of those maps $S \to X(A)$ whose precomposition with $\rho$ factors over $U(A)$. But this is clear since $U(A) \subseteq X(A)$ is a full subanima.
\end{proof}

We now specialize the general mapping stack to the stack of representations by plugging in classifying stacks of groups.

\begin{defn} \label{def:rep-stack}
Fix a ring $\Lambda$, a group $H$ in $\Stk_\Lambda$ and a group $\Gamma$ in $\Ani$ (e.g. a classical abstract group). We define
\begin{align*}
    \intRep_H(\Gamma) := \intMap(*/\Gamma, */H) \in \Stk_\Lambda.
\end{align*}
and call it the \emph{stack of $H$-representations of $\Gamma$}. Precomposition with the canonical map $* \to */\Gamma$ provides a canonical projection $\tau\colon \intRep_H(\Gamma) \to */H$ and we define
\begin{align*}
    \intRep_H^\square(\Gamma) := \intRep_H(\Gamma) \times_{*/H} *
\end{align*}
and call it the \emph{stack of framed $H$-representations of $\Gamma$}. By construction $\intRep_H^\square(\Gamma)$ comes equipped with an action of $H$ and we have $\intRep_H(\Gamma) = \intRep_H^\square(\Gamma)/H$. We denote by $\tau'\colon \intRep_H^\square(\Gamma) \to */H$ the composition of $\tau$ with the projection to $\intRep_H(\Gamma)$.
\end{defn}

The general description of $\intMap(S, X)$ implies that for a $\Lambda$-algebra $A$ the $A$-valued points of $\intRep_H(\Gamma)$ parametrize $H$-torsors over $\Spec A$ that are equipped with a $\Gamma$-action. The $A$-valued points of $\intRep_H^\square(\Gamma)$ additionally parametrize a trivialization of the $H$-torsor. The following result summarizes the basic general properties of $\intRep_H(\Gamma)$ and in particular provides a better description for it: It classifies group homomorphisms $\Gamma \to H$ up to conjugation.

\begin{prop} \label{rslt:basic-properties-of-rep-stack}
Fix a ring $\Lambda$, a group $H$ in $\Stk_\Lambda$ and a group $\Gamma$ in $\Ani$.
\begin{propenum}
    \item \label{rslt:framed-intrep-equals-group-hom} For every $\Lambda$-algebra $A$ there is a natural equivalence
    \begin{align*}
        \intRep_H^\square(\Gamma)(A) = \Hom_{\Mon(\Ani)}(\Gamma, H(A)),
    \end{align*}
    where $\Mon(\Ani)$ is the category of monoids in $\Ani$. In case that both $\Gamma$ and $H(A)$ are static (i.e. classical groups), the above $\Hom$ is the set of group homomorphisms $\Gamma \to H(A)$. Under this equivalence, the action of $H$ on $\intRep_H^\square(\Gamma)$ is given by conjugation.

    \item \label{rslt:framed-intrep-is-affine} If $H$ is represented by an affine scheme, then so is $\intRep_H^\square(\Gamma)$. In particular, if $H$ is smooth and affine then $\intRep_H(\Gamma)$ is geometric.

    \item \label{rslt:cotangent-complex-of-rep-stack} Suppose that $H$ is represented by a smooth affine scheme over $\Lambda$ and let $\Ad \in \QCoh(*/H)$ denote the adjoint representation (see \cref{def:adjoint-representation}). Then
    \begin{align*}
        L_{\intRep_H(\Gamma)/\Lambda} &= (\tau^* \Ad^\vee)_\Gamma[-1],\\
        L_{\intRep_H^\square(\Gamma)/\Lambda} &= \cofib(\tau'^* \Ad^\vee \to (\tau'^* \Ad^\vee)_\Gamma)[-1].
    \end{align*}
    Here $(-)_\Gamma$ denotes $\Gamma$-coinvariants, i.e. the colimit over $*/\Gamma$.

    More explicitly, suppose we are given a $\Lambda$-algebra $A$ and a group homomorphism $\rho\colon \Gamma \to H(A)$, i.e. a point $\rho\colon \Spec A \to \intRep_H^\square(\Gamma)$. We also denote by $\rho$ its composition with the projection to $\intRep_H(\Gamma)$. Then
    \begin{align*}
        \rho^* L_{\intRep_H(\Gamma)/\Lambda} &= (\Ad^\vee_\rho)_\Gamma[-1],\\
        \rho^* L_{\intRep_H^\square(\Gamma)/\Lambda} &= \cofib(\Ad^\vee_\rho \to (\Ad_\rho^\vee)_\Gamma)[-1].
    \end{align*}
    Here $\Ad^\vee_\rho$ denotes the $\Gamma$-representation obtained by inflating $\Ad^\vee \tensor_\Lambda A$ along $\rho$.
\end{propenum}
\end{prop}
\begin{proof}
All of these claims are also claimed in \cite[\S2]{Zhu}, but we provide more (and slightly different) justification for some of the claims.

We start with the proof of (i), so fix $A$. Then it follows immediately from the definitions and the above description of $\intMap(S, X)(A)$ that
\begin{align*}
    \intRep_H^\square(\Gamma)(A) = \Hom_{\Ani_{*/}}(*/\Gamma, (*/H)(A)).
\end{align*}
Note furthermore that $\Hom_{\Ani_{*/}}(*/\Gamma \to */H(A)) = \Hom_{\Mon(\Ani)}(\Gamma, H(A))$ (this is true for any group in place of $H(A)$, see \cite[Example~5.2.6.13]{lurie-higher-algebra}). Hence it is enough to show that the natural map $*/H(A) \to (*/H)(A)$ is a monomorphism (intuitively, this says that the automorphisms of the trivial $H$-torsor are given by $H$). Here is a quick argument for that: $*/H$ is the sheafification of the presheaf $A \mapsto */H(A)$ and by the explicit description of sheafification in the proof of \cite[Proposition~6.2.2.7]{lurie-higher-topos-theory} it is enough to show that this presheaf is separated, i.e. its value on $A$ sits monomorphically in the descent data along any cover of $A$ -- but this immediately reduces to the fact that $H$ descends along that cover.

The identification of the $H$-action in (i) follows by going through the definitions. Generally, if $X \to */H$ is a map then the $H$-action on $X' := X \times_{*/H} *$ is given by the map
\begin{align*}
    X' \times H = X \times_{*/H} * \times_{*/H} * \to X \times_{*/H} * = X'
\end{align*}
coming from the projection to the first and third coordinate.

We now prove (ii), so assume that $H$ is represented by an affine scheme. Recall that the forgetful functor $\Mon(\Ani) \to \Ani$ has a left adjoint sending any anima $S$ to the free monoid $F(S)$ on it (see \cite[Example~3.1.3.14]{lurie-higher-algebra}). Since $\Ani$ has compact projective generators given by finite sets $I$, it follows that $\Mon(\Ani)$ has compact projective generators given by $F(I)$. In particular, $\Gamma$ can be written as a sifted colimit of $F(I)$'s. By (i) we deduce that $\intRep_H^\square(\Gamma)$ is a cosifted limit of $\intRep_H^\square(F(I))$'s, so the claim reduces to the case that $\Gamma = F(I)$ (because affine schemes are stable under limits). But by the universal property of free monoids and (i) we immediately see that $\intRep_H^\square(F(I)) = H^I$, which is indeed affine.

We now prove (iii), so assume that $H$ is represented by a smooth affine scheme. Then by \cref{rslt:cotangent-of-BG-is-adjoint-representation} the cotangent complex of $*/H$ is given by $\Ad^\vee[-1]$. By \cref{rslt:cotangent-complex-of-mapping-stack} the cotangent complex of $\intRep_H(\Gamma)$ is given as the colimit over $*/\Gamma$ of the pullback of $\Ad^\vee[-1]$ along the evaluation map $\intRep_H(\Gamma) \times */\Gamma \to */H$. But this pullback has as underlying sheaf $f^* \Ad^\vee$ on $\intRep_H(\Gamma)$ together with an induced $\Gamma$-action, which proves the claimed description of the cotangent complex of $\intRep_H(\Gamma)$. For $\intRep_H^\square(\Gamma)$ we use $\intRep_H^\square(\Gamma) = \intRep_H(\Gamma) \times_{*/H} *$ and \cref{rslt:cotangent-complex-of-limit} in order to deduce that the cotangent complex of $\intRep_H^\square(\Gamma)$ is the pushout of the pullbacks of the cotangent complexes of $\intRep_H(\Gamma)$, $*/H$ and $*$, which exactly amounts to the claimed cofiber sequence.

We now prove the second part of (iii), so let $\rho\colon \Gamma \to H(A)$ be given. Consider the induced diagram
\begin{equation*}\begin{tikzcd}
    & \Spec A \times */\Gamma \arrow[r,"\pr"] \arrow[d,"\rho \times \id"] & \Spec A \arrow[d,"\rho"]\\
    */H & \intRep_H(\Gamma) \times */\Gamma \arrow[l,"\ev"] \arrow[r,"\pr"] & \intRep_H(\Gamma)
\end{tikzcd}\end{equation*}
where the square is cartesian. Then
\begin{align*}
    \rho^* L_{\intRep_H(\Gamma)/\Lambda} = \rho^* \pr_\natural \ev^* \Ad^\vee = \pr_\natural (\rho \times \id)^* \ev^* \Ad^\vee = ((\rho \times \id)^* \ev^* \Ad^\vee)_\Gamma.
\end{align*}
It remains to show that $((\rho \times \id)^* \ev^* \Ad^\vee = \Ad^\vee_\rho$. By pulling back the bottom row of the above diagram along $\Spec A \to \Spec \Lambda$ and replacing $A$ by $\Lambda$ we can assume that $A = \Lambda$. Then $\ev \comp (\rho \times \id)$ is a map $*/\Gamma \to */H$ and one sees immediately that it is induced by the group homomorphism $\rho$, as desired. The claim for $\intRep_H^\square(\Gamma)$ is proved similarly.
\end{proof}

Let us introduce a few more standard definitions for representations, reinterpreted in terms of the stack of representations. We refer the reader to \cref{def:intIsom-stack} and \cref{rslt:basic-properties-of-intIsom-stack} for the definition and basic properties of the isomorphism stack $\intIsom$ appearing below.

\begin{defn}
Fix a ring $\Lambda$, a group $H$ in $\Stk_\Lambda$ and a group $\Gamma$ in $\Ani$. Let $\rho\colon \Gamma \to H(\Lambda)$ be a group homomorphism, which we identify with a $\Lambda$-valued point in $\intRep_H^\square(\Gamma)$ by \cref{rslt:framed-intrep-equals-group-hom}.
\begin{defenum}
    \item \label{def:orbit-map} The \emph{orbit map} for $\rho$ is the map
    \begin{align*}
        \orb_\rho\colon H \to \intRep_H^\square(\Gamma), \qquad g \mapsto g^{-1} \rho g.
    \end{align*}
    More precisely, $\rho$ induces a map $* \to \intRep_H(\Gamma)$ and the orbit map is the pullback of this map along $* \to */H$.

    \item \label{def:centralizer} Suppose that $H$ is affine and both $H$ and $\intRep_H^\square(\Gamma)$ are classical. Then we define the \emph{centralizer} of $\rho$ as
    \begin{align*}
        Z_H(\rho) := (\Isom_{\intRep_H(\Gamma)}(\rho))^\cl = (H \times_{\intRep^\square_H(\Gamma)} \{ \rho \})^\cl,
    \end{align*}
    where the implicit map $H \to \intRep^\square_H(\Gamma)$ is $\orb_\rho$. It is clear from the second description that $Z_H(\rho)$ is a classical closed subgroup of $H$.
\end{defenum}
\end{defn}

\begin{rmk}
It might appear unnatural to pass to the underlying classical scheme in the definition of the centralizer. We do this in order to avoid the behavior in \cref{rmk:intIsom-often-derived}. In fact, many classical group theoretic operations (like centralizers and kernels) do not seem to have a good analog in the derived world, much like kernels behave very differently in the derived category compared to its abelian heart.
\end{rmk}

With the above general results and definitions at hand, we now come to specific examples of representation stacks, motivated by the stack of L-parameters. 

\begin{lem} \label{rslt:rep-stack-for-finite-group}
Let $\Lambda$ be a ring, $H$ an affine smooth group scheme over $\Lambda$ and $\Gamma$ a finite group whose order is invertible in $\pi_0 \Lambda$. Then $\intRep_H(\Gamma)$ is a smooth geometric stack over $\Lambda$ and for every $\Lambda$-valued point $\rho$ in $\intRep_H^\square(\Gamma)$ the orbit map $\orb_\rho\colon H \to \intRep_H^\square(\Gamma)$ is smooth. If $\Lambda = k$ is an algebraically closed field then
\begin{align*}
    \intRep_H(\Gamma) = \bigsqcup_\rho */Z_H(\rho),
\end{align*}
where $\rho$ ranges over the set of conjugacy classes of group homomorphisms $\Gamma \to H(k)$.
\end{lem}
\begin{proof}
For special cases of $\Lambda$, this is claimed in \cite[Proposition~2.3.2]{Zhu} and is essentially proved in \cite[Lemma~A.1]{DHKM-ParG}. We give a quick self-contained proof. Since the order of $\Gamma$ is invertible in $\pi_0 \Lambda$, the trivial $\Gamma$-representation over $\Lambda$ is a direct summand of the regular representation, i.e. $\Lambda$ is retract of $\Lambda[\Gamma]$ in the category of $\Gamma$-representations (since $\Lambda$ is a $\Z[1/\abs \Gamma]$-algebra, this reduces to the case of $\Lambda = \Z[1/\abs\Gamma]$, which is clear). This implies that for every $\Gamma$-representation $V$, $V_\Gamma = V \tensor_{\Lambda[\Gamma]} \Lambda$ is a retract of $V$. Thus from \cref{rslt:cotangent-complex-of-rep-stack} we deduce that $L_{\intRep_H^\square(\Gamma)/\Lambda}$ is a vector bundle.

By \cref{rslt:framed-intrep-is-affine} we know that $\intRep^\square_H(\Gamma)$ is affine and admits an explicit description in terms of group homomorphisms $\Gamma \to H(A)$ for $\Lambda$-algebras $A$. This implies that the underlying classical scheme of $\intRep^\square_H(\Gamma)$ is a closed subscheme of $H^\Gamma$ and in particular of finite type over $\pi_0\Lambda$. Using that its cotangent complex is a vector bundle (as seen above) and \cite[Proposition~3.2.18]{DAG-Lurie} we deduce that $\intRep^\square_H(\Gamma)$ is lafp over $\Lambda$ and then \cref{rslt:cotangent-complex-for-geometric-map-and-smoothness} implies that $\intRep^\square_H(\Gamma)$ is smooth over $\Lambda$.

Now fix a representation $\rho\colon \Gamma \to H(\Lambda)$ and consider the induced map $f\colon * \to \intRep_H(\Gamma)$. Using the fiber sequence \cref{rslt:fiber-sequence-for-cotangent-complex} for $* \to \intRep_H(\Gamma) \to *$ and \cref{rslt:cotangent-complex-of-rep-stack} we deduce that $L_f = (\Ad_\rho^\vee)_\Gamma$, which is free and concentrated in degree $0$. Since $f$ is an affine map (its base-change is $\orb_\rho$), it follows from \cref{rslt:cotangent-complex-for-geometric-map-and-smoothness} that $f$ is smooth. Hence $\orb_\rho$ is smooth as well.

Finally, assume that $\Lambda = k$ is an algebraically closed field. Since $\orb_\rho$ is smooth, it has open image (see \cite[Proposition~00I1]{stacks-project}). Therefore, by varying $\rho$ we obtain a decomposition of $\intRep_H^\square(\Gamma)$ into disjoint open subschemes stable under $H$, which are of the form $H/Z_H(\rho)$ by definition of the centralizer. Modding out the $H$-action yields the desired description of $\intRep_H(\Gamma)$.
\end{proof}

With the case of finite $\Gamma$ understood, we now consider the following special type of $\Gamma$, which is closely related to the Weil group.

\begin{defn} \label{def:Gamma-q}
For a prime power $q = p^r$, we consider the following (classical) group:
\begin{align*}
    \Gamma_q := \langle \sigma, \tau \ | \ \sigma \tau \sigma^{-1} = \tau^q \rangle.
\end{align*}
In other words, $\Gamma_q = \tau^{\Z[1/p]} \rtimes \sigma^{\Z}$, where the action of $\sigma^{\Z}$ on $\tau^{\Z[1/p]}$ is given by $\sigma * \tau = \tau^q$.
\end{defn}

\begin{lem} \label{rslt:rep-stack-for-Gamma-q}
Let $\Lambda$ be a ring, $q = p^r$ a prime power invertible in $\pi_0 \Lambda$ and $H$ a smooth affine group scheme over $\Lambda$.
\begin{lemenum}
    \item $\intRep_H(\Gamma_q)$ is a local complete intersection over $\Lambda$ and its dualizing complex is trivial.
    \item If $\Lambda$ is a Cohen-Macaulay classical noetherian ring and $H$ is reductive then $\intRep_H(\Gamma_q)$ is classical.
\end{lemenum}
\end{lem}
\begin{proof}
This is proved in \cite[Proposition~2.3.7]{Zhu} in the case that $H$ is reductive and under stronger assumptions on $\Lambda$, but the argument for (i) works in general. Let us first assume that $\Lambda = \Z[1/p]$. The key observation is that $\Lambda_0 := \Z[1/p]$ has the following resolution as a $\Lambda_0[\Gamma_q]$-module:
\begin{align*}
    0 \to \Lambda_0[\Gamma_q] \xto{(1 - (\sum_{j<q} \tau^j)\sigma,\ \tau - 1)} \Lambda_0[\Gamma_q] \oplus \Lambda_0[\Gamma_q] \xto{(1-\tau,\ 1-\sigma)} \Lambda_0[\Gamma_q] \to \Lambda_0 \to 0.
\end{align*}
This allows an easy computation of $\Gamma_q$-homology: For every $\Lambda_0$-algebra $A$ and every $A[\Gamma_q]$-module $V$ the $\Gamma_q$-homology of $V$ is obtained by tensoring the above resolution over $\Lambda_0[\Gamma_q]$ with $V$, producing a similar resolution for $V_{\Gamma_q}$. In particular, if $V$ is a vector bundle over $A$ then $\cofib(V \to V_{\Gamma_0})[-1]$ has Tor amplitude in $[-1, 0]$. Using \cref{rslt:cotangent-complex-of-rep-stack} and applying the above observations to $V = \tau^* \Ad^\vee$ we deduce that the cotangent complex of $\intRep^\square_H(\Gamma_q)$ is perfect and has Tor amplitude in $[-1,0]$. As in the proof of \cref{rslt:rep-stack-for-finite-group} we deduce that $\intRep_H(\Gamma_q)$ is lafp and from the observation on Tor amplitude we deduce that $\intRep_H^\square(\Gamma_q)$ is a local complete intersection; then the same is true for $\intRep_H(\Gamma_q)$. In particular the dualizing complex of $\intRep_H(\Gamma_q)$ is the determinant of the cotangent complex by \cref{rslt:QCoh-solid-on-proper-and-smooth-maps}, which is trivial for any complex of the shape $V \to V^2 \to V$ by basic properties of the determinant. This finishes the proof of (i).

To prove (ii) we now assume that $\Lambda$ is a Cohen-Macaulay classical noetherian ring and $H$ is reductive. We want to show that the lafp affine scheme $\intRep_H^\square(\Gamma_q)$ is classical. By the proof of (i) we know that it is a local complete intersection of virtual dimension $\dim H$ (the dimension of the adjoint representation). By \cref{rslt:classicality-result-for-lci-maps} it is therefore enough to show that the relative dimension of this affine scheme is (at most) $\dim H$. This statement only depends on the fibers of the underlying classical scheme $\mathcal R := \intRep_H^\square(\Gamma_q)^\cl$, so we can reduce to the case that $\Lambda$ is an algebraically closed field. We now refer the reader to \cite[Proposition~2.3.7]{Zhu} for the rest of the proof.
\end{proof}

If $H$ is not reductive then $\intRep_H(\Gamma_q)$ is often not classical and it is important to consider the derived version; we refer the reader to \cite[Remark~2.3.8]{Zhu} for a discussion. Let us now put \cref{rslt:rep-stack-for-finite-group} and \cref{rslt:rep-stack-for-Gamma-q} together:

\begin{prop} \label{rslt:rep-stack-for-semidirect-product-of-Gamma-q-with-finite}
Let $\Lambda$ be a ring, $q = p^r$ a prime power invertible in $\pi_0\Lambda$, $H$ a smooth affine group over $k$ and $\Gamma = Q \rtimes \Gamma_q$, where $Q$ is a finite $p$-group.
\begin{propenum}
    \item $\intRep_H(\Gamma)$ is a local complete intersection over $\Lambda$ and its dualizing complex is trivial.
    \item If $\Lambda$ is a Cohen-Macaulay classical noetherian ring and $H$ is reductive then $\intRep_H(\Gamma)$ is classical.
\end{propenum}
\end{prop}
\begin{proof}
This is proved in \cite[Proposition~2.3.9]{Zhu} in the case that $H$ is reductive (and under stronger assumptions on $\Lambda$), but the argument also proves (i). In fact, (i) can be proved simply by combining the above ideas. Namely, using the projection $\Gamma \surjto \Gamma_q$ we see that for every $\Gamma$-representation $V$ over $\Lambda$ we have $V_\Gamma = (V_Q)_{\Gamma_q}$. Now $V_Q$ is a retract of $V$ by the proof of \cref{rslt:rep-stack-for-finite-group} and $(-)_{\Gamma_q}$ can be computed using an explicit resolution, see the proof of \cref{rslt:rep-stack-for-Gamma-q}. Thus the proof of \cref{rslt:rep-stack-for-Gamma-q} can be applied to $\Gamma$ in order to prove (i).

We sketch the proof of (ii). Using \cref{rslt:classicality-result-for-lci-maps} we are reduced to showing that the fiberwise dimension of $\intRep_H^\square(\Gamma)$ is at most $\dim H$, so we can assume that $\Lambda = k$ is an algebraically closed field. Now fix a representation $\rho_0\colon Q \to H(k)$ and consider the closed substack $\intRep^\square_H(\Gamma)^{\rho_0}$ of those maps $\Gamma \to H(k)$ that restrict to $\rho_0$; we implicitly pass to the underlying classical stack here, as this does not change the dimension. By \cref{rslt:rep-stack-for-finite-group} it is enough to show that the dimension of $\intRep^\square_H(\Gamma)^{\rho_0}$ is $\dim Z_H(\rho_0)$. But note that one can realize this stack as an open substack
\begin{align*}
    \intRep^\square_H(\Gamma)^{\rho_0} \subseteq \intRep^\square_{N_H(\rho_0)}(\Gamma_q),
\end{align*}
where $N_H(\rho_0)$ is the normalizer of $\rho_0$, which is still reductive. We conclude by the proof of \cref{rslt:rep-stack-for-Gamma-q} that the dimension of $\intRep^\square_H(\Gamma)^{\rho_0}$ is $\dim N_H(\rho_0)$, which is also $\dim Z_H(\rho_0)$ because these two groups have the same neutral component.
\end{proof}

In preparation for the definition of the L-parameter stack we now extend the above definitions to topological (i.e. condensed) groups $\Gamma$.

\begin{defn}
Fix a $\Z_\ell$-algebra $\Lambda$ and an affine group $H$ over $\Lambda$.
\begin{defenum}
    \item We view any $\Lambda$-algebra $A$ as a condensed $\Lambda$-algebra via identifying it with $A^{\rm disc} \tensor_{\Z_\ell^{\rm disc}} \Z_\ell$, where $\Z_\ell = \varprojlim_n \Z/\ell^n$ is the condensed ring induced by the $\ell$-adic topology. We view $H(A)$ as a group in condensed anima via $S \mapsto H(A(S))$ for profinite sets $S$. More generally, for any stack $X$ over $\Lambda$, we view $X(A)$ as a condensed anima via $S \mapsto X(A(S))$ for extremally disconnected sets $S$.

    \item \label{def:condensed-rep-stack} Let $\Gamma$ be a group in condensed anima. We denote by $\intRep_H(\Gamma)$ the stack over $\Lambda$ given by
    \begin{align*}
        \intRep_H(\Gamma)\colon A \mapsto \Hom_{\Cond(\Ani)}(*/\Gamma, (*/H)(A)),
    \end{align*}
    for all $\Lambda$-algebras $A$, where we use the implicit condensed structure on $(*/H)(A)$ from (i). Note that if $\Gamma$ is discrete (i.e. the condensed structure is trivial) then $\intRep_H(\Gamma)$ agrees with the definition from \cref{def:rep-stack}, so there is no ambiguity of definitions. We denote $\tau\colon \intRep_H(\Gamma) \to */H$ the natural projection induced by the unit map $1 \to \Gamma$ and we let
    \begin{align*}
        \intRep^\square_H(\Gamma) := \intRep_H(\Gamma) \times_{*/H} *
    \end{align*}
    as in \cref{def:rep-stack}. Then $H$ acts on $\intRep^\square_H(\Gamma)$ and $\intRep^\square_H(\Gamma)/H = \intRep_H(\Gamma)$.
\end{defenum}
\end{defn}

Let us make the definition of the condensed structure on $A$ and $H(A)$ more explicit. For simplicity assume that $A$ is a classical (i.e. non-derived) $\Z_\ell$-algebra. We write $A$ as a filtered colimit of finitely generated $\Z_\ell$-submodules $M_i \subseteq A$. Then the condensed structure on $A$ is obtained by taking the colimit of the $M_i$'s, each equipped with the $\ell$-adic condensed structure. In other words, for a profinite set $S$ the ring $A(S)$ is the ring of maps $S \to A$ that factor over some finite-type $\Z_\ell$-submodule $M \subset A$ and are continuous for the $\ell$-adic topology on $M$. With the condensed structure on $A$ at hand, the condensed structure on $H(A)$ admits a nice description if $H$ is embedded into some $\GL_n$ -- the condensed structure on $\GL_n(A)$ is induced from the condensed structure on $A$ via the matrix coefficient embedding $\GL_n(A) \subseteq A^{n^2}$, and the condensed structure on $H(A)$ then comes from the restriction along $H(A) \injto \GL_n(A)$.

The next result proves the analog of \cref{rslt:basic-properties-of-rep-stack} in the case of condensed groups $\Gamma$. The statements are mostly the same, except that $\intRep^\square_H(\Gamma)$ may not be affine anymore. We therefore keep the statement brief. For the second part of the claim, recall the definition of derivations from \cref{def:derivation-of-stacks}.

\begin{lem}
Fix a $\Z_\ell$-algebra $\Lambda$, an affine group scheme $H$ over $\Lambda$ and a group $\Gamma$ in $\Cond(\Ani)$.
\begin{lemenum}
    \item \label{rslt:framed-intrep-equals-group-hom-cond} For every $\Lambda$-algebra $A$ there is a natural equivalence
    \begin{align*}
        \intRep^\square_H(\Gamma)(A) = \{ \text{condensed group homomorphisms $\Gamma \to H(A)$} \}.
    \end{align*}
    Under this equivalence, the action of $H$ on $\intRep^\square_H(\Gamma)$ is given by conjugation.

    \item \label{rslt:Der-for-conde-rep-stack} Suppose that $H$ is smooth. Then for every $\Lambda$-algebra $A$, connective $A$-module $M$ and point $\rho\colon \Spec A \to \intRep_H(\Gamma)$ there is a natural isomorphism
    \begin{align*}
        \Der_\Lambda(\intRep_H(\Gamma), M)_\rho = C^*(\Gamma, \Ad_\rho \tensor M),
    \end{align*}
    where $\Ad_\rho$ is defined as in \cref{rslt:cotangent-complex-of-rep-stack} and for every $\Gamma$-representation on an $A$-module $V$ we define
    \begin{align*}
        C^*(\Gamma, V) := \varprojlim_{n\in\Delta} \Hom(\Gamma^n, V \tensor_{\Z_\ell^{\rm disc}} \Z_\ell),
    \end{align*}
    where the $\Hom$ is in condensed anima.
\end{lemenum}
\end{lem}
\begin{proof}
For (i), the same argument as in the proof of \cref{rslt:framed-intrep-equals-group-hom} works here as well, replacing $\Ani$ by $\Cond(\Ani)$ where appropriate (note that this is still a topos, hence \cite[Example~5.2.6.13]{lurie-higher-algebra} applies).

For (ii), we introduce the notation $\intMap(Y, */H)$ for a condensed anima $Y$ in a similar way as in \cref{def:condensed-rep-stack}, so that $\intRep_H(\Gamma) = \intMap(*/\Gamma, */H)$. Then $\intMap(-, */H)$ sends colimits to limits, so since $*/\Gamma = \varinjlim_{n\in\Delta} \Gamma^n$, we are reduced to showing that
\begin{align*}
    \Der_\Lambda(\intMap(\Gamma, */H), M)_\eta = \Hom(\Gamma^n, \Ad_\rho \tensor M),
\end{align*}
where $\eta$ is obtained by composing $\rho$ with $\intRep_H(\Gamma) \to \intMap(\Gamma^n, */H)$, induced by the projection $\Gamma^n \to * \to */\Gamma$. Since both sides send colimits in $\Gamma^n$ to limits, we can replace $\Gamma^n$ by an extremally disconnected set $S$. Note that
\begin{align*}
    \intMap(S, */H)(A \oplus M) = (*/H)((A \oplus M)(S)) = (*/H)(A(S) \oplus M(S)),
\end{align*}
where we implicitly use the induced condensed structure on $A$. Using $\Ad_{\rho'} \tensor N = \Der_\Lambda(*/H, N)_{\rho'}$ for any $\Lambda$-algebra $B$, any connective $B$-module $N$ and any map $\rho'\colon \Spec B \to */H$, the desired identity follows easily by considering $B = A(S)$.
\end{proof}

\Cref{rslt:framed-intrep-equals-group-hom-cond} shows that our definition of the $\intRep_H(\Gamma)$ is compatible with the definition in \cite[\S VIII]{FS}. In order to make the above results available, we need to relate $\intRep_H(\Gamma)$ to a representation stack with \emph{discrete} $\Gamma$. This is achieved by the next result, that is also used implicitly in the proof of \cite[Theorem~VIII.1.3]{FS} (in the non-derived setting) and is partially proved in \cite[Lemma~3.9]{Zhu}. Let us first introduce some notation.

\begin{defn}
Fix a ring $\Lambda$. A \emph{linear algebraic group} over $\Lambda$ is an affine flat group scheme $H$ over $\Lambda$ such that there is some $n \ge 0$ and a group homomorphism $H \to \GL_n$ over $\Lambda$ that is a closed immersion.
\end{defn}

\begin{defn} \label{def:discretization-of-W-F}
Let $P_F \subset W_F$ denote the wild inertia subgroup, so that $W_F / P_F \isom \hat \Z^p \rtimes \Z$, where $\hat \Z^p := \prod_{\ell \ne p} \Z_\ell$. A \emph{discretization} of $W_F$ is an embedding
\begin{align*}
    \iota\colon \Gamma_q \injto W_F/P_F
\end{align*}
with $\Gamma_q = \langle \sigma, \tau \rangle$ as in \cref{def:Gamma-q} (for $q$ the number of elements in the residue field of $F$) such that $\iota(\tau)$ is a generator of the tame inertia and $\iota(\sigma)$ is a lift of Frobenius. Given $\iota$, we denote
\begin{align*}
    W_{F,\iota} := W_F \times_{W_F/P_F} \Gamma_q.
\end{align*}
This is a locally profinite group (equipped with the fiber product topology).
\end{defn}

\begin{prop} \label{rslt:Rep-W-F-colimit-of-mapping-stacks}
Fix a $\Z_\ell$-algebra $\Lambda$, a smooth linear algebraic group $H$ over $\Lambda$ and a discretization $\iota$ of the Weil group $W_F$. 
\begin{propenum}
    \item \label{rslt:Rep-W-F-equal-to-discretized-version} The map $W_{F,\iota} \to W_F$ induces an equivalence
    \begin{align*}
        \intRep_H(W_F) \isoto \intRep_H(W_{F,\iota})
    \end{align*}
    of stacks over $\Lambda$. The same is true for $W_{F,\iota}/P \to W_F/P$ for any open subgroup $P \subseteq P_F$.

    \item \label{rslt:decompose-Rep-W-F-as-colimit} We have
    \begin{align*}
        \intRep_H(W_F) = \varinjlim_{P \subseteq P_F} \intRep_H(W_F/P),
    \end{align*}
    where $P \subseteq P_F$ runs through the open normal subgroups of the wild inertia group. Moreover for $P \subseteq P' \subseteq P_F$, the map $\intRep_H(W_F/P) \injto \intRep_H(W_F/P')$ is an open and closed immersion.
\end{propenum}
\end{prop}
\begin{proof}
To prove (i) we follow the argument in \cite[Lemma~3.1.8]{Zhu}: We first check the claimed isomorphism on the underlying classical stacks, i.e. on $A$-valued points for classical $\Lambda$-algebras $A$. We can also reduce to the framed representation stacks. By choosing a closed embedding $H \to \GL_n$ we can further reduce to the case $H = \GL_n$ (note that this closed embedding induces injections on classical points). We are thus left with proving that the restriction map
\begin{align*}
    \{ \text{cond. group hom. $W_F \to \GL_n(A)$} \} \to \{ \text{cond. group hom. $W_{F,\iota} \to \GL_n(A)$} \}
\end{align*}
is a bijection of sets. This result appears in several places in the literature, see e.g. \cite[Lemma~3.1.8]{Zhu} or the proof of \cite[Theorem~VIII.1.3]{FS}.

To prove the derived version of (i), we follow the strategy in the last part of \cite[Lemma~3.1.8]{Zhu}. We can check the isomorphism on framed representation stacks and then both sides of the isomorphism commute with limits in the ring that we plug in. Using the Postnikov limit $A = \varprojlim_n \tau_{\le n} A$ for any $\Lambda$-algebra $A$, we reduce to the case that $A$ is truncated. We now induct on the number $n$ such that $A$ is $n$-truncated. Using \cref{rslt:truncation-of-rings-via-square-zero-extension} we are then reduced to showing that for a truncated ring $A$ and any connective $A$-module $M$, the map
\begin{align*}
    \intRep_H(W_F)(A \oplus M) \to \intRep_H(W_{F,\iota})(A \oplus M)
\end{align*}
is an isomorphism, where we can assume that the same is true on $A$-points. This isomorphism can be checked fiberwise over a fixed $A$-point $\rho$ and then reduces to the claim that the natural map
\begin{align*}
    \Der_\Lambda(\intRep_H(W_F), M) \to \Der_\Lambda(\intRep_H(W_{F,\iota}), M)
\end{align*}
is an isomorphism for every static $A$-module $M$. By \cref{rslt:Der-for-conde-rep-stack} this reduces to showing that the map $C^*(W_F, \Ad_\rho \tensor M) \to C^*(W_{F,\iota}, \Ad_\rho \tensor M)$ is an isomorphism. Let us denote $V := \Ad_\rho \tensor M$; we can then immediately reduce to the case that $A = \Lambda$ and $V$ is a static $\Z_\ell$-modules with continuous $W_F$-action (for the topology on $V \tensor_{\Z_\ell^{\rm disc}} \Z_\ell$). Then $C^*(W_F, V)$ and $C^*(W_{F,\iota}, V)$ compute continuous group cohomology. To prove the desired isomorphisms, we can first pass to $P_F$-cohomology in order to reduce the claim to showing that the map $C^*(W_F/P_F, V) \to C^*(W_{F,\iota}, V)$ is an isomorphism. This further reduces to showing that the map $C^*(\hat \Z^p, V) \to C^*(\Z[1/p], V)$ is an isomorphism. Now $\hat \Z^p$ is compact and $\Z[1/p]$ has finite cohomological dimension (being a sequential colimit of copies of $\Z$), which shows that we can pull out filtered colimits in $V$ in both cohomologies. This reduces us to the case that $V$ is a finite $\Z_\ell$-module, which we leave to the reader. This finishes the proof of (i).

To prove the colimit description in (ii) we use the same strategy as in the proof of (i): Firstly, it is easily verified on classical rings; secondly, since both sides of the claimed equivalence commute with finite limits of rings, we can reduce to checking an isomorphism on derivations and hence on group cohomology, which is easily verified using the finite cohomological dimension of $W_{F,\iota}$ on $\Z_\ell$-modules (cf. the proof of \cite[Proposition~2.42]{Zhu} for a similar argument). It remains to verify the isomorphism on $A$-valued points for a $\Lambda$-algebra $A$ that is not truncated; this is not trivial because we need to show that the filtered colimit on the right-hand side commutes with the limit induced by $A = \varprojlim_n \tau_{\le n} A$. Denoting $X := \intRep_H(W_F)$ it is enough to show that there is some constant $c \ge 0$ such that for all $n \ge 0$ the map $X(\tau_{\le n+1} A) \to X(\tau_{\le n} A)$ of anima is $(n-c)$-connective (i.e. surjective on $\pi_0$ and an isomorphism on $\pi_i$ for $i < n-c$ and surjective on $\pi_{n-c}$). Using \cref{rslt:truncation-of-rings-via-square-zero-extension} and the fact that connectiveness is preserved under base-change we reduce to showing that $X(\tau_{\le n} A) \to X(\tau_{\le n} A \oplus \pi_{n+1} A[n+1])$ is $(n-c)$-connective. This map is a section to the canonical map $X(\tau_{\le n} A \oplus \pi_{n+1}A[n+1]) \to X(\tau_{\le n} A)$, so by \cite[Corollary~04H9]{Kerodon} it is enough to show that the latter map is $(n-c+1)$-connective. By \cite[Corollary~051F]{Kerodon} this can be checked fiberwise and thus reduces to a similar statement about $\Der(X, \pi_{n+1} A[n+1])_\eta$ for varying $\eta$. By the proof of (i) this follows from the finite cohomological dimension of $W_F$ on $\Z_\ell$-modules.

Now fix open normal subgroups $P \subseteq P' \subseteq P_F$. It remains to show that the map $f\colon \intRep_H(W_F/P') \to \intRep_H(W_F/P)$ is an open and closed embedding. By (i) we can replace $W_F$ by $W_{F,\iota}$ and can also pass to framed representation stacks, which are affine schemes by \cref{rslt:framed-intrep-is-affine}. On the underlying classical schemes the claim is for example shown in the first part of the proof of \cite[Theorem~VIII.1.3]{FS} (this result is only stated for $\Lambda = \Z_\ell$, but using a closed embedding $H \to \GL_n$ one can reduce to the case $H = \GL_n$ and then the situation is base-changed from $\Lambda = \Z_\ell$). This implies that $f$ is a closed immersion, and to show that it is an open immersion it is enough to show that $f$ is flat. By \cref{rslt:cotangent-complex-for-geometric-map-and-smoothness} it is enough to show that $L_f$ is trivial. Denoting $\Gamma := W_{F,\iota}/P$ and $\Gamma' := W_{F,\iota}/P'$ and using the fiber sequence \cref{rslt:fiber-sequence-for-cotangent-complex} we are reduced to showing that $f^* L_{\intRep_H(\Gamma)/\Lambda} = L_{\intRep_H(\Gamma')/\Lambda}$. By \cref{rslt:cotangent-complex-of-rep-stack} this is equivalent to $V_\Gamma = V_{\Gamma'}$ for every $\Lambda$-algebra $A$ and map $\rho\colon \Spec A \to \intRep_H(\Gamma')$, where $V = \Ad_\rho^\vee$. The identity $V_\Gamma = V_{\Gamma'}$ further reduces to showing that the natural map $V_{P'/P} \to V$ is an isomorphism, where the $P'/P$-action is trivial. But $P'/P$ is a finite $p$-group and $p$ is invertible in $\pi_0 \Lambda$, which implies the desired isomorphism (note that the claim immediately reduces to $\Lambda = \Z_\ell$).
\end{proof}

\subsection{The stack of L-parameters}

In the beginning of this section we have fixed a finite extension $F$ of $\Q_p$ and denoted by $W_F$ the Weil group of $F$. The main geometric object on the spectral side of the categorical local Langlands correspondence is the stack $\Par_G$ of $L$-parameters of $G$. This stack can be defined over $\Z[1/p]$, where it requires the choice of a discretization of $W_F$. Over $\Z_\ell$ it admits a completely canonical definition, which we will focus on (in fact, we usually even restrict to $\Qellbar$). In the case that $G$ is split, $\Par_G$ is the stack that sends a $\Z_\ell$-algebra $A$ to the anima of continuous maps $W_F \to \hat G(A)$ modulo $G(A)$-conjugation, where $\hat G$ is the dual group of $G$ over $\Z_\ell$. If $G$ is not split then $\hat G$ comes equipped with a $W_F$-action, which has to be taken into account.

\begin{rmk}
The $L$-parameter stack has been studied extensively in the literature, see e.g. \cite{Zhu}, \cite{DHKM-ParG}, \cite[\S VIII]{FS}. The first two reference also discuss the stack over $\Z[1/p]$ instead of $\Z_\ell$, while the last reference provides the cleanest description over $\Z_\ell$ using condensed methods. Recently Scholze \cite{FS-motivic} presented a new construction of $\Par_G$ over $\Z[1/p]$ that avoids any auxiliary choices by letting $\Par_G$ be defined over a certain \enquote{motivic base stack}.
\end{rmk}

Before we define $\Par_G$, let us first introduce the following more general definition of the \enquote{stack of local systems}. We make use of the continuous representation stack from \cref{def:condensed-rep-stack}.

\begin{defn}
Fix a ring $\Lambda$ and let $H$ be a smooth linear algebraic group scheme over $\Lambda$, equipped with a smooth action by $W_F$. Let $Q$ be a finite quotient of $W_F$ through which this action factors.
\begin{defenum}
    \item \label{def:LocSys} We define
    \begin{align*}
        \LocSys_{H,\Lambda} := \intRep_{H \rtimes Q}(W_F) \times_{\intRep_Q(W_F)} *,
    \end{align*}
    where the map $* \to \intRep_Q(W_F)$ is induced by the projection $W_F \surjto Q$.

    \item Restriction along $1 \to W_F$ induces a map $\LocSys_{H,\Z_\ell} \to */H$ and we define
    \begin{align*}
        \LocSys^\square_{H,\Lambda} := \LocSys_{H,\Lambda} \times_{*/H} *
    \end{align*}
    One sees easily that this stack can also be described as $\intRep^\square_{H\rtimes Q}(W_F) \times_{\intRep^\square_Q(W_F)} *$.
\end{defenum}
\end{defn}

We now come to the definition of the $L$-parameter stack. We follow the construction in \cite[\S VIII]{FS}, but also allow parabolic subgroups, in which case the derived structure on the stack becomes important; here we use the results from the previous subsection, which rely in Zhus's work \cite{Zhu}. The following definition recalls some basics from \cite[\S1, \S2]{Borel-LFunctions}.

\begin{defn} \label{def:L-parameter-stack}
Let $G$ be a connected reductive group over $F$ and $\Lambda$ a $\Z_\ell$-algebra. Moreover, fix a (not necessarily split) maximal torus $T \subseteq G$ and a Borel $B_{\overline F} \subset G_{\overline F}$ containing $T_{\overline F}$. By the classification of reductive groups, the triple $(G_{\overline F}, B_{\overline F}, T_{\overline F})$ corresponds to a based root datum $\Psi_0$ and the group $G$ induces an action of $\Gamma_F = \Gal(\overline F/F)$ on $\Psi_0$.
\begin{defenum}
    \item \label{def:ParG} By dualizing $\Psi_0$ and passing to the corresponding (split) reductive group over $\Lambda$ (base-changed from $\Z$), we obtain a pinned reductive group scheme $\hat G$ over $\Lambda$ together with an action of $\Gamma_F$. We refer to $\hat G$ as the \emph{dual group} of $G$. Using the restricted action of $W_F \subset \Gamma_F$ on $\hat G$ we define the \emph{stack of $L$-parameters for $G$} as
    \begin{align*}
        \Par_{G,\Lambda} := \LocSys_{\hat G,\Lambda}.
    \end{align*}
    There is a natural map $\Par_{G,\Lambda} \to */\hat G$ and the \emph{framed $L$-parameter stack} $\Par_{G,\Lambda}^\square = \LocSys_{\hat G,\Lambda}^\square$ is the pullback of this map along $* \to */\hat G$.

    \item \label{def:ParP} Parabolic subgroups $P \subseteq G$ containing $T$ correspond to parabolic subgroups $\hat P \subseteq \hat G$ which are stable under the $\Gamma_F$-action. We refer to $\hat P$ as the \emph{dual parabolic subgroup} and define
    \begin{align*}
        \Par_{P,\Lambda} := \LocSys_{\hat P,\Lambda}.
    \end{align*}
    As before, we also define the framed version $\Par^\square_{P,\Lambda}$.
\end{defenum}
In the case $\Lambda = \Qellbar$ we drop $\Lambda$ from the notation, so e.g. $\Par_G = \Par_{G,\Qellbar}$.
\end{defn}

Let us unpack the definition. Fix a finite quotient $Q$ of $W_F$ through which the action on $\hat G$ factors. By \cref{rslt:framed-intrep-equals-group-hom-cond} for every $\Lambda$-algebra $A$ the anima $\Par^\square_{G,\Lambda}(A)$ consists of condensed group homomorphisms $W_F \to \hat G(A) \rtimes Q$ together with an identification of their projection to $Q$ with the canonical map $W_F \to Q$. The group $\hat G$ acts via conjugation on $\Par^\square_{G,\Lambda}$ and we have $\Par_{G,\Lambda} = \Par^\square_{G,\Lambda}/\hat G$.

\begin{rmk} \label{rmk:Par-G-only-depends-on-qsplit-form}
It follows immediately from the construction of the $\Gamma_F$-action on $\Psi_0$ that two inner forms $G$ and $G'$ induce the same Galois action on $\hat G$ (see \cite[\S1.3]{Borel-LFunctions}) and thus $\Par_G = \Par_{G'}$. By \cite[Proposition~7.2.12]{Conrad-ReductiveGroupSchemes} we are therefore free to assume that $G$ is quasisplit, i.e. contains a Borel defined over $F$.
\end{rmk}

\begin{exmpl}
In the setting of \cref{def:L-parameter-stack} suppose that $G$ is split, e.g. $G = \GL_n$. Then the action of $W_F$ on $\hat G$ is trivial so we can take $Q = 1$ and hence $\Par_{G,\Lambda} = \intRep_{\hat G}(W_F)$ is the stack of $\hat G$-representations of $W_F$.
\end{exmpl}

Using the results in the previous subsection we can now easily deduce the following general results on the geometry of the $L$-parameter stack:

\begin{thm}
Fix a $\Z_\ell$-algebra $\Lambda$ and a connected reductive group $G$ over $F$ with parabolic subgroup $P \subseteq G$.
\begin{thmenum}
    \item $\Par^\square_{G,\Lambda}$ and $\Par^\square_{P,\Lambda}$ are schemes and a local complete intersections of dimension $\dim G$ resp. $\dim P$ over $\Lambda$. They can be written as an infinite disjoint union of affine schemes, each of which admits a closed immersion into $\hat G^n$ resp. $\hat P^n$ for large enough $n$.
    
    \item \label{rslt:connected-components-of-Par-G-are-affines-mod-G} $\Par_{G,\Lambda}$ and $\Par_{P,\Lambda}$ are geometric stacks and a local complete intersection of relative dimension $0$ over $\Lambda$. They can be written as a disjoint union of stacks of the form $X/\hat G$ respectively $X/\hat P$ for varying affine schemes $X$ over $\Lambda$, each of which admits a closed immersion into $\hat G^n$ resp. $\hat P^n$ for large enough $n$.

    \item \label{rslt:Par-G-is-classical} Suppose that $\Lambda$ is a Cohen-Macaulay classical noetherian ring. Then $\Par^\square_{G,\Lambda}$ and $\Par_{G,\Lambda}$ are classical.
\end{thmenum}
\end{thm}
\begin{proof}
Since $\Par_{G,\Lambda} = \Par^\square_{G,\Lambda} / \hat G$ (and similarly for $P$), the claims (i) and (ii) are equivalent, so we only prove (ii) and (iii). Let $Q$ be a finite quotient of $W_F$ through which the $W_F$-action on $\hat G$ factors. Note that $\intRep^\square_Q(W_F)$ is a disjoint union of points (it is a constant stack). In particular the map $* \to \intRep^\square_Q(W_F)$ is a clopen immersion and hence the claims about $\Par_{G,\Lambda}$ and $\Par_{P,\Lambda}$ immediately reduce to similar claims about $\intRep_{\hat G \rtimes Q}(W_F)$ and $\intRep_{\hat P \rtimes Q}(W_F)$. The second part of (ii) (apart from the closed immersions) now follows immediately from \cref{rslt:decompose-Rep-W-F-as-colimit} and \cref{rslt:framed-intrep-is-affine}. By \cref{rslt:Rep-W-F-equal-to-discretized-version} the remaining claims reduce to $\intRep_H(W_{F,\iota}/P)$, where $W_{F,\iota}$ is a discretization of $W_F$, $H$ is either $\hat G \rtimes Q$ or $\hat P \rtimes Q$, and $P \subseteq P_F$ is a normal open subgroup. But $W_{F,\iota}/P$ is of the form in \cref{rslt:rep-stack-for-semidirect-product-of-Gamma-q-with-finite}, so we immediately conclude by that result. We also obtain the desired closed immersions into some power of $\hat G$ or $\hat P$ using the description of the stacks in terms of group homomorphisms and the fact that $W_{F,\iota}/P$ is finitely generated.
\end{proof}

From now on we will focus on the case $\Lambda = \Qellbar$ and study the geometry of $\Par_G$. Most of the results below generalize to more general $\Lambda$, but often with additional technical subtleties; we refer the reader to \cite{DHKM-ParG} for a comprehensive treatment. Let us fix the setting for the rest of the section:

\begin{setting} \label{set:reductive-group-over-F}
Let $G$ be a quasisplit connected reductive group over $F$ with dual group $\hat G$ over $\Qellbar$. Let $Q$ be a finite quotient of $W_F$ over which the action of $W_F$ on $\hat G$ factors.
\end{setting}

By \cref{rslt:Par-G-is-classical} $\Par_G$ is a classical algebraic stack and in particular all current results in the literature (which mostly work in the non-derived setting) apply to it. By \cref{rslt:Par-G-is-classical} $\Par^\square_G$ is a classical scheme over $\Qellbar$ whose $\Qellbar$-points are given by $\Qellbar$-valued $L$-parameters, i.e. continuous group homomorphisms $W_F \to \hat G(\Qellbar) \rtimes Q$ whose projection to $Q$ is the projection $W_F \to Q$. To study the geometry of $\Par_G$, it is very useful to single out a special class of $L$-parameters. We first need the following auxiliary definition:

\begin{defn}
Let $H$ be a (not necessarily connected) reductive group over $\Qellbar$. An \emph{R-parabolic subgroup} of $H$ is a closed subgroup $P \subseteq H$ such that $P$ is parabolic (i.e. $H/P$ is proper) and the induced map $\pi_0 P \to \pi_0 H$ is surjective. An \emph{R-Levi subgroup} of $H$ is the Levi group of an R-parabolic subgroup.
\end{defn}

Our definition of R-parabolic subgroups is taken from \cite[\S VIII.3.1]{FS} and \cite[\S 3.3]{Borel-LFunctions}, although the terminology is from \cite[\S 6]{CompleteReducibility}. The definitions in these two references agree by \cite[Proposition 6.1]{CompleteReducibility}. Note that if $H = \hat G \rtimes Q$ then the R-parabolic subgroups of $H$ are up to conjugation of the form $\hat P \rtimes Q$ for a parabolic subgroup $\hat P \subseteq \hat G$. We now have the following class of $L$-parameters:

\begin{lem} \label{rslt:def-of-semisimple-L-param}
In \cref{set:reductive-group-over-F}, let $\phi\colon W_F \to \hat G(\Qellbar) \rtimes Q$ be an $L$-parameter. Then the following are equivalent:
\begin{lemenum}
    \item The image of $\phi$ in $\hat G(\Qellbar) \rtimes Q$ is \emph{completely reducible}, i.e. for every R-parabolic subgroup $P \subseteq \hat G \rtimes Q$ that contains the image of $\phi$, there is an R-Levi subgroup of $P$ that contains the image of $\phi$.

    \item The kernel of $\phi$ is open in $W_F$ and $\phi$ is \emph{Frobenius-semisimple}, i.e. for some (equivalently any) lift of the Frobenius $\sigma \in W_F$ the element $\phi(\sigma) \in \hat G(\Qellbar) \rtimes Q$ is semisimple.

    \item The $\hat G$-orbit of $\phi$ in the scheme $\Par_G^\square$ is closed.
\end{lemenum}
\end{lem}
\begin{proof}
The equivalence of (i) and (iii) is proved in \cite[Theorem~4.13]{DHKM-ParG} and \cite[Proposition~VIII.3.2]{FS}. By \cite[Remark~6.9(ii)]{DHKM-ParG} (ii) implies (iii). By \cite[Corollary~4.16]{DHKM-ParG} (i) implies the first part of (ii) and then \cite[Remark~6.9(ii)]{DHKM-ParG} shows that (iii) implies the remaining part of (ii).
\end{proof}

\begin{defn}
In \cref{set:reductive-group-over-F} an $L$-parameter $\phi\colon W_F \to \hat G(\Qellbar) \rtimes Q$ is called \emph{semisimple} if it satisfies the equivalent conditions of \cref{rslt:def-of-semisimple-L-param}.
\end{defn}

We next introduce a notation for the centralizer of an $L$-parameter, which plays an important role in the local Langlands correspondence. Recall the definition of the $\intIsom$-stack from \cref{def:intIsom-stack}.

\begin{defn}
In \cref{set:reductive-group-over-F} let $\phi\colon W_F \to \hat G(\Qellbar) \rtimes Q$ be an $L$-parameter, viewed as a point of $\Par^\square_G$ and $\Par_G$. We denote by
\begin{align*}
    S_\phi := (\intIsom_{\Par_G}(\phi))^\cl = (\hat G \times_{\Par^\square_G} \{ \phi \})^\cl
\end{align*}
the \emph{centralizer} of $\phi$ and by $i_\phi\colon */S_\phi \to \Par_G$ the associated closed immersion with image $\phi$ from \cref{rslt:intIsom-is-group-and-classifying-stack-embeds}. Here the implicit map $\hat G \to \Par^\square_G$ is the orbit map (as in \cref{def:orbit-map}).
\end{defn}

By construction $S_\phi$ is the closed subgroup of $\hat G$ whose $\Qellbar$-valued points are given by those elements $g \in \hat G(\Qellbar)$ such that $g \phi(x) g^{-1} = \phi(x)$ for all $x \in W_F$.

Recall that $\Par_G = \Par^\square_G/\hat G$ for the conjugation action by $\hat G$. Here the quotient is a stack quotient and the resulting stack is quite far from being a scheme (e.g. it is lci of dimension $0$). For the study of the geometry of $\Par_G$ it is very convenient to also consider the scheme-theoretic quotient:

\begin{defn} \label{def:coarse-moduli-space}
In \cref{set:reductive-group-over-F}, we denote by
\begin{align*}
    X^\spec_G := \Par^\square_G \sslash \hat G
\end{align*}
the scheme-theoretic quotient of $\Par^\square_G$ by the conjugation action of $\hat G$. In other words, if $C$ is a connected component of $\Par_G$ then in $X^\spec_G$ we have the corresponding connected component $\Spec(H^0(C, \calO_C))$. Since the map $\Par^\square_G \to X^\spec_G$ is $\hat G$-equivariant (for the trivial $\hat G$-action on the target), it induces a natural cover
\begin{align*}
    q\colon \Par_G \surjto X^\spec_G.
\end{align*}
We denote by $\mathfrak Z_G^\spec := \calO(X_G^\spec) = H^0(\Par_G, \calO_{\Par_G})$ the global sections and call it the \emph{spectral Bernstein center} of $G$.
\end{defn}

By \cite[Example~12.9]{AlperGood} $q\colon \Par_G \to X^\spec_G$ is a good moduli space in the sense of \cite[Definition~4.1]{AlperGood}. Moreover, by \cite{Haines} and \cite{DHKM-ParG} there is an explicit description of $X^\spec_G$. We summarize its basic properties in the following result.

\begin{prop} \label{rslt:properties-of-X-spec-G}
Fix \cref{set:reductive-group-over-F}.
\begin{propenum}
    \item The $\Qellbar$-points of $X^\spec_G$ correspond to the semisimple $L$-parameters of $G$ up to conjugation.

    \item \label{rslt:explicit-description-of-connected-components-in-X-spec-G} Let $\phi\colon W_F \to \hat G(\Qellbar) \rtimes Q$ be a semisimple $L$-parameter of $G$ and $C^\coarse \subseteq X^\spec_G$ the associated connected component. Choose a minimal R-parabolic subgroup $M \subseteq \hat G \rtimes Q$ over which $\phi$ factors. Then the $W_F$-action on $\hat G \rtimes Q$ restricts to a $W_F$-action on $Z(M^\circ)$ and for any lift of Frobenius $\sigma \in W_F$ we can define the torus $Y_\phi := (Z(M^\circ)^{I_F})^\circ_\sigma$. There is a finite group $K_\phi$ acting on $Y_\phi$ and a natural isomorphism
    \begin{align*}
        C^\coarse = Y_\phi \sslash K_\phi.
    \end{align*}
\end{propenum}
\end{prop}
\begin{proof}
For (i) note that the $\Qellbar$-points of $X^\spec_G$ are the closed $\hat G$-orbits in $\Par^\square_G$ (see \cite[Theorem~4.16(iv)]{AlperGood}), so we conclude by \cref{rslt:def-of-semisimple-L-param}. For part (ii) see the discussion in \cite[\S 6.3]{DHKM-ParG}.
\end{proof}

\begin{defn}
In \cref{set:reductive-group-over-F}, given an $L$-parameter $\phi$ of $G$, viewed as a $\Qellbar$-point in $\Par^\square_G$, we call its image in $X^\spec_G$ the \emph{semisimplification} of $\phi$.
\end{defn}

\Cref{rslt:properties-of-X-spec-G} gives us a good description of $X^\spec_G$ and thus allows to study $\Par_G$ via its fibers along $q\colon \Par_G \to X^\spec_G$.

In general the geometry of $\Par_G$ is subtle, but the situation is much better around the following class of $L$-parameters:

\begin{defn}
In \cref{set:reductive-group-over-F} fix an $L$-parameter $\phi\colon W_F \to \hat G(\Qellbar) \rtimes Q$.
\begin{defenum}
    \item We say that $\phi$ is \emph{supercuspidal} if it does not factor through any R-parabolic subgroup of $\hat G \rtimes Q$.

    \item We say that $\phi$ is \emph{generous} if it is semisimple and $\Par_G \times_{X^\spec_G} \{ \phi \} = */S_\phi$. 
\end{defenum}
\end{defn}

\begin{defn} \label{def:torus-of-unramified-central-1-cocycles}
We denote 
\begin{align*}
    X^{nr}(G) := (Z(M^\circ)^{I_F})^\circ_\sigma,
\end{align*}
hence $X^{nr}(G) = Y_\phi$ (notation as in \cref{rslt:explicit-description-of-connected-components-in-X-spec-G}) for any supercuspidal $L$-parameter $\phi$. Note that this torus is also identified with the group of unramified central 1-cocycles valued in $\hat G \rtimes Q$, i.e. 1-cocycles $W_F/I_F \to Z(\hat G(\Qellbar)) \rtimes Q$ (see \cite[Eq.~7.2]{Mishra-Bernstein}).
\end{defn}

Generous $L$-parameters were first introduced in \cite[\S2.1]{Beijing} where they are shown to have excellent properties. We summarize some of them as follows.

\begin{prop}
Fix \cref{set:reductive-group-over-F} and let $\phi\colon W_F \to \hat G(\Qellbar) \rtimes Q$ be an $L$-parameter.
\begin{propenum}
    \item \label{rslt:connected-component-of-supercuspidal-rep} If $\phi$ is supercuspidal then it is generous and the connected component $C$ of $\phi$ in $\Par_G$ is smooth and contains only supercuspidal $L$-parameters. Moreover, if $C^\coarse$ denotes the connected component of $\phi$ in $X^\spec_G$ then $C^\coarse = X^{\rm nr}(G)/K_\phi$ for a finite subgroup $K_\phi \subseteq X^{\rm nr}(G) = Y_\phi$, and there is a natural cartesian square in $\Stk_{\Qellbar}$
    \begin{equation*}\begin{tikzcd}
        X^{\rm nr}(G) \times */S_\phi \arrow[r] \arrow[d] & C \arrow[d] \\
        X^{\rm nr}(G) \arrow[r] & C^\coarse
    \end{tikzcd}\end{equation*}
    In particular $C$ is an étale gerbe over the torus $C^\coarse$, banded by $S_\phi$.
    
    \item If $\phi$ is generous then $\Par_G$ is smooth in an open neighborhood of $\phi$ and $q\colon \Par_G \to X^\spec_G$ is flat in an open neighborhood $U \subseteq X^\spec_G$ of $\phi$.
\end{propenum}
\end{prop}
\begin{proof}
Part (ii) is proved in \cite[Lemma~3.3]{HJnote} and \cite[Proposition~3.5]{HJnote}. To prove (i), assume now that $\phi$ is supercuspidal. By the explicit description of the connected components of $X^\spec_G$ in \cref{rslt:explicit-description-of-connected-components-in-X-spec-G} (more precisely, see \cite[\S 5.3]{Haines}) we see that $C^\coarse = Y_\phi \sslash K_\phi$, where $Y_\phi = X^{\rm nr}(G)$ and $K_\phi$ is a finite group acting on $Y_\phi$. Moreover, in the notation of \cite[\S5.3]{Haines} the group $W^{\hat G}_{[\phi]_{\hat G}}$ is trivial, so $K_\phi = \mathrm{stab}_\phi$ and this is by construction a subgroup of $X^{\rm nr}(G)$. Therefore $X^{\rm nr}(G) \sslash K_\phi =X^{\rm nr}(G) / K_\phi$ and this is a torus. In particular the $L$-parameters in $C^\coarse$ differ from $\phi$ only via conjugation by a central element in $\hat G$, hence are supercuspidal.

Note that $\phi$ is strongly semisimple in the sense of \cite[Definition~3.2(1)]{HJnote} (there are no other $L$-parameters with the same semisimplification). The same is then true for any $L$-parameter in $C^\coarse$, as they are all supercuspidal. In particular the fibers of $\pi\colon C \to C^\coarse$ all contain only one point. This shows that all $L$-parameters in $C$ are supercuspidal.

It remains to show the existence of the cartesian square in (i); the other claims of (i) follow immediately from it. The main idea is to construct a canonical action of $X^{\rm nr}(G)$ on $\Par_G$ as follows. The $Q$-action on $\hat G$ restricts to a $Q$-action on $Z(\hat G)$, so we can form $\LocSys_{Z(\hat G)}$. The assignment $H \mapsto \LocSys_H$ is functorial in linear algebraic groups over $\Qellbar$ equipped with $Q$-action (i.e. smooth affine groups over $*/Q$) and commutes with products. Moreover, $Z(\hat G)$ acts on $\hat G$ in that category given by multiplication. Altogether this implies that $\LocSys_{Z(\hat G)}$ is a commutative group stack acting on $\Par_G$ by pointwise multiplication. Note that $\LocSys^\square_{Z(\hat G)} \to \LocSys_{Z(\hat G)}$ is a map of commutative group stacks (because it is a pullback of $* \to */Z(\hat G)$). We can further replace $Z(\hat G)$ by $(Z(\hat G)^{I_F})^\circ$ in order to obtain an action
\begin{align*}
    \LocSys^\square_{(Z(\hat G)^{I_F})^\circ} \times \Par_G \to \Par_G.
\end{align*}
Now $\LocSys^\square_{(Z(\hat G)^{I_F})^\circ}$ is the classical scheme of 1-cocycles $W_F \to (Z(\hat G)^{I_F})^\circ$. As observed in \cite[Eq.~7.2]{Mishra-Bernstein} the closed subspace where the 1-cocycle is trivial on $I_F$ is isomorphic to $X^{\rm nr}(G)$. Altogether we deduce that there is a natural action of $X^{\rm nr}(G)$ on $\Par_G$. By passing to coarse moduli spaces we also obtain an action of $X^{\rm nr}(G)$ on $X^\spec_G$ together with a commuting square
\begin{equation*}\begin{tikzcd}
    X^{\rm nr}(G) \times \Par_G \arrow[r] \arrow[d] & \Par_G \arrow[d]\\
    X^{\rm nr}(G) \times X^\spec_G \arrow[r] & X^\spec_G
\end{tikzcd}\end{equation*}
of stacks over $\Qellbar$. This square is necessarily cartesian: By precomposing the horizontal maps with the compatible \emph{isomorphisms} $(g, x) \mapsto (g, g^{-1} x)$ we reduce to the case that they are given by projection to the second coordinate, where the claim is clear.

With the above cartesian square at hand, it is now easy to finish the argument. First note that if $\phi$ is generous then we conclude by observing that the left vertical map along $(\id, \phi)\colon X^{\rm nr}(G) \to X^{\rm nr}(G) \times X^\spec_G$ yields $X^{\rm nr}(G) \times */S_\phi$. We then easily deduce from the obtained description of $C$ that all $L$-parameters in $C$ are generous. Thus even if we do not a priori know that $\phi$ is generous, it is enough to show that \emph{some} $\phi' \in C$ is generous -- but this follows from \cite[Lemma~0578]{stacks-project} (cf. \cite[Remark~3.7(1)]{HJnote}).
\end{proof}

\begin{exmpl} \label{exmpl:Par-G-for-GL-n}
Suppose that $G = \GL_n$. Then $L$-parameters of $G$ are the same as continuous homomorphisms $\phi\colon W_F \to \GL_n(\Qellbar)$, i.e. $n$-dimensional continuous representations of $W_F$. Such a $\phi$ is supercuspidal if and only if it is irreducible, and it is semisimple if and only if it decomposes into a direct sum of irreducible representations. If $\phi$ is irreducible then $S_\phi = Y_\phi = Z(\GL_n) \isom \mathbb G_m$ and $K_\phi = 1$, hence by \cref{rslt:connected-component-of-supercuspidal-rep} the connected component of $\Par_G$ that contains $\phi$ is isomorphic to $\mathbb G_m \times */\mathbb G_m$ (see \cite[Proposition~3.6]{Nguyen} for an explicit computation).

In fact \cite[Proposition~3.6]{Nguyen} (together with \cite[Lemma~3.7]{Nguyen}) proves something stronger: Suppose $\phi$ decomposes as a direct sum $\phi = \phi_1 \oplus \dots \oplus \phi_r$ of irreducible representations such that for $i \ne j$ there is no unramified character $\chi\colon W_F \to W_F/I_F \to \Qellbar^\times$ with $\phi_i \isom \phi_j \tensor \chi$; then $\phi$ is generous and the connected component $C \subseteq \Par_G$ containing $\phi$ is isomorphic to $\mathbb G_m^r \times */\mathbb G_m^r$. Here an element $(t_i)_i \in \mathbb G_m^r$ corresponds to a family $(\chi_i)_i$ of unramified characters (since $W_F/I_F \isom \Z$) and the corresponding point of $C$ is given by $\chi_1 \tensor \phi_1 \oplus \dots \oplus \chi_r \tensor \phi_r$.
\end{exmpl}

\begin{exmpl} \label{exmpl:Par-G-for-G-m}
Suppose that $G = \mathbb G_m$. This is a special case of \cref{exmpl:Par-G-for-GL-n}. An $L$-parameter for $\mathbb G_m$ is a continuous character $\phi\colon W_F \to \Qellbar^\times$. Every $\phi$ is supercuspidal and its connected component is of the form $\mathbb G_m \times */\mathbb G_m$, where the points in that connected component correspond to twists of $\phi$ by unramified characters of $W_F$.
\end{exmpl}

\begin{exmpl} \label{exmpl:Par-G-for-GL-2}
Suppose that $G = \GL_2$, which is again a special case of \cref{exmpl:Par-G-for-GL-n}. Given an $L$-parameter $\phi\colon W_F \to \GL_2(\Qellbar)$ with associated connected component $C$, there are three cases to consider:
\begin{enumerate}[(a)]
    \item $\phi$ is irreducible: Then $C \isom \mathbb G_m \times */\mathbb G_m$.
    \item $\phi = \chi_1 \oplus \chi_2$ for two characters $\chi_1$ and $\chi_2$ that are not an unramified twist of each other. Then $C \isom \mathbb G_m^2 \times */\mathbb G_m^2$.
    \item $\phi = \chi_1 \oplus \chi_2$ for two characters such that $\chi_1 \chi_2^{-1}$ is unramified, i.e. trivial on $I_F$.
\end{enumerate}
The cases (a) and (b) have been handled in \cref{exmpl:Par-G-for-GL-n} and in this case $C$ is always smooth. Case (c) is analyzed in detail in Bertoloni Meli's notes \cite{BMnotes}, and will also be discussed carefully in our forthcoming Seoul lecture notes. In fact, a hands-on analysis of this example was critical to the genesis of this entire paper, as we will explain there.
\end{exmpl}

\subsection{Ind-coherent sheaves}

In this subsection we study some basic properties of quasicoherent and ind-coherent sheaves on the stack of $L$-parameters. Throughout this section we fix a reductive group $G$ over $F$ and a finite quotient $Q$ of $W_F$ through which the $W_F$-action on $\hat G$ factors (as in \cref{set:reductive-group-over-F}). We furthermore assume without loss of generality that $G$ is quasisplit (see \cref{rmk:Par-G-only-depends-on-qsplit-form}) and we fix closed subgroups $T \subseteq B \subseteq G$ such that $T$ is a (not necessarily split) maximal torus and $B$ is a Borel in $G$.

From the results in \cref{sec:alggeo}, in particular the 6-functor formalism of ind-coherent sheaves from \cref{rslt:6ff-for-ICoh}, we can deduce the following fundamental properties.

\begin{thm}
\begin{thmenum}
    \item \label{rslt:QCoh-and-ICoh-on-ParG-compactly-generated} We have natural isomorphisms
    \begin{align*}
        \QCoh(\Par_G) = \Ind(\Perf^\qc(\Par_G)), \qquad \ICoh(\Par_G) = \Ind(\Coh^\qc(\Par_G)),
    \end{align*}
    where $\Perf^\qc \subseteq \Perf$ and $\Coh^\qc \subseteq \Coh$ denote the full subcategories of those sheaves that are supported on finitely many connected components.

    \item The functor $\gamma$ from \cref{rslt:gamma-from-QCoh-to-ICoh-on-stacks} induces a fully faithful symmetric monoidal embedding
    \begin{align*}
        \gamma\colon \QCoh(\Par_G) \injto \ICoh(\Par_G)
    \end{align*}
    that is obtained as the $\Ind$-extension of the inclusion $\Perf^\qc(\Par_G) \subseteq \Coh^\qc(\Par_G)$ via (i).
    
    \item \label{rslt:ParG-is-lci-with-trivial-dualizing-sheaf} $\Par_G$ is a local complete intersection with trivial dualizing complex and in particular cohomologically smooth. Moreover, it is $\ICoh$-fine and hence admits $!$-functors in the $\ICoh$-formalism.
\end{thmenum}
If $P \subseteq G$ is a parabolic subgroup containing $T$ then all of the above claims are also true for $\Par_P$ in place of $\Par_G$.
\end{thm}
\begin{proof}
Since $P = G$ is a special case for $P$, it is enough to prove all claims for $\Par_P$. We first prove (i). Note that both sides of the claimed isomorphisms commute with colimits along increasing unions of connected components of $\Par_P$ (for the left-hand sides this follows from descent, for the right-hand sides it follows from \cite[Corollary~A.5.9]{heyer-mann-6ff}). We can thus replace $\Par_P$ by one of its connected components, which by \cref{rslt:connected-components-of-Par-G-are-affines-mod-G} is of the form $X/\hat P$ for some affine scheme $X = \Spec A$ with an action of $\hat P$. Then the isomorphism for $\QCoh$ follows from \cite[Corollary~3.22]{BenZviFrancisNadler-IntegralTransforms} and the isomorphism for $\ICoh$ was proved in \cref{rslt:ICoh-equals-IndCoh-on-nice-stacks}.

We now prove (iii). By \cref{rslt:connected-components-of-Par-G-are-affines-mod-G} $\Par_P \to *$ is a local complete intersection morphism with trivial dualizing complex, and in particular cohomologically smooth by \cref{rslt:QCoh-solid-on-proper-and-smooth-maps}. Its connected components are affine schemes modulo $\hat P$ and hence QCA stacks, which by \cref{rslt:QCA-maps-are-ICoh-fine} implies that the connected components (and hence $\Par_P$) are $\ICoh$-fine.

It remains to prove (ii). By (iii) we know that $\Par_P$ is geometric and suave, hence \cref{rslt:gamma-fully-faithful-for-suave} implies that $\gamma\colon \QCoh(\Par_P) \injto \ICoh(\Par_P)$ is fully faithful. To compare $\gamma$ with the ind-extension of the functor in (i), we observe that by \cref{rmk:comparison-of-Coh-and-Perf-in-ICoh} the functors $\nu\colon \Coh^\qc(\Par_P) \to \ICoh(\Par_P)$ and $\gamma\colon \Perf^\qc(\Par_P) \to \ICoh(\Par_P)$ use the same functor $\eta$, but differ a priori as they use suave duality (i.e. $\intHom(-, \omega_{\Par_P})$) and natural duality respectively. However, since $\omega_{\Par_P} = \calO_{\Par_P}$ by (iii), $\nu$ and $\gamma$ are compatible on $\Par_P$ in the naive way.
\end{proof}

We now analyze ind-coherent sheaves on $\Par_G$ further. An important tool to study sheaves on $\Par_G$ is via the central grading, which we define next.

\begin{defn}
As in the proof of \cref{rslt:connected-component-of-supercuspidal-rep} there is an action of $\LocSys_{Z(\hat G)}$ on $\Par_G$ induced by the multiplication map $Z(\hat G) \times \hat G \to \hat G$. Restricting to the point of $\LocSys_{Z(\hat G)^{W_F}}$ given by the identity (i.e. the group homomorphism $W_F \to Z(\hat G) \rtimes Q$ with trivial image in $Z(\hat G)$), whose centralizer is $Z(\hat G)^{W_F}$, we obtain the \emph{central action}
\begin{align*}
    */Z(\hat G)^{W_F} \times \Par_G \to \Par_G.
\end{align*}
\end{defn}

Explicitly, the central action is trivial on objects, but on morphisms (in the groupoid $\Par_G(A)$ for a $\Lambda$-algebra $A$) it induces the multiplication of an element in the centralizer of an $L$-parameter by an element in $Z(\hat G)^{W_F}$.

The central action induces an action in the correspondence category of $\ICoh$-fine stacks (by \cite[Proposition~2.3.7(ii)]{heyer-mann-6ff}) and transferring it along $\ICoh$ we obtain an action of categories\footnote{Using \cite[Lemma~3.2.5]{heyer-mann-6ff} one can upgrade this to an action in the category of $\Qellbar$-linear presentable categories, but we do not need this.}
\begin{align*}
    *\colon \ICoh(*/Z(\hat G)^{W_F}) \times \ICoh(\Par_G) \to \ICoh(\Par_G).
\end{align*}
The symmetric monoidal structure on $\ICoh(*/Z(\hat G)^{W_F})$ is the one coming from convolution, using the abelian group structure of $*/Z(\hat G)^{W_F}$. We further observe by \cref{rslt:gamma-fully-faithful-for-suave} that
\begin{align*}
    \ICoh(*/Z(\hat G)^{W_F}) = \QCoh(*/Z(\hat G)^{W_F}).
\end{align*}
We warn the reader that this equivalence is not compatible with the convolution monoidal structures, only up to invertible twist (because $\gamma$ does not transfer the $\QCoh$ lower-$*$ to the $\ICoh$ lower-$!$ functor). One can instead use the equivalence induced by $\nu$ via \cref{rslt:ICoh-equals-IndCoh-on-nice-stacks}.

\begin{prop} \label{rslt:ParG-central-action-Cartier-dualities}
\begin{propenum}
    \item There is a natural symmetric monoidal equivalence
    \begin{align*}
        \ICoh(*/Z(\hat G)^{W_F}) = \prod_{\chi \in X^*(Z(\hat G)^{W_F})} \D(\Qellbar),
    \end{align*}
    where the symmetric monoidal structure on the left is convolution and on the right is componentwise tensor product. We denote by $\mathcal O(\chi) \in \ICoh(*/Z(\hat G)^{W_F})$ the sheaf corresponding to $1$ in the $\chi$-part and $0$ everywhere else.

    \item \label{rslt:central-grading-on-ParG} The central action induces a decomposition
    \begin{align*}
        \ICoh(\Par_G) = \prod_{\chi \in X^*(Z(\hat G)^{W_F})} \ICoh(\Par_G)^\chi,
    \end{align*}
    where $\ICoh(\Par_G)^\chi \subseteq \ICoh(\Par_G)$ is the full subcategory of those $\mathcal F$ where $\mathcal O(\chi) * \mathcal F = \mathcal F$.
\end{propenum}
\end{prop}
\begin{proof}
Part (i) is a form of Cartier duality and follows in a very similar way as \cref{rslt:Cartier-duality-for-torus}; we leave the details to the reader. 

Part (ii) follows formally from (i). Indeed, since $\mathcal O(\chi)$ is idempotent in $\ICoh(*/Z(\hat G)^{W_F})$, the functor $(-)^\chi := \mathcal O(\chi) * -$ sends $\ICoh(\Par_G)$ to $\ICoh(\Par_G)^\chi$. By combining these functors we obtain a functor from $\ICoh(\Par_G)$ to the product over $\ICoh(\Par_G)^\chi$. It remains to show that this functor is an equivalence, for which we consider its left adjoint $(\mathcal F^\chi)_\chi \mapsto \bigoplus_\chi \mathcal F^\chi$. We need to show that for $\mathcal F \in \ICoh(\Par_G)$ and $(\mathcal E^\chi)_\chi \in \prod_\chi \ICoh(\Par_G)^\chi$ the following maps are isomorphisms:
\begin{align*}
    \bigoplus_\chi \mathcal F^\chi \to \mathcal F, \qquad \mathcal E^\chi \to (\bigoplus \mathcal E^\chi)^\chi
\end{align*}
The first one follows from the fact that the convolution unit in $\ICoh(*/Z(\hat G)^{W_F})$ is the direct sum of $\mathcal O(\chi)$, and the second follows from $\mathcal O(\chi) * \mathcal O(\chi') = 0$ for $\chi \ne \chi'$; both of these observations are immediate from  (i).
\end{proof}

\begin{defn} \label{def:central-grading-on-ParG}
Given $M \in \ICoh(\Par_G)$, we denote by $M^\chi \in \ICoh(\Par_G)^\chi$ its \emph{$\chi$-graded summand} according to \cref{rslt:central-grading-on-ParG}.
\end{defn}

\begin{rmk}
As noted in the proof, part (i) of \cref{rslt:ParG-central-action-Cartier-dualities} is the Cartier duality of the stacks $*/Z(\hat G)^{W_F}$ and $X^*(Z(\hat G)^{W_F})$. In fact, part (ii) is also a form of Cartier duality, but one categorical level higher: The stacks $B^2(Z(\hat G)^{W_F})$ and $X^*(Z(\hat G)^{W_F})$ are 2-Cartier dual for $\ICoh$, meaning that
\begin{align*}
    \ICoh_2(B^2(Z(\hat G)^{W_F})) = \ICoh_2(X^*(Z(\hat G)^{W_F})) = \prod_\chi \ICoh_2(*),
\end{align*}
where $\ICoh_2$ is an appropriate 2-categorical version of $\ICoh$; it is a presentable 2-category defined via descent from the affine case and on an affine scheme $X$ it is roughly given by the 2-category of $\ICoh(X)$-linear presentable categories. By descent, $\ICoh_2(B^2(Z(\hat G)^{W_F}))$ consists of $\Qellbar$-linear presentable categories with an action by $\ICoh(*/Z(\hat G)^{W_F})$. Clearly $\ICoh(\Par_G)$ is an element of this 2-category, hence corresponds to an element in $\prod_\chi \ICoh(*)$, i.e. a family of $\Qellbar$-linear presentable categories, given by $\ICoh(\Par_G)^\chi$. We refer the reader to \cite[\S8]{GLC5} for more details on 2-Cartier duality. See also \cite[\S10]{Scholze-Gestalten} for a very general form of Cartier duality (and all its higher versions).
\end{rmk}

\subsection{Supercuspidal components}

We maintain the setup of the previous subsection, where we introduced ind-coherent sheaves on $\Par_G$ and studied their basic properties. In order to get a deeper understanding of ind-coherent sheaves, we now employ the following strategy: First understand ind-coherent sheaves on the supercuspidal components of $\Par_G$ and then attempt an inductive procedure to understand the non-supercuspidal locus via the supercuspidal components of Levi subgroups of $G$. In the present subsection we study the supercuspidal part and in the following subsections we discuss the induction process.


Let $C$ be a connected component of $\Par_G$ containing a supercuspidal $L$-parameter $\phi$, let $S_{\phi}$ be the centralizer of $\phi$ and let $C^\coarse$ be corresponding connected component in $X_G^\spec$. Recall from \cref{rslt:connected-component-of-supercuspidal-rep} that $C$ is an étale gerbe over the torus $C^\coarse$ banded by $S_\phi$ and that $C^\coarse = X^{\rm nr}(G)/K_\phi$ for the torus $X^{\rm nr}(G) = Y_\phi$ of unramified central 1-cocycles of $G$ (see \cref{def:torus-of-unramified-central-1-cocycles}) and a finite subgroup $K_\phi\subseteq X^{\rm nr}(G)$. In particular $C$ is smooth and hence $\ICoh(C) = \QCoh(C)$. We will need the following characterization of coherent sheaves on $C$ (recall the definition of the central grading from \cref{def:central-grading-on-ParG}).

\begin{prop} \label{rslt:coh-criterion-supercuspidal}
Let $C \subseteq \Par_G$ be a supercuspidal connected component and $\mathcal{F} \in \QCoh(C)$ any object. Then $\mathcal{F}$ is coherent (equivalently perfect) if and only if the following two conditions hold.
\begin{propenum}
    \item The $\chi$-graded summand $\mathcal{F}^\chi$ is nonzero for only finitely many $\chi \in X^{\ast}(Z(\hat{G})^{W_F})$.
    
    \item For all $\mathcal{G} \in \Perf(C)$, $\Gamma(C,\mathcal{G} \otimes \mathcal{F})$ is a bounded complex whose cohomologies are finitely generated $\mathcal{O}(C^{\coarse})$-modules.
\end{propenum}
\end{prop}
\begin{proof}
Suppose first that $\mathcal F$ is coherent, i.e. perfect (by $\ICoh(C) = \QCoh(C)$ and \cref{rslt:QCoh-and-ICoh-on-ParG-compactly-generated}). Then (i) is clear because $\mathcal F$ is compact and (ii) follows because $\mathcal G \tensor \mathcal F$ is still perfect and the pushforward along $C \to C^\coarse$ preserves perfect sheaves (using \cref{rslt:connected-component-of-supercuspidal-rep} and the fact that the same is true for the pushforward along $*/S_\phi \to *$).

Now assume that $\mathcal F$ satisfies (i) and (ii). Let $\pi\colon C' := X^{\mathrm{nr}}(G) \times BS_{\phi} \to C$ be the finite \'etale covering constructed in \cref{rslt:connected-component-of-supercuspidal-rep}. It suffices to show that $\pi^{\ast} \mathcal{F}$ is coherent. If $\mathcal{G}' \in \Perf(C')$ is given, then
\begin{align*}
    \Gamma(C',\mathcal{G}' \otimes \pi^{\ast} \mathcal{F}) = \Gamma(C,\pi_{\ast}(\mathcal{G}' \otimes \pi^{\ast}\mathcal{F})) = \Gamma(C,\pi_{\ast}\mathcal{G}' \otimes \mathcal{F})
\end{align*}
by the projection formula (note that $\pi$ is QCA and use \cref{rslt:QCA-maps-have-base-change-and-projection-formula}). Since $\pi_{\ast} \mathcal{G}'$ is perfect, condition (ii) implies that these cohomologies are bounded complexes of finitely generated $\mathcal{O}(C^{\coarse})$-modules. Now $\pi^{\ast}\mathcal{F}$ is an $\mathrm{Irr}(S_{\phi})$-graded object of $\QCoh( X^{\mathrm{nr}}(G) )$, so it suffices to show that each $\rho$-graded summand is coherent, for varying $\rho \in \mathrm{Irr}(S_{\phi})$. Condition (i) shows that these graded summands vanish for all but finitely many $\rho$, since the $Z(\hat{G})^{W_F}$-grading is a coarsening of the $S_{\phi}$-grading along the finite-index inclusion of groups $Z(\hat{G})^{W_F} \subset S_{\phi}$. Moreover, for a fixed $\rho$, the $\rho$-graded summand is given by $\Gamma(C', p^{\ast}\rho^{\vee} \otimes \pi^{\ast} \mathcal{F})$, where $p\colon X^{\mathrm{nr}}(G) \times BS_{\phi} \to BS_{\phi}$ is the obvious projection. This is a bounded complex of finitely generated $\mathcal{O}(C^{\coarse})$-modules by our observations above. Since $\mathcal{O}(C^{\coarse}) \to \mathcal{O}( X^{\mathrm{nr}}(G) )$ is a finite ring map, this gives the desired result.
\end{proof}

Later, we will also need the following criterion for identifying the structure sheaf of a supercuspidal component. Recall from the proof of \cref{rslt:connected-component-of-supercuspidal-rep} that $X^{\rm nr}(G)$ acts on $C$.

\begin{prop} \label{prop:strsheafcriterionsupercuspidal}
Let $C \subseteq \Par_G$ be a supercuspidal connected component and $\mathcal{F} \in \Coh(C)$ any coherent (equivalently perfect) sheaf. Then the following conditions are equivalent.
\begin{propenum}
    \item There is an isomorphism $\mathcal{F} \simeq \mathcal{O}_C$.
    \item For every $\nu \in X^{\mathrm{nr}}(G)$ with associated closed immersion $i_{\nu}: BS_{\nu \cdot \phi} = BS_{\phi} \to C$, there is an isomorphism $i_{\nu}^{!}\mathcal{F} \simeq \mathcal{O}_{BS_{\nu \cdot \phi}}[-\dim X^{\mathrm{nr}}(G)]$. Here $i_{\nu}^{!}$ denotes the right adjoint to the naive pushforward $i_{\nu,\ast}$ on $\QCoh$.
\end{propenum}
\end{prop}
\begin{proof}
It is clear that (i) implies (ii): By \cref{rslt:connected-component-of-supercuspidal-rep} the statement reduces to the similar claim about the closed immersion $i_\nu\colon * \to C^\coarse$; but $C^\coarse$ is a torus and hence the claim is an easy computation. For the converse, fix $\mathcal{F}$ satisfying (ii). By an easy manipulation with Grothendieck-Serre duality, the condition in (ii) is equivalent to asking that $i_{\nu}^{*}\Dgs \mathcal{F} \simeq \mathcal{O}_{BS_{\nu \cdot \phi}}$ for all $\nu$. This easily implies that $\Dgs \mathcal{F}$ is a line bundle, which moreover has trivial stabilizer action at all closed points. Since $\pi\colon C \to C^{\coarse}$ is a good moduli space, \cite[Theorem 10.3]{AlperGood} now implies that $\Dgs \mathcal{F}$ is the pullback of a line bundle on $C^{\coarse}$. But $C^{\coarse}$ is a torus, so it has no nontrivial line bundles. Therefore $\Dgs\mathcal{F} \simeq \mathcal{O}_C$, so also $\mathcal{F} \simeq \mathcal{O}_C$ using that $C$ has trivial dualizing complex.
\end{proof}

\subsection{Spectral Eisenstein series and constant term functors}

We continue with the setup from the previous subsection, i.e. \cref{set:reductive-group-over-F}; in particular $G$ is a quasisplit reductive group over $F$. We further fix the maximal torus $T \subseteq G$ that is used in the definition of the dual group. In the following we introduce the spectral version of the Eisenstein series and constant term functors on $\Par_G$, which allow inductive arguments for studying $\Par_G$ via $\Par_M$ for Levi subgroups $M$ of $G$.

\begin{defn} \label{def:spectral-Eis-and-CT}
Let $P = MU \subseteq G$ be a parabolic subgroup with Levi $M$ and assume that $P$ contains $T$. We obtain natural group homomorphisms $\hat M \from \hat P \to \hat G$ and hence a diagram of stacks
\begin{equation*}\begin{tikzcd}
    & \Par_P \arrow[dl,"f",swap] \arrow[dr,"g"]\\
    \Par_M && \Par_G
\end{tikzcd}\end{equation*}
We denote
\begin{align*}
    \Eis_P^\spec &:= g_! f^*\colon \ICoh(\Par_M) \to \ICoh(\Par_G),\\
    \CT_P^\spec &:= f_! g^*\colon \ICoh(\Par_G) \to \ICoh(\Par_M),
\end{align*}
and call them the \emph{spectral Eisenstein functor} and \emph{spectral constant term functor} respectively.
\end{defn}

\begin{rmk}
In terms of the classical notation for the pullback and pushforward functors on $\ICoh$, we have
\begin{align*}
    \Eis_P^\spec = g_*^{\ICoh} f^{!,\ICoh}, \qquad \CT_P^\spec = f_*^{\ICoh} g^{!,\ICoh}.
\end{align*}
\end{rmk}

In the following we will summarize the basic properties of $\Eis_P^\spec$ and $\CT_P^\spec$. They rely on a good understanding of the geometric properties of $f$ and $g$:

\begin{lem} \label{rslt:geometric-properties-of-Eis-and-CT-maps}
In the setting of \cref{def:spectral-Eis-and-CT}, we have:
\begin{lemenum}
    \item $f$ is a local complete intersection with trivial dualizing sheaf.
    \item $g$ is schematic and proper.
\end{lemenum}
\end{lem}
\begin{proof}
Both statements are claimed in \cite[Lemma~3.3.1]{Zhu}. For the first part of (i) one argues with the cotangent complex, similarly to the proof that $\Par_P$ is a local complete intersection; in fact the statement holds more generally for the map $\LocSys_{H_1} \to \LocSys_{H_2}$ for \emph{any} surjective map $H_1 \to H_2$ of algebraic groups over $\Qellbar$. The triviality of the dualizing sheaf follows from the triviality of the dualizing sheaf of $\Par_P$ and $\Par M$ (see \cref{rslt:ParG-is-lci-with-trivial-dualizing-sheaf}). Part (ii) is not proved in the reference, so we provide a quick argument: Factor $g$ as
\begin{align*}
    g\colon \Par_P = \Par^\square_P / \hat P \xto{g_1} \Par^\square_G/\hat P \xto{g_2} \Par^\square_G / \hat G = \Par_G.
\end{align*}
The map $g_2$ is proper and schematic because it is a base-change of the map $*/\hat P \to */\hat G$, which has fiber $\hat G/\hat P$, a proper scheme. The map $g_1$ is a closed immersion, e.g. because both source and target admit a closed immersion into $\hat G^n$ for large enough $n$ (see \cref{rslt:connected-components-of-Par-G-are-affines-mod-G}).
\end{proof}

\begin{cor}
In the setting of \cref{def:spectral-Eis-and-CT}, the following is true:
\begin{corenum}
    \item $\Eis_P^\spec$ is left adjoint to $\CT_P^\spec$.
    \item $\Eis_P^\spec$ preserves coherent and $\CT_P^\spec$ preserves admissible ind-coherent sheaves.
\end{corenum}
\end{cor}
\begin{proof}
By \cref{rslt:geometric-properties-of-Eis-and-CT-maps} we see that $f$ is $\ICoh$-prim with $\delta_f^{\ICoh} \isom 1$ and $g$ is $\ICoh$-suave with $\omega_g^{\ICoh} = 1$. Therefore we have $f_! \isom f_*$ and $g^! = g^*$ on $\ICoh$. This immediately implies (i). Part (ii) follows either by direct computation (using that coherent sheaves are the compact objects in $\ICoh$ by \cref{rslt:QCoh-and-ICoh-on-ParG-compactly-generated} and that admissible sheaves are characterized via \cref{rslt:admissible-sheaf-criterion-via-coherent}) or by using the identification of coherent sheaves with $\ICoh$-prim sheaves (see \cref{rslt:ICoh-equals-IndCoh-on-nice-stacks}) and admissible sheaves with $\ICoh$-suave sheaves (by definition) and applying the standard stability properties of suave and prim objects (see \cite[Lemma~4.5.16]{heyer-mann-6ff}).
\end{proof}


Typically, the functor $\CT_P^\spec$ does not preserve coherent sheaves, but as we will see in \cref{rslt:finiteness-for-constant-terms} this is true on graded pieces for the central grading. An important ingredient of that theorem is an inductive principle for understanding (ind-)coherent sheaves on $\Par_G$ via sheaves on the supercuspidal locus and Eisenstein functors, as hinted above. The following result due to Takaya provides a precise formulation of this idea:

\begin{prop}[Takaya] \label{rslt:takaya-generation}
The category $\Coh^\qc(\Par_G)$ is generated under finite colimits and retracts by sheaves of the form $\Eis_P^\spec(\mathcal{G})$ where $P=MU$ is any parabolic containing a fixed Borel of $G$ and $\mathcal{G} \in \Perf^{\mathrm{qc}}(\Par_M)$ is supported on the supercuspidal locus of $\Par_M$.
\end{prop}
\begin{proof}
This is a slight restatement of \cite[Proposition 3.4]{Takaya}, which is stated in terms of $\ICoh$ rather than $\Coh^\qc$.
\end{proof}

\subsection{The spectral geometric lemma}

In this section we develop our main technical tool on the spectral side. We continue with the setup from the previous subsections, so in particular $G$ is a quasisplit connected reductive group over $F$. We also fix a maximal torus and Borel $T \subseteq B \subseteq G$ and we denote the relative Weyl group of $G$ by $W$. Given a parabolic $P = MU \subseteq G$ containing $B$, we get an associated stack $\Par_P$ (see \cref{def:ParP}) and in \cref{def:spectral-Eis-and-CT} we have introduced the functors $\Eis_P^{\spec}$ and $\CT_P^{\spec}$ between $\ICoh(\Par_M)$ and $\ICoh(\Par_G)$. In this subsection we are interested in computing the composition of two such functors.

More precisely, suppose $P_1, P_2 \subseteq G$ are two parabolic subgroups containing the fixed Borel $B$, with Levi quotients $M_1$ and $M_2$. We want to compute the composition
\begin{align*}
    \CT_{P_2}^\spec \comp \Eis_{P_1}^\spec\colon \ICoh(\Par_{M_1}) \to \ICoh(\Par_{M_2}).
\end{align*}
In the spirit of the geometric lemma from representation theory, one might hope this composition is glued from functors which are closely related to compositions $\Eis^\spec_{Q_2} \comp w \comp \CT_{Q_1}^\spec$ for certain auxiliary parabolics $Q_i \subseteq M_i$ and Weyl elements $w$. We will show that something like this is true, but the statement is significantly more complicated than the usual geometric lemma. The essential issue here is that excision works very differently for ind-coherent sheaves than it does for constructible sheaves. Some of the analysis in this section was inspired by \cite[Sections 15.3 and 15.5]{GLC3}, but the arguments there are missing many details. In particular, no proof is provided in loc. cit. for the critical statements concerning the weight gradings of various sheaves. We will give a careful accounting of these issues, ultimately appealing to the contraction principle as developed in Section \ref{sec:contraction}.

By a quick base-change argument (cf. the proof of Theorem \ref{rslt:spectral-geometric-lemma} below), we see that the composition $\CT_{P_2}^\spec \comp \Eis_{P_1}^\spec$ is given by pull-push along the evident correspondence induced by the fiber product $\Par_{P_1} \times_{\Par_G} \Par_{P_2}$, hence we start by analyzing the geometry of that fiber product.

\begin{prop} \label{rslt:geometry-of-Par-of-double-coset-stack}
Let $P_1, P_2 \subseteq G$ be parabolic subgroups containing $B$ with (standard) Levi subgroups $M_1,M_2$ and Weyl groups $W_1$, $W_2$. Then the stack $X := \Par_{P_1} \times_{\Par_G} \Par_{P_2}$ admits a topological stratification
\begin{align*}
    \abs X = \bigsqcup_{w \in W_1 \backslash W / W_2} \abs{X_w}
\end{align*}
into locally closed subsets such that the stratum indexed by $w \in W_1 \backslash W / W_2$ is the image of a natural locally closed immersion $X_w \to X$ with the following properties:
\begin{propenum}
    \item Lift $w$ to an element in $W$ of minimal length and define the classical algebraic closed subgroups $H_w := P_1 \cap w^{-1}(P_2) \subset G$ and $\hat H_w := \hat P_1 \cap w^{-1}(\hat P_2) \subseteq \hat G$. Then
    \begin{align*}
        X_w = \Par_{H_w} := \LocSys_{\hat H_w}.
    \end{align*}

    \item \label{rslt:geometry-of-Par-of-double-coset-stack-diagram} For $i = 1, 2$, let $M_i$ be the Levi quotient of $P_i$, $U_i$ the unipotent radical and define
    \begin{align*}
        Q_1 &:= (P_1 \cap w^{-1}(P_2)) / (U_1 \cap w^{-1}(P_2)) = M_1 \cap w^{-1}(P_2),\\
        Q_2 &:= (w(P_1) \cap P_2) / (w(P_1) \cap U_2) = w(P_1) \cap M_2.
    \end{align*}
    Then $Q_i$ is a standard parabolic subgroup of $M_i$ (i.e. containing $B \cap M_i$) and we denote its Levi quotient by $L_i$. Finally we denote $K_w := (P_1 \cap w^{-1}(P_2)) / (U_1 \cap w^{-1}(U_2))$ and similarly for $\hat K_w$, with the associated stack
    \begin{align*}
        \Par_{K_w} := \LocSys_{\hat K_w}.
    \end{align*}
    There is a commuting diagram
    \begin{equation*}\begin{tikzcd}
        && \Par_{H_w} \isom^w \Par_{w(H_w)} \arrow[d,"h_w"]\\
        && \Par_{K_w} \isom^w \Par_{w(K_w)} \arrow[dl] \arrow[dr] \arrow[ddll,bend right,"g_w",swap] \arrow[ddrr,bend left,"f_w"]\\
        & \Par_{Q_1} \arrow[dl,"g_1"] \arrow[dr,"f_1"] && \Par_{Q_2} \arrow[dl,"f_2"] \arrow[dr,"g_2"]\\
        \Par_{M_1} && \Par_{L_1} \isom^w \Par_{L_2} && \Par_{M_2}
    \end{tikzcd}\end{equation*}
    where all maps are induced by the natural maps between the algebraic groups. In this diagram the center trapezoid is cartesian and the compositions $g_w h_w$ and $f_w h_w$ coincide with the natural maps $X_w \to X \to \Par_{M_i}$. Moreover, $h_w$ is lci with trivial dualizing sheaf.

    \item The center of $\hat K_w \rtimes Q$ is
    \begin{align*}
        Z_w := Z(\hat M_1)^{W_F} \cap w^{-1}(Z(\hat M_2)^{W_F}).
    \end{align*}
    It induces an action of $*/Z_w$ on $\Par_{K_w}$ and hence a natural grading
    \begin{align*}
        \ICoh(\Par_{K_w}) = \prod_{\chi' \in X^*(Z_w)} \ICoh(\Par_{K_w})^{\chi'}.
    \end{align*}
    This grading is compatible with the central gradings on $\Par_{M_1}$ and $\Par_{M_2}$ in the following sense. For $\chi \in X^*(Z(\hat M_1)^{W_F})$ we denote $\chi|_{Z_w} \in X^*(Z_w)$ its restriction to $Z_w$. Then for $\mathcal F \in \ICoh(\Par_{K_w})$, $\mathcal E \in \ICoh(\Par_{M_1})$, $\chi' \in X^*(Z_w)$ and $\chi \in X^*(Z(\hat M_1)^{W_F})$ we have
    \begin{align*}
        (g_w^* \mathcal E)^{\chi|_{Z_w}} = g_w^* \mathcal E^\chi, \qquad g_{w!} \mathcal F^{\chi'} = \bigoplus_{\chi: \; \chi|_{Z_w} = \chi'} (g_{w!} \mathcal F)^\chi.
    \end{align*}
    A similar statement holds for $\hat{M}_2$ and $f_w$, up to a twist by $w$.
\end{propenum}
\end{prop}
We will refer to the diagram appearing in (ii) above as \emph{the huge diagram}.
\begin{proof}
Denote by $W_F \surjto Q$ a finite quotient through which the $W_F$-action on $\hat G$ factors, as in \cref{set:reductive-group-over-F}. By \cref{rslt:Rep-W-F-colimit-of-mapping-stacks} we can fix some discretization $W_{F,\iota}$ of $W_F$ and some open normal subgroup $P \subseteq P_F$ contained in the kernel of $W_F \surjto Q$ such that if $\Gamma := W_{F,\iota}/P$ then for any of the algebraic groups $H$ appearing in the claim, we can replace $\Par_H$ by
\begin{align*}
    \intRep_{\hat H \rtimes Q}(\Gamma) \times_{\intRep_Q(\Gamma)} *.
\end{align*}
Note that $\Gamma$ is a classical discrete group, hence $\intRep_{\hat H \rtimes Q}(\Gamma) = \intMap(*/\Gamma, */(\hat H \rtimes Q))$ is an honest mapping stack. Note that $*/(\hat P_1 \rtimes Q) \times_{*/(\hat G \rtimes Q)} */(P_2 \rtimes Q) = (\hat P_1 \backslash \hat G / \hat P_2) / Q$, so with the above reduction the stack $X$ becomes
\begin{align*}
    \intMap(*/\Gamma, (\hat P_1 \backslash \hat G / \hat P_2) / Q) \times_{\intMap(*/\Gamma, */Q)} *.
\end{align*}

Now let $W^{\mathrm{abs}}$ resp. $\hat{W}^{\mathrm{abs}}$ be the absolute Weyl group of $G$ resp. $\hat{G}$. By the arguments in \cite[\S0.4.3]{KMSW}, there is a natural $Q$-equivariant isomorphism $W^{\mathrm{abs}} \cong \hat{W}^{\mathrm{abs}}$ inducing an isomorphism $W \cong \hat{W}$ on $Q$-invariants, and similarly for Weyl groups of standard Levis. By the Bruhat decomposition (see e.g. \cite[\S4.1.2]{Schieder-HNStrat}) there is a topological stratification of the stack $\hat P_1 \backslash \hat G / \hat P_2$ indexed by $\hat{W}^{\mathrm{abs}}_1 \backslash \hat{W}^{\mathrm{abs}} / \hat{W}^{\mathrm{abs}}_2$ such that the stratum indexed by any \[ w \in W_1 \backslash W / W_2 \cong  \hat{W}_1 \backslash \hat{W} / \hat{W}_2 \subset \hat{W}^{\mathrm{abs}}_1 \backslash \hat{W}^{\mathrm{abs}} / \hat{W}^{\mathrm{abs}}_2 \] is the image of a locally closed immersion $*/\hat H_w \to \hat P_1 \backslash \hat G / \hat P_2$. Note that the action of $Q$ preserves the individual strata indexed by $w \in W_1 \backslash W / W_2$, and $Q$ acts on $\hat{H}_w$ for such $w$, although it permutes the other strata. Now, passing to the quotient stack by $Q$ and noting that $\intMap(*/\Gamma, -)$ preserves locally closed immersions by \cref{rslt:mapping-stack-preserves-locally-closed-immersions},  we obtain a topological stratification indexed by $(\hat{W}^{\mathrm{abs}}_1 \backslash \hat{W}^{\mathrm{abs}} / \hat{W}^{\mathrm{abs}}_2)/Q$ of the subset $Z \subseteq \abs X$ consisting of those points where the induced map $*/\Gamma \to (\hat P_1 \backslash \hat G / \hat P_2)/Q$ factors over one of the strata of $(\hat P_1 \backslash \hat G / \hat P_2)/Q$. But clearly $Z = \abs X$ because $*/\Gamma$ has a single object $*$. Therefore we obtain a $(\hat{W}^{\mathrm{abs}}_1 \backslash \hat{W}^{\mathrm{abs}} / \hat{W}^{\mathrm{abs}}_2)/Q$-indexed topological stratification of $X$ such that the stratum indexed by any \[w \in W_1 \backslash W / W_2 = (\hat{W}^{\mathrm{abs}}_1 \backslash \hat{W}^{\mathrm{abs}} / \hat{W}^{\mathrm{abs}}_2)^Q \subset (\hat{W}^{\mathrm{abs}}_1 \backslash \hat{W}^{\mathrm{abs}} / \hat{W}^{\mathrm{abs}}_2)/Q\]
is the image of the locally closed immersion $X_w = \LocSys_{\hat{H}_w} \to X$ induced by the evident immersion $* / (\hat{H}_w \rtimes Q) \to (\hat P_1 \backslash \hat G / \hat P_2)/Q$ via passing to mapping stacks.

It therefore remains to see that \emph{only} the strata of $X$ indexed by $W_1 \backslash W / W_2$ are nonempty. For this, we will need a slightly finer analysis of the stack $(\hat P_1 \backslash \hat G / \hat P_2)/Q$. By standard Bruhat theory, the quotient $\hat P_1 \backslash \hat G / \hat P_2$ is stratified by the classifying stacks $\ast /( \hat{P}_1 \cap w^{-1}(\hat{P}_2))$ for $w \in \hat{W}^{\mathrm{abs}}_1 \backslash \hat{W}^{\mathrm{abs}} / \hat{W}^{\mathrm{abs}}_2$. A small further argument shows that the quotient $(\hat P_1 \backslash \hat G / \hat P_2)/Q$ is stratified by the classifying stacks $\ast / ((\hat{P}_1 \cap w^{-1}(\hat{P}_2)) \rtimes Q_w)$ where $w$ runs over (a set of representatives of) the $Q$-orbits in $\hat{W}^{\mathrm{abs}}_1 \backslash \hat{W}^{\mathrm{abs}} / \hat{W}^{\mathrm{abs}}_2$ and $Q_w < Q$ is the stabilizer of $w$. But given a point of $X$ mapping to the stratum of $\intMap(*/\Gamma, (\hat P_1 \backslash \hat G / \hat P_2) / Q)$ indexed by some $w$, the associated homomorphism $\Gamma \to (\hat{P}_1 \cap w^{-1}(\hat{P}_2)) \rtimes Q_w$ must have the property that its composition with the evident map $(\hat{P}_1 \cap w^{-1}(\hat{P}_2)) \rtimes Q_w \twoheadrightarrow Q_w \to Q$ is the canonical surjective map $\Gamma \to Q$, by the definition of $X$. This is only possible if $Q_w = Q$, or equivalently if $w$ lies in \[W_1 \backslash W / W_2 = (\hat{W}^{\mathrm{abs}}_1 \backslash \hat{W}^{\mathrm{abs}} / \hat{W}^{\mathrm{abs}}_2)^Q \subset (\hat{W}^{\mathrm{abs}}_1 \backslash \hat{W}^{\mathrm{abs}} / \hat{W}^{\mathrm{abs}}_2)/Q.\]   This completes the proof of (i).

Most of part (ii) immediately reduces to claims about the classifying stacks $*/\hat H$ in place of $\Par_H$ for the various groups $H$ appearing. It is then a simple exercise in algebraic groups, see \cite[\S4.1, \S4.2]{Schieder-HNStrat}. For the claim that $Q_i$ is a standard parabolic in $M_i$, see \cite[Lemma 2.11.(b)]{BZinduced}. The claim that $h_w$ is lci with trivial dualizing sheaf is proved as in \cref{rslt:geometric-properties-of-Eis-and-CT-maps}.

It remains to prove (iii). For the description of the center of $\hat K_w \rtimes Q$, we note that this group is the fiber product $(\hat Q_1 \rtimes Q) \times_{\hat L_1 \rtimes Q} (\hat Q_2 \rtimes Q)$ and the center of $\hat Q_i \rtimes Q$ is the same as the center of $\hat M_i \rtimes Q$, i.e. $Z(\hat M_i)^{W_F}$. To construct the grading on $\ICoh(\Par_{K_w})$, we proceed as in \cref{def:central-grading-on-ParG}: Multiplication yields an action of $Z_w$ on $\hat K_w$, which induces an action of $\LocSys_{Z_w}$ on $\Par_{K_w}$, which restricts to an action of $*/Z_w$ on $\Par_{K_w}$ and hence the desired grading as in \cref{rslt:central-grading-on-ParG}. The claim about the compatibility of the gradings with $f_w$ and $g_w$ follows easily from the fact that the maps $\hat K_w \to \hat M_i$ are compatible with the multiplication actions by $Z_w$ and $Z(\hat M_i)^{W_F}$.
\end{proof}

Before we can come to the statement of the spectral geometric lemma, we need a better understanding of how the central grading interacts with certain operations.

\begin{prop} \label{rslt:A1-retract-and-central-action}
Let $H$ be a linear algebraic group over $\Qellbar$ equipped with a $Q$-action and let $\mu\colon \mathbb G_m \to H$ be a cocharacter of $H$ such that $Q$ acts trivially on its image and the induced conjugation action of $\mathbb G_m$ on $H$ extends to an $\mathbb A^1$-action. Let $K$ be the associated retract and note that the $Q$-action restricts to $K$. We thus obtain the maps
\begin{align*}
    \LocSys_K \xto{i} \LocSys_H \xto{\pi} \LocSys_K,
\end{align*}
whose composition is the identity.
\begin{propenum}
    \item The $\mathbb A^1$-action on $H$ induces an $\mathbb A^1$-action on $\LocSys_H$ such that $(\LocSys_H, \LocSys_K)$ is a QCA $\mathbb A^1$-retract.
    
    \item Let $Z$ be the center of $K \rtimes Q$. Then there is a natural action of $*/Z$ on $\LocSys_K$ and hence a central grading
    \begin{align*}
        \QCoh(\LocSys_K) = \bigoplus_{\chi \in X^*(Z)} \QCoh(\LocSys_K)^\chi.
    \end{align*}
    Furthermore, $\mu$ factors over $Z$ and hence defines a cocharacter of $Z$.

    \item \label{rslt:A1-retract-weight-bounds} Let $\mathcal F \in \Perf(\LocSys_H)$. Then the characters of $Z$ appearing in $\pi_* \mathcal F$ (i.e. those $\chi$ such that $(\pi_* \mathcal F)^\chi \ne 0$) are those that can be written as a sum of a character appearing in $i^* \mathcal F$ and a character $\chi$ with $\chi = 0$ or $\langle \mu, \chi \rangle < 0$.

    \item \label{rslt:A1-retract-LocSys-perfect-pushforward} For all $\mathcal F \in \Perf(\LocSys_H)$ and all $\chi \in X^*(Z)$, the sheaf $(\pi_* \mathcal F)^\chi \in \QCoh(\LocSys_K)$ is perfect.
\end{propenum}
\end{prop}
\begin{proof}
We first prove (i). Note that the $\mathbb A^1$-action on $H$ induces a $\mathbb A^1$-action on $*/(H \rtimes Q)$ and thus, for every group $\Gamma$, an action of $\intMap(*/\Gamma, \mathbb A^1)$ on $\intMap(*/\Gamma, */(H \rtimes Q)) = \intRep_{H \rtimes Q}(\Gamma)$. This action is compatible with the projection to $\intRep_Q(\Gamma)$ and via the natural map of monoids $\mathbb A^1 \to \intMap(*/\Gamma, \mathbb A^1)$ (induced by the projection $*/\Gamma \to *$), so we obtain an action of $\mathbb A^1$ on $\intRep_{H \rtimes Q}(\Gamma) \times_{\intRep_Q(\Gamma)} *$. By \cref{rslt:Rep-W-F-colimit-of-mapping-stacks} these actions glue to an $\mathbb A^1$-action on $\LocSys_H$ with retract $\LocSys_K$. It follows easily from  \cref{rslt:basic-properties-of-rep-stack} that both $\pi$ and $i$ are QCA (after restricting to some quotient $\Gamma$ of a discretization of $W_F$, they are maps between QCA stacks). This finishes the proof of (i).

Part (ii) is done in the same way as in \cref{def:central-grading-on-ParG} and \cref{rslt:central-grading-on-ParG}.

We now prove (iii). In the same way as in (i), the conjugation action of $Z$ on $H$ and $K$ induces an action of $Z$ on $\LocSys_H$ and $\LocSys_K$. The group homomorphism $H \rtimes Z \to H$, $(h, z) \mapsto hz$ induces a commuting diagram
\begin{equation*}\begin{tikzcd}
    \LocSys_H / Z \arrow[r] \arrow[rr,"a'",bend left] \arrow[d,"\pi"] & \LocSys_{H \rtimes Z} \arrow[r] \arrow[d] & \LocSys_H \arrow[d,"\pi"] \\
    \LocSys_K / Z \arrow[r] \arrow[d] & \LocSys_{K \rtimes Z} \arrow[r] \arrow[d] & \LocSys_K\\
    */Z \arrow[r] & \LocSys_Z
\end{tikzcd}\end{equation*}
The right-hand square is cartesian because $\LocSys$ preserves fiber products. The bottom horizontal map is given by the trivial L-parameter, i.e. induced by the projection $W_F \to *$. The left-hand squares are cartesian, where $Z$ acts on $\LocSys_H$ and $\LocSys_K$ as described above---this follows from the construction of these actions. In particular we see that every sheaf on $\LocSys_H$ naturally upgrades to a $Z$-equivariant sheaf via pullback along $a'$.

Let us check how the central grading interacts with the natural $Z$-equivariant structure from above. Note that since $Z$ acts trivially on $\LocSys_K$, we have $\LocSys_K/Z = */Z \times \LocSys_K$. Under this equivalence, the composition of the middle horizontal maps above is exactly the action map $a\colon */Z \times \LocSys_K \to \LocSys_K$. Moreover, \cref{rslt:Cartier-duality-for-torus} induces a $X^*(Z)$-indexed grading on $\QCoh(\LocSys_K/Z)$. An easy computation using base-change in the above diagram shows that for every $\mathcal F \in \QCoh(\LocSys_H)$ and every $\chi \in X^*(Z)$ we have
\begin{align*}
    (\pi_* a'^* \mathcal F)_\chi = \pr_{2*}(\pi_* a'^* \mathcal F \tensor \calO(-\chi)) = \pr_{2*}(a^* \pi_* \mathcal F \tensor \calO(-\chi)) =  a^-_* (\pr_2^* \pi_* \mathcal F \tensor \calO(-\chi)) = (\pi_* \mathcal F)^\chi,
\end{align*}
where $a^-, \pr_2\colon */Z \times \LocSys_K \to \LocSys_K$ denote the \emph{inverse} of the central action and the projection to the second factor, respectively. Altogether we see that the $\chi$-graded piece of the $Z$-equivariant structure on $\pi_* a'^* \mathcal F$ is the $\chi$-graded piece of the central grading on $\pi_* \mathcal F$.

We can now finish the argument using our results on the contraction principle, so fix $\mathcal F \in \Perf(\LocSys_H)$. Restricting the above $Z$-equivariant structure to a $\mathbb G_m$-equivariant structure via $\mu$, we deduce from \cref{rslt:A1-retract-pushforward-in-terms-of-pullback} that $\pi_* \mathcal F$ admits a natural $\mathbb G_m$-equivariant structure whose graded pieces are built using countable limits out of the $\mathbb G_m$-graded pieces of the sheaves $i^* \mathcal F \tensor (i^* i_* 1)^{\tensor n}$ for $n \ge 0$. Since $\pi_* \mathcal F$ is the direct sum of its $\mathbb G_m$-graded pieces, the central characters appearing in $\pi_* \mathcal F$ are at most the ones appearing in $i^* \mathcal F \tensor (i^* i_* 1)^{\tensor n}$ and these themselves are sums of a central character in $i^* \mathcal F$ and arbitrarily many central characters in $i^* i_* 1$.

To finish the proof of (iii) it remains to see that if $\chi$ appears in $i^* i_* 1$ then either $\chi = 0$ or $\langle \mu, \chi \rangle < 0$. But note that by the above discussion, for every $n \in \Z$ and every sheaf $\mathcal F \in \QCoh(\LocSys_H)$, the $n$-th $\mathbb G_m$-graded piece of $\pi_* \mathcal F$ computes as
\begin{align}
    (\pi_* \mathcal F)_n = \bigoplus_{\langle \mu, \chi \rangle = n} (\pi_* \mathcal F)^\chi \label{eq:A1-central-grading-Gm-in-terms-of-Z}
\end{align}
Apply this to $\mathcal F = i_* (i^* i_* 1)$, so that $\pi_* \mathcal F = i^* i_* 1$. By \cref{rslt:A1-retract-pushforward-in-terms-of-pullback} we have $(i^* i_* 1)_0 = 1$ and $(i^* i_* 1)_n = 0$ for $n > 0$. This immediately proves the claim.

It remains to prove (iv), so fix $\mathcal F \in \Perf(\LocSys_H)$. By \cref{eq:A1-central-grading-Gm-in-terms-of-Z} it is enough to show that for all $n \in \Z$, the $n$-th $\mathbb G_m$-graded piece $(\pi_* \mathcal F)_n$ (for the natural $\mathbb G_m$-equivariant structure) is perfect. We can work with $\intRep_{(-)}(\Gamma)$ instead of $\LocSys_{(-)}$, where $\Gamma = W_{F,\iota} / P$ is a quotient of a discretization of $W_F$, as in \cref{rslt:Rep-W-F-colimit-of-mapping-stacks}. After choosing a finite set of generators of $\Gamma$ and considering the induced projection $F_r \to \Gamma$ from the free group anima on $r$-generators, we obtain a natural commuting diagram
\begin{equation*}\begin{tikzcd}
    \intRep_H(\Gamma) \arrow[r,"f'"] \arrow[d,"\pi"] & \intRep_H(F_r) = H^r/H \arrow[d,"\pi'"] \\
    \intRep_K(\Gamma) \arrow[r,"f"] & \intRep_K(F_r) = K^r/K
\end{tikzcd}\end{equation*}
Here the quotient $H^r/H$ is formed by the simultaneous conjugation action of $H$ on $H^r$, similarly for $K^r/K$. We now show that $(\pi_* \mathcal F)_n$ is coherent for all $n \in \Z$. Since the horizontal maps in the above square are closed immersions, it suffices to show that $f_* (\pi_* \mathcal F)_n$ is coherent. On the other hand, the square is compatible with the $\mathbb G_m$-action, hence $f_* (\pi_* \mathcal F)_n = (f_* \pi_* \mathcal F)_n = (\pi'_* (f'_* \mathcal F))_n$. Note that $\mathcal F' := f'_* \mathcal F$ is a coherent and hence perfect (by \cref{rslt:Coh-on-smooth-stack-equals-Perf}) $\mathbb G_m$-equivariant sheaf on $H^r/H$. By applying \cref{rslt:A1-retract-pushforward-preserves-perfect} to the $\mathbb A^1$-retract $(H^r/H, K^r/K)$ we are left to show that $(i'^* i'_* 1)_n$ is perfect for all $n$, where $i'\colon K^r/K \to H^r/H$ is the natural map.

We observe that $i'^* i'_* 1 = g_* 1$, where $g\colon (K^r/K) \times_{H^r/H} K^r/K \to K^r/K$ is the projection to the first factor. Let $U := \ker(H \surjto K)$ and let $g'\colon K^r \times_{H^r/U} K^r \to K^r$ be the base-change of $g$ along $K^r \to K^r/K$; this is still a $\mathbb G_m$-equivariant map and it is enough to show that $(g'_* 1)_n$ is perfect for all $n$. Since $K^r \to H^r/U$ factors over $H^r$, we obtain the following $\mathbb G_m$-equivariant factorization of $g'$:
\begin{align*}
    K^r \times_{H^r/U} K^r \xto{g'_1} K^r \times_{H^r/U} H^r = K^r \times U \xto{g'_2} K^r.
\end{align*}
Note that $g'_1$ is an lci closed immersion and thus $g'_{1*}$ preserves perfect sheaves. We are thus reduced to showing that $g'_{2*}$ sends perfect $\mathbb G_m$-equivariant sheaves to sheaves with perfect graded pieces. We apply \cref{rslt:A1-retract-pushforward-preserves-perfect} to the $\mathbb A^1$-retract $(K^r \times U, K^r)$, which reduces the problem to showing that $(i''^* i''_* 1)_n$ is perfect for all $n$, where $i''\colon K^r \to K^r \times U$ is the trivial section. But $i''$ is an lci closed immersion (by \cref{rslt:map-between-smooth-stacks-is-lci}), hence $i''_*$ preserves perfect sheaves and we conclude.

Above we have shown that $(\pi_* \mathcal F)_n$ is coherent for all $n$. It remains to show that it is perfect. After passing to some affine cover and using \cref{rslt:pseudo-coherent-plus-fin-tor-dim-implies-perfect} we are reduced to showing that $\pi_* \mathcal F$ has finite Tor dimension. Since $\pi$ is QCA and hence $\pi_*$ satisfies the projection formula (by \cref{rslt:QCA-maps-have-base-change-and-projection-formula}), this follows easily from perfectness of $\mathcal F$.
\end{proof}

With \cref{rslt:geometry-of-Par-of-double-coset-stack}, \cref{rslt:A1-retract-and-central-action} and our work on ind-coherent sheaves at hand, it is now easy to deduce the following version of the geometric lemma. It is based on a similar discussion in the context of the geometric Langlands correspondence, cf. \cite[\S15.3]{GLC3}. In the statement, we use the following standard notation: Given an integer $n \ge 1$ and an exact functor $F\colon \cat C_1 \to \cat C_2$ between stable categories, a \emph{filtration} on $F$ of length $n$ is a chain $0 = F_0 \to F_1 \to \dots \to F_n = F$ of functors $\cat C_1 \to \cat C_2$. The \emph{graded pieces} of this filtration are the functors $\cofib(F_{i-1} \to F_i)$ for $i = 1, \dots, n$. We similarly define $\N$-indexed filtrations and their graded pieces.

\begin{thm}[Spectral geometric lemma] \label{rslt:spectral-geometric-lemma}
Let $P_1, P_2 \subseteq G$ be parabolic subgroups containing $B$ with Levi quotients $M_1$, $M_2$ and Weyl groups $W_1$, $W_2$. Maintain the notation from \cref{rslt:geometry-of-Par-of-double-coset-stack}. Then the functor
\begin{align*}
    \CT_{P_2}^\spec \comp \Eis_{P_1}^\spec\colon \ICoh(\Par_{M_1}) \to \ICoh(\Par_{M_2})
\end{align*}
admits a filtration of length $\abs{W_1 \backslash W / W_2}$ whose graded pieces $F_w$, for $w \in W_1 \backslash W / W_2$, satisfy the following properties:
\begin{thmenum}
    \item \label{rslt:spectral-geometric-lemma-graded-pieces} For every $n \ge 0$ denote
    \begin{align*}
        \mathcal K_{w,n} := \gamma(h_{w*} \Sym^n \Nm_w) \in \ICoh(\Par_{K_w}),
    \end{align*}
    where $\Nm_w \in \Perf(\Par_{H_w})$ denotes the normal bundle of $\Par_{H_w} \to X$ (as defined in \cref{rslt:ICoh-excision-normal-bundle}). Then the functor $F_w$ admits a natural exhaustive $\N$-indexed filtration whose $n$-th graded piece $F_{w,n}\colon \ICoh(\Par_{M_1}) \to \ICoh(\Par_{M_2})$ is given by
    \begin{align*}
        F_{w,n} = (f_w \comp w)_! (\mathcal K_{w,n} \tensor g_w^*(-))[1].
    \end{align*}
    
    \item \label{rslt:spectral-geometric-lemma-weight-bounds} For each character $\chi \in X^*(Z_w)$ let $\mathcal K_{w,n}^\chi$ be the $\chi$-graded piece of $\mathcal K_{w,n}$ for the central grading on $\Par_{K_w}$. Then we get decompositions
    \begin{align*}
        \mathcal K_{w,n} = \bigoplus_\chi \mathcal K_{w,n}^\chi, \qquad F_{n,w} = \bigoplus_\chi F_{n,w}^\chi.
    \end{align*}
    Assume that $M_1 = Q_1$ or $M_2 = Q_2$.\footnote{This assumption may be unnecessary for the claimed bounds on $\chi$, see \cref{wtnecessity} below.} Then:
    \begin{enumerate}[(a)]
        \item $\mathcal K_{w,n}^\chi$ is perfect, i.e. it is the image of a perfect sheaf under $\gamma$.

        \item For fixed $w$ and $\chi$, we have $\mathcal K_{w,n}^\chi = 0$ for $n \gg 0$.
    \end{enumerate}
\end{thmenum}
\end{thm}

This theorem perhaps looks somewhat forbidding. In Example \ref{exmpl:SGLrankone} below, we unpack the statement in the simplest nontrivial case.

\begin{proof}
We first prove (i), which is a consequence of \cref{rslt:geometry-of-Par-of-double-coset-stack} and excision for $\ICoh$. More precisely, consider the diagram
\begin{equation*}\begin{tikzcd}
    && X \arrow[ddll,bend right,"g'",swap] \arrow[ddrr,bend left,"f'"] \arrow[dl] \arrow[dr]\\
    & \Par_{P_1} \arrow[dl] \arrow[dr] && \Par_{P_2} \arrow[dl] \arrow[dr]\\
    \Par_{M_1} && \Par_G && \Par_{M_2}
\end{tikzcd}\end{equation*}
Here $X$ is defined so that the trapezoid is a pullback diagram. By base-change we see that $F := \CT_{P_2}^\spec \comp \Eis_{P_1}^\spec = f'_! g'^*$. By \cref{rslt:geometry-of-Par-of-double-coset-stack} we know that $X$ admits a canonical topological stratification indexed by $W_1 \backslash W / W_2$. Pick a chain $\emptyset = V_0 \subseteq V_1 \subseteq \dots \subseteq V_N = X$ of open substacks such that for all $i = 1, \dots, N$ we have $\abs{V_i \setminus V_{i-1}} = \abs{X_{w_i}}$ for some $w_i \in W_1 \backslash W / W_2$ (in particular $N = \abs{W_1 \backslash W / W_2}$). For each $i$ let $\Fil^i F := f'_{i!} g_i'^*$, where $f'_i$ and $g'_i$ are the restrictions of $f'$ and $g'$ to $V_i$. Denoting $j_i\colon V_{i-1} \injto V_i$ the open immersion (which behaves like a proper map in the $\ICoh$-formalism), we see that there is a natural map
\begin{align*}
    \Fil^i F = f'_{i!} \id g'^*_i \to f'_{i!} j_{i!} j_i^* g'^*_i = f'_{i-1,!} g'^*_{i-1} = \Fil^{i-1} F,
\end{align*}
inducing a filtration on $F$. By \cref{rslt:ICoh-excision-triangles}, the $i$-th graded piece of this filtration, i.e. the cofiber of the above map, is given by $\hat f'_{i!} \hat g'^*_i[1]$, where $\hat f'_i$ and $\hat g'_i$ are the restriction of $f'$ and $g'$ to the completion $\hat X_{w_i} \subseteq X$ of $X_{w_i}$ inside $V_i$. From now on we drop the $i$ from the notation and denote this graded piece by $F_w = \hat f'_{w!} \hat g'^*_w[1]$ with the maps defined as in the following diagram:
\begin{equation*}\begin{tikzcd}
    & \hat X_w \subseteq X \arrow[dl,"\hat g'_w",swap] \arrow[dr,"\hat f'_w"]\\
    \Par_{M_1} && \Par_{M_2}
\end{tikzcd}\end{equation*}
We apply \cref{rslt:ICoh-excision-filtration} and \cref{rslt:ICoh-excision-normal-bundle} in order to obtain an exhaustive $\N$-indexed filtration on $F_w$. Let $f'_w$ and $g'_w$ denote the precomposition of $f'$ and $g'$ with the closed immersion $X_w \to X$. Using the projection formula we see that the $n$-th graded piece $F_{w,n}$ of the filtration on $F_w$ is given by
\begin{align*}
    &F_{w,n}[-1] = f'_{w!} (\gamma(\Sym^n \Nm_w) \tensor g'^*_w(-)) = f_{w!} h_{w!} (\gamma(\Sym^n \Nm_w) \tensor h_w^* g_w^*(-)) =\\&\qquad= f_{w!}(h_{w!} \gamma(\Sym^n \Nm_w) \tensor g_w^*(-)).
\end{align*}
By \cref{rslt:geometry-of-Par-of-double-coset-stack-diagram} $h_w$ is lci (hence suave by \cref{rslt:QCoh-solid-on-proper-and-smooth-maps}) with trivial dualizing sheaf, hence by \cref{rslt:ICoh-codualizing-sheaf-for-QCA-maps} we have $h_{w*} = h_{w!}$ on $\ICoh$. Together with \cref{rslt:gamma-commutes-with-pushforward-along-suave-QCA-maps} we obtain $\mathcal K_{w,n} = h_{w!} \gamma(\Sym^n \Nm_w)$. This finishes the proof of (i).

We now prove (ii), so fix $w \in W_1 \backslash W / W_2$ and $n \ge 0$. First observe that $\gamma$ is compatible with the central grading on $\ICoh(\Par_{K_w})$ and a similarly defined central grading on $\QCoh(\Par_{K_w})$. We can thus argue in $\QCoh(\Par_{K_w})$ and we identify $\mathcal K_{w,n}$ with its value as a quasicoherent sheaf, i.e. $\mathcal K_{w,n} = h_{w*} \Sym^n \Nm_w$.

We can without loss of generality assume $M_2 = Q_2$, as the case $M_1 = Q_1$ is symmetric. Then also $\hat M_2 = \hat Q_2 = \hat L_2$ and thus
\begin{align*}
    \hat K_w = \hat Q_1 = \hat M_1 \cap w^{-1}(\hat P_2), \qquad \hat H_w = \hat P_1 \cap w^{-1}(\hat P_2).
\end{align*}
Moreover, it is clear that $Z_w = Z(\hat{M}_1)^{W_F}$. Now pick a strictly $\hat P_1$-dominant cocharacter $\mu\colon \mathbb G_m \to Z(\hat M_1)$, i.e. such that $\langle \mu, \chi \rangle > 0$ for all $\chi \in X^*(Z(\hat M_1))$ lying in the cone spanned by the weights of $Z(\hat M_1)$ acting on $\Lie(\hat P_1)$. Then the induced conjugation action of $\mathbb G_m$ on $\hat P_1$ extends to an $\mathbb A^1$-action with retract $\hat M_1$. The same is true after taking the intersection with the closed subgroup $w^{-1}(\hat P_2)$, i.e. $\mu$ induces an $\mathbb A^1$-action on $\hat H_w$ with retract $\hat K_w$. Since $\hat{P}_{1}$ is $Q$-stable, by replacing $\mu$ with the product of its $Q$-conjugates if necessary we can further assume that $\mu$ factors over $Z_w = Z(\hat M_1)^Q \subseteq Z(\hat M_1)$. We are thus in the setting of \cref{rslt:A1-retract-and-central-action} with $H = \hat H_w$ and $K = \hat K_w$. From \cref{rslt:A1-retract-LocSys-perfect-pushforward} we immediately deduce that all $\mathcal K_{w,n}^\chi$ are perfect, proving (ii).(a).

It remains to prove (ii).(b). Denote $i\colon \Par_{K_w} \to \Par_{H_w}$ the natural map induced by the inclusion $\hat K_w \subseteq \hat H_w$. By \cref{rslt:A1-retract-weight-bounds} the non-trivial characters $\chi \in X^*(Z_w)$ appearing in $\mathcal K_{w,n}$ are those that can be written as a sum of a character appearing in $i^* \Sym^n \Nm_w$ and a character $\chi$ such that $\chi = 0$ or $\langle \mu, \chi \rangle < 0$. Let us first analyze the characters in $i^* \Nm_w$, for which we compute $\Nm_w$. By definition we have $\Nm_w = L_{X_w/X}^\vee[1]$. By dualizing \cref{rslt:fiber-sequence-for-cotangent-complex} we get a fiber sequence
\begin{align*}
    L_{X_w/\Qellbar}^\vee \to L_{X/\Qellbar}^\vee|_{X_w} \to \Nm_w
\end{align*}
in $\QCoh(X_w)$. Using \cref{rslt:Der-for-conde-rep-stack} we compute, for every $\Qellbar$-algebra $A$ and $A$-point $\rho\colon \Spec A \to X_w$,
\begin{align*}
    \rho^*(L_{X_w/\Qellbar}^\vee) = \Hom_A(A, \rho^*(L_{X_w/\Qellbar}^\vee)) = \Der_{\Qellbar}(\Par_{H_w}, A)_\rho = C^*(W_F, (\hat{\mathfrak h}_w)_\rho),
\end{align*}
where $\hat{\mathfrak h}_w := \Lie(\hat H_w)$. Using \cref{rslt:cotangent-complex-of-limit} we compute similarly
\begin{align*}
    \rho^*(L_{X/\Qellbar}^\vee) = C^*(W_F, (\hat{\mathfrak p}_1 \times_{\hat{\mathfrak g}} \hat{\mathfrak p}_2)_\rho),
\end{align*}
where $\hat{\mathfrak p}_1 = \Lie(\hat P_1)$ and so on. Altogether we deduce
\begin{align*}
    \rho^* \Nm_w = \cofib\big( C^*(W_F, (\hat{\mathfrak h}_w)_\rho) \to C^*(W_F, (\hat{\mathfrak p}_1 \times_{\hat{\mathfrak g}} \hat{\mathfrak p}_2)_\rho) = C^*(W_F, (\hat{\mathfrak g} / (\hat{\mathfrak p}_1 + w^{-1}(\hat{\mathfrak p}_2)))_\rho).
\end{align*}
Now $C^*(W_F, -)$ can be computed using a finite resolution (see \cref{rslt:rep-stack-for-Gamma-q}), hence the central characters appearing in $i^* \Nm_w$ are the weights of the adjoint representation of $Z_w$ on $\hat{\mathfrak g} / (\hat{\mathfrak p}_1 + w^{-1}(\hat{\mathfrak p}_2))$; let us denote this set by $D_w \subseteq X^*(Z_w)$. Then the characters appearing in $i^* \Sym^n \Nm_w$ are those of the form $\sum_i n_i \chi_i$ for $\chi_i \in D_w$, $n_i \ge 1$ and $\sum_i n_i = n$; we denote this set by $D_w^n \subseteq X^*(Z_w)$.

It remains to analyze the complimentary set of characters, namely those with  $\langle \mu, \chi \rangle < 0$. Note that by the proof of \cref{rslt:A1-retract-weight-bounds} the relevant characters $\chi$ are those that satisfy $\langle \mu, \chi \rangle < 0$ for \emph{every} strictly $\hat P_1$-dominant cocharacter $\mu$ into $Z_w$ (namely, they are the central characters appearing in $i^* i_* 1$, and this is independent of the chosen $\mu$). By an easy exercise, this implies that $\chi$ lies in the submonoid $C_w^+ \subseteq X^*(Z_w)$ generated by the negatives of the weights appearing in the adjoint action of $Z_w$ on $\hat{\mathfrak u}_1$.\footnote{In the analogous step in \cite[Sections 15.3 and 15.5]{GLC3}, it is claimed without proof that one only needs the negatives of the weights appearing in the adjoint action of $Z_w$ on $\hat{\mathfrak u}_1 \cap w^{-1}(\hat{\mathfrak u}_2)$.}

Putting things together, we see that every $\chi$ appearing in $\mathcal K_{w,n}$ is of the form $\chi_1 + \chi_2$ for $\chi_1 \in D_w^n$ and $\chi_2 \in C_w^+$. One checks easily that this implies (ii).(b), because for $n \gg 0$, the set $D_w^n + C_w^+$ is disjoint from any fixed finite subset of $X^*(Z_w)$.
\end{proof}

\begin{rmks}
\begin{rmksenum}
    \item In the proof of \cref{rslt:spectral-geometric-lemma}, for the (unique) maximal $w$ corresponding to the open stratum, it is understood that $F_{w}^{n} := 0$ for $n > 0$. It is also understood that $\Sym^0 \Nm_w = \calO_{\Par_{H_w}}$ for all $w$. Moreover, in the classical notation of ind-coherent sheaves we have
    \begin{align*}
        F_{w,n} = (f_w \comp w)_*^{\ICoh} (h_{w*}(\Sym^n \Nm_w) \tensor g_w^{!,\ICoh}(-))[1].
    \end{align*}
    Here $h_{w*}(\Sym^n \Nm_w) \tensor -$ denotes the natural $\QCoh$-action on $\ICoh$.
    
    \item \label{wtnecessity}The assumption $M_1 = Q_1$ or $M_2 = Q_2$ is probably required for part (a) of \cref{rslt:spectral-geometric-lemma-weight-bounds}, but may be unnecessary for the weight bounds in (b). Without this assumption, the main issue is that $H_w$ may not contract to $K_w$ anymore. For example, consider $\hat G = \GL_4$ with Weyl group $W = S_4$, let $\hat P_1 = \hat P_2 \subseteq \GL_4$ be the standard parabolic of type $(2, 2)$ and let $w = (23)$ be the shortest representative of the unique nontrivial coset in $W_1 \backslash W / W_2$. Then
    \begin{equation*}
        \hat H_w = \left\{ \begin{pmatrix}
            * & * & * & *\\
            0 & * & 0 & *\\
            0 & 0 & * & *\\
            0 & 0 & 0 & *
        \end{pmatrix} \right\} \subseteq \GL_4, \quad
        \hat U_1 \cap w^{-1}(\hat U_2) = \left\{ \begin{pmatrix}
            1 & 0 & 0 & *\\
            0 & 1 & 0 & 0\\
            0 & 0 & 1 & 0\\
            0 & 0 & 0 & 1
        \end{pmatrix} \right\} \isom \mathbb G_a.
    \end{equation*}
    The projection $\hat H_w \to \hat K_w \isom \hat H_w/\mathbb G_a$ does not seem to have a section. Note that in this example, the center of $\hat H_w$ agrees with the center of $\hat K_w$, so the conclusion of part (b) in \cref{rslt:spectral-geometric-lemma-weight-bounds} is still correct: the central characters appearing in $h_{w*} \mathcal F$ for \emph{any} $\mathcal F \in \QCoh(\Par_{H_w})$ are the same as those appearing in $\mathcal F$.

    \item \label{rslt:precise-weight-bounds-in-spectral-geometric-lemma} From the proof of \cref{rslt:spectral-geometric-lemma-weight-bounds} we see that for $w$ such that $M_2=Q_2$, the only $\chi$ such that $\mathcal K_{w,n}^\chi \ne 0$ are those in $C_w^+ + D_w^n$, where $C_w^+, D_w^n \subseteq X^*(Z_w) = X^*(Z(\hat{M}_1)^{W_F})$ are the following subsets:
    \begin{itemize}
        \item $C_w^+$ is the submonoid generated by the opposites of characters appearing in the adjoint action of $Z_w$ on $\Lie(\hat U_1)$.
        
        \item $D_w^1$ is the set of $\chi \ne 0$ occurring in the adjoint action of $Z_w$ on $\Lie(\hat G) / (\Lie(\hat P_1) + \Lie(w^{-1}(\hat P_2)))$. $D_w^n$ is the set of characters of the form $\sum_i n_i \chi_i$ for $\chi_i \in D_w^1$ and $n_i > 0$ with $\sum_i n_i = n$.
    \end{itemize}
    However, for our applications below the essential case is those $w$ such that $M_1=Q_1$, in which case it is more natural to consider the grading of $w_{\ast}\mathcal K_{w,n} \in \QCoh( \Par_{w(K_w)}) = \QCoh( \Par_{Q_2})$ for the action of $w(Z_w) = Z(\hat{M}_2)^{W_F}$. In this case, the same analysis as above shows that the only $\chi$ such that $(w_{\ast}\mathcal K_{w,n})^{\chi} \neq 0$ are those in $C_w^+ + D_w^n$, where $C_w^+, D_w^n \subseteq X^*(Z(\hat{M}_2)^{W_F})$ are the following subsets:
    \begin{itemize}
        \item $C_w^+$ is the submonoid generated by the opposites of characters appearing in the adjoint action of $Z_w$ on $\Lie(\hat U_2)$.
        
        \item $D_w^1$ is the set of $\chi \ne 0$ occurring in the adjoint action of $Z(\hat{M}_2)^{W_F}$ on $\Lie(\hat G) / (\Lie(w(\hat P_1)) + \Lie(\hat P_2))$. $D_w^n$ is the set of characters of the form $\sum_i n_i \chi_i$ for $\chi_i \in D_w^1$ and $n_i > 0$ with $\sum_i n_i = n$.
    \end{itemize}
\end{rmksenum}
\end{rmks}

\begin{exmpl}\label{exmpl:SGLrankone}
Let us discuss the simplest nontrivial example of \cref{rslt:spectral-geometric-lemma}. Suppose $G$ is split with semisimple rank one, and $P_1=P_2=B$ is the Borel. Then $W_1 \backslash W / W_2 = W = \{1,w_0 \}$. Here $\Par_{H_{1}} = \Par_{B}$ is closed and $\Par_{H_{w_0}} = \Par_{T}$ is open inside $X$. Now the filtration in \cref{rslt:spectral-geometric-lemma} reduces to a fiber sequence
\begin{align*}
    F_1 \to \CT_{B}^{\spec} \circ \Eis_{B}^{\spec} \to F_{w_0}.
\end{align*}
Moreover:
\begin{itemize}
    \item For $w=w_0$, every single group appearing in the huge diagram in \cref{rslt:geometry-of-Par-of-double-coset-stack-diagram} is literally just $T$, and $F_{w_0}=w_0(-)$ is the functor of twisting by the evident automorphism.

    \item For $w=1$, every group appearing in the huge diagram in \cref{rslt:geometry-of-Par-of-double-coset-stack-diagram} is $T$ \emph{except} $H_1 = B$. In this case $h_w$ is just the canonical map $q:\Par_{B} \to \Par_{T}$, and we get $F_{1}^{n}(-)=q_{\ast} (\mathrm{Sym}^{n} \mathrm{Nm}_1) \otimes (-)$ where $\mathrm{Nm}_1$ is the normal bundle of the relative diagonal $\Par_B \to \Par_B \times_{\Par_G} \Par_B$.
\end{itemize}
\end{exmpl}

\begin{rmk}
In the preceding example, the reader may find it odd that we've constructed a canonical map $\CT_{B}^{\spec} \comp \Eis_{B}^{\spec} \to w_0(-)$, since this looks rather close to a counit map but $\CT$ is \emph{right} adjoint to $\Eis$. In fact, this map is a reflection of something else, namely the existence of a canonical natural transformation $\Eis_{B}^{\spec} \to \Eis_{B\ast}^{\spec}$ where $\Eis_{B\ast}^{\spec}$ is the (non-colimit-preserving) right adjoint of $\CT_{B^-}^{\spec}$. This follows immediately from the involutivity of $w_0(-)$ and the isomorphism $w_0 \circ \CT_{B}^{\spec} \simeq \CT_{B^-}^{\spec}$. This transformation is not a freak of nature: it is a spectral incarnation of the canonical transformation $\Eis_{P!} \to \Eis_{P \ast}$, taking into account that $\Eis_{P\ast}$ is right adjoint to $\CT_{P!} \cong \CT_{P^- \ast}$.

In fact, we've also secretly constructed the unit $\mathrm{id} \to \CT_{B}^{\spec} \circ \Eis_{B}^{\spec}$: it is just the composite map $\mathrm{id} \to F_{1}^{0} \to F_{1} \to \CT_{B}^{\spec} \circ \Eis_{B}^{\spec}$, where the first arrow is induced by the canonical map $\mathcal{O}_{\Par_T} \to q_{\ast} (\mathcal{O}_{\Par_B})$ and our description of $F_{1}^{n}$ above.
\end{rmk}

\subsection{Finiteness properties of spectral constant terms}

We keep the notation from the previous subsection and in particular fix a Borel $B \subseteq G$. Our goal in this section is to prove important finiteness properties for constant term functors and deduce finiteness properties for the $\Hom$ between coherent sheaves on the stack of L-parameters. The main work towards these results has been done in previous subsections, most importantly the Spectral Geometric Lemma in \cref{rslt:spectral-geometric-lemma}.

Let us start with the finiteness properties for constant terms. In the following, for any standard parabolic $P \subseteq G$ (i.e. such that $B \subseteq P$) with Levi quotient $M$ and any central character $\chi \in X^*(Z(\hat M)^{W_F})$, we denote by $\CT_P^{\spec,\chi}$ the composition
\begin{align*}
    \ICoh(\Par_G) \xto{\CT_P^{\spec}} \ICoh(\Par_M) \xto{(-)^\chi} \ICoh(\Par_M)^\chi,
\end{align*}
where the second functor comes from the central grading in \cref{def:central-grading-on-ParG}. The functor $\CT_P^{\spec,\chi}$ enjoys excellent properties:

\begin{thm}\label{rslt:finiteness-for-constant-terms}
Let $P \subseteq G$ be a standard parabolic subgroup with Levi quotient $M$. Then:
\begin{thmenum}
    \item For every $\chi \in X^*(Z(\hat M)^{W_F})$ the functor $\CT_P^{\spec,\chi}$ preserves compact objects, i.e. restricts to a functor
    \begin{align*}
        \CT_P^{\spec,\chi}\colon \Coh^\qc(\Par_G) \to \Coh^\qc(\Par_M).
    \end{align*}

    \item \label{rslt:finiteness-for-CT-vanishing} Given $\mathcal{F} \in \Coh^\qc(\Par_G)$ we have $\CT_{P}^{\spec,\chi}(\mathcal{F}) = 0$ for all but finitely many $P$-dominant $\chi$.
\end{thmenum}
\end{thm}

\begin{rmk}
We caution the reader not to misread \cref{rslt:finiteness-for-CT-vanishing}: for $\chi$ which are \emph{not} $P$-dominant, we are not claiming any vanishing properties at all. In fact, it is easy to find examples where $\CT_{P}^{\spec,\chi} (\mathcal{F}) \neq 0$ for \emph{every} $\chi$ which is not $P$-dominant.
\end{rmk}

\begin{proof}[Proof of \cref{rslt:finiteness-for-constant-terms}]
Fix $P$ and $\chi$ as in the statement of the theorem. We want to prove that $\mathcal{F}\mapsto \CT_{P}^{\spec,\chi}(\mathcal{F})$ has the claimed finiteness and vanishing properties. By Takaya's generation theorem (see \cref{rslt:takaya-generation}), we can assume $\mathcal{F}$ is of the form $\Eis_{P'}(\mathcal{G})$ for some standard parabolic $P'=M'U'$ and some perfect complex $\mathcal{G}$ with quasicompact support in the supercuspidal locus. Now applying the spectral geometric lemma with $P_1 := P'$ and $P_2 := P$, we get a finite filtration on $\CT_{P}^{\spec}(\Eis_{P'}^{\spec}(\mathcal{G}))$ with graded pieces $F_w(\mathcal{G})$. It thus suffices to prove the same coherence and vanishing statements for the $\chi$-graded summands $F_w^\chi(\mathcal{G})$.

Now observe that by our assumptions on $\mathcal{G}$, $F_w(\mathcal{G})$ vanishes identically for any $w$ such that $Q_1$ is a proper parabolic in $M_1=M'$, because for any such $w$ the image of the (proper) map $g_1$ in the huge diagram of \cref{rslt:geometry-of-Par-of-double-coset-stack-diagram} is disjoint from the supercuspidal locus in $\Par_{M'}$. Thus the only elements $w$ which contribute nontrivially to the spectral geometric lemma are those with $M_1 = Q_1 = L_1$. Let us call such $w$ \emph{relevant}. For any relevant $w$, it is clear that $w(K_w) = Q_2$ and $w(Z_w) = Z(\hat{M}_2)^{W_F}$. Moreover, for any relevant $w$ the huge diagram in \cref{rslt:geometry-of-Par-of-double-coset-stack-diagram} simplifies to the diagram
\[\xymatrix{
    & \mathrm{Par}_{H_{w}}\simeq^{w}\mathrm{Par}_{w(H_{w})}\ar[d]^{h_{w}}\\
    & \mathrm{Par}_{K_{w}}\simeq^{w}\mathrm{Par}_{Q_{2}}\ar[dl]_{g_w}\ar[dr]^{f_w}\\
    \mathrm{Par}_{M_{1}} &  & \mathrm{Par}_{M_{2}}
}\]
where moreover $g_w$ is quasismooth and $f_w$ is proper and schematic. By \cref{rslt:spectral-geometric-lemma-graded-pieces} the functor $F_w$ admits an exhaustive $\N$-filtration whose graded pieces $F_{w,n}$ satisfy
\begin{align*}
    F_{w,n}(\mathcal{G})[-1] = (f_w \comp w)_!(\mathcal{K}_{w,n} \tensor  g_w^*\mathcal{G})
\end{align*}
which we can simplify further by noting that $\mathcal H_w := g_w^*\mathcal{G}$ is actually a perfect complex (here we implicitly use that $\mathcal G$ lies in the image of $\gamma$ and $\gamma$ is compatible with pullback).

Now let $S_w$ be the finite set of characters occurring in the $w(Z_w) = Z(\hat{M}_2)^{W_F}$-grading of $w_{\ast}\mathcal{H}_w$. Putting everything together, we get that $F_w^\chi(\mathcal{G})$ admits a filtration with graded pieces
\begin{align*}
    F_{w,n}^{\chi}(\mathcal{G})[-1] 
    & = f_{w!}\big(\bigoplus_{\xi \in S_w} (w_{\ast}\mathcal{K}_{w,n})^{\chi-\xi} \tensor (w_{\ast}\mathcal{H}_{w})^{\xi}\big).
\end{align*} 
For any fixed $n$, \cref{rslt:spectral-geometric-lemma-weight-bounds} implies that $\mathcal K_{w,n}^\chi$ is perfect, hence $\bigoplus_{\xi \in S_w} (w_{\ast}\mathcal{K}_{w,n})^{\chi-\xi} \tensor (w_{\ast}\mathcal{H}_{w})^{\xi}$ is perfect and in particular coherent. Since $f_w$ is proper, the functor $f_{w!}$ preserves coherent sheaves ($f$ is $\ICoh$-suave by \cref{rslt:ICoh-on-lafp-schemes}, hence $f_{w!}$ preserves $\ICoh$-prim sheaves, which by \cref{rslt:ICoh-equals-IndCoh-on-nice-stacks} are exactly the coherent sheaves). This shows that each $F_{w,n}^{\chi}(\mathcal{G})$ is coherent. But now \cref{rslt:spectral-geometric-lemma-weight-bounds} shows that $F_{w,n}^{\chi}(\mathcal{G})$ vanishes for all but finitely many $n$. Therefore $F_{w}^{\chi}(\mathcal{G})$ is coherent.

Finally, from the second half of \cref{rslt:precise-weight-bounds-in-spectral-geometric-lemma} we immediately see that if $F_{w}^{\chi}(\mathcal{G}) \neq 0$, then $\chi$ must satisfy
\[ \chi \in S_w + C_{w}^+ + \cup_{n} D_{w}^{n}.\]
But it is easy to see that only finitely many $P$-dominant characters have this property. 
\end{proof}

Using \cref{rslt:finiteness-for-constant-terms} we can now prove the following miraculous finiteness property for coherent sheaves on the stack of $L$-parameters. In this result, recall that $X_G^\spec$ denotes the coarse moduli space of $L$-parameters for $G$ and $\mathfrak Z_G^\spec := \calO(X_G^\spec)$ its global sections (see \cref{def:coarse-moduli-space}).

\begin{thm} \label{rslt:miracle-bound}
For any $\mathcal{F},\mathcal{G} \in \Coh^\qc(\Par_G)$, $\RHom(\mathcal{F},\mathcal{G}) \in \D(\Qellbar)$ is a bounded complex whose cohomology groups are finitely generated $\mathfrak{Z}_{G}^\spec$-modules.
\end{thm}

This theorem is the exact spectral analogue of a famous theorem due to Bernstein: if $A,B \in \D(G(F),\Qellbar)^\omega$ are any two compact objects, then $\RHom(A,B)$ is a bounded complex whose cohomology groups are finitely generated $\mathfrak{Z}_G$-modules. In fact the connection is deeper than mere analogy, as Bernstein's theorem together with categorical local Langlands implies that \cref{rslt:miracle-bound} \emph{must} be true.

It is not hard to see that $\RHom(\mathcal{F},\mathcal{G})$ is left-bounded with individual cohomologies finitely generated over $\mathfrak{Z}_{G}^{\mathrm{spec}}$, but the right-boundedness is severely non-obvious because $\Par_G$ is not smooth. Indeed, such a right-boundedness property is completely false for any singular variety. However, the stackiness of $\Par_G$ gives enough room for it to be true.

\begin{proof}
We first show that $\RHom(\mathcal{F},\mathcal{G})$ is left-bounded with cohomologies finitely generated over $\mathfrak{Z}_{G}^{\mathrm{spec}}$. To see this, write
\[
    \RHom(\mathcal{F},\mathcal{G}) = \Gamma(X_{G}^{\mathrm{spec}}, q_{\ast} \intHom(\mathcal{F},\mathcal{G})),
\]
where $q\colon \Par_G \to X_G^\spec$ is the projection to the coarse moduli space from \cref{def:coarse-moduli-space} and all notation is in terms of quasi-coherent sheaves. Note that $\intHom(\mathcal{F},\mathcal{G})$ is left-bounded with coherent cohomology sheaves. Since $q:\Par_G \to X_{G}^{\mathrm{spec}}$ is a good moduli space, $q_\ast$ is exact by definition and preserves coherence by \cite[Theorem~4.16(x)]{AlperGood}, so $q_* \intHom(\mathcal{F},\mathcal{G})$ is left-bounded with coherent cohomology sheaves and quasicompact support. Since $X_G^\spec$ is a disjoint union of affine varieties with global sections $\mathfrak{Z}_{G}^\spec$, this gives the claim.

It remains to see that $\RHom(\mathcal{F},\mathcal{G})$ is right-bounded. For this, we will argue by induction on the semisimple rank of $G$, the base case of $G$ a torus being clear. When $\mathcal{F}$ is perfect, the argument in the previous paragraph easily shows that $\RHom(\mathcal{F},\mathcal{G})$ is bounded, using the resulting boundedness of $\intHom(\mathcal{F},\mathcal{G})$ together with the exactness of $q_\ast$ and the affineness of $X_{G}^{\spec}$. By \cref{rslt:takaya-generation}, it thus suffices to handle the complementary case where $\mathcal{F}=\Eis_{P}^{\spec}(\mathcal{F}')$ for some proper parabolic $P=MU$ and some $\mathcal{F}' \in \Coh^\qc(\Par_M)$. Without loss of generality we can assume that $\mathcal{F}'$ is $\chi$-graded for some $\chi \in X^*(Z(\hat M)^{W_F})$, so
\[
    \RHom(\Eis_{P}^\spec(\mathcal{F}'),\mathcal{G}) \cong \RHom(\mathcal{F}',\CT_{P}^{\spec,\chi}(\mathcal{G})).
\]
By \cref{rslt:finiteness-for-constant-terms}, $\CT_{P}^{\spec,\chi}$ carries $\Coh^\qc(\Par_G)$ into $\Coh^\qc(\Par_M)$. By induction on the semisimple rank this implies that $\RHom(\mathcal{F}',\CT_{P}^{\spec,\chi}(\mathcal{G}))$ is bounded, as desired.
\end{proof}

\subsection{Finite and restricted sheaves} \label{ss:finrestrcoh}

Maintain the setup from previous subsections and let \[q\colon \Par_G \to X_{G}^{\spec}\] be the map to the coarse moduli from \cref{def:coarse-moduli-space}. In this subsection we define an analogue of $\Par_G$ with \emph{restricted variation}, following the ideas of \cite{AGKRRV}. In fact, we will not really need this stack explicitly, but only the categories of ind-coherent and quasicoherent sheaves on it, which it turns out can be constructed directly inside the analogous categories for $\Par_G$. We begin with the following definition.

\begin{defn}
A coherent sheaf $\mathcal{F}\in \Coh(\Par_G)$ is \emph{finite} if it is set-theoretically supported on finitely many fibers of the map $q$ over closed points of $X_{G}^{\spec}$. We write $\Coh(\Par_G)_{\fin} \subset \Coh(\Par_G)$ for the resulting full subcategory.

We also write $\Perf(\Par_G)_{\fin} = \Perf(\Par_G) \cap \Coh(\Par_G)_{\fin}$ for the analogous category of finite perfect complexes.
\end{defn}

It is clear that $\Coh(\Par_G)_{\fin}$ is stable under finite colimits and retracts and under (twisted) Grothendieck-Serre duality, and similarly for $\Perf(\Par_G)_{\fin}$.

\begin{rmk}The name is meant to suggest that finite sheaves have finite length in some sense. As we will see later, there is a matching notion on the automorphic side where this intuition is even stronger. However, we caution the reader that objects of $\Coh(\Par_G)_{\fin}$ don't have finite length in any literal sense, and closed fibers of the map $q$ can have infinitely many $\Qellbar$-points.\footnote{We thank Teruhisa Koshikawa for showing us an example of such a fiber, with $G = \GL_{21}$!}
\end{rmk}

As a consequence of our cohomological bounds in the previous section, we get the following classification of finite coherent sheaves in terms of the a priori unrelated notion of admissibility (see \cref{def:admissibleIcoh}).

\begin{thm} \label{thm:cohfinclassification}
There is an identification
\[
    \Coh(\Par_G)_{\fin} = \Adm(\Par_G) \cap \Coh^\qc(\Par_G)
\]
as full subcategories of $\IndCoh(\Par_G)$. In particular, any finite coherent sheaf is admissible. Moreover, for any $\mathcal{F} \in \Coh(\Par_G)_{\fin}$, $\Dadm \mathcal{F}$ lies in $\QCoh(\Par_G) \subseteq \IndCoh(\Par_G)$ and satisfies $\mathbf{D}_{\adm}^2 \mathcal{F} \cong \mathcal{F}$.
\end{thm}
\begin{proof}
The second part follows from the first part together with Theorem \ref{rslt:spectral-temperization}.

To prove the first part, note that it is easy to see that $\Adm(\Par_G) \cap \Coh^\qc(\Par_G) \subset \Coh(\Par_G)_{\mathrm{fin}}$. Indeed, if $\mathcal{F} \in \Adm \cap \Coh^\qc(\Par_G)$, then $\End(\mathcal{F})$ is an Artinian $\Qellbar$-algebra (note that $\RHom(\mathcal F, \mathcal F)$ is a perfect $\Qellbar$-algebra by \cref{rslt:admissible-sheaf-criterion-via-coherent}), so the map $\mathfrak{Z}_{G}^{\spec} \to \End(\mathcal{F})$ factors over an Artinian quotient of $\mathfrak{Z}_{G}^{\spec}$, and thus $\mathcal{F}$ is supported set-theoretically on the preimage under $q$ of the associated finite-length subscheme $Z \subset X_{G}^{\spec}$.

Thus we just need to prove that
\[
    \Coh(\Par_G)_{\fin} \subset \Adm(\Par_G).
\]
Fix any $\mathcal{F} \in \Coh(\Par_G)_{\fin}$. By the criterion stated in \cref{rslt:admissible-sheaf-criterion-via-coherent}, we need to show that for any $\mathcal{G} \in \Coh^\qc(\Par_G)$, the complex $\RHom(\mathcal{G},\mathcal{F})$ is bounded with finite-dimensional cohomology groups. By \cref{rslt:miracle-bound}, it is a bounded complex whose cohomologies are finitely generated $\mathfrak{Z}_{G}^{\spec}$-modules. Since also $\mathcal{F} \in \Coh(\Par_G)_{\fin}$ by assumption, the map $\mathfrak{Z}_{G}^{\spec} \to \End(\mathcal{F})$ factors over some Artinian quotient $\mathfrak{Z}_{G}^{\spec}/I$. But then each cohomology of $\RHom(\mathcal{G},\mathcal{F})$ is a finitely generated $\mathfrak{Z}_{G}^{\spec}/I$-module, and hence is finite-dimensional as a vector space.
\end{proof}

We now introduce a mechanism to pass between arbitrary sheaves and finite sheaves, at least after ind-completion. 

\begin{defn}
Given a closed point $\phi \in X_{G}^{\spec}(\Qellbar)$, let $\mathcal{O}^{\restr}_{\phi} \in \QCoh(\Par_G)$ be the local cohomology sheaf of $\mathcal{O}_{\Par_G}$ relative to the closed substack $q^{-1}(\phi)$ (see \cref{def:local-cohomology-sheaf}). Let
\[
    \mathcal{O}^{\restr} \in \QCoh(\Par_G)
\]
be the direct sum of $\mathcal{O}^{\restr}_{\phi}$ over all closed points $\phi$.
\end{defn}

In the notation of \cref{rslt:excision-for-QCoh}, the local cohomology functor $\hat i_\natural \hat i^*$ is oplax monoidal and hence preserves coalgebras. Since also $\hat i^* \hat i_\natural = \id$, it even preserves idempotent coalgebras. Applying this to the evident closed immersion $i: q^{-1}(\phi) \to \Par_G$, we deduce that each $\mathcal{O}^{\restr}_{\phi}$ is an idempotent coalgebra object of $\QCoh(\Par_G)$. Similarly, it is easy to see that
\[
    \calO^{\restr}_{\phi} \tensor \calO^{\restr}_{\phi'} = 0
\]
for distinct $\phi \neq \phi'$, which immediately implies that $\calO^{\restr}$ is an idempotent coalgebra object. 

\begin{rmk} \label{rslt:explicit-description-of-O-restr}
By \cref{rslt:computation-of-local-cohomology} each $\calO^{\restr}_{\phi}$ can be presented very explicitly: choosing elements $f_1,\dots,f_n \in \mathfrak{Z}_{G}^{\spec}$ generating the associated maximal ideal $\mathfrak{m}_{\phi}$, the $\phi$-summand $\calO^{\restr}_{\phi}$ identifies with the (alternating) Cech complex
\[
    \calO_{\Par_G} \to \bigoplus_i \calO_{\Par_G}[\tfrac{1}{f_i}] \to \bigoplus_{i<j} \calO_{\Par_G}[\tfrac{1}{f_i f_j}] \to \cdots.
\]
For later use, we also note that this complex is naturally the direct limit of the system of Koszul complexes $K(\calO_{\Par_G};f_{1}^j,\dots,f_{n}^{j})$, each of which lies in $\Perf(\Par_G)_{\fin}$ (see \cite[Tag 0913]{stacks-project}). Putting these observations together, we see that $\calO^{\restr}$ is naturally an object of $\Ind(\Perf(\Par_G)_{\fin}) \subseteq \QCoh(\Par_G)$.
\end{rmk}

\begin{defn}\label{def:Crestr}
If $\mathcal{C}$ is any presentable $\QCoh(\Par_G)$-linear category, we define $\mathcal{C}^{\restr} \subset \mathcal{C}$ as the full subcategory spanned by objects $A$ such that the canonical map $\calO^{\restr} \tensor A \to A$ is an isomorphism.
\end{defn}

Note that by definition, $\calO^{\restr}$ has a canonical direct sum decomposition indexed by semisimple $L$-parameters $\phi$, so any object of $\mathcal{C}^{\restr}$ inherits such a decomposition. We also observe for later use that if $F:\mathcal{C}\to \mathcal{D}$ is any colimit-preserving $\QCoh(\Par_G)$-linear functor of $\QCoh(\Par_G)$-module categories, it automatically restricts to a functor $\mathcal{C}^{\restr} \to \mathcal{D}^{\restr}$. Note that $\mathcal{C}^{\restr} \subset \mathcal{C}$ is stable under all colimits, so this inclusion has a right adjoint $\mathcal{C} \to \mathcal{C}^{\restr}$. Under mild conditions, this right adjoint has a very simple description.

\begin{prop} \label{prop:restrnaiveadjoint}
If $\mathcal{C}^{\restr}$ is compactly generated and the inclusion $\mathcal{C}^{\restr} \subset \mathcal{C}$ preserves compact objects, the right adjoint $R: \mathcal{C}\to \mathcal{C}^{\restr}$ to the inclusion is the functor $A\mapsto \mathcal{O}^{\restr} \otimes A$.
\end{prop}
\begin{proof}First note that the idempotent coalgebra structure on $\mathcal{O}^{\restr}$ implies that for any $A\in \mathcal{C}$, $\mathcal{O}^{\restr} \otimes A$ is naturally an object of $\mathcal{C}^{\restr}$. Now, our assumptions guarantee that $R$ preserves colimits and is $\QCoh(\Par_G)$-linear. Then for any $A \in \mathcal{C}$, we compute that
\[R(A) \overset{\sim}{\leftarrow} \mathcal{O}^{\restr} \otimes R(A) = R(\mathcal{O}^{\restr} \otimes A) = \mathcal{O}^{\restr} \otimes A,\]
where in the middle equality we used the $\QCoh$-linearity of $R$ and in the final equality we used the previous observation that $\mathcal{O}^{\restr} \otimes A$ already lies in $\mathcal{C}^{\restr}$.    
\end{proof}

Of course, in practice the non-formal task is to describe the category $\mathcal{C}^{\restr}$ explicitly. The following result takes the first steps in this direction. Note that part (i) treats the ``universal'' case where $\mathcal{C} = \QCoh(\Par_G)$ itself.

\begin{prop}\label{prop:restrsheavesParG}
Maintain the notation above.
\begin{propenum}
\item We have an identification
\[\QCoh(\Par_G)^{\restr} = \Ind(\Perf(\Par_G)_{\fin})\]
as full subcategories of $\QCoh(\Par_G)$.
\item We have an identification
\[\IndCoh(\Par_G)^{\restr} = \Ind(\Coh(\Par_G)_{\fin})\]
as full subcategories of $\IndCoh(\Par_G)$.
\end{propenum}
In both cases, the functor $\mathcal{C} \to \mathcal{C}^{\restr}$ is given by $A\mapsto \mathcal{O}^{\restr} \otimes A$.
\end{prop}
\begin{proof}

The proofs of the two statements are parallel, so we only prove the (slightly harder) second statement. First note that the tensor action of $\Perf$ on $\Coh$ clearly restricts to a functor $\Perf(\Par_G)_{\fin} \otimes \Coh^\qc(\Par_G) \to \Coh(\Par_G)_{\fin}$. Since $\mathcal{O}^{\restr} \in \Ind(\Perf(\Par_G)_{\fin})$, say with $\mathcal{O}^{\restr} = \varinjlim_i C_i$ and all $C_i$ finite perfect, then for any $A \in \IndCoh(\Par_G)^{\restr}$ we can write $A = \varinjlim_j A_j$ with $A_j \in \Coh^\qc(\Par_G)$, in which case
\begin{align*}
    A & = \mathcal{O}^{\restr} \otimes A = \varinjlim_{i,j} C_i \otimes A_j,
\end{align*}
where the first equality follows from the assumption that $A$ is restricted. Then $C_i \otimes A_j \in \Coh(\Par_G)_{\fin}$ for all $i,j$, showing that $\IndCoh(\Par_G)^{\restr}$ is contained in $\Ind(\Coh(\Par_G)_{\fin})$. 

For the reverse inclusion, it suffices to see that any $A \in \Coh(\Par_G)_{\fin}$ is restricted. Without loss, we can assume $A$ is supported on a single fiber of $q$ corresponding to a semisimple $L$-parameter $\phi$. Let $f_1,\dots,f_n$ generate the associated maximal ideal $\mathfrak{m}_{\phi} \subset \mathfrak{Z}_{G}^{\spec}$, so the $\phi$-summand of $\mathcal{O}^{\restr}$ is represented by the (alternating) Cech complex
\[ C:= \mathcal{O}_{\Par_G} \to \bigoplus_i \mathcal{O}_{\Par_G}[\tfrac{1}{f_i}] \to \bigoplus_{i<j}\mathcal{O}_{\Par_G}[\tfrac{1}{f_i f_j}] \to \cdots \]
as noted above. Then $\mathcal{O}^{\restr} \otimes A = C \otimes A$. Note that the canonical map $C \to \mathcal{O}_{\Par_G}$ has cone $K$ given by a complex whose terms are finite direct sums of localizations $\mathcal{O}_{\Par_G}[\tfrac{1}{f_{i_1} \cdots f_{i_m}}]$. In particular, this cone has a finite filtration such that some $f_i$ acts invertibly on each graded piece. On the other hand, each $f_i$ acts nilpotently on $A$. Therefore $K \otimes A = 0$, so the canonical map $\mathcal{O}^{\restr} \otimes A = C \otimes A \to A$ is an isomorphism, showing that $A$ is restricted.

The final claim follows from the explicit descriptions of the restricted categories together with \cref{prop:restrnaiveadjoint}.
\end{proof}

\begin{rmk}In most cases of interest, there is an identification
\[\mathcal{C}^{\restr} = \QCoh(\Par_G)^{\restr} \otimes_{\QCoh(\Par_G)} \mathcal{C},\]
but this identification does not actually seem to give any additional insight into the construction. As such, we decided to stick with the more naive approach above.
\end{rmk}

Finally, we note that these constructions can be interpreted directly in terms of sheaves on a stack of parameters with restricted variation. Here is the precise definition.

\begin{defn}
The stack $\Par_{G}^{\restr}$ is the subfunctor of $\Par_G$ whose $A$-points consist exactly of those $L$-parameters $\phi:W_F \to \phantom{}^L G(A)$ such that for any algebraic representation $V\in \mathrm{Rep}(\phantom{}^L G)$ and any element $v \in V\otimes_{\Qellbar} A$, the $W_F$-orbit $(V\circ \phi)(W_F)\cdot v$ spans a finite-dimensional $\Qellbar$-vector space.
\end{defn}

\begin{rmk}
We expect that the functor $\Par_G^\restr$ is a \enquote{formal algebraic stack} with the following description. As a subfunctor of $\Par_G$, it is canonically identified with the disjoint union of the formal completions of $\Par_G$ along all fibers $q^{-1}(x)$ for $x \in X_{G}^{\spec}(\Qellbar)$. The induced map to $X_G^{\spec}$ naturally makes $\Par_G^{\restr}$ into a relative algebraic stack over the formal scheme $X_G^{\spec,\restr}$ given as the disjoint union of the formal completions of $X_G^{\spec}$ at all closed points.

Moreover, pushforward along the natural morphism $\Par_{G}^{\restr} \to \Par_G$ induces natural equivalences $\QCoh(\Par_{G}^{\restr}) \cong \QCoh(\Par_{G})^{\restr}$ and $\ICoh(\Par_{G}^{\restr}) \cong \ICoh(\Par_{G})^{\restr}$.

As we do not need these statements, we do not attempt to prove them here.
\end{rmk}

\section{Automorphic results} \label{sec:automorphic-results}

We now turn to the automorphic side of the picture. We fix a finite extension $F$ of $\Q_p$ and a reductive group $G$ over $F$. Our main concern is the category $\D_{\lis}(\Bun_G,\Lambda)$ defined in \cite{FS}, with its action of Hecke operators and the associated spectral action. As defined originally by Fargues-Scholze, this category is embedded in a much larger category of solid sheaves which admit a certain five-functor formalism. One of our main goals in this chapter is to illustrate how to work with $\D_{\lis}(\Bun_G,\Lambda)$ using the motivic sheaves introduced by Scholze \cite{ScholzeBerkovich}. More precisely, for a fixed $\Z_\ell$-algebra $\Lambda$, we make systematic use of the sheaf theory $\D_{\rel}(-,\Lambda)$ on small v-stacks over $\Spd \overline{\F_p}$ defined in Section \ref{sec:def-of-D-rel} below. Roughly speaking, this sheaf theory takes motivic sheaves on $X/\Spd \overline{\F_p}$ and kills the ``noise'' coming from motivic sheaves on the base geometric point by applying a suitable $\ell$-adic realization functor on the base point. The notation ``rel'' is meant to indicate motivic sheaves realized \emph{relatively} along the base geometric point.  This sheaf theory has already appeared (with slightly different notation) in the papers \cite{HHS} and \cite{LHcourbes}, and will be systematically studied in the forthcoming paper \cite{GH}. 

The first key result here (which is one of the main theorems of \cite{GH}, and which was already used implicitly in \cite{HHS}) is that on $\Bun_G$, $\D_{\rel}$ agrees canonically with $\D_{\lis}$; see Theorem \ref{rslt:D-rel-on-BunG} for a more precise statement. Taking this theorem as our starting point, we explain how to reconstruct the Hecke action and spectral action entirely using the $\D_{\rel}$-formalism (Theorem \ref{rslt:spectral-action-for-D-rel}). From this point of view, the Hecke action is very natural, and does not require the dual embedding of the Satake category into solid sheaves or the exotic $\natural$-pushforward on solid sheaves used in \cite{FS}. Perhaps surprisingly, we show how to canonically realize the Satake category with $\Z_\ell$-coefficients defined in \cite{FS} inside $\D_{\rel}(-,\Z_\ell)$ evaluated on the local Hecke stacks. In particular, despite the essentially motivic nature of $\D_{\rel}$, we make no use of the motivic Satake equivalence proved in \cite{FS-motivic}.\footnote{In fact, we do not directly use any results from \cite[Section 3-6]{FS-motivic}, although our results and arguments at several points are closely inspired by that paper.} We also give a close analysis of $\D_{\rel}(\Div^1,\Lambda)$, and in particular show that this category embeds fully faithfully into a concrete category of Weil group representations.

Using this reconstruction of the spectral action, we prove a very general duality theorem for the spectral action (Theorem \ref{thm:spectraldualityuseful}). This duality theorem is a critical ingredient in our analysis of the categorical local Langlands conjecture.

\subsection{Background on motivic sheaves} \label{sec:def-of-D-rel}

In this section we give some background and general technical results on the six functor formalism $\D_{\rel}$ needed in our arguments. In the following we freely use the notation from \cite{etale-cohomology-of-diamonds,FS} and work with the category of small v-stacks over $\Spd \overline{\F}_p$. A few of the following definitions and results also appear in \cite[\S B]{LHcourbes}. The 6-functor formalism $\D_\rel$ is a variant of the presentable $\Z[1/p]$-linear 6-functor formalism $\D_\mot^\oc$ constructed by Scholze in \cite[\S2]{FS-motivic}.

Let us briefly recall how $\D_\mot^\oc$ is defined. For an affinoid perfectoid space $S = \Spa(R, R^+)$ over $\overline{\F}_p$, $\D_\mot^\oc(S)$ is obtained from the category of $\Z[1/p]$-linear ball-invariant finitary arc-sheaves over $S$ (see \cite[Definitions~4.10,~5.2]{ScholzeBerkovich}) by inverting the Tate twist (see \cite[Definition~9.1]{ScholzeBerkovich}). Here we implicitly choose a norm $\abs{\cdot}$ on $R$ and work with Banach rings over $(R, \abs{\cdot})$. However, we observe that $\D_\mot(S)$ does not depend on the chosen norm.

\begin{lem} Let $R$ be a uniform Tate-Huber ring equipped with a Banach algebra norm $|-|_R$ defining its topology. The category of uniform Banach rings over $R$ (for the fixed norm on $R$) is equivalent to the category of uniform Tate-Huber rings over $R$. In other words, for any map $R \to R'$ of uniform Huber-Tate rings there is a unique norm (up to equivalence) on $R'$ that makes this map into map of Banach rings.
\end{lem}
\begin{proof}By \cite[Definition~2.8.1]{relative-p-adic-hodge-1} every norm on $R$ or $R'$ is equivalent to a unique power-multiplicative norm. In particular we can replace $|-|_R$ with a power-multiplicative Banach algebra norm. We then need to prove that there is a unique power-multiplicative Banach algebra norm $|-|_{R'}$ on $R'$ defining its topology and making the given map $R \to R'$ into a map of Banach rings. 

We first prove uniqueness. More generally, we claim any two $R$-Banach algebra norms $|-|_1, |-|_2$ on $R'$ defining the given topology are equivalent. (If both are power-multiplicative, they are then equal by the usual trick of extracting $n$th roots for all $n\geq 1$.) For this, note that also $|-|'=\max(|-|_1,|-|_2)$ is an $R$-Banach algebra norm defining the given topology on $R'$, so $(R',|-|') \to (R',|-|_i)$ is a bounded (i.e. continuous) bijection of $R$-Banach modules, hence it is open and strict by the open mapping theorem \cite[Theorem 2.2.8]{relative-p-adic-hodge-1}. This shows that both norms $|-|_1, |-|_2$ are equivalent to $|-|'$, hence they are equivalent to each other.

For existence, let $y \in \Spa R'$ be any rank one point, with image $x \in \Spa R$. Let $K_x$ resp. $L_y$ be the completed residue fields of these points, so we have a map of affinoid fields $K_x \to L_y$. Via the homeomorphism $|\Spa R|^h \cong \mathcal{M}(R)$, the point $x$ is associated with a unique norm $|-|_x$ on $K_x$. Let $|-|_y$ be the unique multiplicative norm on $L_y$ which extends $|-|_x$ and defines the topology on $L_y$ (the existence and uniqueness of this extension is a standard exercise). Then by the uniformity of $R'$, the power-multiplicative seminorm $|-|_{R'} \overset{\mathrm{def}}{=} \sup_{y} |-|_y$ is in fact a norm on $R'$ making $R \to R'$ into a map of Banach rings.
\end{proof}
We thus obtain a well-defined functor
\begin{align*}
    (\Perfd_{\overline{\F}_p}^\aff)^\op \to \PrL, \qquad S \mapsto \D_\mot(S)
\end{align*}
on the category of affinoid perfectoid spaces over $\overline{\F}_p$. By construction this functor satisfies descent for all arc-covers, i.e.\ maps $\Spa(R', R'^+) \to \Spa(R, R^+)$ such that the induced map $\Spa(R', \Z) \to \Spa(R, \Z)$ is a v-cover. We thus obtain a sheaf of categories
\begin{align*}
    \vStk^\op \to \PrL, \qquad X \mapsto \D_\mot^\oc(X)
\end{align*}
by sending a v-stack $X$ to its arc-sheafification and then to the descended value of $\D_\mot$ from affinoid perfectoid spaces. In \cite{ScholzeBerkovich} (see also \cite[\S2]{FS-motivic}) a full 6-functor formalism is constructed on $\D_\mot^\oc$.

By construction we have $\D_\mot^\oc(\Spa(R, R^+)) = \D_\mot(\mathcal M_{\mathrm{arc}}(R), \Z[1/p])$ for every perfectoid Huber pair $(R, R^+)$ over $\overline{\F}_p$, where the right-hand side is the notation from \cite{ScholzeBerkovich}. However, one has to be careful with the meaning of $\D_\mot^\oc(*)$, as the final object in v-stacks does not agree with the final object in the (normed) arc-stacks of \cite{ScholzeBerkovich}. To remedy the situation, consider the v-stack $\mathcal N$ that assigns to every affinoid perfectoid $(R, R^+)$ the set of norms on $R$. Then:

\begin{lem} \label{rslt:D-mot-of-point-vs-N}
There are natural symmetric monoidal equivalences
\begin{align*}
    \D_\mot^\oc(*) = \D_\mot^\oc(\mathcal N) = \D_\mot(\overline{\F}_p),
\end{align*}
where the first equivalence comes from the pullback along $\mathcal N \to *$ and the last category is the one considered in \cite[\S11]{ScholzeBerkovich}.
\end{lem}
\begin{proof}
The second equivalence is tautological, because arc-stacks over $\mathcal N$ are the same as the normed arc-stacks over $\overline{\F}_p$ considered in \cite{ScholzeBerkovich}. To get the first equivalence, we first show that the pullback map $\D_\mot^\oc(*) \to \D_\mot^\oc(\mathcal N)$ is fully faithful. Since norms are unique up to scaling, $f\colon \mathcal N \to *$ is a torsor under the group $\R_{>0}$. But $\R_{>0}$ is $\D_\mot^\oc$-smooth with trivial cohomology (by comparison with the 6-functor formalism of sheaves on condensed anima, cf. the proof of \cite[Proposition~2.2]{FS-motivic}), which implies the same for the map $\mathcal N \to *$ by testing after base-change along itself. Thus for all sheaves $\mathcal F \in \D_\mot^\oc(*)$ we have $f_\natural f^* \mathcal F = f_\natural \mathbf 1 \tensor \mathcal F = \mathcal F$, proving full faithfulness of $f^*$.

It remains to show that $f^*$ is essentially surjective. Since it is fully faithful and preserves colimits, it is enough to show that the image of $f^*$ contains generators of $\D_\mot^\oc(\overline{\F}_p)$. By \cite[Theorem~11.1]{ScholzeBerkovich} a set of generators of $\D_\mot^\oc(\overline{\F}_p)$ can be constructed out of the homologies $\Z[X]$ of algebraic varieties over $\overline{\F}_p$. But these make sense over $*$ instead of $\mathcal N$, as they are purely algebraic and do not require norms to be constructed. Hence the relevant homologies can be realized in $\D_\mot^\oc(*)$ by $\natural$-pushing along the maps $X \to *$.
\end{proof}

With \cref{rslt:D-mot-of-point-vs-N} we have great control over the base-category $\D_\mot(*)$. In particular we make the following observation:

\begin{lem} \label{rslt:ell-adic-comparison-for-D-mot}
There is a natural symmetric monoidal and colimit-preserving functor $\D_\mot^\oc(*) \to \D(\Z_\ell)$.
\end{lem}
\begin{proof}
By \cite[Proposition~12.3]{ScholzeBerkovich} there are natural symmetric monoidal isomorphisms
\begin{align*}
    \Mod_{\Z/\ell^n}(\D^\oc_\mot(*)) = \D(\mathbf{Z}/\ell^n)
\end{align*}
for all $n \ge 1$. By passing to the limit of the base-change functors $- \tensor \Z/\ell^n$ we obtain a symmetric monoidal functor
\begin{align*}
    \D^\oc_\mot(*) \to \varprojlim_n \Mod_{\Z/\ell^n}(\D^\oc_\mot(*)) = \varprojlim_n \D(\mathbf{Z}/\ell^n).
\end{align*}
The right-hand side identifies with the full subcategory of $\D(\Z)$ spanned by the $\ell$-complete modules. The full subcategory of dualizable objects agrees with the category $\Perf(\Z_\ell)$ of perfect modules over $\Z_\ell$. Now by \cref{rslt:D-mot-of-point-vs-N} and \cite[Remark~11.2]{ScholzeBerkovich}, $\D^\oc_\mot(*)$ agrees with the usual category of Voevodsky motives over $\overline{\F}_p$ and in particular it is compactly generated by dualizable objects (see \cite[Theorem~11.1]{ScholzeBerkovich}). Thus the restriction of the above functor to dualizable objects induces a symmetric monoidal functor $\D^\oc_\mot(*)^\omega \to \Perf(\Z_\ell)$ and by passing to $\Ind$-categories we arrive at the desired symmetric monoidal functor $\D_\mot^\oc(*) \to \D(\Z_\ell)$.
\end{proof}

\begin{rmk}
We warn the reader that the functor in \cref{rslt:ell-adic-comparison-for-D-mot} does not agree with $\ell$-completion on non-compact objects.
\end{rmk}

Using the symmetric monoidal functor $\D_\mot^\oc(*) \to \D(\Z_\ell)$ from \cref{rslt:ell-adic-comparison-for-D-mot} we now base-change $\D_\mot^\oc$ in order to arrive at the desired 6-functor formalism $\D_\rel$.

\begin{defn} \label{def:D-rel}
Fix a $\Z_\ell$-algebra $\Lambda$ and consider the symmetric monoidal functor
\begin{align*}
    r\colon \D_\mot^\oc(*) \to \D(\Z_\ell) \xto{- \tensor_{\Z_\ell} \Lambda} \D(\Lambda),
\end{align*}
where the first functor is the one from \cref{rslt:ell-adic-comparison-for-D-mot}. Recall from \cite[Lemma~3.2.5]{heyer-mann-6ff} that $\D_\mot^\oc$ automatically factors over $\Mod_{\D_\mot^\oc(*)}(\PrL)$, hence composing $\D_\mot^\oc$ with the base-change along $r$ yields the $\Lambda$-linear presentable 6-functor formalism
\begin{align*}
    \D_\rel(-, \Lambda) := \D_\mot^\oc(-) \tensor_{\D_\mot^\oc(*), r} \D(\Lambda).
\end{align*}
By construction $\D_\rel(-, \Lambda)$ has the same $!$-able maps as $\D_\mot^\oc$ and there is a natural transformation of 6-functor formalisms $\D_\mot^\oc \to \D_\rel$.
\end{defn}

The notation \enquote{rel} indicates motivic sheaves realized \emph{relatively} over the base geometric point. If $\Lambda$ is killed by a power of $\ell$ then $\D_\rel(X,\Lambda)$ agrees with the overconvergent subcategory of $\D_\et(X,\Lambda)$ (see \cref{rslt:compare-D-rel-to-D-et} below). For a general non-torsion ring, this sheaf category is typically huge and completely inexplicit. In the following subsections we establish control over this category for some of the v-stacks appearing in \cite{FS}. Here is an easy example (see also \cite[\S B.4]{LHcourbes}):

\begin{exmpl} \label{rslt:D-rel-on-classifying-stack-is-smooth-reps}
Let $G$ be a $p$-adic Lie group (viewed as a v-stack). Then $\D_\rel(*/G, \Lambda)$ is the (derived) category of smooth $G$-representations on $\Lambda$-modules. Indeed, by descent (cf. \cite[Proposition~5.3.10]{heyer-mann-6ff}) this reduces to showing that $\D_\rel(S, \Lambda) = \D(\mathcal C(S, \Lambda))$ for all profinite sets $S$, where $\mathcal C(S, \Lambda)$ denotes the ring of continuous maps $S \to \Lambda$ (for the discrete topology on $\Lambda$). But this is clear if $S$ is finite and then follows for profinite $S$ by \cite[Lemma~10.4]{ScholzeBerkovich}.
\end{exmpl}

The following result summarizes the main abstract properties of $\D_\rel$.

\begin{thm}
Fix a $\Z_\ell$-algebra $\Lambda$.
\begin{thmenum}
    \item \label{rslt:!-able-maps-for-D-rel} $\D_\rel$ is a presentable 6-functor formalism on $\vStk$ which is sheafy, i.e. descends along all v-covers. The class $E$ of $!$-able maps satisfies the following properties: $E$ is closed under composition and base-change, it is right cancellative and it contains all partially proper maps that are representable in locally spatial diamonds of locally finite $\dimtrg$.    

    \item \label{rslt:compare-D-rel-to-D-et} For every $n \ge 1$ and every small v-stack $X$ there is a natural equivalence
    \begin{align*}
        \D_\rel(X, \Z/\ell^n) = \D_\et(X, \Z/\ell^n)^\oc \subseteq \D_\et(X, \Z/\ell),
    \end{align*}
    where $\D_\et(X, \Z/\ell^n)^\oc$ denotes the full subcategory spanned by the overconvergent étale sheaves. The above embedding is compatible with tensor products, pullbacks and lower-$!$ functors along the maps in (i), i.e. it upgrades to a natural transformation of 6-functor formalisms $\D_\rel(-, \Z/\ell^n) \injto \D_\et(-, \Z/\ell^n)$.
\end{thmenum}
\end{thm}
\begin{proof}
We first prove (i). We know that $\D_\mot^\oc$ satisfies v-descent (and even arc descent) by construction. To deduce v-descent for $\D_\rel$ it is enough to show that base-change $- \tensor_{\D_\mot^\oc(*)} \D(\Lambda)$ commutes with limits. But by \cite[Proposition~10.3]{ScholzeBerkovich} $\D_\mot^\oc(*)$ is rigid, hence $\D(\Lambda)$ is dualizable as a module over $\D_\mot^\oc(*)$ (see \cite[Ch. 1 Proposition 9.4.4]{gaitsgory-rozenblyum-vol1}). The claim about $!$-able maps and part (ii) essentially follow from \cite[Proposition~12.3]{ScholzeBerkovich} together with some arguments sketched in \cite{FS-motivic}; further details will be given in \cite{GH}. 
\end{proof}

\begin{rmk}
In the following, whenever we talk about cohomological properties of maps of small v-stacks (like \enquote{$!$-able} or \enquote{cohomologically smooth}), we implicitly refer to the 6-functor formalism $\D_\rel(-, \Z_\ell)$.    
\end{rmk}

In \cref{rslt:compare-D-rel-to-D-et} we stated a comparison of the inexplicit 6-functor formalism $\D_\rel(-, \Z/\ell^n)$ to the much better understood étale 6-functor formalism $\D_\et(-, \Z/\ell^n)$. This immediately leads to the question of whether a similar comparison is possible for general $\Z_\ell$-algebras $\Lambda$ in place of $\Z/\ell^n$. We refer the reader to the forthcoming paper \cite{GH} for a discussion of this question in great detail and content ourselves with a brief summary of their results, as this is all we will need.

For general $\Z_\ell$-algebras $\Lambda$ we have to replace $\D_\et(-, \Z/\ell^n)$ by a suitable candidate, for which we take the 6-functor formalism of \emph{nuclear} sheaves defined by the second author in \cite{mann-nuclear}. If $X$ is a spatial diamond with finite $\ell$-cohomological dimension, then $\D_{\nuc}(X, \Z_{\ell})$ is the full subcategory of $\D(X_{v},\Z_{\ell})$ generated by $\D_{\et}(X,\Z_{\ell})$ under filtered colimits.\footnote{Here $\D_{\et}(X,\Z_{\ell})$ denotes $\ell$-complete étale sheaves as in \cite[Section 26]{etale-cohomology-of-diamonds}. Slightly different notation is used in \cite{mann-nuclear}.} By the main results of \cite{mann-nuclear} (see \cite[Theorem~1.2]{mann-nuclear}), the assignment $X \mapsto \D_{\nuc}(X,\Z_{\ell})$ extends uniquely to a hypercomplete v-sheaf on the category of all small v-stacks. For a general $\Z_{\ell}$-algebra $\Lambda$, we set
\begin{align*}
    \D_{\nuc}(X,\Lambda) = \Mod_{\Lambda}\D_{\nuc}(X,\Z_{\ell}) = \D_{\nuc}(X,\Z_{\ell}) \otimes_{\D(\Z_{\ell})} \D(\Lambda).
\end{align*}
This also defines a hypercomplete v-sheaf, and \cite[Theorem~1.3]{mann-nuclear} constructs a full 6-functor formalism on this category. Note that $\D_\nuc(X, \Z/\ell^n) = \D_\et(X, \Z/\ell^n)$ by \cite[Proposition~3.20]{mann-nuclear}, so $\D_\nuc$ generalizes $\D_\et$ in an explicit way. Here is the promised comparison between $\D_\rel$ and $\D_\nuc$:

\begin{prop} \label{rslt:compare-D-rel-to-D-nuc}
Fix a $\Z_\ell$-algebra $\Lambda$. For every small v-stack $X$ there is a natural realization functor
\begin{align*}
    \rho\colon \D_\rel(X, \Lambda) \to \D_{\nuc}(X,\Lambda),
\end{align*}
which preserves all small colimits, is symmetric monoidal and is compatible with all pullbacks and base-change in $\Lambda$. It is uniquely determined by these conditions together with the requirement that for dualizable $A \in \D_\rel(X,\Z_\ell)$ we have $\rho(A) = \hat A$; here $\hat{A} = \varprojlim_n A/\ell^n \in \D_{\et}(X,\Z_{\ell})$ denotes the naive $\ell$-completion of $A$, where $A/\ell^n$ is viewed as an étale sheaf via \cref{rslt:compare-D-rel-to-D-et}.
\end{prop}
\begin{proof}
We refer to \cite{GH} for details and only sketch the proof. First note that the case of general $\Lambda$ reduces to the case of $\Lambda = \Z_\ell$ by applying $- \tensor_{\D(\Z_\ell)} \D(\Lambda)$, so from now on we assume $\Lambda = \Z_\ell$. Since both $\D_\rel$ and $\D_\nuc$ satisfy v-descent, we are further reduced to the case that $X$ is a strictly totally disconnected perfectoid space. In this case $\D_\nuc(X, \Z_\ell)^\oc$ is just the category of nuclear modules over the condensed $\Z_\ell$-algebra $\mathcal C(X, \Z_\ell)$ (see \cite[Lemma~6.7(iv)]{mann-nuclear}) and one can perform a similar construction as in \cref{rslt:ell-adic-comparison-for-D-mot} in order to obtain a comparison functor $\D_\mot^\oc(X) \to \D_\nuc(X, \Z_\ell)^\oc$, which is compatible with the functor $r$ from \cref{def:D-rel} by construction and hence provides the desired $\rho$. The main nontrivial input here is that $\D_\mot^\oc(X)$ is compactly generated by dualizable objects for strictly totally disconnected spaces, for which we refer to \cite{GH}.
\end{proof}

\begin{rmk}
The comparison functor $\rho$ from \cref{rslt:compare-D-rel-to-D-nuc} is not compatible with lower-$!$ functors in general, but there is a \emph{lax} compatibility. We refer the reader to \cite{GH} for details.
\end{rmk}

\subsection{Sheaves on \texorpdfstring{$\Bun_G$}{BunG}}

Fix a $\Z_\ell$-algebra $\Lambda$. In this section we quickly review the main geometric player, the stack $\Bun_G$, and its sheaf theory. By definition, $\Bun_G$ is the functor sending an affinoid perfectoid space $S$ over $\overline{\F}_q$ to the groupoid of $G$-torsors on the relative Fargues-Fontaine curve $X_S$. This turns out to be an Artin v-stack, which is smooth of dimension zero in any reasonable sense. At the level of topological spaces, we have a canonical identification $|\Bun_G| = B(G)$, where $B(G)$ is the Kottwitz set of $F$-isocrystals with $G$-structure \cite{Viehmann}. Given $ b \in B(G)$, we write $G_b$ for the $\sigma$-centralizer of $b$, and $i_b: \Bun_{G}^{b} \to \Bun_G$ for the inclusion of the associated locally closed stratum in $\Bun_G$. There is a natural map $s_b: \ast / G_b(F) \to \Bun_{G}^{b}$, which is a torsor for a smooth group diamond (see \cite[Proposition~III.5.3]{FS}).

Fix a $\mathbf{Z}_{\ell}$-algebra $\Lambda$. Our primary interest on the automorphic side is the sheaf category $\D_{\mathrm{lis}}(\Bun_G,\Lambda)$ defined by Fargues-Scholze in \cite[\S VII.7]{FS}. One can define $\D_{\mathrm{lis}}(X,\Lambda)$ is the full subcategory of $\D_{\nuc}(X,\Lambda)$ generated under colimits by sheaves of the form $f_!f^{!}\Lambda$ for all maps $f:Y \to X$ which are representable in locally spatial diamonds, locally compactifiable of finite dim.trg, and $\ell$-cohomologically smooth. This agrees with the definition as written in \cite[\S VII.7]{FS}, using \cite[Proposition 8.10]{mann-nuclear}. When $X=\ast$ is a point, $\D_{\mathrm{lis}}(\ast,\Lambda)=\D(\Lambda)$ is just the usual category of $\Lambda$-modules.

This approach to $\D_{\mathrm{lis}}(\Bun_G)$ has the technical drawback of not underlying a full six functor formalism. However, it turns out that $\D_{\rel}$ gives an alternative approach to this category which avoids this issue:

\begin{thm}\label{rslt:D-rel-on-BunG}
For any reductive group $G$ over $F$, the map $\Bun_G \to *$ is $!$-able and the functor $\rho$ from \cref{rslt:compare-D-rel-to-D-nuc} induces a symmetric monoidal equivalence of categories
    \begin{align*}
        \D_\rel(\Bun_G,\Lambda) \isoto \D_{\mathrm{lis}}(\Bun_G,\Lambda).
    \end{align*}
\end{thm}
This equivalence is one of the main results of \cite{GH}, and we do not prove it here. In what follows, we will simply write $\D(\Bun_G,\Lambda)$ as shorthand for $\D_{\rel}(\Bun_G,\Lambda) = \D_{\lis}(\Bun_G,\Lambda)$.

The key features of $\D(\Bun_G,\Lambda)$ are its compact generation and semiorthogonal gluing structure. To explain this, recall from \cite[Proposition~VII.7.1]{FS} that for any $b$, pullback along $s_b$ induces a canonical t-exact tensor equivalence
\[\D_\rel(\Bun_{G}^{b},\Lambda)\cong \D_\rel(\ast / G_b(F),\Lambda) = \D(G_b(F), \Lambda),\]
where the right-hand side is the derived category of smooth $G_b(F)$-representation on $\Lambda$-modules. We will always identify $\D_\rel(\Bun_{G}^{b},\Lambda)$ with $\D(G_b(F),\Lambda)$ in this way. Under this identification, pullback along $i_{b}$ induces a functor $i_{b}^{\ast}\colon D(\Bun_G,\Lambda) \to \D(G_b(F),\Lambda)$, which admits both a right adjoint $i_{b \ast}$ and a left adjoint $i_{b \sharp}$ (see \cite[Proposition~VII.7.2]{FS}). We also have a functor $i_{b!}$, which admits a right adjoint $i_{b}^{!}$. As in \cite{Beijing}, we can renormalize these functors after fixing a choice of $\sqrt{q} \in \Lambda$, and we write $i_{b}^{\ast \mathrm{ren}}$, $i_{b!}^{\ren}$, etc., for the renormalized functors.

The following result is a combination of several results from \cite{FS}.

\begin{prop} \label{rslt:D-Bun-G-is-compactly-generated}
The category $\D(\Bun_G,\Lambda)$ is compactly generated, and a sheaf $A$ is compact if and only if $i_{b}^{\ast}A$ is compact for all $b$ and vanishes for all but finitely many $b$. A set of compact generators is given by the sheaves $i_{b!}\cind_{K}^{G_b(F)}\Lambda$. Here $b$ runs over $B(G)$ and $K$ runs over pro-$p$ open compact subgroups of $G_b(F)$.
A second set of compact generators is given by the sheaves $i_{b\sharp}\cind_{K}^{G_b(F)}\Lambda$.
\end{prop}
\begin{proof}
The compact generation and the description of compact objects is proved in \cite[Proposition~VII.7.4]{FS}. It is then clear that the objects $i_{b\sharp} \cind_K^{G_b(F)} \Lambda$ form a set of compact generators. The generators $i_{b!} \cind_K^{G_b(F)} \Lambda$ are discussed in the proof of \cite[Proposition~VII.7.6]{FS}.
\end{proof}

Note that we could equally well take the renormalized $!$- and $\sharp$-pushforwards here. 

Next we briefly recall the dualities. Let $\pi_{G}:\Bun_G \to \ast$ be the structure map. It is well-known that $\Bun_G$ has trivial dualizing complex (see e.g. \cite[Proposition 2.6]{FS-motivic}), so the functor
\[ \Gamma_c(\Bun_G,-) := \pi_{G!}(-): \D(\Bun_G,\Lambda) \to \D(\Lambda) \]
is left adjoint to $\pi_{G}^{\ast}$.
By \cite[Proposition~VII.7.6]{FS} there is an involutive self-equivalence
\begin{align*}
    \Dbz\colon \D(\Bun_G,\Lambda)^{\omega,\op} \isoto \D(\Bun_G,\Lambda)^\omega
\end{align*}
given by Bernstein-Zelevinsky duality. This functor is characterized by the formula
\[\RHom(\Dbz(A),B) = \Gamma_c(\Bun_G,A \otimes B)\]
for compact $A$ and arbitrary $B$. In particular, $\D(\Bun_G,\Lambda)$ is self-dual for the evaluative self-duality associated with the functor $\Gamma_c(\Bun_G,-)$.

We also have Verdier duality, which is defined by the formula $\Dverd(A) = \intHom(A,\pi_{G}^{!}\Lambda)$ on all sheaves. (As discussed above, $\pi_{G}^{!}\Lambda \simeq \Lambda$ non-canonically.) This is an involutive self-equivalence on the subcategory of ULA sheaves, see \cite[Definition~VII.7.8]{FS}.

The 6-functor formalism $\D_\rel$ allows a conceptual interpretation of the above dualities. For the following result, recall the notion of prim and suave objects from \cite[\S 4.4]{heyer-mann-6ff}.

\begin{prop} \label{rslt:duality-on-BunG}
Fix a sheaf $A \in \D(\Bun_G, \Lambda)$.
\begin{propenum}
    \item \label{rslt:prim-objects-on-BunG} $A$ is prim iff it is compact iff it has finite support and all $i_b^* A$ are compact. Moreover, $\Dbz$ coincides with the prim duality functor.
    
    \item $A$ is suave iff it is ULA iff all $i_b^* A$ are admissible representations. Moreover $\Dverd$ coincides with the suave duality functor.
\end{propenum}
\end{prop}
\begin{proof}
We first prove (i). It follows from \cite[Proposition~3.2]{FS-motivic} that $\D(\Bun_G \times \Bun_G, \Lambda) = \D(\Bun_G) \tensor_{\D(\Lambda)} \D(\Bun_G)$, which implies that the 2-functor $\Psi_{\D_\rel}\colon\cat K_{\D_\rel} \to \Mod_{\D(\Lambda)}(\PrL)$ from \cite[Proposition~4.1.5]{heyer-mann-6ff} is fully faithful on the full subcategory spanned by $*$ and $\Bun_G$. By definition, prim objects in $\D(\Bun_G,\Lambda)$ coincide with right adjoint maps $\Bun_G \to *$ in $\cat K_{\D_\rel}$, so by the above fully faithfulness of $\Psi_{\D_\rel}$ they correspond to right adjoint maps $\D(\Bun_G,\Lambda) \to \D(\Lambda)$ in $\Mod_{\D(\Lambda)}(\PrL)$. Said differently, $A \in \D(\Bun_G,\Lambda)$ is prim if and only if the functor $\pi_{G!}(A \tensor -)$ admits a $\D(\Lambda)$-linear left adjoint. By \cite[Proposition~VII.7.6]{FS} this is equivalent to $A$ being compact. The comparison of dualities is immediate.

We now prove (ii). Since the prim objects form a set of compact generators of $\D(\Bun_G,\Lambda)$, it follows from \cite[Corollary~4.4.15]{heyer-mann-6ff} that $A$ is suave if and only if $\RHom(P, A) \in \D(\Lambda)$ is dualizable for all compact sheaves $P$. By \cref{rslt:D-Bun-G-is-compactly-generated} this reduces immediately to $P$ of the form $P = i_{b\sharp} \cind_K^{G_b(F)} \Lambda$, so we see that $A$ is suave if and only if for all $b$, $i_b^* A$ is an admissible $G_b(F)$-representation. But this is the same characterization as for ULA sheaves in \cite[Proposition~VII.7.9]{FS}.
\end{proof}

We have a duality exchange formula relating Bernstein-Zelevinsky duality and Verdier duality, exactly analogous to the duality exchange formula on $\ICoh$. In fact, the proof of parts of \cref{prop:admissiblefacts} is formal enough to adapt directly:

\begin{prop} 
\begin{propenum}
    \item \label{rslt:dualityexchangeBunG} For all $A\in \D(\Bun_G,\Lambda)^\omega$ and $B\in \D(\Bun_G,\Lambda)$, we have a natural isomorphism
    \begin{align*}
        \RHom(A, \Dverd(B)) = \RHom(\Dbz(A), B)^\vee
    \end{align*}
    in $\D(\Lambda)$.

    \item \label{prop:ULAreflect}Suppose that $(-)^\vee$ reflects perfect complexes in $\D(\Lambda)$, which is true in particular if $\Lambda$ is a field. Let $A\in \D(\Bun_G,\Lambda)$ be any sheaf. If $\Dverd(A)$ is ULA, then also $A$ is ULA.
\end{propenum}
\end{prop}
Here, as usual, $(-)^\vee$ denotes the naive dual $\RHom(-,\Lambda)$ in the derived category $\D(\Lambda)$.
\begin{proof}
Using \cref{rslt:prim-objects-on-BunG} the same arguments as in \cref{rslt:duality-exchange-formula} and \cref{rslt:admissible-duality-reflects-suaveness} work (they are formal in any 6-functor formalism).
\end{proof}



\subsection{Sheaves on \texorpdfstring{$\Div^1$}{Div1}}

In the previous subsections we introduced the 6-functor formalism $\D_\rel(-,\Lambda)$ on small v-stacks and computed it (without proof) on $\Bun_G$. In the present subsection we study $\D_\rel$ on the stack $\Div^1$ of Cartier divisors of degree $1$ (see \cite[\S II.1]{FS}). The naive expectation is that this category is closely related to a representation category of the Weil group $W_F$, but the precise formulation is a bit subtle: for instance, the pullback functor $\D_\rel(*, \Lambda) \to \D_\rel(\Spd \C_p, \Lambda)$ is not an equivalence in general (it is not even fully faithful, cf. \cite[Corollary~11.10]{ScholzeBerkovich}), so the analog of \cite[Proposition~IV.7.1]{FS} does not help us. Before we salvage the situation, let us introduce the relevant notation:

\begin{defn}
\begin{defenum}
    \item Fix a condensed ring $\Lambda$ with underlying discrete ring $\Lambda^\disc = \Lambda(*)$. Then the assignment $M \mapsto M \tensor_{\Lambda^\disc} \Lambda$ defines a colimit-preserving symmetric monoidal fully faithful embedding $\D(\Lambda^\disc) \injto \Cond(\Lambda)$, where $\Cond(\Lambda)$ denotes the (derived) category of condensed $\Lambda$-modules. We denote by
    \begin{align*}
        \D(\Lambda) \subseteq \Cond(\Lambda)
    \end{align*}
    the essential image of that functor, so that we get an isomorphism $\D(\Lambda^{\mathrm{disc}}) \isoto \D(\Lambda)$. Similarly, for a condensed group $H$ we denote by
    \begin{align*}
        \Rep_\Lambda(H) \subseteq \Cond(\Lambda[H])
    \end{align*}
    the full subcategory spanned by the condensed $H$-representations whose underlying condensed $\Lambda$-module lies in $\D(\Lambda)$.

    \item \label{def:smooth-and-continuous-reps-of-cond-group} Given a $\Z_\ell$-algebra $\Lambda$ we implicitly equip $\Lambda$ with a condensed structure via $\Lambda = \Lambda^\disc \tensor_{\Z_\ell^\disc} \Z_\ell$. For a condensed group $H$ we thus obtain the category $\Rep_\Lambda(H)$ as above. We also define $\Rep^\sm_\Lambda(H) := \Rep_{\Lambda^\disc}(H)$ and call it the category of \emph{smooth $H$-representations}. The base-change along $\Lambda^\disc \to \Lambda$ induces a symmetric monoidal functor
    \begin{align*}
        \Rep^\sm_\Lambda(H) \to \Rep_\Lambda(H).
    \end{align*}
\end{defenum}
\end{defn}

Explicitly, for a $\Z_\ell$-algebra $\Lambda$, $\Rep_\Lambda(H)$ is the category of continuous representations on $\Lambda$-modules $M$, where we implicitly equip $M$ with the colimit topology of the $\ell$-adic topologies on its finitely generated $\Z_\ell$-submodules. Similarly, $\Rep^\sm_\Lambda(H)$ is the category of continuous $H$-representations on $\Lambda$-modules, where we equip the modules with the discrete topology. The functor $\Rep^\sm_\Lambda(H) \to \Rep_\Lambda(H)$ takes a smooth representation and views it as a continuous one for the colimit topology. We warn the reader that this functor is not fully faithful in general (because it may not preserve $\Ext$-groups).

The main goal of the present subsection is to construct and analyze a canonical fully faithful comparison functor $\varepsilon\colon \D_\rel(\Div^1, \Lambda) \to \Rep_\Lambda(W_F)$. In fact, with the preparations from the previous subsections, the construction of the functor is straightforward (see \cref{rslt:compare-D-rel-Div1-to-Rep-W-F} below), but to get the full faithfulness of the comparison, we need the following lemma.

\begin{lem} \label{rslt:ell-completeness-of-suave-sheaves}
Let $f\colon X \to S$ be a $!$-able map of small v-stacks such that the tensor unit $\mathbf 1 \in \D_\rel(S, \Z_\ell)$ is $\ell$-complete. Then every $f$-suave object $M \in \D_\rel(X, \Z_\ell)$ is $\ell$-complete. In particular, if $f$ is suave (e.g. cohomologically smooth) then all dualizable objects in $\D_\rel(X, \Z_\ell)$ are $\ell$-complete.
\end{lem}
\begin{proof}
Note that $\omega_f := f^! \mathbf 1 \in \D_\rel(X, \Z_\ell)$ is $\ell$-complete, because $f^!$ preserves limits. This implies that for any $N \in \D_\rel(X, \Z_\ell)$ the object $\SD_f(N) := \intHom(N, \omega_f)$ is $\ell$-complete. But by suaveness of $M$ we have $M = \SD_f(\SD_f(M))$, hence $M$ is $\ell$-complete. For the second claim, if $f$ is suave then all dualizable objects on $X$ are $f$-suave (see \cite[Corollary~4.5.18(i)]{heyer-mann-6ff}).
\end{proof}
\begin{rmk}Despite its simplicity, this lemma turns out to be a crucial source of control over the $\D_{\rel}$-formalism. The point is that $\ell$-complete objects in $\D_{\rel}(-,\Z_\ell)$ come from the much better understood $\D_{\et}$-formalism, and this lemma gives a soft criterion for $\ell$-completeness. We will see this mechanism in action several times below.
\end{rmk}

\begin{prop} \label{rslt:compare-D-rel-Div1-to-Rep-W-F}
Fix a $\Z_\ell$-algebra $\Lambda$. Then there are natural colimit-preserving symmetric monoidal functors
\begin{align*}
    \Rep^\sm_\Lambda(W_F) \to \D_\rel(\Div^1, \Lambda) \overset{\varepsilon}{\injto} \Rep_\Lambda(W_F),
\end{align*}
whose composition is the natural map from \cref{def:smooth-and-continuous-reps-of-cond-group}. Moreover, the second above functor is fully faithful.
\end{prop}
\begin{proof}
We first observe that there is a natural symmetric monoidal embedding $\D_\nuc(\Div^1, \Lambda) \subseteq \Cond(\Lambda[W_F])$, identifying $\D_\nuc(\Div^1, \Lambda)$ with the full subcategory of condensed $W_F$-representations whose underlying $\Lambda$-module is nuclear (i.e. can be written as a filtered colimit of $\ell$-complete $\Z_\ell$-modules). Indeed, this follows by descent using that $\D_\nuc(\Spd \C_p, \Lambda) = \D_\nuc(\Lambda) \subseteq \Cond(\Lambda)$ and $\Div^1 = \Spd \C_p / W_F$, see \cite[Lemma~10.5]{mann-nuclear} for the same proof for $*$ in place of $\Spd \C_p$.

From the previous paragraph we see that the functor $\rho\colon \D_\rel(\Div^1, \Lambda) \to \D_\nuc(\Div^1, \Lambda)$ from \cref{rslt:compare-D-rel-to-D-nuc} induces a colimit-preserving symmetric monoidal functor
\begin{align*}
    \D_\rel(\Div^1, \Lambda) \to \Cond(\Lambda[W_F]).
\end{align*}
Since $\rho$ is compatible with pullback and its image on $\Spd \C_p$ is $\D(\Lambda) \subseteq \Cond(\Lambda)$, we immediately see that the above functor factors over $\Rep_\Lambda(W_F)$. Furthermore, by \cref{rslt:D-rel-on-classifying-stack-is-smooth-reps} the pullback along $\Div^1 \to */W_F$ induces a functor $\Rep^\sm_\Lambda(W_F) \to \D_\rel(\Div^1, \Lambda)$ and since $\rho$ is compatible with pullbacks we see that the composition of this functor with $\rho$ does indeed give the natural map $\Rep^\sm_\Lambda(W_F) \to \Rep_\Lambda(W_F)$.

It remains to see that $\rho$ is fully faithful on $\Div^1$, for which we can and will assume that $\Lambda = \Z_\ell$. By descent along the cover $\Spd \Breve{F} \surjto \Div^1$ it is enough to show that $\rho$ is fully faithful on $\Spd \Breve{F}$ (note that all terms in the \v{C}ech nerve of that cover are of the form $\Spd \Breve{F} \times \Z^n$), i.e. that the natural functor
\begin{align*}
    \D_\rel(\Spd \Breve{F}, \Z_\ell) \to \Rep_{\Z_\ell}(I_F)
\end{align*}
is fully faithful. By the proof of \cite[Proposition~10.3]{ScholzeBerkovich} (the case $X = \mathcal M_{\mathrm{arc}}(C)/G$) we see that $\D_\mot^\oc(\Spd \Breve{F})$ is compactly generated by dualizable objects. The same argument works directly for $\D_\rel(\Spd \Breve{F}, \Z_\ell)$, hence also this category is compactly generated by dualizable objects. Since $\rho$ is symmetric monoidal, it preserves dualizable objects and one checks that the dualizable objects in $\Rep_{\Z_\ell}(I_F)$ are compact (this reduces to $\mathbf 1 \in \Rep_{\Z_\ell}(I_F)$ being compact, which follows from the finite $\ell$-cohomological dimension of $I_F$, see \cite[Lemma~10.3(ii)]{mann-nuclear}). Altogether we see that $\rho$ restricts to a functor on compact objects
\begin{align*}
    \D_\rel(\Spd \Breve{F}, \Z_\ell)^\omega \to \Rep_{\Z_\ell}(I_F)^\omega
\end{align*}
and it is enough to show that this functor is fully faithful. Now the map $\Spd \Breve{F} \to *$ is cohomologically smooth (see \cite[Proposition~2.2]{FS-motivic} and note that $\Spd \Breve{F} \to \Div^1$ is étale), and we've already noted that all objects of $\D_\rel(\Spd \Breve{F}, \Z_\ell)^\omega$ are dualizable, so applying \cref{rslt:ell-completeness-of-suave-sheaves} we see that all objects in $\D_\rel(\Spd \Breve{F}, \Z_\ell)^\omega$ are in fact $\ell$-complete. This implies that the natural map $\D_\rel(\Spd \Breve{F}, \Z_\ell)^\omega \to \varprojlim_n \D_\rel(\Spd \Breve{F}, \Z/\ell^n)$ is fully faithful. A similar statement holds for $\Rep_{\Z_\ell}(I_F)^{\omega}$. This reduces the full faithfulness claim for $\rho$ to the case of coefficients in $\Z/\ell^n$ for $n \ge 1$. But then it is even an equivalence by \cref{rslt:compare-D-rel-to-D-et}.
\end{proof}

\begin{rmks}
\begin{rmksenum}
    \item In this paper we do not require the full faithfulness of the functor $\varepsilon$ constructed in \cref{rslt:compare-D-rel-Div1-to-Rep-W-F}. We warn the reader that the analogue of this full faithfulness is \emph{not} true for $\Spd \C_p$: Even though $\D_\rel(\Spd \C_p, \Z_\ell)$ is compactly generated by dualizable objects, not all of them are $\ell$-complete. It is somewhat of a miracle that the full faithfulness holds for $\Spd \Breve{F}$.

    \item In \cite[\S B.10]{LHcourbes} there is the explicit description $\D_\rel(\Div^1, \Lambda) \simeq \QCoh(\Spec \Lambda / \mathrm{WD}_F)$ for a certain group scheme $\mathrm{WD}_F$ over $\Lambda$ that sits in a short exact sequence
    \begin{align*}
        1 \to \varprojlim_n \mathbb G_a \to \mathrm{WD}_F \to W_F^\sm \to 1
    \end{align*}
    (cf. \cite[\S4.1]{FS-motivic}). Here $W_F^\sm$ denotes the group scheme over $\Lambda$ underlying space $\bigsqcup_\Z I_F^\sm$ with $I_F^\sm = \Spec \operatorname{LC}(I_F, \Lambda)$, where $\operatorname{LC}$ denotes locally constant maps. Let us denote $I_F = \Spec \mathcal C(I_F, \Lambda)$, where $\mathcal C$ denotes continuous maps (in the condensed sense) and let $W_F$ denote the induced group scheme over $\Lambda$. Then we have the following maps of group schemes over $\Lambda$:
    \begin{align*}
        W_F \to \mathrm{WD}_F \to W_F^\sm.
    \end{align*}
    Namely, the second map was already introduced above and the first map can be constructed in the case $\Lambda = \Z_\ell$ using that $\mathcal C(I_F, \Z_\ell)$ is $\ell$-complete, which reduces the construction modulo $\ell$ (where it is even an isomorphism because $\varprojlim_n \mathbb G_a$ vanishes modulo $\ell$). We expect the functors in \cref{rslt:compare-D-rel-Div1-to-Rep-W-F} to be compatible, via Li-Huerta's equivalence, with the pullback functors on quasicoherent sheaves on the classifying stacks associated to the above groups. However, the equivalence proved in \cite[\S B.10]{LHcourbes} depends on a number of noncanonical choices which are avoided by our construction presented above.
\end{rmksenum}
\end{rmks}


We do not expect the embedding $\D_\rel(\Div^1, \Lambda) \injto \Rep_\Lambda(W_F)$ from \cref{rslt:compare-D-rel-Div1-to-Rep-W-F} to be an equivalence, but we also do not know a completely explicit description of its essential image. The following results provide an upper and a lower bound on the image.

\begin{prop} \label{rslt:D-rel-upper-bound}
Fix a $\Z_\ell$-algebra $\Lambda$. Then:
\begin{propenum}
    \item For every open subgroup $P \subseteq P_F$ of the wild inertia subgroup $P_F \subset W_F$, the inflation functor
    \begin{align*}
        \Rep_\Lambda(W_F/P) \injto \Rep_\Lambda(W_F)
    \end{align*}
    is fully faithful.

    \item \label{rslt:D-rel-Div1-maps-to-colim-W-F-mod-P} The embedding $\D_\rel(\Div^1, \Lambda) \injto \Rep_\Lambda(W_F)$ factors naturally as
    \begin{align*}
        \D_\rel(\Div^1, \Lambda) \to \varinjlim_P \Rep_\Lambda(W_F/P) \to \Rep_\Lambda(W_F),
    \end{align*}
    where the colimit ranges over all $P \subseteq W_F$ as in (i) and is computed in $\PrL$.
\end{propenum}
\end{prop}
\begin{proof}
We can assume $\Lambda = \Z_\ell$, as the general statement follows by passing to $\Lambda$-modules (i.e. applying $- \tensor_{\D(\Z_\ell)} \D(\Lambda)$).

For (i), we argue as in the first part of the proof of \cite[Proposition~IX.5.1]{FS}: It is enough to prove the claim on the ambient categories of solid representations, i.e. we need to show that the pullback along $f\colon */W_F \to */(W_F/P)$ is fully faithful on $\D_\solid(-, \Z_\ell)$; but this map is pro-étale, hence the pullback has a left adjoint $f_\natural$ that satisfies the projection formula, and we are left with proving the standard fact $f_\natural \Z_\ell = \Z_\ell$, i.e. that $P$ has vanishing $\Z_\ell$-homology.

We now prove (ii). Given any $P \subseteq W_F$ as in (i), let $\Rep'_{\Z_\ell}(W_F/P) \subseteq \Rep_{\Z_\ell}(W_F/P)$ be the full subcategory generated under colimits by the compact objects. The embedding $\Rep_{\Z_\ell}(W_F/P) \injto \Rep_{\Z_\ell}(W_F)$ from (i) preserves compact objects, because its right adjoint computes $P$-cohomology, which commutes with colimits as $P$ is pro-$p$ (see \cite[Proposition~10.9(ii)]{mann-nuclear}). By \cite[Corollary~A.5.9]{heyer-mann-6ff} the colimit $\varinjlim_P \Rep'_{\Z_\ell}(W_F/P)$ is the $\Ind$-category of the union of the compact objects in $\Rep_{\Z_\ell}(W_F/P)$ for all $P$, and since all of these are compact in $\Rep_{\Z_\ell}(W_F)$, we see that the map
\begin{align}
    \varinjlim_P \Rep'_{\Z_\ell}(W_F/P) \injto \Rep_{\Z_\ell}(W_F) \label{eq:colim-of-Rep-W-F-mod-P}
\end{align}
is fully faithful. It now remains to show that the embedding $\D_\rel(\Div^1, \Z_\ell) \injto \Rep_{\Z_\ell}(W_F)$ factors over the image of the functor in \cref{eq:colim-of-Rep-W-F-mod-P}. By the proof of \cref{rslt:compare-D-rel-Div1-to-Rep-W-F} the category $\D_\rel(\Spd \Breve{F}, \Z_\ell)$ is generated by compact dualizable objects, which formally implies that $\D_\rel(\Div^1, \Z_\ell)$ is compactly generated by the objects $h_! \mathcal P$ for $h\colon \Spd \Breve{F} \surjto \Div^1$ the projection and $\mathcal P \in \D_\rel(\Spd \Breve{F}, \Z_\ell)$ compact. It is therefore enough to show that every $h_! \mathcal P$ lies in the image of \cref{eq:colim-of-Rep-W-F-mod-P}. We observe that the functor $\rho$ from \cref{rslt:compare-D-rel-to-D-nuc} commutes with $h_!$: Since $h$ is étale, we have $h_! = h_\natural$, so the commutativity is a condition and not a datum, and by compatibility of $\rho$ with base-change along $h$, one easily reduces to showing that $\rho$ commutes with $!$-push along $\Spd \Breve{F} \times \Z \to \Spd \Breve{F}$, which is clear (we refer to \cite{GH} for a much more general commutation of $\rho$ with $!$-functors). Using the commutation of $\rho$ with $h_!$ we can replace $\Div^1$ by $\Spd \Breve{F}$ and $W_F$ by $I_F$, and using that the compact objects in $\D_\rel(\Spd \Breve{F}, \Z_\ell)$ are dualizable (hence the same is true for their image in $\Rep_{\Z_\ell}(I_F)$), we are reduced to the following claim:
\begin{itemize}
    \item[($*$)] Every dualizable object $V \in \Rep_{\Z_\ell}(I_F)$ is trivial on some $P \subseteq I_F$ as in (i), i.e. lies in the image of the embedding $\Rep_{\Z_\ell}(I_F/P) \injto \Rep_{\Z_\ell}(I_F)$.
\end{itemize}
To prove this claim, fix $V$, i.e. $V$ is a continuous $I_F$-representation on a perfect $\Z_\ell$-module (equipped with the $\ell$-adic topology). In particular $V$ is $\ell$-complete, i.e. $V = \varprojlim_n V/\ell^n$. Note that $V/\ell$ is a smooth representation on a perfect $\F_\ell$-module and hence comes from a $I_F/P$-representation for some $P$ (e.g. because this is true for the compact smooth representations). By using the fiber sequences $V/\ell^n \to V/\ell^{n+1} \to V/\ell$ and the full faithfulness in (i) (for $I_F$ in place of $W_F$) we see inductively that our chosen $P$ acts trivially on all $V/\ell^n$, hence also on their limit, as desired.
\end{proof}

\begin{prop} \label{rslt:image-of-Div1-in-Rep-W-F}
Fix a $\Z_\ell$-algebra $\Lambda$ and assume that $\Lambda \tensor_{\Z_\ell} \Q_\ell$ is either $0$ or a finite extension of $\Q_\ell$. Then the essential image of the embedding $\D_\rel(\Div^1, \Lambda) \injto \Rep_\Lambda(W_F)$ from \cref{rslt:compare-D-rel-Div1-to-Rep-W-F} contains every representation whose underlying $\Lambda$-module is perfect.
\end{prop}
\begin{proof}
Fix $V \in\Rep_\Lambda(W_F)$ such that $V$ has perfect underlying $\Lambda$-module (i.e. it is dualizable). Recall that we have a factorization $\Rep^\sm_\Lambda(W_F) \to \D_\rel(\Div^1, \Lambda) \to \Rep_\Lambda(W_F)$ from \cref{rslt:compare-D-rel-Div1-to-Rep-W-F}, so in order to show that $V$ comes from $\D_\rel(\Div^1, \Lambda)$, it is enough to show that $V$ can be built using colimits from the image of $\Rep^\sm_\Lambda(W_F) \to \Rep_\Lambda(W_F)$. Consider the fiber sequence
\begin{align*}
    V \to V \tensor_{\Z_\ell} \Q_\ell \to V \tensor \Q_\ell/\Z_\ell = \varinjlim_n (V \tensor \Z/\ell^n).
\end{align*}
We know that modulo $n$, $\Rep_{\Lambda/\ell^n}^\sm(W_F) = \Rep_{\Lambda/\ell^n}(W_F)$, hence all $V \tensor \Z/\ell^n$ lie in the image of $\D_\rel(\Div^1, \Lambda)$. It remains to show that the same is true for $V \tensor_{\Z_\ell} \Q_\ell$, so from now on we can assume that $\Lambda$ is a finite extension of $\Q_\ell$, which then immediately reduces to the case $\Lambda = \Q_\ell$. Note that the t-structure on $\Cond(\Q_\ell[W_F])$ restricts to a t-structure on the subcategory of dualizable objects, because they are characterized by the underlying $\Lambda$-module being perfect, i.e. being bounded with all cohomologies finite-dimensional. We can build $V$ using cofiber sequences from its cohomologies, thus reducing to the case that $V$ is static. Thus $V$ is now a classical object: It is a finite-dimensional $\Q_\ell$-vector space equipped with the induced topology and a continuous $\Q_\ell$-linear $W_F$-action. We can further assume that $V$ is irreducible: Otherwise there is some non-trivial subrepresentation $V' \subseteq V$ and hence a fiber sequence $V' \to V \to V/V'$ and it is enough to show that both $V'$ and $V/V'$ lie in the image of $\D_\rel(\Div^1, \Q_\ell)$; but they have lower dimension. Finally, if $V$ is irreducible then by classical results (see \cite[Corollary~4.2.3]{Tate-NumberTheoreticBackground}) $V$ is smooth, as desired.
\end{proof}

\begin{rmk}
In \cref{rslt:image-of-Div1-in-Rep-W-F} it is crucial to work with $\Div^1$ and $W_F$. The analogous statement for $\Spd \Breve{F}$ and $I_F$ likely fails: It is the action of Frobenius in $W_F$ that forces an automatic smoothness on irreducible representations, see \cite[Corollary~4.2.3]{Tate-NumberTheoreticBackground}.
\end{rmk}

\begin{rmk}
As mentioned above, we don't have a concrete description of the essential image of the embedding $\D_\rel(\Div^1, \Lambda) \injto \Rep_\Lambda(W_F)$. It is plausible that (in the notation of the proof of Proposition \ref{rslt:D-rel-upper-bound}) this embedding factors over an equivalence $\D_{\rel}(\Div^1,\Lambda) \cong  \varinjlim_P \Rep'_{\Lambda}(W_F/P)$, but we have not seriously tried to prove this.
\end{rmk}

Having studied sheaves on $\Div^1$, we now turn our attention to $\Div^I := (\Div^1)^I$ for a finite set $I$. An important input in the construction of the spectral action in \cite{FS} is Drinfeld's lemma \cite[Proposition~IV.7.3]{FS} which characterizes dualizable étale sheaves on $\Div^I$. We have the following analog for $\D_\rel$ (but see the warning in \cref{rmk:warning-about-drinfelds-lemma} below):

\begin{prop}
Fix a $\Z_\ell$-algebra $\Lambda$.
\begin{propenum}
    \item \label{rslt:D-rel-Div-tensor-I-fully-faithful} The symmetric monoidal functor
    \begin{align*}
        \boxtimes\colon \D_\rel(\Div^1, \Lambda)^{\tensor I} \injto \D_\rel(\Div^I, \Lambda)
    \end{align*}
    is fully faithful. Here the tensor product on the left is over $\D(\Lambda)$.

    \item \label{rslt:Drinfelds-lemma-for-D-rel} Suppose that $\Lambda \tensor_{\Z_\ell} \Q_\ell$ is either $0$ or a finite extension of $\Q_\ell$. Then there is a natural equivalence
    \begin{align*}
        (\D_\rel(\Div^1, \Lambda)^{\tensor I})^\dbl = \Rep_\Lambda(W_F^I)^\dbl,
    \end{align*}
    where $(-)^\dbl$ denotes the full subcategory of dualizable objects. If $\Lambda$ is finite over $\Z_\ell$ then these categories also agree with $\D_\rel(\Div^I, \Lambda)^\dbl$.
\end{propenum}
\end{prop}
\begin{proof}
We first prove (i). Recall the stack of norms $\mathcal N$ from \cref{sec:def-of-D-rel}. We base-change everything to $\mathcal N$, i.e. we consider the commuting diagram of functors
\begin{equation*}\begin{tikzcd}
    \D_\rel(\Div^1, \Lambda)^{\tensor I} \arrow[r] \arrow[dr] & \D_\rel(\Div^1 \times \mathcal N, \Lambda)^{\tensor_{\D_\rel(\mathcal N, \Lambda)} I} \arrow[r] & \D_\rel(\Div^I \times \mathcal N, \Lambda)\\
    & \D_\rel(\Div^I, \Lambda) \arrow[ur]
\end{tikzcd}\end{equation*}
The left diagonal map is the one we want to be fully faithful. By the commutativity of the diagram, it is enough to show that all the other maps are fully faithful. For the right diagonal map this follows from the fact that the pullback $\D_\rel(X, \Lambda) \to \D_\rel(X \times \mathcal N, \Lambda)$ is fully faithful for every small v-stack $X$, by the same argument as in \cref{rslt:D-mot-of-point-vs-N}. For the left horizontal map, we first observe that $\D_\rel(\mathcal N, \Lambda) = \D(\Lambda)$ by \cref{rslt:D-mot-of-point-vs-N}, hence by \cref{rslt:dualizable-ff-trick} the claim follows from the full faithfulness of the pullback $\D_\rel(\Div^1, \Lambda) \to \D_\rel(\Div^1 \times \mathcal N, \Lambda)$ and the fact that this pullback admits a $\D(\Lambda)$-linear left adjoint by cohomological smoothness of $\mathcal N$. It remains to prove the full faithfulness of the right horizontal map. Applying \cref{rslt:dualizable-ff-trick} to $- \tensor_{\D_\mot^\oc(*)} \D(\Lambda)$ we can reduce to $\D_\mot^\oc$ instead of $\D_\rel$. Now since $\Div^1 \to *$ is proper, $\Div^1 \times \mathcal N$ is qcqs and we can apply \cite[Corollary~10.6]{ScholzeBerkovich} (using that over $\mathcal N$ we are in the normed situation of \cite{ScholzeBerkovich}, cf. the proof of \cref{rslt:D-mot-of-point-vs-N}). This finishes the proof of (i).

The proof of part (ii) uses the same ideas as in the case $\abs I = 1$ discussed above, but with $\D_\rel(\Div^1, \Lambda)^{\tensor I}$ instead of $\D_\rel(\Div^I, \Lambda)$. We first show that the natural functor
\begin{align}
    \D_\rel(\Div^1, \Lambda)^{\tensor I} \to \Rep_\Lambda(W_F)^{\tensor I} \to \Rep_\Lambda(W_F^I) \label{eq:D-rel-Div-1-tensor-I-fully-faithful-in-Rep}
\end{align}
is fully faithful, where the first functor is induced from the embedding from \cref{rslt:compare-D-rel-Div1-to-Rep-W-F}. By \cref{rslt:dualizable-ff-trick} this reduces to the case $\Lambda = \Z_\ell$ by applying $- \tensor_{\D(\Z_\ell)} \D(\Lambda)$; hence from now on we assume $\Lambda = \Z_\ell$. To prove the full faithfulness, note that since $\D_\rel(\Div^1, \Z_\ell)$ is compactly generated and hence dualizable over $\D(\Z_\ell)$ (cf. \cref{rmk:how-to-apply-dualizable-ff-trick}), we can pull out the descent limit along the cover $\Spd \Breve{F} \surjto \Div^1$, whose terms are of the form $\Breve{F} \times \Z^k$. The $\Z^k$-factor just adds an additional product over $\Z^k$ and can be ignored. We can thus replace $\Div^1$ by $\Breve{F}$ and $W_F^I$ by $I_F^I$ and need to show that the functor
\begin{align*}
    \D_\rel(\Breve{F}, \Z_\ell)^{\tensor I} \to \Rep_{\Z_\ell}(I_F^I)
\end{align*}
is fully faithful. To see this, we note that the pullback functor $\D_\rel(\Div^1, \Z_\ell)^{\tensor I} \to \D_\rel(\Breve{F}, \Z_\ell)^{\tensor I}$ admits a $\D(\Z_\ell)$-linear left adjoint: Since $- \tensor_{\D(\Z_\ell)} -$ is a 2-functor, this reduces to the case $\abs I = 1$, where it follows from cohomological étaleness of the map $\Breve{F} \to \Div^1$.  The existence of the left adjoint implies that the above functor preserves limits. Now observe that the tensor unit in $\D_\rel(\Div^I, \Z_\ell)$ is $\ell$-complete (e.g. by \cref{rslt:ell-completeness-of-suave-sheaves}), hence by (i) the same is true for the tensor unit in $\D_\rel(\Div^1, \Z_\ell)^{\tensor I}$ and then follows for the tensor unit in $\D_\rel(\Breve{F}, \Z_\ell)^{\tensor I}$ by the discussed limit-preservation of the pullback. It follows formally that all dualizable objects in $\D_\rel(\Breve{F}, \Z_\ell)^{\tensor I}$ are $\ell$-complete. On the other hand, by \cref{rslt:compact-generation-of-tensor-over-Lambda} this category is compactly generated by dualizable objects (because the similar statement for $\abs I = 1$ is true) and these objects get sent to dualizable (hence compact and $\ell$-complete) objects in $\Rep_{\Z_\ell}(I_F^I)$. This formally reduces the claimed full faithfulness to the full subcategory of compact objects, and then via passage to a limit over $\Mod_{\Z/\ell^n}(-)$ to the case of coefficients in $\Z/\ell^n$ instead of $\Z_\ell$. Now $\D_\rel(\Breve{F}, \Z/\ell^n) = \D_\et(\Breve{F}, \Z/\ell^n) = \Rep_{\Z/\ell^n}(I_F)$, so the claim is that $\Rep_{\Z/\ell^n}(I_F)^{\tensor I} \to \Rep_{\Z/\ell^n}(I_F^I)$ is fully faithful. This is even an equivalence: Using dualizability of $\Rep_{\Z/\ell^n}(I_F)$ (this category is compactly generated) we can descend along $* \to */I_F$ in order to reduce to $\D_\et(I_F^k, \Z/\ell^n)^{\tensor I} = \D_\et((I_F^k)^I, \Z/\ell^n)$, which is clear by identifying both sides as module categories. This finishes the proof that \cref{eq:D-rel-Div-1-tensor-I-fully-faithful-in-Rep} is fully faithful.

We observe that $\Rep_\Lambda^\sm(W_F)^{\tensor I} = \Rep_\Lambda^\sm(W_F^I)$ because these categories are compactly generated, hence dualizable, hence we can pull the descent limit out (as in the previous paragraph). Altogether we obtain the factorization
\begin{align*}
    \Rep_\Lambda^\sm(W_F^I) \to \D_\rel(\Div^1, \Lambda)^{\tensor I} \injto \Rep_\Lambda(W_F^I),
\end{align*}
analogous to \cref{rslt:compare-D-rel-Div1-to-Rep-W-F}. We can now argue as in \cref{rslt:image-of-Div1-in-Rep-W-F} to show that the image of the second arrow contains all dualizable objects in $\Rep_\Lambda(W_F^I)$: By the arguments there, this reduces to showing that every irreducible classical continuous representation of $W_F^I$ on a finite-dimensional $\Q_\ell$-vector space lies in the image of the composite functor from smooth representations; this reduces to the case $|I|=1$ which was already handled in the proof of \cref{rslt:image-of-Div1-in-Rep-W-F}.\footnote{More precisely, if $V$ is an irreducible continuous representation of $W_F^I$ on a finite-dimensional $\Q_\ell$-vector space, we can find some finite Galois extension $E/\Q_\ell$ such that $V \otimes E = \oplus_k V_k$ where each $V_k$ is an absolutely irreducible $E$-linear continuous representation of $W_{F}^{I}$. Since each $V_k$ is absolutely irreducible, the usual argument shows that $V_k =\boxtimes_{i \in I} W_{k,i}$ where each $W_{k,i}$ is an irreducible continuous $E$-linear representation of $W_F$. But then all $W_{k,i}$ are smooth by \cite[Corollary 4.2.3]{Tate-NumberTheoreticBackground}, so then each $V_k$ is smooth, so then $V \subset V \otimes E = \oplus_k V_k$ is smooth as desired.} This proves the first equivalence in (ii). To get the second equivalence, we now assume that $\Lambda$ is finite over $\Z_\ell$ and can thus reduce to the case $\Lambda = \Z_\ell$. By (i) we only need to show that every dualizable object $\mathcal P \in \D_\rel(\Div^I, \Z_\ell)$ lies in the image of $\D_\rel(\Div^1, \Z_\ell)^{\tensor I}$. But by \cref{rslt:ell-completeness-of-suave-sheaves} we know that $\mathcal P$ is $\ell$-complete, hence lies in $\varprojlim_n \D_\et(\Div^I, \Z/\ell^n)$, so we conclude by the first equivalence in (ii) and Drinfeld's lemma \cite[Proposition~IV.7.3]{FS}.
\end{proof}

\begin{rmk} \label{rmk:warning-about-drinfelds-lemma}
We warn the reader that even though \cref{rslt:Drinfelds-lemma-for-D-rel} characterizes the dualizable objects of $\D_\rel(\Div^I, \Lambda)$ in terms of representations of $W_F^I$, it is not directly a generalization of the étale version in \cite[Proposition~IV.7.3]{FS}: Unless $\Lambda$ is torsion, the pullback functor $\D_\rel(*/W_F^I, \Lambda) \to \D_\rel(\Div^I, \Lambda)$ is never an equivalence on dualizable objects; this already fails for $\abs I = 1$ because the functor $\Rep^\sm_\Lambda(W_F) \to \Rep_\Lambda(W_F)$ is not fully faithful.
\end{rmk}

While \cref{rslt:Drinfelds-lemma-for-D-rel} is interesting in its own right, we only need the following corollary of it, which will ultimately allow us to reduce geometric Satake to the torsion case.

\begin{cor} \label{rslt:characterize-dualizable-sheaves-on-Div-I-via-completeness}
Fix a finite $\Z_\ell$-algebra $\Lambda$. Then a sheaf $\mathcal F \in \D_\rel(\Div^I, \Lambda)$ is dualizable if and only if it is $\ell$-complete and $\mathcal F/\ell^n \in \D_\et(\Div^I, \Lambda/\ell^n)$ is dualizable for all $n \ge 1$.
\end{cor}
\begin{proof}
The \enquote{only if} part is easy: By \cref{rslt:ell-completeness-of-suave-sheaves} every dualizable sheaf in $\D_\rel(\Div^I, \Lambda)$ is $\ell$-complete, and it is obviously also dualizable modulo $\ell^n$.

We now prove the \enquote{if} part of the claim, so assume that $\mathcal F$ is $\ell$-complete and all $\mathcal F/\ell^n$ are dualizable. Then $\mathcal F$ lies in $\varprojlim_n \D_\et(\Div^I, \Lambda/\ell^n)^\dbl$, which by Drinfeld's lemma \cite[Proposition~IV.7.3]{FS} agrees with $\Rep_\Lambda(W_F^I)^\dbl$. It follows from \cref{rslt:Drinfelds-lemma-for-D-rel} that $\mathcal F$ is contained in $\D_\rel(\Div^I, \Lambda)^\dbl$, as desired.
\end{proof}

\subsection{Sheaves on the local Hecke stacks}

In preparation for the Satake isomorphism for $\D_\rel$, in the following we provide some basic results on sheaves on the local Hecke stacks $\locHck_G^I$ defined in \cite[\S VI]{FS}. We start with the following variant of \cite[Proposition~VI.4.1]{FS}. It is implicitly also claimed in \cite[\S5]{FS-motivic}.

\begin{lem} \label{rslt:unipotent-invariance-for-D-rel}
Let $H$ be a group small v-sheaf over a small v-sheaf $S$ that admits an $\N$-indexed decreasing complete filtration $H^{\ge m} \subseteq H$ by closed normal subgroups such that, v-locally on $S$, for each $m \ge 1$ each quotient $H^{\ge m}/H^{\ge m+1}$ admits a further finite filtration with graded pieces given by $\mathbb A^1_{S^\sharp}$ for some untilt $S^\sharp$ of $S$ (that may vary). Let $X$ be a small v-stack over $S$ with an action of $H$ that factors over $H^{<m} = H/H^{\ge m}$ for some $m > 0$. Then the pullback functor
\begin{align*}
    \D_\rel(X/H^{<m}, \Lambda) \isoto \D_\rel(X/H, \Lambda)
\end{align*}
is an equivalence, for every $\Z_\ell$-algebra $\Lambda$.
\end{lem}
\begin{proof}
It is enough to prove the claim for $\D_\mot^\oc$ in place of $\D_\rel(-, \Lambda)$. By definition of Berkovich motives (see \cite[Definition~9.1]{ScholzeBerkovich}) we can further reduce to $\D_\mot^\eff(-, \Z[1/p])$ in place of $\D_\mot^\oc$. By the same reduction steps as in the proof of \cite[Proposition~VI.4.1]{FS} we can assume that $H^{<m} = 1$ and that each quotient $H^{\ge m'} / H^{\ge m'+1}$ admits a finite filtration with graded pieces $\mathbb A^1_{S^\sharp}$ for some untilts $S^\sharp$ of $S$, and we need to show that for every $\mathcal F \in \D_\mot^\eff(X, \Z[1/p])$ the natural map
\begin{align*}
    \mathcal F \to f_* f^* \mathcal F
\end{align*}
is an isomorphism, where $f\colon X \times H \to X$ is the projection. Further following the proof of \cite[Proposition~VI.4.1]{FS} we write $H$ as a filtered colimit of subgroups $H_j$ such that each $H_j$ is a successive extension as before, but now the quotients are closed balls inside $\mathbb A^1_{S^\sharp}$.\footnote{ We briefly sketch how to construct the subgroups $H_j$. For brevity set $H^{n}=H / H^{\geq n}$. After possibly passing to a finer filtration and renumbering, we can assume the successive quotients are $\mathbb{A}_{S^\sharp}^{1}$ (with the untilt possibly varying with $n$). We will inductively construct a sequence of subgroups $(H^{n}_{j} \subset H^{n})_{j \in J_n}$ indexed by a sequence of posets $J_n$ with compatible surjective maps $J_{n} \to J_{n-1}$. Setting $J = \lim_n J_n$ with $\pi_n: J \to J_n$ the evident map, the groups $H_{j} = \lim_{n} H_{\pi_n(j)}^{n}$ give the desired family.
When $n=1$, $H^{1}$ is an affine space. We set $J_{1} = \mathbf{N}$ and let $H_{j}^{n}$ be a sequence of balls of increasing unbounded radius. For the inductive step, suppose given $J_{n}$ and the associated family of groups $H_{j}^{n}$. The group $H^{n+1}$ sits in an extension $1 \to \mathbb{A}^{1}_{S^{\sharp}} \to H^{n+1} \to H^{n} \to 1$, and in particular is classified by a class $c \in H^2( B H^{n}, \mathbb{A}^{1}_{S^{\sharp}})$. For any $j \in J_{n}$, the restriction of this class along $H^2( B H^{n}, \mathbb{A}^{1}_{S^{\sharp}}) \to H^2( B H^{n}_{j}, \mathbb{A}^{1}_{S^{\sharp}})$ is the image of a class $c_{j,r} \in H^2( B H^{n}_{j}, \mathbb{A}^{1}_{S^{\sharp}}[r]) \to H^2( B H^{n}_{j}, \mathbb{A}^{1}_{S^{\sharp}})$ for all sufficiently large $r \geq r_{j}$, where $\mathbb{A}^{1}_{S^{\sharp}}[r] \subset \mathbb{A}^{1}_{S^{\sharp}}$ denotes a ball of radius $r$. Here we used the fact that $H_{j}^{n}$ is qcqs by constriction, so $H^2(B H_{j}^{n},-)$ commutes with filtered colimits (as is the case for any coherent topos). Define $J_{n+1}$ to be the set of pairs $(j,r)$ with $j \in J_{n}$ and $r \geq r_{j}$ a rational number, and for every such pair let $H^{n+1}_{(j,r)}$ be the evident quasicompact subgroup of $H^{n+1}$ classified by $c_{j,r}$. By design, $H^{n+1}_{(j,r)}$ surjects onto $H^{n}_{j}$ via the ambient surjection $H^{n+1} \to H^{n}$, with kernel a ball of radius $r$. Moreover, if $j \leq j'$ and $r' \geq \max (r,r_{j})$ there is a natural inclusion $H^{n+1}_{(j,r)} \subset H^{n+1}_{(j',r')}$. Giving $J_{n+1}$ the evident partial order, this completes the inductive construction of the desired family of subgroups.} Denoting $f_j\colon X \times H_j \to X$ the projection, we are reduced to showing that the natural map
\begin{align*}
    \mathcal F \to f_{j*} f_j^* \mathcal F
\end{align*}
is an isomorphism for all $j$. We now leave the realm of small v-stacks and instead work with arc stacks via the comparison $a'^*$ from \cite[\S12]{ScholzeBerkovich} ($\D_\mot^\eff$ factors through $a'^*$ anyway). To prove the above identity, it is enough to prove that for all arc-stacks $S'$ over $S$ the above map becomes an isomorphism on global sections. Replacing $S$ by $S'$, we are reduced to showing that the natural map
\begin{align*}
    \mathcal F(S) \to \mathcal F(S \times H_j)
\end{align*}
is an isomorphism. By \cite[Corollary~6.9]{ScholzeBerkovich} the category $\D_\mot^\eff(S, \Z[1/p])$ admits a left-complete t-structure. Thus by writing $\mathcal F = \varprojlim_n \tau^{\ge -n} \mathcal F$ and pulling the limit out of the claimed identity above, we reduce to the case that $\mathcal F$ is left-bounded. Now by design $H_j = \varprojlim_{m'} H_j^{<m'}$ is a cofiltered limit of qcqs arc stacks. Since $\mathcal F$ is finitary and hence commutes with filtered colimits of strictly totally disconnected spaces (see \cite[Definition~4.10]{ScholzeBerkovich}), the same argument as in \cite[Proposition~3.3.9(i)]{mann-p-adic-6-functors} shows that
\begin{align*}
    \mathcal F(S \times H_j) = \varinjlim_{m'} \mathcal F(S \times H_j^{<m'}).
\end{align*}
We are thus reduced to showing that each map $\mathcal F \to \mathcal F(S \times H_j^{<m'})$ is an isomorphism. But $H_j^{<m'}$ admits a finite filtration whose graded pieces are balls inside $\mathbb A^1_{S^\sharp}$, so the claim follows from ball-invariance of $\mathcal F$.
\end{proof}

As an important corollary of \cref{rslt:unipotent-invariance-for-D-rel} we obtain some control on $\D_\rel$ on the stack $\Div^I \!\! / L^+_I G$, which is crucial to control sheaves on the local Hecke stack.

\begin{prop} \label{rslt:ell-completeness-on-DivI-mod-L+G}
Fix a finite set $I$ and a reductive group $G$ over $F$. Then the monoidal unit in $\D_{\rel}(\Div^I \!\! / L^+_I G,\Z_\ell)$ is $\ell$-complete.
\end{prop}
\begin{proof}
By \cite[Proposition~VI.1.11]{FS} the group $L^+_I G$ admits an exhaustive decreasing $\N$-indexed filtration by closed subgroups $(L^+_I G)^{\ge m} \subseteq L^+_I G$ such that $H = (L^+_I G)^{\ge 1}$ satisfies the assumptions of \cref{rslt:unipotent-invariance-for-D-rel}. We deduce from \cref{rslt:unipotent-invariance-for-D-rel} that the pullback induces an equivalence
\begin{align*}
    \D_\rel(\Div^I \!\! / (L^+_I G)^{<1}, \Z_\ell) \isoto \D_\rel(\Div^I \!\! / L^+_I G, \Z_\ell),
\end{align*}
where $(L^+_I G)^{<1} = L_+^I G / (L^+_I G)^{\ge 1}$. We are therefore reduced to prove that the monoidal unit in $\D_\rel(\Div^I \!\! / (L^+_I G)^{<1}, \Z_\ell)$ is $\ell$-complete, which by \cref{rslt:ell-completeness-of-suave-sheaves} reduces to showing that $\Div^I \!\! / (L^+_I G)^{<1} \to *$ is cohomologically smooth. But $\Div^I \to *$ is cohomologically smooth by \cite[Proposition~2.2]{FS-motivic} and $\Div^I \!\! / (L^+_I G)^{<1} \to \Div^I$ is cohomologically smooth because $(L^+_I G)^{<1} \to \Div^I$ is cohomologically smooth by \cite[Proposition~VI.1.12]{FS} (whose proof adapts to $\D_\rel$) and cohomological smoothness is cohomologically smooth local on the source (see \cite[Lemma~4.5.8(i)]{heyer-mann-6ff}).
\end{proof}

An important tool in the proof of the Satake isomorphism for étale sheaves in \cite[\S VI]{FS} is the constant term functor. We need its basic results also for $\D_\rel$, cf. \cite[Proposition~VI.4.2]{FS} and \cite[Proposition~VI.6.4]{FS}. We first introduce the relevant definitions:

\begin{defn} \label{def:CT-on-Gr-for-D-rel}
Fix a $\Z_\ell$-algebra $\Lambda$, a finite set $I$, a split connected reductive group $G$ over $F$ and a cocharacter $\lambda\colon \mathbb G_m \to G$ with associated parabolic $P_\lambda \subseteq G$ and Levi quotient $M_\lambda$. Consider the associated diagram
\begin{equation*}\begin{tikzcd}
    & \Gr_{P_\lambda}^I \! / \, \mathbb G_m \arrow[dl,"g",swap] \arrow[dr,"f"]\\
    \Gr_G^I / \, \mathbb G_m && \Gr_{M_\lambda}^I \! / \, \mathbb G_m
\end{tikzcd}\end{equation*}
of small v-stacks over $\Div^I$ from \cite[Corollary~VI.3.5]{FS}, where $\mathbb G_m$ acts via the group homomorphism $\mathbb G_m \to L^+_I \mathbb G_m$ and the canonical action of $L^+_I H$ on $\Gr_H^I$ for $H \in \{ G, P_\lambda, M_\lambda \}$. We denote
\begin{align*}
    \CT_{P_\lambda} := f_! g^*\colon \D_\rel^\bdd(\Gr_G^I / \, \mathbb G_m, \Lambda) \to \D_\rel^\bdd(\Gr_{M_\lambda}^I \! / \, \mathbb G_m, \Lambda) 
\end{align*}
the associated \emph{constant term functor}. Here $\D_\rel^\bdd$ is defined as in \cref{rslt:3ff-on-bounded-sheaves-on-vstacks} for the stratification of the affine Grassmannian into closed Schubert cells (see \cite[Definition~VI.2.6]{FS}).

We will sometimes abuse notation and also denote by $\CT_{P_\lambda}$ the same functor on $\Gr_G^I$ instead of $\Gr_G^I / \, \mathbb G_m$; note that this functor only depends on the parabolic $P_\Lambda$ and not on $\lambda$.
\end{defn}

\begin{defn} \label{def:suave-bounded-sheaves-on-Hck}
Fix a ring $\Lambda$, a finite set $I$ and a connected reductive group $G$ over $F$. We say that a sheaf $\mathcal F \in \D_\rel^\bdd(\Gr_G^I, \Lambda)$ is \emph{suave over $\Div^I$} if it suave in the usual sense (see \cite[Definition~4.4.1]{heyer-mann-6ff}) when restricted to a large enough union of closed Schubert strata containing its support. We similarly define the notion of suave sheaves in $\D_\rel^\bdd(\locHck_{G}^{I}, \Lambda)$ over $\Div^I \!\! / L^+_I G$.
\end{defn}

\begin{rmk}A priori, there is an asymmetry in last part of this definition, as there are two natural projections $\locHck_{G}^{I} \to \Div^I / L^+_I G$. However, the proof of \cite[Proposition VI.6.2]{FS} goes through verbatim in the $\D_{\rel}$ formalism, showing that a bounded sheaf on $\locHck_{G}^{I}$ is suave for one of these projections iff it is suave for the other. Said differently, suave objects in $\D_\rel^\bdd(\locHck_{G}^{I}, \Lambda)$ are stable under pullback along the switching isomorphism $\mathrm{sw}$ induced by inversion on $L_{I}G$. We will use these symmetries without particular comment.
\end{rmk}


\begin{prop} \label{rslt:CT-on-Gr-properties}
Fix a $\Z_\ell$-algebra $\Lambda$, a finite set $I$, a split connected reductive group $G$ over $F$ and a Borel $B \subseteq G$ with maximal torus $T$. Consider the functor
\begin{align*}
    \CT'_B\colon \D_\rel^\bdd(\locHck_G^I, \Lambda) \to \D_\rel^\bdd(\Gr_G^I, \Lambda) \xto{\CT_B} \D_\rel^\bdd(\Gr_T^I, \Lambda),
\end{align*}
where the first functor is pullback along the projection $\Gr_G^I \to \locHck_G^I$. Then:
\begin{propenum}
    \item $\CT'_B$ is conservative and preserves all small limits and colimits (on a fixed closed stratum).
    \item Given $\mathcal F \in \D_\rel^\bdd(\locHck_G^I, \Lambda)$, then $\mathcal F$ is suave over $\Div^I \!\! / L^+_I G$ if and only if $\CT'_B(\mathcal F)$ is suave over $\Div^I$.
\end{propenum}
\end{prop}
\begin{proof}
We will use without proof that hyperbolic localization works for $\D_\rel(-, \Lambda)$, i.e. the results in \cite[\S IV.6]{FS} apply to $\D_\rel$. This is not entirely obvious as the proofs make use of quasicompact open subsets, which do not behave nicely in $\D_\rel$, but it is possible work with partially proper open immersions instead. Similar claims have been made in \cite[\S2]{FS-motivic} without proof. We refer the reader to \cite{HM-Hyperbolic} for an upcoming independent proof of hyperbolic localization.

With hyperbolic localization at hand, we can now argue similarly as in \cite[\S VI]{FS}. First, the proof of \cite[Proposition~VI.4.2]{FS} adapts without change to $\D_\rel$, proving the conservativity claim in (i). Now restrict to a closed Schubert stratum. It is clear that $\CT'_B$ preserves small colimits. To prove the claim about limits, note that by \cite[Proposition~VI.2.8]{FS} and \cref{rslt:unipotent-invariance-for-D-rel} the pullback along $\Gr_G^I \to \locHck_G^I$ preserves limits on the chosen Schubert cell. Moreover, this pullback factors through $\Gr_G^I / \, \mathbb G_m$ (by definition of the $\mathbb G_m$-action) and hence by \cite[Theorem~IV.6.5]{FS} $\CT_B$ preserves limits as well. This finishes the proof of (i).

The proof of (ii) works in the same way as in \cite[Proposition~VI.6.4]{FS}. In \cite{HM-Hyperbolic} a slightly different variation of this proof will be discussed, which avoids the question whether certain \enquote{natural} maps are the expected ones.
\end{proof}

The above results culminate in the following characterization of suave objects in $\D_\rel^\bdd(\locHck_G^I, \Lambda)$, which seems surprisingly subtle to attack directly.

\begin{cor} \label{rslt:characterize-suave-sheaves-on-locHck-via-ell-completeness}
Fix a finite $\Z_\ell$-algebra $\Lambda$, a finite set $I$, a connected reductive group $G$ over $F$ and some $\mathcal F \in \D_\rel^\bdd(\locHck, \Lambda)$. Then $\mathcal F$ is suave over $\Div^I \!\! / L^+_I G$ if and only if $\mathcal F$ is $\ell$-complete and $\mathcal F/\ell^n \in \D_\et(\locHck_G^I, \Lambda/\ell)^n$ is suave over $\Div^I \!\! / L^+_I G$ for all $n \ge 1$.
\end{cor}
\begin{proof}
We first prove the \enquote{only if} part, so suppose that $\mathcal F$ is suave. Then clearly the same is true for $\mathcal F/\ell^n$ because suaveness is stable under the base-change $- \tensor_\Lambda \Lambda/\ell^n$ (as this is a morphism of 6-functor formalisms). Moreover, by \cref{rslt:ell-completeness-of-suave-sheaves} and \cref{rslt:ell-completeness-on-DivI-mod-L+G} we see that $\mathcal F$ is $\ell$-complete.

We now prove the \enquote{if} part, so assume that $\mathcal F$ is $\ell$-complete and all $\mathcal F/\ell^n$ are suave. Proving that $\mathcal F$ is suave seems to be surprisingly subtle, as the pullback and tensor product on $\D_\rel$ do not preserve $\ell$-completeness in general. We argue as follows. By finite étale descent we can assume that $G$ is split. We fix a Borel $B$ with maximal torus $T$ and we consider the constant term functor $\CT_B$, or more precisely the functor
\begin{align*}
    \CT'_B\colon \D_\rel^\bdd(\locHck_G^I, \Lambda) \to \D_\rel^\bdd(\Gr_T^I, \Lambda)
\end{align*}
from \cref{rslt:CT-on-Gr-properties}. By that result, $\CT'_B$ detects and preserves suaveness and also preserves $\ell$-completeness. We can thus replace $\mathcal F$ by $\CT'_B(\mathcal F)$ in order to reduce to the following claim: If $\mathcal F' \in \D_\rel^\bdd(\Gr_T^I, \Lambda)$ is $\ell$-complete such that all $\mathcal F'/\ell^n$ are suave over $\Div^I$, then $\mathcal F'$ is suave over $\Div^I$. We have $\Gr_T^I \isom X_*(T) \times \Div^I$, which allows us to replace $\Gr_T^I$ by $\Div^I$. The claim now becomes: If $\mathcal F' \in \D_\rel(\Div^I, \Lambda)$ is $\ell$-complete and all $\mathcal F'/\ell^n$ are dualizable, then $\mathcal F'$ is dualizable. This is \cref{rslt:characterize-dualizable-sheaves-on-Div-I-via-completeness}.
\end{proof}

\subsection{Hecke action and spectral action} \label{sec:spectral-action}

With the preparations from the previous subsections at hand, we now review the construction of the Hecke action and the spectral action, which are the central results in \cite{FS}. In contrast to previous constructions of the spectral action, we have taken some extra care regarding the functoriality of the intermediate constructions.

Throughout this subsection we fix a finite quotient $Q$ of $W_F$ over which the action on $\hat{G}$ factors. In \cref{sec:abstract-Hecke-action} we have established general abstract results about the functoriality of the Hecke stack and bounded sheaves on it. We apply it now to $\D_\rel$. For the following result, recall that for a finite set $I$, the Hecke stack $\locHck_G^I = L^+_I G \backslash L^I_G / L^+_I G$ is a stack over $\Div^I$ that is equipped with the convolution monoidal structure. It is stratified by closed Schubert cells (see \cite[Proposition~VI.2.7]{FS}), so by \cref{def:ind-vstacks-and-bounded-sheaves} we obtain an associated category $\D_\rel^\bdd(\locHck_G^I) \subseteq \D_\rel(\locHck_G^I)$ of bounded sheaves on $\locHck_G^I$. The convolution monoidal structure on $\locHck_G^I$ induces a convolution structure on $\D_\rel^\bdd(\locHck_G^I)$. The following result shows that this monoidal structure is natural in $I$ and that there is an induced action on $\D(\Bun_G)$ via the global Hecke stack (similar to the introduction of \cite[\S IX]{FS}).

\begin{prop} \label{rslt:functoriality-of-D-rel-on-locHck}
Fix a $\Z_\ell$-algebra $\Lambda$. There is a natural $\D_\rel(\Div^I,\Lambda)$-linear monoidal functor
\begin{align*}
    \D_\rel^{\bdd}(\locHck_{G}^{I}, \Lambda) \to \End_{\D(\Lambda)}(\D(\Bun_G,\Lambda)) \otimes_{\D(\Lambda)} \D_\rel(\Div^I, \Lambda)
\end{align*}
which is functorial in $I$.\footnote{We refer the reader to \cref{rslt:full-functoriality-of-D-bdd-Hck} for a precise formulation of what we mean by a \enquote{functorial in $I$}.} Here the monoidal structure on the right-hand side is induced by composition of endofunctors in the first factor and tensor product in the second factor.
\end{prop}
\begin{proof}
Since $\D(\Bun_G,\Lambda) = \D_\rel(\Bun_G, \Lambda)$ and $\D(\Lambda) = \D_\rel(*)$, the claim is exactly the one in \cref{rslt:full-functoriality-of-D-bdd-Hck} for $\D = \D_\rel(-,\Lambda)$. Hence we only need to verify that $\D_\rel(-,\Lambda)$ satisfies condition ($*$) of that result. The existence of enough $!$-able maps is shown in \cref{rslt:!-able-maps-for-D-rel}. As for the categorical Künneth formula $\D_\rel(\Bun_G \times S, \Lambda) = \D_\rel(\Bun_G, \Lambda) \tensor_{\D(\Lambda)} \D_\rel(S, \Lambda)$, it reduces immediately to $\D_\mot$, where it is shown in \cite[Proposition~3.2]{FS-motivic}. 
\end{proof}

Now we bring geometric Satake into the mix. There is a motivic version of the Satake isomorphism in \cite[\S5]{FS-motivic}, but we will not need it. Instead, our following definitions and results rely purely on the Satake isomorphism for torsion coefficients in \cite[\S VI]{FS}. Let us start by defining the Satake category; for now we work with $\Z_\ell[\sqrt q]$-coefficients.

\begin{defn}
Fix a finite flat\footnote{Flatness is needed so that $\Lambda/\ell^n$ is a static ring, as \cite{FS} is only stated for static coefficients.} $\Z_\ell$-algebra $\Lambda$ and a finite set $I$. We define
\begin{align*}
    \Sat^I_G(\Lambda) \subseteq \D_\rel^\bdd(\locHck_G^I, \Lambda)
\end{align*}
to be the full subcategory spanned by those objects $A \in \D_\rel^\bdd(\locHck_G^I, \Lambda)$ that are $\ell$-complete and modulo all $\ell^n$ lie in the Satake category $\Sat^I_G(\Lambda/\ell^n) \subseteq \D_\et^\bdd(\locHck_G^I, \Lambda/\ell^n)$ from \cite[Definition~VI.9.1]{FS}, where we implicitly use the étale comparison of $\D_\rel$ from \cref{rslt:compare-D-rel-to-D-et}.
\end{defn}

\begin{defn}
Fix a $\Z_\ell$-algebra $\Lambda$ and a finite set $I$ and recall the category $\Rep_\Lambda(W_F^I)$ of continuous $W_F^I$-representations on relatively discrete $\Lambda$-modules from \cref{def:smooth-and-continuous-reps-of-cond-group}. The action of $W_F$ on $\hat G$ allows us to view $\calO(\hat G^I)$ as a Hopf-algebra in $\Rep_\Lambda(W_F^I)$ (pulled back from a Hopf algebra in $\Rep_\Lambda(Q^I)$). We denote by
\begin{align*}
    \Rep_\Lambda((\hat G \rtimes W_F)^I)^\fpj \subseteq \Rep_\Lambda((\hat G \rtimes W_F)^I)
\end{align*}
the category of comodules under $\calO(\hat G^I)$ in $\Rep_\Lambda(W_F^I)$ and the full subcategory spanned by those objects where the underlying $\Lambda$-module is finite projective.
\end{defn}

More explicitly, $\Rep_\Lambda((\hat G \rtimes W_F)^I)$ consists of (derived) $\Lambda$-modules equipped with compatible actions of $\hat G^I$ and $W_F^I$, where the first one is algebraic and the second one is continuous for the colimit topology on the module (taking the colimit of the $\ell$-adic topologies on finitely generated  $\Z_\ell$-submodules).

With the above definitions at hand, we can bootstrap the following version of the Satake isomorphism from its torsion analogs:

\begin{thm} \label{rslt:Satake-for-D-rel}
Fix a finite flat $\Z_\ell$-algebra $\Lambda$.
\begin{thmenum}
    \item For every finite set $I$, $\Sat^I_G(\Lambda)$ is an exact ordinary category. It is stable under the monoidal structure on $\D_\rel^\bdd(\locHck_G^I, \Lambda)$, and the induced monoidal structure on $\Sat_I^G(\Lambda)$ upgrades to a symmetric monoidal one.

    \item For every map $I \to J$ of finite sets, the transition functor $\D_\rel^\bdd(\locHck_G^I, \Lambda) \to \D_\rel^\bdd(\locHck_G^J, \Lambda)$ in \cref{rslt:functoriality-of-D-rel-on-locHck} sends $\Sat_G^I(\Lambda)$ to $\Sat_G^J(\Lambda)$ and upgrades to a symmetric monoidal functor.

    \item Suppose that $\Lambda$ contains $\sqrt q$. Then there is a symmetric monoidal equivalence of exact categories
    \begin{align*}
        \Sat_G^I(\Lambda) \isom \Rep_\Lambda((\hat G \rtimes W_F)^I)^\fpj.
    \end{align*}
    It is compatible with the symmetric monoidal transition functors along maps $I \to J$ of finite sets.
\end{thmenum}
\end{thm}
\begin{proof}
The base-changes modulo $\ell^n$ induce a functor
\begin{align*}
    \Sat_G^I(\Lambda) \to \varprojlim_n \Sat_G^I(\Lambda/\ell^n)
\end{align*}
with right adjoint $(A_n)_n \mapsto \varprojlim_n A_n$. Since the objects of $\Sat_G^I(\Lambda)$ are $\ell$-complete, the above functor is an equivalence. Since we understand the case of torsion coefficients (see e.g. \cite[Theorem~VI.0.2]{FS}) this proves everything in (i) except the stability of $\Sat_G^I(\Lambda)$ under convolution, which reduces to showing that for $A, A' \in \Sat_G^I(\Lambda)$ the object $A * A'$ is $\ell$-complete. For this we observe that both $A$ and $A'$ are suave over $\Div^I\!\!/L^+_I G$ by \cref{rslt:characterize-suave-sheaves-on-locHck-via-ell-completeness} (this is a surprisingly non-trivial result, because the pullback and tensor product functors on $\D_\rel$ do not preserve $\ell$-completeness). By the stability properties of suave objects (see \cite[Lemma~4.4.8(i)]{heyer-mann-6ff} and \cite[Lemma~4.5.16(ii)]{heyer-mann-6ff}), this implies that $A * A'$ is suave over $\Div^I \!\! /L^+_I G$, as the relevant pushforward is proper (after restriction to suitably large Schubert cells). But then $A * A'$ is $\ell$-complete by \cref{rslt:characterize-suave-sheaves-on-locHck-via-ell-completeness}. This finishes the proof of (i).

In (ii), the symmetric monoidal upgrade reduces by (i) to the torsion case, where it is discussed in \cite[\S VI.9]{FS}. It only remains to prove that the Satake category is stable under the transition maps associated to a map $I \to J$ of finite sets. Here we argue as in the case of the convolution in (i): It is enough to show that the transition maps preserve suaveness, which follows easily from the explicit description in \cref{rslt:explicit-description-of-transition-functors-on-locHck}.

Part (iii) is an immediate consequence of (i) and (ii) together with the torsion Satake isomorphism in \cite[Theorem~VI.0.2]{FS}.
\end{proof}

\begin{rmk}
In \cite[\S5]{FS-motivic} Scholze defines a motivic version of the Satake category $\Sat_{G}^{I,\mot} \subseteq \D_{\mot}^{\bdd}(\locHck_G^I)$ and proves a version of the Satake isomorphism for this category. By unravelling the definitions (and making some choices as in \cite[\S4]{FS-motivic} that ultimately get cancelled out) this should lead to a different proof of \cref{rslt:Satake-for-D-rel}. Our approach bypasses the motivic Satake isomorphism and avoids all choices.
\end{rmk}

Recall that we fixed a finite quotient $Q$ of $W_F$ through which the action of $\hat G$ factors. We now want to combine the above results in order to obtain natural $\Rep_\Lambda(Q^I)$-linear maps $\Rep_\Lambda((\hat G \rtimes Q)^I) \to \End(\D_\rel(\Bun_G, \Lambda))^{BW_F^I}$ that we can plug into \cite[Theorem~X.0.1]{FS} to construct the spectral action. Here $(-)^{BW_F^I}$ makes use of a condensed structure on $\End(\D_\rel(\Bun_G, \Lambda))$, which we want to avoid. We therefore formulate our result slightly differently:

\begin{prop} \label{rslt:functorial-Hecke-action-for-D-rel}
Fix a $\Z_\ell[\sqrt q]$-algebra $\Lambda$. Then there is a natural $\Rep_\Lambda(Q^I)$-linear monoidal functor
\begin{align*}
    \Rep_\Lambda((\hat G \rtimes Q)^I)^\dbl \to \End_{\D(\Lambda)}(\D(\Bun_G, \Lambda)) \tensor_{\D(\Lambda)} \Rep_\Lambda(W_F^I), 
\end{align*}
functorial in $I$. Here $(-)^\dbl$ refers to the full subcategory of dualizable objects (i.e. those representations where the underlying $\Lambda$-module is perfect).
\end{prop}
\begin{proof}
First consider the case $\Lambda = \Z_\ell[\sqrt q]$. We take the functor in \cref{rslt:functoriality-of-D-rel-on-locHck}, first forget the $\D_\rel(\Div^I, \Z_\ell[\sqrt q])$-linear structure to a $\Rep_{\Z[\sqrt q]}(Q^I)^\fpj$-linear structure via the pullback along $\Div^I \to */Q^I$ and then replace the source $\D_\rel^\bdd(\locHck_G^I, \Z_\ell[\sqrt q])$ by the full monoidal subcategory $\Sat_G^I(\Z_\ell[\sqrt q])$. Using \cref{rslt:Satake-for-D-rel} we can functorially in $I$ replace this category by $\Rep_{\Z_\ell[\sqrt q]}((\hat G \rtimes W_F)^I)^\fpj$ and then by $\Rep_{\Z_\ell[\sqrt q]}((\hat G \rtimes Q)^I)^\fpj$ --- here we use that these categories are ordinary monoidal categories and hence there is no higher functoriality to check. This gives us natural $\Rep_{\Z_\ell[\sqrt q]}(Q^I)^\fpj$-linear monoidal functors
\begin{align*}
    \Rep_{\Z_\ell[\sqrt q]}((\hat G \rtimes Q)^I)^\fpj \to \End_{\D(\Z_\ell[\sqrt q])}(\D(\Bun_G, \Z_\ell[\sqrt q])) \tensor_{\D(\Z_\ell[\sqrt q])} \D_\rel(\Div^I, \Z_\ell[\sqrt q]),
\end{align*}
functorial in $I$. As explained in \cite[\S IX.2]{FS}, it follows from highest weight theory that for a general $\Z_\ell[\sqrt q]$-algebra $\Lambda$ the natural map
\begin{align*}
    \Rep_{\Z_\ell[\sqrt q]}((\hat G \rtimes Q)^I)^\fpj \tensor_{\Rep_{\Z_\ell[\sqrt q]}(Q^I)^\fpj} \Rep_\Lambda(Q^I)^\dbl \isoto \Rep_\Lambda((\hat G \rtimes Q)^I)^\dbl
\end{align*}
is an isomorphism; here the tensor product is to be interpreted as a tensor product in $\PrL$ via passing to $\Ind$-categories and a posteriori passing to compact objects (cf. \cref{rslt:compact-generation-of-relative-tensor} for why the tensor product is automatically compactly generated). Using this isomorphism we can upgrade our functors to $\Rep_\Lambda(Q^I)^\dbl$-linear monoidal functors
\begin{align*}
    \Rep_\Lambda((\hat G \rtimes Q)^I)^\dbl \to \End_{\D(\Lambda)}(\D(\Bun_G, \Lambda)) \tensor_{\D(\Lambda)} \D_\rel(\Div^I, \Lambda),
\end{align*}
functorial in $I$. By \cref{rslt:D-rel-Div-tensor-I-fully-faithful} the natural symmetric monoidal functor $\D_\rel(\Div^1, \Lambda)^{\tensor I} \injto \D_\rel(\Div^I, \Lambda)$ is fully faithful and by \cref{rslt:dualizable-ff-trick} the same is true after tensoring with $\End_{\D(\Lambda)}(\D(\Bun_G, \Lambda))$ because the latter category is compactly generated (e.g. since it identifies with $\D(\Bun_{G\times G}, \Lambda)$). We claim that the above functor factors through this full subcategory: Indeed, this follows immediately from the fact that every object in $\Rep_\Lambda((\hat G \rtimes Q)^I)^\dbl$ admits a resolution by exterior tensor products of objects in $\Rep_\Lambda(\hat G \rtimes Q)^\dbl$ (see the proof of \cite[Corollary~IX.2.3]{FS}) which by functoriality in $I$ reduces us to the case $\abs I = 1$, which is trivial. We thus arrived at natural $\Rep_\Lambda(Q^I)^\dbl$-linear functors
\begin{align*}
    \Rep_\Lambda((\hat G \rtimes Q)^I)^\dbl \to \End_{\D(\Lambda)}(\D(\Bun_G, \Lambda)) \tensor_{\D(\Lambda)} \D_\rel(\Div^1, \Lambda)^{\tensor I},
\end{align*}
functorial in $I$. We conclude by composing the above functors with the natural symmetric monoidal functor
\begin{align*}
    \D_\rel(\Div^1, \Lambda)^{\tensor I} \to \Rep_\Lambda(W_F)^{\tensor I} \to \Rep_\Lambda(W_F^I),
\end{align*}
where the first functor is the one from \cref{rslt:compare-D-rel-Div1-to-Rep-W-F}.
\end{proof}

\begin{rmks}
\begin{rmksenum}
    \item Amusingly, it is \emph{not} clear from this construction that the Hecke action specified in \cref{rslt:functorial-Hecke-action-for-D-rel} coincides with the Hecke action constructed in \cite{FS}, because the relation of the 6-functor formalism $\D_\rel$ and the 5-functor formalism $\D_\solid$ from \cite{FS} is unclear. The agreement of the two Hecke actions will be proved in \cite{GH}, using nuclear sheaves as an intermediary.

    \item Even though the construction of the Hecke action in \cref{rslt:functorial-Hecke-action-for-D-rel} passes through the \emph{ordinary} categories $\Rep_{\Z_\ell[\sqrt q]}((\hat G \rtimes Q)^I)^\fpj$, the claimed functoriality in $I$ is an honest $\infty$-categorical claim, as the target categories $\End_{\D(\Lambda)}(\D(\Bun_G, \Lambda)) \tensor \Rep_\Lambda(W_F^I)$ are inherently $\infty$-categorical. For example, one needs to fill higher homotopies for all finite chains of compositions of finite sets $I$.
\end{rmksenum}
\end{rmks}

With \cref{rslt:functorial-Hecke-action-for-D-rel} at hand, we can now proceed as in \cite[\S X]{FS}, in fact we can reduce to the abstract results proved there. For the following result, recall the definition of the L-parameter stack $\Par_G$ from \cref{def:ParG}.

\begin{thm} \label{rslt:spectral-action-for-D-rel}
Fix a $\Z_\ell[\sqrt q]$-algebra $\Lambda$ and assume either that $\ell$ is invertible in $\Lambda$ or $\ell \nmid \abs{\pi_0(Z(G))}$. Then the functors in \cref{rslt:functorial-Hecke-action-for-D-rel} induce an action of $\Perf(\Par_{G,\Lambda})$ on $\D(\Bun_G, \Lambda)$.
\end{thm}
\begin{proof}
We want to apply \cite[Theorem~X.0.1]{FS} to \cref{rslt:functorial-Hecke-action-for-D-rel}, but it is unclear to us how to get $\End(\D(\Bun_G))^{W_F^I}$ into the picture, as this requires some condensed structure on $\End(\D(\Bun_G))$. Instead, we pursue the same approach as in \cite[\S X]{FS} and reduce the construction to discrete subquotients of $W_F^I$.

By (the proof of) \cref{rslt:functorial-Hecke-action-for-D-rel} we have natural $\Rep(Q^I)^\dbl$-linear monoidal functors
\begin{align*}
    \Rep((\hat G \rtimes Q)^I)^\dbl \to \End_{\D(\Lambda)}(\D(\Bun_G, \Lambda)) \tensor_{\D(\Lambda)} \D_\rel(\Div^1, \Lambda)^{\tensor I},
\end{align*}
functorial in $I$ (in the following we need this more refined statement instead of the one in \cref{rslt:functorial-Hecke-action-for-D-rel}). By \cref{rslt:D-rel-Div1-maps-to-colim-W-F-mod-P} there is a natural symmetric monoidal functor $\D_\rel(\Div^1, \Lambda) \to \varinjlim_P \Rep_\Lambda(W_F/P)$, where $P$ runs through the open subgroups of the wild inertia group $P_F \subset W_F$. We obtain a natural symmetric monoidal functor
\begin{align*}
    \D_\rel(\Div^1, \Lambda)^{\tensor I} \to \varinjlim_P \Rep_\Lambda(W_F/P)^{\tensor I} \to \varinjlim_P \Rep_\Lambda((W_F/P)^I).
\end{align*}
We fix a discretization $W_{F,\iota}$ of $W_F$ as in \cref{def:discretization-of-W-F}. Note that for every anima $S$ and every dualizable $\Lambda$-linear category $\cat C$ there is a natural equivalence $\cat C \tensor_{\D(\Lambda)} \Fun(S, \D(\Lambda)) = \cat C^S$; indeed, since $\cat C$ is dualizable over $\D(\Lambda)$, both sides send colimits in $S$ to limits and the statement is clearly true for the generator $S = *$. Thus, by composing the above functors with the natural symmetric monoidal functors
\begin{align*}
    \Rep_\Lambda((W_F/P)^I) \to \Rep_\Lambda((W_{F,\iota}/P)^I) = \Fun(B(W_{F,\iota}/P)^I, \D(\Lambda)),
\end{align*}
we arrive at the $\Rep_{\Lambda}(Q^I)^\dbl$-linear monoidal functors
\begin{align*}
    \Rep_\Lambda((\hat G \rtimes Q)^I)^\dbl \to \varinjlim_P \End_{\D(\Lambda)}(\D(\Bun_G, \Lambda))^{B(W_{F,\iota}/P)^I},
\end{align*}
functorially in $I$. Note that the colimit over $P$ is taken in $\PrL$, i.e. by \cite[Corollary~A.5.9]{heyer-mann-6ff} it is the $\Ind$-category of the colimit (in $\Cat$) on compact objects. For fixed open $P_0 \subseteq P_F$ in the kernel of $W_F \surjto Q$, let $\End_{\D(\Lambda)}^{P_0}(\D(\Bun_G, \Lambda)) \subseteq \End_{\D(\Lambda)}(\D(\Bun_G, \Lambda))$ be the full monoidal subcategory of those functors that preserve the full subcategory $\D^{P_0}(\Bun_G, \Lambda)^\omega$ from \cite[\S IX.5]{FS} (defined in terms of the Hecke action we constructed above). Then by definition the Hecke functors factor through $\End_{\D(\Lambda)}^P(\D(\Bun_G, \Lambda))$. By composing with the canonical projection to $\End_{\D(\Lambda)}(\D^{P_0}(\Bun_G, \Lambda))$ we obtain $\Rep_\Lambda(Q^I)^\dbl$-linear monoidal functors
\begin{align*}
    \Rep_\Lambda((\hat G \rtimes Q)^I)^\dbl \to \varinjlim_{P \supseteq P_0} \End_{\D(\Lambda)}(\D(\Bun_G,\Lambda)^P)^{B(W_{F,\iota}/P)^I}.
\end{align*}
Note that the transition functors in the colimit diagram are fully faithful and preserve compact objects (as this is true for the corresponding functors on representation categories, it follows after tensoring with the endofunctor category by \cref{rslt:dualizable-ff-trick}). Therefore each of the terms in the colimit embeds fully faithfully into the colimit. By definition of $\D^{P_0}(\Bun_G, \Lambda)$ the Hecke action factors through the term for $P = P_0$, i.e. we get $\Rep_\Lambda(Q^I)$-linear monoidal functors
\begin{align*}
    \Rep_\Lambda((\hat G \rtimes Q)^I)^\dbl \to \End_{\D(\Lambda)}(\D(\Bun_G,\Lambda)^{P_0})^{B(W_{F,\iota}/P_0)^I},
\end{align*}
functorially in $I$. Plugging this into \cite[Theorem~X.0.2]{FS} (and tensoring with $\D(\Lambda)$) produces an action of $\Perf(\intRep_{\hat G \rtimes Q}(W_{F,\iota}/P_0) \times_{\intRep_Q(W_{F,\iota}/P_0)} *)$ on $\D(\Bun_G, \Lambda)^{P_0}$.

We sketch how to glue the above actions along different $P_0$'s in order to arrive at the desired action of $\Par_G$ (this is implicitly claimed in the proof of \cite[Corollary~X.1.3]{FS}). For fixed $P$, denote by $\End'^I_P \subseteq \End_{\D(\Lambda)}(\D(\Bun_G, \Lambda)^P)^{B(W_{F,\iota}/P)^I}$ the full subcategory spanned by those endofunctors of $\D(\Bun_G, \Lambda)^P$ that preserve $\D(\Bun_G, \Lambda)^{P'}$ for all $P' \supseteq P$ and whose restriction to $\D(\Bun_G, \Lambda)^{P'}$ factors over the full subcategory
\begin{align*}
    \End_{\D(\Lambda)}(\D(\Bun_G, \Lambda)^{P'})^{B(W_{F,\iota}/P')^I} \subseteq \End_{\D(\Lambda)}(\D(\Bun_G, \Lambda)^{P'})^{B(W_{F,\iota}/P)^I}.
\end{align*}
By definition of $\D(\Bun_G, \Lambda)^P$, the Hecke action factors through $\End'^I_P$. But by construction, $\End'^I_P$ is functorial in $P$, i.e. for $P \subseteq P'$ we have a natural restriction functor $\End'^I_P \to \End'^I_{P'}$. Taking the limit over the induced actions of representation stacks on $\D(\Bun_G, \Lambda)^P$ and using \cref{rslt:decompose-Rep-W-F-as-colimit}, we arrive at the desired action of $\Perf(\Par_G)$ on $\D(\Bun_G, \Lambda)$.
\end{proof}

\begin{rmk}
The construction in \cref{rslt:spectral-action-for-D-rel} is quite ad-hoc as it passes to discrete subquotients of $W_F$. We will pursue a more conceptual construction using more advanced categorical techniques in future work.
\end{rmk}

\begin{defn}
In the case $\Lambda = \Qellbar$ we denote $\Par_G := \Par_{G,\Lambda}$ and $\D(\Bun_G) := \D(\Bun_G, \Lambda)$ as usual. Since $\QCoh(\Par_G) = \Ind(\Perf^\qc(\Par_G))$ by \cref{rslt:QCoh-and-ICoh-on-ParG-compactly-generated}, we obtain from \cref{rslt:spectral-action-for-D-rel} an action
\begin{align*}
    \QCoh(\Par_G) \times \D(\Bun_G) \to \D(\Bun_G), \qquad (\mathcal F, A) \mapsto \mathcal F *  A.
\end{align*}
This is called the \emph{spectral action}.
\end{defn}

With the spectral action at hand, we can now define the spectral support of a sheaf $A \in \D(\Bun_G)$. This is slightly subtle as the scheme $X_G^\spec$ defined in \cref{def:coarse-moduli-space} is not affine or quasicompact. The correct definition turns out to be the following.

\begin{defn} \label{def:spectral-support}
Fix $A \in \D(\Bun_G)$. For each connected component $X_i \subseteq X_G^\spec$ with corresponding connected component $X'_i \subseteq \Par_G$ and idempotent $e_i \in \mathfrak{Z}_{G}^{\spec}$, denote $e_i A := \calO_{X'_i} * A$. Then there is a natural ring map
\begin{align*}
    e_i \mathfrak{Z}_{G}^{\spec} = \End(\calO_{X'_i}) \to \End(e_i A).
\end{align*}
Writing $C_i \subset X_i$ for the reduced closed subscheme cut out by the radical of the kernel of this map, we define the \emph{spectral support} of $A$ as $\mathrm{spec.supp}(A) =\cup_i C_i$. Note that if $A$ is compact then $e_i A = 0$ for all but finitely many $i$.
\end{defn}

\subsection{Duality for the spectral action}

In the previous subsection we constructed the spectral action of $\QCoh(\Par_G)$ on $\D(\Bun_G)$. In the following we will show that this action satisfies a certain compatibility with respect to duality on $\Bun_G$. A version of this duality is shown in \cite[Theorem~IX.0.1(i)]{FS} for the Hecke actions. Using the functoriality of Hecke actions shown in the previous subsection and the self-duality of convolution stacks, we recover \cite[Theorem~IX.0.1(i)]{FS} and glue it to a similar duality statement for the spectral action.

Without further ado, let us first state the compatibility of Hecke actions with duality. Recall that for a monoidal category $\cat E$ we denote by $\cat E^\rev$ the same category but with opposite monoidal structure, i.e. $X \tensor_{\cat E^\rev} Y = Y \tensor_{\cat E} X$. Recall also that the Chevalley involution $c_G$ on $\hat{G} \rtimes W_F$ induces an involution on $\Par_{G,\Lambda}$, which we also denote by $c_G$ or simply by $c$ if $G$ is clear from context.

\begin{prop} \label{rslt:self-duality-of-Hecke-action}
Fix a $\Z_\ell[\sqrt{q}]$-algebra $\Lambda$. Then there is a commuting square
\begin{equation*}\begin{tikzcd}
    \Rep_\Lambda((\hat G \rtimes Q)^I)^\dbl \arrow[r] \arrow[d,isomb,"c^*"] & \End_{\D(\Lambda)}(\D(\Bun_G, \Lambda))^\rev \tensor_{\D(\Lambda)} \Rep_\Lambda(W_F^I) \arrow[d,isom] \\
    \Rep_\Lambda((\hat G \rtimes Q)^I)^\dbl \arrow[r] & \End_{\D(\Lambda)}(\D(\Bun_G, \Lambda)) \tensor_{\D(\Lambda)} \Rep_\Lambda(W_F^I)
\end{tikzcd}\end{equation*}
of $\Rep_\Lambda(Q^I)$-linear monoidal functors, natural in the finite set $I$. These functors have the following explicit description:
\begin{propenum}
    \item The horizontal maps come from the Hecke action in \cref{rslt:functorial-Hecke-action-for-D-rel} (apply $(-)^\rev$ for the top map).
    
    \item The left vertical map is induced by pullback along the Chevalley involution.
    
    \item \label{rslt:self-duality-of-Hecke-action-map-on-End-Bun-G} The right vertical map is induced by the isomorphism
    \begin{align*}
        \End_{\D(\Lambda)}(\D(\Bun_G,\Lambda))^\rev \isoto \End_{\D(\Lambda)}(\D(\Bun_G,\Lambda)^\vee) \isoto \End_{\D(\Lambda)}(\D(\Bun_G,\Lambda)),
    \end{align*}
    where the second map comes from the canonical isomorphism $\D(\Bun_G)^\vee \isoto \D(\Bun_G)$ (cf. \cref{rslt:duality-of-D-bdd-Hck-map-on-End-Bun-G}).
\end{propenum}
\end{prop}
\begin{proof}
By \cref{rslt:functoriality-of-D-rel-on-locHck} and \cref{rslt:duality-of-D-bdd-Hck} we obtain a commuting square
\begin{equation*}\begin{tikzcd}
    \D_\rel^\bdd(\locHck_G^I, \Lambda)^\rev \arrow[r] \arrow[d,isom] & \End_{\D(\Lambda)}(\D(\Bun_G,\Lambda))^\rev \otimes_{\D(\Lambda)} \D_\rel(\Div^I, \Lambda) \arrow[d,isom] \\
    \D_\rel^\bdd(\locHck_G^I, \Lambda) \arrow[r] & \End_{\D(\Lambda)}(\D(\Bun_G,\Lambda)) \otimes_{\D(\Lambda)} \D_\rel(\Div^I, \Lambda)
\end{tikzcd}\end{equation*}
of $\D_\rel(\Div^I, \Lambda)$-linear monoidal functors, natural in $I$. Here the left vertical map is induced by swapping the two $L^+_I G$-quotients in the definition of $\locHck_G^I$ and the right vertical map is induced by the isomorphism in (iii). We now take this square and go through the construction in \cref{rslt:functorial-Hecke-action-for-D-rel}. The construction works in the same way and produces the desired square if we can show that the restriction of the swapping isomorphism $\D^\bdd(\locHck^I_G, \Lambda)^\rev \isoto \D^\bdd(\locHck^I_G)$ preserves the full subcategory $\Sat^I_G(\Lambda)$. But this is shown in \cite[\S VI.12]{FS}, which also proves that the induced isomorphism $\Sat^I_G(\Lambda) \isoto \Sat^I_G(\Lambda)$ is given by Chevalley involution. Note that this identification of the isomorphism can be made functorially in $I$, as $\Sat_G^I(\Lambda)$ is an ordinary category.
\end{proof}

\begin{rmk}
The self-duality $\D(\Bun_G,\Lambda)^\vee \isoto \D(\Bun_G,\Lambda)$ in \cref{rslt:self-duality-of-Hecke-action-map-on-End-Bun-G} comes from the self-duality of $\Bun_G$ in the category of kernels $\cat K_{\D_\rel(-,\Lambda)}$, as shown in the proof of \cref{rslt:duality-of-D-bdd-Hck-map-on-End-Bun-G}. In the case of $\D_\rel$ it can be made explicit: We know that $\D(\Bun_G,\Lambda)$ is compactly generated by \cref{rslt:D-Bun-G-is-compactly-generated}, hence
\begin{align*}
    \D(\Bun_G,\Lambda) = \Ind(\D(\Bun_G,\Lambda)^\omega), \qquad \D(\Bun_G,\Lambda) = \Ind(\D(\Bun_G,\Lambda)^{\omega,\op}),
\end{align*}
so the desired self-duality reduces to an isomorphism $\D(\Bun_G,\Lambda)^{\omega,\op} \isoto \D(\Bun_G,\Lambda)^{\omega}$. But $\D(\Bun_G,\Lambda)^\omega$ identifies with the category of prim objects in $\D(\Bun_G,\Lambda)$, hence prim duality (i.e., Bernstein-Zelevinsky duality) provides the desired isomorphism. To check that this is the correct isomorphism, one can use the fact that the evaluation map of the self-duality is given by
\begin{align*}
    \D(\Bun_G,\Lambda) \tensor_{\D(\Lambda)} \D(\Bun_G,\Lambda) \to \D(\Lambda), \qquad (A_1, A_2) \mapsto \Gamma_c(\Bun_G , A_1 \tensor A_2),
\end{align*}
as this is how the self-duality is given in the category of kernels.
\end{rmk}

With \cref{rslt:self-duality-of-Hecke-action} at hand, we can now carry the isomorphism therein through the construction of the spectral action and arrive at the following result. We formulate it in terms of 2-categories, which gives us direct access to suave and prim duality below. We refer the reader to \cref{eq:self-duality-of-spectral-action-via-monoidal-functors} and the discussion below for a more down-to-earth formulation of the same result. In the following result, for a monoidal category $\cat V$ we denote by $*_{\cat V}$ the 2-category with a single object and endomorphisms $\cat V$ (cf. \cite[Example~C.1.14]{heyer-mann-6ff}). We also make use of the 2-category of kernels $\cat K_{\D}$ from \cite[Definition~4.1.3(a)]{heyer-mann-6ff}.

\begin{thm} \label{rslt:2-cat-duality-for-spectral-action}
Fix a $\Z_\ell[\sqrt q]$-algebra $\Lambda$ and assume either that $\ell$ is invertible in $\Lambda$ or $\ell \nmid \abs{\pi_0(Z(G))}$. Write $\cat K_{\cat D} := \cat K_{\D_\rel(-,\Lambda)}$. Then there is a commuting square of 2-categories
\begin{equation*}\begin{tikzcd}
    *_{\Perf(\Par_{G,\Lambda})} \arrow[r] \arrow[d,"c^*",swap] & \cat K^\op_{\D} \arrow[d]\\
    *_{\Perf(\Par_{G,\Lambda})} \arrow[r] & \cat K_{\D}
\end{tikzcd}\end{equation*}
Here the horizontal maps send $*$ to $\Bun_G$ and induce the spectral action on endomorphisms, the left vertical map is induced by the Chevalley involution on $\hat G \rtimes Q$ and the right vertical map is the equivalence from \cite[Proposition~4.1.4]{heyer-mann-6ff}.
\end{thm}
\begin{proof}
As explained in \cite[Example~C.1.14]{heyer-mann-6ff}, the bottom horizontal map corresponds to a monoidal functor
\begin{align*}
    \Perf(\Par_{G,\Lambda}) \to \End_{\cat K_{\D}}(\Bun_G) = \End_{\D(\Lambda)}(\D(\Bun_G,\Lambda)),
\end{align*}
where the equality comes from the fully faithful functor $\Psi_{\D}\colon \cat K \injto \Mod_{\D(\Lambda)}(\PrL)$ in the proof of \cref{rslt:full-functoriality-of-D-bdd-Hck}. Similarly, the upper horizontal map corresponds to a monoidal functor $\Perf(\Par_{G,\Lambda}) \to \End_{\D(\Lambda)}(\D(\Bun_G,\Lambda))^\rev$. Altogether we see that the claimed commuting diagram of 2-categories is equivalent to a commuting diagram of monoidal functors
\begin{equation}\begin{tikzcd}
    \Perf(\Par_{G,\Lambda}) \arrow[r] \arrow[d,isom,"c^*",swap] & \End_{\D(\Lambda)}(\D(\Bun_G,\Lambda))^\rev \arrow[d,isom] \\ 
    \Perf(\Par_{G,\Lambda}) \arrow[r] & \End_{\D(\Lambda)}(\D(\Bun_G,\Lambda))
    \label{eq:self-duality-of-spectral-action-via-monoidal-functors}
\end{tikzcd}\end{equation}
Here the right-hand vertical map is the one from \cref{rslt:self-duality-of-Hecke-action-map-on-End-Bun-G} (namely, the isomorphism $\cat K_{\D}^\op \isoto \cat K_{\D}$ is induced by passing to duals, and the proof of \cref{rslt:duality-of-D-bdd-Hck-map-on-End-Bun-G} shows that this induces the claimed isomorphism on $\End$). The left-hand map vertical map is induced by Chevalley involution and the horizontal maps are induced by the spectral action in \cref{rslt:spectral-action-for-D-rel} (resp. the $(-)^\rev$-version of it).

The square \cref{eq:self-duality-of-spectral-action-via-monoidal-functors} is obtained by passing the square in \cref{rslt:self-duality-of-Hecke-action} through the construction of the spectral action in \cref{rslt:spectral-action-for-D-rel}. There is a little complication as we pass through a limit of discrete subquotients of $W_F^I$ in the proof of \cref{rslt:spectral-action-for-D-rel}. However, this problem can be overcome by noticing that for each fixed $P \subseteq W_F$ as in the proof of \cref{rslt:spectral-action-for-D-rel}, we get a similar commuting diagram with $\Perf(\Par_{G,\Lambda})$ replaced by perfect sheaves on the corresponding qcqs clopen subset and $\D(\Bun_G,\Lambda)$ replaced by $\D(\Bun_G,\Lambda)^P$; observe that $\D(\Bun_G,\Lambda)^P$ is still self-dual (via the same duality map), as the compact objects in $\D(\Bun_G,\Lambda)^P$ are stable under prim duality.
\end{proof}

While \cref{rslt:2-cat-duality-for-spectral-action} may be hard to grasp by itself, we obtain the following immediate consequence by looking at adjoint morphisms in the category of kernels. This is the main duality theorem we will actually apply later in the paper.

\begin{thm}\label{thm:spectraldualityuseful}
Fix a $\Z_\ell[\sqrt q]$-algebra $\Lambda$ and assume either that $\ell$ is invertible in $\Lambda$ or $\ell \nmid \abs{\pi_0(Z(G))}$. Denote the spectral action from \cref{rslt:spectral-action-for-D-rel} by $*$ and fix some $\mathcal F \in \Perf(\Par_{G,\Lambda})$. Then $\mathcal F * -$ preserves compactness and ULAness. Moreover:
\begin{thmenum}
    \item \label{rslt:prim-duality-for-spectral-action} For all $A \in \D(\Bun_G,\Lambda)^\omega$ and $\mathcal{F} \in \Perf(\Par_{G,\Lambda})$, we have a natural isomorphism
    \begin{align*}
        \Dbz(\mathcal{F} \ast A) = c^* \mathcal{F}^\vee \ast \Dbz A
    \end{align*}
    in $\D(\Bun_G,\Lambda)^{\omega}$.

    \item For all $A \in \D(\Bun_G,\Lambda)$ and $\mathcal{F} \in \Perf(\Par_{G,\Lambda})$, we have a natural isomorphism
    \begin{align*}
        \Dverd(\mathcal{F} \ast A) = c^* \mathcal F^\vee \ast \Dverd A
    \end{align*}
    in $\D(\Bun_G,\Lambda)$.
\end{thmenum}
\end{thm}

Starting in the next chapter, we will need to consider duals of coherent (not just perfect) sheaves. However, since $\Par_{G}$ has trivial dualizing sheaf, the functor $\Dgs$ restricts to the naive duality $(-)^\vee$ along the evident embedding $\Perf(\Par_G) \subset \Coh(\Par_G)$, and similarly with the Chevalley twist. As such, we will systematically write $\Dtwgs \mathcal{F}$ instead of $c^{\ast} \mathcal{F}^{\vee}$ beginning in the next chapter.

\begin{proof}
Fix $\mathcal F \in \Perf(\Par_{G,\Lambda})$. Via the horizontal maps in \cref{rslt:2-cat-duality-for-spectral-action}, $\mathcal F$ induces a map $\mathcal F_*\colon \Bun_G \to \Bun_G$ in $\cat K_{\D}$. Since $\mathcal F$ is dualizable, it is both left and right adjoint in $*_{\Perf(\Par_{G,\Lambda})}$. Hence the same is true for $\mathcal F_*$ and both its left and right adjoint are given by $(\mathcal F^\vee)_*$.

Given $A \in \D(\Bun_G,\Lambda)^\omega$, we view it as a left adjoint map $* \to \Bun_G$ in $\cat K_{\D}$. Then the left-hand side in (i) evaluates to the right adjoint of the composition $\mathcal F_* \comp A\colon * \to \Bun_G$. This left adjoint is thus given by $\Dbz(A) \comp (\mathcal F^\vee)_*\colon \Bun_G \to *$. We now use the equivalence $\cat K_{\D}^\op = \cat K_{\D}$ to reinterpret this as a map $* \to \Bun_G$; by the square in \cref{rslt:2-cat-duality-for-spectral-action} this procedure replaces $(\mathcal F^\vee)_*$ by $(c^* \mathcal F^\vee)_*$, so we arrive at the map $(c^{\ast} \mathcal F^\vee)_* \comp \Dbz(A)\colon * \to \Bun_G$. This is the right-hand side in (i).

Part (ii) can be proved similarly by using right adjoints instead of left adjoints, but this will only recover the claimed identity in case that $A$ is suave (i.e. ULA), which is not sufficient for our purposes. In order to get the result for all $A$, we instead argue as follows. Fix any compact $B$. Then we compute
\begin{align*}
    \RHom(B,\Dverd(\mathcal{F}\ast A)) & = \RHom(\Dbz B, \mathcal{F} \ast A)^\vee \\
     & = \RHom(\mathcal{F}^\vee \ast \Dbz B,  A)^\vee \\
     & = \RHom(\Dbz (c^{\ast}\mathcal{F} \ast B), A)^\vee \\
     & = \RHom(c^{\ast}\mathcal{F} \ast B, \Dverd A) \\
     & = \RHom(B, (c^{\ast}\mathcal{F})^\vee \ast \Dverd A).
\end{align*}
Here we used the duality exchange formula (\cref{rslt:dualityexchangeBunG}) in the first and fourth lines; the adjunction between $\mathcal{G}^\vee \ast -$ and $\mathcal{G} \ast -$ in the second and fifth lines; and part (i) in the third line. Passing to colimits in $B$, we get the same identity for any (not necessarily compact) $B$, so the result follows by Yoneda.
\end{proof}

\begin{rmk}\label{rmk:dualityGammaform}
There is a different way to prove \cref{rslt:prim-duality-for-spectral-action} that avoids $\cat K_{\D}$. Namely, we first show that for all $A, B \in \D(\Bun_G,\Lambda)$ and $\mathcal F \in \Perf(\Par_{G,\Lambda})$ there is a natural isomorphism
\begin{align*}
   \Gamma_c(\Bun_G, (\mathcal F * A) \tensor B) = \Gamma_c(\Bun_G, A \tensor (c^* \mathcal F * B)).
\end{align*}
Indeed, this identity can be obtained from \cref{eq:self-duality-of-spectral-action-via-monoidal-functors} by composing with the functor
\begin{align*}
    \D(\Bun_G, \Lambda) \times \D(\Bun_G,\Lambda) \times \End_{\D(\Lambda)}(\D(\Bun_G,\Lambda)) &\to \D(\Lambda),\\
    (A, B, f) &\mapsto \Gamma_c(\Bun_G, f(A) \tensor B),
\end{align*}
we leave the details to the reader. With this identity at hand we can now prove \cref{rslt:prim-duality-for-spectral-action} by a direct computation. More precisely, fix any sheaf $B$. Then we compute
\begin{align*}
    \RHom(\Dbz(\mathcal{F} \ast A),B) & = \Gamma_c(\Bun_G,\mathcal{F} \ast A \otimes B) \\
    &= \Gamma_c(\Bun_G,A\otimes c^{\ast}\mathcal{F} \ast B) \\
    &= \RHom(\Dbz A, c^{\ast}\mathcal{F} \ast B) \\
    &= \RHom((c^\ast \mathcal{F})^\vee \ast \Dbz A,B)
\end{align*}
Here we used the defining relation of Bernstein-Zelevinsky duality in the first and third lines; the previous identity in the second line; and the adjunction between $\mathcal{G}^\vee \ast -$ and $\mathcal{G} \ast - $ in the fourth line. Since $B$ is arbitrary, the result follows by Yoneda.
\end{rmk}

\subsection{Eisenstein series and non-basic strata} \label{sec:nonbasic}

In this section we recall the precise relationship between the functors $i_{b!}^{\ren}$ and $i_{b \sharp}^{\ren}$ and Eisenstein series. This is just a recitation of results from \cite{HHS}, but we recall it in slightly different notation more suited to our needs.

For the remainder of this section, we assume $G$ is quasisplit, with $A \subset T \subset B \subset G$ as usual, and our coefficient ring is $\Qellbar$ which we drop from the notation following our usual convention. Let $U_0$ be the unipotent radical of $B$. For any standard parabolic $P=MU \subset G$, there is an associated Eisenstein series functor \[\Eis_{P!}:\D(\Bun_M) \to \D(\Bun_G),\] together with a pair of constant term functors \[\CT_{P\ast}, \CT_{P!} : \D(\Bun_G) \to \D(\Bun_M),\]
all of which were introduced and studied at length in \cite{HHS}. These are geometrizations of the usual operations of parabolic induction and Jacquet module. By the main results in \cite{HHS}, $\Eis_{P!}$ preserves compact objects, and on compact sheaves there is a duality isomorphism $\Dbz \Eis_{P!} = \Eis_{P^- !} \mathbf{D}_{\mathrm{BZ}}^{M}$. The latter result is a geometrization of Bernstein's famous second adjointness theorem, and can be equivalently formulated as an isomorphism $\CT_{P \ast} \cong \CT_{P^- !}$. 

Let $\mathcal{P}$ denote the set of pairs $(M,b^M)$ where $M \subset G$ is a standard Levi and $b^M \in B(M)_{\mathrm{bas}}$ is a basic element whose Newton point $\nu_{b^M}$ lies in $X_{\ast}(A)_{\mathbf{Q}}^{+}$. Let $\mathcal{P}^+ \subset \mathcal{P}$ denote the subset where $\mathrm{Cent}_{G}(\nu_{b^M}) = M$, or equivalently where $MU_0$ is exactly the dynamic parabolic of $\nu_{b^M}$.
By a classic result of Kottwitz \cite[\S 6]{KottwitzIso}, there is a canonical identification $\mathcal{P}^+ = B(G)$ induced by sending a pair $(M,b^M)$ to the image of $b^M$ in $B(G)$. Under this bijection, the $\sigma$-centralizer $G_{b}$ identifies canonically with the inner form $M_{b^M}$ of $M$.

The following result is part of \cite[Corollary 2.2.5]{HHS}, up to a small change of notation.

\begin{prop}For any $b \in B(G)$ with corresponding pair $(M,b^M) \in \mathcal{P}^+$, there are canonical identifications
\[ i_{b!}^{\ren} = \Eis_{M U_{0}^- !}i_{b^M!}^{M} \]
and
\[ i_{b \sharp}^{\ren} = \Eis_{M U_{0} !}i_{b^M!}^{M} \]
as functors $\D(G_b(F),\Qellbar) \to \D(\Bun_G)$.
\end{prop}

Next, observe that the inclusion $\mathcal{P}^+ \subset \mathcal{P}$ has a canonical retraction $r$, sending a pair $(L,b^L)$ to the pair $(M,b^M)$ where $M=\mathrm{Cent}_{G}(\nu_{b^L})$ and $b^M$ is the image of $b^L$ under $B(L) \to B(M)$. Given a fixed element $(M,b^M) \in \mathcal{P}^+$, the fiber $r^{-1}(M,b^M)$ is in canonical bijection with standard Levi subgroups $L\subset M$ with the property that $b^M$ lies in the image of $B(L) \to B(M)$: for $L$ such that $b^M$ does lie in this image, it automatically admits a unique preimage $b^L \in B(L)$, and this preimage is automatically basic. Up to some change of notation, this follows from the discussion in \cite[Lemmas 2.2.10 and 2.2.11]{HHS}.

Now suppose given some element $b \in B(G)$ corresponding to a pair $(M,b^M) \in \mathcal{P}^+$, together with an element $(L,b^L) \in r^{-1} (M,b^M)$. Let $Q= (LU_0) \cap M$ be the standard parabolic in $M$ with Levi $L$. Then the $\sigma$-centralizer $Q_{b^L}$ is a standard parabolic in $M_{b^L}=M_{b^M} = G_{b}$ with Levi $L_{b^L}$.

\begin{prop}For any $b \in B(G)$ with corresponding pair $(M,b^M) \in \mathcal{P}^+$, and any $(L,b^L) \in \mathcal{P}$ retracting to it, there are 
canonical identifications
\[ i_{b!}^{\ren} i_{Q_{b^L}^{-}}^{G_b} = \Eis_{L U_{0}^- !}i_{b^L!}^{L} \]
and
\[ i_{b \sharp}^{\ren} i_{Q_{b^L}}^{G_b} = \Eis_{L U_{0} !}i_{b^L!}^{L} \]
as functors $\D(L_{b^L}(F),\Qellbar) \to \D(\Bun_G)$.
\end{prop}
\begin{proof}
This follows from \cite[Proposition 2.2.12]{HHS}, up to small changes in notation.
\end{proof}

\begin{cor}\label{cor:Eisreducetosupercuspidal} For any $b \in B(G)$ with corresponding pair $(M,b^M) \in \mathcal{P}^+$, the full subcategory
\[  i_{b!}^{\ren}\D(G_b(F),\Qellbar)^\omega \subset \D(\Bun_G) \]
is equal to the full subcategory generated under finite colimits and retracts by the images of the functors $\Eis_{L U_{0}^- !}i_{b^L!}^{L}$ restricted to compact supercuspidal objects in $\D(L_{b^L}(F),\Qellbar)$, as $(L,b^L) \in \mathcal{P}$ runs over all pairs retracting to $(M,b^M)$.
\end{cor}
\begin{proof}
This is an immediate consequence of the previous proposition. The only point to observe is that as $(L,b^L)$ varies over all pairs retracting to $(M,b^M)$, the images of the functors $i_{Q_{b^L}}^{G_b}$ on compact supercuspidals generate $\D(G_b(F),\Qellbar)^{\omega}$ under finite colimits and retracts, which is an immediate consequence of Bernstein's basic structure theory.
\end{proof}

\subsection{Finite and restricted sheaves}\label{ss:restrBunGside}

In this section we fix our coefficient ring as $\Qellbar$.

\begin{prop} \label{prop:finite4ways}
The following four conditions on a sheaf $A\in \D(\Bun_G)$ are equivalent.
\begin{propenum}
    \item $\sum_{b\in B(G), n\in \mathbf{Z}} \mathrm{length} H^n(i_{b}^{\ast} A) < +\infty$
    
    \item $A$ is obtained from sheaves of the form $i_{b!}\pi$, for $\pi$ an irreducible smooth $G_b(F)$-representation, via finitely many shifts and cones.
    
    \item $A$ is compact and ULA.
    
    \item $A$ is compact and the natural map $\mathfrak{Z}_{G}^{\mathrm{spec}} \to \mathrm{End}(A)$ factors over an Artinian $\Qellbar$-algebra.
\end{propenum}
\end{prop}
\begin{proof}
The implications (i) $\Leftrightarrow$ (ii) $\Rightarrow$ (iii) are trivial. For (iii) $\Rightarrow$ (iv), note that $\Hom(B,C)$ is a finite-dimensional vector space if $B$ is compact and $C$ is ULA. Therefore $\mathrm{End}(A)$ is an Artinian $\Qellbar$-algebra for any $A$ which is compact and ULA. 

Finally, we show (iv) $\Rightarrow$ (ii). Fix $A$ as in (iv). Since compact sheaves have finite support, it's enough to show that $i_{b}^{\ast}A$ has finite length for all $b$. Since $i_{b}^{\ast}A$ is compact, it is admissible relative to the Bernstein center of $G_b(F)$ \cite{Ber}. Now using the second condition in (iv) and arguing as in the proof of \cite[Theorem 1.6.3]{Beijing} shows that $i_{b}^{\ast}A$ is admissible in the absolute sense, so it is both compact and admissible, hence of finite length.
\end{proof}

\begin{defn}
A sheaf $A \in \D(\Bun_G)$ is \emph{finite} if it satisfies the equivalent conditions of the previous proposition. We write $\D(\Bun_G)_{\fin}$ for the full subcategory of finite sheaves.
\end{defn}

\begin{prop}The full subcategory $\D(\Bun_G)_{\fin}$ is stable under finite colimits and retracts, Hecke operators, and Bernstein-Zelevinsky duality.
\end{prop}
\begin{proof}
Condition (iii) in Proposition \ref{prop:finite4ways} gives stable under finite colimits and retracts. Stability under Hecke operators follows also follows from (iii), since Hecke operators preserve compactness and ULAness individually (in fact this even holds for the whole spectral action by \cref{thm:spectraldualityuseful}). 

Stability under Bernstein-Zelevinsky duality follows from the fourth condition in Proposition \ref{prop:finite4ways}. See also \cite[Theorem 1.6.3]{Beijing} for another (related) proof.
\end{proof}

Recall that for any presentable $\QCoh(\Par_G)$-linear category $\cat C$ we defined a full subcategory of \enquote{restricted} objects $\mathcal{C}^{\restr} \subset \mathcal{C}$ in \cref{def:Crestr}. This paradigm applies to $\D(\Bun_G)$ equipped with the spectral action of $\QCoh(\Par_G)$, and our next goal is to compute $\D(\Bun_G)^{\restr}$ in explicit terms. For this, let $\D(\Bun_G)_{\mathrm{indfin}}$ denote the ind-completion of the category of finite sheaves on $\Bun_G$. We refer to objects of this category as ind-finite sheaves. Since finite sheaves are compact, the canonical inclusion $\D(\Bun_G)_{\mathrm{fin}} \subset \D(\Bun_G)$ extends to a fully faithful embedding $\D(\Bun_G)_{\mathrm{indfin}} \subset \D(\Bun_G)$.

\begin{thm}\label{thm:restrsheavesBunG}
We have an identification
\[\D(\Bun_G)^{\restr} = \D(\Bun_G)_{\mathrm{indfin}}\]
as full subcategories of $\D(\Bun_G)$. Moreover, the projection $\D(\Bun_G) \to \D(\Bun_G)^{\restr}$ is given by $A\mapsto \mathcal{O}^{\restr} \ast A$.
\end{thm}
\begin{proof}
First we show the inclusion of $\D(\Bun_G)^{\restr}$ in $\D(\Bun_G)_{\mathrm{indfin}}$. For this, recall from \cref{rslt:explicit-description-of-O-restr} that $\mathcal{O}^{\restr}$ can be written as a filtered colimit $\varinjlim_i C_i$ with all $C_i \in \Perf(\Par_G)_{\fin}$. Given any $A \in \D(\Bun_G)^{\restr}$, we can write $A = \varinjlim_j A_j$ with $A_j \in \D(\Bun_G)^\omega$, in which case we get
\begin{align*}
    A = \mathcal{O}^{\restr} \ast A = \varinjlim_{i,j} C_i \ast A_j
\end{align*}
where in the first line we used the assumption that $A$ is restricted. Now observe that the spectral action restricts to a functor $\Perf(\Par_G)_{\fin} \times \D(\Bun_G)^{\omega} \to \D(\Bun_G)_{\fin}$: if $\mathcal{F}$ is finite perfect and $A$ is compact, then $\mathcal{F} \ast A$ satisfies the condition in Proposition \ref{prop:finite4ways}.(iv), hence is finite. Therefore $C_i \ast A_j$ is a finite sheaf for all $i,j$, so $A$ is ind-finite.

For the opposite inclusion, it suffices to show that any finite sheaf is restricted. Let $A$ be a finite sheaf. Without loss, we can assume $A$ is supported at a single semisimple $L$-parameter $\phi$ (in the sense of \cref{def:spectral-support}). Let $f_1,\dots,f_n$ generate the associated maximal ideal $\mathfrak{m}_{\phi} \subset \mathfrak{Z}_{G}^{\spec}$, so the $\phi$-summand of $\mathcal{O}^{\restr}$ is isomorphic to the (alternating) Cech complex
\[ C:= \mathcal{O}_{\Par_G} \to \oplus_i \mathcal{O}_{\Par_G}[\tfrac{1}{f_i}] \to \oplus_{i<j}\mathcal{O}_{\Par_G}[\tfrac{1}{f_i f_j}] \to \cdots. \] Then $\mathcal{O}^{\restr} \ast A = C \ast A$. Note that the canonical map $C \to \mathcal{O}_{\Par_G}$ has cone $K$ given by a complex whose terms are finite direct sums of localizations $\mathcal{O}_{\Par_G}[\tfrac{1}{f_{i_1} \cdots f_{i_m}}]$. In particular, this cone has a finite filtration such that some $f_i$ acts invertibly on each graded piece. On the other hand, each $f_i$ acts nilpotently on $A$. Therefore $K \ast A = 0$, so the canonical map $\mathcal{O}^{\restr} \ast A = C \ast A \to A$ is an isomorphism, showing that $A$ is restricted.

The final claim follows from the explicit description of $\D(\Bun_G)^{\restr}$ together with Proposition \ref{prop:restrnaiveadjoint}.
\end{proof}

Although we will not need the following result, we include it for completeness.

\begin{thm}\label{thm:ULAisrestricted}Any ULA sheaf on $\Bun_G$ is restricted. Equivalently, by Theorem \ref{thm:restrsheavesBunG}, any ULA sheaf on $\Bun_G$ is a filtered colimit of finite sheaves.
\end{thm}
\begin{proof}
Fix a ULA sheaf $A$. By \cite{HanAppendix}, there is a canonical direct sum decomposition $A=\oplus_{\phi} A_{\phi}$ indexed by semisimple $L$-parameters, where all irreducible subquotients of all renormalized stalks of $A_{\phi}$ have Fargues-Scholze parameter $\phi$. Replacing $A$ by $A_{\phi}$, it is enough to see that any ULA sheaf $A$ spectrally supported at a single parameter $\phi$ is restricted. 

To see that $\mathcal{O}^{\restr} \ast A \to A$ is an isomorphism, it suffice to show that for any compact sheaf $B$, the canonical map
\[\RHom(B,\mathcal{O}^{\restr} \ast A) \to \RHom(B,A)\]
is an isomorphism. Fix such a $B$. We first observe that $\RHom(B,\mathcal{O}^{\restr} \ast A)$ is a perfect complex. To see this, we may replace $\mathcal{O}^{\restr}$ with its $\phi$-summand as in the previous proof, which can be written as a colimit of Koszul complexes $C_n$ as explained in \cref{rslt:explicit-description-of-O-restr}. Then 
\[\RHom(B,\mathcal{O}^{\restr}_{\phi} \ast A) = \varinjlim_n \RHom(B,C_n \ast A)\]
where each term on the right side is a perfect complex of $\Qellbar$-vector spaces whose amplitude and total dimension in each degree are uniformly bounded in $n$. Therefore, the colimit is still a perfect complex. Now, writing $K$ for the cone of $\mathcal{O}^{\restr}_{\phi} \to \mathcal{O}_{\Par_G}$, we deduce that $\RHom(B,K\ast A)$ is a perfect complex, and the same argument as in the previous proof then implies that $\RHom(B,K\ast A)=0$.
\end{proof}
Summarizing, we have inclusions of full subcategories
\[\D(\Bun_G)_{\fin} \subset \D(\Bun_G)^{\mathrm{ULA}} \subset \D(\Bun_G)_{\mathrm{indfin}}=\D(\Bun_G)^{\restr} \subset \D(\Bun_G),\]
all stable under Hecke operators.

\section{Generalities on the categorical conjecture}

 We fix some notation which will remain in play for the rest of the paper. Until further notice our coefficient ring will be $\Qellbar$ for a fixed prime $\ell \neq p$, which we omit from the notation. Given a connected reductive group $G$ over a finite extension $F/\mathbf{Q}_p$, we write $\mathfrak{Z}_{G}$ for the Bernstein center of $G$, and $\mathfrak{B}(G)$ for the set of primitive idempotents in $\mathfrak{Z}_G$. We identify $\mathfrak{B}(G)$ with the set of Bernstein components in the usual way, and write $\mathfrak{Z}_{G,\mathfrak{s}}=e_{\mathfrak{s}} \mathfrak{Z}_{G}$ for the quotient ring corresponding to an idempotent $\mathfrak{s} \in \mathfrak{B}(G)$. We write $X_G = \coprod_{\mathfrak{s}} \Spec \mathfrak{Z}_{G,\mathfrak{s}}$ for the Bernstein variety of $G$, so $\mathcal{O}(X_G) = \mathfrak{Z}_G$.

On the spectral side, we write $X_{G}^{\spec}$ for the coarse quotient of the stack $\Par_G$, and we set $\mathfrak{Z}_{G}^{\spec} := \mathcal{O}(X_{G}^{\spec})$. We write $\mathfrak{B}^{\spec}(G)$ for the set of connected components of $X_{G}^{\spec}$. The Fargues-Scholze construction gives a ring map $\mathfrak{Z}_{G}^{\spec} \to \mathfrak{Z}_G$, which induces maps $X_{G} \to X_{G}^{\spec}$ and $\mathfrak{B}(G) \to \mathfrak{B}^{\spec}(G)$.

Throughout this chapter, $G$ will denote a fixed \emph{quasisplit} connected reductive group over a fixed finite extension $F/\mathbf{Q}_p$.

\subsection{Whittaker data and the functor \texorpdfstring{$a_\psi$}{a\textunderscore psi}}

Recall that a \emph{Whittaker datum} for $G$ is a pair $(B,\psi)$ where $B\subset G$ is a Borel subgroup and $\psi:U(F)\to \Qellbar^\times$ is a non-degenerate character of the unipotent radical of $B$. For any such pair, we may form the smooth $G(F)$-representation \[W_\psi=\mathrm{c-ind}_{U(F)}^{G(F)} \psi.\]
This is the \emph{Whittaker representation} attached to the given datum. Since $U(F) \subset G(F)$ is not compact open, $W_{\psi}$ is not a compact representation of $G(F)$. Nevertheless, this still turns out to be an extremely nice object. 
\begin{thm}\label{thm:Whittakerbasics} Let $W_\psi$ be the Whittaker representation as defined above.
\begin{propenum}
\item The representation $W_\psi$ is projective.
\item The endomorphism ring $\mathrm{End}_{G(F)}(W_\psi)$ is commutative.
\item For any Bernstein component $\mathfrak{s} \in \mathfrak{B}(G)$, the corresponding summand $W_{\psi,\mathfrak{s}}:= e_{\mathfrak{s}} W_\psi$ is finitely generated. Moreover, either $W_{\psi,\mathfrak{s}}=0$ or the induced map
\[ \mathfrak{Z}_{G,\mathfrak{s}} \to \mathrm{End}_{G(F)}(W_{\psi,\mathfrak{s}}) \]
is an isomorphism.
\item For any Bernstein component $\mathfrak{s}$, there is an isomorphism $\mathbf{D}_{\mathrm{coh}}(W_{\psi,\mathfrak{s}}) \simeq W_{\psi^{-1},\mathfrak{s}^\vee}$.
\end{propenum}
\end{thm}
\begin{proof}
Part (i) is a theorem of Chan-Savin \cite{CS}, with a different argument due to the first author \cite[Appendix A]{Beijing}. Part (ii) is classical and goes back in some form to Gelfand-Graev. Part (iii) is a theorem of Bushnell-Henniart \cite{BH}. Part (iv) is \cite[Theorem A.0.1.iii]{Beijing}.
\end{proof}

For us, a choice of Whittaker representation will play the role of a basepoint sheaf which normalizes the categorical local Langlands correspondence. The first key functor in this paradigm is the following.

\begin{defn}Let $G$ be a quasisplit group with a fixed choice of Whittaker data $(B,\psi)$. We define a functor $a_\psi: \QCoh(\Par_G)  \to \D(\Bun_G)$ by the formula $a_{\psi}(\mathcal{F})=\mathcal{F} \ast i_{1!}W_{\psi}$.
\end{defn}
In other words, $a_\psi$ is the functor of \textbf{a}cting spectrally on $i_{1!}W_{\psi}$. 
\begin{prop}\label{prop:apsibasics}The functor $a_{\psi}$ is colimit-preserving and $\QCoh(\Par_G)$-linear.  For any $\chi \in X^{\ast}(Z(\hat{G})^{W_F}) = \pi_0 \Bun_G$, $a_{\psi}$ carries the $\chi$-graded summand of $\QCoh(\Par_G)$ towards sheaves supported on the $-\chi$-component of $\Bun_G$.
\end{prop}
Recall that we identify components of $\Bun_G$ with basic isocrystals via the $\kappa$ map, which in turn are identified with $X^{\ast}(Z(\hat{G})^{W_F})$ via the Kottwitz isomorphism. If $G=\mathbf{G}_m$, for instance, the connected component given by the bundle $\mathcal{O}(n)$ corresponds to the weight $-n$ character of $\mathbf{G}_m$ on the spectral side. 
\begin{proof}
Immediate from the construction of the spectral action, and equivalent to Lemma 5.3.3 of \cite{ZouTori}. Note that $c_1$ should be replaced by $\kappa$ in the statement of \cite[Lemma 5.3.2]{ZouTori} (and everywhere else in that paper); this is clear if one reads the proof.
\end{proof}

\begin{prop}\label{prop:apsifinitebasic}The following conditions are equivalent.
\begin{propenum}
\item The functor $a_{\psi}$ preserves compact objects.
\item For any connected component $\beta \in \mathfrak{B}^{\spec}(G)$ with associated idempotent $e_\beta \in \mathfrak{Z}^{\spec}_{G}$, the summand $e_\beta W_\psi$ is finitely generated.
\end{propenum}
\end{prop}
Note that in our current state of knowledge, we cannot prove in general that the map $X_G \to X_{G}^{\spec}$ has finite fibers, so the image of the idempotent $e_\beta$ under the map $\mathfrak{Z}^{\spec}_{G} \to \mathfrak{Z}_{G}$ could in principle cut out an infinite set of Bernstein components for $G$. Conversely, if $X_G \to X_{G}^{\spec}$ has finite fibers, it is easy to see that the equivalent conditions of Proposition \ref{prop:apsifinitebasic} both hold for any choice of Whittaker datum. 
\begin{proof} If $a_\psi$ preserves compact objects, then since $e_\beta \mathcal{O}_{\Par_G} \in \Perf^{\mathrm{qc}}(\Par_G)$ is compact, we get that also $e_\beta W_\psi = i_{1}^{\ast} a_\psi(e_\beta \mathcal{O}_{\Par_G})$ is compact.

Conversely, if $\mathcal{F} \in \Perf^{\mathrm{qc}}(\Par_G)$ is a given compact object, then $\mathcal{F} = \oplus_{\beta} e_\beta \mathcal{F}$ where the sum runs over the finite set of those $\beta$ for which $e_\beta \mathcal{F}\neq 0$. Then $a_\psi(\mathcal{F}) = \oplus_{\beta} \mathcal{F} \ast i_{1!} e_\beta W_{\psi}$, and the hypothesis in (ii) guarantees that each $i_{1!} e_\beta W_{\psi}$ is compact. Since the spectral action of $\Perf(\Par_G)$ preserves compact sheaves on $\Bun_G$, we get the desired result.
\end{proof}

Under the conditions of Proposition \ref{prop:apsifinitebasic}, we can prove a duality theorem for $a_{\psi}$.

\begin{thm}\label{thm:apsiduality}
Assume the functor $a_{\psi}$ preserves compact objects. Then the functor $a_{\psi^{-1}}$ also preserves compact objects, and for any $\mathcal{F}\in\mathrm{Perf^{qc}}(\mathrm{Par}_{G})$, we have
a natural duality
\[
\mathbf{D}_{\mathrm{BZ}}a_{\psi}(\mathcal{F})\cong a_{\psi^{-1}}\mathbf{D}_{\mathrm{tw.GS}}(\mathcal{F}).
\]
\end{thm}

\begin{proof}
For a given $\mathcal{F}$, let $e \in \mathfrak{Z}^{\spec}_{G}$ be an idempotent cutting out finitely many components containing the support of $\mathcal{F}$. Then $\mathcal{F}=e\mathcal{F}$, and on the other hand $e W_{\psi}$ is compact by the previous proposition. Then we compute
\begin{align*}
    \mathbf{D}_{\mathrm{BZ}}a_{\psi}(\mathcal{F}) & = \Dbz( \mathcal{F} \ast i_{1!}W_{\psi}) \\
    & = \Dbz(\mathcal{F} \ast i_{1!} e W_{\psi})\\
    & = \Dtwgs(\mathcal{F}) \ast \Dbz(i_{1!} e W_{\psi}) \\
    & = \Dtwgs(\mathcal{F}) \ast i_{1!} e^{\vee}W_{\psi^{-1}}\\
    & = \Dtwgs(\mathcal{F}) \ast i_{1!} W_{\psi^{-1}}\\
    & = a_{\psi^{-1}}\mathbf{D}_{\mathrm{tw.GS}}(\mathcal{F}).
\end{align*}
Here the first and last lines follow immediately from the definitions. The second and fifth lines follow from the assumption on $e$ and the linearity of the spectral action over the action of the spectral Bernstein center. The third line follows from the duality theorem for the spectral action, Theorem \ref{thm:spectraldualityuseful}.(i). The fourth line follows from the duality theorem for the Whittaker representation, Theorem \ref{thm:Whittakerbasics}.(iv), together with the easy identification $\Dbz i_{1!} = i_{1!} \mathbf{D}_{\mathrm{coh}}$.
\end{proof}

\subsection{Change of Whittaker datum}\label{ss:changeofWhittaker}

Recall that two Whittaker data $(B_1, \psi_1)$ and $(B_2,\psi_2)$ are isomorphic if there is an element $g \in G(F)$ such that $B_{1}^{g} = B_2$ and $\psi_{1}^{g} = \psi_2$. It is well-known that two Whittaker data $(B_i,\psi_i)$ are isomorphic if and only if there exists an isomorphism $W_{\psi_1} \simeq W_{\psi_2}$, and that the set $Wh_G$ of isomorphism classes of Whittaker data is finite. 

If $G$ has connected center, there is a unique isomorphism class of Whittaker data. However, we are typically not so lucky, and in general there are multiple isomorphism classes of Whittaker data for $G$. For instance, if $G=\mathrm{SL}_2$ and $p>2$, there are 4 isomorphism classes of Whittaker data. In general, there is a natural \emph{simply transitive} action of the finite abelian group $\mathbf{A}_G := G_{\mathrm{ad}}(F) / G(F)$ on the set $Wh_G$. Note that if we fix a Borel and maximal torus $T \subset B$, the natural map $T_{\mathrm{ad}}(F) / T(F) \to \mathbf{A}_G$ is an isomorphism, so in particular we can obtain \emph{all} Whittaker data by fixing the Borel and only varying the character on the unipotent radical.\footnote{This is some justification for our decision to notate most objects depending on the pair $(B,\psi)$ in terms of $\psi$ alone.}

The following important conjecture is related to some speculations of Zhu \cite[Section 3]{Zhu}, and is expected to be proved in forthcoming joint work of the first author with Eugen Hellmann and Yifei Zhao. 
\begin{conjecture}\label{conj:ZhaoWhittaker} There is a natural group homomorphism
\begin{align*}
\mathbf{A}_G &  \to \mathrm{Pic}(\Par_G)\\
 \theta & \mapsto \mathcal{L}_{\theta}
\end{align*}
with the following property: If $(B,\psi)$ and $(B',\psi')$ are two Whittaker data and $\theta \in \mathbf{A}_G$ is the unique element such that $(B',\psi')=\theta \cdot (B,\psi)$, there is an isomorphism $i_{1!}W_{\psi'} \simeq \mathcal{L}_{\theta} \ast i_{1!} W_\psi$.
\end{conjecture}
In particular, maintaining this notation, we get a natural equivalence of functors $a_{\psi'}(-) = a_{\psi}(\mathcal{L}_{\theta} \otimes -)$. Since $\mathcal{L}_{\theta}$ is a line bundle, tensoring by $\mathcal{L}_{\theta}$ preserves compact objects, so we immediately see that the conditions of Proposition \ref{prop:apsifinitebasic} hold for \emph{one} Whittaker datum if and only if they hold for \emph{all} Whittaker data.

More generally, in the remainder of this paper, we will encounter a series of functors between $\D(\Bun_G)$ and some category $\mathcal{C}$ equipped with a $\otimes$-action of $\Perf(\Par_G)$, denoted (in order of appearance) $a_{\psi}$, $c_{\psi}$, $c_{\psi}^{\mathrm{cts}}$, $\mathbf{L}_{\psi}$, $t_\psi$, $\mathbf{R}_{\psi}$. Each of these functors is $\Perf(\Par_G)$-linear and depends on a choice of Whittaker datum, and for each functor it will be transparent from Conjecture \ref{conj:ZhaoWhittaker} that varying the Whittaker datum changes the functor via either pre- or post-composing with the self-equivalence of $\mathcal{C}$ induced by $\mathcal{L}_{\theta} \otimes -$ for the appropriate choice of $\theta$. In particular, if $F_{\psi}$ is any one of these functors, and \textbf{P} is any of the properties ``is colimit-preserving'', ``is conservative'', ``preserves compact objects'', ``is fully faithful'', ``is essentially surjective'', ``is an equivalence'', or any other qualitative categorical property, then it follows from Conjecture \ref{conj:ZhaoWhittaker} that $F_{\psi}$ has \textbf{P} for \emph{some} Whittaker datum if and only if it has \textbf{P} for \emph{all} Whittaker data. We will assume the truth of Conjecture \ref{conj:ZhaoWhittaker} throughout the remainder of this paper. The reader who is uncomfortable with this assumption can easily keep track of the various dependencies on choices of Whittaker data in our subsequent analysis.

\subsection{The functor \texorpdfstring{$c_{\psi}$}{c\textunderscore psi} and its colimit-preserving cousin}
In this section we introduce two (extremely) closely related functors from the automorphic side to the spectral side.

\begin{defn}The functor $c_\psi : \D(\Bun_G) \to \QCoh(\Par_G)$ is the right adjoint of $a_\psi$.
\end{defn}
Note that $a_\psi$ is colimit-preserving, so $c_\psi$ exists by the adjoint functor theorem. This is the functor of enhanced Whittaker \textbf{c}oefficient. It is also the \textbf{c}oarse Langlands functor. From first principles, it is easy to see that $c_{\psi}$ is $\Perf(\Par_G)$-linear. 

Without more conditions on the group, it is not clear that $c_\psi$ is colimit-preserving or $\QCoh(\Par_G)$-linear. However, it turns out there is a natural variant of $c_\psi$ with both of these properties. To explain this, we briefly recall the self-duality of $\QCoh(\Par_G)$. Let $X$ be a disjoint union of perfect QCA stacks over some characteristic
zero field $k$. Then $\mathrm{QCoh}(X)$ is canonically self-dual,
since it is compactly generated by $\mathrm{Perf^{qc}}(X)$, and the
latter has an obvious self-duality. The evaluation pairing sends $\mathcal{F\boxtimes G}\in\mathrm{QCoh}(X)\otimes\mathrm{QCoh}(X)$
to $\Gamma_{!}(X,\mathcal{F}\otimes\mathcal{G})$, where $\Gamma_{!}(X,-):\mathrm{QCoh}(X)\to \D(k)$
is the unique colimit-preserving extension of the functor $R\mathrm{Hom}(\mathcal{O}_{X},-)$
from $\mathrm{Perf^{qc}}(X)$. In what follows we apply this to $X= \Par_G$, but see Example \ref{exmpl:QCohselfduality} for a more general discussion.

Since $a_\psi$ is a colimit-preserving functor between dualizable categories, it has a canonical dual $a_{\psi}^{\vee}:\D(\mathrm{Bun}_{G})^{\vee}\to\mathrm{QCoh}(\mathrm{Par}_{G})^{\vee}$.
Pre- and post-composing with the canonical Bernstein-Zelevinsky self-duality of $\D(\mathrm{Bun}_{G})$
and the canonical self-duality of $\mathrm{QCoh}(\mathrm{Par}_{G})$ described above, we regard this as a functor \[ a_{\psi}^{\vee}: \D(\mathrm{Bun}_{G})\to\mathrm{QCoh}(\mathrm{Par}_{G}).\]
A short calculation shows that this functor is characterized by the formula
\begin{align*}
\Gamma_{!}(\mathrm{Par}_{G},\mathcal{F}\otimes a_{\psi}^{\vee}(A)) & =\Gamma_{c}(\mathrm{Bun}_{G},a_{\psi}(\mathcal{F})\otimes A)\\
 & =\Gamma_{c}(\mathrm{Bun}_{G},i_{1!}W_{\psi}\otimes(c_{G}^{\ast}\mathcal{F}\ast A)),
\end{align*}
where the first line is a formal consequence of Lemma \ref{lem:evaluativedualitything} and the second line follows
from the duality theorem in the formulation of \cref{rmk:dualityGammaform}.

\begin{defn}We define $c_{\psi}^{\cts}:\D(\Bun_G) \to \QCoh(\Par_G)$ by the formula $c_{\psi}^{\mathrm{cts}}:=c_{G}^{\ast}a_{\psi^{-1}}^{\vee}$, where $c_G$ denotes the Chevalley involution on $\Par_G$ as usual.
\end{defn}

In particular, by
some slight rearranging of the above formulas, we see that this functor is characterized by the formula
\[
\Gamma_{!}(\mathrm{Par}_{G},\mathcal{F}\otimes c_{\psi}^{\mathrm{cts}}(A))=\Gamma_{c}(\mathrm{Bun}_{G},i_{1!}W_{\psi^{-1}}\otimes(\mathcal{F}\ast A))\;\;\;(\dagger)
\]
for all $\mathcal{F\in\mathrm{QCoh}}(\mathrm{Par}_{G})$ and all $A\in \D(\mathrm{Bun}_{G})$.
Since the functors of tensor product, spectral action, $\Gamma_{c}$
and $\Gamma_{!}$ are all colimit-preserving, this formula shows that $c_{\psi}^{\mathrm{cts}}$
is colimit-preserving and $\mathrm{QCoh}$-linear. 

Now if $A$ is compact
and $\mathcal{F}$ is perfect, then $\mathcal{F}\ast A$ is compact,
so its stalk at $b=1$ only hits finitely many Bernstein components of $G(F)$, and
then we can rewrite the right side of $(\dagger)$ above as
\begin{align*}
\Gamma_{c}(\mathrm{Bun}_{G},i_{1!}W_{\psi^{-1}}\otimes(\mathcal{F}\ast A)) & =R\mathrm{Hom}(i_{1!}W_{\psi},\mathcal{F}\ast A)\\
 & =R\mathrm{Hom}(a_{\psi}(\mathcal{F}^{\vee}),A)\\
 & =R\mathrm{Hom}(\mathcal{F}^{\vee},c_{\psi}(A)).
\end{align*}
Here in the first line we used the defining identity of Bernstein-Zelevinsky duality along with Theorem \ref{thm:Whittakerbasics}.(iv), and in the second line we used that $\mathcal{F}^\vee \ast -$ is both left and right adjoint to $\mathcal{F} \ast -$. On the other hand, if $A$ is compact then $c_{\psi}^{\mathrm{cts}}(A)$
has quasicompact support, so we can rewrite the left side of $(\dagger)$ as
\[
\Gamma_{!}(\mathrm{Par}_{G},\mathcal{F}\otimes c_{\psi}^{\mathrm{cts}}(A))=R\mathrm{Hom}(\mathcal{F}^{\vee},c_{\psi}^{\mathrm{cts}}(A)).
\]
Comparing these calculations and letting $\mathcal{F}$ vary arbitrarily,
we obtain a canonical identification $c_{\psi}=c_{\psi}^{\mathrm{cts}}$ as functors on compact objects.
By abstract nonsense, this induces a natural transformation $c_{\psi}^{\mathrm{cts}}\to c_{\psi}$
which is an equivalence on compact objects. Intuitively, $c_{\psi}^{\mathrm{cts}}$ is the best colimit-preserving approximation of $c_{\psi}$. The following extended version of Proposition \ref{prop:apsifinitebasic} is immediate.
\begin{prop}\label{prop:decentequivalents}The following conditions are equivalent.
\begin{propenum}
   \item For any connected component $\beta \in \mathfrak{B}^{\spec}(G)$ with associated idempotent $e_\beta \in \mathfrak{Z}^{\spec}_{G}$, the summand $e_\beta W_\psi$ is finitely generated.
    \item The functor $a_\psi$ preserves compact objects.
    \item The functor $c_\psi$ is colimit-preserving.
    \item The canonical transformation $c_{\psi}^{\mathrm{cts}}\to c_{\psi}$ constructed above is an equivalence.
\end{propenum}
\end{prop}
Again, the discussion at the end of \S\ref{ss:changeofWhittaker} shows that if these conditions hold for one Whittaker datum, they hold for all Whittaker data.

\begin{defn}\label{def:reasonable}A quasisplit group $G$ is \emph{reasonable} if it satisfies the equivalent conditions of Proposition \ref{prop:decentequivalents} for some (equivalently, for all) Whittaker data.
\end{defn}

Of course, we expect that all quasisplit groups are reasonable, but this seems out of reach with current ideas. However, even for reasonable groups, the functor $c_{\psi}=c_{\psi}^{\mathrm{cts}}$ is very hard to compute in any sense.

Nevertheless, we do have one basic finiteness property in general.

\begin{thm}\label{thm:cpsisupercuspidalcoherent}If $A\in \D(\Bun_G)^{\omega}$ is spectrally supported on the supercuspidal locus of $X_{G}^{\spec}$, then $c_{\psi}(A) \in \Coh^\qc(\Par_G)$.
\end{thm}

\begin{proof}
Since $A$ is compact, it is clear that $\mathcal{F}=c_{\psi}(A)$ is supported on finitely many connected components of the parameter stack. Thus, fixing a single supercuspidal connected component $C\subset \Par_G$ and replacing $A$ with $e_{C}A$, we may assume $\mathcal{F}$ is supported on $C$. By Proposition \ref{rslt:coh-criterion-supercuspidal}, to prove that a sheaf $\mathcal{F} \in \QCoh(C)$ is coherent, it suffices to check two conditions:
\begin{itemize}
    \item Only finitely many characters of $Z(\hat{G})^{W_F}$ occur in the grading of $\mathcal{F}$.
    \item For any $\mathcal{G}\in \Perf(C)$, $\Gamma(C,\mathcal{G}\otimes \mathcal{F})$ is cohomologically bounded, and its cohomologies are finitely generated $\mathcal{O}(C)$-modules.
\end{itemize}
In the scenario at hand, the first condition follows immediately from the compatibility of $c_{\psi}$ with the central grading, cf. Proposition \ref{prop:cpsigrading} below. For the second condition, we compute that
\begin{align*}
    \Gamma(C,\mathcal{G}\otimes \mathcal{F}) & = \Gamma(C,\mathcal{G}\otimes c_{\psi}(A)) \\
    & = \Gamma(C, c_{\psi}(\mathcal{G} \ast A)) \\
    & = \RHom(i_{1!}W_{\psi}, \mathcal{G} \ast A) \\
    & = \RHom(W_{\psi}, i_{1}^{\ast}(\mathcal{G} \ast A))
\end{align*}
using the $\Perf(\Par_G)$-linearity of $c_{\psi}$. Now since $A$ is compact and the spectral action of $\Perf$ preserves compact objects, the restriction $B:=i_{1}^{\ast}(\mathcal{G} \ast A)$ is compact, and in particular is supported on finitely many Bernstein components $\mathfrak{s}$. Fixing any one of these components, 
\[\RHom(W_{\psi}, B_{\mathfrak{s}}) = \RHom(W_{\psi,\mathfrak{s}}, B_{\mathfrak{s}})\]
is a bounded complex whose cohomologies are finitely generated $\mathcal{O}(X_{G,\mathfrak{s}})$-modules, using Theorem \ref{thm:Whittakerbasics}.iv together with Bernstein's basic finiteness theorems. Since $\mathcal{O}(C) \to \mathcal{O}(X_{G,\mathfrak{s}})$ is a finite ring map, this gives the desired result.
\end{proof}

In the preceding proof we used the following easy result.
\begin{prop}\label{prop:cpsigrading}
For any $\chi \in X^{\ast}(Z(\hat{G})^{W_F})\cong \pi_0 \Bun_G$, the functors $c_{\psi}$ and $c_{\psi}^{\cts}$ carry sheaves supported on the $\chi$-component of $\Bun_G$ towards the $-\chi$-graded summand of $\QCoh(\Par_G)$.
\end{prop}
\begin{proof}
Immediate from Proposition \ref{prop:apsibasics} and the basic characterizations of $c_{\psi}$ and $c_{\psi}^{\cts}$.
\end{proof}

\begin{prop}\label{prop:fullyfaithfuleasy}
If $G$ is reasonable, the following conditions are equivalent.
\begin{propenum}
\item The functor $a_\psi$ is fully faithful.
\item The natural map $\mathcal{O}_{\Par_G} \to c_\psi(i_{1!}W_\psi)$ is an isomorphism.
\item The quasicoherent sheaf $c_\psi(i_{1!}W_\psi)$ is a line bundle.
\end{propenum}
\end{prop}
\begin{proof}Since $G$ is reasonable, $c_{\psi} = c_{\psi}^{\cts}$ is colimit-preserving and $\QCoh(\Par_G)$-linear. Now $a_{\psi}$ is fully faithful if and only the adjunction
\[\mathcal{F} \to c_{\psi}(a_{\psi}(\mathcal{F}))\]
is an equivalence for all $\mathcal{F} \in \QCoh(\Par_G)$. By $\QCoh$-linearity of both functors, this holds for all $\mathcal{F}$ if and only if it holds for $\mathcal{F}=\mathcal{O}_{\Par_G}$, in which case it reduces to the statement in (ii). This gives the equivalence of (i) and (ii). Since $c_\psi(i_{1!}W_\psi)$ is a unital associative algebra object in $\QCoh(\Par_G)$, the equivalence of (ii) and (iii) is immediate from \cite[Lemma 17.3.7]{GLC3}.
\end{proof}

\subsection{Compatibility with parabolic induction}
In this section we formulate a precise conjectural compatibility of the coarse Langlands functor with parabolic induction. We warn the reader that the most naive guess for this compatibility is wrong, and the correct statement involves a sign. 

To state the conjecture, fix a quasisplit group $G$ together with a Borel and maximal torus $T \subset B$. There is then a canonical bijection between standard parabolic subgroups of $G$ and standard parabolic subgroups of $\phantom{}^L G$. For any parabolic $P$, we have the associated Eisenstein functor $\Eis_{P !} $ as discussed earlier. On the spectral side, recall from Chapter 3 that we have the canonical diagram \[
\xymatrix{\Par_P \ar[r]^p \ar[d]^q & \Par_G \\
\Par_M & \phantom{}
}
\]
of Artin stacks, which gives rise to the
functors $\Eis_{P}^{\spec} = p_{\ast}^{\IndCoh} q^{*\IndCoh}$ and $\CT_{P}^{\spec}=q_{\ast}^{\IndCoh} p^{!\IndCoh}$ on ind-coherent sheaves, along with their coarse variants $\Eis_{P}^{\spec,\coarse} = p_{\ast} q^{\ast}$ and $\CT_{P}^{\spec,\coarse}=q_{\ast} p^{\ast}$  on quasicoherent sheaves. Recall that $\Eis_{P}^{\spec}$ and $\CT_{P}^{\spec}$ are adjoint, while their coarse variants are dual as shown below. Recall also that $\Eis_{P}^{\spec}$ and $\Eis_{P}^{\spec,\coarse}$ both preserve $\Coh^\qc$, and on $\Coh^\qc$ they are literally the same functor (via the equivalence on copies of $\Coh^\qc$ induced by $\Psi$).

Now extend the given Borel to a choice of Whittaker datum $(B,\psi)$. This choice induces Whittaker data $(B\cap M, \psi_M)$ for all standard Levis $M \subset G$.
\begin{conjecture}\label{conj:cpsiEiscompatible}With notation and assumptions as above, the diagram
\[
\xymatrix{\D(\mathrm{Bun}_{M})\ar[r]^{c_{\psi_{M}}^{\cts}}\ar[d]^{\mathrm{Eis}_{P^{-}!}} & \mathrm{QCoh}(\mathrm{Par}_{M})\ar[d]^{\Eis_{P}^{\spec,\coarse}}\\
\D(\mathrm{Bun}_{G})\ar[r]^{c_{\psi}^{\cts}} & \mathrm{QCoh}(\mathrm{Par}_{G})
}
\]
commutes, i.e. there is an equivalence of functors
\[\alpha: c_{\psi}^{\cts} \circ \Eis_{P^- !} \simeq \Eis_{P}^{\spec,\coarse} \circ c_{\psi_M}^{\cts}.\]
\end{conjecture}
The appearance of the opposite parabolic on the automorphic side is not a convention, but a mathematical necessity, as we will see later in some very concrete examples.\footnote{The necessity of taking the opposite parabolic was already understood by Hellmann when he wrote his paper \cite{Hel}, but D.H. was slow to appreciate this point. The opposite parabolic is also (incorrectly) omitted in \cite{BCHN}. D.H. thanks Teruhisa Koshikawa for some helpful discussions regarding these sign issues.} 

Recall that on compact objects, $c_{\psi}=c_{\psi}^{\cts}$, so by the preservation of compact objects under $\Eis_{P^- !}$, the postulated equivalence induces an equivalence
\[\alpha': c_{\psi} \circ \Eis_{P^- !} \simeq \Eis_{P}^{\spec,\coarse} \circ c_{\psi_M}\]
as functors on compact sheaves. Since $c_{\psi}^{\cts}$ is the ind-completion of the restriction of $c_\psi$ to compact sheaves, and both Eisenstein functors commute with colimits, the datum of $\alpha$ is actually equivalent to the datum of $\alpha'$. More substantially, the datum of such an equivalence can be reformulated in terms of the $a_{\psi}$ and constant term functors.

\begin{prop}\label{prop:cpsiEisVSapsiCT}A commutativity datum for the diagram in Conjecture \ref{conj:cpsiEiscompatible} is logically equivalent to a commutativity datum for the diagram
\[
\xymatrix{\mathrm{QCoh}(\mathrm{Par}_{G})\ar[d]^{\CT_{P^-}^{\spec,\coarse}}\ar[r]^{a_{\psi^{-1}}} & \D(\mathrm{Bun}_{G})\ar[d]^{\mathrm{CT}_{P\ast}}\\
\mathrm{QCoh}(\mathrm{Par}_{M})\ar[r]^{a_{\psi_{M}^{-1}}} & \D(\mathrm{Bun}_{M})
}
\]
i.e. the datum of an equivalence $\alpha$ is logically equivalent to the datum of an equivalence
\[\beta: a_{\psi_{M}^{-1}} \circ \CT_{P^-}^{\spec,\coarse} \simeq \CT_{P \ast} \circ a_{\psi^{-1}}.\]
\end{prop}
The equivalence between these two diagrams will be induced by taking \emph{duals} of functors, rather than adjoints. 

\begin{proof}[Proof of Proposition \ref{prop:cpsiEisVSapsiCT}]All the categories appearing here are canonically self-dual over $D(\Qellbar)$ via evaluative self-dualities as in Examples \ref{exmpl:QCohselfduality} and \ref{exmpl:BunGselfduality}. Now we compute duals of functors.

\begin{enumerate}
    \item The functors $a_{\psi^{-1}}$ and $c_{G}^\ast c_{\psi}^{\cts}$ are dual to each other, and similarly for $M$. This follows immediately from the definition of $c_{\psi}^{\cts}$.

    \item The functors $\Eis_{P}^{\spec,\coarse}$ and $\CT_{P}^{\spec,\coarse}$ are dual to each other. Since the maps $p^{\spec}$ and $q^{\spec}$ are quasicompact, this follows from two applications of the last part of Example \ref{exmpl:QCohselfduality}.

    \item The functors $\Eis_{P^- !}$ and $\CT_{P^- !}\cong \CT_{P \ast}$ are dual to each other. This follows from the calculations performed in the proof of \cite[Theorem 3.2.1]{HHS}, in particular the seven-line equation on p. 21 of loc. cit.
\end{enumerate}
Putting things together, we see that duality of functors induces an equivalence between natural equivalences $\beta$ as in Proposition \ref{prop:cpsiEisVSapsiCT} and natural equivalences
\[\beta^{\vee}: c_{G}^{\ast}c_{\psi}^{\cts} \circ \Eis_{P^- !} \simeq \Eis_{P^-}^{\spec,\coarse} \circ c_{M}^{\ast} c_{\psi_{M}}^{\cts}. \]
Since $c_{G}^{\ast}$ is an involution, we may move it to the other side of this equivalence. Now observing that \[c_{G}^{\ast} \Eis_{P^-}^{\spec,\coarse} c_{M}^{\ast} \cong \Eis_{P}^{\spec,\coarse}\]
gives the desired result.
\end{proof}

\begin{rmk}There is no obvious candidate for an equivalence $\beta$: it must be constructed in the course of proving Conjecture \ref{conj:cpsiEiscompatible}.  However, let us explain one expected compatibility of this transformation, replacing $\psi^{-1}$ with $\psi$ for simplicity. Evaluating $\beta$ on the structure sheaf and noting that $\CT_{P^-}^{\spec,\coarse}(\mathcal{O}_{\Par_G}) = q^{-}_{\ast}\mathcal{O}_{\Par_{P^-}}$ where $q^{-}: \Par_{P^-} \to \Par_M$, we get an induced isomorphism
$a_{\psi_M}(q^{-}_{\ast}\mathcal{O}_{\Par_{P^-}}) \simeq \CT_{P \ast} i_{1!}W_{\psi}$. Precomposing $a_{\psi_M}$ with the canonical map $\mathcal{O}_{\Par_M} \to q^{-}_{\ast}\mathcal{O}_{\Par_{P^-}}$, we get an induced map $i_{1!}^{M} W_{\psi_M} \to \CT_{P \ast} i_{1!}W_{\psi}$, or equivalently a map
\begin{align*}
W_{\psi_M} & \overset{\delta}{\to} i_{1}^{M \ast} \CT_{P \ast} i_{1!}W_{\psi} \\
& = r_{G}^{P^-}W_{\psi} 
\end{align*}
using the identification of functors $i_{1}^{M \ast} \CT_{P \ast} i_{1!} = r_{G}^{P^-}$ proved in \cite{HHS}. Now, by a theorem of Bushnell-Henniart \cite{BH}, there is a \emph{canonical} isomorphism $r_{G}^{P^-}W_{\psi} \cong W_{\psi_M}$, and we expect that $\delta$ agrees with this isomorphism up to scalars. Note the self-consistent appearance of the opposite parabolic.
    
\end{rmk}

\subsection{The categorical conjecture}
We now turn to our central focus. Recall the following soft form of the categorical local Langlands conjecture as stated in the introduction.
\begin{conjecture}The functor $a_\psi$ extends along the full embedding $\Xi:\QCoh(\Par_G)\hookrightarrow \IndCoh(\Par_G)$ to an equivalence of categories $\mathbf{L}_{\psi}:\IndCoh(\Par_G)\simeq \D(\Bun_G)$.
\end{conjecture}
The first key observation is that if this equivalence exists, it is automatically unique, and is completely pinned down by the functor $c_\psi$.

\begin{prop}\label{prop:CLLCuniqueness}There is \emph{at most one} equivalence $\mathbf{L}_\psi : \D(\Bun_G) \overset{\sim}{\to} \IndCoh(\Par_G)$ satisfying either of the following equivalent conditions:
\begin{propenum}
    \item The diagram
\[
\xymatrix{\mathrm{QCoh}(\mathrm{Par}_{G})\ar[r]^{a_{\psi}}\ar[dr]^{\Xi} & \D(\mathrm{Bun}_{G})\ar[d]_{\wr}^{\mathbf{L}_{\psi}}\\
 & \mathrm{IndCoh}(\mathrm{Par}_{G})
}
\]
commutes.

\item The diagram
\[
\xymatrix{\D(\mathrm{Bun}_{G})\ar[d]_{\wr}^{\mathbf{L}_{\psi}}\ar[r]^{c_{\psi}} & \mathrm{QCoh}(\mathrm{Par}_{G})\\
\mathrm{IndCoh}(\mathrm{Par}_{G})\ar[ur]^{\Psi}
}
\]
commutes.
\end{propenum}
Moreover, such an equivalence exists if and only if $c_\psi$ restricts to an equivalence of categories \[c_\psi : \D(\Bun_G)^{\omega} \overset{\sim}{\to} \Coh^\qc(\Par_G),\] in which case $\mathbf{L}_{\psi}$ is obtained from this equivalence by ind-completion.
\end{prop}
\begin{proof}The conditions (i) and (ii) are equivalent by passing to left or right adjoints. If an equivalence $\mathbf{L}_{\psi}$ satisfying the condition in (ii) exists, then for any $A\in \D(\Bun_G)^\omega$, we have \[c_\psi(A) = \Psi (\mathbf{L}_{\psi}(A)) \in \Coh^\qc(\Par_G),\] using the fact that $\Psi$ induces an equivalence on the evident copies of $\Coh^\qc$. Thus $c_{\psi}$ restricts to an equivalence from $\D(\Bun_G)^\omega$ to $\Coh^\qc(\Par_G)$, and one then recovers $\mathbf{L}_{\psi}$ from $c_{\psi}$ by ind-completing, which gives the uniqueness. Conversely, if $c_\psi$ restricts to an equivalence of categories $c_\psi : \D(\Bun_G)^{\omega} \overset{\sim}{\to} \Coh^\qc(\Par_G)$, the ind-completion of the restriction of $c_{\psi}$ clearly satisfies the conditions in (ii).
\end{proof}
Turning this around, we can now state the official form of the categorical local Langlands conjecture.
\begin{conjecture}\label{conj:CLLC}The functor $c_{\psi}$ restricts to an equivalence of categories $c_\psi : \D(\Bun_G)^{\omega} \overset{\sim}{\to} \Coh^\qc(\Par_G)$.

Equivalently, the functor $c_\psi$ carries $\D(\Bun_G)^{\omega}$ into $\Coh^\qc(\Par_G)$, and the resulting ind-completed functor
\[\mathbf{L}_{\psi}: \D(\Bun_G) \to \IndCoh(\Par_G)\]
is an equivalence of categories.
\end{conjecture}

Again, the discussion at the end of \S \ref{ss:changeofWhittaker} applies, showing that this conjecture holds for one Whittaker datum if and only if it holds for all Whittaker data. 

In our subsequent analysis, it will be very helpful to separate the existence of the functor $\mathbf{L}_{\psi}$ from the question of whether it induces an equivalence of categories. The following proposition is an easy variant of the proof of Proposition \ref{prop:CLLCuniqueness}.

\begin{prop}\label{prop:Lpsiexistence}The following conditions are equivalent.
\begin{propenum}
\item The functor $c_{\psi}$ carries $\D(\Bun_G)^{\omega}$ into $\Coh^\qc(\Par_G)$.

\item There is a (necessarily unique) colimit-preserving functor $\mathbf{L}_{\psi}:\D(\Bun_G) \to \IndCoh(\Par_G)$ which preserves compact objects and such that the diagram 
\[
\xymatrix{\D(\mathrm{Bun}_{G})\ar[r]^{\mathbf{L}_{\psi}\;\;}\ar[dr]^{c_{\psi}} & \mathrm{IndCoh}(\mathrm{Par}_{G})\ar[d]^{\Psi}\\
& \mathrm{QCoh}(\mathrm{Par}_{G})
}
\]
commutes.
\end{propenum}
\end{prop}

We will use the shorthand phrase ``$\mathbf{L}_{\psi}$ exists'' to indicate that either of the equivalent conditions of Proposition \ref{prop:Lpsiexistence} holds true. Again, these conditions hold for one Whittaker datum if and only if they hold for all Whittaker data. We also note that if $\mathbf{L}_{\psi}$ exists, it is automatically $\QCoh(\Par_G)$-linear, and it automatically carries sheaves on the $\chi$-component of $\Bun_G$ towards $-\chi$-graded ind-coherent sheaves.

Let us comment on the use of the letter L. As in the geometric Langlands literature, we refer to the functor $\mathbf{L}_{\psi}$ as the \textbf{L}anglands functor. In light of the previous proposition, it is also instructive to view this functor as a canonical \textbf{l}ifting of $c_{\psi}$ along the functor $\Psi$. Finally, in the next chapter it will be very important to regard $\mathbf{L}_{\psi}$ as a \textbf{l}eft adjoint functor.

\begin{rmk}If $\mathbf{L}_{\psi}$ exists, then for compact $A$ there is literally no difference between the coherent sheaves $\mathbf{L}_{\psi}A$ and $c_{\psi}A$, because $\Psi$ mapping $\Coh^\qc \subset \IndCoh$ onto $\Coh^\qc \subset \QCoh$ is literally the identity functor on $\Coh^\qc$. In other words, the only difference between $\mathbf{L}_{\psi}A$ and $c_{\psi}A$ for compact $A$ is the container (namely, $\IndCoh$ or $\QCoh$) we regard them as living inside. As such, we will sometimes be cavalier about the difference between $\mathbf{L}_{\psi}$ and $c_{\psi}$ on compact sheaves.
\end{rmk}

In the rest of this section, we record some automatic compatibilities and consequences of the categorical local Langlands conjecture. Later we will turn things around and use these results as guidance on our path towards a proof of the conjecture. We start with an easy observation.
\begin{prop}If Conjecture \ref{conj:CLLC} is true, the functor $a_\psi$ preserves compact objects and is fully faithful.
\end{prop}

\begin{proof}This is immediate from the diagram in Proposition \ref{prop:CLLCuniqueness}.(i), since the functor $\Xi$ is fully faithful and preserves compact objects.
\end{proof}

\begin{prop}If Conjecture \ref{conj:CLLC} is true, then $\mathbf{L}_\psi$ automatically restricts to equivalences
\[ \D(\Bun_G)^{\mathrm{ULA}}\overset{\sim}{\to} \Adm(\Par_G) \]
and
\[ \D(\Bun_G)_{\fin}\overset{\sim}{\to} \Coh(\Par_G)_{\fin}. \]
\end{prop}
\begin{proof}We know that $A\in \D(\Bun_G)$ is ULA if and only if $\RHom(B,A)$ lies in $\Perf(\Qellbar)$ for all compact $B$. Assuming Conjecture \ref{conj:CLLC}, this shows that an ind-coherent sheaf $\mathcal{F}$ comes from a ULA sheaf on $\Bun_G$ exactly when $\RHom(\mathcal{G},\mathcal{F}) \in \Perf(\Qellbar)$ for all $\mathcal{G} \in \Coh^\qc(\Par_G)$. But this is exactly the definition of admissibility. This gives the first claim.

For the second claim, suppose given some compact $A \in \D(\Bun_G)$. Then $A$ is finite if and only if $A$ satisfies the fourth condition of Proposition \ref{prop:finite4ways}. Under the $\mathfrak{Z}_{G}^{\spec}$-linearity of the functor $\mathbf{L}_{\psi}$, this immediately translates into the condition that $\mathbf{L}_{\psi}(A)$ be supported on finitely many closed fibers of $q$, i.e. that $\mathbf{L}_{\psi}(A)$ is finite.
\end{proof}

\begin{rmk}
In the previous proof, we did not use Theorem \ref{thm:cohfinclassification}. In fact, we first guessed that Theorem \ref{thm:cohfinclassification} must be true after observing (via the argument above) that the categorical conjecture must match compact resp. ULA resp. finite sheaves with coherent resp. admissible resp. finite sheaves. Then translating the identification
\[\D(\Bun_G)_{\fin} = \D(\Bun_G)^{\omega} \cap \D(\Bun_G)^{\mathrm{ULA}}\]
from Proposition \ref{prop:finite4ways} across the categorical conjecture, one arrives at the statement of Theorem \ref{thm:cohfinclassification}. 
\end{rmk}
More surprisingly, we automatically get compatibility between the dualities on both sides.
\begin{thm}\label{thm:automaticduality}Assume that Conjecture \ref{conj:CLLC} is true.  Then the equivalences $\mathbf{L}_{\psi}$ and $\mathbf{L}_{\psi^{-1}}$  automatically satisfy the following compatibilities with duality:
\begin{propenum}
\item $\Dtwgs \circ \mathbf{L}_{\psi^{-1}} = \mathbf{L}_{\psi} \circ \Dbz$ on all compact sheaves.

\item$\Dtwadm \circ \mathbf{L}_{\psi^{-1}} = \mathbf{L}_{\psi} \circ \Dverd$ on all sheaves.
\end{propenum}
\end{thm}

\begin{proof}Part (ii) is a formal consequence of (i), using the duality exchange formulas on both sides. More precisely, for any given $A \in \D(\Bun_G)$ and $\mathcal{F} \in \Coh^\qc(\Par_G)$, we compute
\begin{align*} \RHom(\mathcal{F},\Dtwadm \mathbf{L}_{\psi^{-1}}A) & = \RHom(\Dtwgs \mathcal{F}, \mathbf{L}_{\psi^{-1}}A)^{\vee} \\
& = \RHom(\mathbf{L}_{\psi^{-1}}^{-1} \Dtwgs \mathcal{F}, A)^{\vee} \\
& = \RHom(\Dbz \mathbf{L}_{\psi}^{-1} \mathcal{F}, A)^{\vee} \\
& = \RHom( \mathbf{L}_{\psi}^{-1} \mathcal{F}, \Dverd A) \\
& = \RHom( \mathcal{F},  \mathbf{L}_{\psi} \Dverd A).
\end{align*}
Here we used duality exchange on both the automorphic and spectral sides, along with (i) and the assumption that $\mathbf{L}_{\psi}$ is an equivalence. Since $\mathcal{F}$ is arbitrary, (ii) follows by Yoneda.

For (i), fix any $A\in \D(\Bun_G)^\omega$ and any $\mathcal{F} \in \Perf^{\mathrm{qc}}(\Par_G)$. Then on one hand, we compute that
\begin{align*}
R\mathrm{Hom}(\mathbf{D}_{\mathrm{tw.GS}}\mathbf{L}_{\psi^{-1}}A,\mathcal{F}) & =R\mathrm{Hom}(\mathbf{D}_{\mathrm{tw.GS}}\mathcal{F},\mathbf{L}_{\psi^{-1}}A)\\
 & =R\mathrm{Hom}(\mathbf{D}_{\mathrm{tw.GS}}\mathcal{F},c_{\psi^{-1}}A)\\
 & =R\mathrm{Hom}(a_{\psi^{-1}}\mathbf{D}_{\mathrm{tw.GS}}\mathcal{F},A)\\
 & =R\mathrm{Hom}(\mathbf{D}_{\mathrm{BZ}}a_{\psi}\mathcal{F},A).
\end{align*}
Here we used the involutivity of twisted Grothendieck-Serre duality in the first line, the agreement of $\mathbf{L}_\psi$ and $c_\psi$ on compact sheaves in the second line, the adjunction between $a_\psi$ and $c_\psi$ in the third line, and the duality theorem for $a_{\psi}$ (Theorem \ref{thm:apsiduality}) in the last line. On the other hand, since
we are assuming that $\mathbf{L}_{\psi}$ is an equivalence, we compute
that
\begin{align*}
R\mathrm{Hom}(\mathbf{L}_{\psi}\mathbf{D}_{\mathrm{BZ}}A,\mathcal{F}) & =R\mathrm{Hom}(\mathbf{D}_{\mathrm{BZ}}A,\mathbf{L}_{\psi}^{-1}\mathcal{F})\\
 & =R\mathrm{Hom}(\mathbf{D}_{\mathrm{BZ}}A,a_{\psi}\mathcal{F})\\
 & =R\mathrm{Hom}(\mathbf{D}_{\mathrm{BZ}}a_{\psi}\mathcal{F},A).
\end{align*}
Here we crucially used the fact that $\mathbf{L}_{\psi}^{-1}|_{\mathrm{Perf^{\mathrm{qc}}}}=a_{\psi}$
holds \emph{automatically }if $\mathbf{L}_{\psi}$ is an equivalence of categories, as in the diagram in Proposition \ref{prop:CLLCuniqueness}.(i).
Putting these two calculations together, we see that
\[
R\mathrm{Hom}(\mathbf{D}_{\mathrm{tw.GS}}\mathbf{L}_{\psi^{-1}}A,\mathcal{F})=R\mathrm{Hom}(\mathbf{L}_{\psi}\mathbf{D}_{\mathrm{BZ}}A,\mathcal{F})
\]
functorially in $A\in D(\mathrm{Bun}_{G})^{\omega}$ and $\mathcal{F}\in\mathrm{Perf^{qc}}(\mathrm{Par}_{G})$.
By Lemma \ref{lem:cohdualitytrick} below, this implies that $\mathbf{D}_{\mathrm{tw.GS}}\mathbf{L}_{\psi^{-1}}=\mathbf{L}_{\psi}\mathbf{D}_{\mathrm{BZ}}$
as functors on compact sheaves. 
\end{proof}
\begin{lem}\label{lem:cohdualitytrick}Let $X$ be a disjoint union of QCA stacks which are perfect and Gorenstein. If $\mathcal{F}_1,\mathcal{F}_2 \in \Coh(X)$ corepresent the same functor on $\Perf^{\mathrm{qc}}(X)$, then $\mathcal{F}_1 = \mathcal{F}_2$.
\end{lem}
\begin{proof} By the assumption that $X$ is Gorenstein, $\Dgs$ preserves $\Coh(X)$ and $\Perf^{\mathrm{qc}}(X)$ and induces compatible anti-equivalences on them. Since $\mathcal{F}_1$ and $\mathcal{F}_2$ corepresent the same functor on $\Perf^{\mathrm{qc}}(X)$, this formally implies that $\Dgs \mathcal{F}_1$ and $\Dgs \mathcal{F}_2$ represent the same (contravariant) functor on $\Perf^{\mathrm{qc}}(X)$. Since $\QCoh(X)=\Ind(\Perf^{\mathrm{qc}}(X))$ by assumption, passing to colimits formally upgrades this to the statement that $\Dgs \mathcal{F}_1$ and $\Dgs \mathcal{F}_2$ represent the same functor on $\QCoh(X)$. But then $\Dgs \mathcal{F}_1 = \Dgs \mathcal{F}_2$ by Yoneda, so then $\mathcal{F}_1 = \mathcal{F}_2$ using again the involutivity of $\Dgs$.
\end{proof}

\subsection{A partial right adjoint}
In this section we allow $G$ to be any reasonable (as in Definition \ref{def:reasonable}) quasisplit group.
Fix a Whittaker datum, and let $a_{\psi}:\mathrm{QCoh}(\mathrm{Par}_{G})\to D(\mathrm{Bun}_{G})$
and $c_{\psi}:D(\mathrm{Bun}_{G})\to\mathrm{QCoh}(\mathrm{Par}_{G})$
be the adjoint pair of functors given by spectral action on the Whittaker
sheaf, and the functor of enhanced Whittaker coefficient. Since $G$ is assumed reasonable, $a_{\psi}$ preserves compact objects and $c_{\psi}$ is colimit-preserving and $\QCoh$-linear.

Our goal in this section is to define a new functor from the spectral side to the automorphic side. This definition is one of the key new ideas in this paper. To motivate this construction, recall from Theorem \ref{thm:automaticduality} that if categorical local Langlands is true, the inverse of the equivalence $\mathbf{L}_{\psi}$ automatically satisfies the duality relation $\Dverd \mathbf{L}_{\psi^{-1}}^{-1} = \mathbf{L}_{\psi}^{-1} \Dtwadm$ on all sheaves. Evaluating $\mathbf{L}_{\psi}^{-1}$ on an admissible ind-coherent sheaf $\mathcal{F}$ and using this relation together with biduality for $\Dtwadm$ on admissible sheaves, we calculate that
\begin{align*}
\mathbf{L}_{\psi}^{-1}\mathcal{F} & = \mathbf{L}_{\psi}^{-1} \mathbf{D}^{2}_{\mathrm{tw.adm}} \mathcal{F} \\
& = \Dverd \mathbf{L}_{\psi^{-1}}^{-1} \Dtwadm \mathcal{F}.
\end{align*}
Now we restrict attention further to $\mathcal{F} \in \Coh(\Par_G)_{\fin}$, i.e. to $\mathcal{F}$ which are admissible and coherent. Then Theorem \ref{thm:cohfinclassification} guarantees that $\Dtwadm \mathcal{F}$ is in $\QCoh$, or equivalently that $\Dtwadm \mathcal{F} = \Xi \Psi \Dtwadm \mathcal{F}$, in which case we can rewrite the above expression further as
\[\mathbf{L}_{\psi}^{-1}\mathcal{F} = \Dverd \mathbf{L}_{\psi^{-1}}^{-1} \Xi \Psi \Dtwadm \mathcal{F}.  \]
Finally, recall from Proposition \ref{prop:CLLCuniqueness}.(i) that if categorical local Langlands is true, then $\mathbf{L}_{\psi^{-1}}^{-1} \Xi = a_{\psi^{-1}}$ automatically, so we can rewrite this one more time as
\[\mathbf{L}_{\psi}^{-1}\mathcal{F} = \Dverd a_{\psi^{-1}} \Psi \Dtwadm \mathcal{F}.  \]
We have thus arrived at the following observation.

$(\dagger)$ If categorical local Langlands is true, then the inverse functor $\mathbf{L}_{\psi}^{-1}$ restricted to $\Coh(\Par_G)_{\fin}$ is equivalent to the composite functor $\Dverd a_{\psi^{-1}} \Psi \Dtwadm \mathcal{F}$. 

This composite functor looks complicated at first glance, but the major gain here is that all terms in the expression $\Dverd a_{\psi^{-1}} \Psi \Dtwadm \mathcal{F}$ are \emph{unconditionally defined} for every quasisplit group. Our goal in this section is to study this functor from first principles.
\begin{defn}\label{defn:tpsi}
Let $t_{\psi}:\mathrm{Coh}(\mathrm{Par}_{G})_{\fin}\to \D(\mathrm{Bun}_{G})$
be the functor defined by the formula
\[
t_{\psi}(\mathcal{F})=\mathbf{D}_{\mathrm{Verd}}a_{\psi^{-1}}\Psi \mathbf{D}_{\mathrm{tw.adm}}\mathcal{F}.
\]
\end{defn}
Again, we emphasize that this definition makes sense unconditionally for any quasisplit group. Note that $t_{\psi}$ is compatible with the spectral action of $\Perf(\Par_G)$, but essentially no other properties are obvious from this definition. Note also that as written, the expression defining $t_\psi$ makes sense on any ind-coherent sheaf, but we are quite sure that evaluating it on anything outside $\Coh(\Par_G)_{\fin}$ will produce garbage.

Regarding the notation, we previously chose the notations $a_{\psi}$
for \textbf{a}ction and $c_{\psi}$ for \textbf{c}oefficient. Here we
choose $t_{\psi}$ for \textbf{t}empered, since the good properties of this functor rely
on the spectral temperization (Theorem \ref{rslt:spectral-temperization}), and/or \textbf{t}wisted, since the definition
of $t_{\psi}$ twists $a_{\psi^{-1}}$ by an excessive number of duality
functors.

To analyze this functor, we'll need the following duality result.
\begin{thm}\label{thm:cpsiVerdier}
For any $A\in \D(\mathrm{Bun}_{G})$, we have a natural isomorphism
\[
c_{\psi}\mathbf{D}_{\mathrm{Verd}}A\cong\Psi\mathbf{D}_{\mathrm{tw.adm}}\Xi c_{\psi^{-1}}A.
\]
\end{thm}

\begin{proof}
For any $\mathcal{F}\in\mathrm{Perf^{qc}}(\mathrm{Par}_{G})$, we compute
that
\begin{align*}
R\mathrm{Hom}(\mathcal{F},\Psi\mathbf{D}_{\mathrm{tw.adm}}\Xi c_{\psi^{-1}}A) & \cong R\mathrm{Hom}(\Xi\mathcal{F},\mathbf{D}_{\mathrm{tw.adm}}\Xi c_{\psi^{-1}}A)\\
 & \cong R\mathrm{Hom}(\Dtwgs \Xi \mathcal{F},\Xi c_{\psi^{-1}}A)^{\vee}\\
 & \cong R\mathrm{Hom}(\mathbf{D}_{\mathrm{tw.GS}}\mathcal{F},c_{\psi^{-1}}A)^{\vee}\\
 & \cong R\mathrm{Hom}(a_{\psi^{-1}}\mathbf{D}_{\mathrm{tw.GS}}\mathcal{F},A)^{\vee}\\
 & \cong R\mathrm{Hom}(\mathbf{D}_{\mathrm{BZ}}a_{\psi}\mathcal{F},A)^{\vee}\\
 & \cong R\mathrm{Hom}(a_{\psi}\mathcal{F},\mathbf{D}_{\mathrm{Verd}}A)\\
 & \cong R\mathrm{Hom}(\mathcal{F},c_{\psi}\mathbf{D}_{\mathrm{Verd}}A)
\end{align*}
using the duality exchange formula on both the automorphic and spectral side together with Theorem \ref{thm:apsiduality}. Since $\mathcal{F}$ is arbitrary, the result follows by Yoneda. In passing from the second to the third line, we also used that $\mathbf{D}_{\mathrm{GS}}$
is compatible with $\Xi$ on perfect complexes, since $\omega_{\mathrm{Par}_{G}}\simeq\mathcal{O}_{\mathrm{Par}_{G}}$.
\end{proof}

\begin{thm}\label{thm:tpsiadjunction}
For any $A\in \D(\mathrm{Bun}_{G})$ and any $\mathcal{F}\in\mathrm{Coh}(\mathrm{Par}_{G})_{\fin}$,
there is a natural isomorphism
\[
R\mathrm{Hom}(c_{\psi}A,\Psi \mathcal{F})\cong R\mathrm{Hom}(A,t_{\psi}\mathcal{F}).
\]
\end{thm}

\begin{proof}
We compute as follows:
\begin{align*}
R\mathrm{Hom}(A,t_{\psi}\mathcal{F}) & =R\mathrm{Hom}(A,\mathbf{D}_{\mathrm{Verd}}a_{\psi^{-1}}\Psi\mathbf{D}_{\mathrm{tw.adm}}\mathcal{F})\\
 & \cong R\mathrm{Hom}(a_{\psi^{-1}}\Psi\mathbf{D}_{\mathrm{tw.adm}}\mathcal{F},\mathbf{D}_{\mathrm{Verd}}A)\\
 & \cong R\mathrm{Hom}(\Psi \mathbf{D}_{\mathrm{tw.adm}}\mathcal{F},c_{\psi^{-1}}\mathbf{D}_{\mathrm{Verd}}A)\\
 & \cong R\mathrm{Hom}(\Psi \mathbf{D}_{\mathrm{tw.adm}}\mathcal{F}, \Psi\mathbf{D}_{\mathrm{tw.adm}}\Xi c_{\psi}A)\\
 & \cong R\mathrm{Hom}(\Xi \Psi \mathbf{D}_{\mathrm{tw.adm}}\mathcal{F}, \mathbf{D}_{\mathrm{tw.adm}}\Xi c_{\psi}A)\\ 
 & \cong R\mathrm{Hom}(\mathbf{D}_{\mathrm{tw.adm}}\mathcal{F}, \mathbf{D}_{\mathrm{tw.adm}}\Xi c_{\psi}A)\\
 & \cong R\mathrm{Hom}(\Xi c_{\psi}A,\mathbf{D}_{\mathrm{tw.adm}}^{2}\mathcal{F})\\
 & \cong R\mathrm{Hom}(\Xi c_{\psi}A,\mathcal{F})\\
 & \cong R\mathrm{Hom}(c_{\psi}A,\Psi \mathcal{F}).
\end{align*}
Here, aside from the obvious adjunctions, we used Theorem \ref{thm:cpsiVerdier} in the fourth line, and biduality of $\mathbf{D}_{\mathrm{tw.adm}}$
on admissible sheaves in the eighth line. Most crucially, in the sixth line we used the fact that $\mathbf{D}_{\mathrm{tw.adm}}\mathcal{F}$ is in the essential
image of the fully faithful functor $\Xi$ by Theorem \ref{thm:cohfinclassification}, so $\mathbf{D}_{\mathrm{tw.adm}}\mathcal{F} = \Xi\Psi\mathbf{D}_{\mathrm{tw.adm}}\mathcal{F}$.
\end{proof}
Note that if $A \in \D(\Bun_G)$ is any object with $c_\psi(A) \in \Coh^\qc(\Par_G)$, we can formally promote the isomorphism of the previous theorem to an isomorphism \[
R\mathrm{Hom}(c_{\psi}A, \mathcal{F})\cong R\mathrm{Hom}(A,t_{\psi}\mathcal{F})
\]
functorial in $A$ and $\mathcal{F}$, using the fact that $\Psi : \IndCoh \twoheadrightarrow \QCoh$ induces the identity functor between the two copies of $\Coh^\qc$. In fact, we have the following result.

\begin{prop}\label{prop:tpartialright}If $\mathbf{L}_{\psi}$ exists (in the sense that the equivalent conditions of Proposition \ref{prop:Lpsiexistence} are true), we automatically get an isomorphism
\[
R\mathrm{Hom}(\mathbf{L}_{\psi}A, \mathcal{F})\cong R\mathrm{Hom}(A,t_{\psi}\mathcal{F})
\]
functorially in all $A \in \D(\Bun_G)$ and all $\mathcal{F} \in \Coh(\Par_G)_{\fin}$.
\end{prop}
In other words, if $\mathbf{L}_{\psi}$ exists, the functor $t_{\psi}$ gives an \textbf{explicit} partially defined right adjoint for $\mathbf{L}_{\psi}$ on the subcategory $\Coh(\Par_G)_{\fin} \subset \IndCoh(\Par_G)$.

\begin{proof}Quite generally, we have
\[\RHom(\mathcal{G},\mathcal{F}) = \RHom(\Psi \mathcal{G}, \Psi \mathcal{F})\]
for all $\mathcal{G} \in \IndCoh$ and all $\mathcal{F} \in \Coh^\qc$. Applying this in the situation at hand, we get
\begin{align*}
    \RHom(\mathbf{L}_{\psi}A, \mathcal{F}) & = \RHom(\Psi \mathbf{L}_{\psi}A, \Psi \mathcal{F}) \\
    & = \RHom(c_{\psi} A, \Psi \mathcal{F}) \\
    & = \RHom(A, t_{\psi} \mathcal{F}).
\end{align*}
Here in the second line we used the equation $\Psi \mathbf{L}_{\psi} = c_{\psi}$, and in the final line we used the previous theorem.
\end{proof}
\begin{cor}\label{cor:tpsiULA}
If $\mathbf{L}_{\psi}$ exists,
then $t_{\psi}$ carries $\mathrm{Coh}(\mathrm{Par}_{G})_{\fin}$
towards ULA sheaves.
\end{cor}
\begin{proof}
For any $\mathcal{F}\in\mathrm{Coh}(\mathrm{Par}_{G})_{\fin}$
and any $A\in \D(\mathrm{Bun}_{G})^{\omega}$, we compute that
\[
R\mathrm{Hom}(\mathbf{L}_{\psi}A,\mathcal{F})\cong R\mathrm{Hom}(A,t_{\psi}\mathcal{F})
\]
is a perfect complex by the definition of admissible ind-coherent
sheaves, the classification of Theorem \ref{thm:cohfinclassification}, and the fact that $\mathbf{L}_{\psi}$ carries compact sheaves towards coherent sheaves. But this exactly characterizes the ULA property of $t_{\psi}(\mathcal{F})$.
\end{proof}

We end this section by stating a form of categorical local Langlands
with restricted variation.
\begin{conjecture}
The functor $t_\psi$ induces an equivalence of categories
\[
t_{\psi}:\Coh(\Par_G)_{\fin} \overset{\sim}{\to} \D(\Bun_G)_{\fin}.
\]
\end{conjecture}
The interplay between this restricted form of categorical local Langlands and the full conjecture will be examined in the next chapter.

\section{The strategy for well-understood groups}
In this chapter we formulate a strategy to prove the categorical local Langlands conjecture for certain quasisplit groups $G$. Very roughly, our strategy can be divided into four idealized steps:
\begin{enumerate}
    \item Construct the Langlands functor $\mathbf{L}_{\psi}$ and prove that it has good properties.
    \item Construct a right adjoint $\mathbf{R}_{\psi}$ to the Langlands functor and prove that it has good properties.
    \item Prove that $\mathbf{R}_{\psi}$ is fully faithful.
    \item Prove that $\mathbf{L}_{\psi}$ is conservative.
\end{enumerate}

Note that (3) and (4) together imply that $\mathbf{L}_{\psi}$ is an equivalence of categories, on account of the following simple lemma.
\begin{lem}\label{lem:trivialadjointcondition}Let $F:\mathcal{C} \rightleftarrows \mathcal{D}:G$ be an adjoint pair of exact functors between stable $\infty$-categories. If $F$ (resp. $G$) is fully faithful and $G$ (resp. $F$) is conservative, then $F$ and $G$ are mutually inverse equivalences of categories.
\end{lem}

In our opinion, carrying out any of these steps for arbitrary groups is completely out of reach. As such, we will restrict our attention to \emph{well-understood} groups in the sense of Definition \ref{def:wellunderstoodG}. Intuitively, these are groups for which the Fargues-Scholze construction of $L$-parameters can be provably lifted to a true local Langlands parametrization with good properties. This condition implies some basic finiteness properties of the Fargues-Scholze parametrization: notably, any well-understood group is reasonable (Definition \ref{def:reasonable}), so the basic functors $a_\psi$ and $c_\psi$ have good properties.  We also prove that two naturally defined orthogonal decompositions of $\D(\Bun_G)$ into ``cuspidal and Eisenstein parts'' agree for well-understood groups, cf. Theorem \ref{thm:wellunderstoodSameOrthogonal}.

We will also assume the compatibility of the functor $c_{\psi}$ with Eisenstein series, as formulated in Conjecture \ref{conj:cpsiEiscompatible}. To us it seems completely impossible to make any significant progress on categorical local Langlands without this compatibility.  In any case, we will only need this compatibility as a black box. Moreover, the analogue of Conjecture \ref{conj:cpsiEiscompatible} in classical geometric Langlands is a well-known theorem, and we expect to settle this conjecture (for all groups!) by an adaptation of their arguments in future joint work with Linus Hamann.

In any case, we are able to execute steps (1) and (2) above for all well-understood groups. The construction of $\mathbf{L}_{\psi}$ crucially uses the agreement of the two orthogonal decompositions mentioned above. Notably, we also prove that $\mathbf{R}_{\psi}$ preserves compact objects, which is very important for the later steps of the argument. This is not at all formal and relies on most of our main results on the spectral side.

For step (3), we obtain a significant reduction: we show that full faithfulness of $\mathbf{R}_{\psi}$ is implied by full faithfulness of $a_{\psi}$, which by Proposition \ref{prop:fullyfaithfuleasy} simply amounts to the concrete statement that $c_{\psi}(i_{1!} W_{\psi})$ is a line bundle. The proof of this reduction crucially uses the functor $t_{\psi}$ and the spectral temperization Theorem \ref{rslt:spectral-temperization}, and some ideas inspired by ``geometric Langlands with restricted variation.'' 

Step (4) is the most subtle and least well-understood. In an ideal world, we would have a conceptual proof that $c_{\psi}$ is conservative on compact objects, uniformly in all groups. In the setting of classical geometric Langlands, this conservativity theorem is a recent breakthrough of Faergeman-Raskin. However, their microlocal methods are unavailable to us and a direct adaptation of their proof seems out of reach. Instead, we first show that $\mathbf{L}_{\psi}$ is conservative on sheaves on the individual strata of $\Bun_G$, assuming inductively that CLLC is known for proper Levis in $G$ to deal with the nonbasic strata. We then use a weak semiorthogonality result on the spectral side to upgrade this to the full conservativity of $\mathbf{L}_{\psi}$ on compact objects. Here we rely crucially on the ongoing work of Bertoloni Meli-Koshikawa \cite{BMK}.

The capstone result in this section is Theorem \ref{thm:DreamInductionOnLevisTheorem}. In \cref{sec:generallinear}, we will use this theorem to prove categorical local Langlands for $\mathrm{GL}_n$. In \cref{sec:classical}, we will use this theorem to reduce categorical local Langlands for many classical groups to a concrete pair of auxiliary conjectures.

\subsection{Well-understood groups}\label{ss:wellunderstood}
Let $G$ be a quasisplit group with a fixed Borel and maximal torus $T \subset B$, so we may speak of standard Levi subgroups. We also fix a Whittaker datum $(B,\psi)$, which induces a Whittaker datum for all standard Levi subgroups.

Let $\Phi(G)$ be the set of isomorphism classes of Frobenius-semisimple $L$-parameters valued in the $L$-group of $G$. Given any $\phi \in \Phi(G)$, set $S_\phi^{\natural} = S_{\phi} / (S_{\phi} \cap \hat{G}_{\mathrm{der}})^{\circ}$. There is a natural map $Z(\hat{G})^{W_F} \to S_{\phi}^{\natural}$ which by restriction induces a map
\[ \mathrm{Irr}(S_{\phi}^{\natural}) \to X^{\ast}(Z(\hat{G})^{W_F}) \cong B(G)_{\mathrm{bas}}.\]
We write $\mathrm{Irr}(S_{\phi}^{\natural},b)$ for the fiber of this map over any basic $b$.

\begin{defn}\label{def:wellunderstoodG}The group $G$ is \emph{well-understood} if there is a $B(G)_{\mathrm{bas}}$ local Langlands correspondence
\begin{align*}
\coprod_{b\in B(G)_{\mathrm{bas}}} \Pi(G_b) & \to \Phi(G) \\
 \pi & \mapsto \phi_{\pi}
\end{align*}
for $G$ \emph{and} all its standard Levi subgroups, with the following properties.

\begin{propenum}
\item For each $b \in B(G)_{\mathrm{bas}}$, the induced map $\Pi(G_b) \to \Phi(G)$ has finite fiber $\Pi_{\phi}(G_b)$ over every parameter $\phi$. If $\phi$ is essentially tempered, the chosen Whittaker datum induces a bijection $\Pi_{\phi}(G_b) \cong \mathrm{Irr}(S_{\phi}^{\natural},b)$.
\item A parameter $\phi$ is discrete if and only if $\coprod_b \Pi_{\phi}(G_b)$ consists entirely of discrete series representations, and the discrete series $L$-packets $\Pi_{\phi}(G_b)$ satisfy the endoscopic character identities as formulated in \cite[Conjectures F-G]{KalethaNQS}.

\item Every discrete series representation of some $G_b$ occurs in the fiber over a discrete parameter $\phi$.

\item A parameter $\phi$ is supercuspidal if and only if $\coprod_{b} \Pi_{\phi}(G_b)$ consists entirely of supercuspidal representations. 

\item For any $\pi \in \Pi(G_b)$, there is an equality $\phi_{\pi}^{\mathrm{ss}} = \varphi_{\pi}$, where $\varphi_{\pi}$ is the Fargues-Scholze parameter of $\pi$.
\end{propenum}
\end{defn}

Informally, this definition says the local Langlands correspondence is in good shape for all $B(G)$-twists of $G$, and agrees up to semisimplification with the Fargues-Scholze construction of $L$-parameters. Note that tautologically, all standard Levi subgroups in a well-understood group are well-understood.

\begin{rmk}\label{rmk:WellUnderstoodRoster}
The following groups are known to be well-understood. In each of these examples we fix a base field $F$.

\begin{itemize}
\item $\mathrm{GL}_n$. Here the local Langlands parametrization is due to Harris-Taylor and Henniart for the split form \cite{HarrisTaylor, Henniart} and Deligne-Kazhdan-Vigneras and Badulescu for the inner forms \cite{DKV, Bad1, Bad2}. Compatibility with Fargues-Scholze was shown in \cite{FS} for the split form and in \cite{HKW} for general inner forms.

\item $\mathrm{GSp}_4$. Here the local Langlands correspondence is due to Gan-Takeda \cite{GanTakeda} for the split form and Gan-Tantono \cite{GanTantono} for the unique inner form, and the endoscopic character identities are due to Gan-Chan \cite{ChanGan}. The compatibility with Fargues-Scholze was proved by Hamann \cite{Hamann} under a mild condition, and then by Daniels-van Hoften-Kim-Zhang in general \cite{DvHKZ}.

\item $\mathrm{SO}_{2n+1}$. We recall that the split form of $\mathrm{SO}_{2n+1}$ is in fact the unique quasisplit form, and it admits a unique inner form. Here the local Langlands parametrization is due to Arthur \cite{ArthurBook} for the split form and Ishimoto \cite{Ishimoto} for the unique inner form, with full verification of properties (i)-(iv) due to Moeglin \cite{moeglin2002, moeglin2007, moeglin2014}. The compatibility with Fargues-Scholze was proved by Peng \cite{Peng} under a mild condition, and then by Daniels-van Hoften-Kim-Zhang in general \cite{DvHKZ}.

\item $\mathrm{U}_n$. We recall that quasisplit unitary groups $\mathrm{U}_n = \mathrm{U}_n (E/F)$ are parametrized by quadratic extensions $E/F$, and that $\mathrm{U}_n(E/F)$ has a unique inner form when $n$ is even and no inner forms when $n$ is odd. Here the local Langlands parametrization is due to Mok \cite{MokUnitary} for the quasisplit form and Kaletha-Minguez-Shin-White \cite{KMSW} for inner forms, with full verification of the remaining parts of properties (i)-(iv) due to Moeglin and Chen-Zou \cite{CZortho2, Moeglinmult, Moeglinpac}. The compatibility with Fargues-Scholze was proved by Bertoloni Meli-Hamann-Nguyen for $n$ odd and $F=\mathbf{Q}_p$ \cite{BMHN}, by Peng \cite{Peng} under a mild condition, and by Daniels-van Hoften-Kim-Zhang in general \cite{DvHKZ}.
\end{itemize}

\end{rmk}

\begin{rmk} \label{rmk:typeDsubtlety} There are groups - for instance, $\mathrm{Sp}_{2n}$ with $n\geq 3$ - such that conditions (i)-(iv) are known, but where (v) is out of reach with current technology. On the other hand, the quasisplit groups of type $\mathrm{SO}_{2n}$ are nearly well-understood in the sense that all conditions are met except the endoscopic character identities are only known up to some outer automorphism ambiguity. However, we only use the endoscopic character identities in one place, namely the proof of Theorem \ref{thm:wellunderstoodSameOrthogonal} where we appeal to the main results of \cite{HKW}, and in type $\mathrm{SO}_{2n}$ we can simply use \cite[Theorem 7.1.2]{Peng} instead. In particular, all of the results in this chapter hold for type $D$ classical groups as well.
\end{rmk}

For well-understood groups, we have substantial control over the basic finiteness properties of the Fargues-Scholze correspondence.
\begin{prop} \label{rslt:well-understood-implies-reasonable}
If $G$ is well-understood, then $G$ is reasonable in the sense of Definition \ref{def:reasonable}. Moreover, the Fargues-Scholze map 
\begin{align*}
 \Pi(G_b) & \to \Phi^{\mathrm{ss}}(G) = X_{G}^{\spec}(\Qellbar) \\
 \pi & \mapsto \varphi_{\pi}
\end{align*}
has finite fibers for all $b$, and (equivalently) the natural map $X_{G_b} \to X_{G}^{\spec}$ has finite fibers for all $b$.
\end{prop}
\begin{proof}
If $G$ is well-understood, then all the maps in the sequence \[\Pi(G_b) \overset{\pi \mapsto \phi_\pi}{\longrightarrow} \Phi(G_b) \to \Phi^{\mathrm{ss}}(G_b) \to \Phi^{\mathrm{ss}}(G)\]
have finite fibers. For the first map this is condition (i), and for the second and third maps this is general nonsense. Since the composite map agrees with the Fargues-Scholze map by condition (v), we get the result.
\end{proof}

More substantially, for well-understood groups we can compare the two natural orthogonal decompositions of $\D(\Bun_G)$. To formulate this result, recall from \cite[Theorem 1.3.2]{HHS} that the theory of geometric Eisenstein series induces a natural orthogonal decomposition
\[\D(\Bun_G) = \D(\Bun_G)^{\mathrm{Eis}} \oplus \D(\Bun_G)^{\mathrm{cusp}}.\]
Here $\D(\Bun_G)^{\mathrm{Eis}}$ is the full subcategory generated under colimits by the essential images of the Eisenstein functors $\Eis_{P!}:\D(\Bun_M) \to \D(\Bun_G)$ for all proper parabolics $P=MU \subsetneq G$. Similarly $\D(\Bun_G)^{\mathrm{cusp}}$ is the full subcategory given as the common kernel of all constant term functors $\CT_{P\ast}:\D(\Bun_G) \to \D(\Bun_M)$ for all proper parabolics.\footnote{In classical geometric Langlands, there is an identically defined pair of subcategories inside $D-\mathrm{mod}(\Bun_G)$, but they only form a semiorthogonal decomposition in that setting. See the remarks after \cite[Theorem 1.3.2]{HHS} for some further discussion.}

On the other hand, writing $\Par_G = \Par_G^{\mathrm{irr}} \coprod \Par_G^{\mathrm{red}}$ as the disjoint union of the open-closed loci corresponding to irreducible (i.e., supercuspidal) resp. reducible $L$-parameters, we get a pair of complementary idempotents $e^{\mathrm{irr}},e^{\mathrm{red}} \in \mathcal{O}(\Par_G) = \mathcal{O}(X_{G}^{\spec})$. The action of the spectral Bernstein center on $\D(\Bun_G)$ then induces an orthogonal decomposition
\[\D(\Bun_G) = e^{\mathrm{red}}\D(\Bun_G) \oplus e^{\mathrm{red}}\D(\Bun_G). \]
By some soft compatibility of geometric Eisenstein series with the action of the spectral Bernstein center, it is not difficult to show the dual inclusions 
\[e^{\mathrm{irr}} \D(\Bun_G) \subseteq \D(\Bun_G)^{\mathrm{cusp}}\]
and
\[e^{\mathrm{red}} \D(\Bun_G) \supseteq \D(\Bun_G)^{\mathrm{Eis}}\]
for all quasisplit $G$ (see \cite{Takaya} for a detailed proof). 

It is very natural to guess that these two orthogonal decompositions agree for all quasisplit $G$. In fact, this expectation is forced on us by the categorical local Langlands conjecture together with the expected Eisenstein compatibility of the Langlands functor. Unfortunately, proving the agreement of these decompositions for arbitrary groups seems completely out of reach.
\begin{thm}\label{thm:wellunderstoodSameOrthogonal}If $G$ is well-understood, then 
\[e^{\mathrm{irr}} \D(\Bun_G) = \D(\Bun_G)^{\mathrm{cusp}}\]
and
\[e^{\mathrm{red}} \D(\Bun_G) = \D(\Bun_G)^{\mathrm{Eis}}\]
as full subcategories of $\D(\Bun_G)$.
\end{thm}
\begin{proof}
By the remarks preceding the theorem, it's enough to show the inclusion $e^{\mathrm{red}} \D(\Bun_G) \subseteq \D(\Bun_G)^{\mathrm{Eis}}$. Suppose given a sheaf $A \in \D(\Bun_G)$ with $e^{\mathrm{red}}A = A$. Quite generally, if the support of some sheaf $A'$ is disjoint from the basic locus, then $A'$ lies in $\D(\Bun_G)^{\mathrm{Eis}}$. Since the excision triangles are compatible with the action of $e^{\mathrm{red}}$, we can thus assume $A$ is of the form $i_{b!}B$ for some basic $b$ and some $B\in \D(G_b(F),\Qellbar)$. Splitting $B$ into direct summands indexed by Bernstein components, we can assume $B$ is supported on one component. If that component is parabolically induced, there is nothing to prove. Thus we can assume $B$ is supported on a supercuspidal Bernstein component $\mathfrak{s}$ such that the $L$-parameters of the irreducible representations in this component are \emph{reducible}. Letting $\sigma$ be any irreducible supercuspidal in the component $\mathfrak{s}$, the universal unramified twist $\sigma \otimes \chi^{\mathrm{univ}}$ is a compact projective generator for the component $\mathfrak{s}$. By some further easy compatibilities of everything with twisting, we're reduced to proving that if $b$ is basic and $\sigma$ is an irreducible supercuspidal $G_b(F)$-representation with reducible Fargues-Scholze parameter, then $i_{b!}\sigma \in \D(\Bun_G)^{\mathrm{Eis}}$.

By conditions (ii)-(v), the $L$-parameter $\phi:=\phi_{\sigma}$ is discrete but non-supercuspidal, and we may choose some other basic $b'$ together with a discrete \emph{non-supercuspidal} irreducible representation $\pi \in \Pi_{\phi}(G_{b'})$. By the main theorem of \cite{HKW}, we can now choose some $V\in \mathrm{Rep} \hat{G}$ such that there is a nonzero morphism $f:i_{b!}\sigma \to T_{V} i_{b'!}\pi$.\footnote{This is the only place we make use of the endoscopic character identities.} Since $\pi$ is non-supercuspidal, it lives in a parabolically induced Bernstein block, so $i_{b'!}\pi \in \D(\Bun_G)^{\Eis}$. By \cite{HHS2}, $\D(\Bun_G)^{\Eis}$ is stable under Hecke operators, so also $T_V i_{b'!}\pi \in \D(\Bun_G)^{\mathrm{Eis}}$. As the cuspidal and Eisenstein subcategories form an orthogonal decomposition, the existence of the nonzero map $f$ then implies that $i_{b!}\sigma \in \D(\Bun_G)^{\mathrm{Eis}}$, as desired.
\end{proof}

\subsection{The running assumptions}
Throughout the remainder of this entire paper, \emph{we assume unless explicitly stated otherwise} that $G$ is well-understood and that $c_{\psi}$ is compatible with Eisenstein series in the sense of Conjecture \ref{conj:cpsiEiscompatible}.

\subsection{Construction of the Langlands functor}

\begin{thm} \label{rslt:c-psi-sends-compact-to-coherent}
Let $G$ be a well-understood group. Then the functor $c_\psi$ carries $\D(\Bun_G)^\omega$ into $\Coh^\qc(\Par_G)$, so we may define the Langlands functor 
\[\mathbf{L}_\psi:\D(\Bun_G) \to \IndCoh(\Par_G)\]
by ind-completing.
\end{thm}
\begin{proof}
Let $A \in \D(\Bun_G)^\omega$ be any compact sheaf. By the cuspidal-Eisenstein decomposition, we can write $A= A^{\cusp} \oplus A^{\Eis}$. By Theorem \ref{thm:wellunderstoodSameOrthogonal}, $A^{\cusp}= e^{\mathrm{irr}}A$ is supported on the supercuspidal locus, so $c_{\psi}(A^{\cusp})$ is coherent with quasicompact support by Theorem \ref{thm:cpsisupercuspidalcoherent}. On the other hand, $A^{\Eis}$ can be obtained by finitely many colimits and retracts from compact sheaves of the form $\Eis_{P^- !}(B)$, where $P=MU \subsetneq G$ is a proper standard parabolic and $B \in \D(\Bun_M)^{\omega}$ is compact. Now clearly
\[c_{\psi}(\Eis_{P^- !}(B)) \simeq \Eis_{P}^{\spec,\coarse}(c_{\psi_{M}}(B))\]
by Eisenstein compatibility. But then $c_{\psi_{M}}(B) \in \Coh^\qc(\Par_M)$ by induction on the semisimple rank, and $\Eis_{P}^{\spec,\coarse}$ preserves coherent sheaves with quasicompact support, so $c_{\psi}(\Eis_{P^- !}(B))$ is coherent with quasicompact support.
\end{proof}
Note that by construction, $\mathbf{L}_{\psi}$ is colimit-preserving and $\QCoh(\Par_G)$-linear, preserves compact objects, and satisfies $c_\psi = \Psi \circ \mathbf{L}_{\psi}$. We also get the following compatibility for free.

\begin{prop}\label{prop:wellunderstoodLpsiEis}The functor $\mathbf{L}_{\psi}$ is compatible with Eisenstein series, in the sense that there is a commutative diagram\[
\xymatrix{\D(\mathrm{Bun}_{M})\ar[r]^{\mathbf{L}_{\psi_{M}}}\ar[d]^{\mathrm{Eis}_{P^{-}!}} & \IndCoh(\mathrm{Par}_{M})\ar[d]^{\Eis_{P}^{\spec}}\\
\D(\mathrm{Bun}_{G})\ar[r]^{\mathbf{L}_{\psi}} & \IndCoh(\mathrm{Par}_{G})
}
\]
for all standard parabolics $P \subset G$.
\end{prop}
\begin{proof}By construction, we have $\mathbf{L}_{\psi} \circ \Eis_{P^- !} \simeq \Eis_{P}^{\mathrm{spec}} \circ \mathbf{L}_{\psi_M}$ on compact objects of $\D(\Bun_M)$, so the result follows formally by ind-completing.
\end{proof}

\subsection{The right adjoint}

Note that $\mathbf{L}_\psi$ commutes with colimits and preserves compact objects, so by the $\infty$-categorical adjoint functor theorem it admits a colimit-preserving right adjoint \[\mathbf{R}_\psi: \IndCoh(\Par_G) \to \D(\Bun_G)\]
which is automatically $\QCoh(\Par_G)$-linear.
\begin{prop}\label{prop:rightadjointbasics}The functor $\mathbf{R}_{\psi}$ carries $\Adm(\Par_G)$ into $\D(\Bun_G)^{\mathrm{ULA}}$, carries $\chi$-graded ind-coherent sheaves towards sheaves supported on the $-\chi$-component of $\Bun_G$, and is compatible with constant term functors in the sense that there is a natural equivalence
\[ \CT_{P^- \ast} \circ \mathbf{R}_{\psi} \simeq \mathbf{R}_{\psi_M} \circ \CT_{P}^{\spec}.\]
Moreover, there is a natural isomorphism \[\mathbf{R}_{\psi}|_{\Coh(\Par_G)_{\fin}} = t_{\psi}\] with $t_{\psi}$ as in Definition \ref{defn:tpsi}.
\end{prop}
\begin{proof}If $\mathcal{F}$ is admissible and $A\in \D(\Bun_G)^{\omega}$ is compact, then
\[\RHom(\mathbf{L}_{\psi}(A), \mathcal{F}) = \RHom(A,\mathbf{R}_{\psi}(\mathcal{F}))\]
is a perfect complex of vector spaces, by the definition of admissibility and the coherence of $\mathbf{L}_{\psi}(A)$. This shows that $\mathbf{R}_{\psi}(A)$ is ULA. Compatibility with the central grading follows formally by adjunction from the same compatibility for $\mathbf{L}_{\psi}$. The compatibility with constant term functors is immediate by passing to right adjoints in Proposition \ref{prop:wellunderstoodLpsiEis}.  Finally, the comparison with $t_{\psi}$ is exactly the content of Proposition \ref{prop:tpartialright}.
\end{proof}

Much less formally, we have the following result.

\begin{thm}\label{thm:rightadjointcompactness}The functor $\mathbf{R}_\psi$ preserves compact objects, i.e. it carries $\Coh^\qc(\Par_G)$ towards $\D(\Bun_G)^\omega$.
\end{thm}

\begin{cor}\label{cor:tpsifin}The functor $t_{\psi}$ carries $\Coh(\Par_G)_{\fin}$ towards $\D(\Bun_G)_{\fin}$.
\end{cor}

\begin{proof}[Proof of Corollary \ref{cor:tpsifin}] By Corollary \ref{cor:tpsiULA}, we know that $t_{\psi}$ carries $\Coh(\Par_G)_{\fin}$ towards ULA sheaves. Combining Theorem \ref{thm:rightadjointcompactness} and Proposition \ref{prop:rightadjointbasics}, we also see that $t_{\psi}$ carries $\Coh(\Par_G)_{\fin}$ towards compact sheaves. Since \[\D(\Bun_G)_{\fin} = \D(\Bun_G)^{\omega} \cap \D(\Bun_G)^{\mathrm{ULA}},\] this gives the claim.
\end{proof}

Thus we now have access to an adjoint pair of functors $\mathbf{L}_\psi$ and $\mathbf{R}_\psi$ between $\D(\Bun_G)$ and $\IndCoh(\Par_G)$, both colimit-preserving and $\QCoh(\Par_G)$-linear and \emph{both} preserving compact objects.

For the proof of Theorem \ref{thm:rightadjointcompactness}, we will need the following compactness criterion for sheaves on $\Bun_G$.
\begin{thm}\label{thm:compactcriterion}Let $G$ be a well-understood group. For any sheaf $A\in \D(\Bun_G)$, the following are equivalent.

\begin{enumerate}
    \item $A$ is compact.

 \item $A$ has quasicompact support and quasicompact spectral support, and for all $B \in \D(\Bun_G)^\omega$, $\RHom(B,A)$ is a bounded complex whose cohomology groups are finitely generated $\mathfrak{Z}_{G}^{\mathrm{spec}}$-modules.
\end{enumerate}
\end{thm}

We note that the well-understood condition is overkill for this theorem, and the following argument works under the weaker condition that the Fargues-Scholze map $\Pi(G_b) \to \Phi^{\mathrm{ss}}(G)$ has finite fibers for all $b$.

\begin{proof}First we prove that 2) implies 1). It is clear from the hypotheses in 2) and the basic properties of well-understood groups that for each fixed $b$, $\Pi=i_{b}^{\ast}A$ is supported on finitely many Bernstein components for $G_b$, so its enough to see that for each Bernstein component $\mathfrak{s}$, the corresponding summand $\Pi_{\mathfrak{s}}$ is compact. Let $\Sigma_{\mathfrak{s}}$ be Bernstein's canonical compact projective generator of the component, so the functor $\RHom_{G_b(F)}(\Sigma_{\mathfrak{s}},-)$ is an exact equivalence from the component $\D(G_b(F))_\mathfrak{s}$ towards the derived category of right $\mathcal{H}_{\mathfrak{s}}:=\mathrm{End}(\Sigma_{\mathfrak{s}})$-modules, with essential inverse $M \mapsto M\otimes_{\mathcal{H}_{\mathfrak{s}}} \Sigma_{\mathfrak{s}}$. Since $\mathcal{H}_{\mathfrak{s}}$ is noetherian of finite global dimension, it is easy to see that $\Pi_{\mathfrak{s}}$ is compact if and only if $\RHom(\Sigma_{\mathfrak{s}},\Pi_{\mathfrak{s}})$ is a bounded complex whose cohomology groups are finitely generated $\mathcal{H}_{\mathfrak{s}}$-modules. But the hypotheses in 2) guarantee that \[\RHom(\Sigma_{\mathfrak{s}},\Pi_{\mathfrak{s}})\simeq \RHom(i_{b\sharp}\Sigma_{\mathfrak{s}}, A)\] is a bounded complex whose cohomology groups are finitely generated $\mathfrak{Z}_{G}^{\mathrm{spec}}$-modules, and this module structure is compatible with the $\mathcal{H}_{\mathfrak{s}}$-module structure via the canonical maps
\[\mathfrak{Z}_{G}^{\mathrm{spec}} \to \mathcal{Z}(\mathcal{H}_{\mathfrak{s}}) \to \mathcal{H}_{\mathfrak{s}}.\]
Therefore the cohomology groups are finitely generated $\mathcal{H}_{\mathfrak{s}}$-modules, which gives the desired implication.

For the reverse implication, it is clear that if $A$ is compact then it is has quasicompact support and quasicompact spectral support. Now pick any compact $B \in \D(\Bun_G)$. We need to see that $\RHom(B,A)$ is a bounded complex whose cohomology groups are finitely generated finitely generated $\mathfrak{Z}_{G}^{\mathrm{spec}}$-modules.  We first choose a finite filtration on $B$ whose graded pieces are $\sharp$-pushfowards of compact objects on strata. This reduces us to checking that if $\Pi, \Pi' \in \D(G_b(F))$ are compact, then $\RHom(\Pi,\Pi')$ is a bounded complex whose cohomology groups are finitely generated $\mathfrak{Z}_{G}^{\mathrm{spec}}$-modules. Arguing as in the previous step, we may assume $\Pi$ and $\Pi'$ lie in a single Bernstein component $\mathfrak{s}$, in which case they each can be obtained from the compact projective generator $\Sigma_{\mathfrak{s}}$ by finitely many shifts, cones and retracts. Filtering $\RHom(\Pi,\Pi')$ further via this observation, it is finally enough to see that $\mathcal{H}_\mathfrak{s}=\mathrm{End}(\Sigma_\mathfrak{s})$ is a finitely generated $\mathfrak{Z}_{G}^{\mathrm{spec}}$-module. But $\mathcal{H}_\mathfrak{s}$ is finitely generated as a module over its center $\mathcal{Z}(\mathcal{H}_\mathfrak{s})$ \cite{Bernstein} and the map $\mathfrak{Z}_{G}^{\mathrm{spec}} \to \mathcal{Z}(\mathcal{H}_\mathfrak{s})$ is module-finite by \cite[Lemma 1.6.4]{Beijing}.
\end{proof}

\begin{proof}[Proof of Theorem \ref{thm:rightadjointcompactness}] Fix any $\mathcal{F} \in \Coh^\qc(\Par_G)$. We will show that $\mathbf{R}_{\psi}(\mathcal{F})$ is compact by verifying the criterion of Theorem \ref{thm:compactcriterion}. First, it is clear that $\mathbf{R}_{\psi}(\mathcal{F})$ has quasicompact spectral support. Next, if $B \in \D(\Bun_G)^{\omega}$ is any compact object, then
\[\RHom(B,\mathbf{R}_{\psi}(\mathcal{F})) \simeq \RHom(\mathbf{L}_{\psi}(B),\mathcal{F})\]
by the obvious adjointness. Since $\mathbf{L}_{\psi}(B) \in \Coh^\qc(\Par_G)$, we deduce that $\RHom(\mathbf{L}_{\psi}(B),\mathcal{F})$ is a bounded complex whose cohomology groups are finitely generated $\mathfrak{Z}_{G}^{\mathrm{spec}}$-modules by Theorem \ref{rslt:miracle-bound}.

Finally, we need to see that $\mathbf{R}_{\psi}(\mathcal{F})$ has quasicompact support. This immediately reduces to the statement that for all but finitely many $b$, $i_{b}^{\ast \ren} \mathbf{R}_{\psi}(\mathcal{F})=0$. By the compatibility of $\mathbf{R}_{\psi}$ with the central grading, it is clear that $\mathbf{R}_{\psi}(\mathcal{F})$ is supported on finitely many connected components of $\Bun_G$, and in particular $i_{b}^{\ast}\mathbf{R}_{\psi}(\mathcal{F})$ can only be nonzero for finitely many \emph{basic} $b$. Therefore we can restrict our attention to nonbasic $b$. By the constant term compatibility of $\mathbf{R}_{\psi}$ together with compatibility of $\mathbf{R}_{\psi_M}$ with the central grading, we compute that
\[\CT_{P^- \ast}^{\chi} \mathbf{R}_{\psi}(\mathcal{F})= \mathbf{R}_{\psi_M}\CT_{P}^{\spec,-\chi}(\mathcal{F})\] for any standard parabolic $P=MU$ and any $\chi \in X^\ast (Z(\hat{M})^{W_F}) = \pi_0 (\Bun_M)$. By Theorem \ref{rslt:finiteness-for-constant-terms}, we deduce that for every proper standard parabolic $P$, $\CT_{P^- \ast}^{\chi} \mathbf{R}_{\psi}(\mathcal{F})=0$ for all but finitely many $P$-antidominant characters.

Now choose any non-basic element $b \in B(G)$. Let $P_{ \{ b \} }=M_{ \{ b \} } U_{ \{ b \}} $ be the repelling dynamic parabolic of the Newton point of $b$, and let $\chi_{b} \in X^{\ast}(Z(\widehat{M_{ \{b \}}})^{W_F}) $ be the Kottwitz point of the obvious basic $M_{ \{b \} }$-reduction of $b$. Replacing $b$ by a $\sigma$-conjugate representative if necessary, we can assume that $P_{ \{ b \} }$ is standard. Then $i_{b \sharp}^{\ren}= \Eis_{P_{ \{ b \} }^{-}!} i_{b!}^{ M_{ \{ b \} }}$ by \cite[Corollary 2.2.5.(5)]{HHS}, so passing to right adjoints gives $i_{b}^{\ast \ren} = i_{b}^{ \ast M_{ \{ b \} } } \CT_{P_{ \{ b \} }^{-} \ast}$. Since $i_{b}^{ \ast M_{ \{ b \} } }$ vanishes on sheaves outside the $\chi_{b}$-component of $\Bun_{ M_{ \{ b \}} }$, this expression is unchanged by passing to the $\chi_{b}$-graded part of the constant term.
In particular, we see that $i_{b}^{\ast \ren} \mathbf{R}_{\psi}(\mathcal{F})$ is a stalk of $\CT_{P_{ \{b \} }^- \ast}^{\chi_b} \mathbf{R}_{\psi}(\mathcal{F})$. Since $\chi_{b}$ is $P_{ \{ b \} }$-antidominant, the argument in the previous paragraph shows that $\CT_{P_{ \{b \} }^- \ast}^{\chi_b} \mathbf{R}_{\psi}(\mathcal{F})=0$ for all but finitely many $b$. More precisely, we know from the previous paragraph that for any (proper standard) $P$, $\CT_{P^- \ast}^{\chi} \mathbf{R}_{\psi}(\mathcal{F})=0$ for all but finitely many $P$-antidominant $\chi$; but the map from the set of nonbasic $b$ towards such pairs $(P,\chi)$ sending $b$ to $(P_{ \{b \} }, \chi_b)$ is injective. Since $i_{b}^{\ast \ren} \mathbf{R}_{\psi}(\mathcal{F})$ is a stalk of $\CT_{P_{ \{b \} }^- \ast}^{\chi_b} \mathbf{R}_{\psi}(\mathcal{F})$, this implies that $i_{b}^{\ast \ren} \mathbf{R}_{\psi}(\mathcal{F})=0$ for all but finitely many $b$, as desired.
\end{proof}

\subsection{With or without restricted variation}
In this section we prove the following theorem.

\begin{thm}\label{thm:restrvsfull}The following statements are equivalent.

\begin{enumerate}
    \item The functor $\mathbf{L}_{\psi}: \D(\Bun_G) \to \IndCoh(\Par_G)$ is an equivalence.

\item The functor $c_\psi:\D(\Bun_G)^\omega \to \Coh^\qc(\Par_G)$ is an equivalence.

\item The functor $c_\psi: \D(\Bun_G)_{\fin} \to \Coh(\Par_G)_{\fin}$ is an equivalence.

\item The functor $t_\psi: \Coh(\Par_G)_{\fin} \to \D(\Bun_G)_{\fin}$ is an equivalence.
\end{enumerate}
\end{thm}
Here the implications (1) $\Leftrightarrow$ (2) $\Rightarrow$ (3) $\Leftrightarrow$ (4) are straightforward consequences of our results so far. However, the implication (4) $\Rightarrow$ (1) is much more interesting.\footnote{A superficially similar result is proved in \cite{GLC1}, but the theorem proved there assumes the analogue of (4) for all coefficient fields $k/\Qellbar$ (see \cite[\S 5.1.9]{GLC1}), which makes the argument much easier.}

For the proof, recall that in Sections \ref{ss:finrestrcoh} and \ref{ss:restrBunGside} we defined full subcategories of ``restricted'' objects
\[\IndCoh(\Par_G)^{\restr} \subset \IndCoh(\Par_G)\]
and
\[\D(\Bun_G)^{\restr} \subset \D(\Bun_G)\]
such that the displayed inclusions admit a right adjoint given by the action of $\mathcal{O}^{\restr}$.
Moreover, by Proposition \ref{prop:restrsheavesParG} and Theorem \ref{thm:restrsheavesBunG}, we have identifications
\[\IndCoh(\Par_G)^{\restr} = \Ind(\Coh(\Par_G)_{\fin})\]
and
\[\D(\Bun_G)^{\restr} = \D(\Bun_G)_{\mathrm{indfin}}.\]

Now we can exploit these observations in combination with the basic fact that $\mathbf{L}_{\psi}$ and $\mathbf{R}_{\psi}$ are colimit-preserving and $\QCoh(\Par_G)$-linear. In light of the above results, we immediately deduce that they restrict to an adjoint pair of functors $\mathbf{L}_{\psi}: \D(\Bun_G)_{\mathrm{indfin}} \rightleftarrows \IndCoh(\Par_G)^{\restr}: \mathbf{R}_{\psi}$ such that the diagrams
\[
\xymatrix{\D(\mathrm{Bun}_{G})_{\mathrm{indfin}}\ar@<-.5ex>[d]\ar[r]^{\mathbf{L}_{\psi}} & \mathrm{IndCoh}(\mathrm{Par}_{G})^{\restr}\ar@<-.5ex>[d]\\
\D(\mathrm{Bun}_{G})\ar@<-.5ex>[u]\ar[r]^{\mathbf{L}_{\psi}} & \mathrm{IndCoh}(\mathrm{Par}_{G})\ar@<-.5ex>[u]
}
\]
and
\[
\xymatrix{\D(\mathrm{Bun}_{G})_{\mathrm{indfin}}\ar@<-.5ex>[d] & \mathrm{IndCoh}(\mathrm{Par}_{G})^{\restr}\ar@<-.5ex>[d]\ar[l]_{\mathbf{R}_{\psi}}\\
\D(\mathrm{Bun}_{G})\ar@<-.5ex>[u] & \mathrm{IndCoh}(\mathrm{Par}_{G})\ar@<-.5ex>[u]\ar[l]_{\mathbf{R}_{\psi}}
}
\]
commute. Here the downwards arrows are the evident full embeddings, and the upwards arrows are the projection functors given by the action of $\mathcal{O}^{\restr}$.

With these preparations, we can finish the proof of Theorem \ref{thm:restrvsfull}. Assume that $t_\psi$ is an equivalence as in Theorem \ref{thm:restrvsfull}.(4). Using that $\mathbf{R}_{\psi}|_{\Coh(\Par_G)_{\fin}} = t_{\psi}$ and passing to ind-completions, this implies that \[\mathbf{R}_{\psi}:\Ind(\Coh(\Par_G)_{\fin}) \cong \IndCoh(\Par_{G})^{\restr} \to \D(\Bun_G)_{\mathrm{indfin}}\]
is an equivalence of categories. Therefore the adjoint functors \[\mathbf{L}_{\psi}: \D(\Bun_G)_{\mathrm{indfin}} \rightleftarrows \IndCoh(\Par_G)^{\restr}: \mathbf{R}_{\psi}\] are mutually inverse equivalences.

Now, since $\mathbf{L}_{\psi}$ and $\mathbf{R}_{\psi}$ preserve compact objects and both the automorphic and spectral categories are compactly generated, it's enough to see that for any $A \in \D(\Bun_G)^\omega$ resp. $\mathcal{F} \in \Coh^\qc(\Par_G)$, the adjunction maps $A \to \mathbf{R}_{\psi} \mathbf{L}_{\psi}A$ resp. $\mathbf{L}_{\psi} \mathbf{R}_{\psi}\mathcal{F} \to \mathcal{F}$ are isomorphisms. Let $C$ resp. $\mathcal{K}$ be the cone of the first resp. second map. Then $C$ and $\mathcal{K}$ are compact objects, and since $\mathbf{L}_{\psi}$ and $\mathbf{R}_{\psi}$ are equivalences on the categories with restricted variation, we deduce from the commutativity of the diagrams above that $\mathcal{O}^{\restr} \ast C =0$ and $\mathcal{O}^{\restr} \otimes \mathcal{K}=0$. We now conclude by the following lemma.

\begin{lem}\label{lem:conservativetrick} Let $G$ be any reductive group.
\begin{lemenum}
    \item The endofunctor $\mathcal{O}^{\restr} \ast - \circlearrowright \D(\Bun_G)$ is conservative on $\D(\Bun_G)^\omega$.
    \item The endofunctor $\mathcal{O}^{\restr} \otimes - \circlearrowright \IndCoh(\Par_G)$ is conservative on $\Coh^\qc(\Par_G)$.
\end{lemenum}
\end{lem}
\begin{proof}
Using Theorem \ref{thm:restrsheavesBunG} and unwinding definitions, (1) reduces to the claim that if $A \in \D(\Bun_G)$ is a nonzero compact object, then there is a nonzero map $B \to A$ for some $B \in \D(\Bun_G)_{\fin}$. Indeed, writing $i:\D(\Bun_G)_{\mathrm{indfin}} \to \D(\Bun_G)$ for the obvious full inclusion and $R$ for its right adjoint, we have
\begin{align*}
    \RHom(B,A) = & \RHom(i(R(B)),A) \\
     & = \RHom(R(B),R(A)) \\
     & = \RHom(i(R(B)),i(R(A))) \\
     & = \RHom(B,i(R(A))) \\
     & = \RHom(B, \mathcal{O}^{\restr} \ast A)
\end{align*}
where we used the full faithfulness of $i$ several times, along with the identification $i(R(-)) = \mathcal{O}^{\restr} \ast -$ and the fact that $B$ is already in the essential image of $i$. A nonzero map $B \to A$ with $B$ finite then obviously induces a nonzero map $B\to \mathcal{O}^{\restr} \ast A $, which clearly implies that $\mathcal{O}^{\restr} \ast A \neq 0$.

Since $\Dbz$ induces compatible anti-equivalences of compact and finite sheaves on $\Bun_G$, dualizing shows that the existence of the desired map $B \to A$ is equivalent to the claim that for any nonzero compact sheaf $A$, there is a nonzero map $A \to B$ for some $B \in \D(\Bun_G)_{\fin}$. Fix such an $A$, let $b$ be a maximally special point in the support of $A$, let $n$ be the maximal degree for which $H^n(i_{b}^{\ast}A)\neq 0$, and let $s:H^n(i_{b}^{\ast}A) \twoheadrightarrow \pi$ be any surjective map towards some nonzero irreducible $G_b(F)$-representation. The existence of such a surjection $s$ follows from the finite generation of $H^n(i_{b}^{\ast}A)$. Indeed, $H^n(i_{b}^{\ast}A)$ is admissible over the Bernstein center $\mathfrak{Z}_{G_b}$, so choosing any maximal ideal $\mathfrak{m} \subset \mathfrak{Z}_{G_b}$ in the support of $H^n(i_{b}^{\ast}A)$, the (nonzero) quotient $H^n(i_{b}^{\ast}A) / \mathfrak{m}$ is finitely generated and admissible in the absolute sense, hence of finite length, so we may take $\pi$ to be any irreducible quotient of this representation. Then the induced composite map 
\[A \to i_{b!}i_{b}^{\ast}A \to i_{b!} \tau^{\geq n}i_{b}^{\ast}A \to i_{b!}\pi[-n]\] is nonzero, and the target is a finite sheaf.

Similarly, for (2) one shows that 
\[\RHom(\mathcal{G},\mathcal{F}) = \RHom(\mathcal{G},\mathcal{O}^{\restr} \otimes \mathcal{F})\]
for all $\mathcal{G} \in \Coh(\Par_G)_{\fin}$ and all $\mathcal{F} \in \Coh^\qc(\Par_G)$. The argument is formally identical to the case of $\Bun_G$, so we omit it. Therefore it is enough to show that for any nonzero coherent sheaf $\mathcal{F}$, there is a nonzero map $\mathcal{F}'\to\mathcal{F}$ with $\mathcal{F}' \in \Coh(\Par_G)_{\fin}$. We will prove more generally that if $X$ is a QCA stack and $\{i_t: Z_t \subset X \}_{t \in T}$ is a collection of closed substacks whose union contains all closed points of $|X|$, then for any nonzero $\mathcal{F} \in \Coh(X)$ there is a nonzero map $i_{t \ast} \mathcal{G} \to \mathcal{F}$ for some $t\in T$ and some $\mathcal{G} \in \Coh(Z_t)$. Applying this general claim to the case where $X=\Par_G$, $T$ is the set of closed points in $X_G^{\spec}$, and $\{ Z_t \}$ is the set of closed fibers of the map $q$, we get our desired result. 

Now, since Grothendieck-Serre duality preserves coherence and commutes with pushforward along closed immersions, producing a nonzero map $i_{t \ast} \mathcal{G} \to \mathcal{F}$ as above is equivalent (after replacing the sheaves by their duals) to producing a nonzero map $\mathcal{F}\to i_{t \ast} \mathcal{G}$, or equivalently a nonzero map $f:i_{t}^{\ast} \mathcal{F} \to \mathcal{G}$ where the pullback is the usual pullback in $\QCoh$. Since $\mathcal{F}$ is coherent, $\mathrm{supp}\,\mathcal{F}$ is closed in $|X|$, so by our assumption on closed points we may choose some $t$ such that $Z_t$ meets the support of $\mathcal{F}$, in which case $i_{t}^{\ast} \mathcal{F} \neq 0$. Now we can simply take $\mathcal{G} = \tau^{\geq -N} i_{t}^{\ast}\mathcal{F}$ for some sufficiently large $N$, with $f:i_{t}^{\ast} \mathcal{F} \to \tau^{\geq -N} i_{t}^{\ast}\mathcal{F}$ the canonical map.
\end{proof}

\subsection{Variants on full faithfulness}
In this section we prove the following two results.
\begin{thm}\label{thm:atff}If $a_\psi$ is fully faithful, then $t_\psi$ is fully faithful.
\end{thm}
\begin{thm}\label{thm:tRff}If $t_\psi$ is fully faithful, then $\mathbf{R}_\psi$ is fully faithful.
\end{thm}
Of course, as an immediate consequence we get the following corollary. While this corollary is probably the statement of primary interest, we emphasize that we do not know an argument for it which avoids the functor $t_\psi$!
\begin{cor}\label{cor:apsiRpsiffmagic}If $a_\psi$ is fully faithful, then $\mathbf{R}_\psi$ is fully faithful.
\end{cor}

\begin{proof}[Proof of Theorem \ref{thm:atff}] By definition, $t_{\psi}$ is the composition
\[
\mathrm{Coh}(\mathrm{Par}_{G})_{\fin}\overset{\mathbf{D}_{\mathrm{tw.adm}}}{\longrightarrow}\mathrm{QCoh}(\mathrm{Par}_{G})\overset{a_{\psi^{-1}}}{\longrightarrow}\D(\mathrm{Bun}_{G})\overset{\mathbf{D}_{\mathrm{Verd}}}{\longrightarrow}\D(\mathrm{Bun}_{G})
\]
where the first arrow comes from Theorem \ref{thm:cohfinclassification}.
The first two functors in this sequence are fully faithful, so the
only possible problem is the final $\mathbf{D}_{\mathrm{Verd}}$.
However, by Corollary \ref{cor:tpsifin}, the total functor $t_{\psi}$ has image contained in finite sheaves, thus in ULA sheaves. By Proposition \ref{prop:ULAreflect}, $\mathbf{D}_{\mathrm{Verd}}$ reflects
the ULA property, so the partial composition
\[
\mathrm{Coh}(\mathrm{Par}_{G})_{\fin} \overset{\mathbf{D}_{\mathrm{tw.adm}}}{\longrightarrow}\mathrm{QCoh}(\mathrm{Par}_{G})\overset{a_{\psi^{-1}}}{\longrightarrow}\D(\mathrm{Bun}_{G})
\]
must already have its image contained in ULA sheaves. Since this partial
composition is fully faithful, and $\mathbf{D}_{\mathrm{Verd}}$ is
fully faithful on ULA sheaves, we get the desired result.
\end{proof}
\begin{proof}[Proof of Theorem \ref{thm:tRff}] Since $\mathbf{R}_{\psi}|_{\Coh(\Par_G)_{\fin}} = t_\psi$, we immediately get that $\mathbf{R}_\psi$ is fully faithful on $\Coh(\Par_G)_{\fin}$ and sends it towards $\D(\Bun_G)_{\fin}$. Since finite sheaves are compact and $\mathbf{R}_{\psi}$ preserves compact objects, passing to ind-categories we then get that $\mathbf{R}_\psi$ is fully faithful on $\Ind(\Coh(\Par_G)_{\fin}) \cong \IndCoh(\Par_{G})^{\mathrm{restr}}$. Now, since $\mathbf{R}_{\psi}$ preserves compact objects, it is enough to see that it is fully faithful on $\Coh^\qc(\Par_G)$. Picking any $\mathcal{F} \in \Coh^\qc(\Par_G)$, we need to see that the adjunction map
\[\mathbf{L}_{\psi} \mathbf{R}_\psi \mathcal{F} \to \mathcal{F} \] is an isomorphism. Let $\mathcal{K}$ be the cone of this map. By the finiteness theorems proved in the previous sections, we know that $\mathcal{K}$ is coherent. Moreover, arguing as in the proof of Theorem \ref{thm:restrvsfull}, we get that $\mathcal{O}^{\restr} \otimes \mathcal{K}$ is the cone of the adjunction 
\[\mathbf{L}_{\psi} \mathbf{R}_\psi (\mathcal{O}^{\restr} \otimes\mathcal{F}) \to \mathcal{O}^{\restr} \otimes\mathcal{F}, \]
and thus $\mathcal{O}^{\restr} \otimes \mathcal{K}=0$ since $\mathcal{O}^{\restr} \otimes \mathcal{F} \in \IndCoh(\Par_{G})^{\mathrm{restr}}$ and $\mathbf{R}_\psi$ is fully faithful on this category by the analysis above. But $\mathcal{O}^{\restr} \otimes -$ is conservative on $\Coh^\qc(\Par_G)$ by Lemma \ref{lem:conservativetrick}, so we deduce that $\mathcal{K}=0$ as desired.
\end{proof}
\begin{rmk}\label{rmk:aRcomparison}If $a_\psi$ is fully faithful, it is easy to see that $\mathbf{L}_\psi a_\psi$ coincides with the canonical embedding $\Xi:\QCoh(\Par_G)\to \IndCoh(\Par_G)$. By adjunction we then get a canonical natural transformation
\[a_\psi \to \mathbf{R}_\psi \mathbf{L}_\psi a_\psi \cong \mathbf{R}_{\psi} \Xi.\]
It is natural to conjecture that this transformation $a_\psi \to \mathbf{R}_{\psi} \Xi$ is an equivalence, but this doesn't seem to be a formal consequence of our assumptions. By $\QCoh(\Par_G)$-linearity of all functors in sight, this reduces to showing the natural map
\[i_{1!}W_{\psi} \to \mathbf{R}_{\psi} \mathcal{O}_{\Par_G}\]
obtained by evaluating on the structure sheaf is an isomorphism. If categorical local Langlands is true, it is easy to see that this is the case, since this map becomes an isomorphism after applying the conservative functor $\mathbf{L}_{\psi}$. However, much less is needed: we will see in the next section that $\mathbf{R}_{\psi} \mathcal{O}_{\Par_G}$ is automatically of the form $i_{1!} B$ for some sheaf $B$, whence it is enough to know that $\mathbf{L}_{\psi} \circ i_{1!}$ is conservative.
\end{rmk}

\subsection{Basic and cuspidal sheaves}

In this section we single out a natural collection of sheaves on the spectral side which should match the sheaves on the automorphic side which are $!$-extended from the basic locus in $\Bun_G$, along with a smaller collection corresponding to $!$-extensions of supercuspidals from the basic locus. The latter definition is due entirely to Bertoloni Meli-Koshikawa.

Until further notice we allow any quasisplit $G$, with a fixed Borel and torus as usual. 

\begin{lem}\label{lem:nonbasiccoincidence} Let $\Lambda$ be any coefficient ring. The following two collections of sheaves generate the same subcategory of $\D(\Bun_G,\Lambda)$ under colimits.

\begin{enumerate}
    \item Sheaves of the form $i_{b\sharp}A$ for non-basic $b$ and $A\in \D(G_b(F),\Lambda)^\omega$.

    \item Sheaves of the form $\Eis_{P^-! }(A)$ where $P$ is a proper standard parabolic and $A\in \D(\Bun_{M}^{\chi},\Lambda)^\omega$ is supported on a component whose slope vector $\chi \in B(M)_{\mathrm{bas}} \cong X^{\ast}(Z(\hat{M})^{W_F})$ is strictly $P$-antidominant.
\end{enumerate}

\end{lem}
Here again we warn the reader of the sign change in passing from isocrystals to bundles.
\begin{proof}It is clear that every sheaf of type (1) is of type (2). The more interesting statement is that every sheaf of type (2) has a finite filtration whose graded peices are sheaves of type (1). Moreover, noting that pairs $(P,\chi)$ as in (2) are in canonical bijection with the nonbasic elements in $B(G)$, we will see that the graded pieces of this filtration are all of the form $i_{b'\sharp}A'$ for $b'$ with $b_{P,\chi} \preceq b'$.

To show this, keep the notation as in (2), but consider first the Eisenstein series $\Eis_{P!}(A)$ for the opposite parabolic $P$. By Remark 3.5.(3) in \cite{Viehmann}, $\Eis_{P!}(A)$ is supported only at points which are specializations of $b_{P,\chi}$ in $B(G)$. Since it is moreover compact, we deduce that it has a finite filtration with graded pieces of the form $i_{b'!}A'$ for some compact $A'$'s and some elements $b'$ such that $b_{P,\chi}\preceq b'$. But now, replacing $A$ with its Bernstein-Zelevinsky dual, we can apply second adjointness \cite[Theorem 1.2.1(3)]{HHS} in the form $\Eis_{P!} \mathbf{D}_{\mathrm{BZ}}^{M} \simeq \Dbz \Eis_{P^- !}$. Using that Bernstein-Zelevinsky exchanges $!$- and $\sharp$-pushforwards, this gives the desired result.
\end{proof}
\begin{defn}\label{defn:basiccoherent} (1) A sheaf $\mathcal{F}\in \IndCoh(\Par_G)$ is \emph{basic} if $\mathrm{CT}_{P}^{\mathrm{spec},\chi}(\mathcal{F})=0$ for all proper standard parabolics $P$ and all strictly $P$-dominant $\chi \in X^{\ast}(Z(\hat{M})^{W_F})$. We write $\IndCoh(\Par_G)_{\mathrm{bas}}$ for the full subcategory of basic sheaves.

(2) (Bertoloni Meli-Koshikawa)  A sheaf $\mathcal{F}\in \IndCoh(\Par_G)$ is \emph{cuspidal} if $\mathrm{CT}_{P}^{\mathrm{spec},\chi}(\mathcal{F})=0$ for all proper standard parabolics $P$ and all $P$-dominant $\chi \in X^{\ast}(Z(\hat{M})^{W_F})$. We write $\IndCoh(\Par_G)_{\mathrm{cusp}}$ for the full subcategory of cuspidal sheaves.
\end{defn}

It is clear from this definition that cuspidal sheaves are basic, but not conversely. Note also that basic and cuspidal sheaves are still graded by characters of $Z(\hat{G})^{W_F}$ in the usual way.

\begin{prop}\label{prop:basicisright}  Assume that $G$ is well-understood. Let $\mathcal{F}\in \IndCoh(\Par_G)$ be any sheaf.
\begin{enumerate}
    \item If $\mathcal{F}$ is basic, $\mathbf{R}_{\psi}(\mathcal{F})$ is $!$-extended from the basic locus in $\Bun_G$. 
    
    \item If $\mathcal{F}$ is cuspidal, $\mathbf{R}_{\psi}(\mathcal{F})$ is $!$-extended from the basic locus and $i_{b}^{\ast}\mathbf{R}_{\psi}(\mathcal{F})$ lies in a union of supercuspidal Bernstein components for every basic $b$. 
    
    \item If $\mathbf{L}_{\psi_M}$ is an equivalence of categories for all standard Levis $M \subseteq G$, the converses to (1)-(2) are true.
\end{enumerate}
\end{prop}
\begin{proof}
Quite generally, a sheaf $A\in \D(\Bun_G)$ is $!$-extended from the basic locus if and only if it lies in the right orthogonal of the subcategory generated by sheaves of the form $i_{b\sharp}B$ for non-basic $b$ and arbitrary $B$. By Lemma \ref{lem:nonbasiccoincidence} and the adjunction between Eisenstein series and constant terms, this is equivalent to the condition that $\mathrm{CT}_{P^- \ast}^{\chi}(A)=0$ for all proper parabolics $P$ and all strictly $P$-antidominant $\chi \in X^{\ast}(Z(\hat{M})^{W_F})$. Using the compatibility of $\mathbf{R}_{\psi}$ with constant terms (Proposition \ref{prop:rightadjointbasics}), we see that if $\mathcal{F}$ is basic, then
\[\mathrm{CT}_{P^- \ast}^{\chi}(\mathbf{R}_{\psi}(\mathcal{F})) = \mathbf{R}_{\psi_M}(\mathrm{CT}_{P}^{\mathrm{spec},-\chi}(\mathcal{F}))=0\]
upon noting that $-\chi$ is $P$-dominant, so $\mathbf{R}_{\psi}(\mathcal{F})$ is $!$-extended from the basic locus. This gives (1). If moreover $\mathcal{F}$ is cuspidal, then for every basic $b$ the Jacquet modules of $i_{b}^{\ast} \mathbf{R}_{\psi}(\mathcal{F})$ vanish for all proper parabolics $P \subsetneq G_b$, which gives (2). The final converse statement is clear.
\end{proof}

Although we will not need it, we sketch a proof of the following instructive result.

\begin{thm}\label{thm:structuresheafbasic} The structure sheaf $\mathcal{O}_{\Par_G}$ is basic.
\end{thm}
\begin{proof}Unwinding the definition of basic sheaves, this amounts to the claim that for any parabolic $P$ with associated map $q:\Par_P \to \Par_M$, the sheaf $q_{\ast}^{\IndCoh}\mathcal{O}_{\Par_P}$ has vanishing $\chi$-graded summand for all strictly $P$-dominant $\chi$.

This follows from a stronger quantitative statement: if $\chi$ occurs in the grading of $q_{\ast}^{\IndCoh}\mathcal{O}_{\Par_P}$, then $-\chi$ is a nonegative linear combination of characters occurring in the $Z(\hat{M})^{W_F}$-representation $\Lie(\hat{U})$. This in turn follows from the contraction principle.
\end{proof}

With this language in place, we want to state a conjecture purely on the spectral side. Here we freely use the sets $\mathcal{P}$ and $\mathcal{P}^+$ and the attendant notation introduced in Section \ref{sec:nonbasic}. If $M$ is a standard Levi and $b\in B(M)_{\mathrm{bas}}$ is any element, we write $\chi_b \in X^{\ast}(Z(\hat{M})^{W_F})$ for the associated character.

\begin{conjecture}\label{conj:spectralsemiorthweak}Fix the notation as above.
\begin{enumerate}
    \item For any $(M_1,b_1) \in \mathcal{P}$ and $(M_2,b_2) \in \mathcal{P}$ such that $r(M_2,b_2) \npreceq r(M_1,b_1)$, we have
    \[\RHom(\Eis_{M_2 U_0}^{\mathrm{spec}}(\mathcal{F}_2), \Eis_{M_1 U_0}^{\mathrm{spec}}(\mathcal{F}_1))=0\]
    for all $\mathcal{F}_2 \in \IndCoh(\Par_{M_2})^{\chi_{b_2}}$ and all $\mathcal{F}_1 \in \IndCoh(\Par_{M_1})^{\chi_{b_1}}_{\mathrm{bas}}$.

    \item   For any $(M_1,b_1) \in \mathcal{P}$ and $(M_2,b_2) \in \mathcal{P}$ such that $r(M_2,b_2) \npreceq r(M_1,b_1)$, we have
    \[\RHom(\Eis_{M_1 U_0}^{\mathrm{spec}}(\mathcal{F}_1), \Eis_{M_2 U_0}^{\mathrm{spec}}(\mathcal{F}_2))=0\]
    for all $\mathcal{F}_2 \in \IndCoh(\Par_{M_2})^{-\chi_{b_2}}$ and all $\mathcal{F}_1 \in \IndCoh(\Par_{M_1})^{-\chi_{b_1}}_{\mathrm{bas}}$.

    \item For any $(M,b) \in \mathcal{P}^+$, $\Eis_{P}^{\mathrm{spec}}$ is fully faithful on $\Coh^\qc(\Par_{M})^{\chi_{b}}_{\mathrm{bas}}$.
    
\end{enumerate}
\end{conjecture}

Under categorical local Langlands, part (1) resp. (2) corresponds to (a slight generalization of) the geometrically evident fact that $\RHom(i_{b_2 \sharp}-,i_{b_1 \sharp} -)=0$ resp. $\RHom(i_{b_1 !}-, i_{b_2 !} -)=0$ whenever $b_1$ is not a specialization of $b_2$. Part (3) corresponds to the full faithfulness of $i_{b\sharp}$. Note that (1)-(2) resp. (3) only has content when $P_2$ resp. $P$ is a proper parabolic. Moreover, parts (1) and (2) would be equivalent if we knew that basic sheaves are preserved by Grothendieck-Serre duality; this preservation is natural to expect, as it follows from categorical local Langlands and Eisenstein/CT compatibility using Theorem \ref{thm:automaticduality}.

Granted this conjecture, one would be able to construct a full $B(G)$-indexed semiorthogonal decomposition on $\IndCoh(\Par_G)$. Unfortunately, the above conjecture seems extremely difficult to attack directly. Very luckily, however, the following more restrictive variant will be proved in forthcoming work of Bertoloni Meli-Koshikawa, and this is enough for our purposes.

\begin{thm}[Bertoloni Meli-Koshikawa] \label{thm:BMKsemiorthogonalKey} Let $\mathcal{P}$ be as in Section \ref{sec:nonbasic}.
\begin{enumerate}
    \item For any $(M,b^M) \in \mathcal{P}$ and $(L,b^L) \in \mathcal{P}$ such that $r(L,b^L) \npreceq r(M, b^M)$, we have
    \[\RHom(\Eis_{MU_0}^{\mathrm{spec}}(\mathcal{F}_M), \Eis_{LU_0}^{\mathrm{spec}}(\mathcal{F}_L))=0\]
    for all $\mathcal{F}_L \in \Coh(\Par_{L})^{-\chi_{b^L}}_{\mathrm{cusp}}$ and all $\mathcal{F}_M \in \Coh(\Par_{M})^{-\chi_{b^M}}_{\mathrm{cusp}}$.

    \item For any $(M,b^M) \in \mathcal{P}$ and any  $\mathcal{F}_M \in \Coh(\Par_{M})^{\pm\chi_{b^M}}_{\mathrm{cusp}}$ with $\mathcal{F}_M \neq 0$, then also $\Eis_{MU_0}^{\mathrm{spec}}(\mathcal{F}_M)\neq 0$.
    \end{enumerate}
\end{thm}

\subsection{An induction principle}
 
In  this section we put everything together. Let $G$ be a well-understood group, equipped with a fixed Whittaker datum $(B,\psi)$ and maximal torus $T \subset B$ as usual. Then all standard Levis $M \subset G$ inherit their own Whittaker data $\psi_M$, with respect to which they are also well-understood. 

The following theorem informally says that if we know categorical Langlands for all proper Levis in $G$, and we understand the functor $\mathbf{L}_{\psi} i_{1!}$ sufficiently well, then categorical Langlands for $G$ follows.

\begin{thm}\label{thm:DreamInductionOnLevisTheorem}Let $G$ be as above. Assume the following conditions.

\begin{enumerate}

    \item For each proper standard Levi $M \subsetneq G$, the functor $\mathbf{L}_{\psi_M}$ is an equivalence of categories. 

    \item The functor $\mathbf{L}_{\psi} i_{1!}$ sends $W_{\psi}$ to a line bundle.

    \item The functor $\mathbf{L}_{\psi} i_{1!}$ is conservative on compact supercuspidal representations and sends them towards cuspidal coherent sheaves.
\end{enumerate}
Then $\mathbf{L}_{\psi}$ is an equivalence of categories.
\end{thm}

It is immaterial in (2)-(3) whether we argue with $\mathbf{L}_{\psi}$ or $c_{\psi}$. We also note that in (3) it would be enough to ask for conservativity on finite-length supercuspidal representations.

\begin{proof}

We already know that both categories are compactly generated and both functors are colimit-preserving and preserve compact objects. It's therefore enough to show that $\mathbf{R}_{\psi}$ is fully faithful and that $\mathbf{L}_{\psi}$ is conservative on compact objects. From (2), we get that $a_{\psi}$ is fully faithful by Proposition \ref{prop:fullyfaithfuleasy}, so then $\mathbf{R}_{\psi}$ is fully faithful by Corollary \ref{cor:apsiRpsiffmagic}. 

The conservativity of $\mathbf{L}_{\psi}$ on compact objects is trickier, and will be handled in several steps. Ideally, we would like to use the easy observation that if $F:\mathcal{C}\to \mathcal{D}$ is any functor between categories with semiorthogonal decompositions $\mathcal{C}=\left\langle \mathcal{C}_i \right\rangle_{i \in I}$ resp. $\mathcal{D}=\left\langle \mathcal{D}_i \right\rangle_{i \in I}$ such that $F$ is compatible with the decompositions and $F|_{\mathcal{C}_i}$ is conservative for all $i$, then $F$ is conservative. Since we don't have access to a full $B(G)$-indexed semiorthogonal decomposition of $\Coh^\qc(\Par_G)$, we have to argue much more indirectly. Fortunately, Theorem \ref{thm:BMKsemiorthogonalKey} together with the assumptions (1) and (3) give us just enough leverage to push through an argument in this style.

First we note that the general case reduces to the special case of compact sheaves supported on the neutral component of $\Bun_G$. For this, suppose $A \neq 0$ is a compact sheaf supported on the $\chi$-component of $\Bun_G$ for some $\chi \in X^\ast (Z (\hat{G})^{W_F})$. Pick any $\hat{G}$-representation $V$ whose restriction to $Z (\hat{G})^{W_F}$ is $\chi$-isotypic. Then $T_V A$ is compact and supported on the neutral connected component, and $T_V A \neq 0$ by the conservativity of Hecke operators.\footnote{The identity functor is a summand of $T_{V^\vee} T_V \cong T_{V^{\vee} \otimes V}$ for any $V$.} Therefore if $\mathbf{L}_{\psi}$ is conservative for sheaves on the neutral component, we trivially get $0 \neq \mathbf{L}_{\psi} (T_V A)$, but since $\mathbf{L}_{\psi} (T_V A) = V \otimes \mathbf{L}_{\psi}A$, this implies that $\mathbf{L}_{\psi}A \neq 0$ as desired.

This reduces us to showing the conservativity of $\mathbf{L}_{\psi}$ on compact sheaves on the neutral connected component of $\Bun_G$. We next observe that our hypotheses imply the conservativity of the functors $\mathbf{L}_{\psi} i_{b!}^{\ren}$ on compact objects for all $b$ with $\kappa(b)=0$, i.e. for all points of $\Bun_G$ in the neutral connected component. For $b=1$ on supercuspidal blocks, this is exactly the first half of (3). For $b=1$ on non-supercuspidal blocks, we will give the argument at the very end of the proof. For non-basic $b$, let $P$ be the dynamic parabolic of the Newton point of $b$, with Levi $M$. Then $i_{b!}^{\ren} = \Eis_{P^- !} i_{b!}^{M}$, so 
\begin{align*}
    \mathbf{L}_{\psi} i_{b!}^{\ren} & = \mathbf{L}_{\psi} \Eis_{P^- !} i_{b!}^{M} \\
    & = \Eis_{P}^{\spec} \mathbf{L}_{\psi_M} i_{b!}^{M}.
\end{align*}
Now by (1), $\mathbf{L}_{\psi_M} i_{b!}^{M}$ is fully faithful and carries supercuspidal $G_b(F)$-representations towards cuspidal coherent sheaves, so the composite functor $\Eis_{P}^{\spec} \mathbf{L}_{\psi_M} i_{b!}^{M}$ is conservative on supercuspidal $G_b(F)$-representations by Theorem \ref{thm:BMKsemiorthogonalKey}.(2). For non-basic $b$ on non-supercuspidal blocks, a small further argument is required, which we will give at the end of the proof.

Now let $A$ be any nonzero compact sheaf supported on the neutral component of $\Bun_G$. We need to see that $\mathbf{L}_{\psi}A \neq 0$. Let $b$ be a maximally special point in the support of $A$, so we get a distinguished triangle
\[B \to A \to i_{b!}i_{b}^{\ast} A \to \]
whose third term is nonzero, and such that no point in the support of $B$ is a specialization of $b$. Consider the resulting triangle
\[\mathbf{L}_{\psi} B \to \mathbf{L}_{\psi} A \to \mathbf{L}_{\psi} i_{b!}i_{b}^{\ast} A \to.\]
By the analysis in the previous paragraph, the third term is nonzero. If $\mathbf{L}_{\psi} B=0$, then $\mathbf{L}_{\psi} A \simeq \mathbf{L}_{\psi} i_{b!}i_{b}^{\ast} A \neq 0$ and we're done. Now assume that $\mathbf{L}_{\psi} B \neq 0$. If $\mathbf{L}_{\psi} A =0$, we would get an isomorphism $\mathbf{L}_{\psi} B \simeq \mathbf{L}_{\psi} i_{b!}i_{b}^{\ast} A[-1]$, so in particular
\[ \RHom (\mathbf{L}_{\psi} B, \mathbf{L}_{\psi} i_{b!}i_{b}^{\ast} A) \neq 0.\] Filtering $B$ by the finitely many $!$-pushfowards $i_{b_i!} i_{b_i}^{\ast}B$ with $b_i \in \mathrm{supp}B$, it is clear that one of the groups
\[ \RHom (\mathbf{L}_{\psi} i_{b_i!} i_{b_i}^{\ast}B, \mathbf{L}_{\psi} i_{b!}i_{b}^{\ast} A)\]
must then be nonzero. Up to finite colimits and retracts, we now apply Corollary \ref{cor:Eisreducetosupercuspidal} to further filter each $i_{b_i!}i_{b_i}^{\ast}B$ by terms of the form $\Eis_{LU_0 ^- !} i_{b^L !} B_L$, where $(L,b^L)$ runs over the finite set of pairs $(L,b^L) \in \mathcal{P}$ retracting to $b_i \in B(G) = \mathcal{P}^+$ and $B_L$ is a compact supercuspidal representation of $L_{b^L}(F)$, and analogously for $i_{b!}i_{b}^{\ast} A$.  Thus up to finite colimits and retracts, the $\RHom$ above can be filtered by finitely many groups of the form
\[\RHom(\mathbf{L}_{\psi}\Eis_{LU_0 ^- !} i_{b^L !} B_L,\mathbf{L}_{\psi}\Eis_{MU_0 ^- !} i_{b^M !} A_{M})  \]
where $(L,b^L)$ resp. $(M,b^{M})$ run over elements of $\mathcal{P}$ retracting to some $b_i$ resp. to $b$, and $B_L$ resp. $A_{M}$ run over finitely many compact supercuspidal representations of $L_{b^L}(F)$ resp. $M_{b^M}(F)$. By Eisenstein compatibility, we can rewrite this as
\[\RHom(\Eis_{LU_0 }^{\spec} \mathbf{L}_{\psi_L}i_{b^L !} B_L,\Eis_{MU_0}^{\spec}\mathbf{L}_{\psi_M} i_{b^M !} A_{M}). \] Now combining (1) and the second half of (3) with Proposition \ref{prop:basicisright}, we deduce that the sheaves $\mathbf{L}_{\psi_L}i_{b^L !} B_L$ and $\mathbf{L}_{\psi_M}i_{b^M !} A_M$ are cuspidal and graded for the evident characters $-\chi_{b^L}$ resp. $-\chi_{b^M}$.
But no $b_i$ appearing is a specialization of $b$, so the nonvanishing of any of these $\RHom$'s would violate Theorem \ref{thm:BMKsemiorthogonalKey}.(1). Therefore $\mathbf{L}_{\psi} A \neq 0$.

It remains to show that $\mathbf{L}_{\psi} i_{b!}^{\ren}$ is conservative on compact $G_b(F)$-representations supported on non-supercuspidal Bernstein blocks, for all $b$ with $\kappa(b)=0$. Our argument for this is inspired by the proof of \cite[Corollary IX.7.3]{FS}. We first treat the case where $b=1$.  Fix $\Pi$ a compact $G(F)$-representation supported in a non-supercuspidal Bernstein block, and choose some proper standard parabolic $P=MU$ such that $r_{G}^{P^-}\Pi \neq 0$. Shrinking $P$ if necessary, we can assume $j_{G}^{P^-}\Pi$ is supercuspidal. Now choose some rational cocharacter $\mu: \mathbf{G}_m \to G$ with dynamic parabolic $P$, and set $b=\mu(\pi)$. We can assume that $\kappa(b)=0$, e.g. by requiring that $\mu$ factors over $G^{\mathrm{sc}}\to G$. Let $V_{\mu} \in \mathrm{Irr}(\hat{G})$ be the irreducible representation of highest weight $\mu$. Then $G_{b} = M$, $T_{V_{\mu}}i_{1!}\Pi$ is supported in $\Bun_{G}^{\preceq b}$, and there is a canonical identification $i_{b}^{\ast \ren} T_{V_{\mu}}i_{1!} \Pi = r_{G}^{P^-}\Pi$. In particular, we get a distinguished triangle
\[B \to T_{V_{\mu}}i_{1!} \Pi \to i_{b!}^{\ren}r_{G}^{P^-}\Pi \to \]
where the third term is nonzero and the first term is only supported at points which are strict generizations of $b$. By the supercuspidality of $r_{G}^{P^-}\Pi$, the argument in the fourth paragraph of the proof shows that $\mathbf{L}_{\psi}i_{b!}^{\ren}r_{G}^{P^-}\Pi\neq 0$. Now arguing as in the previous paragraph of the proof, we deduce that $\mathbf{L}_{\psi}( T_{V_{\mu}}i_{1!} \Pi) \neq 0$: otherwise, the nonvanishing of $\mathbf{L}_{\psi}i_{b!}^{\ren}r_{G}^{P^-}\Pi$ would again lead to a contradiction with Theorem \ref{thm:BMKsemiorthogonalKey}.(1), using the above-stated fact about the support of $B$. Therefore 
\[0 \neq \mathbf{L}_{\psi}( T_{V_{\mu}}i_{1!} \Pi) = V_{\mu} \otimes \mathbf{L}_{\psi}i_{1!}\Pi \]
so $\mathbf{L}_{\psi}i_{1!}\Pi \neq 0$ as desired.

Finally, we need to show the conservativity of $\mathbf{L}_{\psi} i_{b!}^{\ren}$ on non-supercuspidal blocks for non-basic $b$ with $\kappa(b)=0$. This is a variant of the argument for $b=1$ in the previous paragraph. Fix such a $b$, and fix $\Pi$ a compact $G_b(F)$-representation supported in a non-supercuspidal Bernstein block. We may choose a rational parabolic $P_b =M_b U_b \subset G_b$ which is the $b$-twist of a standard parabolic $P = MU \subset G$ such that $r_{G_b}^{P_{b}^-}\Pi$ is nonzero and supercuspidal. Pick some $\mu:\mathbf{G}_m \to G$ factoring over $G^{\mathrm{sc}} \to G$ with dynamic parabolic $P$, and set $b'=b\mu(\pi)$.  Then $\kappa(b')=0$, $G_{b'} = M_b$, $T_{V_\mu} i_{b!}^{\ren}\Pi$ is supported in $\Bun_{G}^{\preceq b'}$, and there is a canonical identification $i_{b'}^{\ast \ren} T_{V_{\mu}}i_{b!}^{\ren} \Pi = r_{G_b}^{P_{b}^-}\Pi$. In particular, we get a distinguished triangle
\[B \to T_{V_{\mu}}i_{b!}^{\ren} \Pi \to i_{b'!}^{\ren}r_{G_b}^{P_{b}^-}\Pi \to \]
where the third term is nonzero and the first term is only supported at points which are strict generizations of $b'$. Then $\mathbf{L}_{\psi}$ applied to the third term is again nonzero by the supercuspidality of $r_{G_b}^{P_{b}^-}\Pi$, and the argument proceeds exactly as in the previous paragraph.
\end{proof}

In the remainder of this paper, we will explain how to apply Theorem \ref{thm:DreamInductionOnLevisTheorem} to various families of groups, in particular to $\mathrm{GL}_n$.

\section{General linear groups}\label{sec:generallinear}

In this chapter we specialize the results of Chapter 6 to the case of $\mathrm{GL}_n$. More generally, we will sometimes consider groups of \emph{GL-type}, i.e. groups of the form $G=\mathrm{GL}_{n_{1}}\times\cdots\times\mathrm{GL}_{n_{k}}$.
This slight extra generality helps to facilitate some inductions,
but the apparent increase in scope is an illusion, since all the sheaf
categories for such a group are just Lurie tensor products of the
sheaf categories for the individual factor groups, and all the functors
we care about are compatible with this structure. Unless stated otherwise, we will always pin the group by the (product of) upper triangular Borel(s) and the standard diagonal maximal torus inside it. We also observe that for groups of this type, there is a unique isomorphism class of Whittaker data.

\begin{thm}\label{thm:CLLCforGLn}For $G$ any group of $GL$-type, the functor $\mathbf{L}_{\psi}$ is an equivalence of categories.
\end{thm}

The key feature of $\mathrm{GL}_n$, as opposed to arbitrary well-understood groups, is the work of Ben--Zvi-Chen-Helm-Nadler \cite{BCHN}. More precisely, recall that \cite{BCHN} constructed a fully faithful embedding
\[\mathscr{S}_{\GL_n}:\D(\mathrm{GL}_n(F),\Qellbar)\to \IndCoh(\Par_{\mathrm{GL}_n})\]
using their theory of the coherent Springer sheaf.\footnote{Their construction and all its properties extend trivially to groups of GL-type. We will use this without any comment.} The key connection with our program is the following theorem.
\begin{thm}\label{thm:GLnfullfaithndconsequences}Assume that $G$ is of GL-type. Then there is an isomorphism of functors \[\mathscr{S}_{G} \simeq \mathbf{L}_{\psi} \circ i_{1!}.\]
\end{thm}

\begin{proof}[Proof of Theorem \ref{thm:CLLCforGLn}] We will inductively apply Theorem \ref{thm:DreamInductionOnLevisTheorem}. By induction on the semisimple rank, condition (1) is obviously valid, using the results of \cite{ZouTori} for the base case where the group is a torus. For groups of $\GL$-type, the second half of condition (3) is trivially true. Now by the main results of \cite{BCHN} we know that $\mathscr{S}_{G}$ is conservative, in fact fully faithful, on compact objects, and that it sends $W_{\psi}$ to a line bundle (see Theorem  \ref{thm:BCHNWhittakerkey} for the latter statement on general Bernstein blocks). Then by Theorem \ref{thm:GLnfullfaithndconsequences} we get exactly the same properties for the functor $\mathbf{L}_{\psi} \circ i_{1!}$. This verifies condition (2) and the first half of condition (3), so all conditions are satisfied and Theorem \ref{thm:DreamInductionOnLevisTheorem} applies. The result follows.
\end{proof}

The proof of Theorem \ref{thm:GLnfullfaithndconsequences} will occupy Sections 7.1 and 7.2 below. In Section 7.1 we first handle the case of the principal Bernstein block, where the argument actually works uniformly for any split group. In Section 7.2 we deal with general blocks for $\mathrm{GL}_n$: with the actual idea laid out in Section 7.1, this is just a matter of notation.

\subsection{Iwahori-Hecke magic}\label{sec:IwahoriGLn}
In this section we momentarily allow $G$ to be any split group with fixed choices of $T \subset B$ as usual; we can and do choose extensions of them to smooth group schemes over $\mathcal{O}_F$. Let $W$ be the Weyl group, and let $I \subset G(F)$ be the Iwahori subgroup associated with $B$. 

Let $\D(G(F),\overline{\mathbf{Q}_{\ell}})^{\omega}_{\mathbf{1}}$ denote the compact objects in the principal block of the derived category of smooth representations. Let $\Par_{G,\mathbf{1}}$ denote the principal component of the stack of $L$-parameters for $G$.\footnote{This is a connected component of the whole stack, and it is the unique quasicompact open substack whose $\Qellbar$-points correspond to $L$-parameters whose semisimplification is trivial on the inertia $I_F$.} Note that the ring of global functions on $\mathrm{Par}_{G,\mathbf{1}}$ is canonically identified with the ring $\mathfrak{Z}=\mathcal{O}(\hat T)^W$, which also canonically identifies with the center of the category $\D(G(F),\overline{\mathbf{Q}_{\ell}})^{\omega}_{\mathbf{1}}$.

In this generality, \cite{BCHN} constructed a fully faithful functor
\[\mathscr{S}_{G,\mathbf{1}}: \D(G(F),\overline{\mathbf{Q}_{\ell}})^{\omega}_{\mathbf{1}} \to \Coh(\Par_{G,\mathbf{1}}) \subset \Coh(\Par_G)\]
whose definition we will recall below.

\begin{thm}\label{thm:splitgroupBCHNcomparison} For any split group $G$ as above, the functor
\[
c_{\psi}\circ i_{1!}:\D(G(F),\overline{\mathbf{Q}_{\ell}})^{\omega}_{\mathbf{1}} \to\mathrm{QCoh}(\mathrm{Par}_{G})
\]
is fully faithful with image in $\mathrm{Coh}(\Par_{G,\mathbf{1}})$ and agrees with the
functor $\mathscr{S}_{G,\mathbf{1}}$ constructed by Ben--Zvi-Chen-Helm-Nadler.
\end{thm}

The proof of this requires some preparation. We first recall the construction of the functor $\mathscr{S}_{G,\mathbf{1}}$. Let \[\Pi_{G,\mathbf{1}} = C_c(I \backslash G(F)) \simeq i_{B^{-}}^{G} C_c(T(F)/T(\mathcal{O}_F))\] be the canonical generator for the principal block, and let $\mathcal{H}=\mathrm{End}_{G(F)}(\Pi_{G,\mathbf{1}})$ be the Iwahori-Hecke algebra. Then $\D(G(F),\overline{\mathbf{Q}_{\ell}})^{\omega}_{\mathbf{1}}$ is canonically identified with the bounded derived category of finitely generated right $\mathcal{H}$-modules, via the functor $V \mapsto V^I$. By classic results of Bernstein, $\mathcal{H}$ is finitely generated over its center, which admits a canonical identification $\mathfrak{Z} = \mathcal{O}(\hat{T})^W$.

On the other hand, the stack of parameters for $T$ has an unramified component admitting a canonical decomposition $\mathrm{Par}_{T}^{\mathrm{nr}}\cong\hat{T}\times B\hat{T}$, and we define the coherent sheaf \[\mathcal{S}_{G,\mathbf{1}}=\mathrm{Eis}_{B}^{\mathrm{spec}}(\mathcal{O}_{\mathrm{Par}_{T}^{\mathrm{nr}}})\in\mathrm{Coh}(\mathrm{Par}_{G,\mathbf{1}}).\]
This is the coherent Springer sheaf as defined in \cite{BCHN}. The following difficult theorem summarizes the main results concerning this sheaf.

\begin{thm}[Ben--Zvi-Chen-Helm-Nadler] \label{thm:BCHNsplitGspringer} Maintain the above notation.
\begin{enumerate}
    \item For all $i \neq 0$, $\Hom(\mathcal{S}_{G,\mathbf{1}},\mathcal{S}_{G,\mathbf{1}}[i]) =0$.

    \item There is a $\mathfrak{Z}$-algebra isomorphism $\mathrm{End}(\mathcal{S}_{G,\mathbf{1}}) \simeq \mathcal{H}$. 
\end{enumerate}
\end{thm}

Note that part (1) of this theorem is a hard counterpart to the easy fact that $\Pi_{G,\mathbf{1}}$ has no self-exts in nonzero degrees.

With this theorem in hand, the functor $\mathscr{S}_{G,\mathbf{1}}$ is literally characterized by the requirement that it sends $\Pi_{G,\mathbf{1}}$ to $\mathcal{S}_{G,\mathbf{1}}$ and that the induced map $\mathrm{End}(\Pi_{G,\mathbf{1}}) \to \mathrm{End}(\mathcal{S}_{G,\mathbf{1}})$ is a $\mathfrak{Z}$-algebra isomorphism. Slightly more precisely, after fixing a choice of isomorphism as in (2) of the previous theorem, we may regard $\mathcal{S}_{G,\mathbf{1}}$ as a left $\mathcal{H}$-module, and the functor $\mathscr{S}_{G,\mathbf{1}}$ is defined by
\begin{align*}
    \D(G(F),\overline{\mathbf{Q}_{\ell}})^{\omega}_{\mathbf{1}} \cong \D^{b}_{\mathrm{fg}}(\mathcal{H}) & \to \Coh(\Par_{G,\mathbf{1}}) \\
    M & \mapsto M \otimes_{\mathcal{H}} \mathcal{S}_{G,\mathbf{1}}.
\end{align*}
This discussion reduces Theorem \ref{thm:splitgroupBCHNcomparison} to the following result.

\begin{thm}\label{thm:splitgroupcpsicomputation}Maintain the setup and notation above. Then there is an isomorphism $c_{\psi}(i_{1!} \Pi_{G,\mathbf{1}}) \simeq \mathcal{S}_{G,\mathbf{1}}$, and the induced ring map
\[ \mathcal{H} = \mathrm{End}(\Pi_{G,\mathbf{1}}) \overset{c_{\psi} \circ i_{1!}}{\longrightarrow} \mathrm{End}(\mathcal{S}_{G,\mathbf{1}})\simeq \mathcal{H} \]
is a $\mathfrak{Z}$-algebra isomorphism.
\end{thm}
For this we need several preparatory lemmas.

\begin{lem}\label{lem:tracemagic}
Let $R$ be a noetherian normal domain with fraction field $K$, and
let $A$ be an $R$-order in a central simple $K$-algebra $D$. Then
any injective $R$-algebra endomorphism $f:A\to A$ is an isomorphism.
\end{lem}

\begin{proof}
It is clear that $f$ extends to a $K$-algebra automorphism of $A_{K}=D$,
so taking backwards iterates of $A\subset D$ under this automorphism
we get a chain of $R$-orders
\[
A\subseteq f^{-1}(A)\subseteq f^{-2}(A)\subseteq\cdots\subset D.
\]
Now standard trace magic gives an $R$-module embedding
\[
A^{\ast}=\mathrm{Hom}_{R}(A,R)\hookrightarrow D
\]
such that $B\subseteq A^{\ast}$ for any inclusion of $R$-orders
$A\subseteq B\subset D$ (see e.g. Corollary 3.5 in \cite{Chan}).
Therefore we get a containment
\[
A\subseteq f^{-1}(A)\subseteq f^{-2}(A)\subseteq\cdots\subseteq A^{\ast}.
\]
But $R$ is noetherian and $A^{\ast}$ is a finite $R$-module, so
this chain must stabilize.
\end{proof}
\begin{lem}\label{lem:IHcentralsimple}
Let $\mathcal{H}=\mathrm{End}(C_{c}(I\backslash G(F)))$ be the Iwahori-Hecke
algebra of a split p-adic group, with center $\mathfrak{Z}=\mathcal{O}(\hat{T})^{W}$.
Set $K=\mathrm{Frac}\mathfrak{Z}$. Then $\mathcal{H}_{K}$ is a central
simple $K$-algebra.
\end{lem}

After writing the argument below, we found a paper of Psaromiligkos \cite{Psaro} which proves strictly more: Lemma 2.37 there gives an explicit presentation
of $\mathcal{H}_{K}$ as a twisted group algebra $L\left\langle W\right\rangle $. Later Solleveld proved (at our request) an ultimately general result of this type, which applies to any Bernstein block of any reductive group \cite[Corollary 3.9]{Solleveld}. We decided to retain the proof below for its simplicity.

\begin{proof}
Set $L=\mathrm{Frac}\mathcal{O}(\hat{T})$, so $K\to L$ is a finite
Galois extension with Galois group $W$. It suffices to show that
$\mathcal{H}_{L}$ is a matrix algebra over $L$. Set $M=C_{c}(T(F)/T(\mathcal{O}_F))\simeq\mathcal{O}(\hat{T})$.
Then $C_{c}(I\backslash G(F))\simeq i_{B}^{G}M$, compatibly with the
obvious $\mathfrak{Z}$-module structure on $i_{B}^{G}M$ induced
by $\mathcal{O}(\hat{T})^{W}\to\mathcal{O}(\hat{T})\circlearrowright M$,
so
\[
\mathcal{H}\simeq\mathrm{End}_{G}(i_{B}^{G}M)=\mathrm{End}_{\mathfrak{Z}[G]}(i_{B}^{G}M).
\]
We then get isomorphisms 
\begin{align*}
\mathcal{H}_{L} & \simeq\mathrm{End}_{L[G]}((i_{B}^{G}M)\otimes_{\mathfrak{Z}}L)\\
 & \simeq\mathrm{End}_{L[G]}(i_{B}^{G}(M\otimes_{\mathfrak{Z}}L)).
\end{align*}
Let $\delta:T\to L^{\times}$ be the universal unramified character,
i.e. the character such that $M\otimes_{\mathcal{O}(\hat{T})}L\simeq\delta$
as smooth $L[T]$-modules. Then $M\otimes_{\mathfrak{Z}}L\simeq\oplus_{w\in W}\delta^{w}$.
But the smooth $L[G]$-modules $i_{B}^{G}(\delta^{w})$ are absolutely simple
and pairwise-isomorphic as $w$ varies, so 
\begin{align*}
i_{B}^{G}(M\otimes_{\mathfrak{Z}}L) & \simeq\oplus_{w\in W}i_{B}^{G}(\delta^{w})\\
 & \simeq i_{B}^{G}(\delta)^{\oplus|W|}
\end{align*}
and thus
\begin{align*}
\mathcal{H}_{L} & \simeq\mathrm{End}_{L[G]}\left(i_{B}^{G}(\delta)^{\oplus|W|}\right)\\
 & \simeq M_{|W|}(L)
\end{align*}
is a matrix algebra over $L$ as desired.
\end{proof}

Combining these observations, we get the following result.
\begin{prop}\label{prop:injisomagic}Let $\mathcal{H}=\mathrm{End}(C_{c}(I\backslash G(F)))$ be the Iwahori-Hecke
algebra of a split p-adic group, with center $\mathfrak{Z}=\mathcal{O}(\hat{T})^{W}$. Then any injective
$\mathfrak{Z}$-algebra endomorphism $\mathcal{H}\to\mathcal{H}$
is automatically an isomorphism.
\end{prop}
\begin{proof}
Let $K=\mathrm{Frac}\mathfrak{Z}$. By Lemma \ref{lem:IHcentralsimple}, $\mathcal{H}_K$ is a central simple $K$-algebra. Moreover, by the Bernstein presentation of $\mathcal{H}$, it is easy
to see that $\mathcal{H}$ is a $\mathfrak{Z}$-order in $\mathcal{H}_{K}$. Therefore Lemma \ref{lem:tracemagic} applies.
\end{proof}

\begin{proof}[Proof of Theorem \ref{thm:splitgroupcpsicomputation}]

Consider the sheaf $i_{1!}^{T}C_{c}(T(F)/T(\mathcal{O}_F))$ on $\mathrm{Bun}_{T}$. Under the known categorical equivalence for tori, this corresponds
to the structure sheaf on the unramified component $\Par_{T}^{\mathrm{nr}}$. Now $\Pi_{G,\mathbf{1}} \simeq i_{B^-}^{G}C_{c}(T(F)/T(\mathcal{O}_F))$ as noted above, so
\begin{align*}
 i_{1!}\Pi_{G,\mathbf{1}} & \simeq i_{1!} i_{B^-}^{G}C_{c} (T(F)/T(\mathcal{O}_F)) \\
  & \simeq \Eis_{B^- !} i_{1!}^{T}C_{c}(T(F)/T(\mathcal{O}_F)).
 \end{align*}
By the assumed compatibility of $c_{\psi}$ with Eisenstein series, we then compute 
\begin{align*}
c_{\psi}\left(i_{1!} \Pi_{G,\mathbf{1}} \right) & \simeq c_{\psi}\left(\mathrm{Eis}_{B^-}i_{1!}^{T}C_{c}(T(F)/T(\mathcal{O}_F))\right)\\
 & = \mathrm{Eis}_{B}^{\mathrm{spec}}(\mathcal{O}_{\mathrm{Par}_{T}^{\mathrm{nr}}})\\
 & = \mathcal{S}_{G,\mathbf{1}}.
\end{align*}
Taking endomorphism algebras on both sides, we see that $c_{\psi}$
induces a canonical map
\[
f:\mathrm{End}(\Pi_{G,\mathbf{1}})\overset{c_{\psi} \circ i_{1!}}{\longrightarrow}\mathrm{End}(\mathcal{S}_{G,\mathbf{1}}).
\]
We need to see that this
map is an isomorphism. By the results recalled above, the source and target of $f$ are both isomorphic to $\mathcal{H}$,
and $f$ is a map of $\mathfrak{Z}$-algebras, so by Proposition \ref{prop:injisomagic}
 we just need to check that $f$ is \emph{injective}. 

Now, it is easy to prove that for any maximal ideal $\mathfrak{m}\subset\mathfrak{Z}$ corresponding
to a generous parameter \cite[Section 2.1]{Beijing}, the lower horizontal arrow in the commutative
diagram
\[
\xymatrix{\mathrm{End}(\Pi_{G,\mathbf{1}})\ar[r]^{c_{\psi} \circ i_{1!}}\ar[d] & \mathrm{End}(\mathcal{S}_{G,\mathbf{1}})\ar[d]\\
\mathrm{End}(\Pi_{G,\mathbf{1}}\otimes_{\mathfrak{Z}}\mathfrak{Z}/\mathfrak{m})\ar[r]^{c_{\psi} \circ i_{1!}} & \mathrm{End}(\mathcal{S}_{G,\mathbf{1}}\otimes_{\mathfrak{Z}}\mathfrak{Z}/\mathfrak{m})
}
\]
is an isomorphism. We will justify this in the next paragraph. Admitting this result, we then conclude using the observation that the maps
\[
\mathrm{End}(\Pi_{G,\mathbf{1}})\to\prod_{\mathfrak{m}\,\mathrm{generous}}\mathrm{End}(\Pi_{G,\mathbf{1}}\otimes_{\mathfrak{Z}}\mathfrak{Z}/\mathfrak{m})
\]
and
\[
\mathrm{End}(\mathcal{S}_{G,\mathbf{1}})\to\prod_{\mathfrak{m}\,\mathrm{generous}}\mathrm{End}(\mathcal{S}_{G,\mathbf{1}}\otimes_{\mathfrak{Z}}\mathfrak{Z}/\mathfrak{m})
\]
are injective, which is an easy consequence of the Zariski-density of generous parameters in $\Spec \mathfrak{Z}$ plus the fact that 
$\mathcal{H}$ is $\mathfrak{Z}$-torsion free.

It remains to justify the claim about specialization at generous parameters. Fix such a parameter $\varphi$, with associated maximal ideal $\mathfrak{m}_{\varphi}$. Let $\pi_{\varphi}$ be the unique irreducible principal series representation of $G(F)$ with Fargues-Scholze parameter $\varphi$. Let $k_{\varphi} \in \Coh(\Par_{G,\mathbf{1}})$ be the pushforward of the structure sheaf on the closed substack $BS_{\varphi} \subset \Par_{G,\mathbf{1}}$. Note that both $\pi_{\varphi}$ and $k_{\varphi}$ have endomorphism ring $\Qellbar$. Now it is easy to see that there are isomorphisms \[\Pi_{G,\mathbf{1}}\otimes_{\mathfrak{Z}}\mathfrak{Z}/\mathfrak{m_{\varphi}} \simeq \pi_{\varphi}^{\oplus |W|} \]
and
\[\mathcal{S}_{G,\mathbf{1}} \otimes_{\mathfrak{Z}}\mathfrak{Z}/\mathfrak{m_{\varphi}} \simeq k_{\varphi}^{\oplus |W|},\]
using the assumption that $\varphi$ is generous in both cases.
We then conclude by the observation that $c_{\psi}(i_{1!} \pi_{\varphi}) = k_{\varphi}$, which again follows from the assumed compatibility with Eisenstein series.
\end{proof}

\subsection{General components}\label{generalGLn}
We now return to $\GL_n$. To go beyond the principal block, we will exploit the special fact that every Bernstein block of $\mathrm{GL}_n(F)$ is isomorphic to the principal block for some auxiliary group of the form $\prod_{i} \mathrm{GL}_{n_i}(E_i)$. To implement this idea precisely, we need to fix some convenient notation and terminology.

\begin{defn}A \emph{block datum} for $\mathrm{GL}_n$ is an ordered set of triples $\{ (n_i,m_i,\pi_i) \mid 1 \leq i \leq l \}$ where $n_i$ and $m_i$ are positive integers with $\sum_{1 \leq i \leq l} n_i m_i =n$ and $\pi_i$ is an irreducible supercuspidal representation of $\mathrm{GL}_{n_i}(F)$. We require that for every pair $1\leq i \neq j \leq l$ with $n_i = n_j$, $\pi_i$ is not an unramified twist of $\pi_j$.

Two block data $\{ (n_i,m_i,\pi_i) \mid 1 \leq i \leq l \}$ and $\{ (n_i',m_i',\pi_i') \mid 1 \leq i \leq l' \}$ are \emph{equivalent} if $l=l'$ and there exists $\sigma \in S_l$ such that for all $1\leq i \leq l$, $n_i' = n_{\sigma(i)}$, $m_i'=m_{\sigma(i)}$, and $\pi_i' \simeq \pi_{\sigma(i)} \otimes \eta_i$ for $\eta_i$ some unramified character of $\mathrm{GL}_{n_i'}(F)$.
\end{defn}

The following proposition is well-known.
\begin{prop}There are canonical bijections between the connected components of $\Par_{\GL_n}$, the Bernstein components of $\mathrm{GL}_n(F)$, and equivalence classes of block data for $\mathrm{GL}_n$.
\end{prop}

By a slight abuse, we will now identify the set of Bernstein components $\mathfrak{B}(\mathrm{GL}_n)$ with a fixed set of representatives of all equivalence classes of block data, writing $\mathfrak{s}$ for both a typical Bernstein component and for the corresponding element $\{ (n_i,m_i,\pi_i) \mid 1 \leq i \leq l \}$. Given any such $\mathfrak{s}$, we get standard block-diagonal Levis $L_{\mathfrak{s}} = \prod_i \mathrm{GL}_{n_i m_i}$ and $M_{\mathfrak{s}}= \prod_i \mathrm{GL}_{n_i}^{m_i}$ with  obvious inclusions $M_{\mathfrak{s}} \subset L_{\mathfrak{s}} \subset G$. We also get obvious Bernstein components for these two Levis, and corresponding components of their parameter stacks, which we also all label by $\mathfrak{s}$. Let $P_{\mathfrak{s}}$ be the standard parabolic of $G$ with Levi $M_{\mathfrak{s}}$.

Now, the external tensor product $ \sigma_{\mathfrak{s}}:=\boxtimes_{1 \leq i \leq l} \pi_{i}^{\boxtimes m_i}$ is an irreducible supercuspidal representation of $M_{\mathfrak{s}}(F)$. Choose an irreducible $M_{\mathfrak{s}}(F)^\circ$-representation $\tau'$ occurring in the restriction of $\sigma_{\mathfrak{s}}$ to $M_{\mathfrak{s}}(F)^{\circ}$, and set $\tau_{\mathfrak{s}} = \mathrm{ind}_{M_{\mathfrak{s}}(F)^\circ}^{M_{\mathfrak{s}}(F)} \tau'$. This is a compact projective generator for the supercuspidal block $\mathfrak{s} \in \mathfrak{B}(M_{\mathfrak{s}})$, and up to isomorphism it is independent of the choice of $\tau'$; in fact, there is an identification $\tau_{\mathfrak{s}}\simeq W_{\psi_{M_{\mathfrak{s}}},\mathfrak{s}}$. This identification shows in particular that every irreducible $M_{\mathfrak{s}}(F)$-representation in the block $\mathfrak{s}$ occurs with multiplicity one as a quotient of $\tau_{\mathfrak{s}}$.

\begin{lem}\label{lem:cpsisupercuspidalseedGLn}There is an isomorphism $c_{\psi_{M_{\mathfrak{s}}}}(i_{1!}^{M_{\mathfrak{s}}} \tau_{\mathfrak{s}}) \simeq \mathcal{O}_{\Par_{M_{\mathfrak{s}},\mathfrak{s}}}$.
\end{lem}
\begin{proof}
By Theorem \ref{thm:cpsisupercuspidalcoherent}, we already know $ \mathcal{F}:= c_{\psi_{M_{\mathfrak{s}}}}(i_{1!}^{M_{\mathfrak{s}}} \tau_{\mathfrak{s}}) $ is a coherent sheaf supported on $\Par_{M_{\mathfrak{s}},\mathfrak{s}}$.
The Fargues-Scholze parameter of $\sigma_{\mathfrak{s}}$ determines a canonical basepoint in $\Par_{M_{\mathfrak{s}},\mathfrak{s}}$, which then induces a canonical $Z(\hat{M_{\mathfrak{s}}})$-gerbe
\[ \Par_{M_{\mathfrak{s}},\mathfrak{s}} \to X/K  \]
where $X$ is the torus of unramified characters of $M_{\mathfrak{s}}(F)$ and $K$ is the finite subgroup of characters such that $\sigma_{\mathfrak{s}} \simeq \sigma_{\mathfrak{s}} \otimes \eta$. Note that $X/K$ is still a torus, so it has no nontrivial vector bundles by the Quillen-Suslin theorem, and the only character of $Z(\hat{M_{\mathfrak{s}}})$ occurring in the grading of $\mathcal{F}$ is the trivial character by Proposition \ref{prop:cpsigrading}. By Proposition \ref{prop:strsheafcriterionsupercuspidal}, we're reduced to checking that for every closed point $x \in X/K$ with corresponding closed immersion $i_x : B Z(\hat{M_{\mathfrak{s}}}) \to \Par_{M_{\mathfrak{s}},\mathfrak{s}}$, $i_{x}^{!}\mathcal{F}$ is of length one and concentrated in degree $\dim X$.

We already know the only possibly nontrivial graded pieces can occur for the trivial character, so it's enough to check that \[\RHom(\mathcal{O}_{B Z(\hat{M_{\mathfrak{s}}})}, i_{x}^{!}\mathcal{F}) \simeq \Qellbar[-\dim X].\]
Now we compute
\begin{align*}
    \RHom(\mathcal{O}_{B Z(\hat{M_{\mathfrak{s}}})}, i_{x}^{!}\mathcal{F}) & \simeq \RHom(i_{x \ast} \mathcal{O}_{B Z(\hat{M_{\mathfrak{s}}})}, \mathcal{F}) \\
    & \simeq \RHom(a_{\psi_{M_{\mathfrak{s}}}}(i_{x \ast} \mathcal{O}_{B Z(\hat{M_{\mathfrak{s}}})}), i_{1!}^{M_{\mathfrak{s}}} \tau_{\mathfrak{s}})
\end{align*}
where in the second line we used the definition of $\mathcal{F}$ and the $a_{\psi}$-$c_{\psi}$ adjunction. Now $x$ determines a $K$-coset of unramified characters $K \eta_x \subset X$, and the twist $\sigma_{\mathfrak{s}}\otimes \eta_x$ is well-defined independently of the choice of coset representative by the definition of $K$. It is then easy to see that \[a_{\psi_{M_{\mathfrak{s}}}}(i_{x \ast} \mathcal{O}_{B Z(\hat{M_{\mathfrak{s}}})}) \simeq i_{1!}^{M_{\mathfrak{s}}} \sigma_{\mathfrak{s}}\otimes \eta_x\] by the linearity of the spectral action over the spectral Bernstein center together with the fact that $W_{M_{\mathfrak{s}},\mathfrak{s}} \otimes_{\mathfrak{Z}_{M_{\mathfrak{s}},\mathfrak{s}}}\mathfrak{Z}_{M_{\mathfrak{s}},\mathfrak{s}}/\mathfrak{m}_x \simeq \sigma_{\mathfrak{s}} \otimes \eta_x$. This plus the full faithfulness of $i_{1!}^{M_{\mathfrak{s}}}$ reduces us to verifying that
\[ \RHom(\sigma_{\mathfrak{s}}\otimes \eta_x, \tau_{\mathfrak{s}}) \simeq \Qellbar[-\dim X],\]
which is an easy computation. Indeed, we can take Bernstein-Zelevinsky duals on both terms, which has the effect of replacing $\mathfrak{s}$ by the evident "dual" block datum $\mathfrak{s}^{\vee}$, and sends $\tau_{\mathfrak{s}}$ to $\tau_{\mathfrak{s}^{\vee}}$ and $\sigma_{\mathfrak{s}}\otimes \eta_x$ to $\sigma_{\mathfrak{s}^{\vee}}\otimes \eta_{x}^{-1}[-\dim X]$. This gives
\begin{align*}
    \RHom(\sigma_{\mathfrak{s}}\otimes \eta_x, \tau_{\mathfrak{s}}) & \simeq \RHom(\tau_{\mathfrak{s}^{\vee}},\sigma_{\mathfrak{s}^{\vee}}\otimes \eta_{x}^{-1})[-\dim X] \\
    & \simeq \Qellbar[-\dim X].
\end{align*}
To justify the final line, first note that concentration in degree $\dim X$ is automatic because both terms in the dualized $\RHom$ are in degree zero and the first term is a projective representation. Finally, the multiplicity one property of $\tau_{\mathfrak{s}}$ noted immediately before the lemma (which passes to $\tau_{\mathfrak{s}^{\vee}}$) implies that $\Hom(\tau_{\mathfrak{s}^{\vee}},\sigma_{\mathfrak{s}^{\vee}}\otimes \eta_{x}^{-1}) \simeq \Qellbar$. This gives the result.
\end{proof}

Now set $\Pi_{\mathfrak{s}}=i_{P_{\mathfrak{s}}^{-}}^{G} \tau_{\mathfrak{s}}$. This is a compact projective generator for the block $\mathfrak{s} \in \mathfrak{B}(\mathrm{GL}_n)$. Set $\mathcal{H}_{\mathfrak{s}}=\mathrm{End}_{\GL_n(F)}(\Pi_{\mathfrak{s}})$. It is well-known that $\mathcal{H}_{\mathfrak{s}}$ is a tensor product of type $A$ Iwahori-Hecke algebras, with center $\mathfrak{Z}_{\mathfrak{s}}$ canonically isomorphic to the ring of functions on $\Par_{\GL_n,\mathfrak{s}}$, and also canonically isomorphic to the $\mathfrak{s}$-component of the Bernstein center $\mathfrak{Z}_{\GL_n,\mathfrak{s}}$. Now again following \cite{BCHN}, we consider the coherent sheaf $\mathcal{S}_{\mathfrak{s}} = \Eis_{P_{\mathfrak{s}}}^{\spec}(\mathcal{O}_{\Par_{M_{\mathfrak{s}},\mathfrak{s}}})$. We then have the following variant on Theorem \ref{thm:BCHNsplitGspringer}.
\begin{thm}[Ben--Zvi-Chen-Helm-Nadler] For any block datum $\mathfrak{s} \in \mathfrak{B}(\GL_n)$, we have \[\Hom(\mathcal{S}_{\mathfrak{s}},\mathcal{S}_{\mathfrak{s}}[i])=0\] for all $i\neq 0$, and there is a $\mathfrak{Z}_{\mathfrak{s}}$-algebra isomorphism $\mathrm{End}(\mathcal{S}_{\mathfrak{s}}) \simeq \mathcal{H}_{\mathfrak{s}}$.
\end{thm}

Now just as in the previous section, the functor $\mathscr{S}_{\GL_n}$ restricted to the block $D(\mathrm{GL}_n(F),\Qellbar)_{\mathfrak{s}}^{\omega}$ is characterized by the requirements that a) it send the representation $\Pi_{\mathfrak{s}}$ to the coherent sheaf $\mathcal{S}_{\mathfrak{s}}$ and b) the induced map $\mathrm{End}(\Pi_{\mathfrak{s}}) \to \mathrm{End}(\mathcal{S}_{\mathfrak{s}})$ is a $\mathfrak{Z}_{\mathfrak{s}}$-algebra isomorphism. This reduces Theorem \ref{thm:GLnfullfaithndconsequences} to the following result.

\begin{thm}\label{thm:BCHNcpsiGLnessence}For any block datum $\mathfrak{s} \in \mathfrak{B}(\GL_n)$, there is an isomorphism $c_{\psi}(i_{1!} \Pi_{\mathfrak{s}}) \simeq \mathcal{S}_{\mathfrak{s}}$, and the induced map
\[ \mathcal{H}_{\mathfrak{s}} = \mathrm{End}(\Pi_{\mathfrak{s}}) \overset{c_{\psi} \circ i_{1!}}{\longrightarrow} \mathrm{End}(\mathcal{S}_{\mathfrak{s}})\simeq \mathcal{H}_{\mathfrak{s}} \]
is a $\mathfrak{Z}_{\mathfrak{s}}$-algebra isomorphism.
\end{thm}
\begin{proof}
We check the first claim. By the assumed compatibility of $c_{\psi}$ with Eisenstein series, we  compute that
\begin{align*}
    c_{\psi}(i_{1!} \Pi_{\mathfrak{s}}) & = c_{\psi}(i_{1!} i_{P_{\mathfrak{s}}^{-}}^{G} \tau_{\mathfrak{s}}) \\
    & = c_{\psi}(\Eis_{P_{\mathfrak{s}}^{-} !} i_{1!}^{M_{\mathfrak{s}}} \tau_{\mathfrak{s}}) \\
    & \simeq \Eis_{P_{\mathfrak{s}}}^{\spec}(c_{\psi_{M_{\mathfrak{s}}}}i_{1!}^{M_{\mathfrak{s}}} \tau_{\mathfrak{s}}) \\
    & \simeq \Eis_{P_{\mathfrak{s}}}^{\spec}(\mathcal{O}_{\Par_{M_{\mathfrak{s}},\mathfrak{s}}}) \\
    & = \mathcal{S}_{\mathfrak{s}}.
\end{align*}
Here the first and second lines follow from straightforward manipulations with Eisenstein series, the third line follows from the compatibility assumption, the fourth line follows from Lemma \ref{lem:cpsisupercuspidalseedGLn}, and the last line follows by the definition of $\mathcal{S}_{\mathfrak{s}}$. Now for the claim about the induced map on endomorphism rings, one argues exactly as in the proof of Theorem \ref{thm:splitgroupcpsicomputation}, first reducing to checking injectivity of the map by Proposition \ref{prop:injisomagic} and then verifying the injectivity by specialization at a dense set of generous parameters. We leave the details to the reader.
\end{proof}

Finally, we establish the following result.
\begin{thm}\label{thm:BCHNWhittakerkey}The BCHN functor
\[\mathscr{S}_{\GL_n}: \D(\GL_n(F),\Qellbar) \to \IndCoh(\Par_{\GL_n})\]
sends the Whittaker representation $W_{\psi}$ to $\mathcal{O}_{\Par_{\GL_n}}$.
\end{thm}
On the principal block, this is explicitly proved in \cite{BCHN} for all split groups. For general blocks the idea is already implicit in \cite{BCHN}, and we claim no originality here. 

Clearly we can argue one block at a time. We fix a block datum $\mathfrak{s}$ and all associated notation as above. In particular we write $M=M_{\mathfrak{s}}$ and $L=L_{\mathfrak{s}}$ for the standard Levis defined above. Let $Q$ be the standard parabolic in $G$ with Levi $L$.

The first key observation is that we have a commutative diagram
\[
\xymatrix{\otimes_{i}\D(\mathrm{GL}_{n_{i}m_{i}}(F),\overline{\mathbf{Q}_{\ell}})_{\mathfrak{s}_{i}}\ar[r]^{\sim}\ar[d]^{\otimes_{i}\mathscr{S}_{i}} & \D(L(F),\overline{\mathbf{Q}_{\ell}})_{\mathfrak{s}}\ar[r]_{i_{Q^{-}}^{\mathrm{GL}_{n}}}^{\sim}\ar[d]^{\mathscr{S}_{L}} & \D(\mathrm{GL}_{n}(F),\overline{\mathbf{Q}_{\ell}})_{\mathfrak{s}}\ar[d]^{\mathscr{S}_{\mathrm{GL}_{n}}}\\
\otimes_{i}\mathrm{IndCoh}(\mathrm{Par}_{\mathrm{GL}_{n_{i}m_{i}},\mathfrak{s}_{i}})\ar[r]^{\sim} & \mathrm{IndCoh}(\mathrm{Par}_{L,\mathfrak{s}})\ar[r]_{\mathrm{Eis}_{Q}^{\mathrm{spec}}}^{\sim} & \mathrm{IndCoh}(\mathrm{Par}_{\mathrm{GL}_{n},\mathfrak{s}})
}
\]
with arrows described as follows.  The unlabelled horizontal arrows are the obvious equivalences induced by Lurie tensor product via the decomposition $L=\prod_i \GL_{n_i m_i}$. The subscript $\mathfrak{s}_i$ refers to the evident projection of $\mathfrak{s}\in \mathfrak{B}(L)$ to the $i$th factor in this decomposition, and similarly for $M$. Writing $P_i$ for the standard parabolic in $\GL_{n_i m_i}$ with standard Levi $\GL_{n_i}^{m_i}$, the functor $\mathscr{S}_i$ sends $\Pi_i:=i_{P_{i}^{-}}^{\GL_{n_i m_i}} \sigma_{\mathfrak{s}_i}$ to $\Eis_{P_i}^{\spec}(\mathcal{O}_{\Par_{\GL_{n_i}^{m_i},\mathfrak{s}_i}})$.

Next we reduce to showing that $\mathscr{S}_L$ sends $W_{\psi_L, \mathfrak{s}}$ to $\mathcal{O}_{\Par_{L,\mathfrak{s}}}$. Going around the righthand square, this follows by combining from two observations:

$\bullet$ The functor $i_{Q^-}^{\GL_n}$ sends $W_{\psi_L, \mathfrak{s}}$ to $W_{\psi,\mathfrak{s}}$. Since $i_{Q^-}^{\GL_n}$ and $j_{\GL_n}^{Q^-}$ are mutually inverse equivalences on the $\mathfrak{s}$-blocks, this follows from taking the $\mathfrak{s}$-projection of Bushnell-Henniart's isomorphism $j_{\GL_n}^{Q^-} W_{\psi} \cong W_{\psi_L}$.

$\bullet$ The functor $\Eis_{Q}^{\spec}$ sends $\mathcal{O}_{\Par_{L,\mathfrak{s}}}$ to $\mathcal{O}_{\Par_{\GL_n,\mathfrak{s}}}$. This is immediate.

Now we consider the lefthand square. The lower horizontal arrow sends the external product of the structure sheaves to the structure sheaf, and the upper horizontal arrow sends the external product of $W_{\psi_{\GL_{n_i m_i}},\mathfrak{s}_i}$ to $W_{\psi_L,\mathfrak{s}}$. Therefore it suffices to see that each $\mathscr{S}_i$ sends $W_{\psi_{\GL_{n_i m_i}},\mathfrak{s}_i}$ to $\mathcal{O}_{\mathrm{Par}_{\mathrm{GL}_{n_{i}m_{i}},\mathfrak{s}_{i}}}$.  But now $\mathfrak{s}_i$ is a simple type for $\GL_{n_i m_i}$, and $\mathcal{H}_i=\mathrm{End}(\Pi_i)$ is a type A Iwahori-Hecke algebra; moreover, the main result of \cite{CS} shows that $W_{\psi_{\GL_{n_i m_i}},\mathfrak{s}_i}$ corresponds to the antispherical module of $\mathcal{H}_i$. Now we are reduced to the fact that $\mathscr{S}_i$ sends the antispherical $\mathcal{H}_i$-module to the structure sheaf, which is already proved in \cite{BCHN}.

\section{Classical groups}\label{sec:classical}

In this section we return to a general well-understood group $G$, equipped with a fixed Whittaker datum $(B,\psi)$ and maximal torus $T\subset B$ as usual. However, we will add two more axioms to the list of axioms (i)-(v) introduced in Section 6.1. Namely, we impose the following conditions on $G$ and all its standard Levi subgroups.

(vi) We have \emph{local Langlands in families} (borrowing a piece of terminology from \cite{DHKM2}): the natural composite map
\[\mathfrak{Z}_{G}^{\spec} \to \mathfrak{Z}_{G} \to \mathrm{End}(W_{\psi})\]
is an isomorphism of commutative rings.

(vii) For all supercuspidal $L$-parameters $\phi$, every $\rho \in \mathrm{Irr}(S_{\phi},1)$ is one-dimensional, and the corresponding supercuspidal $L$-packet $\Pi_{\phi}(G)$ over $b=1$ contains a unique $W_{\psi}$-generic member $\pi_{\phi}$.

If $G$ satisfies axioms (i)-(vii), we will say $G$ is \emph{very well-understood}. Again, this tautologically passes to standard Levi subgroups. Notably, all the groups listed in Remark \ref{rmk:WellUnderstoodRoster} are known to be very well-understood, and conjecturally all quasisplit classical groups satisfy these conditions. This explains the title of this chapter.

\begin{rmk} The first part of (vii) amounts to the condition that for all supercuspidal $L$-parameters $\phi$, the finite group $S_{\phi} / Z(\hat{G})^{W_F}$ is abelian. This assumption is admittedly unnatural, but we did not see how to avoid it in the proof of Theorem \ref{thm:CLLCoverirredparsVWE} below. It's certainly possible that some improvement to the argument could eliminate this condition. Note that this one-dimensionality definitely fails for some groups which are ``nearly'' classical groups, e.g. for $\mathrm{Spin}_9$. 
\end{rmk}

We emphasize that if $G$ is very well-understood, the natural (surjective) map $\mathfrak{B}(G)\to \mathfrak{B}^{\spec}(G)$ from Bernstein components to spectral Bernstein components admits a canonical section: by axiom (vi) together with Theorem \ref{thm:Whittakerbasics}.iii, there is a unique $\psi$-generic Bernstein component mapping to any spectral Bernstein component. We will freely identify $\psi$-generic Bernstein components with spectral Bernstein components via this observation. Note that under this identification, Theorem \ref{thm:Whittakerbasics}.iii implies that we have canonical identifications
\[\mathfrak{Z}_{G,\mathfrak{s}}^{\spec} = \mathfrak{Z}_{G,\mathfrak{s}} = \mathrm{End}(W_{\psi,\mathfrak{s}}) \]
for any spectral/$\psi$-generic Bernstein component $\mathfrak{s}$.

Recall that for well-understood $G$, the two natural orthogonal decompositions of $\D(\Bun_G)$ agree by Theorem \ref{thm:wellunderstoodSameOrthogonal}. In particular, $c_{\psi}$ restricts to a functor from $\D(\Bun_G)^{\cusp}=e^{\mathrm{irr}}\D(\Bun_G)$ towards $\IndCoh(\Par_{G}^{\mathrm{irr}}) = \QCoh(\Par_{G}^{\mathrm{irr}})$, where the last equality follows from the smoothness of $\Par_{G}^{\mathrm{irr}}$.

\begin{thm}\label{thm:CLLCoverirredparsVWE}For very well-understood $G$, the functor $c_{\psi}$ restricts to an equivalence of categories
\[c_{\psi}:\D(\Bun_G)^{\cusp} \overset{\sim}{\to} \QCoh(\Par_{G}^{\mathrm{irr}}).\]
For all the groups listed in Remark \ref{rmk:WellUnderstoodRoster}, this equivalence is t-exact for the standard t-structure on $\QCoh$.
\end{thm}

For very well-understood $G$, we also reduce the full faithfulness of $a_{\psi}$ and $\mathbf{R}_{\psi}$ to a precise conjectural analogue of Ben--Zvi-Chen-Helm-Nadler's work for \emph{$\psi$-generic Bernstein blocks} of $G$; the exact result we need is formulated in Conjecture \ref{conj:BCHNpsigenericdreamconjecture} below.

\begin{thm}\label{thm:genericBCHNimpliesFFgeneral}Assume that $G$ is very well-understood, and that Conjecture \ref{conj:BCHNpsigenericdreamconjecture} below holds true for all $\psi$-generic Bernstein blocks of $G$. Then $c_{\psi}i_{1!} |_{\D(G(F),\Qellbar)^{\omega}_{\mathfrak{s}}} \simeq \mathscr{S}_{G,\mathfrak{s}}$ for all $\psi$-generic Bernstein blocks $\mathfrak{s}$, where the right-hand functor is defined as in Conjecture \ref{conj:BCHNpsigenericdreamconjecture}.(2). Moreover, $\mathbf{L}_{\psi}i_{1!}W_{\psi}$ is a line bundle, and $a_{\psi}$ and $\mathbf{R}_{\psi}$ are fully faithful.
\end{thm}

We emphasize that here we require no information about inner forms or nongeneric blocks. We also emphasize that Conjecture \ref{conj:BCHNpsigenericdreamconjecture} is a \emph{consequence} of CLLC, so there is no doubt it must be true.

\subsection{Irreducible \texorpdfstring{$L$}{L}-parameters}
In this section we prove Theorem \ref{thm:CLLCoverirredparsVWE}. 
\begin{prop}\label{prop:cpsisupercuspidalclassical}For any $\psi$-generic supercuspidal block $\mathfrak{s}\in \mathfrak{B}(G)$, with associated component $\mathfrak{s} \in \mathfrak{B}^{\spec}(G)$, there is an isomorphism
\[c_{\psi}(W_{\psi,\mathfrak{s}}) \simeq \mathcal{O}_{\Par_{G,\mathfrak{s}}}.\]
\end{prop}
\begin{proof}This is an elaboration of the proof of Lemma \ref{lem:cpsisupercuspidalseedGLn}, but additional arguments are needed to deal with the nontriviality of $L$-packets. More precisely, we will use the criterion of Proposition \ref{prop:strsheafcriterionsupercuspidal}.

By Theorem \ref{thm:cpsisupercuspidalcoherent}, we know $\mathcal{F}:=c_{\psi}(W_{\psi,\mathfrak{s}})$ is a coherent sheaf supported on $\Par_{G,\mathfrak{s}}$. Choosing a basepoint parameter $\phi$ in this component of the stack, we get a finite \'etale surjection $ X^{\mathrm{nr}}(G) \times BS_{\phi} \to \Par_{G,\mathfrak{s}}$ where $X^{\mathrm{nr}}(G)$ is the torus of unramified characters of $G(F)$. By Proposition \ref{prop:strsheafcriterionsupercuspidal}, it suffices to check that for every closed point $\nu \in \Psi(G)$ with associated immersion $i_\nu: BS_{\phi} = B S_{\nu \cdot \phi} \to \Par_{G,\mathfrak{s}}$, $i_{\nu}^{!}\mathcal{F}$ is of length one, concentrated in degree $\dim X^{\mathrm{nr}}(G)$, and graded for the trivial representation of $S_{\phi}$. Note that the image of this immersion in the coarse moduli corresponds to the unramified twist $\nu \cdot \phi$ of $\phi$.

Letting $\rho$ denote any irreducible representation of $S_{\phi}$, arguing as in the proof of Lemma \ref{lem:cpsisupercuspidalseedGLn} shows that
\begin{align*}
    \RHom(\rho,i_{\nu}^{!}\mathcal{F}) & = \RHom(a_{\psi}(i_{\nu \ast} \rho),i_{1!}W_{\psi,\mathfrak{s}}) \\
    & = \RHom(i_{1!}W_{\psi^{-1},\mathfrak{s}^{\vee}},\Dbz(a_{\psi}(i_{\nu \ast} \rho)))
\end{align*}
where the second line follows from BZ duality. By compatibility of $a_{\psi}$ with the central grading, this $\RHom$ obviously vanishes unless $\rho$ is trivial on $Z(\hat{G})^{W_F}$, so we may assume $\rho$ satisfies this condition. Then by \cite[Theorem 4.5]{HJnote}, $a_{\psi}(i_{\nu \ast} \rho)$ is of the form $i_{1!}\nu \otimes \pi_{\rho}[n_{\rho}]$ for some irreducible supercuspidal representation $\pi_{\rho} \in \Pi_{\phi}(G)$ and some integer $n_{\rho}$, and the map $\rho \mapsto \nu \otimes \pi_{\rho}$ parametrizes the entire $L$-packet over $\nu \cdot \phi$. It's then easy to see that $n_{1}=0$ and $\nu \otimes \pi_{1}=\nu \otimes \pi_{\phi}$ is the unique $\psi$-generic member of the packet, and consequently that $\nu \otimes \pi_{\rho}$ is not $\psi$-generic for any nontrivial $\rho$. Using that $\Dbz$ on supercuspidals agrees with smooth duality up to an obvious shift and thus exchanges $\psi$-genericity with $\psi^{-1}$-genericity, this implies that the above $\RHom$ vanishes for all $\rho \neq 1$, while for $\rho=1$ it is of length one and concentrated in degree $\dim X^{\mathrm{nr}}(G)$ as desired.
\end{proof}
\begin{rmk}The one-dimensionality hypotheses in axiom (vii) was used in the proof above, in order to apply \cite[Theorem 4.5]{HJnote}. Again, we note that this condition will fail for some groups.
\end{rmk}

\begin{proof}[Proof of Theorem \ref{thm:CLLCoverirredparsVWE}]
By the previous proposition, $a_{\psi}$ induces a fully faithful functor from $\QCoh(\Par_{G}^{\mathrm{irr}})$ towards $\D(\Bun_G)^{\cusp}$. For essential surjectivity, let $A$ be any finite automorphic sheaf supported at the supercuspidal parameter $\phi$. By our axioms, we know $A$ will be concentrated at basic points, and at each basic point it will have finite length with irreducible subquotients lying in the packets $\Pi_{\phi}(G_b)$. Now it's enough to see that $c_{\psi}(A)\neq 0$.

Let $\pi_{\phi}$ be the unique $\psi$-generic member in the packet $\Pi_{\phi}(G)$. By the full faithfulness of $a_{\psi}$ and the argument in the previous proof, we know that $i_{1!}\pi_{\phi} = a_{\psi}(i_{\phi \ast} 1)$.  Now by our axioms together with the main theorem in \cite{HKW}, we get a nonzero map $i_{1!}\pi_{\phi}[n] \to T_V A$ for some $V$ and some integer $n$, which by adjunction corresponds to a nonzero map
\[i_{\phi \ast} 1[n] \to c_{\psi}(T_V A) = V \otimes c_{\psi}(A). \]
But the existence of such a map clearly forces $c_{\psi}(A)\neq 0$ as desired.
\end{proof}

\subsection{A coherent Springer conjecture for generic blocks}

Fix a very well-understood group $G$, with a fixed Whittaker datum $(B,\psi)$ and maximal torus $T\subset B$ as usual. Let $\mathfrak{s}=[M,\sigma]_G \in \mathfrak{B}(G)$ be a representative of a $\psi$-generic Bernstein component, with corresponding spectral component $\mathfrak{s}\in \mathfrak{B}^{\spec}(G)$. We also get an evident $\psi_M$-generic supercuspidal component $\mathfrak{s}_M=[M,\sigma]_{M}$ which again lives over a component $\mathfrak{s}_M\in \mathfrak{B}^{\spec}(M)$. Let $P$ be the standard parabolic in $G$ with Levi $M$.

First consider the representation $W_{\psi_M,\mathfrak{s}_M}$. This is a compact projective generator for the component $\mathfrak{s}_M$ which is indecomposable and satisfies $\mathrm{End}(W_{\psi_M,\mathfrak{s}_M})=\mathfrak{Z}_{M,\mathfrak{s}_M}$. Now set $\Pi_{\mathfrak{s}}=i_{P^-}^{G}W_{\psi_M,\mathfrak{s}_M}$. This is a compact projective generator for the component $\mathfrak{s}$. Set $\mathcal{H}_{\mathfrak{s}}=\mathrm{End}(\Pi_{\mathfrak{s}})$. Note that we have natural identifications of rings
\[\mathfrak{Z}_{G,\mathfrak{s}}^{\spec}=\mathcal{O}(\Par_{G,\mathfrak{s}}) = \mathrm{End}(W_{\psi,\mathfrak{s}}) = \mathfrak{Z}_{G,\mathfrak{s}} = Z(\mathcal{H}_{\mathfrak{s}}),  \]
where the first equality is tautological, the second equality comes by axiom (vi), the third equality is a consequence of Theorem \ref{thm:Whittakerbasics}.(iii), and the fourth equality is general nonsense (using that $\Pi_{\mathfrak{s}}$ is a projective generator). 

On the other hand, set $\mathcal{S}_{\mathfrak{s}}=\Eis_{P}^{\spec}(\mathcal{O}_{\Par_{M,\mathfrak{s}_M}})$. This is a coherent sheaf supported on $\Par_{G,\mathfrak{s}}$.
\begin{conjecture}\label{conj:BCHNpsigenericdreamconjecture}Maintain the above notation.
\begin{enumerate}
    \item For all $i\neq 0$, $\Hom(\mathcal{S}_{\mathfrak{s}},\mathcal{S}_{\mathfrak{s}}[i])=0$, and there is a $\mathfrak{Z}_{G,\mathfrak{s}}$-algebra isomorphism
    \[\mathcal{H}_{\mathfrak{s}} \simeq \mathrm{End}(\mathcal{S}_{\mathfrak{s}}). \]
    \item Choosing an isomorphism as in (1), the resulting functor
    \begin{align*}
    \mathscr{S}_{G,\mathfrak{s}}:\D(G(F),\overline{\mathbf{Q}_{\ell}})^{\omega}_{\mathfrak{s}} \cong \D^{b}_{\mathrm{fg}}(\mathcal{H}_{\mathfrak{s}}) & \to \Coh(\Par_{G,\mathfrak{s}}) \\
    M & \mapsto M \otimes_{\mathcal{H}_{\mathfrak{s}}} \mathcal{S}_{\mathfrak{s}}
\end{align*}
sends $W_{\psi,\mathfrak{s}}$ to a line bundle.
    
\end{enumerate}
\end{conjecture}

Of course we expect in (2) that the functor sends $W_{\psi,\mathfrak{s}}$ to $\mathcal{O}_{\Par_{G,\mathfrak{s}}}$, but the weaker statement above is sufficient for our purposes. We note that a proof of Conjecture \ref{conj:BCHNpsigenericdreamconjecture} for all quasisplit classical groups has been announced by Helm-Solleveld-Xu.

\begin{proof}[Proof of Theorem \ref{thm:genericBCHNimpliesFFgeneral}] By Proposition \ref{prop:fullyfaithfuleasy} and Corollary \ref{cor:apsiRpsiffmagic}, it's enough to prove that $c_{\psi}i_{1!}W_{\psi}$ is a line bundle. Choosing any $\mathfrak{s} \in \mathfrak{B}^{\spec}(G)$ which we identify with the associated $\psi$-generic Bernstein component of $G$ as in the discussion above, we're immediately reduced to proving that $c_{\psi}i_{1!}W_{\psi,\mathfrak{s}}$ is a line bundle on $\Par_{G,\mathfrak{s}}$ for every $\mathfrak{s}$. For this, the last condition of Conjecture \ref{conj:BCHNpsigenericdreamconjecture}.(2) shows that it suffices to check that the functor $c_{\psi}i_{1!}$ restricted to compact representations in the block $\mathfrak{s}$ identifies with the functor $\mathscr{S}_{G,\mathfrak{s}}$ defined in Conjecture \ref{conj:BCHNpsigenericdreamconjecture}.(2).

Now the argument is entirely analogous to the proof of Theorem \ref{thm:BCHNcpsiGLnessence}. We freely use the notation introduced in the formulation of Conjecture \ref{conj:BCHNpsigenericdreamconjecture}. By the assumed compatibility of $c_{\psi}$ with Eisenstein series, we compute
\begin{align*}
 c_{\psi}(i_{1!}\Pi_{\mathfrak{s}}) & = c_{\psi}(i_{1!}i_{P^-}^{G}W_{\psi_M, \mathfrak{s}_M}) \\
 & =c_{\psi}(\Eis_{P^-}i_{1!}^M W_{\psi_M, \mathfrak{s}_M}) \\
& = \Eis_{P}^{\mathrm{spec}}(c_{\psi_M} i_{1!}^M W_{\psi_M, \mathfrak{s}_M}) \\
& = \Eis_{P}^{\mathrm{spec}}(\mathcal{O}_{\Par_{M,\mathfrak{s}_M}}) \\
& = \mathcal{S}_{\mathfrak{s}}
\end{align*}
where the fourth line follows from Proposition \ref{prop:cpsisupercuspidalclassical} and the remaining lines are straightforward. 
Now $c_{\psi}i_{1!}$ induces a $\mathfrak{Z}_{G,\mathfrak{s}}$-algebra map $f: \mathrm{End}(\Pi_{\mathfrak{s}}) \to \mathrm{End}(\mathcal{S}_{\mathfrak{s}})$ where the source and target are both isomorphic to $\mathcal{H}_{\mathfrak{s}}$ (using Conjecture \ref{conj:BCHNpsigenericdreamconjecture}.(1) for the target). Writing $K=\mathrm{Frac}(\mathfrak{Z}_{G,\mathfrak{s}})$, a result of Solleveld \cite[Corollary 3.9]{Solleveld} shows that $\mathcal{H}_{\mathfrak{s}}\otimes K$ is a central simple $K$-algebra containing $\mathcal{H}_{\mathfrak{s}}$ as a $\mathfrak{Z}_{G,\mathfrak{s}}$-order, so Lemma \ref{lem:tracemagic} applies, reducing us to checking that $f$ is injective.  Arguing exactly as in the proof of Theorem \ref{thm:splitgroupcpsicomputation}, this injectivity can be verified after specialization at any Zariski-dense set of generous parameters $\varphi \in \Spec \mathfrak{Z}_{G,\mathfrak{s}}$. 

At any such parameter, let $k_{\varphi}$ be the skyscraper at the parameter $\varphi$ corresponding to the trivial representation of $BS_{\varphi}$, and let $\pi_{\varphi}$ be the unique irreducible representation in the block $\mathfrak{s}$ with Fargues-Scholze parameter $\varphi$. One then checks that $\Pi_{\mathfrak{s}} \otimes \mathfrak{Z}_{G,\mathfrak{s}}/\mathfrak{m}_{\varphi} = \pi_{\varphi}^{\oplus n}$ and $\mathcal{S}_{\mathfrak{s}} \otimes \mathfrak{Z}_{G,\mathfrak{s}}/\mathfrak{m}_{\varphi} = k_{\varphi}^{\oplus n}$ for some integer $n$ depending only on $\mathfrak{s}$, and that $c_{\psi}i_{1!}\pi_{\varphi} = k_{\varphi}$ using that $\pi_{\varphi}$ is an irreducible parabolic induction. 
\end{proof}

\subsection{Cuspidal control and endgame}
In applying Theorem \ref{thm:DreamInductionOnLevisTheorem} to very well-understood groups, condition (1) holds inductively, and Theorem \ref{thm:genericBCHNimpliesFFgeneral} reduces condition (2) to the tractable Conjecture \ref{conj:BCHNpsigenericdreamconjecture}. It remains to handle condition (3). In this section we reduce condition (3) to a concrete counting problem. Recall that condition (3) requires that $\mathbf{L}_{\psi} i_{1!}$ is conservative on compact supercuspidal representations and sends them into cuspidal coherent sheaves. The essential subtletly is that for supercuspidal representations whose Fargues-Scholze parameter is \emph{not} supercuspidal, the local geometry of the $L$-parameter stack at the relevant loci is much more subtle than in the case of supercuspidal parameters. 

We take the following approach. Assuming condition (2), the functor $\mathbf{R}_{\psi}$ is fully faithful, and Proposition \ref{prop:basicisright} crucially implies that it carries $\Coh^\qc(\Par_G)_{\cusp}^{0}$ into $i_{1!}\D(G(F),\Qellbar)_{\mathrm{sc}}^{\omega}$. In particular, the functor $i_{1}^\ast \mathbf{R}_{\psi}$ restricts to a fully faithful functor
\[r_{\psi}: \Coh^\qc(\Par_G)_{\cusp}^{0} \hookrightarrow \D(G(F),\Qellbar)_{\mathrm{sc}}^{\omega}.  \]
To verify condition (3), it is enough to see this fully faithful functor is essentially surjective, since then it is an \emph{equivalence} and its left adjoint $\mathbf{L}_{\psi} i_{1!}$ visibly induces an inverse equivalence, which (more than) verifies condition (3). We will reduce the essential surjectivity of $r_{\psi}$ to a simple numerical inequality, using the full power of Bertoloni Meli-Koshikawa's theory of cuspidal coherent sheaves.

To carry this out, we introduce some notation. For a fixed integer $d \geq 0$, set \[\mathcal{A}_d=\RHom_{\Qellbar[X_1,\dots,X_d]}(\Qellbar,\Qellbar)\] where $\Qellbar[X_1,\dots,X_d] \twoheadrightarrow \Qellbar$ is evaluation at zero. This is an associative algebra object in $\Perf(\Qellbar)$. Let $\mathcal{C}_d$ be the category of $\mathcal{A}_d$-modules in $\Perf(\Qellbar)$. Note that $\mathcal{C}_d$ is equivalent to the category $\Coh_{ \{0 \} }(\mathbf{A}^d)$. We also note that if $X$ is any object in a $\Qellbar$-linear stable $\infty$-category $\mathcal{D}$ such that $\RHom(X,X) \simeq \mathcal{A}_d$, then the thick subcategory of $\mathcal{D}$ generated by $X$ under shifts and cones is equivalent to $\mathcal{C}_d$.

\begin{prop}Let $\pi$ be an irreducible supercuspidal representation of $G(F)$, and let $\mathcal{C}_{\pi} \subset \D(G(F),\Qellbar)^{\omega}_{\mathrm{sc}}$ be the full subcategory generated by $\pi$ under shifts and cones. Then $\mathcal{C}_{\pi} \simeq \mathcal{C}_{d}$, where $d$ is the rank of the split center of $G$.

Moreover, if $\phi$ is a discrete $L$-parameter, and $\D(G(F),\Qellbar)^{\omega}_{\mathrm{sc},\phi} \subset \D(G(F),\Qellbar)^{\omega}_{\mathrm{sc}}$ denotes the subcategory of compact ULA objects such that all irreducible subquotients have $L$-parameter $\phi$, then $\D(G(F),\Qellbar)^{\omega}_{\mathrm{sc},\phi} \simeq \mathcal{C}_{d}^{\oplus m_{\phi}^{\mathrm{aut}}}$. Here $m_{\phi}^{\mathrm{aut}}$ denotes the number of supercuspidal representations in the $L$-packet $\Pi_{\phi}(G)$.
\end{prop}
\begin{proof}
The last part follows easily from the first part. Indeed, by definition $\D(G(F),\Qellbar)^{\omega}_{\mathrm{sc},\phi}$ is generated by the mutually orthogonal subcategories $\mathcal{C}_{\pi}$ as $\pi$ varies over the supercuspidal elements of the packet $\Pi_{\phi}(G)$.

The first part follows immediately from the observation that $\RHom(\pi,\pi) \simeq \mathcal{A}_d$.
\end{proof}

On the other side, for a discrete parameter $\phi$, let $\Coh(\Par_G)^{0}_{\cusp,\phi}$ be the category of $0$-graded cuspidal coherent sheaves which are supported set-theoretically on the closure $V_{\phi}$ of $\phi$. We will need the following deep result of Bertoloni Meli-Koshikawa.

\begin{thm}[Bertoloni Meli-Koshikawa] There is an integer $m_{\phi}^{\spec}\geq 0$ such that the category $\Coh(\Par_G)^{0}_{\cusp,\phi}$ is generated by a finite collection of mutually orthogonal objects $\mathcal{F}_i$, $1 \leq i \leq m_{\phi}^{\spec}$, each of which is the pushforward of a vector bundle on $V_{\phi}$, and with $\RHom(\mathcal{F}_i,\mathcal{F}_i) \simeq \mathcal{A}_{d}$ where again $d$ denotes the rank of the split center of $G$. In other words, there is an equivalence $\Coh(\Par_G)^{0}_{\cusp,\phi} \simeq \mathcal{C}_{d}^{\oplus m_{\phi}^{\spec}}$.
\end{thm}

Now we return to our functor $r_{\psi}$. 

\begin{prop}Maintain the notation and assumptions as above. For any discrete parameter $\phi$, $r_{\psi}$ restricts to a fully faithful embedding
\[\Coh(\Par_G)^{0}_{\cusp,\phi} \hookrightarrow  \D(G(F),\Qellbar)^{\omega}_{\mathrm{sc},\phi}.\]
\end{prop}
\begin{proof}
Anything in the source will automatically map to a finite length object in $\D(G(F),\Qellbar)^{\omega}_{\mathrm{sc}}$, all of whose irreducible subquotients have Fargues-Scholze parameter $\phi^{\mathrm{ss}}$. If $\phi' \neq \phi$ is another parameter with $\phi'^{\mathrm{ss}} = \phi^{\mathrm{ss}}$, then $\phi'$ cannot be a discrete parameter, so $\Pi_{\phi'}(G)$ does not contain any supercuspidal representations by our axioms. This implies that any irreducible supercuspidal $G(F)$-representation with Fargues-Scholze parameter $\phi^{\mathrm{ss}}$ automatically lies in the packet $\Pi_{\phi}(G)$. Unwinding definitions, the result follows.
\end{proof}

Maintaining the notation of this proposition, and recalling our prior descriptions of the source and target, we see that for any discrete parameter $\phi$, $r_{\psi}$ induces a fully faithful functor
\[\mathcal{C}_{d}^{\oplus m_{\phi}^{\spec}} \simeq \Coh(\Par_G)^{0}_{\cusp,\phi} \hookrightarrow  \D(G(F),\Qellbar)^{\omega}_{\mathrm{sc},\phi} \simeq \mathcal{C}_{d}^{\oplus m_{\phi}^{\mathrm{aut}}}.\]
Now we have the following crucial lemma.

\begin{lem}\label{lem:fullyfaithfultrick}Let $F:\mathcal{C}_{d} ^ {\oplus m} \hookrightarrow \mathcal{C}_{d} ^ {\oplus n}$ be an exact $\Qellbar$-linear fully faithful functor. Then $ m \leq n$, and if $m=n$ then $F$ is an equivalence of categories.
\end{lem}
\begin{proof}
For brevity write $k=\Qellbar$. We will freely use the equivalence $\mathcal{C}_d = \Coh_{ \{ 0 \} }(\mathbf{A}^{d})$. Fix $F$ as in the Lemma, and consider the object $G=(k[0],0,\dots,0)$ in $\mathcal{C}_{d}^{\oplus m}$. Since $G$ has endomorphism ring $k$ and $F$ is fully faithful, $F(G)$ also has endomorphism ring $k$, which in particular is a local ring. Therefore $F(G)$ projects to zero in all but one direct factor of $\mathcal{C}_{d}^{\oplus n}$. Applying a permutation autoequivalence of the target category, we can assume $F(G)\simeq (N,0,\dots,0)$ for some $N \in \mathcal{C}_{d}$. Let $a \leq b$ be the least resp. greatest degrees in which $H^j(N)$ is nonzero. Note that these cohomology groups are Artinian $\mathcal{O}_{\mathbf{A}^d,0}$-modules, so we can always choose a nonzero map $H^b(N) \to H^a(N)$. This induces a nonzero map
\[ N \to H^b(N)[-b] \to H^a (N)[-b] \to N[a-b], \]
so we get a nonzero element in $\Hom(N,N[a-b])$. But $F$ is fully faithful, so $\RHom(N,N)=\RHom(G,G)=\RHom_{\mathcal{C}_d}(k,k)$ vanishes in negative degrees. This shows that $a-b \geq 0$, so $a=b$ and $N$ is concentrated in degree $a$. (We learned this trick from the proof of \cite[Tag 0A7E]{stacks-project}.) We then also deduce that $\End(H^a(N))=\End(N)=\End(G)=k$, so $H^a(N)$ is a simple $\mathcal{O}_{\mathbf{A}^d,0}$-module, and thus $N=k[-a]$. This shows that $F(G)=(k[-a],0,\dots,0)$ is a generator of the first direct factor of $\mathcal{C}_{d} ^ {\oplus n}$. Passing to the Verdier quotient of the source resp. target by the direct factor subcategory generated by $G$ resp. $F(G)$ and continuing by induction, we see that $m \leq n$. If $m=n$ then after $n-1$ steps we will be left with a fully faithful functor $\mathcal{C}_d \to \mathcal{C}_d$, and the argument above shows that any such functor is automatically an equivalence.
\end{proof} 

From this Lemma, we immediately deduce the following corollary.

\begin{cor}\label{cor:cuspidalcounting}Let $G$ be a very well-understood group. Suppose that the Langlands functor carries $i_{1!}W_{\psi}$ to a line bundle, so $\mathbf{R}_{\psi}$ is fully faithful. Suppose also that for every discrete parameter $\phi$, we have an inequality $m_{\phi}^{\spec} \geq m_{\phi}^{\mathrm{aut}}$. 
Then $r_{\psi}$ induces an equivalence from $\Coh(\Par_G)_{\cusp,\fin}^{0}$ towards $\D(G(F),\Qellbar)_{\mathrm{sc},\fin}$, and the functor $\mathbf{L}_{\psi}i_{1!}$ induces an inverse equivalence.
\end{cor}
\begin{proof}
We have gradings
\[\Coh(\Par_G)_{\cusp,\fin}^{0} = \bigoplus_{\phi} \Coh(\Par_G)_{\cusp,\phi}^{0}\]
and
\[\D(G(F),\Qellbar)_{\mathrm{sc},\fin} = \bigoplus_{\phi} \D(G(F),\Qellbar)_{\mathrm{sc},\phi}^{\omega} \]
indexed by discrete parameters $\phi$, which are compatible under the fully faithful functor $r_{\psi}$. The result now follows from Lemma \ref{lem:fullyfaithfultrick}, using the numerical hypothesis together with our prior descriptions of the individual $\phi$-graded summands.
\end{proof}

\begin{rmk}
In the notation of the previous proof, there are discrete parameters for which both $\phi$-graded summands are empty. Note also that any supercuspidal parameter is discrete, so there is some redundancy at supercuspidal parameters with the results of Theorem \ref{thm:CLLCoverirredparsVWE}. The whole point of this section is to handle the ``mixed'' supercuspidals, i.e. the supercuspidals whose Fargues-Scholze parameter is not supercuspidal. However, the present arguments would not gain clarity if we artificially discarded the supercuspidal parameters.
\end{rmk}

Note that for all the groups listed in Remark \ref{rmk:WellUnderstoodRoster}, given any discrete parameter $\phi$, the integer $m_{\phi}^{\mathrm{aut}}$ is completely explicit and can be read off elementarily from $\phi$ by the work of Moeglin \cite{moeglin2014}. With this observation in hand, the forthcoming work of Bertoloni Meli-Koshikawa should establish the equality $m_{\phi}^{\spec} = m_{\phi}^{\mathrm{aut}}$ for all discrete parameters of all quasisplit classical groups.

Putting everything together, we conclude with the following theorem.

\begin{thm}\label{thm:ultimateVWO}Let $G$ be a very well-understood group with a fixed Whittaker datum $\psi$. Assume that $\mathbf{L}_{\psi_M}$ is an equivalence of categories for all proper standard Levis $M\subsetneq G$. Assume also that Conjecture \ref{conj:BCHNpsigenericdreamconjecture} is true for all $\psi$-generic Bernstein components of $G(F)$, and that $m_{\phi}^{\spec} = m_{\phi}^{\mathrm{aut}}$ for all discrete $L$-parameters into $\phantom{}^L G$.

Then $\mathbf{L}_{\psi}$ is an equivalence of categories.
\end{thm}

Note that we can trivially rearrange the hypotheses in this theorem: instead of assuming full CLLC for proper Levis, it is equivalent to assume that for $G$ \emph{and} all its proper Levis, Conjecture \ref{conj:BCHNpsigenericdreamconjecture} holds for all $\psi$-generic Bernstein blocks and $m_{\phi}^{\spec} = m_{\phi}^{\mathrm{aut}}$ for all discrete $L$-parameters.

\begin{proof}
This follows from a minor variant of the proof of Theorem \ref{thm:DreamInductionOnLevisTheorem}. More precisely, Conjecture \ref{conj:BCHNpsigenericdreamconjecture} implies the full faithfulness of $\mathbf{R}_{\psi}$ via Theorem \ref{thm:genericBCHNimpliesFFgeneral}. Then we may apply Corollary \ref{cor:cuspidalcounting}, which shows in particular that $\mathbf{L}_{\psi}i_{1!}$ is conservative on finite-length supercuspidals and sends them into cuspidal coherent sheaves. Rerunning the argument for conservativity in the proof of Theorem \ref{thm:DreamInductionOnLevisTheorem} with the word ``finite'' everywhere, we deduce from this together with our inductive hypothesis on Levis that $\mathbf{L}_{\psi}$ is conservative on $\D(\Bun_G)_{\fin}$ and (obviously) carries it towards $\Coh(\Par_G)_{\fin}$. On the other hand, we also know that the right adjoint $\mathbf{R}_{\psi}$ restricts to a fully faithful functor $\Coh(\Par_G)_{\fin} \to \D(\Bun_G)_{\fin}$,  so by Lemma \ref{lem:trivialadjointcondition} we see that $\mathbf{L}_{\psi}$ restricts to an equivalence $\D(\Bun_G)_{\fin} \overset{\sim}{\to} \Coh(\Par_G)_{\fin}$. The full equivalence now follows from Theorem \ref{thm:restrvsfull}.
\end{proof}

Combining this Theorem with the forthcoming work of Helm-Solleveld-Xu on Conjecture \ref{conj:BCHNpsigenericdreamconjecture}, we deduce the full categorical local Langlands conjecture for all special orthogonal and unitary groups (modulo our running assumption of Eisenstein compatibility).

\appendix
\section{Abstract nonsense} \label{sec:abstract-nonsense}

In this section we prove several abstract results about higher categories that are needed in the paper, including properties of the category of correspondences and 6-functor formalisms. Among others we discuss the extension of a 6-functor formalism to a 3-functor formalism of \enquote{bounded sheaves} in \cref{sec:Corr-bounded} and study convolution monoidal structures in \cref{sec:convolution-stacks}. The reader is invited to skip the appendix on a first read and only refer to it when required by the main text.

\subsection{Basic lemmas}

This subsection contains several unrelated results in category theory, to be cited throughout the main text and in other parts of the appendix.

We start with two functoriality results about certain constructions in symmetric monoidal categories and 6-functor formalisms. In the following result, $\CMon := \CAlg(\Cat)$ denotes the category of symmetric monoidal categories and symmetric monoidal functors.

\begin{lem} \label{rslt:functoriality-of-dualizable-dual}
Let $F\colon \cat I \to \CMon$ be a functor and assume that for each $i \in \cat I$, all objects in $\cat C_i := F(i)$ are dualizable. Then the duality functor $(-)^\vee\colon \cat C_i^\op \isoto \cat C_i$ induces an isomorphism
\begin{align*}
    (-)^\vee\colon F(-)^\op \isoto F(-)
\end{align*}
of functors $\cat I \to \CMon$.
\end{lem}
\begin{proof}
The claims are roughly that for $P, Q \in \cat C_i$ we have $(P \tensor Q)^\vee = P^\vee \tensor Q^\vee$ and if $f\colon \cat C_i \to \cat C_j$ is a transition functor in the diagram $F$ then $f(P)^\vee = f(P^\vee)$. Both of these claims are easy to check on objects; the non-trivial claim is that these identities are compatible with all higher homotopies in the diagram $F$.

By definition a symmetric monoidal category is a functor $f\colon \Fin_* \to \Cat$ such that $f(\langle n \rangle) = f(\langle 1 \rangle)^n$ via the natural maps; then $\cat C := f(\langle 1 \rangle)$ is the underlying category and the functoriality of $f$ induces the tensor structure. In particular, the functor $F$ induces a functor $F'\colon \cat I \times \Fin_* \to \Cat$ with $F'(i, \langle n \rangle) = \cat C_i^n$. Consider the functor $q\colon \Fin_* \to \Fin_*$ that sends $\langle n \rangle \mapsto \langle 2n \rangle$ and acts on morphisms by repeating the same morphism twice. There is a natural transformation $\mu\colon q \to \id$ given by the maps $\langle 2n \rangle \to \langle n \rangle$ mapping $k \in \langle 2n \rangle$ to $k$ if $k \le n$ and to $k - n$ if $k > n$. By 2-functoriality of $\cat I \times -$ and $F' \comp -$, $\mu$ induces a natural transformation of functors $\alpha\colon F' q \to F'$. Explicitly, $\alpha$ encodes the natural transformation $\tensor\colon F(-) \times F(-) \to F(-)$. We now pass to the associated cocartesian fibrations, translating $\alpha$ into a functor $A\colon \cat E_{F'} \times_{\cat I \times \Fin_*} \cat E_{F'} \to \cat E_{F'}$ of cocartesian fibrations over $\cat I \times \Fin_*$. The category $\Fin_*$ has the initial object $\langle 0 \rangle$ and the fiber of $\cat E_{F'}$ over $\cat I \times \{ \langle 0 \rangle \}$ is $\cat I \times * = \cat I$, i.e. has a unique section. By cocartesian lifting this section extends to a section $s\colon \cat I \times \Fin_* \to \cat E_{F'}$ sending $(i, \langle n \rangle) \mapsto (1, \dots, 1) \in \cat C_i^n$. Plugging $s$ into the first argument of the relative hom-functor $\Hom_{\cat I \times \Fin_*}(-,-)$ produces a functor $\cat E_{F'} \to \Ani$ sending $(P_k)_k \in \cat C_i^n$ to $\prod_k \Hom(1, P_k)$. We precompose this functor with $A$ to obtain the functor
\begin{align*}
    \cat E_{F'} \times_{\cat I \times \Fin_*} \cat E_{F'} \to \Ani, \qquad ((P_k)_k, (Q_k)_k) \mapsto \prod_k \Hom(1, P_k \tensor Q_k).
\end{align*}
As in the proof of \cite[Lemma~A.2.32]{heyer-mann-6ff} (using \cite[Definition~3.77]{Heine.2023}) we obtain an induced functor
\begin{align*}
    B\colon \cat E_{F'} \to \Fun^{\cat I \times \Fin_*}(\cat E_{F'}, \Ani), \qquad (P_k)_k \mapsto \left[(Q_k)_k \mapsto \prod_k \Hom(1, P_k \tensor Q_k) \right].
\end{align*}
over $\cat I \times \Fin_*$, where the right-hand side the is the cartesian fibration classifying the functor $\Fun(F'(-), \Ani)$. It is easy to see from the definition of dualizable objects that in the fiber over $(i, [n])$ the functor $B$ identifies $\cat C_i^n$ with the full subcategory of $\Fun(\cat C_i^n, \Ani)$ spanned by the corepresentable functors; in fact $B((P_k)_k)$ is corepresented by $(P_k^\vee)_k$. By the relative Yoneda embedding (see \cite[Remark~A.2.24]{heyer-mann-6ff} and the proof of \cite[Lemma~A.2.32]{heyer-mann-6ff}) we obtain the desired equivalence of the claim.
\end{proof}

\begin{lem} \label{rslt:functoriality-of-suave-dual}
Let $\D$ be a 6-functor formalism on some category $\cat C$ and denote by $\D_!, \D_*\colon \cat C \to \Cat$ and $\D^*, \D^!\colon \cat C^\op \to \Cat$ the functors sending $X \in \cat C$ to $\D(X)$ and a morphism $f\colon Y \to X$ to $f_!$, $f_*$, $f^*$ and $f^!$, respectively. For each $X \in \cat C$ denote $\omega_X := f^! 1$, where $f\colon X \to *$ is the projection. Then the assignment $M \mapsto \intHom(M, \omega_X)$ defines natural transformations
\begin{align*}
    \D_!(X) \to \D_*(X)^\op, \qquad \D^*(X) \to \D^!(X)^\op
\end{align*}
of functors $\cat C \to \Cat$ and $\cat C^\op \to \Cat$, respectively.
\end{lem}
\begin{proof}
The claimed naturality of the functor $\intHom(-, \omega_X)\colon \D_!(X) \to \D_*(X)^\op$ boils down to the easy identity
\begin{align*}
    \intHom(g_! M, \omega_X) = g_* \intHom(M, g^! \omega_X) = g_* \intHom(M, \omega_Y)
\end{align*}
for any map $g\colon Y \to X$, and a similar computation exhibits the naturality of $\D^*(X) \to \D^!(X)^\op$. The challenge is to show that these identities are natural in $g$.

We thank Claudius Heyer for suggesting the following argument. Let $\cat K = \cat K_{\D}$ be the category of kernels of $\D$ (see \cite[Definition~4.1.3(a)]{heyer-mann-6ff}). Using the Yoneda embedding for 2-categories (see e.g. \cite[Theorem~C.2.10]{heyer-mann-6ff}) we can construct a natural transformation
\begin{align*}
    \Hom_{\cat K}(-, *) \to \Fun(\Hom_{\cat K}(*, -), \Hom_{\cat K}(*, *))
\end{align*}
of functors $\cat K^\op \to \Cat$; indeed, the Yoneda embedding says that such a morphism is given by an element in $\Fun(\Hom_{\cat K}(*, *), \Hom_{\cat K}(*, *))$ and we simply take the identity functor. From now on we view this as a natural transformation of 1-functors $\cat K^\op \to \Cat$ in order to avoid 2-category theory. Explicitly, for $X \in \cat C$ with projection $f\colon X \to *$ the above natural transformation translates to the functor
\begin{align*}
    \D(X) \to \Fun(\D(X), \D(*)), \qquad M \mapsto [N \mapsto f_!(M \tensor N)].
\end{align*}
In particular we see that the image of a given $M \in \D(X)$ is a left adjoint functor $\D(X) \to \D(*)$, so we get an induced natural transformation $\Hom_{\cat K}(-, *) \to \Fun^L(\Hom_{\cat K}(*, -), \Hom_{\cat K}(*, *))$, where $\Fun^L \subseteq \Fun$ denotes the full subcategory of left adjoint functors. Using the equivalence $(-)^\op\colon \Cat \isoto \Cat^\co$ of 2-categories (which follows by tracking the $\Cat$-linear structure on $\Cat$ through $(-)^\op$) and the equivalence $\Cat^L = \Cat^{R,\co,\op}$ from \cite[Example~D.3.18(a)]{heyer-mann-6ff} we obtain natural equivalences
\begin{align*}
    \Fun^L(\Hom_{\cat K}(*, -), \Hom_{\cat K}(*, *)) &= \Fun^R(\Hom_{\cat K}(*, -)^\op, \Hom_{\cat K}(*, *)^\op)^\op\\
    &= \Fun^L(\Hom_{\cat K}(*, *)^\op, \Hom_{\cat K}(*, -)^\op).
\end{align*}
Precomposition with the natural functor $* \to \Hom_{\cat K}(*, *)^\op = \D(*)^\op$ given by the tensor unit yields a natural transformation $\Fun^L(\Hom_{\cat K}(*, *)^\op, \Hom_{\cat K}(*, -)^\op) \to \Hom_{\cat K}(*, -)^\op$. By composing this with the natural transformation from above we obtain the natural transformation
\begin{align*}
    \Hom_{\cat K}(-, *) \to \Hom_{\cat K}(*, -)^\op
\end{align*}
of functors $\cat K^\op \to \Cat$. Tracing through the construction reveals that on a fixed object $X \in \cat C$, the above natural transformation is the functor $\intHom(-, \omega_X)\colon \D(X) \to \D(X)^\op$. Thus, if we precompose the natural transformation with the two natural functors $\cat C \to \cat K^\op$ and $\cat C^\op \to \cat K^\op$ (see \cite[Proposition~1.3.2]{heyer-mann-6ff}) then we arrive at the natural transformations in the claim.
\end{proof}

For the next result we recall the twisted arrow category $\Tw(\cat I)$ of a category $\cat I$. It can be defined as the cocartesian unstraightening of the functor $\cat I^\op \to \cat I \to \Ani$, $(X, Y) \mapsto \Hom(X, Y)$. Explicitly, the objects of $\Tw(\cat I)$ are morphisms $f\colon X \to Y$ in $\cat I$ and a map $f \to f'$ in $\Tw(\cat I)$ is a commuting square
\begin{equation*}\begin{tikzcd}
    X \arrow[d,"f",swap] & X' \arrow[d,"f'"] \arrow[l]\\
    Y \arrow[r] & Y'
\end{tikzcd}\end{equation*}
in $\cat I$. The twisted arrow category can be used to encode the naturality of the unit and counit of an adjunction, as follows.

\begin{lem}
Let $\cat I$, $\cat C$ and $\cat D$ be categories and let $\ell\colon \cat I \to \Fun(\cat C, \cat D)$ be a functor such that for all $S \in \cat I$, $\ell_S := \ell(S)$ admits a right adjoint $r_S\colon \cat D \to \cat C$. Then:
\begin{lemenum}
    \item \label{rslt:naturality-of-adjunction-unit} There are a functor $\Tw(\cat I) \to \Fun(\cat C, \cat C)$, $[S \to T] \mapsto r_S \ell_T$ and a natural transformation $\id_{\cat C} \to r_S \ell_T$ of such functors such that the induced maps $\id_{\cat C} \to r_S \ell_S$ are given by the units of the adjunctions.

    \item \label{rslt:naturality-of-adjunction-counit} There are a functor $\Tw(\cat I)^\op \to \Fun(\cat D, \cat D)$, $[S \to T] \mapsto \ell_S r_T$ and a natural transformation $\ell_S r_T \to \id_{\cat D}$ of such functors such that the induced maps $\ell_S r_S \to \id_{\cat D}$ are given by the counits of the adjunctions.
\end{lemenum}
\end{lem}
\begin{proof}
Part (ii) follows formally from part (i) by applying $(-)^\op$ and using the natural isomorphism $\Fun(\cat D, \cat D)^\op = \Fun(\cat D^\op, \cat D^\op)$. To prove (i), we consider the two functors
\begin{align*}
    &F_1, F_2 \colon \Tw(\cat I) \times \cat C^\op \times \cat C \to \Ani,\\
    &F_1([S \to T], X, Y) = \Hom(X, Y),\\
    &F_2([S \to T], X, Y) = \Hom(\ell_S(X), \ell_T(Y)).
\end{align*}
There is a natural transformation $\alpha\colon F_1 \to F_2$ that acts as the composition
\begin{align*}
    \Hom(X, Y) \to \Hom(\ell_S(X), \ell_S(Y)) \to \Hom(\ell_S(X), \ell_T(Y)).
\end{align*}
More formally, this natural transformation can be constructed on the unstraightenings of $F_1$ and $F_2$ as follows. Consider the commuting square
\begin{equation*}\begin{tikzcd}
    \Tw(\cat I \times \cat C) \arrow[r] \arrow[d] & \Tw(\cat I \times \cat D) \arrow[d]\\
    \Tw(\cat I) \times \cat C^\op \times \cat C \arrow[r] & \Tw(\cat I) \times \cat D^\op \times \cat D
\end{tikzcd}\end{equation*}
where the horizontal maps are induced by the functor $\cat I \times \cat C \to \cat D$ coming from $\ell$; here we implicitly use that $\Tw(-)$ commutes with finite products. Now the left vertical map in the square is the unstraightening of $F_1$, while the fiber product of the lower right part of the square is the unstraightening of $F_2$. Thus the square provides the desired natural transformation $\alpha$.

With $\alpha$ at hand, we now reinterpret $F_1$ and $F_2$ as functors $\Tw(\cat I) \to \Fun(\cat C, \cat P(\cat C))$. We note that both $F_1$ and $F_2$ factor over the full subcategory $\Fun(\cat C, \cat C) \subseteq \Fun(\cat C, \cat P(\cat C))$ and that $F_1$ corresponds to the constant functor $\id_{\cat C}$ and $F_2$ corresponds to the functor $[S \to T] \mapsto r_S \ell_T$. Thus $\alpha$ provides the desired natural transformation.
\end{proof}

The next result provides a comparison result for sheaves on a subsite of a site and is for example used in the paper in order to compare classical stacks to stacks. We refer the reader to \cite[\S A.4]{heyer-mann-6ff} for a quick summary of sheaf theory.

\begin{lem} \label{rslt:embedding-of-sheaves-on-subsite}
Let $\cat C_0$ and $\cat C$ be sites and $f\colon \cat C_0 \injto \cat C$ a fully faithful embedding of categories. Assume that the following is satisfied:
\begin{itemize}
    \item[($*$)] $\cat C_0$ admits finite limits and for every $X \in \cat C_0$, a sieve on $X$ in $\cat C$ is a covering sieve if and only if it is generated by a family of maps $(U_i \to X)_i$ in $ \cat C_0$ such that $f$ preserves all fiber products of the form $U_{i_1} \times_X \dots \times_X U_{i_n}$.
\end{itemize}
Then there is a canonical adjunction
\begin{align*}
    f^*\colon \Shv(\cat C_0) \rightleftarrows \Shv(\cat C) \noloc f_*
\end{align*}
With the following properties:
\begin{lemenum}
    \item Both $f^*$ and $f_*$ preserve all small colimits and $f_*$ preserves all small limits.
    \item $f^*$ is fully faithful.
    \item $f_*$ preserves injective and surjective maps of sheaves.
\end{lemenum}
Assume furthermore that $f$ preserves all finite limits. Then the following is true:
\begin{lemenum}
    \setcounter{lemenumi}{3}
    \item $f^*$ preserves all finite limits.
    \item $f^*$ preserves and detects injective and surjective maps of sheaves.
\end{lemenum}
\end{lem}
\begin{proof}
The embedding $f$ induces an embedding $f^*\colon \cat P(\cat C_0) \injto \cat P(\cat C)$ on presheaf categories via left Kan extension. Composing this functor with the sheafification functor $\cat P(\cat C) \to \Shv(\cat C)$ yields a colimit-preserving functor $\cat P(\cat C_0) \to \Shv(\cat C)$. By \cite[Proposition~5.5.4.20]{lurie-higher-topos-theory} and ($*$) this functor factors uniquely over the sheafication functor $\cat P(\cat C_0) \to \Shv(\cat C_0)$, providing the desired colimit-preserving functor $f^*\colon \Shv(\cat C_0) \to \Shv(\cat C)$. By the adjoint functor theorem $f^*$ admits a right adjoint $f_*$ and this functor can be explicitly computed by restriction along the embedding $\cat C_0 \injto \cat C$.

We now make the following crucial observation: The functor $f_*$ commutes with sheafification on $\cat C$ and $\cat C_0$. To see this, recall the following explicit formula for sheaficiation: Given a presheaf $\mathcal F\colon \cat C^\op \to \Ani$, the sheafification of $\mathcal F$ is obtained by transfinitely applying the construction
\begin{align*}
    \mathcal F \mapsto \mathcal F^\dagger := \big[ X \mapsto \varinjlim_{\cat C_{/X}^{(0)}} \varprojlim_{Y \in \cat C_{/X}^{(0)}} \mathcal F(Y) \big],
\end{align*}
where the colimit is taken over all covering sieves of $X$ in $\cat C$ (see \cite[Remark~6.2.2.12]{lurie-higher-topos-theory} and the proof of \cite[Proposition~6.2.2.7]{lurie-higher-topos-theory}). A similar formula holds for sheafification on $\cat C_0$. Thus it is enough to show that $f_*\colon \cat P(\cat C) \to \cat P(\cat C_0)$ commutes with the functor $\mathcal F \mapsto \mathcal F^\dagger$. Since $f_*$ preserves all small colimits of presheaves, this follows easily from ($*$) and the explicit formula for limits over covering sieves in \cite[Lemma~A.4.6]{heyer-mann-6ff}.

With the above claim for $f_*$ at hand, it is now easy to prove (i)--(v). In (i) the only non-trivial claim is that $f_*$ preserves small colimits. But these are computed by composing the small colimits in $\cat P(\cat C)$ and $\cat P(\cat C_0)$ with the sheafification functor, and $f_*$ commutes with sheafification. In (ii) we need to prove that the adjunction unit $\id \to f_* f^*$ is an equivalence. This is clear on the level of presheaves and hence follows easily from the commutation of $f_*$ with sheafification. Part (iii) follows formally from the fact that $f_*$ preserves (finite) limits and colimits, as both injective maps (i.e. monomorphisms) and surjective maps (i.e. effective epimorphisms) are defined using these operations. Now assume that $f$ preserves finite limits. Since sheafification commutes with finite limits (which is easy to see from the above formula for it, see \cite[Proposition~6.2.2.7]{lurie-higher-topos-theory}) it is easy to deduce that $f^*$ commutes with finite limits (cf. \cite[Proposition~6.2.3.20]{lurie-higher-topos-theory}), proving (iv). Then (v) follows in the same way as (iii).
\end{proof}

The next results study some basic properties of the tensor product of presentable categories (defined in \cite[Proposition~4.8.1.15]{lurie-higher-algebra}) and more specifically the relative version over an algebra in $\PrL$.

\begin{lem} \label{rslt:dualizable-ff-trick}
Fix some $\cat{A} \in \CAlg(\PrL)$, and let $i\colon \cat{M}_1 \injto \cat{M}_2$ be a colimit-preserving fully faithful $\cat{A}$-linear functor of $\cat{A}$-module categories. Let $\cat{N}\in\Mod_{\cat{A}}(\PrL)$ be some other $\cat{A}$-module category. Assume that $\cat{N}$ is dualizable as an $\cat{A}$-module, or that $i$ admits a colimit-preserving and $\cat{A}$-linear left or right adjoint. Then the canonical functor
\begin{align*}
    i \tensor \id\colon \cat{M}_1 \tensor_{\cat{A}} \cat{N} \injto \cat{M}_2 \tensor_{\cat{A}} \cat{N}    
\end{align*}
is fully faithful.
\end{lem}
\begin{proof}
We first treat the case where $\cat{N}$ is dualizable as an $\cat{A}$-module. Under this hypothesis, the tensor product $\cat{M}_i \tensor_{\cat{A}} \cat{N}$ agrees with the functor category $\Fun^L_{\cat{A}}(\cat{N}^\vee, \cat{M}_i)$. Now by general nonsense, $\Fun^{L}_{\cat{A}}(\cat{N}^\vee, \cat{M}_i)$ is computed as the totalization of the cosimplicial object $[n] \mapsto \Fun^L(\cat{N}^\vee \tensor \cat{A}^{\otimes n}, \cat{M}_i)$, and the functor $i \tensor \mathrm{id}_{\cat{N}}$ is the totalization of the obvious functor between these cosimplicial categories. Since $\Fun^L(\cat{C},-)$ preserves full faithfulness for any $\cat{C} \in \Pr^L$, the individual functors
\begin{align*}
    \Fun^{L}(\cat{N}^\vee \tensor \cat{A}^{\tensor n}, \cat{M}_1) \to \Fun^{L}(\cat{N}^\vee \tensor \cat{A}^{\tensor n}, \cat{M}_2)
\end{align*}
induced by $i$ are fully faithful. Since full faithfulness is stable under limits (see \cite[Lemma~A.1.17]{heyer-mann-6ff}), we conclude.

In the case where $i$ has a colimit-preserving $\cat{A}$-linear left resp.\ right adjoint, we observe that the same is then true for $i \tensor \id$. Thus the full faithfulness reduces to checking that the unit resp.\ counit of the adjunction is an equivalence. But this property is preserved under the 2-functor $- \tensor_{\cat A} \cat N$.
\end{proof}

\begin{rmk} \label{rmk:how-to-apply-dualizable-ff-trick}
In the setting of \cref{rslt:dualizable-ff-trick}, note that if $\cat{A}$ is rigid or locally rigid (e.g. $\cat A = \D(\Lambda)$ for some ring $\Lambda$) then the dualizability of $\cat{N}$ over $\cat{A}$ is equivalent to the abstract dualizability of $\cat{N}$ (which holds e.g. if $\cat N$ is compactly generated), see \cite[Chapter~1, Proposition 9.4.4]{gaitsgory-rozenblyum-vol1}. This is the only case we will need.
\end{rmk}

The following results shows that the relative tensor product often preserves the property of being compactly generated. They are a version of \cite[Chapter~1, Corollary~8.7.4]{gaitsgory-rozenblyum-vol1}.

\begin{lem} \label{rslt:compact-generation-of-relative-tensor}
Fix some stable $\cat A \in \Alg(\PrL)$ and let $\cat M$ and $\cat N$ be $\cat A$-linear presentable categories. Assume that $\cat A$, $\cat M$ and $\cat N$ are compactly generated and satisfy the following condition
\begin{itemize}
    \item[($*$)] The tensor unit $\mathbf 1 \in \cat A$ is compact and for compact objects $a_1, a_2 \in \cat A$, $m \in \cat M$ and $n \in \cat N$, also $a_1 \tensor a_2 \in \cat A$, $a_1 \tensor m \in \cat M$ and $a_2 \tensor n \in \cat N$ are compact.
\end{itemize}
Then $\cat M \tensor_{\cat A} \cat N$ is compactly generated by objects of the form $m \boxtimes n$ for compact $m \in \cat M$ and $n \in \cat N$.
\end{lem}
\begin{proof}
We have $\cat M \tensor_{\cat A} \cat N = \varinjlim_{[k] \in \Delta^\op} \cat M \tensor \cat A^{\tensor k} \cat N$. By \cite[Lemma~5.3.2.11]{lurie-higher-algebra} all terms in this colimit are compactly generated, hence by \cite[Lemma~5.3.2.9(3)]{lurie-higher-algebra} it is enough to show that all transition functors in the colimit diagram preserve compact objects. Since all the categories are stable, this is equivalent to having colimit-preserving right adjoints. This then reduces to the claim that the unit map $\Sp \to \cat A$ and the multiplication maps $\cat A \tensor \cat A \to \cat A$, $\cat A \tensor \cat M \to \cat M$, $\cat A \tensor \cat N \to \cat N$ preserve compact objects. By the explicit description of compact generators in these tensor products (see \cite[Lemma~5.3.2.11(1)]{lurie-higher-algebra}) the claim immediately reduces to condition ($*$).

It remains to see that the objects $m \tensor n$ generate $\cat M \tensor_{\cat A} \cat N$. We know this to be true for $\cat M \tensor \cat N$ by \cite[Lemma~5.3.2.11(1)]{lurie-higher-algebra}, hence it is enough to show that the right adjoint of the canonical map $\cat M \tensor \cat N \to \cat M \tensor_{\cat A} \cat N$ is conservative. But by writing the above colimit over $\Delta^\op$ as a limit over $\Delta$ (via passage to right adjoints), we see that the functor $\cat M \tensor_{\cat A} \cat N \to \cat M \tensor \cat N$ is simply the projection to the component of $[0] \in \Delta$, and this is clearly conservative (as any map in $\Delta$ receives a map from $[0]$).
\end{proof}

\begin{cor} \label{rslt:compact-generation-of-tensor-over-Lambda}
Fix a ring $\Lambda$ and let $\cat M$ and $\cat N$ be $\D(\Lambda)$-linear presentable categories. If both $\cat M$ and $\cat N$ are compactly generated then so is $\cat M \tensor_{\D(\Lambda)} \cat N$. Moreover, a set of compact generators is given by $m \boxtimes n$ for compact $m \in \cat M$ and $n \in \cat N$.
\end{cor}
\begin{proof}
We only need to show that condition ($*$) in \cref{rslt:compact-generation-of-relative-tensor} is satisfied. It is clearly true if $a_1 = a_2 = \mathbf 1 = \Lambda \in \D(\Lambda)$. But every compact object in $\D(\Lambda)$ is obtained from $\Lambda$ via finite colimits and retracts, hence the claim follows for all $a_1$ and $a_2$.
\end{proof}

\subsection{Evaluative dualities}

Here we briefly recall some general nonsense around duality. For this, fix a ring $\Lambda$ with $\D(\Lambda)$ the derived category of $\Lambda$-modules. Let $\mathcal{C}$ be a $\Lambda$-linear compactly generated presentably symmetric monoidal stable category. The compact generation of $\mathcal{C}$ implies that it is dualizable. Recall that a self-duality on $\mathcal{C}$ is the datum of an equivalence $\mathbf{D}:\mathcal{C} \simeq \mathcal{C}^{\vee}$, where $\mathcal{C}^{\vee} = \Fun^{L}_{\mathrm{Sp}}(\mathcal{C},\mathrm{Sp})= \Fun^L_{\D(\Lambda)}(\mathcal{C},\D(\Lambda))$ as usual. This is equivalent to the datum of a contravariant self-equivalence $\mathbf{D}:\mathcal{C}^{\omega} \overset{\sim}{\to} \mathcal{C}^{\omega,\op}$, and can also be described in terms of unit and counit functors $\mu:\D(\Lambda) \to \mathcal{C} \otimes_{\D(\Lambda)} \mathcal{C}$ and $\epsilon: \mathcal{C} \otimes_{\D(\Lambda)} \mathcal{C} \to \D(\Lambda)$ in the usual way.

\begin{defn}A self-duality $\mathbf{D}$ on $\mathcal{C}$ is \emph{evaluative} if there is a (necessarily unique) colimit-preserving functor $\gamma\colon \mathcal{C} \to D(\Lambda)$ such that the counit functor $\epsilon$ satisfies $\epsilon(A \boxtimes B) = \gamma(A \otimes B)$.
\end{defn}

A functor $\gamma$ inducing the counit of a self-duality is often called a \emph{Frobenius structure} (as in e.g.\ \cite[Definition 4.6.5.1]{lurie-higher-algebra}).
If a self-duality is evaluative, there is a tight relation between $\mathbf{D}$ and $\gamma$, and in fact they completely determine each other: for $A$ compact and $B$ arbitrary, we have $\RHom(\mathbf{D}A,B) = \gamma(A \otimes B)$. Moreover, for an evaluative self-duality, the equivalence $\mathcal{C} \simeq \mathcal{C}^{\vee}$ is completely transparent: an object $C \in \mathcal{C}$ corresponds to the colimit-preserving functor $X \mapsto \gamma(C \otimes X)$.

Now suppose $\mathcal{C},\mathcal{D}$ are $\Lambda$-linear categories satisfying the above conditions and equipped with self-dualities. Then any colimit-preserving $\Lambda$-linear functor $F:\mathcal{C} \to \mathcal{D}$ induces a dual functor $F^\vee : \mathcal{D}^{\vee} \to \mathcal{C}^{\vee}$. Pre- and post-composing with the given self-dualities on $\mathcal{C}$ and $\mathcal{D}$, we can and will view the dual as a functor $F^\vee: \mathcal{D} \to \mathcal{C}$. The following lemma is then a straightforward computation.

\begin{lem}\label{lem:evaluativedualitything}If $F:\mathcal{C}\to\mathcal{D}$ is any colimit-preserving functor, and the given self-dualities on $\mathcal{C}$ and $\mathcal{D}$ are evaluative with associated functors $\gamma_{\mathcal{C}}$ and $\gamma_{\mathcal{D}}$, then $F^\vee : \mathcal{D}\to \mathcal{C}$ is characterized by the formula
\[\gamma_{\mathcal{C}}(C \otimes F^\vee (D)) = \gamma_{\mathcal{D}}(F(C) \otimes D)\]
bifunctorially in all $C \in \mathcal{C}$ and $D\in \mathcal{D}$.
\end{lem}

Let us discuss the relevant examples for our purposes.

\begin{exmpl}\label{exmpl:rigiddualizable}
    Suppose $\mathcal{C}$ is any category satisfying the general conditions above, which is moreover \emph{rigid dualizable}. In other words, we assume $\mathcal{C}$ is compactly generated by dualizable compact objects and with compact unit $1_{\mathcal{C}}$. (See \cite[Chapter 1, Lemma 9.1.5]{gaitsgory-rozenblyum-vol1}.) Then $\mathcal{C}$ is canonically self-dual via the naive dual $X\mapsto X^\vee$ on dualizable compact objects. This self-duality is evaluative, with $\gamma(-)=\RHom(1_{\mathcal{C}},-)$.
\end{exmpl}

\begin{exmpl}\label{exmpl:QCohselfduality}
    Assume $\Lambda=k$ is a field with $\mathrm{char} k =0$, and let $X$ be a disjoint union of perfect QCA stacks over $k$. Then $\QCoh(X)$ is canonically self-dual via the naive dualizability of all compact objects $X \in \Perf^{\mathrm{qc}}(X)$. This self-duality is evaluative, with $\gamma(-) = \Gamma_!(X,-)$ the unique colimit-preserving functor $\QCoh(X) \to D(k)$ extending the functor $\mathcal{F} \mapsto \Gamma(X,\mathcal{F})$ on $\Perf^\qc(X)$.

    When $X$ is quasicompact then $\Gamma_!(X,-) = \RHom(\mathcal{O}_X,-)$ and this is a special case of the previous example, but in general it is distinct. We will need this extra generality, since $\Par_G$ is not quasicompact.
    
    If $f:X\to Y$ is any quasicompact morphism of such stacks, the functors $f_{\ast}$ and $f^\ast$ are dual to each other: there are canonical identifications $(f^{\ast})^\vee = f_{\ast}$ and $(f_{\ast})^\vee = f^\ast$. This follows easily from Lemma \ref{lem:evaluativedualitything} via a simple projection formula argument.
\end{exmpl}

\begin{exmpl}
    Assume $\Lambda=k$ is a field with $\mathrm{char} k =0$, and let $X$ be any disjoint union of QCA stacks over $k$. Then $\ICoh(X)$ is canonically self-dual via the contravariant self-equivalence $\Dgs\colon \Coh^\qc(X) \overset{\sim}{\to} \Coh^\qc(X)^{\op}$ induced by Grothendieck-Serre duality. This self-duality is evaluative, with $\gamma(-) = \Gamma_!(X,\Psi(-))$.

    For non-smooth $X$, this example is disjoint from the previous two examples, since $\IndCoh(X)$ is never rigid dualizable for non-smooth $X$.

    If $f\colon X\to Y$ is any quasicompact lafp morphism of such stacks, the functors $f_!$ and $f^*$ on $\ICoh$ are dual to each other. This again follows from the projection formula.
\end{exmpl}

\begin{exmpl} \label{exmpl:ICoh-self-duality}
   Assume $\Lambda=k$ is a field with $\mathrm{char} k =0$, and let $X$ be any perfect QCA stack over $k$, so we have the natural colimit-preserving functor $\Psi_{X}: \IndCoh(X) \to \QCoh(X)$ given as the ind-completion of the evident inclusion $\Coh^\qc(X) \to \QCoh(X)$. Then the source and target of $\Psi_{X}$ are canonically self-dual via the previous two examples, and the dual $\Psi_{X}^{\vee}$ is nothing more than the symmetric monoidal embedding $\gamma_X: \QCoh(X) \to \IndCoh(X)$ induced by the action of $\QCoh(X)$ on the monoidal unit $\omega_X \in \IndCoh(X)$. This is an easy exercise, using the definition of the self-dualities together with the formula $\Psi_X(\gamma_X(\mathcal{F}) \otimes \mathcal{G}) = \mathcal{F} \otimes \Psi_X(\mathcal{G})$.
\end{exmpl}

\begin{exmpl}\label{exmpl:BunGselfduality}
    For $\Lambda$ any $\mathbf{Z}_{\ell}$-algebra, $\D(\Bun_G,\Lambda)$ is canonically self-dual via the self-equivalence on compact objects induced by Bernstein-Zelevinsky duality. This self-duality is evaluative, with $\gamma(-) = \Gamma_c(\Bun_G,-)$. 
\end{exmpl}

\begin{rmk}
Both \cref{exmpl:ICoh-self-duality} and \cref{exmpl:BunGselfduality} are a special case of prim duality (for $\ICoh$ and $\D_\rel$, respectively) together with the fact that in these examples the 2-functor $\cat K_{\D} \to \Mod_{\D(*)}(\PrL)$ is fully faithful on the full subcategory spanned by $*$ and $\Par_G$ resp.\ $\Bun_G$. In fact, also \cref{exmpl:QCohselfduality} is a special case of this, for the 6-functor formalism on $\QCoh$ which has $f_! = f_*$ for all qcqs QCA maps.
\end{rmk}

\subsection{Filtered and graded rings} \label{sec:filtered-rings}

In the following we describe a general procedure of constructing the category of \enquote{animated/derived rings} in a base category $\cat C$. In the case $\cat C = \D(\Z)$ we recover the usual category $\Ring$ of animated rings that forms the basis of (derived) algebraic geometry, but the following formalism also captures filtered and graded variants. All definitions and constructions are based on \cite{Raksit-derived-rings} and were inspired by discussions with Zachary Gardner and Jeroen Hekking (their recent work \cite{GardnerHekking-Ideals} contains similar ideas as below).

Let us start by introducing the necessary structure on the base category $\cat C$ in order to define a good category of animated rings in $\cat C$. 

\begin{defn}
Let $\cat A$ be an abelian symmetric monoidal category with an object $X \in \cat A$.
\begin{defenum}
    \item \label{def:symmetric-power} Let $n \ge 0$. The \emph{$n$-th symmetric power} $\Sym_{\cat A}^n(X) \in \cat A$ of $X$ is given by the $S_n$-coinvariants of $X^{\tensor n}$, where the permutation group $S_n$ acts on $X^{\tensor n}$ via permuting the coordinates.

    \item Suppose that $\cat A$ admits countable direct sums. Then the \emph{symmetric algebra} associated to $X$ is the algebra $\Sym_{\cat A}(X) \in \CAlg(\cat A)$ whose underlying module is $\bigoplus_{n\ge0} \Sym^n_{\cat A}(X)$ and where multiplication is given by the obvious maps $\Sym_{\cat A}^n(X) \tensor \Sym_{\cat A}^m(X) \to \Sym_{\cat A}^{m+n}(X)$.
\end{defenum}
\end{defn}

\begin{defn}[{see \cite[Definition~4.2.1]{Raksit-derived-rings}}]
A \emph{derived algebraic context} is a stable presentable symmetric monoidal category $\cat C$ with a t-structure together with a full subcategory $\cat C^0 \subseteq \cat C^\heartsuit$ satisfying the following conditions:
\begin{enumerate}[(i)]
    \item The t-structure on $\cat C$ is right-complete and $\cat C^{\ge0}$ is stable under filtered colimits in $\cat C$. Moreover, $\cat C^{\le0}$ is stable under the tensor product and contains the tensor unit of $\cat C$. In particular $\cat C^\heartsuit$ inherits a tensor product via $X \tensor_{\cat C^\heartsuit} Y := H^0(X \tensor_{\cat C} Y)$.

    \item The subcategory $\cat C^0 \subseteq \cat C$ is stable under the tensor product and under the formation of symmetric powers in $\cat C^\heartsuit$.

    \item $\cat C^0$ is stable under finite coproducts in $\cat C$ and its objects form a set of compact projective generators for $\cat C^{\le0}$.
\end{enumerate}
A \emph{(lax) morphism between derived algebraic contexts} $\cat C$ and $\cat D$ is a colimit-preserving right t-exact (lax) symmetric monoidal functor $F\colon \cat C \to \cat D$ such that $F(\cat C^0) \subseteq \cat D^0$. 
\end{defn}

\begin{rmk} \label{rmk:DAG-morphisms-obtained-from-ordinary-categories}
As observed in \cite[Remark~4.2.2]{Raksit-derived-rings}, a derived algebraic context is completely determined by the subcategory $\cat C^0$. Namely, by the compact projective generation we obtain $\cat C^{\le0}$ from $\cat C^0$ via animation, i.e. $\cat C^{\le0} = \cat P_\Sigma(\cat C^0)$, and we then obtain $\cat C$ from $\cat C^{\le0}$ via passing to spectrum objects. The assignment $\cat C \mapsto \cat C^0$ thus defines a fully faithful embedding of the 2-category of derived algebraic contexts into the 2-category of symmetric monoidal small categories. As the latter category is an ordinary 2-category, this allows us to make many constructions with derived algebraic contexts very explicitly.
\end{rmk}

\begin{rmk}
Unlike \cite[Definition~4.2.1]{Raksit-derived-rings} we also consider \emph{lax} morphisms of derived algebraic contexts, where the morphism is only required to be \emph{lax} symmetric monoidal. We do this because they still induce associated functors on ring objects and we have several important examples in this paper. However, we warn the reader that $\Ring(-)$ is not compatible with compositions of lax morphisms (see \cref{rmk:warning-about-functoriality-of-Ring} below).
\end{rmk}

Let us introduce some important examples of derived algebraic contexts. As mentioned above, we will mostly focus on filtered and derived objects (cf. \cite[\S3]{Raksit-derived-rings}).

\begin{exmpl}
The category $\D(\Z)$ equipped with the standard t-structure and the usual tensor product is a derived algebraic context if we let $\D(\Z)^0 \subseteq \D(\Z)^\heartsuit = \Ab$ be the full subcategory spanned by the objects $\Z^n$ for $n \ge 0$.
\end{exmpl}

\begin{defn} \label{def:filtered-and-graded-algebraic-contexts}
Let $\cat C$ be a derived algebraic context.
\begin{defenum}
    \item We regard $\Z$ as a category via the usual partial order. We define
    \begin{align*}
        \Fil(\cat C) := \Fun(\Z^\op, \cat C)
    \end{align*}
    and call it the category of \emph{filtered objects of $\cat C$}. Explicitly, an object in $\Fil(\cat C)$ is a diagram
    \begin{align*}
        X^* = [\dots \to X^2 \to X^1 \to X^0 \to X^{-1} \to X^{-2} \to \dots]
    \end{align*}
    in $\cat C$. Via Day convolution (cf. \cite[Example~2.2.6.9]{lurie-higher-algebra}) the tensor product on $\cat C$ induces a tensor product on $\Fil(\cat C)$ which is computed as
    \begin{align*}
        (X^* \tensor Y^*)^n = \varinjlim_{i+j \ge n} X^i \tensor Y^j.
    \end{align*}
    We equip $\Fil(\cat C)$ with a t-structure by letting $\Fil(\cat C)^{\le0}$ resp. $\Fil(\cat C)^{\ge0}$ be the full subcategories of those filtered objects $X^*$ such that $X^n \in \cat C^{\le0}$ resp. $X^n \in \cat C^{\ge0}$ for all $n \in \Z$. For $n \in \Z$ we obtain an adjunction
    \begin{align*}
        \ins^n\colon \cat C \rightleftarrows \Fil(\cat C) \noloc \ev^n,
    \end{align*}
    where $\ev^n$ is the restriction along $\{n\} \injto \Z^\op$ and $\ins^n$ is the left Kan extension along the same map. Explicitly, $\ev^n(X^*) = X^n$ and $\ins^n(X)^m$ is equal to $X$ for $m \le n$ and equal to $0$ otherwise. Finally, we let $\Fil(\cat C)^0 \subseteq \Fil(\cat C)^\heartsuit$ be the full subcategory spanned by finite coproducts of the objects $\ins^n(X)$ for $n \in \Z$ and $X \in \cat C^0$. Altogether this equips $\Fil(\cat C)$ with the structure of a derived algebraic context (cf. \cite[Construction~4.3.4]{Raksit-derived-rings}).

    \item Let $\Z^{\mathrm{ds}}$ denote the discrete set $\Z$ viewed as a category with no non-trivial morphisms. We define
    \begin{align*}
        \Gr(\cat C) := \Fun(\Z^{\mathrm{ds}}, \cat C) = \prod_{\Z} \cat C
    \end{align*}
    and call it the category of \emph{graded objects of $\cat C$}. An object in $\Gr(\cat C)$ is a collection $X^* = (X^n)_{n\in\Z}$ of objects in $\cat C$. Day convolution defines a tensor product on $\Gr(\cat C)$ which is computed as
    \begin{align*}
        (X^* \tensor Y^*)^n = \bigoplus_{i+j=n} X^i \tensor Y^j.
    \end{align*}
    We equip $\Gr(\cat C)$ with the pointwise t-structure as in (a) and for $n \in \Z$ we denote
    \begin{align*}
        \ins^n\colon \cat C \rightleftarrows \Gr(\cat C) \noloc \ev^n
    \end{align*}
    the similar adjunction, where $\ev^n(X^*) = X^n$ and $\ins^n(X)$ is the graded object which is $X$ at position $n$ and $0$ everywhere else. We let $\Gr(\cat C)^0 \subseteq \Gr(\cat C)^\heartsuit$ be the full subcategory spanned by finite coproducts of $\ins^n(X)$ for $n \in \Z$ and $X \in \cat C^0$. This equips $\Gr(\cat C)$ with the structure of a derived algebraic context (cf. \cite[Construction~4.3.4]{Raksit-derived-rings}).
\end{defenum}
\end{defn}

\begin{defn}
Let $\cat C$ be a derived algebraic context, let $k \ge 0$ be an integer or $\infty$, and let $I_k = [0, k] \subset \Z$. We define the derived algebraic contexts $\Fil^{[0,k]}(\cat C)$ and $\Gr^{[0,k]}(\cat C)$ in a similar way as in \cref{def:filtered-and-graded-algebraic-contexts} with $\Z$ replaced by $I_k$. In the case $k = \infty$ we also denote the categories by $\Fil^{\ge0}(\cat C)$ and $\Gr^{\ge0}(\cat C)$.
\end{defn}

Let us write out explicitly what $\Fil^{[0,1]}(\cat C)$ looks like. An object in that category is a map $X^1 \to X^0$ in $\cat C$. The tensor product is given by
\begin{align*}
    (X^1 \to X^0) \tensor (Y^1 \to Y^0) = [(X^1 \tensor Y^0) \sqcup_{(X^1 \tensor Y^1)} (X^0 \tensor Y^1) \to X^0 \tensor Y^0].
\end{align*}
The objects of $\Fil^{[0,1]}(\cat C)^0$ are the objects of the form $Y \to Y \oplus X$ for $X, Y \in \cat C^0$, where the map is the natural inclusion. The category $\Fil^{[0,1]}(\cat C)$ plays an important role to us because it allows us to define ideals, as we will see below.

Let us now come to the definition of rings in a derived algebraic context. As we only care about connective rings, the definition is straightforward (cf. \cite[Remark~4.2.24]{Raksit-derived-rings}).

\begin{defn}
Let $\cat C$ be an algebraic context.
\begin{defenum}
    \item Let $\cat R^0 \subseteq \CAlg(\cat C^\heartsuit)$ be the full subcategory spanned by the objects $\Sym_{\cat C^\heartsuit}(X)$ for $X \in \cat C^0$. We define
    \begin{align*}
        \Ring(\cat C) := \cat P_{\Sigma}(\cat R^0),
    \end{align*}
    i.e. $\Ring(\cat C)$ is the free sifted colimit completion of $\cat R^0$. The objects of $\Ring(\cat C)$ are called the \emph{rings in $\cat C$}. There is a fully faithful embedding $\CAlg(\cat C^\heartsuit) \subseteq \Ring(\cat C)$ whose image consists of those rings in $\cat C$ where the underlying object lies in $\cat C^\heartsuit$ (cf. \cite[Remark~4.2.24]{Raksit-derived-rings}).

    \item \label{def:Sym-and-forgetful-functor-on-DAG} The functor $\Sym_{\cat C^\heartsuit}\colon \cat C^0 \to \cat R^0$ upgrades to a colimit-preserving functor $\Sym_{\cat C}$ with a right adjoint $U_{\cat C}$:
    \begin{align*}
        \Sym_{\cat C}\colon \cat C^{\le0} \rightleftarrows \Ring(\cat C) \noloc U_{\cat C}.
    \end{align*}
    We call $U_{\cat C}$ the \emph{forgetful functor}. By adjunction properties $U_{\cat C}$ is conservative and preserves all small limits and sifted colimits.

    \item The composition $\cat R^0 \to \CAlg(\cat C^{\heartsuit}) \to \CAlg(\cat C)$ defines a colimit-preserving functor
    \begin{align*}
        (-)^\circ\colon \Ring(\cat C) \to \CAlg(\cat C), \qquad A \mapsto A^\circ
    \end{align*}
    which factors over $\CAlg(\cat C^{\le0})$. We call $A^\circ$ the \emph{underlying $\mathbb E_\infty$-algebra} of $A$.
\end{defenum}
\end{defn}

\begin{lem} \label{rslt:functoriality-of-Ring}
Fix $\cat C, \cat D \in \DAC$. There is a canonical functor
\begin{align*}
    \Fun^{\DAC,\lax}(\cat C, \cat D) \to \Fun^\Sigma(\Ring(\cat C), \Ring(\cat D)),
\end{align*}
where $\Fun^\Sigma \subseteq \Fun$ denotes the full subcategory spanned by sifted-colimit preserving functors and $\Fun^{\DAC,\lax}(\cat C, \cat D)$ denotes the category of lax morphisms of derived algebraic contexts. Moreover, for every lax morphism $f\colon \cat C \to \cat D$ with associated map $f'\colon \Ring(\cat C) \to \Ring(\cat D)$ there are commuting squares
\begin{equation*}
    \begin{tikzcd}
        \Ring(\cat C) \arrow[r,"f'"] \arrow[d,"U_{\cat C}",swap] & \Ring(\cat D) \arrow[d,"U_{\cat D}"]\\
        \cat C \arrow[r,"f"] & \cat D
    \end{tikzcd}
    \qquad\qquad
    \begin{tikzcd}
        \Ring(\cat C) \arrow[r,"f'"] \arrow[d,"(-)^\circ",swap] & \Ring(\cat D) \arrow[d,"(-)^\circ"]\\
        \CAlg(\cat C) \arrow[r,"f"] & \CAlg(\cat D)
    \end{tikzcd}
\end{equation*}
\end{lem}
\begin{proof}
We first observe that $\CAlg = \Hom(\operatorname{Comm}^\tensor, -)$ defines a 2-functor $\Op \to \Cat$ and hence restricts to a 2-functor $\DAC \to \Cat$. In particular we obtain the functor $\Fun^{\DAC,\lax}(\cat C, \cat D) \to \Fun(\CAlg(\cat C), \CAlg(\cat D))$. By definition of lax morphisms of derived algebraic contexts, this functor factors through the full subcategory spanned by those functors $\CAlg(\cat C) \to \CAlg(\cat D)$ that preserve sifted colimits (note that this can be checked on underlying objects) and send $\cat R^0 \subseteq \CAlg(\cat C)$ to $\CAlg(\cat D^\heartsuit) \subseteq \CAlg(\cat D)$ (because $\cat C^0$ is stable under symmetric powers). We obtain the functor
\begin{align*}
    \Fun^{\DAC,\lax}(\cat C, \cat D) \to \Fun(\cat R^0, \CAlg(\cat D^\heartsuit)) \injto \Fun(\cat R^0, \Ring(\cat D)) = \Fun^\Sigma(\Ring(\cat C), \Ring(\cat D)),
\end{align*}
where in the last step we used the universal property of $\cat P_\Sigma$.

It remains to verify the commuting squares. Since all functors in the squares preserve sifted colimits, we can replace $\Ring(\cat C)$ by $\cat R^0$. It is then clear that the right-hand square commutes. The commutativity of the left square follows from this together with the observation that $U_{\cat C}$ factors as the composition of $(-)^\circ$ and the forgetful functor $\CAlg(\cat C) \to \cat C$ (this can again be checked on $\cat R^0$).
\end{proof}

\begin{rmk} \label{rmk:warning-about-functoriality-of-Ring}
We warn the reader that despite \cref{rslt:functoriality-of-Ring}, the assignment $\cat C \mapsto \Ring(\cat C)$ is not functorial in derived algebraic contexts and \emph{lax} morphisms between them, as it is not compatible with compositions of such morphisms in general (this is similar to the fact that taking left derived functors is not compatible with composition in general). If one restricts to (strict) morphisms of derived algebraic contexts, this issue gets resolved.
\end{rmk}


We can use the above definitions to define the category of pairs $(A, I)$ consisting of a ring $A$ and an ideal $I$ in $A$, as follows:

\begin{defn}
Let $\cat C$ be a derived algebraic context. We define
\begin{align*}
    \Pair(\cat C) := \Ring(\Fil^{[0,1]}(\cat C))
\end{align*}
and call its objects the \emph{pairs in $\cat C$}.
\end{defn}

We refer the reader to \cite[Definition~3.21]{Mao-derived-crystalline-cohom} for a similar definition in the case that $\cat C = \D(\Z)$. Explicitly, a pair in $\cat C$ consists of a map $I \to A$ in $\cat C^{\le0}$ and a multiplication map of the form
\begin{align*}
    (I \to A)^{\tensor2} = ((I \tensor A) \sqcup_{I^{\tensor2}} (A \tensor I) \to A \tensor A) \to (I \to A)
\end{align*}
together with higher associative and strictly commutative coherences. In particular we get a multiplication map $A \tensor A \to A$ which makes $A$ into a ring; more formally the projection $\Fil^{[0,1]}(\cat C) \to \cat C$ of derived algebraic contexts induces a functor $\Pair(\cat C) \to \Ring(\cat C)$ sending $(A, I) \mapsto A$. Classically, to a pair $(A, I)$ one can associate the ring $A/I$ together with a surjective ring map $A \surjto A/I$. The same is true in the derived world (see \cite[Theorem~3.29]{Mao-derived-crystalline-cohom} for a similar result):

\begin{prop} \label{rslt:equivalence-of-pairs-and-ring-maps}
Let $\cat C$ be a derived algebraic context. There is a natural fully faithful embedding
\begin{align*}
    \Pair(\cat C) \injto \Fun([1], \Ring(\cat C)),
\end{align*}
whose essential image consists of those ring maps $A \to B$ such that the induced map $\pi_0 A \to \pi_0 B$ is surjective (equivalently, $\fib(A \to B) \in \cat C^{\le0}$). On underlying objects in $\cat C$, the above embedding sends $(A, I)$ to $A \to A/I$, where $A/I := \cofib(I \to A)$; in particular $A/I$ is naturally a ring in $\cat C$.
\end{prop}
\begin{proof}
We equip the category $\Fun([1], \cat C)$ with the structure of a derived algebraic context as follows. The tensor product is given pointwise. The t-structure is defined by letting $\Fun([1], \cat C)^{\le0}$ be the full subcategory of those maps $M \to N$ such that $M \in \cat C^{\le0}$ and $\fib(M \to N) \in \cat C^{\le0}$, and similarly for $(-)^{\ge0}$. We let $\Fun([1], \cat C)^0$ be the full subcategory spanned by objects $M \to N$ such that the induced map $\fib(M \to N) \to M$ is of the form $X \to X \oplus Y$ for $X, Y \in \cat C^0$. We claim that these definitions indeed induce the structure of a derived algebraic context on $\Fun([1], \cat C)$ and that there is a natural equivalence of derived algebraic contexts
\begin{align*}
    \Fil^{[0,1]}(\cat C) \isoto \Fun([1], \cat C), \qquad [Y \to X] \mapsto [X \to \cofib(Y \to X)].
\end{align*}
In fact, it is enough to show that the above formula defines an equivalence of symmetric monoidal categories, because then it is clearly compatible with the structures of derived algebraic contexts and in particular $\Fun([1], \cat C)$ is one. It is clear that the formula defines an equivalence of categories, so it only remains to construct a symmetric monoidal structure on the functor, for which we refer the reader to \cite[Proposition~3.1]{Mao-derived-crystalline-cohom}. Altogether we deduce that there is a natural isomorphism
\begin{align*}
    \Pair(\cat C) = \Ring(\Fun([1], \cat C)).
\end{align*}
It remains to identify the right-hand side with the expected full subcategory of $\Ring(\cat C)$. By going through the definitions we see that $\Ring(\Fun([1], \cat C))$ is freely generated under sifted colimits by the maps $s_{X,Y}\colon \Sym_{\cat C^\heartsuit}(X \oplus Y) \to \Sym_{\cat C^\heartsuit}(Y)$ for $X, Y \in \cat C^0$, which immediately induces a functor $\Ring(\Fun([1], \cat C)) \to \Fun([1], \Ring(\cat C))$ that commutes with sifted colimits. To see that this functor is fully faithful we have to show that all $s_{X,Y}$ are compact projective in $\Fun([1], \Ring(\cat C))$. By the adjunction of $\Sym_{\cat C}$ and the forgetful functor (which is preserved under $\Fun([1],-)$) we deduce
\begin{align*}
    \Hom_{\Fun([1],\Ring(\cat C)}(s_{X,Y}, [A \to B]) &= \Hom_{\Fun([1],\cat C)}([X \oplus Y \to Y], [A \to B])\\
    &= \Hom_{\cat C}(Y, A) \times \Hom_{\cat C}(X, \fib(A \to B)).
\end{align*}
This immediately shows that $s_{X,Y}$ is compact projective. Moreover, it also shows that the family of functors $\Hom(s_{X,Y}, -)$ is conservative on the category of ring maps $A \to B$ such that $\fib(A \to B) \in \cat C^{\le0}$, which says that the image of the above embedding is exactly this subcategory.
\end{proof}

We now discuss a version of \cref{rslt:equivalence-of-pairs-and-ring-maps} for $\Z_{\ge0}$-filtered rings. Let us first introduce the relevant notation.

\begin{defn}
Given a derived algebraic context $\cat C$, we denote
\begin{align*}
    \FilgeRing(\cat C) := \Ring(\Fil^{\ge0}(\cat C))
\end{align*}
and call its objects the \emph{$\Z_{\ge0}$-filtered rings in $\cat C$}.
\end{defn}

Thus a filtered ring in $\cat C$ is a sequence $\dots \to I^2 \to I^1 \to I^0 =: A$ in $\cat C$ together with a certain commutative algebra structure. By the description of commutative algebras in the Day convolution (see \cite[Example~2.2.6.9]{lurie-higher-algebra}) we see that there are natural maps $I^n \tensor I^m \to I^{n+m}$ for all $n, m \ge 0$. We will often denote objects in $\FilgeRing(\cat C)$ by $(A, I^\bullet)$.

\begin{prop} \label{rslt:quotients-from-filtered-ring}
Let $\cat C$ be a derived algebraic context. There is a canonical functor
\begin{align*}
    \FilgeRing(\cat C) \to \Fun(\Z_{\ge0}^{\triangleright,\op}, \Ring(\cat C)), \qquad (A, I^\bullet) \mapsto (A/I^n)_n
\end{align*}
such that on underlying objects in $\cat C$ we have $A/I^n = \cofib(I^n \to A)$.
\end{prop}
\begin{proof}
We first construct the functor
\begin{align*}
    F\colon \FilgeRing(\cat C) \to \Fun(\Z_{\ge0}^\op, \Pair(\cat C)), \qquad (A, I^\bullet) \mapsto (A, I^n)_{n\ge0}.
\end{align*}
This functor is equivalently a functor $\Z_{\ge0}^\op \to \Fun^\Sigma(\FilgeRing(\cat C), \Pair(\cat C))$, where $\Fun^\Sigma$ denotes the category of functors that preserve sifted colimits. By \cref{rslt:functoriality-of-Ring} it is enough to construct the functor $\Z_{\ge0}^\op \to \Fun^{\DAC,\lax}(\Fil^{\ge0}(\cat C), \Fil^{[0,1]}(\cat C))$. By the functoriality of Day convolution, lax symmetric monoidal functors $\Fil^{\ge0}(\cat C) \to \Fil^{[0,1]}(\cat C)$ are induced by lax symmetric monoidal functors $[1]^\op \to \Z_{\ge0}^\op$, where $[1] = \{ 0, 1 \}$ is equipped with the monoid structure such that multiplication is given by taking the maximum. Considering the functors $[1]^\op \to \Z_{\ge0}^\op$ that send $0 \mapsto 0$ and $1 \mapsto n$ for $n \ge 0$, we arrive at the desired diagram. This finishes the construction of $F$.

We now compose $F$ with the equivalence from \cref{rslt:equivalence-of-pairs-and-ring-maps} in order to get a functor $\FilgeRing(\cat C) \to \Fun(\Z_{\ge0}^\op, \Fun([1], \Ring(\cat C)))$ and then note that the projection to $0 \in [1]$ is the constant functor $(A, I^\bullet) \mapsto A$. By some formal transformations we arrive at the desired functor in the claim.
\end{proof}

\subsection{Extended geometric setups} \label{sec:extended-geom-setups}

In \cite[\S3.1]{heyer-mann-6ff}, 3-functor formalisms are defined as lax monoidal functors $\Corr(\cat C, E)^\tensor \to \Cat^\times$ for a geometric setup $(\cat C, E)$. Such a 3-functor formalism defines a pullback functor $f^*$ for every map in $\cat C$ and a !-functor $f_!$ for every map in $E$. In the following we need a slight generalization of this concept, where pullbacks are only defined for a certain subclass of maps in $\cat C$. This will be important when constructing the 3-functor formalism of bounded sheaves in \cref{sec:Corr-bounded}.

\begin{defn} \label{def:extended-geometric-setup}
An \emph{extendeded geometric setup} is a triple $(\cat C, B, F)$ consisting of a category $\cat C$ and collections of edges $B$ and $F$ in $\cat C$, called the \emph{backward} and \emph{forward} morphisms respectively, such that the following properties are satisfied:
\begin{enumerate}[(i)]
    \item $(\cat C, F)$ is a geometric setup, i.e. $F$ is stable under composition, pullback and diagonals (in particular it is right cancellative, cf. \cite[Lemma~2.1.5]{heyer-mann-6ff}).

    \item $B$ contains all isomorphisms and is stable under composition and pullback.

    \item $\cat C$ admits finite products.
\end{enumerate}
\end{defn}

The definition of extended geometric setups is related to but more restrictive than the definition of adequate triples from \cite[Definition~2.1]{Haugseng-Hebestreit.Spans}. We impose the right cancellativeness hypothesis on $F$ in order to be compatible with the definition of geometric setups from \cite[Definition~2.1.1]{heyer-mann-6ff} (see the remarks and lemmas following that definition) and we impose pullback stability of $B$ along arbitrary maps and the existence of finite products in $\cat C$ in order to construct a well-defined operad structure on the correspondence category (cf. \cref{rmk:subtle-definition-of-extended-Corr-operad} below). Note that an extended geometric setup $(\cat C, B, F)$ with $B = all$ is the same as a geometric setup such that $\cat C$ has finite limits.

\begin{defn}
Let $(\cat C, B, F)$ be an extended geometric setup. We denote by $\Corr(\cat C, B, F)$ the \emph{category of correspondences}, which is the full subcategory of $\Corr(\cat C, F)$ containing the same objects, but only those morphisms where the backward map lies in $B$. Explicitly, its objects are those of $\cat C$ and for $X, Y \in \cat C$ a morphism from $X$ to $Y$ in $\Corr(\cat C, B, F)$ is given by a diagram
\begin{equation*}\begin{tikzcd}[column sep=small]
    & Z \ar[dl,"b",swap] \ar[dr,"f"]\\
    X && Y
\end{tikzcd}\end{equation*}
where $b \in B$ and $f \in F$. Composition is given by pullback diagrams.
\end{defn}

\begin{defn} \label{def:Corr-operad-for-extended-geometric-setup}
Let $(\cat C, B, F)$ be an extended geometric setup. Consider the cocartesian operad $(\cat C^\op)^\sqcup$ (see e.g.\ \cite[Example~B.1.7]{heyer-mann-6ff}). We define the classes of edges $B^+$ and $F^-$ in $(\cat C^\op)^\sqcup$ as follows:
\begin{itemize}
    \item Let $f\colon Y_\bullet \to X_\bullet$ be a map in $(\cat C^\op)^\sqcup$ lying over a map $\alpha\colon \langle n \rangle \to \langle m \rangle$ in $\Fin_*$. Then $f$ lies in $B^+$ if for all $j = 1, \dots, m$ the induced map $X_j \to \prod_{i\in \alpha^{-1}(j)} Y_i$ lies in $B$.

    \item The edges of $F^-$ are the same as in \cite[Lemma~2.3.1]{heyer-mann-6ff}, i.e. they are those morphisms $Y_\bullet \to X_\bullet$ in $(\cat C_F^\op)^\sqcup$ such that $Y_\bullet$ and $X_\bullet$ are tuples of the same size and the map is given by a tuple of maps $Y_i \to X_i$ in $F$.
\end{itemize}
With this notation at hand, we define
\begin{align*}
    \Corr(\cat C, B, F)^\tensor := \Corr((\cat C^\op)^{\sqcup,\op}, B^+, F^-).
\end{align*}
The same argument as in \cite[Proposition~2.3.3]{heyer-mann-6ff} shows that $\Corr(\cat C, B, F)^\tensor$ is an operad and if $\cat C_B$ admits finite products that coincide with those in $\cat C$ then $\Corr(\cat C, B, F)$ is a symmetric monoidal category.
\end{defn}

\begin{rmk} \label{rmk:subtle-definition-of-extended-Corr-operad}
Our definition of $\Corr(\cat C, B, F)^\tensor$ differs slightly from the one in \cite[Proposition~6.1.2]{liu-zheng-enhanced-operations}, as we have a different condition on the edges in $B^+$. Even if $B$ is stable under diagonals (which it is in our applications), our condition differs: For example, a map $\emptyset \to X$ in $(\cat C^\op)^\sqcup$ over $\alpha\colon \langle 0 \rangle \to \langle 1 \rangle$ lies in $B^+$ if and only if the induced map $X \to *$ lies in $B$. We chose our definition in order to ensure that $\Corr(\cat C, B, F)^\tensor$ is a category, for which we need to ensure that morphisms in $B^+$ are stable under pullbacks along morphisms in $E^-$ -- see the proof of \cite[Lemma~2.3.1]{heyer-mann-6ff} for an explicit computation of these pullbacks. In particular it is unclear to us whether \cite[Proposition~6.1.2]{liu-zheng-enhanced-operations} is correct as stated.
\end{rmk}

We warn the reader that in the applications we have in mind, $\Corr(\cat C, B, F)^\tensor$ is not a symmetric monoidal category, because $\cat C_B$ does not have finite products: It does not admit a final object.

\begin{defn} \label{def:3ff-on-extended-geometric-setup}
A \emph{3-functor formalism} on an extended geometric setup $(\cat C, B, F)$ is a map
\begin{align*}
    \Corr(\cat C, B, F)^\tensor \to \Cat^\times
\end{align*}
of operads.
\end{defn}

\noindent A 3-functor formalism $\D\colon \Corr(\cat C, B, F)^\tensor \to \Cat^\times$ encodes the following data:
\begin{itemize}
    \item For every $X \in \cat C$ a category $\D(X)$.
    \item For every $b\colon Y \to X$ in $B$ a pullback functor $b^*\colon \D(X) \to \D(Y)$.
    \item For every $f\colon Y \to X$ in $F$ an exceptional pushforward functor $f_!\colon \D(Y) \to \D(X)$.
    \item For every object $X \in \cat C$ such that $\Delta_X\colon X \to X \times X$ lies in $B$, a non-unital symmetric monoidal structure $\tensor\colon \D(X) \times \D(X) \to \D(X)$ on $\D(X)$. If the map $X \to *$ lies in $B$ then this structure is unital (i.e. a symmetric monoidal structure).
\end{itemize}
Moreover, these data satisfy the following compatibilities:
\begin{itemize}
    \item For every cartesian square
    \begin{equation*} \begin{tikzcd}
        Y' \arrow[r,"b'"] \arrow[d,"f'"] & Y \arrow[d,"f"] \\
        X' \arrow[r,"b"] & X
    \end{tikzcd}\end{equation*}
    in $\cat C$ with $b, b' \in B$ and $f, f' \in F$ there is an isomorphism $b^* f_! = f'_! b'^*$.
    
    \item For every map $f\colon Y \to X$ such that $f$ lies in $E$ and $f$, $\Delta_X$ and $\Delta_Y$ lie in $B$, there is an isomorphism $f_! (f^* M \tensor N) = f_! f^* M \tensor N$ for all $M \in \D(X)$ and $N \in \D(Y)$.
\end{itemize}
These data come together with higher homotopies exhibiting the naturality of $b^*$, $f_!$, $\tensor$, the base-change isomorphism and the projection formula.

\subsection{Sheaves on ind-stacks} \label{sec:Corr-bounded}

Given a 3-functor formalism $\D\colon \Corr(\cat C, E) \to \Cat$ on a geometric setup $(\cat C, E)$, we will extend $\D$ to a 3-functor formalism of \enquote{bounded sheaves} on a certain full subcategory of $\Ind(\cat C)$. We start with the definition of the relevant (extended) geometric setup.

\begin{defn} \label{def:geom-setup-on-Ind-stacks}
Let $(\cat C, E)$ be a geometric setup such that $\cat C$ admits finite limits. We denote by
\begin{align*}
    \Ind_E(\cat C) \subseteq \Ind(\cat C)
\end{align*}
the full subcategory spanned by the image of $\Ind(\cat C_E) \to \Ind(\cat C)$. Note that there is a natural fully faithful embedding $\cat C \subseteq \Ind_E(\cat C)$ and the objects of $\Ind_E(\cat C)$ can be written as filtered diagrams of objects in $\cat C$ along morphisms in $E$.

We denote by $E'$ the class of morphisms in $\Ind_E(\cat C)$ that lie in the image of $\Ind(\cat C_E) \to \Ind(\cat C)$. We denote by $B$ the class of representable morphisms in $\Ind_E(\cat C)$, i.e. those maps $b\colon Y \to X$ such that for every $X' \in \cat C$ with a map $X' \to X$ we have $Y \times_X X' \in \cat C$.
\end{defn}

\begin{lem}
Given a geometric setup $(\cat C, E)$, the triple $(\Ind_E(\cat C), B, E')$ from \cref{def:geom-setup-on-Ind-stacks} defines an extended geometric setup.
\end{lem}
\begin{proof}
Recall that every object $X$ in $\Ind(\cat C)$ admits a description as a filtered diagram $X = (X_i)_i$ in $\cat C$. Given a diagram $(Y_j)_j \to (X_i)_i \from (Z_k)_k$, the fibrer product is computed via a diagram of the form $(Y_{j(\ell)} \times_{X_{i(\ell)}} Z_{k(\ell)})_\ell$. From this description and the properties of $E$ one sees directly that $\Ind_E(\cat C)$ is stable under products and fiber products in $\Ind(\cat C)$ and that $E'$ is stable under composition (this is obvious), fiber products and diagonals. It is clear from the definition that $B$ is stable under composition and fiber products.
\end{proof}

We now show that a 3-functor formalism on $(\cat C, E)$ always extends uniquely and explicitly to a 3-functor formalism on $(\Ind_E(\cat C), B, E')$:

\begin{prop} \label{rslt:extend-3ff-to-Ind-stacks}
Let $\D\colon \Corr(\cat C, E)^\tensor \to \Cat^\times$ be a 3-functor formalism on a geometric setup $(\cat C, E)$ such that $\cat C$ admits finite limits. Then there is a unique extension of $\D$ to a 3-functor formalism
\begin{align*}
    \D\colon \Corr(\Ind_E(\cat C), B, E')^\tensor \to \Cat^\times
\end{align*}
such that for every $X = (X_i)_i \in \Ind_E(\cat C)$ we have
\begin{align*}
    \D(X) = \varinjlim_i \D(X_i).
\end{align*}
\end{prop}
\begin{proof}
We argue similarly as in \cite[Proposition~A.5.16]{mann-p-adic-6-functors} and \cite[Lemma~A.5.11]{mann-p-adic-6-functors}. Let us denote $\cat X^\tensor = \Corr(\cat C, E)^\tensor$ and $\cat X'^\tensor = \Corr(\Ind_E(\cat C), B, E')^\tensor$. Since $\Cat^\times$ is cartesian, the operad map $\D\colon \cat X^\tensor \to \Cat^\times$ is equivalently given by a functor $\D\colon \cat X^\tensor \to \Cat$ such that $\D(X_\bullet) = \prod_i \D(X_i)$ (see \cite[Proposition~2.4.1.7]{lurie-higher-algebra}). We perform a left Kan extension of this functor along $\cat X^\tensor \to \cat X'^\tensor$; let us denote the resulting functor by $\D'\colon \cat X'^\tensor \to \Cat$. We now compute $\D'$, so fix an object $X_\bullet \in \cat X^\tensor$ and denote $\cat K := \cat X^\tensor \times_{\cat X'^\tensor} \cat X'^\tensor_{/X_\bullet}$. By the pointwise formula for left Kan extensions we know that
\begin{align*}
    \D'(X_\bullet) = \varinjlim_{Y_\bullet \in \cat K} \D(Y_\bullet).
\end{align*}
The objects in $\cat K$ are correspondences $Y_\bullet \from Y'_\bullet \to X_\bullet$ such that all $X_i$ and $Y'_i$ are in $\Ind_E(\cat C)$, all $Y_i$ are in $\cat C$, the map $Y'_\bullet \to X'_\bullet$ lies in $E'^-$ and the map $Y'_\bullet \to Y_\bullet$ lies in $B^+$. By definition of $B^+$ and $B$, all $Y'_i$ lie in $\cat C$. Let $\cat K' \subseteq \cat K$ be the full subcategory spanned by those morphisms where $Y_\bullet \from Y'_\bullet$ lies over the identity in $\Fin_*$ and is given by a collection of isomorphisms $Y_i = Y'_i$. As in the proof of \cite[Lemma~A.5.11]{mann-p-adic-6-functors} the inclusion $\cat K' \injto \cat K$ is cofinal, hence
\begin{align*}
    \D'(X_\bullet) = \varinjlim_{Y_\bullet \in \cat K'} \D(Y_\bullet).
\end{align*}
Now write $X_\bullet = (X_1, \dots, X_n)$ and for $i = 1, \dots, n$ let $\cat K'_i$ be the category of maps in $E'$ of the form $Y_i \to X_i$ for some $Y_i \in \cat C$; more precisely $\cat K'_i = \cat C_E \times_{\Ind(\cat C_E)} \Ind(\cat C_E)_{/Y_i}$. Then $\cat K' = \prod_{i=1}^n \cat K'_i$. By definition of $\Ind$-categories, all $\cat K'_i$ and hence also $\cat K'$ are filtered. Since filtered colimits commute with finite limits in $\Cat$ we deduce
\begin{align*}
    \D'(X_\bullet) = \varinjlim_{(Y_i)_i \in \prod_i \cat K'_i} \prod_i \D(Y_i) = \prod_i \varinjlim_{Y_i \in \cat K'_i} \D(Y_i).
\end{align*}
This shows that $\D'$ is still a lax cartesian structure and thus corresponds to an operad map $\cat X'^\tensor \to \Cat^\times$. Moreover, we obtain the claimed description of the extension.
\end{proof}

Let us make the 3-functor formalism constructed in \cref{rslt:extend-3ff-to-Ind-stacks} more explicit. We can describe it as follows:
\begin{itemize}
    \item For an ind-stack $(X_i)_i \in \Ind_E(\cat C)$ we have $\D((X_i)_i) = \varinjlim_i \D(X_i)$, i.e. an object in this category is the same as an object in some $\D(X_i)$ and for any map $i \to i'$ with induced map $\gamma_{ii'}\colon X_i \to X_{i'}$ we identify $M \in \D(X_i)$ with $\gamma_{ii'!} M \in \D(X_{i'})$.
    
    \item Given a map $f\colon (Y_j)_j \to (X_i)_i$ in $B$ and som $M \in \D((X_i)_i)$, the object $f^* M$ is defined as follows. By definition of $B$ we can write $(Y_j)_j = (Y_i)_i$ for $Y_i = X_i \times_{(X_i)_i} (Y_j)_j$ and thus $f$ is induced by a family of maps $f_i\colon Y_i \to X_i$. Now pick some $i$ such that $M \in \D(X_i)$ and let $f^* M := f_i^* M \in \D(Y_i)$. The definition of $B$ guarantees that this construction is well-defined, i.e. independent of the chosen $i$.

    \item Given a map $f\colon (Y_j)_j \to (X_i)_i$ in $E'$ and some $N \in \D((Y_j)_j)$, the object $f_! N$ is defined as follows. Pick some $j$ such that $N \in \D(Y_j)$. Now $f$ induces a map $f_{ij}\colon Y_j \to X_i$ for some $i$ and by definition of $E'$ we can assume that $f_{ij} \in E$. We then let $f_! N := f_{ij!} N$. One checks easily that this definition is independent of the chosen $i$ and $j$.
\end{itemize}

In practice the category $\Ind_E(\cat C)$ is too big to be useful and it is hard to find many maps in $B$, so we will usually restrict to the full subcategory spanned by $\Ind$-systems along \emph{monomorphisms}:

\begin{defn} \label{def:Ind-stacks-with-monomorphisms}
Let $(\cat C, E)$ be a geometric setup such that $\cat C$ admits finite limits. We denote by
\begin{align*}
    \Ind_E'(\cat C) \subseteq \Ind_E(\cat C)
\end{align*}
the full subcategory spanned by the image of $\Ind(\cat C_E') \to \Ind(\cat C)$, where $\cat C_E'$ denotes the category with objects in $\cat C$ and morphisms the monomorphisms in $E$. Note that $\Ind_E'(\cat C)$ is stable under fiber products in $\Ind_E(\cat C)$, hence the classes $B$ and $E'$ from above restrict to an extended geometric setup $(\Ind_E'(\cat C), B, E')$.
\end{defn}

\begin{rmk} \label{rmk:survey-of-Ind-E-mono}
In the setting of \cref{def:Ind-stacks-with-monomorphisms}, assume that $\cat C$ admits colimits along the filtered diagrams $(X_i)_i \in \Ind_E'(\cat C)$ and that filtered colimits commute with fiber products in $\cat C$ (e.g.\ $\cat C$ is a topos). Pick some $(X_i)_i \in \Ind_E'(\cat C)$ with colimit $X$. Then the induced map $\iota\colon (X_i)_i \injto X$ is a monomorphism in $\Ind_E'(\cat C)$ and lies in $E'$. Given a 3-functor formalism $\D$ as in \cref{rslt:extend-3ff-to-Ind-stacks}, the functor
\begin{align*}
    \iota_!\colon \D((X_i)_i) = \varinjlim_i \D(X_i) \injto \D(X)
\end{align*}
is fully faithful. We often call $\D((X_i)_i)$ the full subcategory of \emph{bounded} sheaves on $X$ (with respect to the $\Ind$-structure $(X_i)_i$). Pick some $(Y_j)_j \in \Ind_E'(\cat C)$ with colimit $Y$. Then a map $(Y_j)_j \to (X_i)_i$ is the same as a map $Y \to X$ such that for all $j$ the map $Y_j \to X$ factors over some $X_i \subseteq X$. This map lies in $B$ if and only if we have an equivalence of $\Ind$-systems $(Y_j)_j = (Y \times_X X_i)_i$ (as subobjects of $Y$).
\end{rmk}

As seen in \cref{rmk:survey-of-Ind-E-mono} the category $\Ind_E'(\cat C)$ together with \cref{rslt:extend-3ff-to-Ind-stacks} provide a general way of defining bounded sheaves inside a 6-functor formalism, by explicitly specifying the relevant subobject $(X_i)_i \subset X$ for each $X$. In practice one often starts with a diagram in $\Corr(\cat C)$ and wants to upgrade it to a diagram in $\Corr(\Ind_E'(\cat C))$ by choosing the relevant $\Ind$-subobjects. One then runs into the problem of showing that these choices can be made coherent enough to induce the desired functor $\Corr(\Ind_E'(\cat C))$. In the following we show that the necessary coherence is only a \emph{condition}, no additional data. We first need some preparation:

\begin{lem} \label{rslt:factor-functor-using-subobjects}
Let $F\colon \cat I \to \cat C$ and $\alpha\colon \cat C_0 \to \cat C$ be functors of categories such that $\alpha$ induces monomorphisms on $\Hom$-anima. For each $i \in \cat I$ choose some preimage $X_i \in \cat C_0$ of $F(i)$ under $\alpha$ such that for every map $f\colon i \to j$ in $\cat I$ we have
\begin{align*}
    F(f) \in \operatorname{img}\big( \Hom(X_i, X_j) \injto \Hom(F(i), F(j))\big)
\end{align*}
Then $F$ factors uniquely over a functor $F_0\colon \cat I \to \cat C_0$ which satisfies $F_0(i) = X_i$ for all $i \in \cat I$.
\end{lem}
\begin{proof}
Via the Yoneda embedding $\cat C \injto \cat P(\cat C)$, $F$ corresponds to a functor $\cat I \times \cat C^{\op} \to \Ani$. Precomposition of this functor with $\id_{\cat I} \times \alpha$ produces the functor
\begin{align*}
    F'\colon \cat I \times \cat C_0^{\op} \to \Ani, \qquad (i, X) \mapsto \Hom_{\cat C}(\alpha(X), F(i)).
\end{align*}
Now for each $(i, X) \in \cat I \times \cat C_0^{\op}$ we choose the subobject $\Hom_{\cat C_0}(X, X_i) \subseteq F'(i, X)$. By the conditions on $F$ these subobjects are mapped to each other under morphisms in $\cat I \times \cat C_0^{\op}$. Via unstraightening $F'$, passing to the induced full subcategory, then straightening again, we obtain the functor
\begin{align*}
    F''\colon \cat I \times \cat C_0^{\op} \to \Ani, \qquad (i, X) \mapsto \Hom_{\cat C_0}(X, X_i).
\end{align*}
Using Yoneda for $\cat C_0$ we arrive at the desired functor $F_0\colon \cat I \to \cat C_0$. It is clear from the construction that $F_0$ is unique.
\end{proof}


\begin{prop} \label{rslt:lift-functor-from-Corr-to-Corr-Ind-stacks}
Let $(\cat C, E)$ be a geometric setup such that $\cat C$ admits finite limits and small colimits and they commute. Let $\cat K$ be another category and $F\colon \cat K \to \Corr(\cat C)$ be a functor. For each $k \in \cat K$ we choose an object $(F(k)_i)_{i\in I_k} \in \Ind_E'(\cat C)$ such that $\varinjlim_i F(k)_i = F(k)$ and such that the following condition is satisfied:
\begin{itemize}
    \item[($*$)] For every $\gamma\colon k \to k'$ let $F(k) \from F(\gamma) \to F(k')$ be the correspondence in $\cat C$ that depicts the morphism $F(\gamma)$ in $\Corr(\cat C)$. For each $i \in I_k$ denote $F(\gamma)_i := F(\gamma) \times_{F(k)} F(k)_i$. We require that the map $F(\gamma)_i \to F(k')$ factors over $F(k')_{i'}$ for some $i' \in I_{k'}$ and that the map $F(\gamma)_i \to F(k')_{i'}$ lies in $E$.
\end{itemize}
Then $F$ lifts uniquely to a functor $\cat K \to \Corr(\Ind_{E'}(\cat C), B, E')$ that sends $k \in \cat K$ to $(F(k)_i)_i$.
\end{prop}
\begin{proof}
We apply \cref{rslt:factor-functor-using-subobjects} to $\alpha\colon \Corr(\Ind_E'(\cat C), B, E') \to \Corr(\cat C)$. It only remains to verify that $\alpha$ induces monomorphisms on $\Hom$-anima. Thus we fix two objects $(X_i)_i, (Y_j)_j \in \Ind_E'(\cat C)$ with colimits $X$ and $Y$ and we need to verify that the functor
\begin{align*}
    \Ind_E'(\cat C)_{/(X_i)_i \times (Y_j)_j} \to \cat C_{/X \times Y}
\end{align*}
is fully faithful on the full subcategory of those $(Z_a)_a \to (X_i)_i \times (Y_j)_j$ such that the map $(Z_a)_a \to (X_i)_i$ lies in $B$ and the map $(Z_a)_a \to (Y_j)_j$ lies in $E'$. Pick two such objects $(Z_a)_a$ and $(Z'_{a'})_{a'}$ and let $Z$ and $Z'$ be their colimits. By the description of $B$ in \cref{rmk:survey-of-Ind-E-mono} we have $(Z_a)_a = (Z_i)_i$ and $(Z'_{a'})_{a'} = (Z'_i)_i$ for $Z_i = Z \times_X X_i$ and $Z'_i = Z \times_X X_i$. We are thus reduced to showing that the map
\begin{align}
    \Hom_{(X_i)_i \times (Y_j)_j}((Z_i)_i, (Z'_i)_i) \isoto \Hom_{X \times Y}(Z, Z') \label{eq:lift-to-Corr-of-Ind-stacks-Hom-isom}
\end{align}
is an isomorphism. Note that on both sides of \cref{eq:lift-to-Corr-of-Ind-stacks-Hom-isom} we can pull out the colimit $Z = \varinjlim_i Z_i$ to reduce to the case that $Z = Z_i$. Since $(Z'_i)_i \injto Z'$ is a monomorphism, the map in \cref{eq:lift-to-Corr-of-Ind-stacks-Hom-isom} is a monomorphism, so it only remains to show that it is surjective. This amounts to saying that every map $Z_i \to Z'$ over $X$ factors over $Z'_{i'}$ for some $i' \ge i$. In fact we can choose $i' = i$ because $Z'_i = X_i \times_X Z'$.
\end{proof}

\begin{rmk} \label{rmk:lift-monoid-in-Corr-to-Corr-Ind}
Since the operad $\Corr(\Ind_E'(\cat C), B, E')^\tensor$ is also defined using a correspondence category (see \cref{def:Corr-operad-for-extended-geometric-setup}) one can apply the same argument as in \cref{rslt:lift-functor-from-Corr-to-Corr-Ind-stacks} to the operads as well. For example let $X \in \Alg(\Corr(\cat C))$ with multiplication and unit map given by correspondences
\begin{equation*}
    \begin{tikzcd}
        & Y \arrow[dl] \arrow[dr]\\
        * \arrow[rr,dashed] && X
    \end{tikzcd}
    \qquad\qquad
    \begin{tikzcd}
        & M \arrow[dl] \arrow[dr]\\
        X \times X \arrow[rr,dashed] && X
    \end{tikzcd}
\end{equation*}
Choose a presentation $(X_i)_i \in \Ind_E'(\cat C)$ with $\varinjlim_i X_i = X$. Then the algebra structure on $X$ upgrades (uniquely) to an algebra structure on $(X_i)_i$ in $\Corr(\Ind_E'(\cat C), B, E')$ if and only if the following two conditions are satisfied:
\begin{enumerate}[(a)]
    \item The map $Y \to X$ factors over some $X_i$ such that $Y \to X_i$ lies in $E$.
    \item For every $i$ let $M_i := M \times_{X \times X} (X_i \times X_i)$. Then the map $M_i \to X$ factors over some $X_{i'}$ such that $M_i \to X_{i'}$ lies in $E$.
\end{enumerate}
\end{rmk}

\subsection{Convolution stacks} \label{sec:convolution-stacks}

A crucial step in the construction of the spectral action in \cite{FS} is the convolution monoidal structure on (bounded) sheaves on the local Hecke stack. Some care is required to rigorously define this monoidal structure, as it involves several higher homotopies. In the following we recall the basic definition of convolution stacks and then prove a basic functoriality result, which is needed in the construction of the spectral action. Similar results can probably be obtained from \cite{Goedicke-2Segal}.

\begin{defn} \label{def:convolution-stack}
Let $\cat C$ be a category which has all fiber products and let $X \to S$ be a map in $\cat C$. Recall that the correspondence category $\Corr(\cat C_{/S})$ is enriched over itself, with $\intHom(X, Y) = Y \times_S X$ (see \cite[\S2.4]{heyer-mann-6ff}). We define
\begin{align*}
    \Conv_S(X) := \intHom_S(X, X) \in \Alg(\Corr(\cat C_{/S})).
\end{align*}
Since the forgetful functor $\Corr(\cat C_{/S}) \to \Corr(\cat C)$ is lax symmetric monoidal (cf. \cite[Definition~4.1.3(b)]{heyer-mann-6ff}) we can also view $\Conv_S(X)$ as an algebra in $\Corr(\cat C)$.
\end{defn}

In the setting of \cref{def:convolution-stack}, assume that $\cat C$ also has finite products. Then we can describe $\Conv_S(X)$ as follows: The underlying object is $X \times_S X$, the unit map is given by the correspondence
\begin{equation*}\begin{tikzcd}
    & X \arrow[dl] \arrow[dr,"\Delta_X"]\\
    * \arrow[rr,dashed] && X \times_S X
\end{tikzcd}\end{equation*}
where $\Delta_X$ is the diagonal of $X$ over $S$, and the multiplication is given by
\begin{equation*}\begin{tikzcd}
    & X \times_S X \times_S X \arrow[dl,"{(\id_X, \Delta_X, \id_X)}",swap] \arrow[dr,"\pr_{13}"] \\
    (X \times_S X) \times (X \times_S X) \arrow[rr,dashed] && X \times_S X
\end{tikzcd}\end{equation*}
where $\pr_{13}$ denotes the projection to the first and third entry.

\begin{exmpl} \label{exmpl:convolution-of-classfying-stacks}
Suppose $\cat C$ is a topos and $H \to G$ is a morphism of groups in $\cat C$ (see e.g. \cite[\S B.3]{heyer-mann-6ff}). Denote by $BG = */G$ and $BH = */H$ their classifying stacks, so that we get an induced map $BH \to BG$ of objects in $\cat C$. Then
\begin{align*}
    \Conv_{BG}(BH) = H \backslash G / H
\end{align*}
and the multiplication is given by the correspondence
\begin{equation*}\begin{tikzcd}
    & H \backslash G \times^H  G / H \arrow[dl] \arrow[dr] \\
    H \backslash G / H \times H \backslash G / H && H \backslash G / H
\end{tikzcd}\end{equation*}
where \enquote{$\times^H$} denotes the quotient of the product by the canonical $H$ action (acting by multiplication on the right on the left factor and vice-versa) and the forward map is induced by multiplication in $G$.
\end{exmpl}

We now discuss the functoriality of the assignment $[X \to S] \mapsto \Conv_S(X)$. Since $\Conv_S(X)$ is defined using the enriched endomorphisms in $\Corr(\cat C_{/S})$, the key construction is the functoriality of enriched endomorphisms. In the following, for a monoidal category $\cat V$ (or more generally an operad over the associative operad $\Ass^\tensor$) we denote by $\Enr_{\cat V}$ the category of $\cat V$-enriched categories: Roughly, a $\cat V$-enriched category is a category where all $\Hom$'s are objects in $\cat V$. We refer the reader to \cite[\S C]{heyer-mann-6ff} for an introduction to enriched category theory.

\begin{lem} \label{rslt:functoriality-of-enriched-End}
Let $\cat V^\tensor$ be an operad over $\Ass^\tensor$ (e.g. a monoidal category), $\cat I$ a category and $F\colon \cat I \to \Enr_{\cat V}$ a functor. Denote by $\cat E \to \cat I$ the cocartesian unstraightening of the composition $\cat I \to \Enr_{\cat V} \to \Ani$, where the second functor sends an enriched category to the underlying anima of objects (an object in $\cat E$ is a pair $(i, X)$ consisting of some $i \in \cat I$ and some $X \in F(i)$). Then there is a functor
\begin{align*}
    \End^{\cat V}\colon \cat E \to \Alg(\cat V)
\end{align*}
that sends a pair $(i, X) \in \cat E$ to $\End^{\cat V}_{F(i)}(X) = \Hom^{\cat V}_{F(i)}(X, X) \in \cat V$ equipped with the composition monoidal structure.
\end{lem}
\begin{proof}
We use the Gepner--Haugseng model for enriched categories (cf. \cite[Remark~C.1.4(i)]{heyer-mann-6ff}). Let us briefly recall this model. To an anima $S$ we associate a generalized operad $\Ass_S^\tensor$ over $\Ass^\tensor$ with underlying category $S \times S$. Then a $\cat V$-enriched precategory with object anima $S$ is a map of generalized operads $\Ass_S^\tensor \to \cat V^\tensor$. Denote the category of $\cat V$-enriched precategories on $S$ by $\operatorname{PreEnr}_{\cat V, S}$. One can fit this construction into a cartesian fibration $\operatorname{PreEnr}_\cat V \to \Ani$, see \cite[Definition~4.3.1]{GepnerHaugsengEnriched}.

The crucial observation is that there is another generalized operad associated to an anima $S$, namely $\Ass^\tensor \times S$. We denote the category of operad maps $\Ass^\tensor \times S \to \cat V^\tensor$ by $\Alg^S(\cat V)$ and note that it is equivalent to $\Fun(S, \Alg(\cat V))$. Pullback along the natural map $\Ass^\tensor \times S \to \Ass_S^\tensor$ induces a functor $\operatorname{PreEnr}_{\cat V, S} \to \Alg^S(\cat V)$ and by the construction of Gepner--Haugseng we obtain a map of cartesian fibrations
\begin{align*}
    \operatorname{PreEnr}_{\cat V} \to \Alg^{(-)}(\cat V)
\end{align*}
over $\Ani$. This functor sends a $\cat V$-enriched precategory $\cat C$ to the functor $\cat C^\simeq \to \Alg(\cat V)$ mapping $X \mapsto \End^{\cat V}_{\cat C}(X)$.

By construction, $\Alg^{(-)}(\cat V) \to \Ani$ is the cartesian unstraightening of the functor $\Ani^\op \to \Cat$, $S \mapsto \Fun(S, \Alg(\cat V))$. By \cite[Proposition~7.3]{FreeFibrations}, if $\cat U \to \Ani$ denotes the cocartesian unstraightening of the identity $\id\colon \Ani \to \Ani$ (strictly speaking we have to replace the second $\Ani$ by some large version of it, but we ignore set-theoretic issues here), we obtain the evaluation map
\begin{align*}
    \ev\colon \cat U \times_{\Ani} \Alg^{(-)}(\cat V) \to \Alg(\cat V).
\end{align*}
that takes a pair $(S, X) \in \cat U$ consisting of an anima $S$ and an object $X \in S$ together with a functor $f\colon S \to \Alg(\cat V)$ and sends it to $f(X)$. We are now ready to put everything together. Since $\cat E \to \cat I$ is the cocartesian unstraightening of a functor $\cat I \to \Ani$, we can write it as the pullback of $\cat U \to \Ani$ along $\cat I \to \Ani$. We thus obtain the following pullback square:
\begin{equation*}\begin{tikzcd}
    \cat E \arrow[d] \arrow[rrr] &&& \cat U \arrow[d]\\ 
    \cat I \arrow[r,"F"] & \Enr_{\cat V} \subseteq \operatorname{PreEnr}_{\cat V} \arrow[r] & \Alg^{(-)}(\cat V) \arrow[r] & \Ani.
\end{tikzcd}\end{equation*}
In particular we have a map $\cat E \to \cat U \times_{\Ani} \Alg^{(-)}(\cat V)$. Composing this map with $\ev$ produces the desired functor $\cat E \to \Alg(V)$.
\end{proof}

With the functoriality of enriched endomorphisms at hand, we can now discuss the functoriality of convolution stacks. In the following, for any category $\cat C$ we denote by $\Ar(\cat C) := \Fun([1], \cat C)$ the category of morphisms in $\cat C$. A morphism in $\Ar(\cat C)$ is a commuting square in $\cat C$.

\begin{prop} \label{rslt:functoriality-of-convolution-stacks}
Let $\cat C$ be a category with all fiber products. Let $\Corr'(\Ar(\cat C)) \subseteq \Corr(\Ar(\cat C))$ be the subcategory with the same objects, but where we only allow those morphisms
\begin{equation*}\begin{tikzcd}[column sep=tiny]
    & {[X_2 \to S_2]} \arrow[dl] \arrow[dr] \\
    {[X_1 \to S_1]} && {[X_3 \to S_3]}
\end{tikzcd}\end{equation*}
such that the induced maps $X_2 \isoto X_3$ and $X_2 \isoto X_1 \times_{S_1} S_2$ are isomorphisms. Then there is a canonical functor
\begin{align*}
    \Conv\colon \Corr'(\Ar(\cat C)) \to \Alg(\Corr(\cat C))
\end{align*}
that sends a map $X \to S$ to $\Conv_S(X)$.
\end{prop}
\begin{proof}
We may assume that $\cat C$ has all finite limits by embedding it in a larger category if necessary (and restricting all functors afterwards). In \cite[Lemma~4.2.2(i)]{heyer-mann-6ff} it was shown that the assignment $S \mapsto \Corr(\cat C_{/S})$ is functorial in $S$ and induces a functor $F\colon \Corr(\cat C) \to \Enr_{\Corr(\cat C)}$. We apply \cref{rslt:functoriality-of-enriched-End} to $F$ in order to obtain a functor
\begin{align*}
    \cat E \to \Alg(\Corr(\cat C)),
\end{align*}
where $\cat E \to \Corr(\cat C)$ is the cocartesian unstraightening of the composition $\Corr(\cat C) \to \Enr_{\Corr(\cat C)} \to \Ani$. Note that an object of $\cat E$ is given by a pair $(S, X)$, where $S$ is an object of $\Corr(\cat C)$ and $X$ is an object in $\Corr(\cat C_{/S})^\simeq = (\cat C_{/S})^\simeq$ -- in other words, an object in $\cat E$ is given by a map $X \to S$ in $\cat C$. Moreover, by construction it is sent to $\Conv_S(X) \in \Alg(\Corr(\cat C))$ under the above functor.

It remains to show that there is an isomorphism $\cat E \isom \Corr'(\Ar(\cat C))$. Note that since for $S \in \cat C$ we have $\Hom_{\Corr(\cat C)}(*, S) = (\cat C_{/S})^\simeq$, the functor $\Corr(\cat C) \to \Ani$ considered above is corepresentable by $*$. Similarly, one checks that the target map $\Corr'(\Ar(\cat C)) \to \Corr(\cat C)$ is a cocartesian fibration (as in \cite[Proposition~2.2.14]{heyer-mann-6ff}) and that the corresponding functor $\Corr(\cat C) \to \Ani$ is corepresented by $*$. This shows $\cat E \isom \Corr'(\Ar(\cat C))$, as desired.
\end{proof}

\begin{rmk}
A similar result should be obtainable from \cite[Theorem~1.1]{Goedicke-2Segal} by constructing a functor from $\Corr'(\Ar(\cat C))$ into the category of 2-Segal spaces considered in loc. cit.
\end{rmk}

Let us make the construction in \cref{rslt:functoriality-of-convolution-stacks} more explicit. Suppose we are given a commuting square
\begin{equation*}\begin{tikzcd}
    X' \arrow[d] \arrow[r] & X \arrow[d] \\
    S' \arrow[r] \arrow[r] & S
\end{tikzcd}\end{equation*}
in $\cat C$, i.e. a morphism in $\Ar(\cat C)$. Then \cref{rslt:functoriality-of-convolution-stacks} answers the question under which assumptions we obtain a monoidal map between $\Conv_S(X) = X \times_S X$ and $\Conv_{S'}(X') = X' \times_{S'} X'$. Namely, there are two cases in which such a map exists:
\begin{itemize}
    \item Suppose that the map $X' \to X$ is an isomorphism. Then we get an induced map $\Conv_{S'}(X') = \Conv_{S'}(X) \to \Conv_S(X)$ in $\Alg(\Corr(\cat C))$. The underlying map is the obvious projection $X \times_{S'} X \to X \times_S X$.
    
    \item Suppose that the square is cartesian, i.e. $X' = X \times_S S'$. Then we get an induced map $\Conv_S(X) \to \Conv_{S'}(X')$ in $\Alg(\Corr(\cat C))$. The underlying map in $\Corr(\cat C)$ is given by the backward map $X \times_S X \from X \times_S X' = X' \times_{S'} X'$.
\end{itemize}
It is straightforward to check by hand that the above maps indeed commute with the multiplication map on the convolution stacks, but it requires more work to construct all the higher homotopies (even if $\cat C$ is a category of classical stacks, $\Corr(\cat C)$ is already a $(3,1)$-category, so there are many higher homotopies to consider!). \Cref{rslt:functoriality-of-convolution-stacks} shows this and also proves that the construction is functorial in the square we started with.

An important property of the convolution stack is that it is self-dual: Swapping the two factors in $\Conv_S(X) = X \times_S X$ induces an isomorphism $\Conv_S(X)^\op \isoto \Conv_S(X)$. The following result proves this isomorphism and shows that it is natural in $X$ and $S$:

\begin{prop} \label{rslt:self-duality-of-convolution-stacks}
In the setup of \cref{rslt:functoriality-of-convolution-stacks}, there is a natural isomorphism
\begin{align*}
    \Conv_S(X)^\op \isoto \Conv_S(X)
\end{align*}
of functors $\Corr'(\Ar(\cat C)) \to \Alg(\Corr(\cat C))$ given by swapping the two factors in $\Conv_S(X) = X \times_S X$.
\end{prop}
\begin{proof}
Consider the functor $F\colon \Corr(\cat C) \to \Enr_{\Corr(\cat C)}$, $S \mapsto \Corr(\cat C_{/S})$ from \cref{rslt:functoriality-of-convolution-stacks}. We can pass to opposite enriched categories in order to obtain the functor $F'\colon \Corr(\cat C) \to \Enr_{\Corr(\cat C)}$, $S \mapsto \Corr(\cat C_{/S})^\op$. If we apply \cref{rslt:functoriality-of-enriched-End} to $F'$ then by the definition of opposite enriched categories in terms of opposite quivers (see \cite[Remark~4.43]{Heine.2023}), we obtain the functor $\Corr'(\Ar(\cat C)) \to \Alg(\Corr(\cat C))$, $S \mapsto \Conv_S(X)^\op$. But by \cite[Lemma~4.2.2(ii)]{heyer-mann-6ff} there is an equivalence $F' \isoto F$ of functors $\Corr(\cat C) \to \Enr_{\Corr(\cat C)}$ which induces the identity on underlying anima; applying \cref{rslt:functoriality-of-enriched-End} to this equivalence (i.e. the induced functor $[1] \times \Corr(\cat C) \to \Enr_{\Corr(\cat C)}$) produces the claimed natural equivalence $\Conv_S(X)^\op \isoto \Conv_S(X)$.
\end{proof}

\section{Ingredients for the Hecke action} \label{sec:abstract-Hecke-action}

Fix a finite extension $F$ of $\Q_p$ and an algebraic group $G$ over $F$. For a finite set $I$, let $\Div^I := (\Div^1)^I$ and let $\locHck_G^{I}$ be the local Hecke stack over $\Div^{I}$, in the setting of \cite{FS}. Suppose we are given a reasonable six functor formalism $X\mapsto \D(X)$ on small v-stacks and let $\D^{\bdd}(\locHck_{G}^{I},\Lambda)$ denote the full subcategory of $\D(\locHck_{G}^{I})$ spanned by sheaves with quasicompact support relative to $\Div^{I}$. Our goal in this appendix is to give a precise meaning to the following two statements:
\begin{enumerate}
    \item Convolution naturally gives $\D^{\bdd} (\locHck_{G}^{I})$ the structure of a $\D(\Div^{I})$-linear monoidal category which is functorial in $I$.
    \item There is a natural $\D^{\bdd} (\Div^{I})$-linear monoidal functor \[\D^{\bdd}(\locHck_{G}^{I}) \to \D(\Bun_G \times \Bun_G \times \Div^I)\]
    which is functorial in $I$.
\end{enumerate}
Both of these statements (or minor variants) seem to be regarded as well-known. However, there are several issues involved in making them precise, which in the end all stem from the usual source: all the categories involved are $\infty$-categories, so we cannot just ``write down formulas'' and declare that they work, as this ignores the necessity of exhibiting the higher coherences. Moreover, this issue occurs in two different ``directions'' simultaneously: for a fixed $I$ we need to exhibit a monoidal structure on the $\infty$-category $\D^{\bdd}(\locHck_{G}^{I})$, and then we need to show that this can be done coherently in varying $I$. An additional complication is that the functoriality in $I$ is itself exhibited by a correspondence, as there are no natural maps between the stacks $\locHck_{G}^{I}$ for varying $I$.

We also remark that the arguments in this appendix are quite formal, and should easily adapt e.g. to the versions of $\Bun_G$ and the local and global Hecke stacks used in classical geometric Langlands. It seems to us that the issues we address in this appendix are not dealt with properly in the literature, for any form of geometric Langlands.

\subsection{Interpolating Hecke stacks}

To resolve the functoriality of the Hecke stack in $I$, we need to actually consider some new spaces, and our first goal is to define these. Let $f:I \to J$ denote a map of finite sets. This induces a map $\Div^J \to \Div^I$. Moreover, as in the discussion preceding \cite[Definition VI.1.5]{FS}, we have a sheaf of rings $B^{+}_{\Div^J}$ over $\Div^J$, and for each $j\in J$ we have the associated (locally principal) ideal sheaf $\mathcal{I}_{j} \subset B^{+}_{\Div^J}$. For any affine scheme $Z$ over $F$, we have v-sheaves $L^+_{J}Z$ resp. $L_{J}Z$ whose $S$-points are the $B^{+}_{\Div^J}(S)$-points resp. $B^{+}_{\Div^J}[\frac{1}{\prod_{j \in J} \mathcal{I}_{j}}](S)$-points of $Z$. The following definition interpolates between these two objects.

\begin{defn}
Given a map of finite sets $f:I \to J$, let $B^{I \to J}$ be the sheaf of rings over $\Div^J$ defined as follows. For an affinoid perfectoid $S \to \Div^{J}$ with associated closed Cartier divisors $D_{j} \subset X_{S}$ and ideal sheaves $\mathcal{I}_j \subset \mathcal{O}_{X_S}$, $j \in J$, the sections of $B^{I\to J}$ over $S$ are given by $B^{+}_{\Div^J}[\frac{1}{\prod_{i \in I} \mathcal{I}_{f(i)}}](S)$.
\end{defn}

In other words, $B^{I\to J}$ is the ring of functions on the formal completion of $X_S$ along $\cup_{j \in J} D_j$ punctured along $\cup_{ i \in I}D_{f(i)}$.

\begin{defn} \label{def:interpolating-loop-group}
If $I\to J$ is a map of finite sets and $Z$ is an affine scheme over $F$, we define a v-sheaf over $\Div^J$ by $L_{I \to J}Z (S) = Z(B^{I \to J}(S))$.
\end{defn}

We will only apply this in the case where $Z=G$ is our fixed reductive group, in which case $L_{I \to J} G$ is a sheaf of groups. It is easy to see that $L_{I \to J} G$ is a closed subfunctor of $L_{J}G$ interpolating between two extreme cases: when $I = \emptyset$ it is $L^{+}_{J} G$, and when $I \to J$ is surjective it is $L_{J} G$. Note also that $L^{+}_{J} G$ is a subgroup of $L_{I \to J}G$.

\begin{defn}
Given a map of finite sets $I \to J$, we define the v-stack
\[\locHck_{G}^{I \to J} = L^{+}_{J} G \backslash L_{I \to J} G / L^{+}_{J} G.\]
\end{defn}

Informally, $\locHck_{G}^{I \to J}$ parametrizes a pair of $G$-bundles $\mathcal{E}_1,\mathcal{E}_2$ on the formal completion of $X_{S}$ along $\cup_{j \in J} D_j$, together with an isomorphism on the complement of $\cup_{ i \in I} D_{f(i)}$.

\begin{prop} \label{rslt:properties-of-interpolating-locHck}
For any map of finite sets $I \to J$, we have a natural correspondence
\[
\xymatrix{ & \locHck_{G}^{I \to J}\ar[dl]_{\pi}\ar[dr]^{\iota} & \\
\locHck_{G}^{I} &  & \locHck_{G}^{J}
}
\]
where $\iota$ is a closed immersion and $\pi$ is a base change of the canonical map $[\Div^J / L^{+}_{J} G] \to [\Div^I / L^{+}_{I} G]$. 

Moreover, if $I \to J \to K$ are composable maps of finite sets, there is a natural isomorphism
\[\locHck_{G}^{I \to K} \cong \locHck_{G}^{I \to J} \times_{\locHck_{G}^{J}}\locHck_{G}^{J \to K}.\]
\end{prop}
\begin{proof}
The maps are clearly induced by the definitions. The fact that $\iota$ is a closed immersion is immediate. For the claim about $\pi$, note that by Beauville-Laszlo gluing we have a cartesian diagram
\[
\xymatrix{[\Div^J / L^{+}_{J} G] \ar[r]\ar[d] & [\Div^I / L^{+}_{I} G] \ar[d] \\
[\Div^J / L_{I \to J} G]\ar[r] &  [\Div^I / L_{I} G]
}.
\]
Then we compute that
\begin{align*}
    \locHck_{G}^{I \to J} & = [\Div^J / L^{+}_{J} G] \times_{[\Div^J / L_{I \to J} G]} [\Div^J / L^{+}_{J} G] \\
    & =  [\Div^J / L^{+}_{J} G] \times_{[\Div^J / L_{I \to J} G]} ([\Div^J / L_{I \to J} G] \times_{[\Div^I / L_{I} G]} [\Div^I / L^{+}_{I} G] ) \\
    & = [\Div^J / L^{+}_{J} G] \times_{[\Div^I / L_{I} G]} [\Div^I / L^{+}_{I} G] 
\end{align*}
where the second line follows from the previous cartesian diagram. It is now easy to see that via this isomorphism, the canonical map from $\locHck_{G}^{I \to J}$ towards $\locHck_{G}^{I}$ is induced by the evident base change of the map $[\Div^J / L^{+}_{J} G] \to [\Div^I / L^{+}_{I} G]$ as claimed.

The final claim is straightforward.
\end{proof}

\begin{rmk} \label{rslt:explicit-description-of-transition-functors-on-locHck}
The analysis in the previous proof shows that for any $I\to J$, we have a commutative diagram
\[
    \xymatrix{ & \Gr_G^I \times_{\Div^I} \Div^J \ar[dl]_{\pi'}\ar[dr]^{\iota'}\ar[d]^{t_{I\to J}}\\
    \Gr_G^I \ar[d]^{t_{I}} & \locHck_G^{I\to J}\ar[dl]_{\pi}\ar[dr]^{\iota} & \Gr_G^J\ar[d]^{t_{J}}\\
    \locHck_G^I &  & \locHck_G^J
    }
\]
where $t_{I}$ is an $L^{+}_{I}G$-torsor, $t_{I \to J}$ and $t_J$ are $L^{+}_J G$-torsors, and the right-hand trapezoid is cartesian (so $\iota'$ is a closed immersion). We emphasize that the left-hand trapezoid is not cartesian. If $A$ is any sheaf on $\locHck_{G}^{I}$, an easy calculation shows that $t_{J}^\ast \iota_! \pi^{\ast} A = \iota'_{!} \pi'^{\ast} t_{I}^{\ast}A$. In particular, since the Satake category embeds fully faithfully into sheaves on $\mathrm{Gr}_{G}^{I}$ via pullback along $t_I$ (and similarly for $J$), Fargues-Scholze are able to avoid consideration of $\locHck_{G}^{I \to J}$ and the correspondence it induces by only considering the functor $\iota'_{!} \pi'^{\ast}$ and observing that it preserves (the image of) the Satake category. However, aside from breaking the symmetry of the situation, this approach is not compatible with the convolution structures, and does not give ``functoriality in $I$'' for the full category of bounded sheaves with its convolution monoidal structure.
\end{rmk}

\subsection{Functoriality of the local Hecke stack}

We can now state a slightly more precise version of our goal: For any reasonable six-functor formalism $\D(-)$ on small v-stacks, there is an attached notion of \enquote{bounded sheaves} $\D(-)^{\bdd} \subseteq \D(-)$ and a natural convolution monoidal structure on $\D(\locHck_{G}^{I})^{\bdd}$. Moreover, this structure is functorial in $I$, in the sense that it upgrades to a functor
\begin{align*}
    \Fin \to \Mon, \qquad I \mapsto \D(\locHck_G^I)^{\bdd}
\end{align*}
from the category of finite sets to the category of monoidal ($\infty$-)categories. This functor sends a morphism $I \to J$ to the functor
\begin{align*}
    \iota_!\pi^*: \D(\locHck_{G}^{I})^{\bdd} \to \D(\locHck_{G}^{J})^{\bdd}.
\end{align*}
constructed in \cref{rslt:properties-of-interpolating-locHck}. In particular we claim that this functor is monoidal and natural in $I \to J$. The passage to bounded sheaves is crucial here, as otherwise the necessary $!$-functors do not exist: the local Hecke stack is too big for that. This complicates the situation, as we cannot rely on the functoriality of $\D$ and instead need to construct a 3-functor formalism on $\D(-)^{\bdd}$ directly. We have developed general tools to overcome the challenges of the above goal in \cref{sec:Corr-bounded} and show in the following how they can be applied to our situation.

Let us start with the construction of the geometry. In \cref{sec:convolution-stacks} we discuss a general and abstract definition of convolution monoidal structures on stacks of the form $\Conv_S(X) := X \times_S X$ for a map $X \to S$ and we in particular provide a general functoriality result for this construction. Note that the Hecke stacks $\locHck_G^I$ and even the interpolating Hecke stacks $\locHck_G^{I \to J}$ are convolution stacks in this sense (see \cref{exmpl:convolution-of-classfying-stacks}). This leads to the following result:

\begin{lem} \label{rslt:functoriality-of-locHck}
Fix an algebraic group $G$ over $F$. There is a canonical functor
\begin{align*}
    \locHck_G\colon \Fin \to \Alg(\Corr(\vStk)), \qquad I \mapsto \locHck_G^I
\end{align*}
that sends a map $I \to J$ to the correspondence in \cref{rslt:properties-of-interpolating-locHck}.
\end{lem}
\begin{proof}
We construct this functor in several steps. The first step is the construction of the functor
\begin{align*}
    F_1\colon \Fin \to \Corr(\Ar(\Tw(\Fin))), \qquad I \mapsto \big[ [\emptyset \to I] \to [I \xto{\id_I} I] \big].
\end{align*}
The category $\Corr(\Ar(\Tw(\Fin)))$ may look a bit intimidating, so we give a brief description: An object in $\Tw(\Fin)$ is a map $[I \to J]$ of finite sets; morphisms $[I \to J] \to [I' \to J']$ are given by commuting squares
\begin{equation*}\begin{tikzcd}
    I \arrow[d] \arrow[r] & I' \arrow[d] \\
    J & J' \arrow[l]
\end{tikzcd}\end{equation*}
An object in $\Corr(\Ar(\Tw(\Fin)))$ is thus a map $[I \to J] \to [I' \to J']$ in $\Tw(\Fin)$, and a morphism between such maps is a correspondence of such maps. Now in order to construct $F_1$ we use the adjunction between $\Tw$ and $\Corr$ (see \cite[Theorem~2.18]{Haugseng-Hebestreit.Spans}), which tells us that $F_1$ is equivalently given by the functor
\begin{align*}
    F_1'\colon \Tw(\Fin) \to \Ar(\Tw(\Fin)), \qquad [I \to J] \mapsto \big[ [\emptyset \to J] \to [I \to J] \big].
\end{align*}
This functor can for example be constructed by hand, as both source and target are ordinary 1-categories.

The second step of our construction is the functor
\begin{align*}
    F_2\colon \Tw(\Fin) \to \vStk, \qquad [I \to J] \mapsto \Div^J \!\! / L_{I \to J} G.
\end{align*}
We have defined $L_{I \to J} G$ in \cref{def:interpolating-loop-group} and it is easy to see that it is functorial in $I \to J$ in the claimed way. In order to pass to the quotient stacks, we observe that $L_{(-)} G$ is really a functor $\Tw(\Fin) \to \cat E$, where $\cat E \to \Fin^{\op}$ is the cocartesian unstraightening of the functor $\Fin^{\op} \to \Cat$, $J \mapsto \Grp(\vStk_{/\Div^J})$ (i.e. the fiber of $\cat E$ over $J$ is the category of group stacks over $\Div^J$). One can then fiberwise pass to quotient stacks, inducing a functor $\cat E \to \vStk$.

The third step of our construction is the functor
\begin{align*}
    F_3\colon \Corr'(\Ar(\vStk)) \to \Alg(\vStk), \qquad [X \to S] \mapsto \Conv_S(X) = X \times_S X,
\end{align*}
where $\Corr' \subseteq \Corr$ is the subcategory with the same objects, but conditions on the morphisms; this functor was constructed in \cref{rslt:functoriality-of-convolution-stacks}. Note that the composition
\begin{align*}
    \Fin \xto{F_1} \Corr(\Ar(\Tw(\Fin))) \xto{\Corr(\Ar(F_2))} \Corr(\Ar(\vStk))
\end{align*}
factors through $\Corr'(\Ar(\vStk))$: This follows from Beauville--Laszlo gluing, more precisely the fiber product square in the proof of \cref{rslt:properties-of-interpolating-locHck}. Composing the above functor with $F_3$ yields the desired functor $\Fin \to \Alg(\Corr(\vStk))$.
\end{proof}

\begin{rmk}
We point out that the category $\Corr(\vStk)$ is a $(3,1)$-category, so constructing the functor in \cref{rslt:functoriality-of-locHck} rigorously by hand is cumbersome.
\end{rmk}

We now understand the functoriality of the local Hecke stacks $\locHck_G^I$. In order to construct the spectral action, we also need to understand how the local Hecke stack acts on the stack $\Bun_G$ of $G$-bundles. Following the construction that is sketched in \cite[\S IX]{FS} we use the global Hecke stack $\glbHck_G^I \to \Div^I$ associated to a finite set $I$: It parametrizes two $G$-bundles $\mathcal E_1, \mathcal E_2$ on the Fargues-Fontaine curve and a family $(D_i)_{i\in I}$ of Cartier divisors of degree $1$, together with an isomorphism of $\mathcal E_1$ and $\mathcal E_2$ away from $\bigcup_i D_i$. This stack is a convolution stack and hence comes equipped with a monoidal structure in $\Corr(\vStk)$. We now have the following result:

\begin{lem} \label{rslt:functoriality-of-locHck-to-BunG}
Fix a an algebraic group $G$ over $F$. Then there is a natural transformation
\begin{align*}
    \locHck_G^\bullet \to \Conv_{\Div^\bullet}(\Bun_G \times \Div^\bullet) = \Bun_G \times \Bun_G \times \Div^\bullet
\end{align*}
of functors $\Fin \to \Alg(\Corr(\vStk))$ which for fixed $I \in \Fin$ is realized by the correspondence
\begin{equation*}\begin{tikzcd}
    & \glbHck_G^I \arrow[dl] \arrow[dr]\\
    \locHck_G^I \arrow[rr,dashed] && \Bun_G \times \Bun_G \times \Div^I
\end{tikzcd}\end{equation*}
Here $\glbHck_G^I$ denotes the global Hecke stack of triples $(\mathcal E_1, \mathcal E_2, (D_i)_{i\in I})$ introduced above, the forward map is given by the projection to the three components of the triple and the backward map is given by restricting $\mathcal E_1$ and $\mathcal E_2$ to the completion of $\bigcup_i D_i$.
\end{lem}
\begin{proof}
Consider the category $\Fin \times [1]$. It consists of pairs $(I, a)$ where $a \in \{ 0, 1 \}$. For simplicity we will denote such an object by $I_a$. The desired natural transformation is equivalently a functor
\begin{align*}
    F\colon \Fin \times [1] \to \Alg(\Corr(\vStk)).
\end{align*}
We construct this by upgrading the construction in \cref{rslt:functoriality-of-locHck}. We first define the functor
\begin{align*}
    F_1\colon \Fin \times [1] \to \Corr(\Ar(\Tw(\Fin \times [1]))), \qquad I_a \mapsto \big[ [\emptyset_a \to I_a] \to [I_a \to I_a] \big]
\end{align*}
in a similar fashion as in \cref{rslt:functoriality-of-locHck}. Next we define the functor
\begin{align*}
    F_2\colon \Tw(\Fin \times [1]) \to \vStk
\end{align*}
as follows. Suppose we are given an object $f = [I_a \to J_{a'}]$ in $\Tw(\Fin \times [1])$. Then $F_2(f)$ is the following small v-stack over $\Div^J$:
\begin{itemize}
    \item If $a = a' = 0$ then we let $F_2(f)$ be defined as in \cref{rslt:functoriality-of-locHck}.

    \item If $a = 0$ and $a' = 1$ then for every $S \to \Div^J$ the value of $F_2(f)$ on $S$ is the groupoid of $G$-bundles on $\FF_S \setminus \bigcup_{i\in I} D_{f(i)}$, where $(D_j)_{j\in J}$ denotes the family of Cartier divisors on $\FF_S$ parametrized by $S \to \Div^J$.

    \item If $a = a' = 1$ then for every $S \to \Div^J$ the value of $F'_2(f)$ on $S$ is the groupoid of $G$-bundles on $\FF_S \setminus \bigcup_{i\in I} \FF_S$. In other words, $F_2([\emptyset_1 \to J_1]) = \Bun_G \times \Div^J$ and $F_2([I_1 \to J_1]) = \Div^J$ if $I \ne \emptyset$.
\end{itemize}
The intuition behind the assignment is that $F_2([I_a \to J_{a'}])$ parametrizes $G$-bundles on the completion of the union of the divisors parametrized by $J$ minus the union of the divisors parametrized by $I$, where $a$ and $a'$ determine whether by \enquote{divisor} we mean the given Cartier divisors or instead take copies of the whole curve. One can check that the above description indeed defines a functor $F_2$ (this can be checked by hand because $\Tw(\Fin \times [1])$ is an ordinary category and $\vStk$ is a $(2,1)$-category). A more conceptual definition of this functor would require us to talk about $G$-bundles on the completion of closed subspaces of $\FF_S$, which we avoid here.

We now observe that the composed functor
\begin{align*}
    \Fin \times [1] \xto{F_1} \Corr(\Ar(\Tw(\Fin \times [1]))) \xto{\Corr(Ar(F_2))} \Corr(\Ar(\vStk))
\end{align*}
factors over $\Corr'(\Ar(\vStk))$ (as defined in \cref{rslt:functoriality-of-convolution-stacks}): Similar to the proof of \cref{rslt:functoriality-of-locHck} this reduces to a Beauville--Laszlo gluing to compare the local and global Hecke stacks. We finally compose the above functor with the functor $F_3\colon \Corr'(\Ar(\vStk)) \to \Alg(\Corr(\vStk))$ from \cref{rslt:functoriality-of-locHck} in order to arrive at the desired functor $\Fin \times [1] \to \Alg(\Corr(\vStk))$.
\end{proof}

In the next step we upgrade \cref{rslt:functoriality-of-locHck-to-BunG} slightly: We also need to capture the fact that the natural correspondences $\locHck^I_G \to \Bun_G \times \Bun_G \times \Div^I$ live over $\Div^I$, i.e.\ this is a map in $\Alg(\Corr(\vStk_{/\Div^I}))$. Since this category depends on $I$, the precise meaning of the claimed functoriality gets a bit more complicated:

\begin{prop} \label{rslt:upgraded-functoriality-of-locHck-to-BunG}
Fix a reductive group $G$ over $F$ and let $\cat E \to \Fin$ be the cartesian fibration corresponding to the functor $\Fin^\op \to \Cat$, $I \mapsto \Alg(\Corr(\vStk))_{\Div^I/}$, where we view $\Div^I$ as a (symmetric) monoid in $\Corr(\vStk)$ via the diagonal monoidal structure. Then the natural transformation
\begin{align*}
    \locHck_G^\bullet \to \Conv_{\Div^\bullet}(\Bun_G \times \Div^\bullet)
\end{align*}
lifts to a natural transformation of functors $\Fin \to \cat E$ over $\Fin$.
\end{prop}
\begin{proof}
Recall that there is a functor $\vStk^\op \to \Alg(\Corr(\vStk))$ that sends each v-stack to itself equipped with the diagonal (symmetric) monoidal structure -- this follows from \cite[Proposition~2.3.7(i)]{heyer-mann-6ff} and \cite[Corollary~2.4.3.10]{lurie-higher-algebra}. Moreover, by \cite[Example~A.2.13]{heyer-mann-6ff} the cartesian unstraightening of the functor $\Alg(\Corr(\vStk))^\op \to \Cat$, $X \mapsto \Alg(\Corr(\vStk))_{X/}$ is given by $\Ar(\Alg(\Corr(\vStk)))$. Altogether we see that $\cat E$ sits in the following fiber product square:
\begin{equation*}\begin{tikzcd}
    \cat E \arrow[r] \arrow[d] & \Ar(\Alg(\Corr(\vStk))) \arrow[d] \\
    \Fin \arrow[r] & \Alg(\Corr(\vStk))
\end{tikzcd}\end{equation*}
Here the right vertical arrow is the source map and the bottom horizontal arrow is the composition $\Fin \to \vStk^\op \to \Alg(\Corr(\vStk))$, where the first functor is $I \mapsto \Div^I$ and the second functor is the one discussed above. Thus in order to construct the desired functor $\Fin \times [1] \to \cat E$, we need to construct the functor $\Fin \times [1] \to \Ar(\Alg(\Corr(\vStk)))$ that sends $I_a$ to the map $\Div^I \to F(I)$ for $F$ the functor in the proof of \cref{rslt:functoriality-of-locHck-to-BunG}.

We now go through the construction in \cref{rslt:functoriality-of-locHck-to-BunG} and upgrade it. We first observe that the functor $F_2\colon \Tw(\Fin \times [1]) \to \vStk$ upgrades to the functor
\begin{align*}
    F'_2\colon \Tw(\Fin \times [1]) \to \Ar(\vStk), \qquad [I_a \xto{f} J_{a'}] \mapsto [F_2(f) \to \Div^J].
\end{align*}
Now consider the functor
\begin{align*}
    F_{12}\colon \Fin \times [1] \to \Corr(\Ar(\Ar(\vStk))) \to \Ar(\Corr(\Ar(\vStk))),
\end{align*}
where the first functor is the composition of $F_1$ and $\Corr(\Ar(F'_2))$ and the second functor is defined as follows. Note first that for every category $\cat C$ with finite limits there is a natural functor $\alpha_1\colon \Corr(\Ar(\cat C)) \to \Ar(\Corr(\cat C))$ induced by the composition
\begin{align*}
    \Corr(\Ar(\cat C)) \times [1] \to \Corr(\Ar(\cat C)) \times \Corr([1]) = \Corr(\Ar(\cat C) \times [1]) \xto{\Corr(\ev)} \Corr(\cat C).
\end{align*}
Explicitly, $\alpha_1$ sends a map $f\colon Y \to X$ in $\cat C$ to the correspondence $Y = Y \to X$. There is a similar functor $\alpha_2\colon \Corr(\Ar(\cat C)) \to \Ar(\Corr(\cat C))$ induced by the composition
\begin{align*}
    &\Corr(\Ar(\cat C)) \times [1] \isoto \Corr(\Ar(\cat C)) \times [1]^\op \to \Corr(\Ar(\cat C)) \times \Corr([1]) =\\&\qquad \Corr(\Ar(\cat C) \times [1]) \xto{\Corr(\ev)} \Corr(\cat C).
\end{align*}
Explicitly, $\alpha_2$ sends $f\colon Y \to X$ to $X \from Y = Y$. Using $\alpha_1$, $\alpha_2$ and composition in $\Corr(\Ar(\vStk))$ we can now build the functor
\begin{align*}
    \Corr(\Ar(\Ar(\vStk))) \to \Ar(\Corr(\Ar(\vStk)))
\end{align*}
that sends
\begin{align*}
    \big[[Y' \to X'] \to [Y \to X]\big] \mapsto \qquad
    \left[\begin{tikzcd}[column sep=tiny,row sep=small,ampersand replacement=\&]
        \& {[Y' \to Y']} \arrow[dl] \arrow[dr]\\
        {[Y \to X]} \&\& {[Y' \to X']}
    \end{tikzcd}\right]
\end{align*}
We have finished the construction of the functor $F_{12}$. Explicitly, it sends an object $I_a \in \Fin \times [1]$ to
\begin{equation*}\begin{tikzcd}
    & {[F_2(\emptyset \to I_a) \to F_2(\emptyset \to I_a)]} \arrow[dr] \arrow[dl]\\
    {[\Div^I \xto{\id} \Div^I]} && {[F_2(\emptyset \to I_a) \to F_2(I_a \to I_a)]}
\end{tikzcd}\end{equation*}
One checks that $F_{12}$ factors over $\Ar(\Corr'(\Ar(\vStk)))$. Finally, we compose $F_{12}$ with the functor
\begin{align*}
    \Ar(F_3)\colon \Ar(\Corr'(\Ar(\vStk))) \to \Ar(\Alg(\Corr(\vStk)))
\end{align*}
coming from \cref{rslt:functoriality-of-convolution-stacks} in order to arrive at the desired functor.
\end{proof}

\subsection{Bounded sheaves}

With \cref{rslt:upgraded-functoriality-of-locHck-to-BunG} at hand, we next introduce a suitable notion of bounded sheaves on $\locHck_G^I$ together with the required functoriality. Our idea is to replace the category of small v-stacks with a category of ind-systems $(X_i)_i$ of small v-stacks along closed immersions $X_i \injto X_{i'}$. We view such an ind-system as an additional structure on the small v-stack $X := \varinjlim_i X_i$ which determines the desired boundedness condition on $X$:

\begin{defn} \label{def:ind-vstacks-and-bounded-sheaves}
\begin{defenum}
    \item We denote by $\Ind'(\vStk)$ the full subcategory of $\Ind(\vStk)$ spanned by the filtered diagrams $(X_i)_i$ of small v-stacks such that all transition maps $X_i \to X_{i'}$ are closed immersions. We have natural functors
    \begin{align*}
        \vStk \injto \Ind'(\vStk) \to \vStk
    \end{align*}
    where the first functor takes $X$ to itself (as a constant diagram) and the second functor takes $(X_i)_i$ to $X := \varinjlim_i X_i$. We note that there is a natural map $(X_i)_i \injto X$ in $\Ind'(\vStk)$ and this map is a monomorphism. Moreover, each of the maps $X_i \to X$ is a closed immersion. In particular we see that maps $(Y_j)_j \to (X_i)_i$ are a subanima of maps $Y \to X$ for $Y = \varinjlim_j Y_j$.

    \item Let $(\vStk, E)$ be a geometric setup on $\vStk$ such that $E$ contains all closed immersions. Let $\D$ be a 3-functor formalism on $(\vStk, E)$. For every $(X_i)_i \in \Ind'(\vStk)$ with $X = \varinjlim_i X_i$ we denote by
    \begin{align*}
        \D^\bdd((X_i)_i) \subseteq \D(X)
    \end{align*}
    the full subcategory spanned by those $M \in \D(X)$ such that $M = f_{i!} M_i$ for some $i$, where $f_i\colon X_i \injto X$ denotes the associated closed immersion.

    \item With $(\vStk, E)$ as in (b), let $E'$ be the class of maps $(Y_j)_j \to (X_i)_i$ in $\Ind'(\vStk)$ such that each induced map $Y_j \to X_i$ lies in $E$. Moreover, let $B$ be the class of maps $(Y_j)_j \to (X_i)_i$ such that the induced map $f\colon Y = \varinjlim_j Y_j \to X = \varinjlim_i X_i$ has the following property: For every $i$, the preimage $f^{-1}(X_i)$ is contained in some $Y_j$.
\end{defenum}
\end{defn}

In \cref{sec:Corr-bounded} we prove a general and abstract result about the extension of 3-functor formalisms to bounded sheaves. It tells us the following:

\begin{prop} \label{rslt:3ff-on-bounded-sheaves-on-vstacks}
Let $(\vStk, E)$ be a geometric setup such that $E$ contains all closed immersions and let $\D$ be a 3-functor formalism on $(\vStk, E)$. With the notation from \cref{def:ind-vstacks-and-bounded-sheaves} the assignment $(X_i)_i \mapsto \D^\bdd((X_i)_i)$ upgrades to a 3-functor formalism
\begin{align*}
    \D^{\bdd}\colon \Corr(\Ind'(\vStk), B, E')^\tensor \to \Cat^\times
\end{align*}
in the sense of \cref{def:3ff-on-extended-geometric-setup}.
\end{prop}
\begin{proof}
This is a special case of \cref{rslt:extend-3ff-to-Ind-stacks}.
\end{proof}

We refer the reader to the discussion after \cref{def:3ff-on-extended-geometric-setup} for a more explicit description of the data that goes into the 3-functor formalism $\D^{\bdd}$. To summarize, $\D^{\bdd}$ defines a pullback $f^*$ for every $f \in B$, a lower-$!$ $f_!$ for every $f \in E$ and a non-unital symmetric monoidal structure $\tensor$ on every $\D^\bdd((X_i)_i)$.

We now put the results in this section together. In the following, we equip $\locHck_G^I$ with the structure of an ind-v-stack via its filtration by closed Schubert cells as in \cite[Definition~VI.2.6]{FS}.

\begin{lem} \label{rslt:functoriality-of-D-bdd-Hck}
Let $(\vStk, E)$ be a geometric setup such that $E$ contains all proper maps of locally finite $\dimtrg$ that are representable in spatial diamonds, and let $\D$ be a 3-functor formalism on $(\vStk, E)$. Let $G$ be a reductive group over $F$. Then there is a natural transformation
\begin{align*}
    \D^\bdd(\locHck_G^\bullet) \to \D(\Bun_G \times \Bun_G \times \Div^\bullet)
\end{align*}
of functors $\Fin \to \Mon$, where both $\D^\bdd(\locHck^I_G)$ and $\D(\Bun_G \times \Bun_G \times \Div^\bullet)$ are equipped with the convolution monoidal structures.
\end{lem}
\begin{proof}
We start with the functor $\locHck_G^\bullet\colon \Fin \to \Alg(\Corr(\vStk))$, $I \mapsto \locHck^I_G$ from \cref{rslt:functoriality-of-locHck}. For every finite set $I$ we equip $\locHck^I_G$ with the structure of an ind-stack in the usual way, i.e. with the filtration by finite unions of the closed substacks $(\locHck^I_G)_{\le \mu_\bullet}$ defined in \cite[Definition~VI.2.6]{FS}. By \cref{rslt:lift-functor-from-Corr-to-Corr-Ind-stacks} it is now a \emph{condition} that $\locHck_G$ lifts to a functor
\begin{align*}
    \locHck_G\colon \Fin \to \Alg(\Corr(\Ind'(\vStk), B, E'))
\end{align*}
that sends every finite set $I$ to $\locHck_G^I$ with the chosen ind-structure. Roughly this condition says that certain maps between local Hecke stacks (including the interpolating Hecke stacks $\locHck_G^{I \to J}$) factor through the chosen ind-structures and are built out of maps in $E$. This condition is readily verified using \cref{rslt:properties-of-interpolating-locHck} and the explicit description of the condition in \cref{rmk:lift-monoid-in-Corr-to-Corr-Ind} using basic properties of the local Hecke stacks as in \cite[Proposition VI.2.7]{FS}. Composing the above functor $\locHck_G$ with the 3-functor formalism $\D^{\bdd}$ produces the functor $\D^\bdd(\locHck_G^\bullet)$.

We now apply the above argument to the functor $\Fin \times [1] \to \Alg(\Corr(\vStk))$ from \cref{rslt:functoriality-of-locHck-to-BunG}. We use the ind-structures on $\locHck_G^I$ from above and we put the trivial ind-structure on $\Bun_G \times \Bun_G \times \Div^I$. It is again a condition to see that this works, which is easily verified using \cite[Proposition VI.2.7]{FS} and its obvious global analogue, and by following the above argument we arrive at the desired functor $\Fin \times [1] \to \Mon$.
\end{proof}

The next result upgrades the above construction using \cref{rslt:upgraded-functoriality-of-locHck-to-BunG}, i.e.\ captures the fact that the monoidal maps $\D^\bdd(\locHck^I_G) \to \D(\Bun_G \times \Bun_G \times \Div^I)$ are $\D(\Div^I)$-linear, compatibly in $I$.

\begin{prop} \label{rslt:upgraded-functoriality-of-D-bbd-Hck}
In the setup of \cref{rslt:functoriality-of-D-bdd-Hck}, let $\cat R \to \Fin$ denote the cartesian fibration associated to the functor $\Fin^\op \to \Cat$, $I \mapsto \Alg(\Cat)_{\D(\Div^I)/}$. Then the natural transformation
\begin{align*}
    \D^\bdd(\locHck_G^\bullet) \to \D(\Bun_G \times \Bun_G \times \Div^\bullet)
\end{align*}
upgrades to a natural transformation of functors $\Fin \to \cat R$ over $\Fin$.
\end{prop}
\begin{proof}
We argue as in the proof of \cref{rslt:upgraded-functoriality-of-locHck-to-BunG} to see that there is a fiber product square
\begin{equation*}\begin{tikzcd}
    \cat R \arrow[r] \arrow[d] & \Ar(\Alg(\Cat)) \arrow[d]\\
    \Fin \arrow[r] & \Alg(\Cat)
\end{tikzcd}\end{equation*}
where the right vertical arrow is the source map and the bottom horizontal arrow is the functor $I \mapsto \D(\Div^I)$ (with pullback transition maps). This reduces us to construct the functor $\Fin \times [1] \to \Ar(\Alg(\Cat))$. Now we apply the same construction as in \cref{rslt:functoriality-of-D-bdd-Hck} to the functor $\Fin \times [1] \to \Ar(\Alg(\Corr(\vStk)))$ from the proof of \cref{rslt:upgraded-functoriality-of-locHck-to-BunG}.
\end{proof}

\begin{rmk}
The notation \enquote{$\cat R$} aims to suggest that this category is related to Ran space, which occurs a lot in the literature on geometric Langlands but is subtle to define properly.
\end{rmk}

In order to apply \cref{rslt:upgraded-functoriality-of-D-bbd-Hck} to the construction of the spectral action, we need to perform a slight modification of the target category $\D(\Bun_G \times \Bun_G \times \Div^\bullet)$. This is possible if the 6-functor formalism $\D$ satisfies a certain form of categorical Künneth:

\begin{prop} \label{rslt:full-functoriality-of-D-bdd-Hck}
Let $\D$ be a presentable 6-functor formalism on a geometric setup $(\vStk, E)$ and let $G$ be a reductive group over $F$. Assume that the following condition is satisfied:
\begin{itemize}
    \item[($*$)] $E$ contains the map $\Bun_G \to *$ and all proper maps of locally finite $\dimtrg$ that are representable in spatial diamonds. Moreover, for every small v-stack $S$ the natural map
    \begin{align*}
        \D(\Bun_G) \tensor_{\D(*)} \D(S) \isoto \D(\Bun_G \times S)
    \end{align*}
    is an isomorphism.
\end{itemize}
Let $\cat R \to \Fin$ be the cartesian fibration associated to the functor $\Fin^\op \to \Cat$, $I \mapsto \Alg(\PrL)_{\D(\Div^I)/}$. Then there is a natural transformation
\begin{align*}
    \D^{\bdd}(\locHck_{G}^\bullet) \to \End_{\D(\ast)}(\D(\Bun_G)) \otimes_{\D(\ast)} \D(\Div^\bullet)
\end{align*}
of functors $\Fin \to \cat R$ over $\Fin$. Here $\End_{\D(*)}$ denotes the endomorphism algebra in $\Mod_{\D(*)}(\PrL)$.
\end{prop}
\begin{proof}
By \cref{rslt:upgraded-functoriality-of-D-bbd-Hck} it only remains to construct a natural isomorphism
\begin{align}
    \End(\D(\Bun_G)) \tensor_{\D(*)} \D(\Div^I) \isoto \D(\Bun_G \times \Bun_G \times \Div^I). \label{eq:End-BunG-vs-D-BunG-2-functorial}
\end{align}
This isomorphism needs to be functorial in the finite set $I$ and compatible with the natural algebra structures over $\D(\Div^I)$.

Let us start with the case $I = \emptyset$. We denote by $\cat K_{\D}$ the $\Mod_{\D(*)}(\PrL)$-enriched category of kernels of $\D$, i.e. the transfer of the self-enrichment on $\Corr(\vStk_E)$ along $\D\colon \Corr(\vStk_E) \to \Mod_{\D(*)}(\PrL)$. By definition of convolution stacks, $\Bun_G \times \Bun_G$ is the enriched endomorphism object of $\Bun_G$ in $\Corr(\vStk_E)$, hence $\D(\Bun_G \times \Bun_G)$ is the enriched endomorphism object of $\Bun_G$ in $\cat K_{\D}$. The representable functor $\Psi_{\D} := \Hom(*, -)\colon \cat K_{\D} \to \Mod_{\D(*)}(\PrL)$ is $\Mod_{\D(*)}(\PrL)$-enriched (by enriched Yoneda) and hence defines a natural map of algebras
\begin{align}
    \D(\Bun_G \times \Bun_G) \to \End_{\D(*)}(\D(\Bun_G)). \label{eq:End-BunG-vs-D-BunG-2}
\end{align}
It remains to verify that this map is an isomorphism, for which we can now ignore the algebra structures. Let $\cat K \subseteq \cat K_{\D}$ be the full subcategory spanned by $(\Bun_G)^n$ for $n \ge 0$. Then $\cat K$ is a symmetric monoidal category and by ($*$) the functor $\Psi_{\D}\colon \cat K \to \Mod_{\D(*)}(\PrL)$ is symmetric monoidal. Since every object $X \in \cat K$ is dualizable (even self-dual), the same follows for $\Psi_{\D}(X)$ and hence for $X, Y \in \cat K$ we have
\begin{align*}
    &\Hom_{\D(*)}(\Psi_{\D}(X), \Psi_{\D}(Y)) = \Hom_{\D(*)}(\D(*), \Psi_{\D}(X)^\vee \tensor_{\D(*)} \Psi_{\D}(Y)) =\\&\qquad= \Hom_{\D(*)}(\D(*), \Psi_{\D}(X^\vee \tensor Y)) = \Hom(*, X^\vee \tensor Y) = \Hom(X, Y).
\end{align*}
Thus $\Psi_{\D}\colon \cat K \to \Mod_{\D(*)}(\PrL)$ is fully faithful and in particular \cref{eq:End-BunG-vs-D-BunG-2} is an isomorphism. This settles the case $I = \emptyset$.

In order to handle the case of general $I$ and the required functoriality, we start with the following general observation: Given a functor $G\colon \Fin \to \cat R$ over $\Fin$, there are canonical maps $G(\emptyset) \tensor_{\D(*)} \D(\Div^I) \to G(I)$ for all finite sets $I$; if these maps are isomorphisms then they assemble to an isomorphism $G(\emptyset) \tensor_{\D(*)} \D(\Div^\bullet) \to G$ of functors $\Fin \to \cat R$ over $\Fin$. To see this, first note that $\cat R \to \Fin$ is not only a cartesian fibration, but also a cocartesian fibration: This amounts to the observation that for a map $I \to J$ in $\Fin$, the forgetful functor $\Alg(\PrL)_{\D(\Div^J)/} \to \Alg(\PrL)_{\D(\Div^I)/}$ admits a left adjoint, namely $- \tensor_{\D(\Div^I)} \D(\Div^J)$. The limit of the functor $\Fin \to \Cat$ corresponding to this cocartesian fibration is given by $\Alg(\PrL)_{\D(*)/}$, because $\Fin$ has the initial object $\emptyset$. This limit can also be described as the category of functors $\Fin \to \cat R$ over $\Fin$ that send all maps in $\Fin$ to cocartesian maps. This proves the above observation, as the condition on $G$ guarantees that $G$ sends maps in $\Fin$ to cocartesian maps.

The observation in the previous paragraph together with the isomorphism \cref{eq:End-BunG-vs-D-BunG-2} provide maps as in \cref{eq:End-BunG-vs-D-BunG-2-functorial} and reduce us to showing that these maps are isomorphisms. But this follows easily from condition ($*$).
\end{proof}

In order to prove compatibility of the spectral action with duality, we also need to understand how the family of functors in \cref{rslt:full-functoriality-of-D-bdd-Hck} interacts with the canonical auto-equivalence $\locHck^I_G \isoto \locHck^I_G$ obtained by swapping the two $G$-bundles (i.e. the two quotients by $L^+_I G$). This is handled by the following result. It makes use of some standard facts and notation about dualizable categories: Firstly, for a monoidal category $\cat E$ we denote by $\cat E^\rev$ the same category but with the opposite monoidal structure, i.e. $X \tensor_{\cat E^\rev} Y = Y \tensor_{\cat E} X$. Secondly, if $\cat V$ is a symmetric monoidal presentable category and $\cat C \in \Mod_{\cat V}(\PrL)$ is dualizable, then there is a natural isomorphism of monoidal categories
\begin{align*}
    \End_{\cat V}(\cat C)^\rev \isoto \End_{\cat V}(\cat C^\vee).
\end{align*}
Indeed, this is true for the $\Mod_{\cat V}(\PrL)$-enriched endomorphism objects by \cite[Corollary~C.2.12]{heyer-mann-6ff} and the definition of opposite enriched categories (see \cite[Definition~C.1.16]{heyer-mann-6ff}) and then follows for the underlying categories via the lax symmetric monoidal forgetful functor $\Mod_{\cat V}(\PrL) \to \Cat$.
         
\begin{prop} \label{rslt:duality-of-D-bdd-Hck}
In the setup of \cref{rslt:full-functoriality-of-D-bdd-Hck}, there is a commuting square
\begin{equation*}\begin{tikzcd}
    \D^\bdd(\locHck_G^\bullet)^\rev \arrow[r] \arrow[d,isom] & \End_{\D(*)}(\D(\Bun_G))^\rev \tensor_{\D(*)} \D(\Div^\bullet) \arrow[d,isom]\\
    \D^\bdd(\locHck_G^\bullet) \arrow[r] & \End_{\D(*)}(\D(\Bun_G)) \tensor_{\D(*)} \D(\Div^\bullet)
\end{tikzcd}\end{equation*}
of functors $\Fin \to \cat R$ over $\Fin$. The arrows in this diagram have the following description:
\begin{propenum}
    \item The lower horizontal map is the one from \cref{rslt:full-functoriality-of-D-bdd-Hck} and the upper horizontal map is the $\rev$-version of that.

    \item The left vertical map is induced by swapping the two $L^+_I G$-quotients in $\locHck^I_G$.

    \item \label{rslt:duality-of-D-bdd-Hck-map-on-End-Bun-G} The right vertical map is induced by the isomorphism
    \begin{align*}
        \End_{\D(*)}(\D(\Bun_G))^\rev \isoto \End_{\D(*)}(\D(\Bun_G)^\vee) \isoto \End_{\D(*)}(\D(\Bun_G)),
    \end{align*}
    where the first map is the one discussed above and the second one is induced by a canonical isomorphism $\D(\Bun_G)^\vee \isoto \D(\Bun_G)$.
\end{propenum}
\end{prop}
\begin{proof}
In the following we go through the construction of the natural transformation in \cref{rslt:full-functoriality-of-D-bdd-Hck} and show how the self-duality of convolution stacks from \cref{rslt:self-duality-of-convolution-stacks} induces a version of the desired square in each step.

The construction of the natural transformation $\locHck_G^\bullet \to \Conv_{\Div^\bullet}(\Bun_G \times \Div^\bullet)$ from \cref{rslt:functoriality-of-locHck-to-BunG} factors through $\Corr'(\Ar(\vStk))$, hence by \cref{rslt:self-duality-of-convolution-stacks} there is a diagram
\begin{equation*}\begin{tikzcd}
    (\locHck_G^\bullet)^\op \arrow[r] \arrow[d,isom] & \Conv_{\Div^\bullet}(\Bun_G \times \Div^\bullet)^\op \arrow[d,isom]\\
    \locHck_G^\bullet \arrow[r] & \Conv_{\Div^\bullet}(\Bun_G \times \Div^\bullet)
\end{tikzcd}\end{equation*}
of functors $\Fin \to \Alg(\Corr(\vStk))$, where the vertical maps are induced by swapping the two factors in the convolution stacks. This diagram upgrades to a diagram of functors $\Fin \to \cat E$ over $\Fin$ with $\cat E$ as in \cref{rslt:upgraded-functoriality-of-locHck-to-BunG}, i.e. the above stacks are naturally algebras over $\Div^I$ for $I \in \Fin$; indeed, this upgrade follows by observing that the whole construction in \cref{rslt:upgraded-functoriality-of-locHck-to-BunG} factors through $\Corr'(\Ar(\cat C))$ at the end. We also implicitly use that $(\Div^\bullet)^\op = \Div^\bullet$ (where $\Div^I$ is equipped with the \emph{symmetric} monoidal structure coming from the diagonal).

As in the proof of \cref{rslt:functoriality-of-D-bdd-Hck} we lift the above diagram from $\vStk$ to $\Ind'(\vStk)$, for which we only need to check that the swapping maps $\locHck_G^I \to \locHck_G^I$ preserve the $\Ind$-structure (and similarly for the interpolating and the global Hecke stacks). But this is clear, since the ind-structure induced by \cite[Definition VI.2.6]{FS} is stable under the swapping isomorphism as it simply replaces each bounded piece $\{ \mu_{\bullet} \}$ with $\{ -w_0(\mu_{\bullet}) \}$ and the resulting systems are cofinal in each other. Then the same upgrade as in \cref{rslt:upgraded-functoriality-of-D-bbd-Hck} leads us to a square
\begin{equation*}\begin{tikzcd}
    \D^\bdd(\locHck_G^\bullet)^\rev \arrow[r] \arrow[d,isom] & \D(\Bun_G \times \Bun_G \times \Div^\bullet)^\rev \arrow[d,isom]\\
    \D^\bdd(\locHck_G^\bullet) \arrow[r] & \D(\Bun_G \times \Bun_G \times \Div^\bullet)
\end{tikzcd}\end{equation*}
of functors $\Fin \to \cat R$ over $\Fin$. Using the isomorphism
\begin{align*}
    \End(\D(\Bun_G)) \tensor_{\D(*)} \D(\Div^\bullet) \isoto \D(\Bun_G \times \Bun_G \times \Div^\bullet)
\end{align*}
from the proof of \cref{rslt:full-functoriality-of-D-bdd-Hck} we arrive at the square of functors in the claim. It is clear that the descriptions in (i) and (ii) hold. It remains to prove (iii), so we need to show that there is a commuting diagram of isomorphisms
\begin{equation*}\begin{tikzcd}
    \End(\D(\Bun_G))^\rev \tensor_{\D(*)} \D(\Div^\bullet) \arrow[r,isom] \arrow[d,isom] & \D(\Bun_G \times \Bun_G \times \Div^\bullet)^\rev \arrow[d,isom]\\
    \End(\D(\Bun_G)) \tensor_{\D(*)} \D(\Div^\bullet) \arrow[r,isom] & \D(\Bun_G \times \Bun_G \times \Div^\bullet)
\end{tikzcd}\end{equation*}
where the right vertical isomorphism is induced by the swapping isomorphism on convolution stacks and the left vertical isomorphism is the one from (iii). As observed in the proof of \cref{rslt:full-functoriality-of-D-bdd-Hck}, all the functors $\Fin \to \cat R$ in the above square send all maps in $\Fin$ to cocartesian maps and are thus determined completely by their value on $\emptyset \in \Fin$. Thus, a commuting diagram as above is equivalent to a commuting diagram
\begin{equation*}\begin{tikzcd}
    \End(\D(\Bun_G))^\rev \arrow[r,isom] \arrow[d,isom] & \D(\Bun_G \times \Bun_G)^\rev \arrow[d,isom]\\
    \End(\D(\Bun_G)) \arrow[r,isom] & \D(\Bun_G \times \Bun_G)
\end{tikzcd}\end{equation*}
in $\Alg(\PrL)_{\D(*)/}$. But this square follows immediately from the construction in \cref{rslt:full-functoriality-of-D-bdd-Hck}. Namely, the right-hand isomorphism is induced by the self-duality of $\Bun_G$ in $\cat K \subseteq \cat K_{\D}$, with $\cat K$ as in the proof of \cref{rslt:full-functoriality-of-D-bdd-Hck}. Via the fully faithful and symmetric monoidal functor $\Psi_{\D}\colon \cat K \injto \Mod_{\D(*)}(\PrL)$, $X \mapsto \D(X)$ we see that $\D(\Bun_G)$ is also self-dual and that the above square (whose horizontal maps are induced by $\Psi_{\D}$) does indeed commute.
\end{proof}

\bibliographystyle{alpha}
\bibliography{bibliography}

\end{document}